%% file: PhD.tex
\pdfoutput=1
\documentclass[11pt,a4paper,oneside, openany]{book} 
\usepackage[utf8]{inputenc}
\usepackage[T1]{fontenc}

\input{PhD_file_packages}

\input{PhD_file_theorems}

\input{PhD_file_operators}

\title{Lax structures in 2-category theory}
\author{Miloslav Štěpán}
\date{\today}

\begin{document}

\pdfbookmark[0]{Front cover}{milobookmark1}

\includepdf[pages=-, offset=27 0]{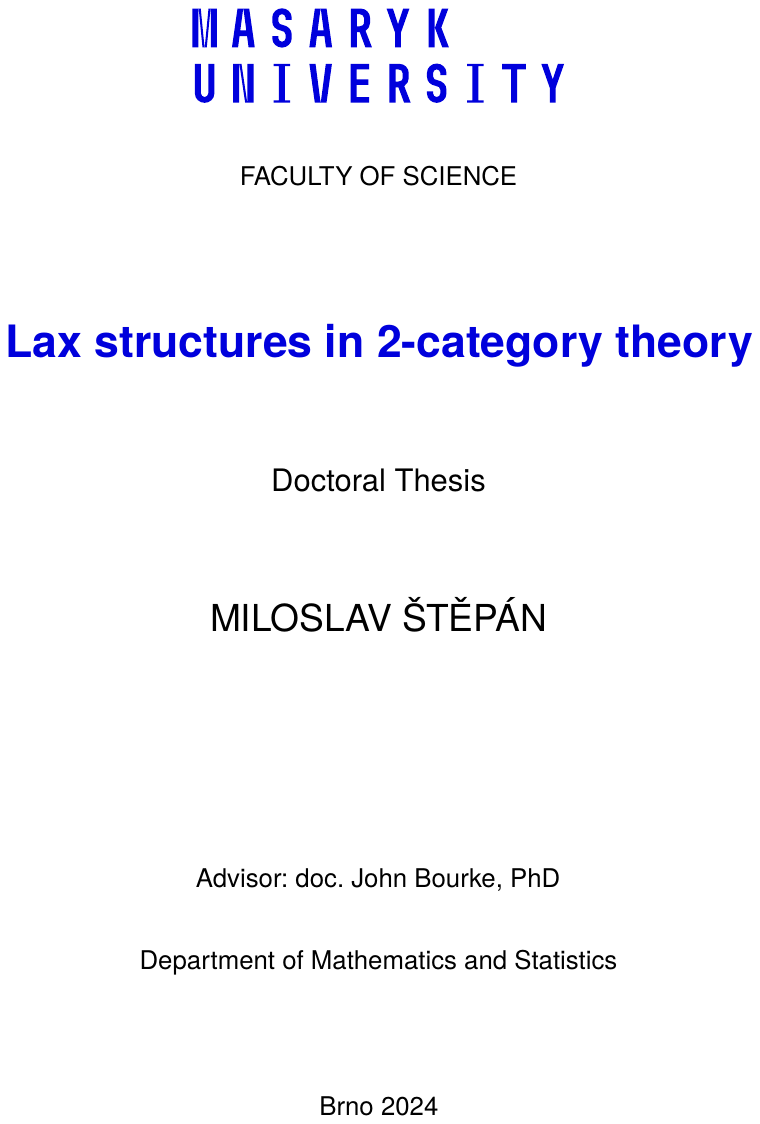}

\tableofcontents

\include{PhD_CH_1_INTRO}

\include{PhD_CH_2_TWOCATS}

\include{PhD_CH_3_FS_DBL}

\include{PhD_CH_4_COMPUTATION}

\include{PhD_CH_5_COLAX_ADJ_LAX_IDEMP_PSMON}

\include{PhD_CH_6_COKZFICATION}

\include{PhD_CH_7_Concluding}

\appendix

\include{PhD_CH_A_Further_remarks}

\include{PhD_CH_A_AUXILIARY_LEMMAS}

\bibliographystyle{plain} 
\bibliography{PhD} 
\addcontentsline{toc}{chapter}{Bibliography}

\end{document}

%% file: PhD_file_packages.tex
\usepackage{latexsym, url}
\usepackage{amsthm}
\usepackage{amsmath}
\usepackage{amsfonts}
\usepackage{amssymb}
\usepackage{adjustbox}
\usepackage{mathabx} 
\usepackage{nameref} 
\usepackage{hyperref}
\hypersetup{
    bookmarksnumbered=true
}
\usepackage{quiver}
\input xy
\xyoption{all}
\usepackage{url}
\usepackage{tikz-cd}

\usetikzlibrary{babel}

\usepackage{fancyhdr}
\usepackage{bm}
\usepackage{caption}
\usepackage{mathtools}
\usepackage{xcolor}
\usepackage{paralist}
\usepackage{opakuj}
\usepackage{graphicx}
\usepackage{booktabs}
\usepackage{graphicx}
\usepackage{tabu}
\usepackage{multirow}
\usepackage{relsize}
\usepackage{titletoc}
\usepackage[normalem]{ulem}

\usepackage{pdfpages}
\usepackage{zref-totpages}

%% file: PhD_file_theorems.tex
\theoremstyle{definition}
\newtheorem{theo}{Theorem}[section]

\newtheorem{lemma}[theo]{Lemma}
\newtheorem{ass}[theo]{Assumption}
\newtheorem{prop}[theo]{Proposition}
\newtheorem{defi}[theo]{Definition}
\newtheorem{cor}[theo]{Corollary}
\newtheorem{pozn}[theo]{Remark}

\newtheorem{con}[theo]{Construction}
\newtheorem{nota}[theo]{Notation}
\newtheorem{termi}[theo]{Terminology}
\newtheorem{pr}[theo]{Example}
\newtheorem{prc}[theo]{Class of examples}
\newtheorem{nonpr}[theo]{Non-example}
\newtheorem*{vari}{Variations}

\newtheorem*{warning}{Warning}

%% file: PhD_file_operators.tex
\DeclareMathOperator{\talg}{T\text{-}Alg}
\DeclareMathOperator{\talgl}{T\text{-}Alg_l}
\DeclareMathOperator{\salgl}{S\text{-}Alg_l}
\DeclareMathOperator{\talgc}{T\text{-}Alg_c}
\DeclareMathOperator{\talgs}{T\text{-}Alg_s}

\DeclareMathOperator{\tgph}{T\text{-}Graph}
\DeclareMathOperator{\tmult}{T\text{-}Multicat}

\DeclareMathOperator{\pstalg}{Ps\text{-}T\text{-}Alg}
\DeclareMathOperator{\pstalgl}{Ps\text{-}T\text{-}Alg_l}
\DeclareMathOperator{\psdalg}{Ps\text{-}D\text{-}Alg}

\DeclareMathOperator{\laxtalgl}{Lax\text{-}T\text{-}Alg_l}

\DeclareMathOperator{\colaxtalgl}{Colax\text{-}T\text{-}Alg_l}
\DeclareMathOperator{\colaxtalg}{Colax\text{-}T\text{-}Alg}

\DeclareMathOperator{\ncolaxtalgl}{NColax\text{-}T\text{-}Alg_l}
\DeclareMathOperator{\ncolaxtalg}{NColax\text{-}T\text{-}Alg}

\DeclareMathOperator{\qlcoalgs}{Q_l\text{-}CoAlg_s}

\DeclareMathOperator{\prof}{Prof}
\DeclareMathOperator{\PROF}{PROF}
\DeclareMathOperator{\cat}{Cat}
\DeclareMathOperator{\CAT}{CAT}
\DeclareMathOperator{\twocat}{2\text{-}Cat}
\DeclareMathOperator{\twoCAT}{2\text{-}CAT}
\DeclareMathOperator{\dbl}{Dbl}
\DeclareMathOperator{\rel}{Rel}

\DeclareMathOperator{\cocone}{Cocone(Cat)}
\DeclareMathOperator{\pair}{Pair}

\DeclareMathOperator{\ob}{ob }
\DeclareMathOperator{\mor}{mor }

\DeclareMathOperator{\comp}{comp}
\DeclareMathOperator{\ids}{ids}

\DeclareMathOperator{\gph}{Graph}
\DeclareMathOperator{\fc}{fc}
\DeclareMathOperator{\set}{Set}

\DeclareMathOperator{\parf}{Par}

\DeclareMathOperator{\spn}{Span}
\DeclareMathOperator{\comnd}{CoMnd}

\newcommand\cmndl[1]{\text{CoMnd}_l(#1)}

\newcommand\mndioo[2]{\text{Mnd}_{#1}(#2)}
\newcommand\cmndioo[2]{\text{CoMnd}_{#1}(#2)}

\newcommand{\ddashv}{\dashv\!\!\!\!\!\!\dashv}
\newcommand{\ttop}{\rotatebox{90}{$\ddashv$}}
\newcommand{\tbot}{\rotatebox{-90}{$\ddashv$}}

\DeclareMathOperator{\xxcnr}{Cnr}
\newcommand\cnr[1]{\xxcnr(#1)}

\newcommand\psd[2]{\text{Psd}[#1,#2]}
\newcommand\lax[2]{\text{Lax}[#1,#2]}
\newcommand\colax[2]{\text{CoLax}[#1,#2]}

\newcommand\Hom[2]{\text{Hom}[#1,#2]}
\newcommand\Homl[2]{\text{Hom}_l[#1,#2]}
\newcommand\Homc[2]{\text{Hom}_c[#1,#2]}
\newcommand\laxhomc[2]{\text{LaxHom}_c[#1,#2]}
\newcommand\laxhoml[2]{\text{LaxHom}_l[#1,#2]}

\DeclareMathOperator{\xxres}{Res}
\newcommand\res[1]{\xxres(#1)}

\newcommand\opres[1]{\text{OpRes}(#1)}

\DeclareMathOperator{\xxcoh}{Coh}
\newcommand\coh[3]{\xxcoh(#1,#2,#3)}

\newcommand\cod[1]{\text{CoDesc}(#1)}

\newcommand\sq[1]{\text{Sq}(#1)}
\newcommand\pbsq[1]{\text{PbSq}(#1)}

\DeclareMathOperator{\pscr}{Crossed}

\newcommand\sstackrel[2]{\stackrel{\mathmakebox[\widthof{#2}]{#1}}{#2}}

\newcommand\spanc{\text{Span}}

\newcommand\coconecat{\text{Cocone}(\cat)}

\newcommand{\sfs}{\mathcal {SFS}}
\newcommand{\ofs}{\mathcal {OFS}}
\newcommand{\coddiscr}{\text{CodDiscr}}
\newcommand{\factdbl}{\text{FactDbl}}

\newcommand\mbbd{\mathbb{D}}

\newcommand\mbbc{\mathbb{C}}

\newcommand\mbba{\mathbb{A}}
\newcommand\mbbl{\mathbb{L}}

\newcommand\ca{\mathcal {A}}
\newcommand\cb{\mathcal {B}}
\newcommand\cc{\mathcal {C}}
\newcommand\cd{\mathcal {D}}
\newcommand\cf{\mathcal {F}}

\newcommand\ce{\mathcal {E}}
\newcommand\cj{\mathcal {J}}
\newcommand\ck{\mathcal {K}}
\newcommand\cl{\mathcal {L}}
\newcommand\cm{\mathcal {M}}
\newcommand\cn{\mathcal {N}}
\newcommand\co{\mathcal {O}}

\newcommand\cs{\mathcal {S}}

\newcommand\cp{\mathcal {P}}
\newcommand\cw{\mathcal {W}}
\newcommand\cx{\mathcal {X}}
\newcommand\cy{\mathcal {Y}}
\newcommand\cv{\mathcal {V}}
\newcommand\cz{\mathcal {Z}}

%% file: PhD_CH_1_INTRO.tex
\chapter{Introduction}\label{CH_intro}

This thesis aims to contribute to the study of low-dimensional higher category theory -- to the study of 2-categories, double categories, pseudomonads and the like. The goal is to understand ``lax” analogues of various phenomena surrounding ``pseudo” structures -- the word ``pseudo” here meaning ``having invertible 2-cells”, with ``lax” meaning ``having not necessarily invertible 2-cells”.

For instance, monoidal categories (with the tensor $\otimes$ being associative only up to an isomorphism) and monoidal functors (the tensor being preserved up to an isomorphism) are a typical example of pseudo structures -- they are also the kind of monoidal structures one typically comes across. A \textbf{lax} monoidal functor $F: (\ca, \otimes) \to (\cb,\odot)$ would be a functor $F: \ca \to \cb$ for which there is only a comparison morphism between the tensors, without assumptions on its invertibility:
\[
\overline{F}_{a,b}: Fa \odot Fb \to F(a \otimes b).
\]
These may feel more rare to come by, but if for instance the tensors $\otimes, \odot$ are given by the coproducts in $\ca, \cb$, \textbf{any} functor $F: \ca \to \cb$ inherits a canonical lax monoidal structure. In other instances, the lax notion is the default one -- one does not consider monads $(t,\mu,\eta)$ in a 2-category for which both 2-cells $\mu$, $\eta$ are invertible. The natural kinds of morphism of monads -- \textit{monad functors} -- likewise do not have the accompanying 2-cell invertible. Going one level higher, it is the \textbf{lax-idempotent} pseudomonads that are relevant, with their pseudo variant being relatively rare.

The motivating questions the author had in mind when writing the thesis were the following:
\begin{itemize}
\item Lack \cite{codobjcoh} gave a general framework under which pseudo algebras are canonically equivalent to strict algebras -- this is called the \textit{coherence theorem for pseudo algebras}, an example of which is the well-known coherence for monoidal categories.

\textbf{Question 1}: What is the coherence theorem for colax algebras? Under what conditions can a colax algebra $A$ be turned into a strict one $A'$ and what is the relationship between $A, A'$? What are the examples of this phenomenon?

\item Blackwell, Kelly, Power \cite{twodim} study biadjunctions related to pseudo morphisms and use these to prove the bicocompleteness of strict algebras and pseudo morphisms for a 2-monad. They also introduce and study an important class of algebras called the \textit{flexible} algebras, an example of which are the flexible limits.

\textbf{Question 2}: What do we obtain if we work with lax morphisms instead? What is the lax analogue of flexible algebras and what properties do they enjoy?
\end{itemize}

\medskip

\noindent \textbf{Question 1} led the author to study codescent objects and in particular how to compute them in certain desirable situations. The conditions under which a colax algebra $A$ can be turned into a strict one $A'$ such that there is an adjunction between $A$ and $A'$ (perhaps the best lax analogue one can hope for) have been known in the literature \cite{moritactxs} -- when they are satisfied for a 2-monad $T$, we will say that \textit{coherence for colax algebras} holds for $T$. Apart from some general classes of 2-monads (those preserving all colimits), the author was not aware of any other examples for which coherence for colax algebras holds. What will be done in this thesis is to slightly refine the result from \cite{moritactxs} (Section \ref{SUBS_lax_coh_THMS}) and then exhibit a class of 2-monads for which this is possible (Section \ref{SEC_application_lax_coh}). This class is wide enough to include colax monoidal categories, colax functors $\cj \to \cat$ (for $\cj$ a 1-category), lax functors between 2-categories, or lax double functors between double categories.

At the same time, the observation that the codescent objects in \cite{internalalg} often admit canonical factorization systems led the author to study the correspondence between double categories and factorization systems, where a double category is turned into a factorization system by regarding it as a diagram in $\cat$ and then computing its codescent object. This work resulted in Chapter \ref{CH_FS_DBLCATS}.

\smallskip

\noindent \textbf{Question 2} led the author to the study of general lax-idempotent pseudomonads on 2-categories and colax adjunctions associated to their Kleisli 2-categories (Chapter \ref{CH_colax_adjs_lax_idemp_psmons}). We provide the answers to this question in Section \ref{SEKCE_applications_to_twodim} and Section \ref{SUBS_multicats_codescent}.

\bigskip
\noindent \textbf{The structure of the thesis is as follows}:

\begin{itemize}

\item Chapter \ref{CH_prelims} primarily serves to gather all the known results in 2-category theory that will be relevant for the later chapters. This includes providing exposition on 2-categorical limits, colax adjunctions, pseudomonads, internal category theory and codescent objects.

\item In Chapter \ref{CH_FS_DBLCATS} we show that orthogonal factorization systems are equivalent to a certain class of double categories, and describe these double categories in elementary terms (Theorem \ref{thm_equiv_ofs_factdbl}). We also do the same for strict factorization systems (Theorem \ref{thmekvivalence_sfs_coddiscr}).

\item In Chapter \ref{CH_COMPUTATION} we use codescent objects to show that both equivalences from the above chapter are in fact restrictions of a common adjunction (Theorem \ref{THM_fundamentaladj}). In later sections, we show that for a certain class of 2-monads, every colax algebra $A$ can be turned into a strict algebra $A'$ and moreover there is a canonical adjunction between $A$ and $A'$ (Theorem \ref{THM_coherence_for_colax_algs_CatT}). We also provide an explicit formula for such an $A'$ and demonstrate the constructions on numerous examples, including colax monoidal categories (Example \ref{PR_freestrictmoncat2monad_coherence}). 

\item In Chapter \ref{CH_colax_adjs_lax_idemp_psmons} we prove a generalization of a theorem of Bunge and Gray about forming colax adjunctions out of relative Kan extensions (Theorem \ref{THM_Bunge}) and apply it to the study of the Kleisli 2-category for a lax-idempotent pseudomonad (Theorem \ref{THM_BIG_lax_adj_thm}). We then show how the duals of these results provide lax analogues of results in two-dimensional monad theory (Theorem \ref{THM_big_colax_adj_THM}).

\item Chapter \ref{CH_coKZfication} consists of two independent sections:
\begin{itemize}
\item Section \ref{SEC_coKZfication} describes a process that turns a 2-monad (that satisfies certain properties) into a colax-idempotent one. This is then applied to the lax morphism classifier 2-comonad that is relevant for two-dimensional monad theory (Theorem \ref{THM_univ_prop_of_Ql}).
\item Section \ref{SEC_lax_pie_algebras} studies a special class of $T$-algebras that we call \textit{lax-pie} ones, and describes how every such algebra has an underlying (generalized) multicategory on which it is free. Examples include strict monoidal categories (Example \ref{PR_freestrictmoncat2monad_laxpie}).
\end{itemize}

\item Chapter \ref{CH_concludingremarks} lists ideas for further development.
\item Appendix contains two chapters:
\begin{itemize}
\item Chapter \ref{CH_generalized_colax_coh} mentions a generalization of the coherence theorem for colax algebras.
\item Chapter \ref{CH_APPENDIX_AUXILIARY} is a technical section containing proofs for some of the results in Chapter \ref{CH_colax_adjs_lax_idemp_psmons}.
\end{itemize}
\end{itemize}

\section{Table of prerequisites}

The thesis is intended to be absolutely self-contained -- the reader is only assumed to be familiar with ordinary category theory (in the scope of \cite{categoriesworking_HIST}, \cite{context}) and with the basic definitions of 2-category theory (2-categories, 2-functors, 2-natural transformation as covered by \cite[\S 1]{review}) and the Yoneda lemma for 2-categories (see \cite[Exercise 7.10.3]{handbook}).

With minor exceptions, Chapters \ref{CH_FS_DBLCATS} - \ref{CH_coKZfication} are independent of each other and can be read in any order. The exception for instance being that Section \ref{SUBS_COMPUTATION_Examples} uses the construction from the previous chapter (Section \ref{subsekce_catofcnrs}) as an example. As already mentioned, Chapter \ref{CH_prelims} primarily serves to gather all the results and definitions that are well known\footnote{There are again exceptions to this rule, see Section \ref{SEC_whats_new} for the details.} in the literature, and Chapters \ref{CH_FS_DBLCATS} - \ref{CH_coKZfication} will freely use the results from this chapter. The reader is welcome to skip to any chapter or section they prefer -- the exact prerequisites required to understand each chapter are in this table:

\begin{flushleft}
\begin{tabu}{|l|l|l|}
\hline
No. & Section in Chapter \ref{CH_prelims} & Prerequisites \\
\hline
\ref{SEC_lax_idemp_psmonads} & \nameref{SEC_lax_idemp_psmonads} & \ref{SEC_psmonads}\\
\ref{SEC_weak_limits} & \nameref{SEC_weak_limits} & \ref{SEC_2cats}\\
\ref{SEC_codescent_objects} & \nameref{SEC_codescent_objects} & \ref{SEC_2cats}, \ref{SEC_psmonads}\\
\hline
No. & Chapter & Prerequisites \\
\hline
\ref{CH_FS_DBLCATS} & \nameref{CH_FS_DBLCATS}            & \ref{SEC_double_cats}       \\
\ref{CH_COMPUTATION} & \nameref{CH_COMPUTATION}          & \ref{SEC_2cats}, \ref{SEC_psmonads}, \ref{SEC_int_cat_theory}, \ref{SEC_double_cats}, \ref{SEC_codescent_objects} \\
\ref{CH_colax_adjs_lax_idemp_psmons} & \nameref{CH_colax_adjs_lax_idemp_psmons}                        & \ref{SEC_lax_transf}, \ref{SEC_lax_adjs}, \ref{SEC_psmonads}, \ref{SEC_lax_idemp_psmonads}, \ref{SEC_weak_limits} \\
\ref{SEC_coKZfication} & \nameref{SEC_coKZfication}        & \ref{SEC_psmonads}, \ref{SEC_lax_idemp_psmonads} \\
\ref{SEC_lax_pie_algebras} & \nameref{SEC_lax_pie_algebras}        & \ref{SEC_Tmulticats}, \ref{SEC_codescent_objects} \\
\hline
\end{tabu}
\end{flushleft}

\section{Notational conventions}

Let us fix our terminological conventions for the thesis:
\begin{itemize}
\item Generally, we use the prefix ``$2$-” to be synonymous with $\cat$-enriched -- so 2-categories, 2-functors, 2-natural transformations, 2-adjunctions\dots will be considered to be strict. For weaker concepts we use a different prefix, e.g. ``bi”.
\item In case $\ck$ is a 2-category and we say ``$\cx$ holds for the category $\ck$”, we mean ``$\cx$ holds for \textbf{the underlying category} of the 2-category $\ck$”.
\item A lax-idempotent pseudomonad has algebra structure left adjoint to the unit, and a lax-idempotent pseudocomonad has coalgebra structure left adjoint to the counit.
\end{itemize}

\noindent Let us also fix the notation regarding some basic concepts:

\begin{itemize}
\item for a 2-category $\ck$, by $\ck^{op}$ we mean the 2-category obtained from $\ck$ by reversing the directions of 1-cells, and by $\ck^{co}$ the 2-category obtained by reversing the directions of 2-cells,
\item $\boldsymbol{2}$: the free category on a graph with a single non-endomorphism arrow,
\item $\cat$: the 2-category of small categories, functors, natural transformations,
\item $\CAT$: the 2-category of locally small categories, functors, natural transformations,
\item $\PROF$: the virtual double category of locally small categories and profunctors,
\item $\prof$: the bicategory of locally small categories and \textit{small} profunctors (Example \ref{PR_presheafpsmonad}),
\item the symbol ``$\dashv$” denotes either an adjunction, 2-adjunction or a biadjunction depending on the context. The symbol ``$\ddashv$\hspace{0.1cm}” will be used for (co)lax adjunctions.
\end{itemize}

\noindent Let us also fix the notation regarding the common hom categories (the definitions for objects, 1-cells, 2-cells below may be found in Section \ref{SEC_lax_transf}):

\begin{nota}\label{NOTA_kl_psdkl_laxkl_colaxkl_hom}
For 2-categories $\ck, \cl$, denote:

\begin{flushleft}
{\tabulinesep=1.2mm
\begin{tabu} {|l|l|l|m{3.5cm}|}
\hline
\textbf{2-category}          & \textbf{functors}           & \textbf{transformations}        & \textbf{2-cells}                    \\ \hline
$[\ck, \cl]$     & \multirow{ 4}{*}{2-functors}                & 2-natural  &  \multirow{ 8}{*}{modifications}                          \\ \cline{1-1}\cline{3-3}
$\psd{\ck}{\cl}$      &                  & pseudonatural  &                            \\ \cline{1-1}\cline{3-3}
$\lax{\ck}{\cl}$     &                   & lax  &     \\ \cline{1-1}\cline{3-3}
$\colax{\ck}{\cl}$     &                   & colax  &     \\ \cline{1-3}
$\Hom{\ck}{\cl}$    & pseudofunctors & pseudonatural  &                            \\ \cline{1-3}
$\laxhoml{\ck}{\cl}$    & lax functors & lax  &                            \\ \cline{1-4}
\end{tabu}}
\end{flushleft}
\end{nota}

\medskip

\begin{defi}\label{DEFI_lali_rali_rari_lari}
Consider an adjunction (with counit $\epsilon$ and unit $\eta$) in a 2-category:
\[\begin{tikzcd}
	{(\epsilon, \eta):} & A && B
	\arrow[""{name=0, anchor=center, inner sep=0}, "f"', curve={height=24pt}, from=1-2, to=1-4]
	\arrow[""{name=1, anchor=center, inner sep=0}, "u"', curve={height=24pt}, from=1-4, to=1-2]
	\arrow["\dashv"{anchor=center, rotate=90}, draw=none, from=0, to=1]
\end{tikzcd}\]

\noindent It will be called:
\begin{itemize}
\item a \textit{reflection} if $\epsilon$ is invertible,
\item a \textit{coreflection} if $\eta$ is invertible,
\item an \textit{adjoint equivalence} if both $\epsilon$, $\eta$ are invertible.
\end{itemize}

\noindent The terminology regarding the 1-cells present will be as follows:
\begin{itemize}
\item in a reflection, $f$ will be called the \textit{reflector} and $u$ the \textit{reflection-inclusion},
\item in a coreflection, $u$ will be called the \textit{coreflector} and $f$ the \textit{coreflection-inclusion},
\item in an equivalence, $f$ will be said to be an \textit{equivalence inverse} to $u$.
\end{itemize}

\noindent Further pieces of terminology will be used for the cases where an invertible 2-cell is replaced by the identity. If, for instance, $\epsilon$ is the identity 2-cell, we refer to $f$ as a \textit{lali} (short for \textit{left adjoint-left inverse}) and to $u$ as a \textit{rari} (\textit{right adjoint-right inverse}).

\end{defi}

\section{Statement of originality}\label{SEC_whats_new}

I declare that the content of this thesis is the outcome of my own effort. To my best knowledge, the Chapters \ref{CH_FS_DBLCATS} - \ref{CH_coKZfication} contain no material previously published by anyone other than myself, except for places where I explicitly state otherwise. Chapter \ref{CH_prelims} contains the material I consider well known in the published literature, with the following exceptions that I am unaware of being published:

\begin{itemize}
\item Remark \ref{POZN_lax_adjs_not_unique_up_to_equiv}: lax adjoints are not unique up to an equivalence.

\item Proposition \ref{THMP_colax_algebras_for_freeinicobj2monad}: description of colax algebras for the free initial object 2-monad,
\item the proof of the Yoneda lemma for lax functors using colax adjunctions -- Theorem \ref{THM_Yonedalemma_lax_funs}.

\item Subsections \ref{SUBS_psmonads_doctrinal} - \ref{SUBS_psmonads_adjsvsalgs} contain the proof of the doctrinal adjunction for pseudomonads and define a 2-functor sending adjunctions to colax algebras (Theorem \ref{THM_2fun_phi_adj_colaxtalgl}).

\item Subsection \ref{SUBS_lax_idemp_props} uses this theorem to prove known things about lax-idem\-potent pseudomonads, but moreover establishes:
\begin{itemize}
\item Theorem \ref{THM_lax_idemp_psmonad_Phi}: a characterization of lax-idempotent pseudomonads in terms of their normal colax algebras,
\item Proposition \ref{THMP_normal_colax_algebra_for_lax_idemp_is_pseudo}: every normal colax algebra for a lax-idempotent pseudomo\-nad is a pseudo algebra.
\end{itemize}

\item Section \ref{SUBS_leftKan2monads} introduces a ``Left Kan” analogue of lax-idempotent 2-monads -- this is used in Chapter \ref{CH_colax_adjs_lax_idemp_psmons}.

\item Proposition \ref{THMP_2Catl_OFS}: There is an orthogonal factorization system on the category of 2-categories and lax functors -- this is used in Example \ref{prcodlax-colaxfunclass}.
\end{itemize}

\subsection{Publication information}
Some contents of the thesis have been published or sent for a publication:
\begin{itemize}
\item The content of Chapter \ref{CH_FS_DBLCATS}, except for Remarks \ref{POZN_bicrossed_is_bicartesian_NEW}, \ref{POZN_repr_defi_bicrossed} that have been added later, has been published in \cite{milo1}. The paper \cite{milo1} moreover contains results on ``codescent objects and double categories” (Propositions \ref{THMP_cohX_cong_pairX}, \ref{THMP_invariance_under_transposition} for $\ce = \set$, Proposition \ref{THMP_cnr_is_cod_obj}) and results on ``lax morphism classifiers” (special case of Theorem \ref{THM_catE_admits_codobj_of_domdiscr_codlax} for strict $T$-algebras and $\ce = \set$ and the special case of Remark \ref{POZN_colax_mor_classifier}).
\item The content of Chapter \ref{CH_colax_adjs_lax_idemp_psmons}, except for Remarks \ref{POZN_almost_a_biadj}, \ref{POZN_idemp_monad_char_algeber}, \ref{POZN_almost_a_biadj_big_l_adj_thm} that have been added later, has been put in a preprint \cite{milo2} that is currently awaiting a review.
\end{itemize}

\section{Information on the ArXiv version}

This is an edited version of the thesis (with original available at the Masaryk University repository \cite{milophd}), taking into account most of the suggestions of the two reviewers.

\newpage
\section{Acknowledgements}

I wish to thank my supervisor John Bourke for his careful guidance and all the feedback he's given me, always suggesting simpler and more conceptual alternatives to my overly technical proofs.

I thank Professor Emily Riehl for letting me visit her at Johns Hopkins University and inviting me to speak about my research at the university's category theory seminar in February 2022. Her book \cite{context} is also what sparked my deep interest in category theory.

I thank the two thesis reviewers for giving a very valuable feedback.

The work in Chapter \ref{CH_FS_DBLCATS} and Chapter \ref{CH_COMPUTATION} has been supported from Operational Programme Research, Development and Education ``Project Internal Grant Agency of Masa\-ryk University'' (No. $\text{CZ}.02.2.69 \symbol{92}$ $0.0\symbol{92} 0.0\symbol{92}19$ $\underline{\enskip} 073\symbol{92}0016943)$. I also acknowledge the support of Grant Agency of the Czech Republic under the grant 22-02964S.

My thanks goes to all my colleagues and friends Nathanael Arkor, Yuriy Dupyn, Jan Jurka, Martin Dzúrik, Kristóf Kanalas, Joanna Ko, Petr Liczman, Mária Šimková, Dominik Trnka for all the great conversations regarding academia, mathematics and non-mathematics.

I also thank my friends Barbora and Remi for their support when the times were bleak.

Last but not least, I thank my partner for discovering my heart in the university library, and I thank my family for believing in me and supporting me my whole life.

\vspace{2cm}
\begin{center}
Když radosti není dosti,\\
raduju se z maličkostí.\\
Představím si třeba kdyby\\
lidi žili jako ryby.\\
\smallskip
\textit{Hodina zpěvu (1988)}
\end{center}

%% file: PhD_CH_2_TWOCATS.tex
\chapter{Structures in 2-category theory}\label{CH_prelims}

\noindent \textbf{The chapter is organized as follows}:
\begin{itemize}
\item We first recall basic notions from 2-category theory including 2-categorical limits and comonads internal to a 2-category (Section \ref{SEC_2cats}), colax functors, colax natural transformations, modifications (Section \ref{SEC_lax_transf}), colax adjunctions (Section \ref{SEC_lax_adjs}).

\item Section \ref{SEC_psmonads} recalls pseudomonads and proves some of their basic properties, with primary focus being on their colax algebras.

\item Section \ref{THM_kz_monad_TFAE} recalls lax-idempotent pseudomonads and proves some of their basic properties.

\item Sections \ref{SEC_int_cat_theory}, \ref{SEC_double_cats}, \ref{SEC_Tmulticats} recall internal categories, their special case -- double categories, and their generalization -- $T$-multicategories.

\item Section \ref{SEC_weak_limits} recalls \textit{enriched weak limits} of \cite{enrichedweakness} and focuses on their special case in 2-category theory.

\item Section \ref{SEC_codescent_objects} defines a 2-categorical colimit called the \textit{codescent object} and describes its usage in two-dimensional monad theory -- in particular the focus is on turning colax algebras into strict ones and the coherence theorem for colax algebras.
\end{itemize}

\section{2-categories}\label{SEC_2cats}

\noindent \textbf{The section is organized as follows}:
\begin{itemize}
\item In subsections \ref{SUBS_limits} -- \ref{SUBS_colimits_in_Cat} we give an exposition on 2-categorical limits (in the scope of \cite{elemobs}) and also give explicit formulas for some simple colimits in $\cat$.

\item In Subsection \ref{SEKCE_monads_in_2cat} we recall some facts about comonads in a 2-category. 
\end{itemize}

\subsection{Limits}\label{SUBS_limits}

\begin{defi}\label{DEFI_2limit_of_2funs}
Let $\cj$ be a small 2-category, $\ck$ a 2-category and let $F: \cj \to \ck$ and $W: \cj \to \cat$ be two 2-functors. A 2-natural transformation $\lambda: W \Rightarrow \ck(L,F-)$ is called a \textit{$W$-weighted cone} for $F$ and $L$ will be called its \textit{apex}.

Such a cone is the \textit{limit of $F$ weighted by $W$} (in this case $L$ is called the \textit{limit object}) if it has the property that for every $A \in \ck$, the whiskering 2-natural transformation $\kappa : \ck(-,L) \Rightarrow [\cj, \cat](W, \ck(-, F?))$, whose component at $A \in \ck$ is the functor:
\begin{align*}
\kappa_A : \ck(A,L) &\to [\cj, \cat](W, \ck(A, F?)),\\
\kappa_A : (\theta: A \to L) &\mapsto (\ck(\theta,F?) \circ \lambda),
\end{align*}
is an isomorphism in $[\ck,\cat]$. The fact that $\kappa_A$ is bijective on objects will be referred to as a \textit{one-dimensional universal property} and the fact that it is fully faithful will be referred to as a \textit{two-dimensional universal property}.
\end{defi}

\begin{vari} We will encounter the following variations of the definition:
\begin{itemize}
\item given 2-functors $F: \cj \to \ck$, $W: \cj^{op} \to \cat$, the \textit{colimit of $F$ weighted by $W$} is defined to be the limit of $F^{op}:J^{op} \to \ck^{op}$ weighted by $W$.
\item (co)limits in which $W$ is the terminal weight (constant functor at the terminal category) will be called \textit{conical},
\item replacing $\lambda$ by a pseudonatural transformation and $[\cj, \cat]$ by $\psd{\cj}{\cat}$, we get the notion of a \textit{pseudo limit},
\item using lax natural $\lambda$ and $\lax{\cj}{\cat}$, this is the notion of a \textit{lax limit},
\item a weaker notion of a pseudolimit in which instead of requiring that $\kappa$ is an isomorphism we only require that it is an equivalence, will be called a \textit{bilimit}.
\end{itemize}
Additionally, both the 2-functors $F,W$ and the 2-category $\ck$ may be replaced by something weaker (bicategories, pseudofunctors\dots). We will not introduce additional terminology for these cases.
\end{vari}

\subsection{Examples of limits}

Here we list various examples of (co)limits that we will meet throughout the thesis. For the more important ones we describe both the one-dimensional and two-dimensional universal properties. In this section, $\ck$ will always be a general 2-category, $W: J \to \cat$ will be a 2-functor that will be the weight, and $F: J \to \ck$ will be the 2-functor whose limit we will be taking.

\begin{termi}
Our terminology regarding the names of limits will be as follows. If the 2-limit of $F$ weighted by $W$ is called a ``foo” in $\ck$, the corresponding colimit notion is called a ``co-foo” in $\ck^{op}$, and the bilimit version is called ``bi-foo”.
\end{termi}

\begin{prc}[Powers]
In case $J = *$, we can identify $F$ with an object $F \in \ob \ck$ and $W$ with a category $\cw \in \cat$.

\begin{itemize}
\item A cone with apex $L$ is a functor $\lambda: \cw \to \ck(L,F)$,
\item The one-dimensional universal property says that for any other cone:
\[
\mu: \cw \to \ck(Y,F),
\]
there is a unique 1-cell $\theta: Y \to L$ such that the following diagram commutes:
\[\begin{tikzcd}
	\cw && {\ck(L,F)} && {\ck(Y,F)}
	\arrow["\lambda", from=1-1, to=1-3]
	\arrow["\mu"', curve={height=24pt}, from=1-1, to=1-5]
	\arrow["{\theta^*}", from=1-3, to=1-5]
\end{tikzcd}\]
\item The two-dimensional universal property says that for any cone morphism\\$\rho: \mu' \to \mu''$, that is, a natural transformation $\rho: \mu' \Rightarrow \mu'' : \cw \to \ck(Y,F)$, there is a unique 2-cell $\overline{\theta}: \theta' \Rightarrow \theta''$ such that:
\[\begin{tikzcd}[column sep=scriptsize]
	\cw && {\ck(L,F)} && {\ck(Y,F)} & {=} & \cw && {\ck(Y,F)}
	\arrow["\lambda", from=1-1, to=1-3]
	\arrow[""{name=0, anchor=center, inner sep=0}, "{(\theta')^*}", curve={height=-24pt}, from=1-3, to=1-5]
	\arrow[""{name=1, anchor=center, inner sep=0}, "{(\theta'')^*}"', curve={height=24pt}, from=1-3, to=1-5]
	\arrow[""{name=2, anchor=center, inner sep=0}, "{\mu''}"', curve={height=24pt}, from=1-7, to=1-9]
	\arrow[""{name=3, anchor=center, inner sep=0}, "{\mu'}", curve={height=-24pt}, from=1-7, to=1-9]
	\arrow["{\overline{\theta}^*}", shorten <=6pt, shorten >=6pt, Rightarrow, from=0, to=1]
	\arrow["\rho", shorten <=6pt, shorten >=6pt, Rightarrow, from=3, to=2]
\end{tikzcd}\]

\end{itemize}

The 2-limit cone is called the \textit{power} (also the \textit{cotensor}) of $F$ by $\cw$ and is typically denoted by $F^\cw$.
\end{prc}

\begin{prop}\label{THMP_powers_with2_implies_1dim_impl_2dim}
Let $\ck$ be a 2-category that admits powers with $\boldsymbol{2}$. Then the two-dimensional universal property of any colimit follows from the one-dimensional one.
\end{prop}
\begin{proof}
See \cite[Page 306]{elemobs}.
\end{proof}

\begin{prc}[Inserters]
Consider the categories and functors:
\[\begin{tikzcd}
	{J:=} & \bullet && \bullet \\
	{W:=} & {*} && {\boldsymbol{2}} \\
	{F:=} & A && B
	\arrow[shift left=2, from=1-2, to=1-4]
	\arrow[shift right=2, from=1-2, to=1-4]
	\arrow["0", shift left=2, from=2-2, to=2-4]
	\arrow["1"', shift right=2, from=2-2, to=2-4]
	\arrow["f", shift left=2, from=3-2, to=3-4]
	\arrow["g"', shift right=2, from=3-2, to=3-4]
\end{tikzcd}\]
\begin{itemize}
\item A cone with apex $L$ is a pair $(i: L \to A, \xi)$ as pictured below:
\[\begin{tikzcd}
	& A \\
	L && B \\
	& A
	\arrow["f", from=1-2, to=2-3]
	\arrow["\xi", shorten <=6pt, shorten >=10pt, Rightarrow, from=1-2, to=3-2]
	\arrow["i", from=2-1, to=1-2]
	\arrow["i"', from=2-1, to=3-2]
	\arrow["g"', from=3-2, to=2-3]
\end{tikzcd}\]
\item The one-dimensional universal property says that given any other cone\\$(k: X \to A, \psi: fk \Rightarrow gk)$, there is a unique 1-cell $\theta: X \to L$ such that:
\begin{align*}
i \theta &= k,\\
i \xi &= \psi.
\end{align*}
\item The two-dimensional universal property says that given a morphism of cones:
\[
\rho: (k',\psi') \to (k'',\psi''),
\]
that is, a 2-cell $\rho: k' \Rightarrow k''$ satisfying:
\[\begin{tikzcd}[column sep=scriptsize]
	&& A &&&&&& A \\
	\\
	X &&&& B & {=} & X &&&& B \\
	\\
	&& A &&&&&& A
	\arrow["f", from=1-3, to=3-5]
	\arrow["{\psi'}", shorten <=15pt, shorten >=15pt, Rightarrow, from=1-3, to=5-3]
	\arrow["f", from=1-9, to=3-11]
	\arrow["{\psi'}", shorten <=15pt, shorten >=15pt, Rightarrow, from=1-9, to=5-9]
	\arrow[""{name=0, anchor=center, inner sep=0}, "{k''}"', from=3-1, to=1-3]
	\arrow[""{name=1, anchor=center, inner sep=0}, "{k'}", curve={height=-30pt}, from=3-1, to=1-3]
	\arrow["{k''}"', from=3-1, to=5-3]
	\arrow["{k'}", from=3-7, to=1-9]
	\arrow[""{name=2, anchor=center, inner sep=0}, "{k'}", from=3-7, to=5-9]
	\arrow[""{name=3, anchor=center, inner sep=0}, "{k''}"', curve={height=30pt}, from=3-7, to=5-9]
	\arrow["g"', from=5-3, to=3-5]
	\arrow["g"', from=5-9, to=3-11]
	\arrow["\rho", shorten <=5pt, shorten >=5pt, Rightarrow, from=1, to=0]
	\arrow["\rho", shorten <=5pt, shorten >=5pt, Rightarrow, from=2, to=3]
\end{tikzcd}\]
there is a unique 2-cell $\overline{\theta}: \theta' \Rightarrow \theta''$ such that:
\[
i \overline{\theta} = \rho.
\]
\end{itemize}
The limit cone is called the \textit{inserter}.
\end{prc}

\begin{prc}[Equifiers]
Let:
\[\begin{tikzcd}
	{J:=} & \bullet && \bullet \\
	\\
	{W:=} & {*} && {\boldsymbol{2}} \\
	\\
	{F:=} & A && B
	\arrow[""{name=0, anchor=center, inner sep=0}, curve={height=24pt}, from=1-2, to=1-4]
	\arrow[""{name=1, anchor=center, inner sep=0}, curve={height=-24pt}, from=1-2, to=1-4]
	\arrow[""{name=2, anchor=center, inner sep=0}, "1"{description}, curve={height=24pt}, from=3-2, to=3-4]
	\arrow[""{name=3, anchor=center, inner sep=0}, "0"{description}, curve={height=-24pt}, from=3-2, to=3-4]
	\arrow[""{name=4, anchor=center, inner sep=0}, "g"', curve={height=24pt}, from=5-2, to=5-4]
	\arrow[""{name=5, anchor=center, inner sep=0}, "f", curve={height=-24pt}, from=5-2, to=5-4]
	\arrow[shift right=3, shorten <=6pt, shorten >=6pt, Rightarrow, from=1, to=0]
	\arrow[shift left=3, shorten <=6pt, shorten >=6pt, Rightarrow, from=1, to=0]
	\arrow[shorten <=6pt, shorten >=6pt, Rightarrow, from=3, to=2]
	\arrow["\alpha"', shift right=3, shorten <=6pt, shorten >=6pt, Rightarrow, from=5, to=4]
	\arrow["\beta", shift left=3, shorten <=6pt, shorten >=6pt, Rightarrow, from=5, to=4]
\end{tikzcd}\]
\begin{itemize}
\item A cone is a 1-cell $i: L \to A$ such that $\alpha i = \beta i$.
\item The one-dimensional universal property says that given any other $k: X \to A$ such that $\alpha k = \beta k$, there exists a unique $\theta: X \to L$ such that:
\[\begin{tikzcd}
	L && A \\
	\\
	X
	\arrow["i", from=1-1, to=1-3]
	\arrow["{\exists ! \theta}", dashed, from=3-1, to=1-1]
	\arrow["k"', from=3-1, to=1-3]
\end{tikzcd}\]
\item The two-dimensional universal property that given any morphism $k' \to k''$ of cones, that is, a 2-cell $\rho: k' \Rightarrow k''$, there is a unique 2-cell $\overline{\theta}: \theta' \Rightarrow \theta''$ such that:
\[
i \overline{\theta}= \overline{k}.
\]
\end{itemize}
\noindent The limit is called the \textit{equifier} of $\alpha, \beta$.
\end{prc}

\begin{prc}[Comma objects]
Let:
\[\begin{tikzcd}
	{J:=} & \bullet && \bullet && \bullet \\
	{W:=} & {*} && {\boldsymbol{2}} && {*} \\
	{F:=} & A && C && B
	\arrow[from=1-2, to=1-4]
	\arrow[from=1-6, to=1-4]
	\arrow["1"', from=2-2, to=2-4]
	\arrow["0", from=2-6, to=2-4]
	\arrow["f"', from=3-2, to=3-4]
	\arrow["g", from=3-6, to=3-4]
\end{tikzcd}\]
\begin{itemize}
\item A cone with apex $L$ is a triple $(\lambda_0, \lambda_1, \overline{\lambda})$ as follows:
\[\begin{tikzcd}
	L && B \\
	\\
	A && C
	\arrow[""{name=0, anchor=center, inner sep=0}, "{\lambda_1}", from=1-1, to=1-3]
	\arrow["{\lambda_0}"', from=1-1, to=3-1]
	\arrow["g", from=1-3, to=3-3]
	\arrow[""{name=1, anchor=center, inner sep=0}, "f"', from=3-1, to=3-3]
	\arrow["{\overline{\lambda}}", shorten <=9pt, shorten >=9pt, Rightarrow, from=0, to=1]
\end{tikzcd}\]
\item The one-dimensional universal property says that given any other cone\\$(\mu_0, \mu_1, \overline{\mu})$, there is a unique 1-cell $\theta: Y \to L$ such that:
\[
\lambda_i \circ \theta = \mu_i, \text{ for }i\in \{0, 1 \},
\]
and:
\[\begin{tikzcd}
	Y \\
	& L && B && Y && B \\
	&&&& {=} \\
	& A && C && A && C
	\arrow["\theta", from=1-1, to=2-2]
	\arrow["{\mu_1}", curve={height=-24pt}, from=1-1, to=2-4]
	\arrow["{\mu_0}"', curve={height=24pt}, from=1-1, to=4-2]
	\arrow[""{name=0, anchor=center, inner sep=0}, "{\lambda_1}", from=2-2, to=2-4]
	\arrow["{\lambda_0}"', from=2-2, to=4-2]
	\arrow["g", from=2-4, to=4-4]
	\arrow[""{name=1, anchor=center, inner sep=0}, "{\mu_1}", from=2-6, to=2-8]
	\arrow["{\mu_0}"', from=2-6, to=4-6]
	\arrow["g", from=2-8, to=4-8]
	\arrow[""{name=2, anchor=center, inner sep=0}, "f"', from=4-2, to=4-4]
	\arrow[""{name=3, anchor=center, inner sep=0}, "f"', from=4-6, to=4-8]
	\arrow["{\overline{\lambda}}", shorten <=9pt, shorten >=9pt, Rightarrow, from=0, to=2]
	\arrow["{\overline{\mu}}", shorten <=9pt, shorten >=9pt, Rightarrow, from=1, to=3]
\end{tikzcd}\]

\item The two-dimensional universal property says that given any morphism of cocones $\rho: \mu' \to \mu''$, that is, a pair of 2-cells:
\[
\rho_0 : \mu_0' \Rightarrow \mu_0'', \hspace{0.5cm}\rho_1 : \mu_1' \Rightarrow \mu_1'',
\]
satisfying:
\[\begin{tikzcd}
	X && B &&& X && B \\
	&&& {=} \\
	A && C &&& A && C
	\arrow[""{name=0, anchor=center, inner sep=0}, "{\mu_1''}"'{pos=0.2}, from=1-1, to=1-3]
	\arrow[""{name=1, anchor=center, inner sep=0}, "{\mu'_1}", curve={height=-30pt}, from=1-1, to=1-3]
	\arrow["{\mu_0''}"', from=1-1, to=3-1]
	\arrow["g", from=1-3, to=3-3]
	\arrow[""{name=2, anchor=center, inner sep=0}, "{\mu_1'}", from=1-6, to=1-8]
	\arrow[""{name=3, anchor=center, inner sep=0}, "{\mu_0'}", from=1-6, to=3-6]
	\arrow[""{name=4, anchor=center, inner sep=0}, "{\mu''_0}"', curve={height=30pt}, from=1-6, to=3-6]
	\arrow["g", from=1-8, to=3-8]
	\arrow[""{name=5, anchor=center, inner sep=0}, "f"', from=3-1, to=3-3]
	\arrow[""{name=6, anchor=center, inner sep=0}, "f"', from=3-6, to=3-8]
	\arrow["{\rho_1}", shorten <=4pt, shorten >=4pt, Rightarrow, from=1, to=0]
	\arrow["{\overline{\mu''}}", shorten <=9pt, shorten >=9pt, Rightarrow, from=0, to=5]
	\arrow["{\rho_0}"', shorten <=6pt, shorten >=6pt, Rightarrow, from=3, to=4]
	\arrow["{\overline{\mu'}}", shorten <=9pt, shorten >=9pt, Rightarrow, from=2, to=6]
\end{tikzcd}\]
there is a unique 2-cell $\overline{\theta}: \theta' \Rightarrow \theta''$ satisfying:
\[
\lambda_i \overline{\theta} = \rho_i \text{ for } i \in \{ 0,1 \}.
\]
\end{itemize}
The limit cone is called the \textit{comma object} and in this case would be denoted by $f \downarrow g$. In case $\ck = \cat$, this is the usual comma category.
\end{prc}

\begin{prc}[Lax limit of an arrow]\label{PR_lax_limit_of_arrow}
Let:
\[\begin{tikzcd}
	{J:=} & \bullet && \bullet \\
	{W:=} & {*} \\
	{F:=} & A && B
	\arrow[from=1-2, to=1-4]
	\arrow["f", from=3-2, to=3-4]
\end{tikzcd}\]
\begin{itemize}
\item A lax cone with apex $L$ is a triple $\lambda := (\lambda_0, \lambda_1, \overline{\lambda})$ as follows:
\[\begin{tikzcd}
	& L \\
	\\
	A && B
	\arrow["{\lambda_0}"', from=1-2, to=3-1]
	\arrow[""{name=0, anchor=center, inner sep=0}, "{\lambda_1}", from=1-2, to=3-3]
	\arrow["f"', from=3-1, to=3-3]
	\arrow["{\overline{\lambda}}"', shorten <=8pt, shorten >=8pt, Rightarrow, from=3-1, to=0]
\end{tikzcd}\]
\item The one-dimensional universal property states that given any other lax cone $\mu$, there is a unique 1-cell $\theta: X \to L$ such that:
\[
\lambda_i \theta = \mu_i \text{ for } i \in \{ 0,1 \},
\]
and also that:
\[\begin{tikzcd}
	& X &&&& X \\
	\\
	& L \\
	&&& {=} \\
	A && B && A && B
	\arrow["{\exists ! \theta}"', dashed, from=1-2, to=3-2]
	\arrow["{\mu_0}"', curve={height=18pt}, from=1-2, to=5-1]
	\arrow["{\mu_1}", curve={height=-18pt}, from=1-2, to=5-3]
	\arrow["{\mu_0}"', from=1-6, to=5-5]
	\arrow[""{name=0, anchor=center, inner sep=0}, "{\mu_1}", from=1-6, to=5-7]
	\arrow["{\lambda_0}"', from=3-2, to=5-1]
	\arrow[""{name=1, anchor=center, inner sep=0}, "{\lambda_1}", from=3-2, to=5-3]
	\arrow["f"', from=5-1, to=5-3]
	\arrow["f"', from=5-5, to=5-7]
	\arrow["{\overline{\lambda}}"', shorten <=8pt, shorten >=8pt, Rightarrow, from=5-1, to=1]
	\arrow["{\overline{\mu}}"', shorten <=11pt, shorten >=11pt, Rightarrow, from=5-5, to=0]
\end{tikzcd}\]
\item the two-dimensional universal property says that given any modification $\rho: \mu' \to \mu'': \Delta X \Rightarrow F$ between lax cones, that is, a pair of 2-cells:
\[
\rho_0 : \mu_0' \Rightarrow \mu_0'', \hspace{0.5cm} \rho_1 : \mu_1' \Rightarrow \mu_1'',
\]
satisfying:
\[\begin{tikzcd}
	& X &&&& X \\
	&&& {=} \\
	A && B && A && B
	\arrow[""{name=0, anchor=center, inner sep=0}, "{\mu_0'}"', curve={height=24pt}, from=1-2, to=3-1]
	\arrow[""{name=1, anchor=center, inner sep=0}, "{\mu_0''}"{pos=0.3}, from=1-2, to=3-1]
	\arrow[""{name=2, anchor=center, inner sep=0}, "{\mu''_1}", curve={height=-12pt}, from=1-2, to=3-3]
	\arrow["{\mu_0'}"', from=1-6, to=3-5]
	\arrow[""{name=3, anchor=center, inner sep=0}, "{\mu''_1}", curve={height=-24pt}, from=1-6, to=3-7]
	\arrow[""{name=4, anchor=center, inner sep=0}, "{\mu_1'}"'{pos=0.2}, from=1-6, to=3-7]
	\arrow["f"', from=3-1, to=3-3]
	\arrow["f"', from=3-5, to=3-7]
	\arrow["{\rho_0}", shorten <=4pt, shorten >=4pt, Rightarrow, from=0, to=1]
	\arrow["{\rho_1}", shorten <=4pt, shorten >=4pt, Rightarrow, from=4, to=3]
	\arrow["{\overline{\mu''}}"', shorten <=10pt, shorten >=10pt, Rightarrow, from=3-1, to=2]
	\arrow["{\overline{\mu'}}"', shorten <=8pt, shorten >=8pt, Rightarrow, from=3-5, to=4]
\end{tikzcd}\]
there is a unique 2-cell $\overline{\theta}: \theta' \Rightarrow \theta''$ satisfying:
\[
\lambda_i \overline{\theta} = \rho_i \text{ for } i \in \{ 0,1 \}.
\]
\end{itemize}
The limit cone is called the \textit{lax limit of the arrow} $f: A \to B$. Note that it can equivalently be described as an ordinary 2-limit of $F$ weighted by the 2-functor $W: \boldsymbol{2} \to \cat$, that sends the walking arrow to:
\[\begin{tikzcd}
	{*} && {\boldsymbol{2}}
	\arrow["0", from=1-1, to=1-3]
\end{tikzcd}\]
Alternatively, the lax limit of $f$ is precisely the comma object $f \downarrow 1_B$.

Another variation is called the \textit{pseudo limit of an arrow} -- it has an invertible 2-cell $\overline{\lambda}$ that has the universal properties only with respect to other triples with an invertible 2-cell, and the \textit{oplax limit of an arrow} where the 2-cell goes in the other direction.
\end{prc}

\begin{pozn}\label{POZN_2limit_is_not_laxlimit}
In the previous example we were able to express a lax limit using an ordinary 2-limit. This can always be done -- given a 2-functor $W: \cj \to \cat$, \textbf{every} $W$-weighted pseudo limit and lax limit can be shown to be a 2-limit if one appropriately changes the weight. To elaborate: by \cite[Theorem 3.16]{twodim}, the inclusions into lax and pseudonatural transformations have left 2-adjoints:
\[\begin{tikzcd}
	{[\cj,\cat]} && {\lax{\cj}{\cat}} && {[\cj,\cat]} && {\psd{\cj}{\cat}}
	\arrow[""{name=0, anchor=center, inner sep=0}, "J"', curve={height=24pt}, hook, from=1-1, to=1-3]
	\arrow[""{name=1, anchor=center, inner sep=0}, "{(-)'}"', curve={height=24pt}, dashed, from=1-3, to=1-1]
	\arrow[""{name=2, anchor=center, inner sep=0}, "{J_p}"', curve={height=24pt}, hook, from=1-5, to=1-7]
	\arrow[""{name=3, anchor=center, inner sep=0}, "{(-)^\dagger}"', curve={height=24pt}, dashed, from=1-7, to=1-5]
	\arrow["\dashv"{anchor=center, rotate=-90}, draw=none, from=1, to=0]
	\arrow["\dashv"{anchor=center, rotate=-90}, draw=none, from=3, to=2]
\end{tikzcd}\]
So, for instance a $W$-weighted lax limit of $F$ is equivalently a $W'$-weighted 2-limit of $F$ since:
\[
\lax{\cj}{\cat}(W, \ck(-, F?)) \cong [\cj, \cat](W', \ck(-, F?)).
\]
The converse is not true -- not every 2-limit can be expressed as a lax limit. A~counterexample will be given in Example \ref{PR_reslan_laxcoh_explicit}.
\end{pozn}

\noindent As is well known, a monad $(A,t,\mu,\eta)$ in a 2-category $\ck$ is equivalently a lax functor\footnote{The reader not familiar with the concept may see Definition \ref{SEC_lax_transf}.} $F: * \to \ck$ from the terminal 2-category to $\ck$. For this next example consider the obvious extension of  Definition \ref{DEFI_2limit_of_2funs} to the case when $F$ is a lax functor:

\begin{prc}[EM-objects]\label{PR_EM_Kleisli_objects}
Let:
\[\begin{tikzcd}
	{J:=} & {*} && {W :=} & {*} && {F:= (A,t,\mu,\eta)}
\end{tikzcd}\]
\begin{itemize}
\item A lax cone with apex $L$ is a pair $(e: L \to A, \xi: te \Rightarrow e)$ satisfying:
\[\begin{tikzcd}
	&& A \\
	L &&&& {=} & {1_e} \\
	&& A
	\arrow[""{name=0, anchor=center, inner sep=0}, "t"'{pos=0.7}, from=1-3, to=3-3]
	\arrow[""{name=1, anchor=center, inner sep=0}, curve={height=-30pt}, Rightarrow, no head, from=1-3, to=3-3]
	\arrow["e", from=2-1, to=1-3]
	\arrow["e"', from=2-1, to=3-3]
	\arrow["\xi"', shorten <=11pt, shorten >=17pt, Rightarrow, from=0, to=2-1]
	\arrow["\eta"', shorten <=6pt, shorten >=6pt, Rightarrow, from=1, to=0]
\end{tikzcd}\]
and also satisfying:
\[\begin{tikzcd}
	&& A &&&& A \\
	\\
	L && A & {=} & L &&& A \\
	\\
	&& A &&&& A
	\arrow["t", from=1-3, to=3-3]
	\arrow["t", from=1-7, to=3-8]
	\arrow[""{name=0, anchor=center, inner sep=0}, "t"{pos=0.7}, curve={height=18pt}, from=1-7, to=5-7]
	\arrow[""{name=1, anchor=center, inner sep=0}, "e", from=3-1, to=1-3]
	\arrow[""{name=2, anchor=center, inner sep=0}, "e"{description}, from=3-1, to=3-3]
	\arrow[""{name=3, anchor=center, inner sep=0}, "e"', from=3-1, to=5-3]
	\arrow["t", from=3-3, to=5-3]
	\arrow["e", from=3-5, to=1-7]
	\arrow["e"', from=3-5, to=5-7]
	\arrow["t", from=3-8, to=5-7]
	\arrow["\xi"', shorten <=8pt, shorten >=8pt, Rightarrow, from=0, to=3-5]
	\arrow["\xi"', shorten <=4pt, shorten >=4pt, Rightarrow, from=1, to=2]
	\arrow["\xi"', shorten <=4pt, shorten >=4pt, Rightarrow, from=2, to=3]
	\arrow["\mu"', shorten <=8pt, shorten >=8pt, Rightarrow, from=3-8, to=0]
\end{tikzcd}\]
\item The one-dimensional universal property says that for any other cone\\$(g: X \to A, \psi: tg \Rightarrow g)$ there is a unique 1-cell $\theta: X \to L$ such that:
\begin{align*}
e \theta &= g,\\
\xi \theta &= \psi.
\end{align*}
\item The two-dimensional universal property says that given a morphism of cones\\ $\rho: (g',\psi') \to (g'',\psi'')$, that is, a 2-cell $\rho: g' \Rightarrow g''$ satisfying:
\[\begin{tikzcd}
	X && A && X && A \\
	&&& {=} \\
	&& A &&&& A
	\arrow["{g'}", from=1-1, to=1-3]
	\arrow[""{name=0, anchor=center, inner sep=0}, "{g''}"', curve={height=30pt}, from=1-1, to=3-3]
	\arrow[""{name=1, anchor=center, inner sep=0}, "{g'}"{description, pos=0.3}, from=1-1, to=3-3]
	\arrow["t", from=1-3, to=3-3]
	\arrow[""{name=2, anchor=center, inner sep=0}, "{g'}", from=1-5, to=1-7]
	\arrow[""{name=3, anchor=center, inner sep=0}, "{g''}"'{pos=0.7}, curve={height=30pt}, from=1-5, to=1-7]
	\arrow[""{name=4, anchor=center, inner sep=0}, "{g''}"', curve={height=30pt}, from=1-5, to=3-7]
	\arrow["t", from=1-7, to=3-7]
	\arrow["\rho"', shorten <=5pt, shorten >=5pt, Rightarrow, from=1, to=0]
	\arrow["{\psi'}"', shorten >=6pt, Rightarrow, from=1-3, to=1]
	\arrow["\rho", shorten <=4pt, shorten >=4pt, Rightarrow, from=2, to=3]
	\arrow["{\psi''}", shorten <=5pt, shorten >=5pt, Rightarrow, from=3, to=4]
\end{tikzcd}\]
there is a unique 2-cell $\overline{\theta}: \theta' \Rightarrow \theta''$ such that:
\[
e \overline{\theta} = \rho.
\]
\end{itemize}
The lax limit cone is called the \textit{Eilenberg-Moore object} (\textit{EM-object} for short) of the monad $(A,t,\mu,\eta)$. Any monad equivalently gives a lax functor $* \to \ck^{op}$, the lax colimit of which is called the \textit{Kleisli object} of the monad in $\ck$. EM-objects and Kleisli objects in $\cat$ are the usual Eilenberg-Moore and Kleisli categories associated to monads. The Kleisli object is also a special case of a 2-colimit called the \textit{codescent object} that will be studied in detail in Section \ref{SEC_codescent_objects}.
\end{prc}

There are many other important classes of examples - for instance, the \textit{inverter} of a 2-cell $\alpha : f \Rightarrow g$ is the universal 1-cell $i$ with the property that $\alpha i$ is invertible. The \textit{identifier} of two parallel 2-cells $\alpha, \beta: f \Rightarrow g$ is the universal 1-cell $j$ with $\alpha j = \beta j$.

\subsection{Some colimits in $\cat$}\label{SUBS_colimits_in_Cat}

The goal of this short sections is to give explicit formulas for coinserters and coequifiers in the 2-category $\cat$. The purpose of this is to motivate the formulas for codescent objects that will be given in Chapter \ref{CH_COMPUTATION}. Recall also that $\cat$ admits all small 2-limits and 2-colimits by \cite[Page 304]{elemobs} and \cite[Page 306]{elemobs}.

\begin{pr}[Coinserters in $\cat$]\label{PR_coinserter_in_Cat}
Consider the following diagram in $\cat$:
\[\begin{tikzcd}
	{X_1} && {X_0}
	\arrow["{d_1}"{description}, shift left=3, from=1-1, to=1-3]
	\arrow["{d_0}"{description}, shift right=3, from=1-1, to=1-3]
\end{tikzcd}\]
We will refer to the morphisms of $X_0$ as \textit{vertical morphisms} and to the objects of $X_1$ as \textit{horizontal morphisms} and denote $g \in \ob X_1$ as $g: a \to b$ if $d_1(g) = a, d_0(g) = b$.

The coinserter $(F: X_0 \to \cc, \xi: Fd_1 \Rightarrow Fd_0)$ is given as follows. The category $\cc$ has objects the objects of $X_0$, while a morphism $a \to b$ is an equivalence class of paths:
\[\begin{tikzcd}
	a & {a_1} & {a_2} & \dots & {a_n} & {b,}
	\arrow["{f_1}", from=1-1, to=1-2]
	\arrow["{f_2}", from=1-2, to=1-3]
	\arrow["{f_3}", from=1-3, to=1-4]
	\arrow["{f_n}", from=1-4, to=1-5]
	\arrow["{f_{n+1}}", from=1-5, to=1-6]
\end{tikzcd}\]
with each $f_i$ being either a vertical or a horizontal morphism. The equivalence relation on morphisms is then generated by the following:
\begin{alignat*}{2}
(f_1, f_2) &\sim (f_2 \circ f_1) &\phantom{....}&\text{if both }f_1, f_2\text{ are vertical},\\
(1_a) &\sim ()_a                 & &\text{if }1_a\text{ is the identity in }X_0,\\
(g,v) &\sim (u,h)                & &\text{if there is }(\alpha: g \to h) \in X_1 \\
&                                & &\text{such that } d_1(\alpha) = u, d_0(\alpha)= v.
\end{alignat*}
\end{pr}

\begin{pr}[Coequifiers in $\cat$]
Consider a pair of parallel natural transformations:
\[\begin{tikzcd}
	\ca && \cb
	\arrow[""{name=0, anchor=center, inner sep=0}, "H"', curve={height=24pt}, from=1-1, to=1-3]
	\arrow[""{name=1, anchor=center, inner sep=0}, "G", curve={height=-24pt}, from=1-1, to=1-3]
	\arrow["\beta", shift left=3, shorten <=6pt, shorten >=6pt, Rightarrow, from=1, to=0]
	\arrow["\alpha"', shift right=3, shorten <=6pt, shorten >=6pt, Rightarrow, from=1, to=0]
\end{tikzcd}\]
The coequifier $\cc$ is a category that is given in terms of generators and relations as follows:
\begin{itemize}
\item the objects are the objects of $\cb$,
\item the morphisms are the equivalence classes of paths of morphisms in $\cb$, with the equivalence relation being generated by the following:
\begin{align*}
(1_B) &\sim ()_B             &\phantom{....}& \text{if }B \in \ob \cb, \\
(f,g) &\sim (g \circ f)      && \text{if } (f,g) \text{ is a composable pair in }\cb,\\
(\alpha_A) &\sim (\beta_A)   && \text{for all }A \in \ob \ca.
\end{align*}
\end{itemize}
\end{pr}

\subsection{Comonads in a 2-category}\label{SEKCE_monads_in_2cat}

The results we provide here have been proven in the classical papers \cite{formaltheory}, \cite[\S 3]{review}. The reason for including them in the thesis is to give a motivation for Section \ref{SEC_codescent_objects} and  later for Chapter \ref{CH_COMPUTATION}, in which the subject of study is the colax algebras for a 2-monad $T$. In the case $T$ is the identity 2-monad, colax algebras are the same thing as comonads and the results from these later sections become the results presented here.

\begin{defi}\label{DEFI_comonad}
A \textit{comonad} $(A,t,\delta,\epsilon)$ consists of an object $A$, 1-cell $t : A \to A$, 2-cells $\delta: t \Rightarrow t^2$, $\epsilon : t \Rightarrow 1_A$ (the \textit{comultiplication} and the \textit{counit}) such that the following 2-cells are equal:
\[\begin{tikzcd}[column sep=scriptsize]
	& A && A &&& A && A \\
	A && A && A & A && A && A \\
	& A && A &&& A && A
	\arrow[Rightarrow, no head, from=1-2, to=1-4]
	\arrow["t", from=1-4, to=2-5]
	\arrow["\delta", shorten <=6pt, shorten >=6pt, Rightarrow, from=1-4, to=3-4]
	\arrow[Rightarrow, no head, from=1-7, to=1-9]
	\arrow["t"', from=1-7, to=2-8]
	\arrow["t", from=1-9, to=2-10]
	\arrow[Rightarrow, no head, from=2-1, to=1-2]
	\arrow[""{name=0, anchor=center, inner sep=0}, Rightarrow, no head, from=2-1, to=2-3]
	\arrow["t"', from=2-1, to=3-2]
	\arrow[Rightarrow, no head, from=2-3, to=1-4]
	\arrow["t"', from=2-3, to=3-4]
	\arrow[Rightarrow, no head, from=2-6, to=1-7]
	\arrow["t"', from=2-6, to=3-7]
	\arrow[""{name=1, anchor=center, inner sep=0}, "t"', from=2-8, to=2-10]
	\arrow["t"', from=3-2, to=3-4]
	\arrow["t"', from=3-4, to=2-5]
	\arrow[Rightarrow, no head, from=3-7, to=2-8]
	\arrow[""{name=2, anchor=center, inner sep=0}, "t"', from=3-7, to=3-9]
	\arrow["t"', from=3-9, to=2-10]
	\arrow["\delta", shift right=5, shorten <=3pt, shorten >=3pt, Rightarrow, from=1-9, to=1]
	\arrow["\delta", shorten <=3pt, shorten >=3pt, Rightarrow, from=0, to=3-2]
	\arrow["\delta", shorten <=3pt, shorten >=3pt, Rightarrow, from=2-8, to=2]
\end{tikzcd}\]
and also satisfying that both of these 2-cells are equal to the identity on $t$:
\[\begin{tikzcd}
	&&&&& A &&&& A \\
	&&& A &&&& A \\
	A && A && A &&& A \\
	&&& A && A &&&& A
	\arrow[Rightarrow, no head, from=1-6, to=1-10]
	\arrow[Rightarrow, no head, from=1-6, to=2-8]
	\arrow["t"', from=1-6, to=4-6]
	\arrow["t", from=1-10, to=4-10]
	\arrow["t", from=2-4, to=3-5]
	\arrow["\delta", shorten <=6pt, shorten >=6pt, Rightarrow, from=2-4, to=4-4]
	\arrow[""{name=0, anchor=center, inner sep=0}, Rightarrow, no head, from=2-8, to=1-10]
	\arrow["t"', from=2-8, to=3-8]
	\arrow[curve={height=-24pt}, Rightarrow, no head, from=3-1, to=2-4]
	\arrow[Rightarrow, no head, from=3-1, to=3-3]
	\arrow[""{name=1, anchor=center, inner sep=0}, curve={height=24pt}, Rightarrow, no head, from=3-1, to=4-4]
	\arrow[Rightarrow, no head, from=3-3, to=2-4]
	\arrow["t"', from=3-3, to=4-4]
	\arrow[""{name=2, anchor=center, inner sep=0}, "t"', from=3-8, to=4-10]
	\arrow["t"', from=4-4, to=3-5]
	\arrow[Rightarrow, no head, from=4-6, to=3-8]
	\arrow[""{name=3, anchor=center, inner sep=0}, Rightarrow, no head, from=4-6, to=4-10]
	\arrow["\delta", shorten <=9pt, shorten >=9pt, Rightarrow, from=0, to=2]
	\arrow["\epsilon"', shorten >=5pt, Rightarrow, from=3-3, to=1]
	\arrow["\epsilon", shorten >=3pt, Rightarrow, from=3-8, to=3]
\end{tikzcd}\]
\end{defi}

\begin{vari}
A \textit{monad} in $\ck$ is a comonad in $\ck^{co}$. We will denote the data by tuple $(A,t,\mu: t^2 \Rightarrow t,\eta: 1_A \Rightarrow t)$.
\end{vari}

\begin{defi}\label{DEFI_comonadfun}
A \textit{comonad functor} $(f,\overline{f}): (A,t,\delta,\epsilon) \to (B,s,\delta', \epsilon')$ between two comonads is a pair of a 1-cell $f: A \to B$ and a 2-cell $\overline{f}: s f \Rightarrow f t$ such that the following two pairs of 2-cells are equal:
\[\begin{tikzcd}
	A && B &&& A && B \\
	& A && B &&&&& B \\
	A &&&& {=} & A && B \\
	& A && B &&& A && B
	\arrow["f", from=1-1, to=1-3]
	\arrow[""{name=0, anchor=center, inner sep=0}, Rightarrow, no head, from=1-1, to=2-2]
	\arrow["t"', from=1-1, to=3-1]
	\arrow[Rightarrow, no head, from=1-3, to=2-4]
	\arrow[""{name=1, anchor=center, inner sep=0}, "f", from=1-6, to=1-8]
	\arrow["t"', from=1-6, to=3-6]
	\arrow[""{name=2, anchor=center, inner sep=0}, Rightarrow, no head, from=1-8, to=2-9]
	\arrow["s"', from=1-8, to=3-8]
	\arrow[""{name=3, anchor=center, inner sep=0}, "f"{description}, from=2-2, to=2-4]
	\arrow["t", from=2-2, to=4-2]
	\arrow["s", from=2-4, to=4-4]
	\arrow["s", from=2-9, to=4-9]
	\arrow[""{name=4, anchor=center, inner sep=0}, "t"', from=3-1, to=4-2]
	\arrow[""{name=5, anchor=center, inner sep=0}, "f"{description}, from=3-6, to=3-8]
	\arrow["t"', from=3-6, to=4-7]
	\arrow[""{name=6, anchor=center, inner sep=0}, "s"', from=3-8, to=4-9]
	\arrow[""{name=7, anchor=center, inner sep=0}, "f"', from=4-2, to=4-4]
	\arrow[""{name=8, anchor=center, inner sep=0}, "f"', from=4-7, to=4-9]
	\arrow["\delta", shorten <=9pt, shorten >=9pt, Rightarrow, from=0, to=4]
	\arrow["{\overline{f}}", shorten <=9pt, shorten >=9pt, Rightarrow, from=1, to=5]
	\arrow["{\delta'}", shorten <=9pt, shorten >=9pt, Rightarrow, from=2, to=6]
	\arrow["{\overline{f}}", shorten <=9pt, shorten >=9pt, Rightarrow, from=3, to=7]
	\arrow["{\overline{f}}"', shorten <=3pt, shorten >=3pt, Rightarrow, from=3-8, to=8]
\end{tikzcd}\]
\[\begin{tikzcd}
	A && B &&& A && B \\
	& A && B &&&&& B \\
	&&&& {=} \\
	& A && B &&& A && B
	\arrow["f", from=1-1, to=1-3]
	\arrow[Rightarrow, no head, from=1-1, to=2-2]
	\arrow[""{name=0, anchor=center, inner sep=0}, curve={height=24pt}, Rightarrow, no head, from=1-1, to=4-2]
	\arrow[Rightarrow, no head, from=1-3, to=2-4]
	\arrow["f", from=1-6, to=1-8]
	\arrow[curve={height=24pt}, Rightarrow, no head, from=1-6, to=4-7]
	\arrow[Rightarrow, no head, from=1-8, to=2-9]
	\arrow[""{name=1, anchor=center, inner sep=0}, curve={height=24pt}, Rightarrow, no head, from=1-8, to=4-9]
	\arrow[""{name=2, anchor=center, inner sep=0}, "f", from=2-2, to=2-4]
	\arrow["t"', from=2-2, to=4-2]
	\arrow["s", from=2-4, to=4-4]
	\arrow["s", from=2-9, to=4-9]
	\arrow[""{name=3, anchor=center, inner sep=0}, "f"', from=4-2, to=4-4]
	\arrow["f"', from=4-7, to=4-9]
	\arrow["\epsilon"', shorten >=6pt, Rightarrow, from=2-2, to=0]
	\arrow["{\overline{f}}", shorten <=9pt, shorten >=9pt, Rightarrow, from=2, to=3]
	\arrow["{\epsilon'}"', shorten >=6pt, Rightarrow, from=2-9, to=1]
\end{tikzcd}\]
\end{defi}

\begin{vari}
We will encounter the following variations:
\begin{itemize}
\item in case the 2-cell $\overline{f}$ goes in the other direction and appropriate analogues of the above equations hold, $(f,\overline{f})$ is called a \textit{comonad opfunctor},
\item a comonad functor in $\ck^{op}$ is called a \textit{monad opfunctor},
\item in case the 1-cell component is the identity (so that $A = B$), we will call $\overline{f}: t \Rightarrow s$ a \textit{comonad morphism}.
\end{itemize}
\end{vari}

\begin{defi}
A \textit{comonad functor transformation} between two comonad functors $(f,\overline{f}), (g,\overline{g}): (A,t,\delta,\epsilon) \to (B,s,\delta',\epsilon')$ is a 2-cell $\alpha: f \Rightarrow g$ satisfying the equation:
\[\begin{tikzcd}
	A && B && A && B \\
	&&& {=} \\
	A && B && A && B
	\arrow[""{name=0, anchor=center, inner sep=0}, "g"', from=1-1, to=1-3]
	\arrow[""{name=1, anchor=center, inner sep=0}, "f", curve={height=-24pt}, from=1-1, to=1-3]
	\arrow["t"', from=1-1, to=3-1]
	\arrow["s", from=1-3, to=3-3]
	\arrow[""{name=2, anchor=center, inner sep=0}, "f", from=1-5, to=1-7]
	\arrow["t"', from=1-5, to=3-5]
	\arrow["s", from=1-7, to=3-7]
	\arrow[""{name=3, anchor=center, inner sep=0}, "g"', from=3-1, to=3-3]
	\arrow[""{name=4, anchor=center, inner sep=0}, "f", from=3-5, to=3-7]
	\arrow[""{name=5, anchor=center, inner sep=0}, "g"', curve={height=24pt}, from=3-5, to=3-7]
	\arrow["\alpha", shorten <=3pt, shorten >=3pt, Rightarrow, from=1, to=0]
	\arrow["{\overline{g}}", shorten <=13pt, shorten >=9pt, Rightarrow, from=0, to=3]
	\arrow["{\overline{f}}", shorten <=9pt, shorten >=13pt, Rightarrow, from=2, to=4]
	\arrow["\alpha", shorten <=3pt, shorten >=3pt, Rightarrow, from=4, to=5]
\end{tikzcd}\]
\end{defi}

\begin{nota}
Denote by:
\begin{itemize}
\item $\cmndl{\ck}$ the 2-category of comonads, comonad functors, comonad functor transformations,
\item $\cmndioo{\ck}{A}$ the 1-category of comonads with an underlying object $A$ and monad morphisms between them (comonad functors whose 1-cell component is the identity).
\end{itemize}
\end{nota}

\begin{pozn}
$\cmndioo{\ck}{A}$ is equivalently the category of comonoids in the strict monoi\-dal category $\ck(A,A)$.
\end{pozn}

\begin{theo}\label{THM_strictification_of_comonads}
Let $\ck$ be a 2-category and denote by $J$ the 2-functor sending an object to the identity comonad on it. Then $\ck$ admits Kleisli objects of comonads if and only if $J$ admits a left 2-adjoint:
\[\begin{tikzcd}
	\ck && {\cmndl{\ck}}
	\arrow[""{name=0, anchor=center, inner sep=0}, "J"', curve={height=24pt}, from=1-1, to=1-3]
	\arrow[""{name=1, anchor=center, inner sep=0}, "{\text{Kl}(-)}"', curve={height=24pt}, dashed, from=1-3, to=1-1]
	\arrow["\dashv"{anchor=center, rotate=-90}, draw=none, from=1, to=0]
\end{tikzcd}\]
\end{theo}
\begin{proof}
Notice that the hom-category:
\[
\cmndl{\ck}((A,t,\delta,\epsilon),(B,1,1,1)),
\]
is isomorphic to the category of cocones for the Kleisli object (dual of Example \ref{PR_EM_Kleisli_objects}). This hom-category being isomorphic to $\ck(A',B)$ for some object $A' \in \ob \ck$ thus happens if and only if this object is the Kleisli object for $(A,t,\delta,\epsilon)$. 
\end{proof}

\noindent The following is the first instance of what we will call a \textbf{coherence theorem for colax algebras}:

\begin{theo}\label{THM_laxcoherencethm_for_identity2monad}
Let $\ck$ be a 2-category that admits Kleisli objects of comonads. In the above 2-adjunction, each component of the unit has a left adjoint in $\ck$.
\end{theo}
\begin{proof}
Denote by $(f_t: A \to A_t, \xi: f_t \Rightarrow f_t t)$ the Kleisli object of $(A,t,\delta,\epsilon)$. It is equivalently a comonad opfunctor $(A,t,\delta,\epsilon) \to (A_t,1,1,1)$. Since the pair $(t, \delta)$ is a cocone, by the one-dimensional universal property there exists a unique 1-cell $u_t: A_t \to A$ with the property that:
\begin{align*}
u_t f_t &= t,\\
u_t \xi &= \delta.
\end{align*}

\noindent Next, notice that the 2-cell $\xi: f_t \Rightarrow f_t t = f_t u_t f_t$ is a cocone morphism:
\[
\xi : (f_t, \xi) \to (f_t u_t f_t, f_t u_t \xi).
\]
By the two-dimensional universal property, there is a unique 2-cell $\widetilde{\eta}: 1_{A_t} \Rightarrow f_t u_t$ with the property that $\widetilde{\eta} f_t = \xi$. We claim that this data gives an adjunction:
\[
(\epsilon, \widetilde{\eta}): u_t \dashv f_t \hspace{0.5cm}\text{ in } \ck.
\]
The second triangle equality:
\[
f_t \epsilon \circ \widetilde{\eta} f_t = 1_A,
\]
follows immediately from the cocone axiom for $(f_t, \xi)$. To prove the first triangle equality:
\[
\epsilon u_t \circ u_t \widetilde{\eta} = 1_{A_t},
\]
it suffices to prove that it holds after pre-composing with $f_t: A \to A_t$. It then follows from the comonad counit axiom.
\end{proof}

\section{Colax functors and transformations}\label{SEC_lax_transf}

The main reference for the definitions from this section is the classical text \cite{introbicats}, or the more modern \cite{2dimensionalcats}. Each of these also contains the definition of a \textit{bicategory} (a weak 2-category), but since our focus will be on weak \textbf{morphisms} of 2-categories, not weak 2-categories, we omit the definition.

\begin{defi}
Let $\ca, \cb$ be 2-categories. A \textit{colax functor} $F: \ca \to \cb$ consists of the following:
\begin{itemize}
\item a function $F_0: \ob \ca \to \ob \cb$,
\item for every pair $A, B$ of objects of $\ca$ a functor $F_{A,B}: \ca(A,B) \to \cb(FA,FB)$,
\item for every $A \in \ca$ a 2-cell $\iota_A : F1_A \Rightarrow 1_{FA}$ called the \textit{counitor},
\item for every $A,B,C \in \ca$ a natural transformation called the \textit{coassociator}:
\[\begin{tikzcd}
	{\ca(A,B)\times \ca(B,C)} && {\cb(FA,FB)\times \cb(FB,FC)} \\
	{\ca(A,C)} && {\cb(FA,FC)}
	\arrow[""{name=0, anchor=center, inner sep=0}, "{F_{A,B} \times F_{B,C}}", from=1-1, to=1-3]
	\arrow[""{name=0p, anchor=center, inner sep=0}, phantom, from=1-1, to=1-3, start anchor=center, end anchor=center]
	\arrow["\comp"', from=1-1, to=2-1]
	\arrow["\comp", from=1-3, to=2-3]
	\arrow[""{name=1, anchor=center, inner sep=0}, "{F_{A,C}}"', from=2-1, to=2-3]
	\arrow[""{name=1p, anchor=center, inner sep=0}, phantom, from=2-1, to=2-3, start anchor=center, end anchor=center]
	\arrow["\gamma"', shorten <=4pt, shorten >=4pt, Rightarrow, from=1p, to=0p]
\end{tikzcd}\]

Such that for each triple $(f,g,h)$ of composable morphisms in $\ca$, the following pairs of 2-cells are equal:
\[\begin{tikzcd}[column sep=small]
	FA && FC && FD & FA && FB && FD \\
	& FB &&&&&&& FC
	\arrow[""{name=0, anchor=center, inner sep=0}, "{F(gf)}", from=1-1, to=1-3]
	\arrow[""{name=1, anchor=center, inner sep=0}, "{F(hgf)}", curve={height=-30pt}, from=1-1, to=1-5]
	\arrow[""{name=2, anchor=center, inner sep=0}, draw=none, from=1-1, to=1-5]
	\arrow["Ff"', from=1-1, to=2-2]
	\arrow["Fh"', from=1-3, to=1-5]
	\arrow["Ff"', from=1-6, to=1-8]
	\arrow[""{name=3, anchor=center, inner sep=0}, "{F(hgf)}", curve={height=-30pt}, from=1-6, to=1-10]
	\arrow[""{name=4, anchor=center, inner sep=0}, draw=none, from=1-6, to=1-10]
	\arrow[""{name=5, anchor=center, inner sep=0}, "{F(hg)}", from=1-8, to=1-10]
	\arrow["Fg"', from=1-8, to=2-9]
	\arrow["Fg"', from=2-2, to=1-3]
	\arrow["Fh"', from=2-9, to=1-10]
	\arrow["{\gamma_{gf,h}}", shorten <=4pt, shorten >=4pt, Rightarrow, from=1, to=2]
	\arrow["{\gamma_{f,g}}", shorten <=3pt, Rightarrow, from=0, to=2-2]
	\arrow["{\gamma_{f,hg}}"', shorten <=4pt, shorten >=4pt, Rightarrow, from=3, to=4]
	\arrow["{\gamma_{g,h}}", shorten <=3pt, Rightarrow, from=5, to=2-9]
\end{tikzcd}\]
and for every morphism $f: A \to B$ in $\ca$, both of the 2-cells below are equal to the identity 2-cell on $Ff$:
\[\begin{tikzcd}[column sep=scriptsize]
	FA && FB && FB & FA && FA && FB
	\arrow["Ff"', from=1-1, to=1-3]
	\arrow[""{name=0, anchor=center, inner sep=0}, "Ff", curve={height=-30pt}, from=1-1, to=1-5]
	\arrow[""{name=1, anchor=center, inner sep=0}, draw=none, from=1-1, to=1-5]
	\arrow[""{name=2, anchor=center, inner sep=0}, curve={height=30pt}, Rightarrow, no head, from=1-3, to=1-5]
	\arrow[""{name=3, anchor=center, inner sep=0}, "{F1_B}"{description}, from=1-3, to=1-5]
	\arrow[""{name=4, anchor=center, inner sep=0}, curve={height=30pt}, Rightarrow, no head, from=1-6, to=1-8]
	\arrow[""{name=5, anchor=center, inner sep=0}, "{F1_A}"{description}, from=1-6, to=1-8]
	\arrow[""{name=6, anchor=center, inner sep=0}, "Ff", curve={height=-30pt}, from=1-6, to=1-10]
	\arrow[""{name=7, anchor=center, inner sep=0}, draw=none, from=1-6, to=1-10]
	\arrow["Ff"', from=1-8, to=1-10]
	\arrow["{\gamma_{f,1_B}}", shorten <=4pt, shorten >=4pt, Rightarrow, from=0, to=1]
	\arrow["{\iota_B}", shorten <=4pt, shorten >=4pt, Rightarrow, from=3, to=2]
	\arrow["{\iota_A}", shorten <=4pt, shorten >=4pt, Rightarrow, from=5, to=4]
	\arrow["{\gamma_{1_A,f}}", shorten <=4pt, shorten >=4pt, Rightarrow, from=6, to=7]
\end{tikzcd}\]

\end{itemize}
\begin{vari} We will encounter the following variations:
\begin{itemize}
\item if $\gamma$, $\iota$ go in the other direction, this is called a \textit{lax functor}.
\item if $\gamma$, $\iota$ are invertible, this is called a \textit{pseudofunctor},
\item if $\gamma$, $\iota$ are the identities, this is what is called a \textit{2-functor}.
\end{itemize}
\end{vari}
\end{defi}

\noindent For simplicity, we will always use the letters $\gamma, \iota$ for the associator and the unitor of a~colax functor, and always omit the index for any of their components.

\begin{defi}\label{DEFI_colax_nat_tr}
Given 2-categories $\ca, \cb$ and two colax functors $F,G: \ca \to \cb$, a \textit{colax natural transformation} $\alpha: F \Rightarrow G$ consists of the following data:
\begin{itemize}
\item for every $A \in \ca$ a 1-cell $\alpha_A : FA \to GA$,
\item for every $f: A \to B \in \ca$ a 2-cell:
\[\begin{tikzcd}
	FA && GA \\
	FB && GB
	\arrow[""{name=0, anchor=center, inner sep=0}, "{\alpha_B}"', from=2-1, to=2-3]
	\arrow[""{name=1, anchor=center, inner sep=0}, "{\alpha_A}", from=1-1, to=1-3]
	\arrow["Ff"', from=1-1, to=2-1]
	\arrow["Gf", from=1-3, to=2-3]
	\arrow["{\alpha_f}"', shorten <=4pt, shorten >=4pt, Rightarrow, from=0, to=1]
\end{tikzcd}\]
\end{itemize}
These must satisfy the \textit{unit}, \textit{composition} axioms for every $A \in \ca$ and every composable pair of morphisms $(f,g)$ in $\ca$:
\[\begin{tikzcd}
	FA && GA &&&& FA && GA \\
	&&&& {=} \\
	FA && GA &&&& FA && GA
	\arrow["{\alpha_A}", from=1-1, to=1-3]
	\arrow[""{name=0, anchor=center, inner sep=0}, curve={height=-24pt}, Rightarrow, no head, from=1-1, to=3-1]
	\arrow[""{name=1, anchor=center, inner sep=0}, "{F1_A}"', curve={height=24pt}, from=1-1, to=3-1]
	\arrow[curve={height=-24pt}, Rightarrow, no head, from=1-3, to=3-3]
	\arrow[""{name=2, anchor=center, inner sep=0}, "{\alpha_A}", from=1-7, to=1-9]
	\arrow["{F1_A}"', curve={height=24pt}, from=1-7, to=3-7]
	\arrow[""{name=3, anchor=center, inner sep=0}, curve={height=-24pt}, Rightarrow, no head, from=1-9, to=3-9]
	\arrow[""{name=4, anchor=center, inner sep=0}, "{G1_A}"{description}, curve={height=24pt}, from=1-9, to=3-9]
	\arrow["{\alpha_A}"', from=3-1, to=3-3]
	\arrow[""{name=5, anchor=center, inner sep=0}, "{\alpha_A}"', from=3-7, to=3-9]
	\arrow["\iota", shorten <=10pt, shorten >=10pt, Rightarrow, from=1, to=0]
	\arrow["\iota"', shorten <=10pt, shorten >=10pt, Rightarrow, from=4, to=3]
	\arrow["{\alpha_{1_A}}", shift left=5, shorten <=9pt, shorten >=9pt, Rightarrow, from=5, to=2]
\end{tikzcd}\]

\[\begin{tikzcd}
	FA && GA && FA && GA \\
	\\
	FB && GB & {=} &&& GB \\
	\\
	FC && GC && FC && GC
	\arrow[""{name=0, anchor=center, inner sep=0}, "{\alpha_A}", from=1-1, to=1-3]
	\arrow["Ff"', from=1-1, to=3-1]
	\arrow[""{name=1, anchor=center, inner sep=0}, "{F(gf)}"', shift right=3, curve={height=30pt}, from=1-1, to=5-1]
	\arrow[""{name=1p, anchor=center, inner sep=0}, phantom, from=1-1, to=5-1, start anchor=center, end anchor=center, shift right=3, curve={height=30pt}]
	\arrow["Gf", from=1-3, to=3-3]
	\arrow[""{name=2, anchor=center, inner sep=0}, "{\alpha_A}", from=1-5, to=1-7]
	\arrow["{F(gf)}"', from=1-5, to=5-5]
	\arrow["Gf", from=1-7, to=3-7]
	\arrow[""{name=3, anchor=center, inner sep=0}, "{G(gf)}"{description, pos=0.3}, shift right=5, curve={height=24pt}, from=1-7, to=5-7]
	\arrow[""{name=3p, anchor=center, inner sep=0}, phantom, from=1-7, to=5-7, start anchor=center, end anchor=center, shift right=5, curve={height=24pt}]
	\arrow[""{name=4, anchor=center, inner sep=0}, "{\alpha_B}"{description}, from=3-1, to=3-3]
	\arrow["Fg"', from=3-1, to=5-1]
	\arrow["Gg", from=3-3, to=5-3]
	\arrow["Gg", from=3-7, to=5-7]
	\arrow[""{name=5, anchor=center, inner sep=0}, "{\alpha_C}"', from=5-1, to=5-3]
	\arrow[""{name=6, anchor=center, inner sep=0}, "{\alpha_C}"', from=5-5, to=5-7]
	\arrow["\gamma", Rightarrow, from=1p, to=3-1]
	\arrow["\gamma", Rightarrow, from=3p, to=3-7]
	\arrow["{\alpha_f}"', shorten <=9pt, shorten >=9pt, Rightarrow, from=4, to=0]
	\arrow["{\alpha_g}"', shorten <=9pt, shorten >=9pt, Rightarrow, from=5, to=4]
	\arrow["{\alpha_{gf}}", shift left=5, shorten <=17pt, shorten >=17pt, Rightarrow, from=6, to=2]
\end{tikzcd}\]
together with the \textit{local naturality condition} for every 2-cell $\alpha: f \Rightarrow g$ in $\ca$:
\[\begin{tikzcd}[column sep=scriptsize]
	FA &&& GA &&&& FA &&& GA \\
	&&&&& {=} \\
	FB &&& GB &&&& FB &&& GB
	\arrow[""{name=0, anchor=center, inner sep=0}, "{\alpha_B}"{description}, from=3-1, to=3-4]
	\arrow[""{name=1, anchor=center, inner sep=0}, "{\alpha_A}", from=1-1, to=1-4]
	\arrow["Ff"', from=1-1, to=3-1]
	\arrow[""{name=2, anchor=center, inner sep=0}, "Gf"', curve={height=24pt}, from=1-4, to=3-4]
	\arrow[""{name=3, anchor=center, inner sep=0}, "Gg", curve={height=-24pt}, from=1-4, to=3-4]
	\arrow[""{name=4, anchor=center, inner sep=0}, "Fg", curve={height=-24pt}, from=1-8, to=3-8]
	\arrow[""{name=5, anchor=center, inner sep=0}, "Ff"', curve={height=24pt}, from=1-8, to=3-8]
	\arrow[""{name=6, anchor=center, inner sep=0}, "{\alpha_A}", from=1-8, to=1-11]
	\arrow["Gg", from=1-11, to=3-11]
	\arrow[""{name=7, anchor=center, inner sep=0}, "{\alpha_B}"', from=3-8, to=3-11]
	\arrow["{\alpha_f}", shorten <=9pt, shorten >=9pt, Rightarrow, from=0, to=1]
	\arrow["G\alpha", shorten <=10pt, shorten >=10pt, Rightarrow, from=2, to=3]
	\arrow["{\alpha_g}"', shorten <=9pt, shorten >=9pt, Rightarrow, from=7, to=6]
	\arrow["F\alpha", shorten <=10pt, shorten >=10pt, Rightarrow, from=5, to=4]
\end{tikzcd}\]
\begin{vari} We will encounter the following variations:
\begin{itemize}
\item if the 2-cells $\alpha_f$ go in the other direction, this is referred to as a \textit{lax natural transformation},
\item if $\alpha_f$ is invertible for all morphisms $f$, $\alpha$ is called a \textit{pseudonatural transformation},
\item if the $\alpha_f$'s are the identities, we use the term \textit{2-natural transformation},
\item if for every object $A$, $\alpha_A$ is the identity (in particular $F,G$ agree on objects), $\alpha$ is what has been called an \textit{icon} in \cite{icons} (it stands for \textit{identity-component oplax natural transformation}).
\end{itemize}
Let us also fix the terminology for the notions categorifying the notion of a ``natural isomorphism” from ordinary category theory:
\begin{itemize}
\item If the 1-cell components $\alpha_A$ of a pseudonatural transformation are the equivalences in $\cb$, $\alpha$ is called a \textit{pseudonatural equivalence}. In case $\alpha$ is 2-natural, we will call this a \textit{2-natural equivalence}.
\item If the 1-cell components $\alpha_A$ of a 2-natural transformation are isomorphisms in $\cb$, $\alpha$ is called a \textit{2-natural isomorphism}.
\end{itemize}
\end{vari}
\end{defi}

\begin{defi}
Given two colax natural transformations $\alpha, \beta$ between pseudofunctors $F,G: \ck \to \cl$, a \textit{modification} $\Gamma: \alpha \to \beta$ consists of a 2-cell $\Gamma_A: \alpha_A \Rightarrow \beta_A$ for every object $A \in \ck$, subject to the \textit{modification axiom} for each 1-cell in $\ck$:
\[\begin{tikzcd}
	FA && GA && FA && GA \\
	&&& {=} \\
	FB && GB && FB && GB
	\arrow[""{name=0, anchor=center, inner sep=0}, "{\alpha_A}"{description}, curve={height=24pt}, from=1-1, to=1-3]
	\arrow[""{name=1, anchor=center, inner sep=0}, "{\beta_A}", curve={height=-24pt}, from=1-1, to=1-3]
	\arrow["Ff"', from=1-1, to=3-1]
	\arrow["Gf", from=1-3, to=3-3]
	\arrow[""{name=2, anchor=center, inner sep=0}, "{\beta_A}", from=1-5, to=1-7]
	\arrow["Ff"', from=1-5, to=3-5]
	\arrow["Gf", from=1-7, to=3-7]
	\arrow[""{name=3, anchor=center, inner sep=0}, "{\alpha_B}"', from=3-1, to=3-3]
	\arrow[""{name=4, anchor=center, inner sep=0}, "{\alpha_B}"', curve={height=24pt}, from=3-5, to=3-7]
	\arrow[""{name=5, anchor=center, inner sep=0}, "{\beta_B}"{description}, curve={height=-24pt}, from=3-5, to=3-7]
	\arrow["{\Gamma_A}"', shorten <=6pt, shorten >=6pt, Rightarrow, from=0, to=1]
	\arrow["{\alpha_f}"', shorten <=5pt, shorten >=5pt, Rightarrow, from=3, to=0]
	\arrow["{\beta_f}"', shorten <=5pt, shorten >=5pt, Rightarrow, from=5, to=2]
	\arrow["{\Gamma_B}"', shorten <=6pt, shorten >=6pt, Rightarrow, from=4, to=5]
\end{tikzcd}\]
\end{defi}

\begin{pr}\label{PR_colax_nat_gives_modif}
Given an endo-2-functor $T: \ca \to \ca$, any colax natural transformation $c: T \Rightarrow 1_\ca$ induces a modification $(cc): c \circ Tc \to c \circ cT$, whose component at $A \in \ca$ is given by:
\[
(cc)_A := c_{c_A} :c_A \circ Tc_A \Rightarrow c_A \circ c_{TA}.
\]
\end{pr}

\begin{pozn}[Pseudofunctors preserve colax natural transformations]

Notice that if $H: \cc \to \cd$ is a pseudofunctor and $\alpha: F \Rightarrow G: \cb \to \cc$ is colax natural, there is an induced colax natural transformation $H\alpha$ whose 1-cell component at $A$ is $H\alpha_A$ and whose 2-cell component at a morphism $f: A \to B$ is the following composite 2-cell that we denote by $(H\alpha)_f$:
\[\begin{tikzcd}
	FA && GA \\
	\\
	\\
	FB && GB
	\arrow["{H\alpha_A}", from=1-1, to=1-3]
	\arrow["HFf"', from=1-1, to=4-1]
	\arrow[""{name=0, anchor=center, inner sep=0}, curve={height=-18pt}, from=1-1, to=4-3]
	\arrow[""{name=1, anchor=center, inner sep=0}, curve={height=18pt}, from=1-1, to=4-3]
	\arrow["HGf", from=1-3, to=4-3]
	\arrow["{H\alpha_B}"', from=4-1, to=4-3]
	\arrow["\gamma", shorten <=5pt, Rightarrow, from=0, to=1-3]
	\arrow["{H\alpha_f}"', shorten <=6pt, shorten >=6pt, Rightarrow, from=1, to=0]
	\arrow["{\gamma^{-1}}"', shorten >=5pt, Rightarrow, from=4-1, to=1]
\end{tikzcd}\]
\end{pozn}

\noindent 2-categories, lax functors and icons form a 2-category that we denote by $\twocat_l$. The following is a generalization of the bijective on objects/fully faithful factorisation of functors between categories:

\begin{prop}\label{THMP_2Catl_OFS}
The underlying 1-category of $\twocat_l$ admits an orthogonal factorization system\footnote{See Definition \ref{DEFI_ofs}.} $(\ce,\cm)$, where:
\begin{itemize}
\item $\ce$ is the class of lax functors that are bijective on objects,
\item $\cm$ is the class of fully faithful 2-functors.
\end{itemize}
\end{prop}
\begin{proof}
Given a lax functor $F: \cc \to \cl$ between 2-categories, the factorization:
\[\begin{tikzcd}
	\cc && \cl \\
	& {\widetilde{\cc}}
	\arrow["{(F,\gamma,\iota)}", squiggly, from=1-1, to=1-3]
	\arrow["{(F_\circ,\gamma_\circ,\iota_\circ)}"', squiggly, from=1-1, to=2-2]
	\arrow["{F_\bullet}"', from=2-2, to=1-3]
\end{tikzcd}\]
is given as follows. The 2-category $\widetilde{\cc}$ has objects of $\cc$ as objects, with hom categories being given by:
\[
\widetilde{\cc}(A,B) := \cl(FA,FB).
\]
The identity on $A$ is given by $1_{FA}$. We will denote the objects of this hom by $(A,f,B)$, where $(f: FA \to FB) \in \cl$. The lax functor $F_\circ: \cc \to \widetilde{\cc}$ is the identity on objects, and sends:
\[
(f: A \to B) \mapsto (A, Ff, B) \in \widetilde{\cc}.
\]
Given $A \in \cc$ and a pair $g: A \to B, f: B \to C$, the associator and unitor for $F_\circ$ are given by those of $F$:
\begin{align*}
\gamma_{\circ} &:= (A, \gamma: Ff Fg \Rightarrow F(fg), C),\\
\iota_{\circ} &:= (A, \iota: 1_{FA} \Rightarrow F1_A, A).
\end{align*}
The lax functor axioms follow from those of $F$. The 2-functor $F_\bullet: \widetilde{\cc} \to \cl$ is given on objects by an assignment $A \mapsto FA$, and on hom categories is given by:
\begin{align*}
\widetilde{\cc}(A,B) &\to \cl(FA,FB),\\
(A, f: FA \to FB, B) &\mapsto f.
\end{align*}
This is clearly a fully faithful 2-functor satisfying $F = F_\bullet \circ F_\circ$. To prove the uniqueness of factorizations up to an isomorphism, assume that $F = G_\bullet \circ G_\circ$ is another $(\ce,\cm)$-factorization of $F$. We must show that there is a unique isomorphism making the following commute:
\[\begin{tikzcd}
	\cc && \cx && \cl \\
	\\
	\cc && {\widetilde{\cc}} && \cl
	\arrow[""{name=0, anchor=center, inner sep=0}, "{G_\circ}", from=1-1, to=1-3]
	\arrow[Rightarrow, no head, from=1-1, to=3-1]
	\arrow[""{name=1, anchor=center, inner sep=0}, "{G_\bullet}", from=1-3, to=1-5]
	\arrow[Rightarrow, no head, from=1-5, to=3-5]
	\arrow[""{name=2, anchor=center, inner sep=0}, "{F_\circ}"', from=3-1, to=3-3]
	\arrow["\Phi"', dashed, from=3-3, to=1-3]
	\arrow[""{name=3, anchor=center, inner sep=0}, "{F_\bullet}"', from=3-3, to=3-5]
	\arrow["{(A)}"{description}, draw=none, from=0, to=2]
	\arrow["{(B)}"{description}, draw=none, from=1, to=3]
\end{tikzcd}\]
The commutativity of the left square forces that $\Phi$ agrees with $G_\circ$ on objects -- $\Phi$ is thus also bijective on objects. The commutativity of the right square forces that for each $A,B \in \cc$ the action of $\Phi$ on hom categories is given by the following isomorphism:
\[\begin{tikzcd}
	{\widetilde{\cc}(A,B)} && {\cl(FA,FB)} && {\cx(G_\circ A, G_\circ B)}
	\arrow["{(F_{\bullet})_{A,B}}", from=1-1, to=1-3]
	\arrow["{(G_{\bullet}^{-1})_{G_\circ A, G_\circ B}}", from=1-3, to=1-5]
\end{tikzcd}\]
Since $G_\bullet \circ \Phi = F_\bullet$ is a 2-functor and $G_\bullet$ is fully faithful, $\Phi$ is a 2-functor. Since it is fully faithful and bijective on objects, it is an isomorphism.
\end{proof}

\begin{pr}
A lax functor $* \to \ck$ from a terminal 2-category into a 2-category is the same thing as a monad in $\ck$. We have mentioned (co)limits of such in Example \ref{PR_EM_Kleisli_objects}.
\end{pr}

\section{Colax adjunctions}\label{SEC_lax_adjs}

A lax adjunction is a categorification of an adjunction between functors where the unit and the counit are replaced by lax natural transformations, and the triangle identities are replaced by modifications.

Lax adjunctions (among with their many variants) have first been extensively studied in \cite[I,7.1]{formalcat} and \cite{bunge}. In the first named paper they have been called \textit{transcendental quasi-adjunctions} and in the second one \textbf{colax} adjunctions. See also \cite{localadj}.

\begin{defi}\label{DEFI_colax_adjunction}
A \textit{colax adjunction} consists of two pseudofunctors $U: \cd \to \cc$ and $F: \cc \to \cd$, two colax natural transformations:
\[
\eta : 1_{\cc} \Rightarrow UF, \hspace{2cm} \epsilon : FU \Rightarrow 1_{\cd},
\]
and two modifications:
\[\begin{tikzcd}
	F && FUF && U && UFU \\
	\\
	&& F &&&& U
	\arrow["F\eta", from=1-1, to=1-3]
	\arrow["{\epsilon F}", from=1-3, to=3-3]
	\arrow[""{name=0, anchor=center, inner sep=0}, Rightarrow, no head, from=1-1, to=3-3]
	\arrow[""{name=1, anchor=center, inner sep=0}, Rightarrow, no head, from=1-5, to=3-7]
	\arrow["{\eta U}", from=1-5, to=1-7]
	\arrow["{U \epsilon}", from=1-7, to=3-7]
	\arrow["\Psi"', shorten <=5pt, shorten >=5pt, Rightarrow, from=1-3, to=0]
	\arrow["\Phi", shorten <=5pt, shorten >=5pt, Rightarrow, from=1, to=1-7]
\end{tikzcd}\]

Before stating the axioms required, let us fix a convention: we will use the symbol $U \Psi$ to denote the modification obtained from $\Psi$ by not just applying $U$, but also by pre- and post-composing it with the associator and the unitor for $U$ so that its domain and codomain are $U\epsilon F \circ UF\eta$, $1_{UF}$. Let us use the same convention for $F \Phi$. The axioms are the \textit{swallowtail identities}, which assert that the two composite modifications below\footnote{The notation $\eta \eta$ is explained in Example \ref{PR_colax_nat_gives_modif}.} are the identities on $\eta$ and $\epsilon$:
\[
\adjustbox{scale=0.9,center}{
\begin{tikzcd}[column sep=scriptsize]
	{1_\cc} && UF &&& FU \\
	&&& {} &&&& {} & {} \\
	UF && UFUF & {} &&& {} & FUFU && FU \\
	& {} & {} &&&& {} \\
	&&&& UF &&& FU && {1_{\cd}}
	\arrow[curve={height=24pt}, Rightarrow, no head, from=3-1, to=5-5]
	\arrow[curve={height=-24pt}, Rightarrow, no head, from=1-3, to=5-5]
	\arrow["\eta"', from=1-1, to=3-1]
	\arrow["\eta", from=1-1, to=1-3]
	\arrow["UF\eta", from=1-3, to=3-3]
	\arrow["{\eta UF}"', from=3-1, to=3-3]
	\arrow["{\eta \eta}"{description}, shorten <=12pt, shorten >=12pt, Rightarrow, from=3-1, to=1-3]
	\arrow["{U\epsilon F}"{description}, from=3-3, to=5-5]
	\arrow["{F\eta U}"{description}, from=1-6, to=3-8]
	\arrow["FU\epsilon", from=3-8, to=3-10]
	\arrow["{\epsilon_{FU}}"', from=3-8, to=5-8]
	\arrow["\epsilon"', from=5-8, to=5-10]
	\arrow["\epsilon", from=3-10, to=5-10]
	\arrow[curve={height=24pt}, Rightarrow, no head, from=1-6, to=5-8]
	\arrow[curve={height=-24pt}, Rightarrow, no head, from=1-6, to=3-10]
	\arrow["{\epsilon \epsilon}", shorten <=12pt, shorten >=12pt, Rightarrow, from=3-10, to=5-8]
	\arrow[""{name=0, anchor=center, inner sep=0}, draw=none, from=3-3, to=3-4]
	\arrow[""{name=0p, anchor=center, inner sep=0}, phantom, from=3-3, to=3-4, start anchor=center, end anchor=center]
	\arrow[""{name=1, anchor=center, inner sep=0}, draw=none, from=4-2, to=4-3]
	\arrow[""{name=1p, anchor=center, inner sep=0}, phantom, from=4-2, to=4-3, start anchor=center, end anchor=center]
	\arrow[""{name=2, anchor=center, inner sep=0}, draw=none, from=2-8, to=2-9]
	\arrow[""{name=2p, anchor=center, inner sep=0}, phantom, from=2-8, to=2-9, start anchor=center, end anchor=center]
	\arrow[""{name=3, anchor=center, inner sep=0}, shift left=3, draw=none, from=3-7, to=3-8]
	\arrow[""{name=3p, anchor=center, inner sep=0}, phantom, from=3-7, to=3-8, start anchor=center, end anchor=center, shift left=3]
	\arrow["U\Psi"', shorten >=10pt, Rightarrow, from=0p, to=2-4]
	\arrow["{\Phi F}"', shorten >=4pt, Rightarrow, from=1p, to=3-3]
	\arrow["F\Phi"', shorten <=4pt, shorten >=2pt, Rightarrow, from=2p, to=3-8]
	\arrow["{\Psi U}"', shorten <=5pt, shorten >=9pt, Rightarrow, from=3p, to=4-7]
\end{tikzcd}
}\]

\begin{nota}
We will denote a colax adjunction as follows and say that $F$ is a \textit{left colax adjoint} to $U$:
\[\begin{tikzcd}
	{(\Psi, \Phi): (\epsilon, \eta):} & \cc && \cd
	\arrow[""{name=0, anchor=center, inner sep=0}, "F", curve={height=-24pt}, from=1-2, to=1-4]
	\arrow[""{name=1, anchor=center, inner sep=0}, "U", curve={height=-24pt}, from=1-4, to=1-2]
	\arrow["\tbot"{description}, draw=none, from=0, to=1]
\end{tikzcd}\]
\end{nota}
\end{defi}

\begin{vari}
There are several important variations or special cases:
\begin{itemize}
\item if $\epsilon, \eta$ are lax natural, $\Psi,\Phi$ go in the other directions and an appropriate dual of the swallowtail identities holds, we will call it a \textit{lax adjunction},
\item in case that $\epsilon, \eta$ are pseudonatural transformations and $\Psi ,\Phi$ are isomorphisms, we will use the term \textit{biadjunction}.
\item if $U,F$ are 2-functors, $\epsilon, \eta$ are 2-natural and $\Psi, \Phi$ are the identities, we will call this a \textit{2-adjunction}.
\end{itemize}
The last two cases are the more usual notion so we will use the usual symbol $\dashv$ instead of $\ddashv$\phantom{.} for them. Let us also fix the terminology for the notions categorifying an ``equivalence between categories” from ordinary category theory:

\begin{itemize}
\item a 2-adjunction in which $\epsilon, \eta$ are 2-natural isomorphisms will be called a \textit{2-equivalence},
\item a biadjunction in which $\epsilon, \eta$ are pseudonatural equivalences will be called a \textit{biequivalence}.
\end{itemize}
\end{vari}

In Theorem \ref{THM_folklore_biadjunctions} we will give a folklore characterization of biadjunctions using methods developed in that section.

\begin{prop}\label{THMP_colax_adj_adjoint_homcats}
Given a colax adjunction $(\Psi, \Phi): (\epsilon, \eta): F \ddashv U: \cd \to \cc$, for any pair of objects $A \in \ob \cc, B \in \ob \cd$ there is an induced adjunction of categories:
\begin{equation}\label{EQ_prop_colax_adj_adjoint_homcats}
\begin{tikzcd}
	{\cd(FA, B)} && {\cc(A,UB)}
	\arrow[""{name=0, anchor=center, inner sep=0}, "{(\eta_A)^* \circ U:}", curve={height=-24pt}, from=1-1, to=1-3]
	\arrow[""{name=1, anchor=center, inner sep=0}, "{(\epsilon_B)_* \circ F}", curve={height=-24pt}, from=1-3, to=1-1]
	\arrow["\dashv"{anchor=center, rotate=90}, draw=none, from=1, to=0]
\end{tikzcd}
\end{equation}
Moreover, the collection of such adjunctions is \textit{lax natural in $A$} and \textit{colax natural in $B$}. By this mean there are adjunctions:
\[\begin{tikzcd}
	{\cd(F-, B)} && {\cc(-,UB)} & {\text{ in }\Homl{\cc^{op}}{\cat}} \\
	\\
	{\cd(FA, -)} && {\cc(A,U-)} & {\text{ in }\Homc{\cd}{\cat}}
	\arrow[""{name=0, anchor=center, inner sep=0}, "{\eta^* \circ U:}", curve={height=-24pt}, from=1-1, to=1-3]
	\arrow[""{name=1, anchor=center, inner sep=0}, "{(\epsilon_B)_* \circ F}"{description}, curve={height=-24pt}, from=1-3, to=1-1]
	\arrow[""{name=2, anchor=center, inner sep=0}, "{\epsilon_* \circ F}"', curve={height=24pt}, from=3-1, to=3-3]
	\arrow[""{name=3, anchor=center, inner sep=0}, "{(\eta_A)^* \circ U:}"{description}, curve={height=24pt}, from=3-3, to=3-1]
	\arrow["\dashv"{anchor=center, rotate=90}, draw=none, from=1, to=0]
	\arrow["\dashv"{anchor=center, rotate=90}, draw=none, from=2, to=3]
\end{tikzcd}\]
\end{prop}
\begin{proof}
The counit and unit 2-cells evaluated at $h: FA \to B$ and $g: A \to UB$ are given as follows:
\[\begin{tikzcd}
	FA &&& A && UFA \\
	FUFA && FA & UB && UFUB \\
	FUB && B &&& UB
	\arrow["{Fy_A}", from=1-1, to=2-1]
	\arrow[""{name=0, anchor=center, inner sep=0}, Rightarrow, no head, from=1-1, to=2-3]
	\arrow[""{name=0p, anchor=center, inner sep=0}, phantom, from=1-1, to=2-3, start anchor=center, end anchor=center]
	\arrow[""{name=1, anchor=center, inner sep=0}, "{F(Uhy_A)}"', shift right=5, curve={height=30pt}, from=1-1, to=3-1]
	\arrow[""{name=1p, anchor=center, inner sep=0}, phantom, from=1-1, to=3-1, start anchor=center, end anchor=center, shift right=5, curve={height=30pt}]
	\arrow[""{name=2, anchor=center, inner sep=0}, curve={height=6pt}, draw=none, from=1-1, to=3-1]
	\arrow[""{name=2p, anchor=center, inner sep=0}, phantom, from=1-1, to=3-1, start anchor=center, end anchor=center, curve={height=6pt}]
	\arrow[""{name=3, anchor=center, inner sep=0}, "{\eta_A}", from=1-4, to=1-6]
	\arrow[""{name=3p, anchor=center, inner sep=0}, phantom, from=1-4, to=1-6, start anchor=center, end anchor=center]
	\arrow["g"', from=1-4, to=2-4]
	\arrow["UFg"', from=1-6, to=2-6]
	\arrow[""{name=4, anchor=center, inner sep=0}, "{U(\epsilon_B Fg)}", shift left=5, curve={height=-30pt}, from=1-6, to=3-6]
	\arrow[""{name=4p, anchor=center, inner sep=0}, phantom, from=1-6, to=3-6, start anchor=center, end anchor=center, shift left=5, curve={height=-30pt}]
	\arrow[""{name=5, anchor=center, inner sep=0}, curve={height=-6pt}, draw=none, from=1-6, to=3-6]
	\arrow[""{name=5p, anchor=center, inner sep=0}, phantom, from=1-6, to=3-6, start anchor=center, end anchor=center, curve={height=-6pt}]
	\arrow[""{name=6, anchor=center, inner sep=0}, "{\epsilon_{FA}}"{description, pos=0.7}, from=2-1, to=2-3]
	\arrow[""{name=6p, anchor=center, inner sep=0}, phantom, from=2-1, to=2-3, start anchor=center, end anchor=center]
	\arrow[""{name=6p, anchor=center, inner sep=0}, phantom, from=2-1, to=2-3, start anchor=center, end anchor=center]
	\arrow["FUh", from=2-1, to=3-1]
	\arrow["h", from=2-3, to=3-3]
	\arrow[""{name=7, anchor=center, inner sep=0}, "{\eta_{UB}}"{description, pos=0.3}, from=2-4, to=2-6]
	\arrow[""{name=7p, anchor=center, inner sep=0}, phantom, from=2-4, to=2-6, start anchor=center, end anchor=center]
	\arrow[""{name=7p, anchor=center, inner sep=0}, phantom, from=2-4, to=2-6, start anchor=center, end anchor=center]
	\arrow[""{name=8, anchor=center, inner sep=0}, Rightarrow, no head, from=2-4, to=3-6]
	\arrow[""{name=8p, anchor=center, inner sep=0}, phantom, from=2-4, to=3-6, start anchor=center, end anchor=center]
	\arrow["{U\epsilon_B}"', from=2-6, to=3-6]
	\arrow[""{name=9, anchor=center, inner sep=0}, "{\epsilon_B}"', from=3-1, to=3-3]
	\arrow[""{name=9p, anchor=center, inner sep=0}, phantom, from=3-1, to=3-3, start anchor=center, end anchor=center]
	\arrow["\gamma"{pos=0.2}, shorten >=17pt, Rightarrow, from=1p, to=2p]
	\arrow["{\gamma^{-1}}"{pos=0.7}, shorten <=17pt, shorten >=3pt, Rightarrow, from=5p, to=4p]
	\arrow["{\Psi_A}", shorten <=2pt, shorten >=2pt, Rightarrow, from=6p, to=0p]
	\arrow["{\eta_g}", shorten <=4pt, shorten >=4pt, Rightarrow, from=7p, to=3p]
	\arrow["{\Phi_B}"', shorten <=2pt, shorten >=2pt, Rightarrow, from=8p, to=7p]
	\arrow["{\epsilon_h}"', shorten <=4pt, shorten >=4pt, Rightarrow, from=9p, to=6p]
\end{tikzcd}\]
The triangle identities can be seen to essentially follow from the swallowtail identities -- we thus obtain the adjunction \eqref{EQ_prop_colax_adj_adjoint_homcats}.

The fact that the adjunctions are (co)lax natural in $A,B$ can be proven directly -- we do however prefer to use a more abstract argument that will appear later in Section \ref{SEC_psmonads}. For the remainder of the proof, we assume familiarity with that section. We will only prove the lax naturality in the first component, since the other component is done in an analogous way. Let us fix $B \in \ob \cd$. With $A$ replaced by ``blank”, we have the following $(\ob \cc)$-graded adjunction of categories:
\begin{equation}\label{EQP_homcats}
\begin{tikzcd}
	{\cd(F-, B)} && {\cc(-,UB)} & {\text{ in } [\ob \cc^{op}, \cat]}
	\arrow[""{name=0, anchor=center, inner sep=0}, "{(\eta_{(-)})^* \circ U}", curve={height=-24pt}, from=1-1, to=1-3]
	\arrow[""{name=1, anchor=center, inner sep=0}, "{(\epsilon_B)_* \circ F}", curve={height=-24pt}, from=1-3, to=1-1]
	\arrow["\dashv"{anchor=center, rotate=90}, draw=none, from=1, to=0]
\end{tikzcd}
\end{equation}

\noindent Consider the reslan 2-monad $T$ from Example \ref{PR_reslan2monad}. The collections:
\[
\cc(F-, B), \cd(-,UB): \ob \cc^{op} \to \cat,
\]
are of course pseudofunctors and so they are pseudo $T$-algebras. Notice also that the collection $(\epsilon_B)_* \circ F: \cd(-,UB) \Rightarrow \cc(F-,B)$ is a pseudonatural transformation (the pseudonaturality squares use the pseudofunctor associators) and so it is a pseudo $T$-morphism. By Corollary \ref{THMC_doct_adj_corollary}, the collection $\eta^* \circ U$ admits a canonical lax transformation structure and the adjunction \eqref{EQP_homcats} lifts to $\Homl{\cc^{op}}{\cat}$ (which is the 2-category pseudo algebras and lax morphisms for $T$).
\end{proof}

\begin{pozn}[Special cases]\label{POZN_colax_adj_adjoint_homcats_corefl_refl}
Given a colax adjunction:
\[
(\Psi, \Phi): (\epsilon, \eta): F \ddashv U: \cd \to \cc,
\]
the hom category adjunction \eqref{EQ_prop_colax_adj_adjoint_homcats} is:
\begin{itemize}
\item a reflection if and only if $\epsilon$ is pseudonatural and $\Psi$ is invertible,
\item a coreflection if and only if $\eta$ is pseudonatural and $\Phi$ is invertible.
\end{itemize}
The reader may also observe that the canonical lax naturality square for the transformation:
\[
\eta^* \circ U: \cd(F-,B) \Rightarrow \cc(-,UB),
\]
only depends on the pseudofunctor associators and the colax naturality square for $\eta$. In particular, the adjunctions \eqref{EQ_prop_colax_adj_adjoint_homcats} are:
\begin{itemize}
\item pseudonatural in $A$ if $\eta$ is pseudonatural,
\item pseudonatural in $B$ if $\epsilon$ is pseudonatural,
\end{itemize}
\end{pozn}

\begin{pr}
An existence of a colax left adjoint to the unique 2-functor $\ck \to *$ to the terminal 2-category implies the existence of an object $I \in \ob \ck$ with the property that for each $B \in \ob \ck$, the category $\ck(I,B)$ has an initial object. This is a special case of a new kind of weak colimit to be introduced in Section \ref{SEC_weak_limits}.
\end{pr}

\begin{pr}
In Chapter \ref{CH_colax_adjs_lax_idemp_psmons} we will show that colax adjunctions naturally appear in two-dimensional algebra whenever \textbf{lax} algebra morphisms are involved. An example of this kind will also be given in Proposition \ref{THMP_adj_colax_alg_base_2cat}.
\end{pr}

\begin{pozn}\label{POZN_lax_adjs_not_unique_up_to_equiv}
Contrary to the case of biadjunctions, left colax adjoints are not unique up to an equivalence, not even when $U$ is a 2-functor, $\eta$ is 2-natural and $\Psi$, $\Phi$ are the identities. An example will be given in Remark \ref{POZN_rali_limits_not_unique_up_to_equiv}.
\end{pozn}

As the reader has likely already noticed, we will primarily use the \textbf{colax} versions of weak functors, weak transformations, weak adjunctions. This is because in Chapter \ref{CH_colax_adjs_lax_idemp_psmons} we are building on the work of Bunge \cite{bunge} which uses colax structures as the base notion\footnote{Note that in the older terminology, they were called \textit{lax}.}.

\section{Pseudomonads}\label{SEC_psmonads}

\noindent \textbf{The section is organized as follows}:
\begin{itemize}
\item In Subsection \ref{SUBS_psmonads_defi} we define pseudomonads and their algebras (\textit{pseudomonads}, \cite{coherentapproach_HIST}, \cite{distrlawspsmonads_HIST}).
\item In Subsection \ref{SUBS_psmonads_doctrinal} we give an exposition on the doctrinal adjunction in two-dimen\-sional monad theory, and list some of its corollaries.
\item In Subsection \ref{SUBS_psmonads_adjsvsalgs} we describe a 2-functor that sends adjunctions to colax algebras and study its properties. The motivation for this is to make the exposiiton in the next section (Section \ref{SEC_lax_idemp_psmonads}) clearer.

\item The final Subsection \ref{SUBS_psmonads_examples} is here to list examples of 2-monads. Some of these examples will be used in later chapters as a demonstration of the developed theory (for instance in Subsection \ref{SEC_examples_lax_coh} and Corollary \ref{THMC_examples_of_rali_cocomplete_2cats}).
\end{itemize}

\begin{pozn}
We believe that it does not cost much to prove all the basic properties of pseudomonads in full generality (and we will do so here), even though the only examples we encounter will be either 2-monads or lax-idempotent pseudomonads (for this latter class, a convenient framework is available -- we will recall it in Section \ref{SECS_left_kan_psmons}).
\end{pozn}

\subsection{Definition}\label{SUBS_psmonads_defi}

\begin{defi}
Let $\ck$ be a 2-category. A \textit{pseudomonad} $(T,m,i,\mu,\eta,\beta)$ on $\ck$ consists of a pseudofunctor $T: \ck \to \ck$, pseudonatural transformations:
\[
m: T^2 \Rightarrow T, \hspace{2cm} i : 1_{\ck} \Rightarrow T,
\]
and invertible modifications:
\[\begin{tikzcd}
	{T^3} && {T^2} && T && {T^2} && T \\
	\\
	{T^2} && T &&&& T
	\arrow["mT", from=1-1, to=1-3]
	\arrow["Tm"', from=1-1, to=3-1]
	\arrow["m", from=1-3, to=3-3]
	\arrow["iT", from=1-5, to=1-7]
	\arrow[""{name=0, anchor=center, inner sep=0}, Rightarrow, no head, from=1-5, to=3-7]
	\arrow["m"{description}, from=1-7, to=3-7]
	\arrow["Ti"', from=1-9, to=1-7]
	\arrow[""{name=1, anchor=center, inner sep=0}, Rightarrow, no head, from=1-9, to=3-7]
	\arrow["\mu"', shorten <=11pt, shorten >=11pt, Rightarrow, from=3-1, to=1-3]
	\arrow["m"', from=3-1, to=3-3]
	\arrow["\beta"', shorten >=6pt, Rightarrow, from=1-7, to=0]
	\arrow["\eta"', shorten <=6pt, Rightarrow, from=1, to=1-7]
\end{tikzcd}\]
such that the following pairs of 2-cells are equal:
\[\begin{tikzcd}[column sep=tiny]
	& {T^3A} && {T^2A} &&& {T^3A} && {T^2A} \\
	{T^4A} && {T^3A} && TA & {T^4A} && {T^2A} && TA \\
	& {T^3A} && {T^2A} &&& {T^3A} && {T^2A}
	\arrow["{Tm_A}", from=1-2, to=1-4]
	\arrow["{m_A}", from=1-4, to=2-5]
	\arrow["{\mu_A}", shorten <=6pt, shorten >=6pt, Rightarrow, from=1-4, to=3-4]
	\arrow["{Tm_A}", from=1-7, to=1-9]
	\arrow["{m_{TA}}", from=1-7, to=2-8]
	\arrow["{m_{m_A}}"', shorten <=6pt, shorten >=6pt, Rightarrow, from=1-7, to=3-7]
	\arrow["{m_A}", from=1-9, to=2-10]
	\arrow["{T^2m_A}", from=2-1, to=1-2]
	\arrow[""{name=0, anchor=center, inner sep=0}, "{Tm_{TA}}"', from=2-1, to=2-3]
	\arrow["{m_{T^2A}}"', from=2-1, to=3-2]
	\arrow["{Tm_A}", from=2-3, to=1-4]
	\arrow["{m_{TA}}"', from=2-3, to=3-4]
	\arrow["{T^2m_A}", from=2-6, to=1-7]
	\arrow["{m_{T^2A}}"', from=2-6, to=3-7]
	\arrow[""{name=1, anchor=center, inner sep=0}, "{m_A}"{description}, from=2-8, to=2-10]
	\arrow[""{name=2, anchor=center, inner sep=0}, "{m_{TA}}"', from=3-2, to=3-4]
	\arrow["{m_A}"', from=3-4, to=2-5]
	\arrow["{Tm_A}"', from=3-7, to=2-8]
	\arrow["{m_{TA}}"', from=3-7, to=3-9]
	\arrow["{m_A}"', from=3-9, to=2-10]
	\arrow["{T\mu_A}", shorten >=3pt, Rightarrow, from=1-2, to=0]
	\arrow["{\mu_A}"', shorten >=3pt, Rightarrow, from=1-9, to=1]
	\arrow["{\mu_{TA}}"', shift right=5, shorten >=3pt, Rightarrow, from=2-3, to=2]
	\arrow["{\mu_A}"', shorten <=3pt, Rightarrow, from=1, to=3-9]
\end{tikzcd}\]
and the following pairs of 2-cells are equal:
\[\begin{tikzcd}[column sep=scriptsize]
	&& {T^2A} &&& {T^3A} \\
	{T^2A} & {T^3A} && TA & {T^2A} && {T^2A} & TA \\
	&& {T^2A} &&& {T^3A}
	\arrow["{m_A}", from=1-3, to=2-4]
	\arrow["{\mu_A}", shorten <=6pt, shorten >=6pt, Rightarrow, from=1-3, to=3-3]
	\arrow["{Tm_A}", from=1-6, to=2-7]
	\arrow["{Ti_{TA}}", from=2-1, to=2-2]
	\arrow["{Tm_A}", from=2-2, to=1-3]
	\arrow["{m_{TA}}"', from=2-2, to=3-3]
	\arrow["{Ti_{TA}}", from=2-5, to=1-6]
	\arrow[""{name=0, anchor=center, inner sep=0}, Rightarrow, no head, from=2-5, to=2-7]
	\arrow["{Ti_{TA}}"', from=2-5, to=3-6]
	\arrow["{m_A}", from=2-7, to=2-8]
	\arrow["{m_A}"', from=3-3, to=2-4]
	\arrow["{m_{TA}}"', from=3-6, to=2-7]
	\arrow["{T\beta_A}", shorten >=3pt, Rightarrow, from=1-6, to=0]
	\arrow["{\eta_{TA}}", shorten <=3pt, Rightarrow, from=0, to=3-6]
\end{tikzcd}\]

\end{defi}

\begin{vari}
In case $T$ is a 2-functor, $m,i$ are 2-natural, and the modifications $\mu, \beta, \eta$ are the identities, this is called a \textit{2-monad}.
\end{vari}

For simplicity, we will usually denote the pseudomonad just by a triple $(T,m,i)$, with symbols $\mu, \eta, \beta$ always being used for its modifications.

\begin{lemma}\label{THML_in_every_psmonad}
For any pseudomonad $(T,m,i)$ on a 2-category $\ck$, the following 2-cells are equal:
\[\begin{tikzcd}[column sep=scriptsize]
	& TA &&&&& {T^2A} \\
	A && {T^2A} & TA & A & TA && TA \\
	& TA &&&&& {T^2A}
	\arrow["{i_{TA}}", from=1-2, to=2-3]
	\arrow["{i_{i_A}}"', shorten <=6pt, shorten >=6pt, Rightarrow, from=1-2, to=3-2]
	\arrow["{m_A}", from=1-7, to=2-8]
	\arrow["{i_A}", from=2-1, to=1-2]
	\arrow["{i_A}"', from=2-1, to=3-2]
	\arrow["{m_A}", from=2-3, to=2-4]
	\arrow["{i_A}", from=2-5, to=2-6]
	\arrow["{i_{TA}}", from=2-6, to=1-7]
	\arrow[""{name=0, anchor=center, inner sep=0}, Rightarrow, no head, from=2-6, to=2-8]
	\arrow["{Ti_A}"', from=2-6, to=3-7]
	\arrow["{Ti_A}"', from=3-2, to=2-3]
	\arrow["{m_A}"', from=3-7, to=2-8]
	\arrow["{\beta_A}"', shorten >=3pt, Rightarrow, from=1-7, to=0]
	\arrow["{\eta_A}"', shorten <=3pt, Rightarrow, from=0, to=3-7]
\end{tikzcd}\]
also, the following 2-cells are equal:
\[\begin{tikzcd}[column sep=tiny]
	{T^2A} && {T^2A} && TA & {T^2A} && TA && TA \\
	& {T^3A} && {T^2A} &&& {T^3A} && {T^2A}
	\arrow[""{name=0, anchor=center, inner sep=0}, Rightarrow, no head, from=1-1, to=1-3]
	\arrow["{T^2i_A}"', from=1-1, to=2-2]
	\arrow["{m_A}", from=1-3, to=1-5]
	\arrow["{m_A}", from=1-6, to=1-8]
	\arrow["{T^2i_A}"', from=1-6, to=2-7]
	\arrow[""{name=1, anchor=center, inner sep=0}, Rightarrow, no head, from=1-8, to=1-10]
	\arrow["{Ti_A}"', from=1-8, to=2-9]
	\arrow["{Tm_A}"', from=2-2, to=1-3]
	\arrow[""{name=2, anchor=center, inner sep=0}, "{m_{TA}}"', from=2-2, to=2-4]
	\arrow["{m_A}"', from=2-4, to=1-5]
	\arrow[""{name=3, anchor=center, inner sep=0}, "{m_{TA}}"', from=2-7, to=2-9]
	\arrow["{m_A}"', from=2-9, to=1-10]
	\arrow["{T\eta_A}"', shorten <=3pt, Rightarrow, from=0, to=2-2]
	\arrow["{\mu_A}", shift left=5, shorten >=3pt, Rightarrow, from=1-3, to=2]
	\arrow["{m^{-1}_{i_A}}"', shift right=5, shorten >=3pt, Rightarrow, from=1-8, to=3]
	\arrow["{\eta_A}", shorten <=3pt, Rightarrow, from=1, to=2-9]
\end{tikzcd}\]
\end{lemma}
\begin{proof}
This is proven in \cite[Proposition 8.1]{doctrinesffadjstring}.
\end{proof}

\begin{defi}
A \textit{colax $T$-algebra} $(A,a,\gamma,\iota)$ consists of an object $A$, 1-cell $a : TA \to A$, 2-cells $\gamma: am_A \Rightarrow aTa$, $\iota : ai_A \Rightarrow 1$ (called the \textit{coassociator} and the \textit{counitor}) such that the following 2-cells are equal:
\[\begin{tikzcd}[column sep=small]
	& {T^2A} && TA &&& {T^2A} && TA \\
	{T^3A} && {T^2A} && A & {T^3A} && TA && A \\
	& {T^2A} && TA &&& {T^2A} && TA
	\arrow[""{name=0, anchor=center, inner sep=0}, "{m_A}", from=1-2, to=1-4]
	\arrow["a", from=1-4, to=2-5]
	\arrow["\gamma", shorten <=6pt, shorten >=6pt, Rightarrow, from=1-4, to=3-4]
	\arrow["{m_A}", from=1-7, to=1-9]
	\arrow["Ta", from=1-7, to=2-8]
	\arrow["{m_{a}^{-1}}"', shorten <=6pt, shorten >=6pt, Rightarrow, from=1-7, to=3-7]
	\arrow["a", from=1-9, to=2-10]
	\arrow["{m_{TA}}", from=2-1, to=1-2]
	\arrow[""{name=1, anchor=center, inner sep=0}, "{Tm_A}", from=2-1, to=2-3]
	\arrow["{T^2a}"', from=2-1, to=3-2]
	\arrow["{m_A}"{description}, from=2-3, to=1-4]
	\arrow["Ta"{description}, from=2-3, to=3-4]
	\arrow["{m_{TA}}", from=2-6, to=1-7]
	\arrow["{T^2a}"', from=2-6, to=3-7]
	\arrow[""{name=2, anchor=center, inner sep=0}, "a"{description}, from=2-8, to=2-10]
	\arrow["Ta"', from=3-2, to=3-4]
	\arrow["a"', from=3-4, to=2-5]
	\arrow["{m_A}", from=3-7, to=2-8]
	\arrow[""{name=3, anchor=center, inner sep=0}, "Ta"', from=3-7, to=3-9]
	\arrow["a"', from=3-9, to=2-10]
	\arrow["{\mu_{A}^{-1}}"', shorten <=3pt, Rightarrow, from=0, to=2-3]
	\arrow["\gamma", shift right=5, shorten <=3pt, shorten >=3pt, Rightarrow, from=1-9, to=2]
	\arrow["{T\gamma}", shift right=4, shorten <=3pt, shorten >=3pt, Rightarrow, from=1, to=3-2]
	\arrow["\gamma", shorten <=3pt, shorten >=3pt, Rightarrow, from=2-8, to=3]
\end{tikzcd}\]
and also satisfying that both of these 2-cells are equal to the identity on $a$:
\[\begin{tikzcd}[column sep=scriptsize]
	&&&&& TA &&&& TA \\
	&&& TA &&&& {T^2A} \\
	TA && {T^2A} && A &&& TA \\
	&&& TA && A &&&& A
	\arrow[""{name=0, anchor=center, inner sep=0}, Rightarrow, no head, from=1-6, to=1-10]
	\arrow[""{name=1, anchor=center, inner sep=0}, "{i_{TA}}"{description}, from=1-6, to=2-8]
	\arrow["a"', from=1-6, to=4-6]
	\arrow["a", from=1-10, to=4-10]
	\arrow["a", from=2-4, to=3-5]
	\arrow["\gamma", shorten <=6pt, shorten >=6pt, Rightarrow, from=2-4, to=4-4]
	\arrow[""{name=2, anchor=center, inner sep=0}, "{m_A}"{description}, from=2-8, to=1-10]
	\arrow["Ta"', from=2-8, to=3-8]
	\arrow[""{name=3, anchor=center, inner sep=0}, curve={height=-24pt}, Rightarrow, no head, from=3-1, to=2-4]
	\arrow["{Ti_A}"{description}, from=3-1, to=3-3]
	\arrow[""{name=4, anchor=center, inner sep=0}, curve={height=24pt}, Rightarrow, no head, from=3-1, to=4-4]
	\arrow["{m_A}"{description}, from=3-3, to=2-4]
	\arrow["Ta"{description}, from=3-3, to=4-4]
	\arrow[""{name=5, anchor=center, inner sep=0}, "a"{description}, from=3-8, to=4-10]
	\arrow["a"', from=4-4, to=3-5]
	\arrow[""{name=6, anchor=center, inner sep=0}, "{i_A}"{description}, from=4-6, to=3-8]
	\arrow[""{name=7, anchor=center, inner sep=0}, Rightarrow, no head, from=4-6, to=4-10]
	\arrow["{\beta^{-1}_A}", shorten <=3pt, Rightarrow, from=0, to=2-8]
	\arrow["{i_a^{-1}}"', shorten <=9pt, shorten >=9pt, Rightarrow, from=1, to=6]
	\arrow["\gamma", shorten <=9pt, shorten >=9pt, Rightarrow, from=2, to=5]
	\arrow["{\eta_A}"', shorten <=5pt, Rightarrow, from=3, to=3-3]
	\arrow["{T\iota}"', shorten >=5pt, Rightarrow, from=3-3, to=4]
	\arrow["\iota", shorten >=3pt, Rightarrow, from=3-8, to=7]
\end{tikzcd}\]
\end{defi}

\begin{vari}
We will encounter the following variations:
\begin{itemize}
\item If $\iota$ is invertible, this is called a \textit{normal} colax $T$-algebra,
\item if $\gamma, \iota$ go in the other direction, this is called a \textit{lax $T$-algebra},
\item if $\gamma, \iota$ are invertible, this is called a \textit{pseudo $T$-algebra},
\item if $\gamma, \iota$ are the identities, this is called a \textit{strict $T$-algebra}\footnote{In case $T$ is a 2-monad, by a \textit{$T$-algebra} we will always mean the strict one.}.
\end{itemize}
\end{vari}

\begin{pr}
For every object $A \in \ck$, there is a \textit{free $T$-algebra} on $A$ that is given by the tuple $(TA, m_A, \mu^{-1}_A, \beta_A)$. For pseudomonads, this is a pseudo algebra. The algebra associativity axiom is given by the pseudomonad axiom. The first algebra unit axiom is given by the second pseudomonad axiom, and the second algebra unit axiom is proven in Lemma \ref{THML_in_every_psmonad}. 
\end{pr}

\begin{defi}
A \textit{lax morphism} $(f,\overline{f}): (A,a,\gamma,\iota) \rightsquigarrow (B,b,\gamma,\iota)$ between two colax $T$-algebras consists of a 1-cell $f: A \to B$ and a 2-cell $\overline{f}: b Tf \Rightarrow f a$ (the \textit{associator}) with the property that the following pair of 2-cells is equal:
\begin{equation}\label{EQD_laxmor_assoc}
\begin{tikzcd}[column sep=scriptsize]
	{T^2A} && {T^2B} && {T^2A} && {T^2B} \\
	& TA && TB &&&& TB \\
	TA &&&& TA && TB \\
	& A && B && A && B
	\arrow[""{name=0, anchor=center, inner sep=0}, "{T^2f}", from=1-1, to=1-3]
	\arrow[""{name=0p, anchor=center, inner sep=0}, phantom, from=1-1, to=1-3, start anchor=center, end anchor=center]
	\arrow[""{name=1, anchor=center, inner sep=0}, "{m_A}"{description}, from=1-1, to=2-2]
	\arrow[""{name=1p, anchor=center, inner sep=0}, phantom, from=1-1, to=2-2, start anchor=center, end anchor=center]
	\arrow["Ta"', from=1-1, to=3-1]
	\arrow["{m_B}", from=1-3, to=2-4]
	\arrow[""{name=2, anchor=center, inner sep=0}, "{T^2f}", from=1-5, to=1-7]
	\arrow[""{name=2p, anchor=center, inner sep=0}, phantom, from=1-5, to=1-7, start anchor=center, end anchor=center]
	\arrow["Ta"', from=1-5, to=3-5]
	\arrow[""{name=3, anchor=center, inner sep=0}, "{m_B}", from=1-7, to=2-8]
	\arrow[""{name=3p, anchor=center, inner sep=0}, phantom, from=1-7, to=2-8, start anchor=center, end anchor=center]
	\arrow["Tb"', from=1-7, to=3-7]
	\arrow[""{name=4, anchor=center, inner sep=0}, "Tf"{description}, from=2-2, to=2-4]
	\arrow[""{name=4p, anchor=center, inner sep=0}, phantom, from=2-2, to=2-4, start anchor=center, end anchor=center]
	\arrow["a"{description}, from=2-2, to=4-2]
	\arrow["b", from=2-4, to=4-4]
	\arrow["b", from=2-8, to=4-8]
	\arrow[""{name=5, anchor=center, inner sep=0}, "a"', from=3-1, to=4-2]
	\arrow[""{name=5p, anchor=center, inner sep=0}, phantom, from=3-1, to=4-2, start anchor=center, end anchor=center]
	\arrow[""{name=6, anchor=center, inner sep=0}, "Tf"{description}, from=3-5, to=3-7]
	\arrow[""{name=6p, anchor=center, inner sep=0}, phantom, from=3-5, to=3-7, start anchor=center, end anchor=center]
	\arrow["a"', from=3-5, to=4-6]
	\arrow[""{name=7, anchor=center, inner sep=0}, "b"{description}, from=3-7, to=4-8]
	\arrow[""{name=7p, anchor=center, inner sep=0}, phantom, from=3-7, to=4-8, start anchor=center, end anchor=center]
	\arrow[""{name=8, anchor=center, inner sep=0}, "f"', from=4-2, to=4-4]
	\arrow[""{name=8p, anchor=center, inner sep=0}, phantom, from=4-2, to=4-4, start anchor=center, end anchor=center]
	\arrow["f"', from=4-6, to=4-8]
	\arrow["{m_f}", shorten <=3pt, Rightarrow, from=0p, to=2-2]
	\arrow["\gamma", shorten <=9pt, shorten >=9pt, Rightarrow, from=1p, to=5p]
	\arrow["{T\overline{f}}", shorten <=9pt, shorten >=9pt, Rightarrow, from=2p, to=6p]
	\arrow["{\gamma'}", shorten <=9pt, shorten >=9pt, Rightarrow, from=3p, to=7p]
	\arrow["{\overline{f}}", shorten <=9pt, shorten >=9pt, Rightarrow, from=4p, to=8p]
	\arrow["{\overline{f}}", shorten <=6pt, Rightarrow, from=6, to=4-6]
\end{tikzcd}
\end{equation}
and also the following pair of 2-cells is equal:
\begin{equation}\label{EQD_laxmor_unit}
\begin{tikzcd}
	A && B && A && B \\
	& TA && TB &&&& TB \\
	\\
	& A && B && A && B
	\arrow["f", from=1-1, to=1-3]
	\arrow["{i_A}", from=1-1, to=2-2]
	\arrow[""{name=0, anchor=center, inner sep=0}, curve={height=24pt}, Rightarrow, no head, from=1-1, to=4-2]
	\arrow["{i_B}", from=1-3, to=2-4]
	\arrow["f", from=1-5, to=1-7]
	\arrow[curve={height=24pt}, Rightarrow, no head, from=1-5, to=4-6]
	\arrow["{i_B}", from=1-7, to=2-8]
	\arrow[""{name=1, anchor=center, inner sep=0}, curve={height=24pt}, Rightarrow, no head, from=1-7, to=4-8]
	\arrow[""{name=2, anchor=center, inner sep=0}, "Tf"{description}, from=2-2, to=2-4]
	\arrow["a"', from=2-2, to=4-2]
	\arrow["b", from=2-4, to=4-4]
	\arrow["b", from=2-8, to=4-8]
	\arrow[""{name=3, anchor=center, inner sep=0}, "f"', from=4-2, to=4-4]
	\arrow["f"', from=4-6, to=4-8]
	\arrow["{i_f}"', shorten <=3pt, shorten >=3pt, Rightarrow, from=1-3, to=2]
	\arrow["\iota"', shorten >=6pt, Rightarrow, from=2-2, to=0]
	\arrow["{\overline{f}}", shorten <=9pt, shorten >=9pt, Rightarrow, from=2, to=3]
	\arrow["{\iota'}"', shorten >=6pt, Rightarrow, from=2-8, to=1]
\end{tikzcd}
\end{equation}
\end{defi}

\begin{vari}
We will encounter the following variations:
\begin{itemize}
\item If $\overline{f}$ goes in the other direction, this is called a \textit{colax morphism},
\item if $\overline{f}$ is invertible, this is called a \textit{pseudo morphism},
\item if $\overline{f}$ is the identity, this is called a \textit{strict morphism}.
\end{itemize}
\end{vari}

\begin{defi}
An \textit{algebra 2-cell} $\alpha: (f,\overline{f}) \Rightarrow (g,\overline{g}): (A,a,\gamma, \iota) \to (B,b,\gamma',\iota')$ between two lax morphisms of colax $T$-algebras is a 2-cell $\alpha: f \Rightarrow g$ in $\ck$ that satisfies:
\[\begin{tikzcd}
	TA && TB && TA && TB \\
	&&& {=} \\
	A && B && A && B
	\arrow[""{name=0, anchor=center, inner sep=0}, "Tg"', from=1-1, to=1-3]
	\arrow[""{name=1, anchor=center, inner sep=0}, "Tf", curve={height=-24pt}, from=1-1, to=1-3]
	\arrow["a"', from=1-1, to=3-1]
	\arrow["b", from=1-3, to=3-3]
	\arrow[""{name=2, anchor=center, inner sep=0}, "Tf", from=1-5, to=1-7]
	\arrow["a"', from=1-5, to=3-5]
	\arrow["b", from=1-7, to=3-7]
	\arrow[""{name=3, anchor=center, inner sep=0}, "g"', from=3-1, to=3-3]
	\arrow[""{name=4, anchor=center, inner sep=0}, "f", from=3-5, to=3-7]
	\arrow[""{name=5, anchor=center, inner sep=0}, "g"', curve={height=24pt}, from=3-5, to=3-7]
	\arrow["T\alpha", shorten <=3pt, shorten >=3pt, Rightarrow, from=1, to=0]
	\arrow["{\overline{g}}", shorten <=13pt, shorten >=9pt, Rightarrow, from=0, to=3]
	\arrow["{\overline{f}}", shorten <=9pt, shorten >=13pt, Rightarrow, from=2, to=4]
	\arrow["\alpha", shorten <=3pt, shorten >=3pt, Rightarrow, from=4, to=5]
\end{tikzcd}\]
\end{defi}

\begin{nota}\label{NOTA_categories_of_algebras}
Given a 2-category $\ck$ and a pseudomonad $(T,m,i)$ on it, there are many notions of algebras and morphisms and thus many different 2-categories of algebras. We will denote by:
\[
\talg_x, \pstalg_x, \colaxtalg_x, \ncolaxtalg_x
\]
the 2-categories of strict, pseudo, colax and normal colax algebras for the pseudomonad. The index $x$ will denote the kind of algebra morphisms we consider: $x=s$ will denote strict morphisms, $x=l$ or $c$ will denote lax/colax morphisms. We follow the convention that if we omit the index $x$, we always mean pseudo morphisms.

\end{nota}

\begin{prop}\label{THMP_adj_colax_alg_base_2cat}
Let $(T,m,i,\mu,\eta,\beta)$ be a pseudomonad on a 2-category $\ck$. There is a ``free-forgetful” colax adjunction between lax algebras and the base 2-category:
\[\begin{tikzcd}
	{(\eta^{-1}, \Phi): (\epsilon, i):} & \laxtalgl && \ck
	\arrow[""{name=0, anchor=center, inner sep=0}, "U"', curve={height=18pt}, from=1-2, to=1-4]
	\arrow[""{name=1, anchor=center, inner sep=0}, "F"', curve={height=18pt}, from=1-4, to=1-2]
	\arrow["\tbot"{description}, draw=none, from=1, to=0]
\end{tikzcd}\]
\end{prop}
\begin{proof}
Let us first define the \textbf{2-functors}. The forgetful functor $U: \laxtalgl \to \ck$ sends algebras to their underlying objects. The 2-functor $F$ sends $A \in \ck$ to the free algebra $(TA, m_A, \mu_A, \beta_A^{-1})$ (regarded as a lax algebra this time) and it sends a 1-cell $f: A \to B$ to the pseudo morphism $(Tf,m_f)$ -- the pseudo morphism axioms follow from the modification axioms of $\mu: m \circ mT \to m \circ Tm$ and $\beta : m \circ iT \to 1_{T}$. Finally, a 2-cell $\alpha$ is sent to $T\alpha$ -- the algebra 2-cell axiom follows from the local naturality of $m: T^2 \Rightarrow T$.

The \textbf{unit} is the unit of the pseudomonad. The \textbf{counit} $\epsilon: FU \Rightarrow 1_{\laxtalgl}$ evaluated at a lax $T$-algebra $(A,a,\gamma,\iota)$ is the lax morphism:
\[
(a, \gamma) : (TA,m_A) \to (A,a,\gamma,\iota).
\]
The counit is colax natural, its 2-cell component at a lax morphism:
\[
(f,\overline{f}): (A,a,\gamma, \iota) \to (B,b,\gamma',\iota'),
\]
is given by the 2-cell $\overline{f}$:
\[\begin{tikzcd}
	{(TA,m_A)} && {(A,a,\gamma,\iota)} \\
	&&& {\text{ in } \laxtalgl} \\
	{(TB,m_B)} && {(B,b,\gamma',\iota')}
	\arrow[""{name=0, anchor=center, inner sep=0}, "{(a,\gamma)}", from=1-1, to=1-3]
	\arrow[""{name=0p, anchor=center, inner sep=0}, phantom, from=1-1, to=1-3, start anchor=center, end anchor=center]
	\arrow["Tf"', from=1-1, to=3-1]
	\arrow["{(f,\overline{f})}", from=1-3, to=3-3]
	\arrow[""{name=1, anchor=center, inner sep=0}, "{(b,\gamma')}"', from=3-1, to=3-3]
	\arrow[""{name=1p, anchor=center, inner sep=0}, phantom, from=3-1, to=3-3, start anchor=center, end anchor=center]
	\arrow["{\overline{f}}"', shorten <=9pt, shorten >=9pt, Rightarrow, from=1p, to=0p]
\end{tikzcd}\]

The \textbf{first modification} $\eta^{-1}: \epsilon F \circ F \eta \to 1_F$ evaluated at an object $B \in \ck$ is the pseudomonad 2-cell $\eta^{-1}_B: m_B Ti_B \cong 1_{TB}$. The \textbf{second modification} $\Phi: 1_U \to U\epsilon \circ \eta U$ evaluated at a lax $T$-algebra $(A,a,\gamma,\iota)$ is given by the 2-cell:
\[
\Phi_{(A,a,\gamma,\iota)} := \iota : 1_A \Rightarrow ai_A.
\]
The 2-cells in the swallowtail identities evaluated at $A \in \ck$ (on the left) and\\$(A,a,\gamma,\iota) \in \laxtalgl$ (on the right) become the following:
\[\begin{tikzcd}[column sep=scriptsize]
	A && TA &&& TA \\
	&&& {} &&&& {} & {} \\
	TA && {T^2A} & {} &&& {} & {T^2A} && TA \\
	& {} & {} &&&& {} \\
	&&&& TA &&& TA && {1_{\cd}}
	\arrow["{i_A}", from=1-1, to=1-3]
	\arrow["{i_A}"', from=1-1, to=3-1]
	\arrow["{Ti_A}", from=1-3, to=3-3]
	\arrow[curve={height=-24pt}, equals, from=1-3, to=5-5]
	\arrow["{Ti_A}"{description}, from=1-6, to=3-8]
	\arrow[curve={height=-24pt}, equals, from=1-6, to=3-10]
	\arrow[curve={height=24pt}, equals, from=1-6, to=5-8]
	\arrow[""{name=0, anchor=center, inner sep=0}, draw=none, from=2-8, to=2-9]
	\arrow[""{name=0p, anchor=center, inner sep=0}, phantom, from=2-8, to=2-9, start anchor=center, end anchor=center]
	\arrow["{i_{i_A}}", shorten <=11pt, shorten >=11pt, Rightarrow, from=3-1, to=1-3]
	\arrow["{i_{TA}}"', from=3-1, to=3-3]
	\arrow[curve={height=24pt}, equals, from=3-1, to=5-5]
	\arrow[""{name=1, anchor=center, inner sep=0}, draw=none, from=3-3, to=3-4]
	\arrow[""{name=1p, anchor=center, inner sep=0}, phantom, from=3-3, to=3-4, start anchor=center, end anchor=center]
	\arrow["{m_A}"{description}, from=3-3, to=5-5]
	\arrow[""{name=2, anchor=center, inner sep=0}, shift left=3, draw=none, from=3-7, to=3-8]
	\arrow[""{name=2p, anchor=center, inner sep=0}, phantom, from=3-7, to=3-8, start anchor=center, end anchor=center, shift left=3]
	\arrow["Ta", from=3-8, to=3-10]
	\arrow["{m_A}"', from=3-8, to=5-8]
	\arrow["\gamma", shorten <=11pt, shorten >=11pt, Rightarrow, from=3-10, to=5-8]
	\arrow["a", from=3-10, to=5-10]
	\arrow[""{name=3, anchor=center, inner sep=0}, draw=none, from=4-2, to=4-3]
	\arrow[""{name=3p, anchor=center, inner sep=0}, phantom, from=4-2, to=4-3, start anchor=center, end anchor=center]
	\arrow["a"', from=5-8, to=5-10]
	\arrow["{T\iota}"', shorten <=4pt, shorten >=2pt, Rightarrow, from=0p, to=3-8]
	\arrow["{\eta^{-1}_A}"', shorten >=9pt, Rightarrow, from=1p, to=2-4]
	\arrow["{\eta^{-1}_{A}}"', shorten <=5pt, shorten >=9pt, Rightarrow, from=2p, to=4-7]
	\arrow["{\beta_A}"', shorten >=4pt, Rightarrow, from=3p, to=3-3]
\end{tikzcd}\]
The 2-cell above left is the identity because of the equality proven in Lemma \ref{THML_in_every_psmonad}. The 2-cell above right is the identity because of the lax algebra unit axiom.
\end{proof}

\begin{pozn}\label{POZN_homcat_adj_laxtalgl_base}
If we restrict the above colax adjunction to pseudo algebras and pseudo morphisms, we see that we obtain a biadjunction. For a 2-monad, restricting to strict algebras and strict morphisms gives a 2-adjunction. Note also that by Proposition \ref{THMP_colax_adj_adjoint_homcats}, for each $A \in \ck$ and a lax $T$-algebra $\mathbb{B} := (B,b,\gamma,\iota)$ we have an adjunction:
\[\begin{tikzcd}
	{\laxtalgl(FA,\mathbb{B})} && {\ck(A,B)}
	\arrow[""{name=0, anchor=center, inner sep=0}, "{(i_A)^* \circ U}", curve={height=-24pt}, from=1-1, to=1-3]
	\arrow[""{name=1, anchor=center, inner sep=0}, "{(\epsilon_{\mathbb{B}})_* \circ F}", curve={height=-24pt}, from=1-3, to=1-1]
	\arrow["\dashv"{anchor=center, rotate=90}, draw=none, from=1, to=0]
\end{tikzcd}\]
The counit uses the modification $\eta^{-1}$ and the 2-cell components of morphisms in the 2-category $\laxtalgl$. The unit uses the coassociator for $\mathbb{B}$ and the pseudonaturality squares for $i: 1_{\ck} \Rightarrow T$. So in case $\mathbb{B}$ is a pseudo algebra, the above becomes an adjoint equivalence if we restrict to pseudo morphisms.
\end{pozn}

\subsection{Doctrinal adjunction}\label{SUBS_psmonads_doctrinal}

The following theorem is perhaps the most important formal result in the two-dimensional monad theory. It has first been observed in \cite[Theorem 1.2]{doctrinaladj} for the case where $T$ is a 2-monad. The name comes from the fact that some authors used to call 2-monads \textit{doctrines}.

\begin{theo}[Doctrinal adjunction]\label{THM_doctrinal_adj}
Let $(A,a,\gamma,\iota)$, $(B,b,\gamma',\iota')$ be two colax $T$-algebras for a pseudomonad $(T,m,i)$ on a 2-category $\ck$. Assume that we have an adjunction in $\ck$:
\[\begin{tikzcd}
	{(\epsilon,\eta):} & B && A
	\arrow[""{name=0, anchor=center, inner sep=0}, "f"', curve={height=24pt}, from=1-2, to=1-4]
	\arrow[""{name=1, anchor=center, inner sep=0}, "u"', curve={height=24pt}, from=1-4, to=1-2]
	\arrow["\dashv"{anchor=center, rotate=90}, draw=none, from=0, to=1]
\end{tikzcd}\]
Then the 2-cell below gives a lax morphism $(A,a,\gamma,\iota) \to (B,b,\gamma',\iota')$:
\[\begin{tikzcd}
	TA && TB \\
	\\
	A && B
	\arrow[""{name=0, anchor=center, inner sep=0}, "Tu", from=1-1, to=1-3]
	\arrow["a"', from=1-1, to=3-1]
	\arrow["b", from=1-3, to=3-3]
	\arrow[""{name=1, anchor=center, inner sep=0}, "u"', from=3-1, to=3-3]
	\arrow["{\overline{u}}", shorten <=9pt, shorten >=9pt, Rightarrow, from=0, to=1]
\end{tikzcd}\]
if and only if its \textit{mate}\footnote{See \cite[Page 87]{review}.} under adjunctions $f \dashv u$, $Tf \dashv Tu$ gives a colax morphism $(B,b,\gamma',\iota') \to (A,a,\gamma,\iota)$:
\[\begin{tikzcd}
	& TB \\
	\\
	& TA && TB \\
	{\overline{f}:=} \\
	& A && B \\
	\\
	&&& A
	\arrow["Tf"', from=1-2, to=3-2]
	\arrow[""{name=0, anchor=center, inner sep=0}, Rightarrow, no head, from=1-2, to=3-4]
	\arrow[""{name=0p, anchor=center, inner sep=0}, phantom, from=1-2, to=3-4, start anchor=center, end anchor=center]
	\arrow[""{name=1, anchor=center, inner sep=0}, "Tu"{description}, from=3-2, to=3-4]
	\arrow[""{name=1p, anchor=center, inner sep=0}, phantom, from=3-2, to=3-4, start anchor=center, end anchor=center]
	\arrow[""{name=1p, anchor=center, inner sep=0}, phantom, from=3-2, to=3-4, start anchor=center, end anchor=center]
	\arrow["a"', from=3-2, to=5-2]
	\arrow["b", from=3-4, to=5-4]
	\arrow[""{name=2, anchor=center, inner sep=0}, "u"{description}, from=5-2, to=5-4]
	\arrow[""{name=2p, anchor=center, inner sep=0}, phantom, from=5-2, to=5-4, start anchor=center, end anchor=center]
	\arrow[""{name=2p, anchor=center, inner sep=0}, phantom, from=5-2, to=5-4, start anchor=center, end anchor=center]
	\arrow[""{name=3, anchor=center, inner sep=0}, Rightarrow, no head, from=5-2, to=7-4]
	\arrow[""{name=3p, anchor=center, inner sep=0}, phantom, from=5-2, to=7-4, start anchor=center, end anchor=center]
	\arrow["f", from=5-4, to=7-4]
	\arrow["{T\eta}"', shorten <=4pt, shorten >=4pt, Rightarrow, from=0p, to=1p]
	\arrow["{\overline{u}}", shorten <=9pt, shorten >=9pt, Rightarrow, from=1p, to=2p]
	\arrow["\epsilon", shorten <=4pt, shorten >=4pt, Rightarrow, from=2p, to=3p]
\end{tikzcd}\]
\end{theo}
\begin{proof}
Notice that for a lax morphism $(u,\overline{u})$, the associativity axiom holds:
\[\begin{tikzcd}[column sep=scriptsize]
	{T^2A} && {T^2B} && {T^2A} && {T^2B} \\
	& TA && TB &&&& TB \\
	TA &&&& TA && TB \\
	& A && B && A && B
	\arrow[""{name=0, anchor=center, inner sep=0}, "{T^2u}", from=1-1, to=1-3]
	\arrow[""{name=0p, anchor=center, inner sep=0}, phantom, from=1-1, to=1-3, start anchor=center, end anchor=center]
	\arrow[""{name=1, anchor=center, inner sep=0}, "{m_A}"{description}, from=1-1, to=2-2]
	\arrow[""{name=1p, anchor=center, inner sep=0}, phantom, from=1-1, to=2-2, start anchor=center, end anchor=center]
	\arrow["Ta"', from=1-1, to=3-1]
	\arrow["{m_B}", from=1-3, to=2-4]
	\arrow[""{name=2, anchor=center, inner sep=0}, "{T^2u}", from=1-5, to=1-7]
	\arrow[""{name=2p, anchor=center, inner sep=0}, phantom, from=1-5, to=1-7, start anchor=center, end anchor=center]
	\arrow["Ta"', from=1-5, to=3-5]
	\arrow[""{name=3, anchor=center, inner sep=0}, "{m_B}", from=1-7, to=2-8]
	\arrow[""{name=3p, anchor=center, inner sep=0}, phantom, from=1-7, to=2-8, start anchor=center, end anchor=center]
	\arrow["Tb"', from=1-7, to=3-7]
	\arrow[""{name=4, anchor=center, inner sep=0}, "Tu"{description}, from=2-2, to=2-4]
	\arrow[""{name=4p, anchor=center, inner sep=0}, phantom, from=2-2, to=2-4, start anchor=center, end anchor=center]
	\arrow["a"{description}, from=2-2, to=4-2]
	\arrow["b", from=2-4, to=4-4]
	\arrow["b", from=2-8, to=4-8]
	\arrow[""{name=5, anchor=center, inner sep=0}, "a"', from=3-1, to=4-2]
	\arrow[""{name=5p, anchor=center, inner sep=0}, phantom, from=3-1, to=4-2, start anchor=center, end anchor=center]
	\arrow[""{name=6, anchor=center, inner sep=0}, "Tu"{description}, from=3-5, to=3-7]
	\arrow[""{name=6p, anchor=center, inner sep=0}, phantom, from=3-5, to=3-7, start anchor=center, end anchor=center]
	\arrow["a"', from=3-5, to=4-6]
	\arrow[""{name=7, anchor=center, inner sep=0}, "b"{description}, from=3-7, to=4-8]
	\arrow[""{name=7p, anchor=center, inner sep=0}, phantom, from=3-7, to=4-8, start anchor=center, end anchor=center]
	\arrow[""{name=8, anchor=center, inner sep=0}, "u"', from=4-2, to=4-4]
	\arrow[""{name=8p, anchor=center, inner sep=0}, phantom, from=4-2, to=4-4, start anchor=center, end anchor=center]
	\arrow["u"', from=4-6, to=4-8]
	\arrow["{m_u}", shorten <=3pt, Rightarrow, from=0p, to=2-2]
	\arrow["\gamma", shorten <=9pt, shorten >=9pt, Rightarrow, from=1p, to=5p]
	\arrow["{T\overline{u}}", shorten <=9pt, shorten >=9pt, Rightarrow, from=2p, to=6p]
	\arrow["{\gamma'}", shorten <=9pt, shorten >=9pt, Rightarrow, from=3p, to=7p]
	\arrow["{\overline{u}}", shorten <=9pt, shorten >=9pt, Rightarrow, from=4p, to=8p]
	\arrow["{\overline{u}}", shorten <=6pt, Rightarrow, from=6, to=4-6]
\end{tikzcd}\]
if and only if the mates of both sides under the adjunctions $f \dashv u$, $T^2f \dashv T^2u$ are equal:
\[\begin{tikzcd}[column sep=scriptsize]
	{T^2B} &&&& {T^2B} \\
	{T^2A} && {T^2B} && {T^2A} && {T^2B} \\
	& TA && TB &&&& TB \\
	TA &&&& TA && TB \\
	& A && B && A && B \\
	&&& A &&&& A
	\arrow["{T^2f}"', from=1-1, to=2-1]
	\arrow[""{name=0, anchor=center, inner sep=0}, Rightarrow, no head, from=1-1, to=2-3]
	\arrow[""{name=0p, anchor=center, inner sep=0}, phantom, from=1-1, to=2-3, start anchor=center, end anchor=center]
	\arrow["{T^2f}"', from=1-5, to=2-5]
	\arrow[""{name=1, anchor=center, inner sep=0}, Rightarrow, no head, from=1-5, to=2-7]
	\arrow[""{name=1p, anchor=center, inner sep=0}, phantom, from=1-5, to=2-7, start anchor=center, end anchor=center]
	\arrow[""{name=2, anchor=center, inner sep=0}, "{T^2u}"{description}, from=2-1, to=2-3]
	\arrow[""{name=2p, anchor=center, inner sep=0}, phantom, from=2-1, to=2-3, start anchor=center, end anchor=center]
	\arrow[""{name=2p, anchor=center, inner sep=0}, phantom, from=2-1, to=2-3, start anchor=center, end anchor=center]
	\arrow[""{name=3, anchor=center, inner sep=0}, "{m_A}"{description}, from=2-1, to=3-2]
	\arrow[""{name=3p, anchor=center, inner sep=0}, phantom, from=2-1, to=3-2, start anchor=center, end anchor=center]
	\arrow["Ta"', from=2-1, to=4-1]
	\arrow["{m_B}", from=2-3, to=3-4]
	\arrow[""{name=4, anchor=center, inner sep=0}, "{T^2u}"{description}, from=2-5, to=2-7]
	\arrow[""{name=4p, anchor=center, inner sep=0}, phantom, from=2-5, to=2-7, start anchor=center, end anchor=center]
	\arrow[""{name=4p, anchor=center, inner sep=0}, phantom, from=2-5, to=2-7, start anchor=center, end anchor=center]
	\arrow["Ta"', from=2-5, to=4-5]
	\arrow[""{name=5, anchor=center, inner sep=0}, "{m_B}", from=2-7, to=3-8]
	\arrow[""{name=5p, anchor=center, inner sep=0}, phantom, from=2-7, to=3-8, start anchor=center, end anchor=center]
	\arrow["Tb"', from=2-7, to=4-7]
	\arrow[""{name=6, anchor=center, inner sep=0}, "Tu"{description}, from=3-2, to=3-4]
	\arrow[""{name=6p, anchor=center, inner sep=0}, phantom, from=3-2, to=3-4, start anchor=center, end anchor=center]
	\arrow["a"{description}, from=3-2, to=5-2]
	\arrow["b", from=3-4, to=5-4]
	\arrow["b", from=3-8, to=5-8]
	\arrow[""{name=7, anchor=center, inner sep=0}, "a"', from=4-1, to=5-2]
	\arrow[""{name=7p, anchor=center, inner sep=0}, phantom, from=4-1, to=5-2, start anchor=center, end anchor=center]
	\arrow[""{name=8, anchor=center, inner sep=0}, "Tu"{description}, from=4-5, to=4-7]
	\arrow[""{name=8p, anchor=center, inner sep=0}, phantom, from=4-5, to=4-7, start anchor=center, end anchor=center]
	\arrow["a"', from=4-5, to=5-6]
	\arrow[""{name=9, anchor=center, inner sep=0}, "b"{description}, from=4-7, to=5-8]
	\arrow[""{name=9p, anchor=center, inner sep=0}, phantom, from=4-7, to=5-8, start anchor=center, end anchor=center]
	\arrow[""{name=10, anchor=center, inner sep=0}, "u"{description}, from=5-2, to=5-4]
	\arrow[""{name=10p, anchor=center, inner sep=0}, phantom, from=5-2, to=5-4, start anchor=center, end anchor=center]
	\arrow[""{name=10p, anchor=center, inner sep=0}, phantom, from=5-2, to=5-4, start anchor=center, end anchor=center]
	\arrow[""{name=11, anchor=center, inner sep=0}, Rightarrow, no head, from=5-2, to=6-4]
	\arrow[""{name=11p, anchor=center, inner sep=0}, phantom, from=5-2, to=6-4, start anchor=center, end anchor=center]
	\arrow["f", from=5-4, to=6-4]
	\arrow[""{name=12, anchor=center, inner sep=0}, "u"{description}, from=5-6, to=5-8]
	\arrow[""{name=12p, anchor=center, inner sep=0}, phantom, from=5-6, to=5-8, start anchor=center, end anchor=center]
	\arrow[""{name=13, anchor=center, inner sep=0}, Rightarrow, no head, from=5-6, to=6-8]
	\arrow[""{name=13p, anchor=center, inner sep=0}, phantom, from=5-6, to=6-8, start anchor=center, end anchor=center]
	\arrow["f", from=5-8, to=6-8]
	\arrow["{T^2\eta}"'{pos=0.4}, shift right=5, shorten <=2pt, shorten >=2pt, Rightarrow, from=0p, to=2p]
	\arrow["{T^2\eta}"'{pos=0.4}, shift right=5, shorten <=2pt, shorten >=2pt, Rightarrow, from=1p, to=4p]
	\arrow["{m_u}", shorten <=3pt, Rightarrow, from=2p, to=3-2]
	\arrow["\gamma", shorten <=9pt, shorten >=9pt, Rightarrow, from=3p, to=7p]
	\arrow["{T\overline{u}}", shorten <=9pt, shorten >=9pt, Rightarrow, from=4p, to=8p]
	\arrow["{\gamma'}", shorten <=9pt, shorten >=9pt, Rightarrow, from=5p, to=9p]
	\arrow["{\overline{u}}", shorten <=9pt, shorten >=9pt, Rightarrow, from=6p, to=10p]
	\arrow["{\overline{u}}", shorten <=6pt, Rightarrow, from=8, to=5-6]
	\arrow["\epsilon", shift left=5, shorten <=2pt, shorten >=2pt, Rightarrow, from=10p, to=11p]
	\arrow["\epsilon", shift left=5, shorten <=2pt, shorten >=2pt, Rightarrow, from=12p, to=13p]
\end{tikzcd}\]

Using the local naturality of $m: T^2 \Rightarrow T$, the triangle identity for $Tf \dashv Tu$ and the definition of $\overline{f}$, it can be seen that they are equal if and only if the following 2-cells are equal, which is the associativity axiom for the colax morphism $(f,\overline{f})$:
\begin{equation}\label{EQ_PF_doct_adj_assoc}
\begin{tikzcd}[column sep=scriptsize]
	{T^2A} && TA && {T^2A} && TA \\
	& TA && A &&&& A \\
	{T^2B} &&&& {T^2B} && TB \\
	& TB && B && TB && B
	\arrow["Ta", from=1-1, to=1-3]
	\arrow[""{name=0, anchor=center, inner sep=0}, "{m_A}", from=1-1, to=2-2]
	\arrow["a", from=1-3, to=2-4]
	\arrow[""{name=1, anchor=center, inner sep=0}, "Ta", from=1-5, to=1-7]
	\arrow[""{name=2, anchor=center, inner sep=0}, "a", from=1-7, to=2-8]
	\arrow[""{name=3, anchor=center, inner sep=0}, "a"', from=2-2, to=2-4]
	\arrow["{T^2f}", from=3-1, to=1-1]
	\arrow[""{name=4, anchor=center, inner sep=0}, "{m_B}"', from=3-1, to=4-2]
	\arrow["{T^2f}", from=3-5, to=1-5]
	\arrow[""{name=5, anchor=center, inner sep=0}, "Tb", from=3-5, to=3-7]
	\arrow["{m_B}"', from=3-5, to=4-6]
	\arrow["Tf"{description}, from=3-7, to=1-7]
	\arrow[""{name=6, anchor=center, inner sep=0}, "b"{description}, from=3-7, to=4-8]
	\arrow["Tf"{description}, from=4-2, to=2-2]
	\arrow[""{name=7, anchor=center, inner sep=0}, "b"', from=4-2, to=4-4]
	\arrow["f"', from=4-4, to=2-4]
	\arrow["b"', from=4-6, to=4-8]
	\arrow["f"', from=4-8, to=2-8]
	\arrow["\gamma", shorten <=3pt, Rightarrow, from=3, to=1-3]
	\arrow["{m_f^{-1}}", shorten <=9pt, shorten >=9pt, Rightarrow, from=4, to=0]
	\arrow["{T\overline{f}}"', shorten <=9pt, shorten >=9pt, Rightarrow, from=5, to=1]
	\arrow["{\overline{f}}"', shorten <=9pt, shorten >=9pt, Rightarrow, from=6, to=2]
	\arrow["{\overline{f}}"', shorten <=9pt, shorten >=9pt, Rightarrow, from=7, to=3]
	\arrow["{\gamma'}"', shorten <=3pt, shorten >=3pt, Rightarrow, from=4-6, to=5]
\end{tikzcd}
\end{equation}
Likewise for the morphism unit axioms for $(u,\overline{u})$, $(f,\overline{f})$.
\end{proof}

\begin{pozn}\label{POZN_mixed_algebra_2cells}
Note that if $(u,\overline{u})$ is a lax morphism and $\overline{f}$ its mate as in the above theorem, each of the 2-cells $\epsilon$, $\eta$ becomes a kind of a ``mixed algebra 2-cell” i.e. the following 2-cells are equal:
\[\begin{tikzcd}
	TB \\
	\\
	TA && TB & TB && B \\
	\\
	A && B & TA && A && B
	\arrow["Tf"', from=1-1, to=3-1]
	\arrow[""{name=0, anchor=center, inner sep=0}, Rightarrow, no head, from=1-1, to=3-3]
	\arrow[""{name=0p, anchor=center, inner sep=0}, phantom, from=1-1, to=3-3, start anchor=center, end anchor=center]
	\arrow[""{name=1, anchor=center, inner sep=0}, "Tu"{description}, from=3-1, to=3-3]
	\arrow[""{name=1p, anchor=center, inner sep=0}, phantom, from=3-1, to=3-3, start anchor=center, end anchor=center]
	\arrow["a"', from=3-1, to=5-1]
	\arrow["b", from=3-3, to=5-3]
	\arrow[""{name=2, anchor=center, inner sep=0}, "b", from=3-4, to=3-6]
	\arrow["Tf"', from=3-4, to=5-4]
	\arrow["f"{description}, from=3-6, to=5-6]
	\arrow[""{name=3, anchor=center, inner sep=0}, Rightarrow, no head, from=3-6, to=5-8]
	\arrow[""{name=3p, anchor=center, inner sep=0}, phantom, from=3-6, to=5-8, start anchor=center, end anchor=center]
	\arrow[""{name=4, anchor=center, inner sep=0}, "u"', from=5-1, to=5-3]
	\arrow[""{name=5, anchor=center, inner sep=0}, "a"', from=5-4, to=5-6]
	\arrow[""{name=6, anchor=center, inner sep=0}, "u"', from=5-6, to=5-8]
	\arrow[""{name=6p, anchor=center, inner sep=0}, phantom, from=5-6, to=5-8, start anchor=center, end anchor=center]
	\arrow["{T\eta}"', shorten <=4pt, shorten >=4pt, Rightarrow, from=0p, to=1p]
	\arrow["{\overline{u}}", shorten <=9pt, shorten >=9pt, Rightarrow, from=1, to=4]
	\arrow["{\overline{f}}", shorten <=9pt, shorten >=9pt, Rightarrow, from=2, to=5]
	\arrow["\eta"', shorten <=4pt, shorten >=4pt, Rightarrow, from=3p, to=6p]
\end{tikzcd}\]
and the following 2-cells are equal\footnote{These tell us that $\eta, \epsilon$ really are squares in the double category of algebras for a pseudomonad, for the details see \cite{shulman2011comparing_HIST}. We will not pursue this direction in the thesis.}:
\[\begin{tikzcd}
	TA && TB \\
	\\
	A && B & TA && TB && B \\
	\\
	&& A &&& TA && A
	\arrow[""{name=0, anchor=center, inner sep=0}, "Tu"{description}, from=1-1, to=1-3]
	\arrow["a"', from=1-1, to=3-1]
	\arrow["b", from=1-3, to=3-3]
	\arrow[""{name=1, anchor=center, inner sep=0}, "u"{description}, from=3-1, to=3-3]
	\arrow[""{name=1p, anchor=center, inner sep=0}, phantom, from=3-1, to=3-3, start anchor=center, end anchor=center]
	\arrow[""{name=2, anchor=center, inner sep=0}, Rightarrow, no head, from=3-1, to=5-3]
	\arrow[""{name=2p, anchor=center, inner sep=0}, phantom, from=3-1, to=5-3, start anchor=center, end anchor=center]
	\arrow["f", from=3-3, to=5-3]
	\arrow[""{name=3, anchor=center, inner sep=0}, "Tu", from=3-4, to=3-6]
	\arrow[""{name=3p, anchor=center, inner sep=0}, phantom, from=3-4, to=3-6, start anchor=center, end anchor=center]
	\arrow[""{name=4, anchor=center, inner sep=0}, Rightarrow, no head, from=3-4, to=5-6]
	\arrow[""{name=4p, anchor=center, inner sep=0}, phantom, from=3-4, to=5-6, start anchor=center, end anchor=center]
	\arrow[""{name=5, anchor=center, inner sep=0}, "b", from=3-6, to=3-8]
	\arrow["Tf"', from=3-6, to=5-6]
	\arrow["f"{description}, from=3-8, to=5-8]
	\arrow[""{name=6, anchor=center, inner sep=0}, "a"', from=5-6, to=5-8]
	\arrow["{\overline{u}}", shorten <=9pt, shorten >=9pt, Rightarrow, from=0, to=1]
	\arrow["\epsilon", shorten <=4pt, shorten >=4pt, Rightarrow, from=1p, to=2p]
	\arrow["{T\epsilon}", shorten <=4pt, shorten >=4pt, Rightarrow, from=3p, to=4p]
	\arrow["{\overline{f}}", shorten <=9pt, shorten >=9pt, Rightarrow, from=5, to=6]
\end{tikzcd}\]
\end{pozn}

Informally put, these equations tell us that the ``only” obstruction preventing the doctrinal adjunction from Theorem \ref{THM_doctrinal_adj} from living in $\colaxtalgl$ is the fact that the 2-cell $\overline{f}$ goes in the other direction. If $\overline{f}$ is invertible, the adjunction readily lifts to $\colaxtalgl$ to an adjunction $(\epsilon, \eta): (f,\overline{f}^{-1}) \dashv (u,\overline{u})$. Another reformulation gives us two categorifications of the statement from ordinary category theory that ``monadic functors reflect isomorphisms”:

\begin{cor}\label{THMC_doct_adj_corollary}
Let $(T,m,i)$ be a pseudomonad on a 2-category $\ck$.
\begin{itemize}
\item Let $(f,\overline{f}): \mathbb{A} \to \mathbb{B}$ be a pseudo morphism between colax $T$-algebras such that there is an adjunction $(\epsilon, \eta): f \dashv u$ in $\ck$. Then the adjunction lifts to an adjunction in $\colaxtalgl$:
\begin{equation}\label{EQ_eps_eta_fovf_uovu}
(\epsilon, \eta): (f,\overline{f}) \dashv (u,\overline{u}),
\end{equation}
where $\overline{u}$ is the mate of $(\overline{f})^{-1}$.

\item Moreover, if the adjunction $(\epsilon, \eta): f \dashv u$ in $\ck$ is an adjoint equivalence, $\overline{u}$ is invertible and so \eqref{EQ_eps_eta_fovf_uovu} is an adjoint equivalence in $\colaxtalg$.
\end{itemize}
\end{cor}
\begin{proof}
Both statements follow immediately from Theorem \ref{THM_doctrinal_adj} and Remark \ref{POZN_mixed_algebra_2cells}.
\end{proof}

Every lax morphism that is a left adjoint is a pseudo morphism:

\begin{cor}\label{THMC_left_adjoint_lax_mor_je_pseudo}
Given any adjunction in $\colaxtalgl$:
\[\begin{tikzcd}
	{(\epsilon,\eta):} & {(A,a,\gamma,\iota)} && {(B,b,\gamma',\iota')}
	\arrow[""{name=0, anchor=center, inner sep=0}, "{(f,\overline{f})}"', curve={height=24pt}, from=1-2, to=1-4]
	\arrow[""{name=1, anchor=center, inner sep=0}, "{(u,\overline{u})}"', curve={height=24pt}, from=1-4, to=1-2]
	\arrow["\dashv"{anchor=center, rotate=90}, draw=none, from=0, to=1]
\end{tikzcd}\]
The morphism $(f,\overline{f})$ is a pseudo morphism, i.e. $\overline{f}$ is invertible. In particular, every equivalence in $\colaxtalgl$ is a pseudo morphism.
\end{cor}
\begin{proof}
Denote by $\widetilde{f}$ the mate of $\overline{u}$ that gives $f$ a colax morphism structure. We claim that $\widetilde{f}$ is inverse to the 2-cell $\overline{f}$. For instance, consider the composite $\overline{f} \circ \widetilde{f}$ as portrayed below left. It reduces to the composite below right since $\eta$ is an algebra 2-cell:
\[\begin{tikzcd}
	&& TA &&&&& TA \\
	\\
	A && TB && TA & A &&&& TA \\
	\\
	&& B && A &&& B && A \\
	&&&& B &&&&& B
	\arrow["a"', from=1-3, to=3-1]
	\arrow["Tf", from=1-3, to=3-3]
	\arrow[""{name=0, anchor=center, inner sep=0}, Rightarrow, no head, from=1-3, to=3-5]
	\arrow[""{name=0p, anchor=center, inner sep=0}, phantom, from=1-3, to=3-5, start anchor=center, end anchor=center]
	\arrow["a"', from=1-8, to=3-6]
	\arrow[Rightarrow, no head, from=1-8, to=3-10]
	\arrow[""{name=1, anchor=center, inner sep=0}, "f"', from=3-1, to=5-3]
	\arrow[""{name=2, anchor=center, inner sep=0}, "Tu"{description}, from=3-3, to=3-5]
	\arrow[""{name=2p, anchor=center, inner sep=0}, phantom, from=3-3, to=3-5, start anchor=center, end anchor=center]
	\arrow[""{name=2p, anchor=center, inner sep=0}, phantom, from=3-3, to=3-5, start anchor=center, end anchor=center]
	\arrow["b", from=3-3, to=5-3]
	\arrow["a", from=3-5, to=5-5]
	\arrow["f"', from=3-6, to=5-8]
	\arrow[""{name=3, anchor=center, inner sep=0}, Rightarrow, no head, from=3-6, to=5-10]
	\arrow[""{name=3p, anchor=center, inner sep=0}, phantom, from=3-6, to=5-10, start anchor=center, end anchor=center]
	\arrow["a", from=3-10, to=5-10]
	\arrow[""{name=4, anchor=center, inner sep=0}, "u", from=5-3, to=5-5]
	\arrow[""{name=4p, anchor=center, inner sep=0}, phantom, from=5-3, to=5-5, start anchor=center, end anchor=center]
	\arrow[""{name=4p, anchor=center, inner sep=0}, phantom, from=5-3, to=5-5, start anchor=center, end anchor=center]
	\arrow[""{name=5, anchor=center, inner sep=0}, Rightarrow, no head, from=5-3, to=6-5]
	\arrow[""{name=5p, anchor=center, inner sep=0}, phantom, from=5-3, to=6-5, start anchor=center, end anchor=center]
	\arrow["f", from=5-5, to=6-5]
	\arrow[""{name=6, anchor=center, inner sep=0}, "u", from=5-8, to=5-10]
	\arrow[""{name=6p, anchor=center, inner sep=0}, phantom, from=5-8, to=5-10, start anchor=center, end anchor=center]
	\arrow[""{name=7, anchor=center, inner sep=0}, Rightarrow, no head, from=5-8, to=6-10]
	\arrow[""{name=7p, anchor=center, inner sep=0}, phantom, from=5-8, to=6-10, start anchor=center, end anchor=center]
	\arrow["f", from=5-10, to=6-10]
	\arrow["{T\eta}"', shorten <=4pt, shorten >=4pt, Rightarrow, from=0p, to=2p]
	\arrow["{\overline{f}}"', shorten >=6pt, Rightarrow, from=3-3, to=1]
	\arrow["{\overline{u}}", shorten <=9pt, shorten >=9pt, Rightarrow, from=2p, to=4p]
	\arrow["\eta", shorten <=3pt, Rightarrow, from=3p, to=5-8]
	\arrow["\epsilon", shorten <=2pt, shorten >=2pt, Rightarrow, from=4p, to=5p]
	\arrow["\epsilon", shorten <=2pt, shorten >=2pt, Rightarrow, from=6p, to=7p]
\end{tikzcd}\]
By the first triangle identity for $f \dashv u$, the above right 2-cell is the identity. Similarly, the composite $\widetilde{f} \circ \overline{f}$ can be shown to be the identity using the fact that $\epsilon$ is an algebra 2-cell and the first triangle identity for $Tf \dashv Tu$.

\end{proof}

\begin{pr}
Applied to the free cocompletion pseudomonad that will be introduced later (in Example \ref{PR_presheafpsmonad}), Theorem \ref{THM_doctrinal_adj} implies the statement from ordinary category theory that says that ``left adjoints preserve colimits”.
\end{pr}

\subsection{Adjunctions versus algebras}\label{SUBS_psmonads_adjsvsalgs}

\begin{lemma}\label{THML_Tpsmonad_adjunkce_gives_colaxalg}
Let $(T,m,i)$ be a pseudomonad. Let $(\iota, \sigma): a \dashv i_A: A \to TA$ be an adjunction. Then there exists a unique 2-cell $\gamma: a m_A \Rightarrow a Ta$ such that $(A,a,\gamma,\iota)$ is a colax $T$-algebra.
\end{lemma}
\begin{proof}
Due to taking mates, the algebra unit axiom asserting the equality of the following 2-cells:
\[\begin{tikzcd}[column sep=scriptsize]
	&&& TA \\
	TA && {T^2A} && A & TA & {T^2A} & TA & A \\
	&&& TA
	\arrow["a", from=1-4, to=2-5]
	\arrow["\gamma", shorten <=6pt, shorten >=6pt, Rightarrow, from=1-4, to=3-4]
	\arrow["{Ti_A}", from=2-1, to=2-3]
	\arrow[""{name=0, anchor=center, inner sep=0}, curve={height=24pt}, Rightarrow, no head, from=2-1, to=3-4]
	\arrow["{m_A}", from=2-3, to=1-4]
	\arrow["Ta"{description}, from=2-3, to=3-4]
	\arrow["{Ti_A}", from=2-6, to=2-7]
	\arrow[""{name=1, anchor=center, inner sep=0}, curve={height=30pt}, Rightarrow, no head, from=2-6, to=2-8]
	\arrow[""{name=1p, anchor=center, inner sep=0}, phantom, from=2-6, to=2-8, start anchor=center, end anchor=center, curve={height=30pt}]
	\arrow["{m_A}", from=2-7, to=2-8]
	\arrow["a", from=2-8, to=2-9]
	\arrow["a"', from=3-4, to=2-5]
	\arrow["{T\iota}"', shorten >=5pt, Rightarrow, from=2-3, to=0]
	\arrow["{\eta_A^{-1}}", shorten >=3pt, Rightarrow, from=2-7, to=1p]
\end{tikzcd}\]

holds if and only if:
\[\begin{tikzcd}
	& {T^2A} \\
	\\
	{\gamma \sstackrel{!}{=}} & TA && {T^2A} && TA & A
	\arrow["Ta"', from=1-2, to=3-2]
	\arrow[""{name=0, anchor=center, inner sep=0}, Rightarrow, no head, from=1-2, to=3-4]
	\arrow[""{name=0p, anchor=center, inner sep=0}, phantom, from=1-2, to=3-4, start anchor=center, end anchor=center]
	\arrow[""{name=1, anchor=center, inner sep=0}, "{Ti_A}"{description}, from=3-2, to=3-4]
	\arrow[""{name=1p, anchor=center, inner sep=0}, phantom, from=3-2, to=3-4, start anchor=center, end anchor=center]
	\arrow[""{name=2, anchor=center, inner sep=0}, curve={height=30pt}, Rightarrow, no head, from=3-2, to=3-6]
	\arrow["{m_A}", from=3-4, to=3-6]
	\arrow["a", from=3-6, to=3-7]
	\arrow["{T\sigma}"', shorten <=4pt, shorten >=4pt, Rightarrow, from=0p, to=1p]
	\arrow["{\eta^{-1}_A}", shorten <=3pt, shorten >=3pt, Rightarrow, from=3-4, to=2]
\end{tikzcd}\]
This forces the definition of the 2-cell $\gamma$. Let us verify that the second algebra unit axiom holds. In other words, this composite is the identity on $a:TA \to A$:
\[\begin{tikzcd}
	TA \\
	& {T^2A} \\
	\\
	& TA && {T^2A} && TA & A \\
	A
	\arrow[""{name=0, anchor=center, inner sep=0}, "{i_{TA}}"{description}, from=1-1, to=2-2]
	\arrow[""{name=0p, anchor=center, inner sep=0}, phantom, from=1-1, to=2-2, start anchor=center, end anchor=center]
	\arrow[""{name=1, anchor=center, inner sep=0}, curve={height=-24pt}, Rightarrow, no head, from=1-1, to=4-6]
	\arrow["a"', from=1-1, to=5-1]
	\arrow["Ta"{description}, from=2-2, to=4-2]
	\arrow[""{name=2, anchor=center, inner sep=0}, Rightarrow, no head, from=2-2, to=4-4]
	\arrow[""{name=2p, anchor=center, inner sep=0}, phantom, from=2-2, to=4-4, start anchor=center, end anchor=center]
	\arrow[""{name=3, anchor=center, inner sep=0}, "{Ti_A}"{description}, from=4-2, to=4-4]
	\arrow[""{name=3p, anchor=center, inner sep=0}, phantom, from=4-2, to=4-4, start anchor=center, end anchor=center]
	\arrow[""{name=4, anchor=center, inner sep=0}, shift right=5, curve={height=24pt}, Rightarrow, no head, from=4-2, to=4-6]
	\arrow["{m_A}", from=4-4, to=4-6]
	\arrow["a", from=4-6, to=4-7]
	\arrow[""{name=5, anchor=center, inner sep=0}, "{i_A}"', from=5-1, to=4-2]
	\arrow[""{name=5p, anchor=center, inner sep=0}, phantom, from=5-1, to=4-2, start anchor=center, end anchor=center]
	\arrow[""{name=6, anchor=center, inner sep=0}, shift right=5, curve={height=30pt}, Rightarrow, no head, from=5-1, to=4-7]
	\arrow["{i^{-1}_a}"', shorten <=13pt, shorten >=13pt, Rightarrow, from=0p, to=5p]
	\arrow["{\beta_A^{-1}}", shorten <=8pt, Rightarrow, from=1, to=4-4]
	\arrow["{T\sigma}"', shorten <=4pt, shorten >=4pt, Rightarrow, from=2p, to=3p]
	\arrow["\iota", shorten <=3pt, shorten >=3pt, Rightarrow, from=4, to=6]
	\arrow["{\eta^{-1}_A}", shorten <=3pt, shorten >=3pt, Rightarrow, from=4-4, to=4]
\end{tikzcd}\]
Using the local naturality of $i: 1_\ck \Rightarrow T$ at the 2-cell $\sigma$, this reduces to:
\[\begin{tikzcd}
	&&& TA \\
	TA &&&& {T^2A} \\
	& A && TA & {} & TA \\
	&&&&&& A
	\arrow["{i_{TA}}", from=1-4, to=2-5]
	\arrow["{i_{i_A}}"', shorten <=6pt, shorten >=6pt, Rightarrow, from=1-4, to=3-4]
	\arrow[""{name=0, anchor=center, inner sep=0}, shift left=2, curve={height=-30pt}, Rightarrow, no head, from=1-4, to=3-6]
	\arrow[""{name=1, anchor=center, inner sep=0}, Rightarrow, no head, from=2-1, to=1-4]
	\arrow["a"', from=2-1, to=3-2]
	\arrow["{m_A}", from=2-5, to=3-6]
	\arrow["{i_A}"{description}, from=3-2, to=1-4]
	\arrow["{i_A}"', from=3-2, to=3-4]
	\arrow[""{name=2, anchor=center, inner sep=0}, curve={height=24pt}, Rightarrow, no head, from=3-2, to=4-7]
	\arrow["{Ti_A}"{description}, from=3-4, to=2-5]
	\arrow[""{name=3, anchor=center, inner sep=0}, draw=none, from=3-4, to=3-5]
	\arrow[""{name=4, anchor=center, inner sep=0}, Rightarrow, no head, from=3-4, to=3-6]
	\arrow["a", from=3-6, to=4-7]
	\arrow["{\beta^{-1}_A}"', shorten <=4pt, Rightarrow, from=0, to=2-5]
	\arrow["\sigma"', shorten <=6pt, shorten >=6pt, Rightarrow, from=1, to=3-2]
	\arrow["{\eta^{-1}_A}", shorten >=3pt, Rightarrow, from=2-5, to=4]
	\arrow["\iota", shorten <=5pt, shorten >=5pt, Rightarrow, from=3, to=2]
\end{tikzcd}\]
The 2-cell composite involving $\beta^{-1}_A$, $i_{i_A}$, $\eta^{-1}_A$ now equals the identity on $i_A$ because of Lemma \ref{THML_in_every_psmonad}. The final composite involving $\sigma, \iota$ is equal to the identity because of the first triangle identity for $a \dashv i_A$.

To verify the \textbf{associativity}, we need to show that this composite:
\[\begin{tikzcd}[column sep=scriptsize]
	{T^3A} &&&& {T^3A} && {T^2A} \\
	& {T^2A} && {T^2A} && {T^2A} && TA \\
	&& TA
	\arrow[""{name=0, anchor=center, inner sep=0}, Rightarrow, no head, from=1-1, to=1-5]
	\arrow["{T^2a}"', from=1-1, to=2-2]
	\arrow["{m_{TA}}", from=1-5, to=1-7]
	\arrow["{Tm_A}", from=1-5, to=2-6]
	\arrow["{m_A}", from=1-7, to=2-8]
	\arrow[""{name=1, anchor=center, inner sep=0}, Rightarrow, no head, from=2-2, to=2-4]
	\arrow["Ta"', from=2-2, to=3-3]
	\arrow["{T^2i_A}", from=2-4, to=1-5]
	\arrow[""{name=2, anchor=center, inner sep=0}, Rightarrow, no head, from=2-4, to=2-6]
	\arrow[""{name=3, anchor=center, inner sep=0}, "{m_A}"', from=2-6, to=2-8]
	\arrow["{Ti_A}"', from=3-3, to=2-4]
	\arrow["{T^2\sigma}", shorten <=6pt, shorten >=4pt, Rightarrow, from=0, to=1]
	\arrow["{T\eta^{-1}_A}", shorten >=3pt, Rightarrow, from=1-5, to=2]
	\arrow["{\mu^{-1}_A}"', shorten <=3pt, shorten >=3pt, Rightarrow, from=1-7, to=3]
	\arrow["{T\sigma}", shorten <=3pt, Rightarrow, from=1, to=3-3]
\end{tikzcd}\]
equals this composite:
\[\begin{tikzcd}
	{T^3A} &&& {T^2A} && {T^2A} & TA \\
	&&&& TA \\
	{T^2A} && {T^2A} & {} \\
	TA
	\arrow[""{name=0, anchor=center, inner sep=0}, "{m_{TA}}", from=1-1, to=1-4]
	\arrow[""{name=0p, anchor=center, inner sep=0}, phantom, from=1-1, to=1-4, start anchor=center, end anchor=center]
	\arrow["{T^2a}"', from=1-1, to=3-1]
	\arrow[""{name=1, anchor=center, inner sep=0}, Rightarrow, no head, from=1-4, to=1-6]
	\arrow["Ta"', from=1-4, to=2-5]
	\arrow["{m_A}", from=1-6, to=1-7]
	\arrow["{Ti_A}"', from=2-5, to=1-6]
	\arrow[""{name=2, anchor=center, inner sep=0}, Rightarrow, no head, from=3-1, to=3-3]
	\arrow[""{name=2p, anchor=center, inner sep=0}, phantom, from=3-1, to=3-3, start anchor=center, end anchor=center]
	\arrow[""{name=3, anchor=center, inner sep=0}, draw=none, from=3-1, to=3-4]
	\arrow[""{name=3p, anchor=center, inner sep=0}, phantom, from=3-1, to=3-4, start anchor=center, end anchor=center]
	\arrow["Ta"', from=3-1, to=4-1]
	\arrow["{m_A}", from=3-3, to=2-5]
	\arrow[""{name=4, anchor=center, inner sep=0}, curve={height=30pt}, Rightarrow, no head, from=4-1, to=2-5]
	\arrow[""{name=5, anchor=center, inner sep=0}, "{Ti_A}"', from=4-1, to=3-3]
	\arrow[""{name=5p, anchor=center, inner sep=0}, phantom, from=4-1, to=3-3, start anchor=center, end anchor=center]
	\arrow["{m_{a}^{-1}}"', shorten <=9pt, shorten >=9pt, Rightarrow, from=0p, to=3p]
	\arrow["{T\sigma}", shorten <=3pt, Rightarrow, from=1, to=2-5]
	\arrow["{T\sigma}"', shorten <=2pt, shorten >=2pt, Rightarrow, from=2p, to=5p]
	\arrow["{\eta^{-1}_A}", shorten >=3pt, Rightarrow, from=3-3, to=4]
\end{tikzcd}\]
Using the local naturality of $m: T^2 \Rightarrow T$, the first mentioned composite changes to:
\[\begin{tikzcd}
	& {T^2A} && TA && {T^2A} \\
	{T^3A} &&&& {T^3A} && TA \\
	& {T^2A} && {T^2A} && {T^2A} \\
	&& TA
	\arrow[""{name=0, anchor=center, inner sep=0}, "Ta", from=1-2, to=1-4]
	\arrow[""{name=1, anchor=center, inner sep=0}, shift left=2, curve={height=-30pt}, Rightarrow, no head, from=1-2, to=1-6]
	\arrow[""{name=1p, anchor=center, inner sep=0}, phantom, from=1-2, to=1-6, start anchor=center, end anchor=center, shift left=2, curve={height=-30pt}]
	\arrow[""{name=2, anchor=center, inner sep=0}, draw=none, from=1-2, to=1-6]
	\arrow[""{name=2p, anchor=center, inner sep=0}, phantom, from=1-2, to=1-6, start anchor=center, end anchor=center]
	\arrow[""{name=3, anchor=center, inner sep=0}, "{Ti_A}", from=1-4, to=1-6]
	\arrow["{m_A}", from=1-6, to=2-7]
	\arrow["{\mu^{-1}_A}", shorten <=6pt, shorten >=6pt, Rightarrow, from=1-6, to=3-6]
	\arrow["{m_{TA}}", from=2-1, to=1-2]
	\arrow["{T^2a}"', from=2-1, to=3-2]
	\arrow["{m_{TA}}"{description}, from=2-5, to=1-6]
	\arrow["{Tm_A}", from=2-5, to=3-6]
	\arrow[""{name=4, anchor=center, inner sep=0}, Rightarrow, no head, from=3-2, to=3-4]
	\arrow["Ta"', from=3-2, to=4-3]
	\arrow["{m_A}", from=3-4, to=1-4]
	\arrow["{T^2i_A}"{description}, from=3-4, to=2-5]
	\arrow[""{name=5, anchor=center, inner sep=0}, Rightarrow, no head, from=3-4, to=3-6]
	\arrow["{m_A}"', from=3-6, to=2-7]
	\arrow["{Ti_A}"', from=4-3, to=3-4]
	\arrow["{T\sigma}", shorten <=2pt, shorten >=5pt, Rightarrow, from=1p, to=2p]
	\arrow["{m_a^{-1}}"', shorten <=9pt, shorten >=9pt, Rightarrow, from=0, to=4]
	\arrow["{m_{i_A}^{-1}}"', shorten <=3pt, Rightarrow, from=3, to=2-5]
	\arrow["{T\eta^{-1}_A}", shorten >=3pt, Rightarrow, from=2-5, to=5]
	\arrow["{T\sigma}", shorten <=3pt, Rightarrow, from=4, to=4-3]
\end{tikzcd}\]
We see that these composites will be equal if we prove that the 2-cells below are equal:
\[\begin{tikzcd}
	TA && {T^2A} && TA && {T^2A} && TA \\
	\\
	&&&& TA && {T^2A} && TA \\
	TA && {T^2A} && {T^3A} && {T^2A}
	\arrow["{Ti_A}", from=1-1, to=1-3]
	\arrow[""{name=0, anchor=center, inner sep=0}, curve={height=30pt}, Rightarrow, no head, from=1-1, to=1-5]
	\arrow[""{name=0p, anchor=center, inner sep=0}, phantom, from=1-1, to=1-5, start anchor=center, end anchor=center, curve={height=30pt}]
	\arrow[""{name=1, anchor=center, inner sep=0}, draw=none, from=1-1, to=1-5]
	\arrow[""{name=1p, anchor=center, inner sep=0}, phantom, from=1-1, to=1-5, start anchor=center, end anchor=center]
	\arrow["{m_A}", from=1-3, to=1-5]
	\arrow["{Ti_A}", from=1-5, to=1-7]
	\arrow["{m_A}", from=1-7, to=1-9]
	\arrow["{Ti_A}", from=3-5, to=3-7]
	\arrow["{m_{i_A}^{-1}}", Rightarrow, from=3-5, to=4-5]
	\arrow["{m_A}", from=3-7, to=3-9]
	\arrow["{\mu_{A}^{-1}}", Rightarrow, from=3-7, to=4-7]
	\arrow["{Ti_A}", from=4-1, to=4-3]
	\arrow["{m_A}", from=4-3, to=3-5]
	\arrow["{T^2i_A}"', from=4-3, to=4-5]
	\arrow[""{name=2, anchor=center, inner sep=0}, curve={height=30pt}, Rightarrow, no head, from=4-3, to=4-7]
	\arrow["{m_{TA}}"{description}, from=4-5, to=3-7]
	\arrow["{Tm_A}"{description}, from=4-5, to=4-7]
	\arrow["{m_A}"', from=4-7, to=3-9]
	\arrow["{\eta^{-1}_A}", shorten <=4pt, shorten >=4pt, Rightarrow, from=1p, to=0p]
	\arrow["{T\eta^{-1}_A}"{pos=0.3}, shorten >=3pt, Rightarrow, from=4-5, to=2]
\end{tikzcd}\]
But this equality follows from Lemma \ref{THML_in_every_psmonad}.
\end{proof}

\begin{lemma}\label{THML_Tpsmonad_1cell_mezi_adjs_laxmor}
Let $(T,m,i)$ be a pseudomonad and let $(\iota, \sigma): a \dashv i_A: A \to TA$ and $(\iota', \sigma'): b \dashv i_B$ be two adjunctions.
\begin{itemize}
\item For any 1-cell $f: A \to B$ in $\ck$ there is a \textbf{unique} 2-cell $\overline{f}: b Tf \Rightarrow f a$ making $(f,\overline{f})$ a lax morphism between the induced colax $T$-algebras $(A,a,\gamma, \iota)$, $(B,b,\gamma', \iota')$ from the previous lemma,
\item any 2-cell $\alpha: f \Rightarrow g$ in $\ck$ is an algebra 2-cell $(f,\overline{f}) \Rightarrow (g,\overline{g})$ for the canonical lax morphism structures from the previous point.
\end{itemize}
\end{lemma}
\begin{proof}
To prove the \textbf{first statement}, note that due to taking mates, the lax morphism unit axiom:
\[\begin{tikzcd}[column sep=scriptsize]
	A & TA && TB && A && TA && TB \\
	&&&& {=} &&& B \\
	& A && B &&& A &&& B
	\arrow["{i_A}", from=1-1, to=1-2]
	\arrow[""{name=0, anchor=center, inner sep=0}, curve={height=18pt}, Rightarrow, no head, from=1-1, to=3-2]
	\arrow[""{name=1, anchor=center, inner sep=0}, "Tf", from=1-2, to=1-4]
	\arrow["a"', from=1-2, to=3-2]
	\arrow["b", from=1-4, to=3-4]
	\arrow["{i_A}", from=1-6, to=1-8]
	\arrow["f", from=1-6, to=2-8]
	\arrow[curve={height=18pt}, Rightarrow, no head, from=1-6, to=3-7]
	\arrow["Tf", from=1-8, to=1-10]
	\arrow["{i_f^{-1}}", Rightarrow, from=1-8, to=2-8]
	\arrow["b", from=1-10, to=3-10]
	\arrow["{i_B}", from=2-8, to=1-10]
	\arrow[""{name=2, anchor=center, inner sep=0}, Rightarrow, no head, from=2-8, to=3-10]
	\arrow[""{name=3, anchor=center, inner sep=0}, "f"', from=3-2, to=3-4]
	\arrow["f"', from=3-7, to=3-10]
	\arrow["\iota"', shorten <=3pt, shorten >=6pt, Rightarrow, from=1-2, to=0]
	\arrow["{\overline{f}}", shorten <=9pt, shorten >=9pt, Rightarrow, from=1, to=3]
	\arrow["{\iota'}"', shorten <=7pt, shorten >=7pt, Rightarrow, from=1-10, to=2]
\end{tikzcd}\]
holds if and only if:
\[\begin{tikzcd}
	TA && TB && TA && TA && TB \\
	&&& {\sstackrel{!}{=}} \\
	A && B && A && B && B
	\arrow[""{name=0, anchor=center, inner sep=0}, "Tf", from=1-1, to=1-3]
	\arrow[""{name=0p, anchor=center, inner sep=0}, phantom, from=1-1, to=1-3, start anchor=center, end anchor=center]
	\arrow["a"', from=1-1, to=3-1]
	\arrow["b", from=1-3, to=3-3]
	\arrow[""{name=1, anchor=center, inner sep=0}, Rightarrow, no head, from=1-5, to=1-7]
	\arrow[""{name=1p, anchor=center, inner sep=0}, phantom, from=1-5, to=1-7, start anchor=center, end anchor=center]
	\arrow["a"', from=1-5, to=3-5]
	\arrow["Tf", from=1-7, to=1-9]
	\arrow["{i_f^{-1}}", shorten <=6pt, shorten >=6pt, Rightarrow, from=1-7, to=3-7]
	\arrow["b", from=1-9, to=3-9]
	\arrow[""{name=2, anchor=center, inner sep=0}, "f"', from=3-1, to=3-3]
	\arrow[""{name=2p, anchor=center, inner sep=0}, phantom, from=3-1, to=3-3, start anchor=center, end anchor=center]
	\arrow[""{name=3, anchor=center, inner sep=0}, "{i_A}"', from=3-5, to=1-7]
	\arrow[""{name=3p, anchor=center, inner sep=0}, phantom, from=3-5, to=1-7, start anchor=center, end anchor=center]
	\arrow["f", from=3-5, to=3-7]
	\arrow[""{name=4, anchor=center, inner sep=0}, "{i_B}", from=3-7, to=1-9]
	\arrow[""{name=4p, anchor=center, inner sep=0}, phantom, from=3-7, to=1-9, start anchor=center, end anchor=center]
	\arrow[""{name=5, anchor=center, inner sep=0}, Rightarrow, no head, from=3-7, to=3-9]
	\arrow[""{name=5p, anchor=center, inner sep=0}, phantom, from=3-7, to=3-9, start anchor=center, end anchor=center]
	\arrow["{\overline{f}}", shorten <=9pt, shorten >=9pt, Rightarrow, from=0p, to=2p]
	\arrow["\sigma"', shorten <=4pt, shorten >=4pt, Rightarrow, from=1p, to=3p]
	\arrow["{\iota'}", shorten <=4pt, shorten >=4pt, Rightarrow, from=4p, to=5p]
\end{tikzcd}\]
So the definition of $\overline{f}$ is forced upon us. To verify the lax morphism associativity axiom \eqref{EQD_laxmor_assoc}, note that the right-hand side of it becomes the composite:
\[\begin{tikzcd}
	{T^2A} && {T^2A} & {} && {T^2B} \\
	\\
	TA & TA && TB && TB && {T^2B} \\
	A && B &&&&& TB \\
	&&&&&&& B
	\arrow[""{name=0, anchor=center, inner sep=0}, Rightarrow, no head, from=1-1, to=1-3]
	\arrow[""{name=0p, anchor=center, inner sep=0}, phantom, from=1-1, to=1-3, start anchor=center, end anchor=center]
	\arrow["Ta"', from=1-1, to=3-1]
	\arrow["{T^2f}", from=1-3, to=1-6]
	\arrow["{Ti_f^{-1}}"', Rightarrow, from=1-4, to=3-4]
	\arrow["Tb"', from=1-6, to=3-6]
	\arrow[""{name=1, anchor=center, inner sep=0}, Rightarrow, no head, from=1-6, to=3-8]
	\arrow[""{name=1p, anchor=center, inner sep=0}, phantom, from=1-6, to=3-8, start anchor=center, end anchor=center]
	\arrow[""{name=2, anchor=center, inner sep=0}, "{Ti_A}"', from=3-1, to=1-3]
	\arrow[""{name=2p, anchor=center, inner sep=0}, phantom, from=3-1, to=1-3, start anchor=center, end anchor=center]
	\arrow[""{name=3, anchor=center, inner sep=0}, Rightarrow, no head, from=3-1, to=3-2]
	\arrow[""{name=3p, anchor=center, inner sep=0}, phantom, from=3-1, to=3-2, start anchor=center, end anchor=center]
	\arrow["a"', from=3-1, to=4-1]
	\arrow[""{name=4, anchor=center, inner sep=0}, "Tf", from=3-2, to=3-4]
	\arrow[""{name=5, anchor=center, inner sep=0}, "{Ti_B}", from=3-4, to=1-6]
	\arrow[""{name=5p, anchor=center, inner sep=0}, phantom, from=3-4, to=1-6, start anchor=center, end anchor=center]
	\arrow[""{name=6, anchor=center, inner sep=0}, Rightarrow, no head, from=3-4, to=3-6]
	\arrow[""{name=6p, anchor=center, inner sep=0}, phantom, from=3-4, to=3-6, start anchor=center, end anchor=center]
	\arrow[""{name=7, anchor=center, inner sep=0}, "{Ti_B}"{description}, from=3-6, to=3-8]
	\arrow[""{name=7p, anchor=center, inner sep=0}, phantom, from=3-6, to=3-8, start anchor=center, end anchor=center]
	\arrow[""{name=7p, anchor=center, inner sep=0}, phantom, from=3-6, to=3-8, start anchor=center, end anchor=center]
	\arrow[""{name=8, anchor=center, inner sep=0}, Rightarrow, no head, from=3-6, to=4-8]
	\arrow[""{name=8p, anchor=center, inner sep=0}, phantom, from=3-6, to=4-8, start anchor=center, end anchor=center]
	\arrow["{m_B}", from=3-8, to=4-8]
	\arrow[""{name=9, anchor=center, inner sep=0}, "{i_A}"', from=4-1, to=3-2]
	\arrow[""{name=9p, anchor=center, inner sep=0}, phantom, from=4-1, to=3-2, start anchor=center, end anchor=center]
	\arrow["f"', from=4-1, to=4-3]
	\arrow[""{name=10, anchor=center, inner sep=0}, "{i_B}", from=4-3, to=4-8]
	\arrow[""{name=10p, anchor=center, inner sep=0}, phantom, from=4-3, to=4-8, start anchor=center, end anchor=center]
	\arrow[""{name=11, anchor=center, inner sep=0}, Rightarrow, no head, from=4-3, to=5-8]
	\arrow[""{name=11p, anchor=center, inner sep=0}, phantom, from=4-3, to=5-8, start anchor=center, end anchor=center]
	\arrow["b", from=4-8, to=5-8]
	\arrow["{T\sigma}"', shorten <=4pt, shorten >=4pt, Rightarrow, from=0p, to=2p]
	\arrow["{T\sigma'}"', shift right=5, shorten <=4pt, shorten >=4pt, Rightarrow, from=1p, to=7p]
	\arrow["\sigma"', shorten <=2pt, shorten >=2pt, Rightarrow, from=3p, to=9p]
	\arrow["{i^{-1}_f}", shorten <=3pt, Rightarrow, from=4, to=4-3]
	\arrow["{T\iota'}", shorten <=4pt, shorten >=4pt, Rightarrow, from=5p, to=6p]
	\arrow["{\eta_B^{-1}}", shift left=5, shorten <=2pt, Rightarrow, from=7p, to=8p]
	\arrow["{\iota'}"{pos=0.3}, shorten <=2pt, shorten >=2pt, Rightarrow, from=10p, to=11p]
\end{tikzcd}\]
After using the triangle equality for $Tb \dashv Ti_B$, this becomes:
\[\begin{tikzcd}
	{T^2A} &&& {T^2A} \\
	TA & TA && {} && TB && {T^2B} \\
	\\
	A & {} & B &&& {} && TB \\
	&&&&&&& B
	\arrow[""{name=0, anchor=center, inner sep=0}, Rightarrow, no head, from=1-1, to=1-4]
	\arrow[""{name=0p, anchor=center, inner sep=0}, phantom, from=1-1, to=1-4, start anchor=center, end anchor=center]
	\arrow["Ta"', from=1-1, to=2-1]
	\arrow["{Ti_f^{-1}}", Rightarrow, from=1-4, to=2-4]
	\arrow["{T^2f}", from=1-4, to=2-8]
	\arrow[""{name=1, anchor=center, inner sep=0}, "{Ti_A}"', from=2-1, to=1-4]
	\arrow[""{name=1p, anchor=center, inner sep=0}, phantom, from=2-1, to=1-4, start anchor=center, end anchor=center]
	\arrow[""{name=2, anchor=center, inner sep=0}, Rightarrow, no head, from=2-1, to=2-2]
	\arrow[""{name=2p, anchor=center, inner sep=0}, phantom, from=2-1, to=2-2, start anchor=center, end anchor=center]
	\arrow["a"', from=2-1, to=4-1]
	\arrow[""{name=3, anchor=center, inner sep=0}, "Tf"{description}, from=2-2, to=2-6]
	\arrow[""{name=4, anchor=center, inner sep=0}, "{Ti_B}"{description}, from=2-6, to=2-8]
	\arrow[""{name=4p, anchor=center, inner sep=0}, phantom, from=2-6, to=2-8, start anchor=center, end anchor=center]
	\arrow[""{name=5, anchor=center, inner sep=0}, Rightarrow, no head, from=2-6, to=4-8]
	\arrow[""{name=5p, anchor=center, inner sep=0}, phantom, from=2-6, to=4-8, start anchor=center, end anchor=center]
	\arrow["{m_B}", from=2-8, to=4-8]
	\arrow[""{name=6, anchor=center, inner sep=0}, "{i_A}"', from=4-1, to=2-2]
	\arrow[""{name=6p, anchor=center, inner sep=0}, phantom, from=4-1, to=2-2, start anchor=center, end anchor=center]
	\arrow["f"', from=4-1, to=4-3]
	\arrow[""{name=7, anchor=center, inner sep=0}, draw=none, from=4-2, to=4-6]
	\arrow[""{name=8, anchor=center, inner sep=0}, "{i_B}", from=4-3, to=4-8]
	\arrow[""{name=8p, anchor=center, inner sep=0}, phantom, from=4-3, to=4-8, start anchor=center, end anchor=center]
	\arrow[""{name=9, anchor=center, inner sep=0}, Rightarrow, no head, from=4-3, to=5-8]
	\arrow[""{name=9p, anchor=center, inner sep=0}, phantom, from=4-3, to=5-8, start anchor=center, end anchor=center]
	\arrow["b", from=4-8, to=5-8]
	\arrow["{T\sigma}"', shorten <=2pt, shorten >=2pt, Rightarrow, from=0p, to=1p]
	\arrow["\sigma"', shorten <=4pt, shorten >=4pt, Rightarrow, from=2p, to=6p]
	\arrow["{i^{-1}_f}", shorten <=9pt, shorten >=9pt, Rightarrow, from=3, to=7]
	\arrow["{\eta_B^{-1}}", shorten <=4pt, shorten >=4pt, Rightarrow, from=4p, to=5p]
	\arrow["{\iota'}"{pos=0.3}, shorten <=2pt, shorten >=2pt, Rightarrow, from=8p, to=9p]
\end{tikzcd}\]
The left-hand side of the associativity axiom \eqref{EQD_laxmor_assoc} for a lax morphism is the composite:
\[\begin{tikzcd}
	& {T^2A} &&&&& {T^2B} \\
	TA \\
	& {T^2A} \\
	& TA & TA & {} &&& TB \\
	& A && B &&& B
	\arrow[""{name=0, anchor=center, inner sep=0}, "{T^2f}", from=1-2, to=1-7]
	\arrow["Ta"', from=1-2, to=2-1]
	\arrow[""{name=1, anchor=center, inner sep=0}, Rightarrow, no head, from=1-2, to=3-2]
	\arrow["{m_B}", from=1-7, to=4-7]
	\arrow[""{name=2, anchor=center, inner sep=0}, "{Ti_A}"', from=2-1, to=3-2]
	\arrow[""{name=3, anchor=center, inner sep=0}, curve={height=18pt}, Rightarrow, no head, from=2-1, to=4-2]
	\arrow["{m_A}", from=3-2, to=4-2]
	\arrow[""{name=4, anchor=center, inner sep=0}, Rightarrow, no head, from=4-2, to=4-3]
	\arrow[""{name=4p, anchor=center, inner sep=0}, phantom, from=4-2, to=4-3, start anchor=center, end anchor=center]
	\arrow[""{name=5, anchor=center, inner sep=0}, draw=none, from=4-2, to=4-7]
	\arrow["a"', from=4-2, to=5-2]
	\arrow["Tf", from=4-3, to=4-7]
	\arrow["{i^{-1}_f}"', Rightarrow, from=4-4, to=5-4]
	\arrow["b", from=4-7, to=5-7]
	\arrow[""{name=6, anchor=center, inner sep=0}, "{i_A}"', from=5-2, to=4-3]
	\arrow[""{name=6p, anchor=center, inner sep=0}, phantom, from=5-2, to=4-3, start anchor=center, end anchor=center]
	\arrow["f"', from=5-2, to=5-4]
	\arrow[""{name=7, anchor=center, inner sep=0}, "{i_B}", from=5-4, to=4-7]
	\arrow[""{name=7p, anchor=center, inner sep=0}, phantom, from=5-4, to=4-7, start anchor=center, end anchor=center]
	\arrow[""{name=8, anchor=center, inner sep=0}, Rightarrow, no head, from=5-4, to=5-7]
	\arrow[""{name=8p, anchor=center, inner sep=0}, phantom, from=5-4, to=5-7, start anchor=center, end anchor=center]
	\arrow["{T\sigma}"', shorten <=4pt, shorten >=4pt, Rightarrow, from=1, to=2]
	\arrow["{m_f}", shorten <=13pt, shorten >=13pt, Rightarrow, from=0, to=5]
	\arrow["{\eta^{-1}_A}", shorten >=5pt, Rightarrow, from=3-2, to=3]
	\arrow["\sigma"', shorten <=2pt, shorten >=2pt, Rightarrow, from=4p, to=6p]
	\arrow["{\iota'}"{pos=0.3}, shorten <=2pt, shorten >=2pt, Rightarrow, from=7p, to=8p]
\end{tikzcd}\]
We see that it remains to verify that the following 2-cells are equal:
\[\begin{tikzcd}[column sep=scriptsize]
	&&&&& {T^2A} && {T^2B} \\
	{T^2A} & {T^2A} && {T^2B} & TA \\
	&&&&& {T^2A} \\
	TA && TB & TB && TA && TB
	\arrow[""{name=0, anchor=center, inner sep=0}, "{T^2f}", from=1-6, to=1-8]
	\arrow[""{name=0p, anchor=center, inner sep=0}, phantom, from=1-6, to=1-8, start anchor=center, end anchor=center]
	\arrow["Ta"', from=1-6, to=2-5]
	\arrow[""{name=1, anchor=center, inner sep=0}, Rightarrow, no head, from=1-6, to=3-6]
	\arrow["{m_B}", from=1-8, to=4-8]
	\arrow[""{name=2, anchor=center, inner sep=0}, Rightarrow, no head, from=2-1, to=2-2]
	\arrow[""{name=2p, anchor=center, inner sep=0}, phantom, from=2-1, to=2-2, start anchor=center, end anchor=center]
	\arrow["{m_A}"', from=2-1, to=4-1]
	\arrow["{T^2f}", from=2-2, to=2-4]
	\arrow["{m_B}", from=2-4, to=4-4]
	\arrow[""{name=3, anchor=center, inner sep=0}, "{Ti_A}"', from=2-5, to=3-6]
	\arrow[""{name=4, anchor=center, inner sep=0}, curve={height=18pt}, Rightarrow, no head, from=2-5, to=4-6]
	\arrow["{m_A}", from=3-6, to=4-6]
	\arrow[""{name=5, anchor=center, inner sep=0}, "{Ti_A}"', from=4-1, to=2-2]
	\arrow[""{name=5p, anchor=center, inner sep=0}, phantom, from=4-1, to=2-2, start anchor=center, end anchor=center]
	\arrow[""{name=6, anchor=center, inner sep=0}, "Tf"', from=4-1, to=4-3]
	\arrow[""{name=6p, anchor=center, inner sep=0}, phantom, from=4-1, to=4-3, start anchor=center, end anchor=center]
	\arrow[""{name=7, anchor=center, inner sep=0}, "{Ti_B}", from=4-3, to=2-4]
	\arrow[""{name=7p, anchor=center, inner sep=0}, phantom, from=4-3, to=2-4, start anchor=center, end anchor=center]
	\arrow[""{name=8, anchor=center, inner sep=0}, Rightarrow, no head, from=4-3, to=4-4]
	\arrow[""{name=8p, anchor=center, inner sep=0}, phantom, from=4-3, to=4-4, start anchor=center, end anchor=center]
	\arrow[""{name=9, anchor=center, inner sep=0}, "Tf"', from=4-6, to=4-8]
	\arrow[""{name=9p, anchor=center, inner sep=0}, phantom, from=4-6, to=4-8, start anchor=center, end anchor=center]
	\arrow["{T\sigma}"', shorten <=4pt, shorten >=4pt, Rightarrow, from=1, to=3]
	\arrow["{m_f}", shorten <=13pt, shorten >=13pt, Rightarrow, from=0p, to=9p]
	\arrow["{T\sigma}"', shorten <=4pt, shorten >=4pt, Rightarrow, from=2p, to=5p]
	\arrow["{Ti_f^{-1}}", shift left=5, shorten <=4pt, shorten >=7pt, Rightarrow, from=2-2, to=6p]
	\arrow["{\eta^{-1}_A}", shorten >=5pt, Rightarrow, from=3-6, to=4]
	\arrow["{\eta^{-1}_B}", shorten <=4pt, shorten >=4pt, Rightarrow, from=7p, to=8p]
\end{tikzcd}\]
This follows from the modification axiom of $\eta^{-1}: m \circ Ti \to 1_{1_T}$ for the 1-cell $f: A \to B$. Let us now prove the \textbf{second statement}. Let $\alpha: f \Rightarrow g : A \to B$ be a 2-cell. The algebra 2-cell axiom amounts to verifying the equality of the 2-cells:
\[\begin{tikzcd}
	TA && TA && TB & TA && TA && TB \\
	\\
	A && B && B & A && B && B
	\arrow[""{name=0, anchor=center, inner sep=0}, Rightarrow, no head, from=1-1, to=1-3]
	\arrow["a"', from=1-1, to=3-1]
	\arrow["Tf", from=1-3, to=1-5]
	\arrow["{i_f^{-1}}", shorten <=6pt, shorten >=6pt, Rightarrow, from=1-3, to=3-3]
	\arrow["b", from=1-5, to=3-5]
	\arrow[""{name=1, anchor=center, inner sep=0}, Rightarrow, no head, from=1-6, to=1-8]
	\arrow["a"', from=1-6, to=3-6]
	\arrow[""{name=2, anchor=center, inner sep=0}, "Tg"', from=1-8, to=1-10]
	\arrow[""{name=3, anchor=center, inner sep=0}, "Tf", curve={height=-30pt}, from=1-8, to=1-10]
	\arrow["{i_g^{-1}}", shorten <=6pt, shorten >=6pt, Rightarrow, from=1-8, to=3-8]
	\arrow["b", from=1-10, to=3-10]
	\arrow[""{name=4, anchor=center, inner sep=0}, "{i_A}"', from=3-1, to=1-3]
	\arrow[""{name=5, anchor=center, inner sep=0}, "f", from=3-1, to=3-3]
	\arrow[""{name=6, anchor=center, inner sep=0}, "g"', curve={height=30pt}, from=3-1, to=3-3]
	\arrow[""{name=7, anchor=center, inner sep=0}, "{i_B}", from=3-3, to=1-5]
	\arrow[""{name=8, anchor=center, inner sep=0}, Rightarrow, no head, from=3-3, to=3-5]
	\arrow[""{name=9, anchor=center, inner sep=0}, "{i_A}"', from=3-6, to=1-8]
	\arrow["g"', from=3-6, to=3-8]
	\arrow[""{name=10, anchor=center, inner sep=0}, "{i_B}", from=3-8, to=1-10]
	\arrow[""{name=11, anchor=center, inner sep=0}, Rightarrow, no head, from=3-8, to=3-10]
	\arrow["\sigma"', shorten <=4pt, shorten >=4pt, Rightarrow, from=0, to=4]
	\arrow["\sigma"', shorten <=4pt, shorten >=4pt, Rightarrow, from=1, to=9]
	\arrow["{T\alpha}", shorten <=4pt, shorten >=4pt, Rightarrow, from=3, to=2]
	\arrow["\alpha", shorten <=4pt, shorten >=4pt, Rightarrow, from=5, to=6]
	\arrow["{\iota'}", shorten <=4pt, shorten >=4pt, Rightarrow, from=7, to=8]
	\arrow["\iota", shorten <=4pt, shorten >=4pt, Rightarrow, from=10, to=11]
\end{tikzcd}\]
This equality follows from the local naturality of $i: 1_\ck \Rightarrow T$ at $\alpha$.
\end{proof}

\begin{defi}
Let $(T,m,i)$ be a pseudomonad. Denote by $\text{Adj}(i)$ the 2-category whose objects are pairs $(A, (\iota, \sigma): a \dashv i_A))$ of an object $A \in \ck$ and an adjoint to the pseudomonad unit evaluated at $A$. The 1-cells and 2-cells of $\text{Adj}(i)$ are those of $\ck$:
\[
\text{Adj}(i)((A, (\iota, \sigma): a \dashv i_A),(B, (\iota', \sigma'): b \dashv i_B)) := \ck(A,B).
\]
\end{defi}

From the previous two lemmas we obtain:
\begin{theo}\label{THM_2fun_phi_adj_colaxtalgl}
For any pseudomonad $(T,m,i)$ on a 2-category $\ck$ there is a 2-fully faithful 2-functor\footnote{A 2-functor $F$ is said to be \textit{2-fully faithful} if each hom-functor $F_{A,B}$ is an isomorphism of categories.} that is injective on objects:
\[
\Phi: \text{Adj}(i) \to \colaxtalgl.
\]
\end{theo}
\begin{proof}
Define $\Phi$ so that it sends $(A,(\iota, \sigma))$ to the colax algebra $(A,a,\gamma,\iota)$, where $\gamma$ is the unique coassociator from Lemma \ref{THML_Tpsmonad_adjunkce_gives_colaxalg}. Likewise, send the 1-cell $f$ to the lax morphism $(f,\overline{f})$ from Lemma \ref{THML_Tpsmonad_1cell_mezi_adjs_laxmor} and send the 2-cell $\alpha$ to $\alpha$ (by the last cited lemma, it is an algebra 2-cell). Again by this second lemma, we see that $\Phi$ is 2-fully faithful, and injectivity on objects is obvious.
\end{proof}

\begin{pozn}\label{POZN_Phi_sends_LAs_to_psmors}
Note that by Proposition \ref{THMC_left_adjoint_lax_mor_je_pseudo}, the 2-functor $\Phi$ sends left adjoints in $\ck$ to pseudo morphisms.
\end{pozn}

\subsection{Examples}\label{SUBS_psmonads_examples}

\begin{pr}[The identity 2-monad]\label{PR_identity2monad}

Consider the identity 2-monad $T := 1_{\ck}$ on a 2-category $\ck$. Strict algebras are trivial -- they are just the objects of $\ck$. Colax algebras are comonads in $\ck$ (Definition \ref{DEFI_comonad}), and lax morphisms are comonad functors (Definition \ref{DEFI_comonadfun}). In fact:
\[
\colaxtalgl = \cmndl{\ck}.
\]
\end{pr}

Since for every 2-monad $(T,m,i)$ on a 2-category $\ck$ there is a 2-monad morphism $i: 1_{\ck} \to T$ from the identity 2-monad to $T$, one may expect some kind of a relationship between colax $T$-algebras and comonads:

\begin{pozn}\label{POZN_laxalg_underlying_monad}
Fix a 2-monad $(T,m,i)$ on a 2-category $\ck$. Any colax algebra $(A,a,\gamma,\iota)$ has an underlying comonad in $\ck$:
\[
\mathbb{U}(A,a,\gamma,\iota) := (A, ai_A, \gamma i_{TA} i_A, \iota).
\]
This extends to a 2-functor $\colaxtalgl \to \cmndl{\ck}$ that sends a morphism:
\[
(F,\overline{F}): (A,a,\gamma, \iota) \to (B,b,\gamma', \iota'),
\]
to the comonad functor $(F,\overline{F}i_A)$, and sends an algebra 2-cell $\alpha$ to $\alpha$. There is also a 2-functor $\mathbb{F}: \cmndl{\ck} \to \colaxtalgl$ going the other way, defined on objects as:
\[
\mathbb{F}(A,t,\delta, \epsilon) := (TA, Tt m_A, T\delta m_A Tm_A, T\epsilon),
\]
that sends a comonad functor $(F,\overline{F}): (A,t,\delta, \epsilon) \to (B,s,\delta', \epsilon')$ to the lax morphism $(TF, T\overline{F}m_A)$. There also seems to be a unit $\widetilde{\eta}: 1_{\cmndl{\ck}} \Rightarrow \mathbb{UF}$ and a counit $\widetilde{\epsilon}: \mathbb{FU} \Rightarrow 1$, both colax natural, and also two modifications $\Psi: 1_{\mathbb{F}} \to \widetilde{\epsilon} \mathbb{F} \circ \mathbb{F} \widetilde{\eta}$, $\Phi: 1_{\mathbb{U}} \to \mathbb{U}\widetilde{\epsilon} \circ \widetilde{\eta}\mathbb{U}$. This does not however assemble into a (co)lax adjunction since one of the modifications always goes in the wrong way.
\end{pozn}

\begin{pr}(The terminal 2-monad)\label{PR_terminal2monad}
Fix a 2-category $\ck$ with a terminal object $*$ and consider the endo-2-functor $T$ on $\ck$ that is constant on $*$. This 2-functor has a 2-monad structure -- the multiplication is the identity, while the unit is given by the unique map $A \to *$. An object $A \in \ck$ is a $T$-algebra if and only if it is the terminal object in $\ck$.

It can be seen that a colax $T$-algebra is an ``object with an initial object”, i.e. a triple $(A,a: * \to A, \iota: a !_A \Rightarrow 1_A)$ such that there is an adjunction in $\ck$:
\[\begin{tikzcd}
	{(\iota,1):} & A && {*}
	\arrow[""{name=0, anchor=center, inner sep=0}, "{!_A}"', curve={height=24pt}, from=1-2, to=1-4]
	\arrow[""{name=1, anchor=center, inner sep=0}, "a"', curve={height=24pt}, from=1-4, to=1-2]
	\arrow["\dashv"{anchor=center, rotate=-90}, draw=none, from=1, to=0]
\end{tikzcd}\]
Next, notice that any 1-cell $f:A \to B$ between colax-$T$-algebras $(A,a,\iota)$, $(B,b,\iota')$ admits a unique lax morphism structure $(f,\overline{f})$, since the component $\overline{f}: b \Rightarrow fa$ is necessarily forced to equal the following:
\[\begin{tikzcd}
	& {*} \\
	&&& {*} &&&& {*} && {*} \\
	{\overline{f}=} & A &&&&& {=} & A \\
	&&& A && B && B && B
	\arrow[Rightarrow, no head, from=1-2, to=2-4]
	\arrow["a"', from=1-2, to=3-2]
	\arrow["a"{description}, from=2-4, to=4-4]
	\arrow[""{name=0, anchor=center, inner sep=0}, "b", from=2-4, to=4-6]
	\arrow[Rightarrow, no head, from=2-8, to=2-10]
	\arrow["a"', from=2-8, to=3-8]
	\arrow["b", from=2-10, to=4-10]
	\arrow["{!_A}"{description}, from=3-2, to=2-4]
	\arrow[""{name=1, anchor=center, inner sep=0}, Rightarrow, no head, from=3-2, to=4-4]
	\arrow["f"', from=3-8, to=4-8]
	\arrow[""{name=2, anchor=center, inner sep=0}, "f"', from=4-4, to=4-6]
	\arrow["{!_B}", from=4-8, to=2-10]
	\arrow[""{name=3, anchor=center, inner sep=0}, Rightarrow, no head, from=4-8, to=4-10]
	\arrow["\iota", shorten <=7pt, shorten >=7pt, Rightarrow, from=2-4, to=1]
	\arrow["{\overline{f}}"', shorten <=4pt, shorten >=4pt, Rightarrow, from=0, to=2]
	\arrow["{\iota'}"', shorten <=13pt, shorten >=9pt, Rightarrow, from=2-10, to=3]
\end{tikzcd}\]
Here the first equality follows from the first triangle identity for $a \dashv !_A$ and the second equality uses the lax morphism unit axiom.
\end{pr}

\begin{pr}[Free pointed category 2-monad]\label{PR_free_pointed_cat_2monad}
The 2-monad $(-)+1$ on the 2-category $\cat$ that adds an object to a category. The unit $i_\ca: \ca \to \ca + 1$ is given by the inclusion, while the multiplication $m_\ca : \ca + 1 + 1 \to \ca + 1$ concatenates the two added objects. One can verify that:
\begin{itemize}
\item a strict $T$-algebra is a pair $(\ca, a_0)$ of a small category and an object $a_0 \in \ob \ca$,
\item a strict morphism $(\ca, a_0) \to (\cb, b_0)$ is a functor $F: \ca \to \cb$ such that $Fa_0 = b_0$,
\item a lax morphism $(\ca, a_0) \to (\cb, b_0)$ is a pair $(F,f)$ of a functor $F: \ca \to \cb$ and a morphism $f: b_0 \to Fa_0$ in $\cb$.
\item a colax $T$-algebra consists of a comonad $(\ca, t, \delta, \epsilon)$ with a choice of a $t$-coalgebra $(a_0, \gamma_\bullet: a_0 \to ta_0)$.
\end{itemize}
\end{pr}

\begin{pr}[Free initial object 2-monad]\label{PR_freeinitialobject}
This 2-monad $T$ on $\cat$ adds a new object to a category $\ca$ together with a unique arrow from this object to every object of $\ca$. We have:
\begin{itemize}
\item a category $\ca$ admits the structure of a pseudo $T$-algebra if and only if it has an initial object,
\item a strict morphism is a functor preserving the chosen initial objects,
\item a pseudo/colax morphism is a functor that preserves initial objects,
\item any functor between pseudo $T$-algebras admits a unique lax morphism structure.
\end{itemize}

We will also provide the description of colax algebras for this 2-monad. We will use the following notation: in $T^2A$ we will denote by $*_2$ the second added initial object and by $*_1$ the first added initial object. We will use the following observation: if $\cc$ is a category and $\lambda: \Delta i \Rightarrow 1_{\cc}$ is a cone with the property that $\lambda_i = 1_i$, then $i$ is the initial object of $\cc$.

\begin{prop}\label{THMP_colax_algebras_for_freeinicobj2monad}
A colax $T$-algebra structure on a category $\ca$ is equivalently a comonad $(\ca, t, \delta, \epsilon)$ and an initial object $a_0$ in $\ca$.
\end{prop}
\begin{proof}
``$\Leftarrow$”: Since there is a 2-monad morphism $(-)+1 \Rightarrow T$ from the free pointed category 2-monad to $T$, every colax $T$-algebra $(\ca,A,\gamma,\iota)$ is in particular a colax $((-)+1)$-algebra, i.e. there is an underlying comonad $(\ca, A_0, \gamma i_{T\ca} i_{\ca}, \iota)$ (with $A_0 := A i_{\ca}$) and an $A_0$-coalgebra $(a_0, \gamma_{*_1}: a_0 \to Aa_0)$ (with $a_0 := A*$). 

\textbf{Our goal} will be to show that $a_0$ is the initial object of $\ca$. First, denote:

\[
\mu_b := A(!_b : * \to b) : a_0 \to Ab.
\]
Notice that this data assembles into a cone $\mu: \Delta a_0 \Rightarrow A_0: \ca \to \ca$. Also, because of the naturality of $\gamma: A \circ m_\ca \Rightarrow A \circ TA$ at $!_b: *_1 \to b$, the following commutes:
\[\begin{tikzcd}
	{Aa_0} && {A^2b} \\
	\\
	{a_0} && Ab
	\arrow["{A\mu_b}", from=1-1, to=1-3]
	\arrow["{\gamma_{*_1}}", from=3-1, to=1-1]
	\arrow["{\mu_b}"', from=3-1, to=3-3]
	\arrow["{\gamma_b}"', from=3-3, to=1-3]
\end{tikzcd}\]
In other words, for every $b \in \ca$, $\mu_b$ is an $A_0$-coalgebra morphism:
\[
\mu_b: (a_0, \gamma_{*_1}) \to (Ab, \gamma_b).
\]
If we denote the free $A_0$-coalgebra functor by $F_{A_0}: \ca \to \ca^{A_0}$, the collection of $\mu_b$'s is a cone:
\[
\mu: \Delta (a_0, \gamma_{*_1}) \Rightarrow F_{A_0}: \ca \to \ca^{A_0}.
\]
Transpose this cone along the free-forgetful adjunction $U_{A_0} \dashv F_{A_0}$ (take a mate) to obtain a cone:
\[
\widehat{\mu} : \Delta a_0 \Rightarrow 1_\ca,
\]
whose component at $b \in \ob \ca$ is the composite:
\[\begin{tikzcd}
	{a_0} && Ab && b
	\arrow["{\mu_b}", from=1-1, to=1-3]
	\arrow["{\iota_b}", from=1-3, to=1-5]
\end{tikzcd}\]

\noindent Consider again the naturality of $\gamma$, this time at the morphism $!: *_2 \to *_1$:
\[\begin{tikzcd}
	{a_0} && {a_0} \\
	\\
	{a_0} && {Aa_0}
	\arrow["{\gamma_{*_2}}", from=1-1, to=1-3]
	\arrow[Rightarrow, no head, from=1-1, to=3-1]
	\arrow["{\mu_{a_0}}", from=1-3, to=3-3]
	\arrow["{\gamma_{*_1}}"', from=3-1, to=3-3]
\end{tikzcd}\]

\noindent The colax $T$-algebra unit axiom says that $\gamma_{*_2}$ is the identity. From this we obtain that $\mu_{a_0} = \gamma_{*_1}$, and in particular $\widehat{\mu}_{a_0}$ is the identity because of the $A_0$-coalgebra counit axiom for $(a_0, \gamma_{*_1})$. By the observation above this proposition, $a_0$ is the initial object of $\ca$.

``$\Leftarrow$”: Conversely, given a comonad $(A,t,\delta,\iota)$ and an intial object $\emptyset \in \ca$, define the functor $A: T\ca \to \ca$ on objects and morphisms as follows:
\[\begin{tikzcd}
	{*} & b & c & \mapsto & \emptyset & tb & tc
	\arrow["{!_b}", from=1-1, to=1-2]
	\arrow["f", from=1-2, to=1-3]
	\arrow["{!_{tb}}", from=1-5, to=1-6]
	\arrow["tf", from=1-6, to=1-7]
\end{tikzcd}\]
Also, define the natural transformation $\gamma: A m_\ca \Rightarrow A TA$ on components by:
\begin{align*}
\gamma_{*_2} &:= 1_{\emptyset},\\
\gamma_{*_1} &:=\hspace{0.13cm}!_{t\emptyset} : \emptyset \to t\emptyset,\\
\gamma_{x} &:= \delta_x : tx \to t^2x \hspace{1cm} \text{for } x \in \ca.
\end{align*}
It is now a routine exercise to verify that $(\ca,A,\gamma,\iota)$ is a colax $T$-algebra so we omit the rest of the proof.
\end{proof}
\end{pr}

\begin{pr}[The 2-monad for weights]\label{PR_reslan2monad}
Let $\cj$ be a small 2-category and denote by $\iota: \ob \cj \to \cj$ the inclusion 2-functor. Recall that there is a 2-adjunction with the left 2-adjoint being given by the left Kan extension along $\iota$ and the right 2-adjoint being the restriction:
\[\begin{tikzcd}
	{[\ob \cj, \cat]} && {[\cj, \cat]}
	\arrow[""{name=0, anchor=center, inner sep=0}, "{\text{Lan}_\iota}"', curve={height=24pt}, from=1-1, to=1-3]
	\arrow[""{name=1, anchor=center, inner sep=0}, "{\text{res}}"', curve={height=24pt}, from=1-3, to=1-1]
	\arrow["\dashv"{anchor=center, rotate=90}, draw=none, from=0, to=1]
\end{tikzcd}\]
Given a collection $W: \ob \cj \to \cat$ of small categories, the 2-functor $\text{Lan}_\iota W : \cj \to \cat$ is the assignment:
\[
(b \in \cj) \mapsto \sum_{j \in \cj} \cj(j,b) \times Wj.
\]
This 2-adjunction generates a 2-monad $T$ on $[ \ob \cj, \cat]$. Since colimits are computed pointwise, $\text{res}: [\cj, \cat] \to [\ob \cj, \cat]$ creates them -- in particular it creates coequalizers of $\text{res}$-split pairs. It is thus 2-monadic and we have $\talgs \simeq [ \cj, \cat]$. Moreover, it can be verified that:
\begin{itemize}
\item a pseudo/lax/colax $T$-algebra is a pseudo/lax/colax functor $\cj \to \cat$,
\item a pseudo/lax/colax morphism is a pseudo/lax/colax natural transformation,
\item an algebra 2-cell is a modification.
\end{itemize}
For a better memorability, we will call this 2-monad the \textit{reslan 2-monad} in the thesis. In case $\cj = *$, we obtain the identity 2-monad of Example \ref{PR_identity2monad}. Let us also point out that this 2-monad may be used to prove the lax version of the Yoneda Lemma:

\begin{theo}[Yoneda lemma for lax functors]\label{THM_Yonedalemma_lax_funs}
For any 2-category $\cj$, $j \in \cj$ and a lax functor $W: \cj \to \cat$ there is an adjunction:
\[\begin{tikzcd}
	{\laxhoml{\cj(j,-)}{W}\phantom{.................}} && Wj
	\arrow[""{name=0, anchor=center, inner sep=0}, "\Phi", curve={height=-24pt}, from=1-1, to=1-3]
	\arrow[""{name=1, anchor=center, inner sep=0}, "\Psi", curve={height=-24pt}, from=1-3, to=1-1]
	\arrow["\dashv"{anchor=center, rotate=90}, draw=none, from=1, to=0]
\end{tikzcd}\]
The right adjoint is the assignment $(\alpha: \cj(j,-) \Rightarrow W) \mapsto \alpha_j(1_j) \in Wj$. The left adjoint sends an object $x \in Wj$ to a lax natural transformation:
\begin{align*}
\Psi^x &: \cj(j,-) \Rightarrow W, \\
(\Psi^x)_k &: (f: j \to k) \mapsto Wf(x) \in Wk.
\end{align*}
\end{theo}
\begin{proof}
Consider a family $X: \ob \cj \to \cat$ given by:
\[   
Xk = 
     \begin{cases}
       * &\quad\text{if } k = j,\\
       \emptyset &\quad \text{otherwise.} \\
     \end{cases}
\]
We can see that the reslan 2-monad sends $X$ to the 2-functor $TX = \cj(j,-)$. Also, for this choice of $X$, we have $[\ob \cj, \cat](X,W) = Wj$. Theorem now follows from Remark \ref{POZN_homcat_adj_laxtalgl_base}.
\end{proof}

Analogous arguments produce Yoneda lemmas for colax functors, pseudofunctors, 2-functors and the like.
\end{pr}

\begin{pr}[Free strict monoidal category 2-monad]\label{PR_freestrictmoncat2monad}
There is a 2-monad $T$ on $\cat$ given by the assignment $\ca \mapsto \sum_{n \geq 0} \ca^n$.
\begin{itemize}
\item A strict $T$-algebra is a strict monoidal category: a category $\ca$ equipped with a functor $\otimes: \ca \times \ca \to \ca$ and an object $I \in \ca$ such that $\otimes$ is strictly associative and unital,
\item A strict/pseudo/lax/colax morphism is an appropriate version of a monoidal functor,
\item A (co)lax $T$-algebra is what is called (an unbiased) \textit{(co)lax monoidal category}: it is a category $\ca$ equipped with:
\begin{itemize}
\item a functor $\otimes_n: \ca^n \to \ca$ for every $n \geq 0$. We denote:
\begin{align*}
a_1 \otimes \dots \otimes a_n &:= \otimes_n (a_1, \dots, a_n),\\
[a] &:= \otimes_1 (a),\\
I &:= \otimes_0 (*).
\end{align*}
\item for every object of $T^2\ca$ the \textit{coassociator} 1-cell that ``adds brackets”, for instance:
\begin{align*}
\gamma_{((a_1, a_2),(a_3,a_4))} &: a_1 \otimes a_2 \otimes a_3 \otimes a_4 \to (a_1 \otimes a_2) \otimes (a_3 \otimes a_4),\\
\gamma_{((a_1),(),(a_2, a_3))} &: a_1 \otimes a_2 \otimes a_3 \to [a_1] \otimes I \otimes (a_2 \otimes a_3).
\end{align*}
\item for every object of $\ca$ the \textit{counitor 1-cell}:
\[
\iota_a : [a] \to a.
\]
\end{itemize}
These are subject to the coassociativity and counit laws. Lax monoidal categories have been mentioned in a greater detail in the book \cite[Section 3.1]{higheroperads}. They have also been studied in \cite{algsofhigheroperads_HIST} under the name \textit{multitensors} -- see therein for their relationship with operads.

\item A pseudo $T$-algebra is typically called an \textit{unbiased monoidal category}. The 2-category $\pstalg$ can be shown to be 2-equivalent to the 2-category of ordinary monoidal categories, see \cite[Appendix B]{higheroperads},
\end{itemize}
\end{pr}

\begin{pr}[Free symmetric strict monoidal category 2-monad]\label{PR_freesymstrictmoncat2monad}
There is a 2-monad $S$ on $\cat$ sending a category $\ca$ to the category $S\ca$ whose:
\begin{itemize}
\item objects are tuples of objects of $\ca$,
\item a morphism $(a_1, \dots, a_n) \to (b_1, \dots, b_m)$ is a pair $(\phi, (f_1, \dots, f_n))$ of a bijection\\$\phi: n \to m$ and a tuple of morphisms $f_i : a_{\phi(i)} \to b_i$.
\end{itemize}
We have:
\begin{itemize}
\item Strict $S$-algebras are \textit{symmetric strict monoidal categories}: a strict monoidal category $(\ca, \otimes, I)$ equipped with an isomorphisms $\tau_{a,b}: a \otimes b \cong b \otimes a$ subject to axioms.
\item As has been observed in \cite{algsofhigheroperads_HIST}, a normal lax $S$-algebra is what has been called a \textit{functor-operad} in \cite[Definition 4.1]{cosimplicialobjectslittlencubes_HIST}.
\item We do not have a proof, but presumably the 2-category of pseudo $S$-algebras can be shown to be equivalent to ordinary symmetric monoidal categories, much like in the non-symmetric case in Example \ref{PR_freestrictmoncat2monad}.
\end{itemize}

\end{pr}

\begin{pr}[The 2-category 2-monad]\label{PR_2category_2monad}
Fix a set $X$. There is a 2-monad $T$ on $\cat^{X \times X}$ that sends an $(X \times X)$-indexed collection of categories:
\[
\cc := (\cc(A,B))_{(A,B) \in X \times X},
\]
to the collection of paths:
\[
(T\cc)_{A,B} = \text{Path}_\cc(A,B) = \{ (f_1,\dots, f_m) \mid m \in \mathbb{N}, \text{cod}(f_{i})= \text{dom}(f_{i+1}) \forall i < m \}.
\]
Here, we denote $\text{cod}(f) = B$, $\text{dom}(f) = A$ if $f \in \cc(A,B)$
It is readily verified that:
\begin{itemize}
\item a strict $T$-algebra is a 2-category with object-set $X$,
\item a pseudo $T$-algebra is (an unbiased) bicategory,
\item a strict morphism is an identity-on-objects 2-functor,
\item a lax morphism is an identity-on-objects lax functor.
\item an algebra 2-cell is an icon (see \textbf{Variations} for Definition \ref{DEFI_colax_nat_tr}).
\end{itemize}

\noindent Notice that for $X = *$ this reduces to Example \ref{PR_freestrictmoncat2monad} -- ``bicategories are many-objects monoidal categories”.
\end{pr}

\begin{pr}[Fibration 2-monad]\label{PR_fibration_2monad}
Let $\ck$ be a 2-category with comma objects and fix an object $C \in \ck$. There is a 2-monad $P$ on $\ck / C$ that sends a morphism $f: A \to B$ to the morphism $\pi_f : Pf \to C$ which is a projection of the following comma object in $\ck$:
\begin{equation}\label{EQ_comma_object_lonely}
\begin{tikzcd}
	Pf && A \\
	C && C
	\arrow[""{name=0, anchor=center, inner sep=0}, "{\rho_f}", from=1-1, to=1-3]
	\arrow[""{name=0p, anchor=center, inner sep=0}, phantom, from=1-1, to=1-3, start anchor=center, end anchor=center]
	\arrow["{\pi_f}"', from=1-1, to=2-1]
	\arrow["f", from=1-3, to=2-3]
	\arrow[""{name=1, anchor=center, inner sep=0}, Rightarrow, no head, from=2-1, to=2-3]
	\arrow[""{name=1p, anchor=center, inner sep=0}, phantom, from=2-1, to=2-3, start anchor=center, end anchor=center]
	\arrow["\chi", shorten <=4pt, shorten >=4pt, Rightarrow, from=0p, to=1p]
\end{tikzcd}
\end{equation}
We have:
\begin{itemize}
\item a pseudo $P$-algebra is a fibration internal to $\ck$ (see \cite[Proposition 9]{streetfib}),
\item a \textit{split fibration} internal to a 2-category has been defined in \cite{streetfib} to be a strict $P$-algebra. For the case $\ck = \cat$, these are the usual split fibrations (fibrations equipped with a choice of cartesian lifts that contains identities and is closed under composition), the full proof is offered in \cite[Proposition 9.3.17, Proposition 9.4.22]{2dimensionalcats}.
\end{itemize}
\end{pr}

\medskip

\noindent The following examples will be introduced later:
\begin{itemize}
\item The double category 2-monad (Example \ref{PR_doublecat_2monad}),
\item The free cocompletion pseudomonad (Example \ref{PR_presheafpsmonad}),
\item Cauchy completion pseudomonad (Example \ref{PR_ps_idemp_psmonads}).
\end{itemize}

\section{Lax-idempotent pseudomonads}\label{SEC_lax_idemp_psmonads}

The process of categorification leads one to replace ``properties” by ``property-like structure” -- structure unique up to an invertible 2-cell. A particular class of pseudomonads giving an object a property-like structure has been isolated by Anders Kock \cite{monadsforwhichstradjtounit_HIST} and Volker Zöberlein \cite{zoberleindoctrines_HIST} -- such pseudomonads have later been called by an adjective \textit{lax-idempotent} or \textit{KZ}. They are a categorification of \textit{idempotent monads} \cite[Proposition 4.2.3]{handbook2} and reduce to them if the 2-category on which they are defined is locally discrete. For a comprehensive study of them we also refer the reader to \cite{onpropertylike}.

\medskip

\noindent \textbf{The section is organized as follows}:
\begin{itemize}
\item Section \ref{SUBS_lax_idemp_defi_pr}: We recall the definition and relevant examples.
\item Section \ref{SUBS_lax_idemp_props}: We use the construction from Section \ref{SUBS_psmonads_adjsvsalgs} to derive various properties of lax-idempotent pseudomonads.
\item Section \ref{SECS_left_kan_psmons} recalls \textit{Left kan pseudomonads} of \cite{kanext} and in Section \ref{SUBS_leftKan2monads} we focus on their special simple case that we call \textit{left Kan 2-monads}.
\end{itemize}

\subsection{Definition and examples}\label{SUBS_lax_idemp_defi_pr}

\begin{theo}\label{THM_kz_monad_TFAE}
The following are equivalent for a pseudomonad $(T,m,i,\mu,\eta,\beta)$ on a 2-category $\ck$:
\begin{itemize}
\item there exists a modification $\Phi: 1_{T^2} \to iT \circ m$ so that $(\beta, \Phi): m \dashv iT$ is an adjunction in $\Hom{\ck}{\ck}$,
\item there exists a modification $\Psi: Ti \circ m \to 1_{T^2}$ so that $(\Psi, \eta): Ti \dashv m$ is an adjunction in $\Hom{\ck}{\ck}$,
\item there is a modification $\lambda : Ti \to iT$ satisfying:
\begin{equation}\label{EQ_KZ_psmonad_lambdaAiA}
\begin{tikzcd}[column sep=scriptsize]
	&&&&&&& TA \\
	A && TA && {T^2A} & {=} & A && {T^2A} \\
	&&&&&&& TA
	\arrow["{Ti_A}", from=1-8, to=2-9]
	\arrow["{i^{-1}_{i_A}}", shorten <=6pt, shorten >=6pt, Rightarrow, from=1-8, to=3-8]
	\arrow["{i_A}", from=2-1, to=2-3]
	\arrow[""{name=0, anchor=center, inner sep=0}, "{Ti_A}", curve={height=-24pt}, from=2-3, to=2-5]
	\arrow[""{name=1, anchor=center, inner sep=0}, "{i_{TA}}"', curve={height=24pt}, from=2-3, to=2-5]
	\arrow["{i_A}", from=2-7, to=1-8]
	\arrow["{i_A}"', from=2-7, to=3-8]
	\arrow["{i_{TA}}"', from=3-8, to=2-9]
	\arrow["{\lambda_A}", shorten <=6pt, shorten >=6pt, Rightarrow, from=0, to=1]
\end{tikzcd}
\end{equation}
and:
\begin{equation}\label{EQ_KZ_psmonad_mAlambdaA}
\begin{tikzcd}[column sep=scriptsize]
	&&&&&&& {T^2A} \\
	TA && {T^2A} && TA & {=} & TA && TA \\
	&&&&&&& {T^2A}
	\arrow["{m_A}", from=1-8, to=2-9]
	\arrow[""{name=0, anchor=center, inner sep=0}, "{Ti_A}", curve={height=-24pt}, from=2-1, to=2-3]
	\arrow[""{name=1, anchor=center, inner sep=0}, "{i_{TA}}"', curve={height=24pt}, from=2-1, to=2-3]
	\arrow["{m_A}", from=2-3, to=2-5]
	\arrow["{Ti_A}", from=2-7, to=1-8]
	\arrow[""{name=2, anchor=center, inner sep=0}, Rightarrow, no head, from=2-7, to=2-9]
	\arrow["{i_{TA}}"', from=2-7, to=3-8]
	\arrow["{m_A}"', from=3-8, to=2-9]
	\arrow["{\eta_A^{-1}}", shorten >=3pt, Rightarrow, from=1-8, to=2]
	\arrow["{\lambda_A}", shorten <=6pt, shorten >=6pt, Rightarrow, from=0, to=1]
	\arrow["{\beta^{-1}_A}", shorten <=3pt, Rightarrow, from=2, to=3-8]
\end{tikzcd}
\end{equation}

\end{itemize}
\end{theo}
\begin{proof}
The equivalence of $1)$ and $2)$ has been shown for instance in \cite[Lemma 1.1.12]{elephant}. For ``$1) \Rightarrow 3)$”, one takes the 2-cell $\lambda_A$ to be the following:
\[\begin{tikzcd}
	TA && {T^2A} && {T^2A} \\
	&&& TA
	\arrow["{Ti_A}", from=1-1, to=1-3]
	\arrow[""{name=0, anchor=center, inner sep=0}, curve={height=18pt}, Rightarrow, no head, from=1-1, to=2-4]
	\arrow[""{name=1, anchor=center, inner sep=0}, Rightarrow, no head, from=1-3, to=1-5]
	\arrow["{m_A}"', from=1-3, to=2-4]
	\arrow["{i_{TA}}"', from=2-4, to=1-5]
	\arrow["{\Phi_A}", shorten <=3pt, Rightarrow, from=1, to=2-4]
	\arrow["{\eta_A^{-1}}"', shorten >=4pt, Rightarrow, from=1-3, to=0]
\end{tikzcd}\]
For ``$3) \Rightarrow 1)$”, one takes the 2-cell $\Phi_A$ to be the 2-cell:
\[\begin{tikzcd}
	{T^2A} && {T^3A} && {T^2A} \\
	\\
	&& TA
	\arrow[""{name=0, anchor=center, inner sep=0}, "{i_{T^2A}}"', curve={height=24pt}, from=1-1, to=1-3]
	\arrow[""{name=1, anchor=center, inner sep=0}, "{Ti_{TA}}", from=1-1, to=1-3]
	\arrow[""{name=2, anchor=center, inner sep=0}, shift left, curve={height=-30pt}, Rightarrow, no head, from=1-1, to=1-5]
	\arrow["{m_A}"', curve={height=18pt}, from=1-1, to=3-3]
	\arrow["{Tm_A}"', from=1-3, to=1-5]
	\arrow["{i^{-1}_{m_A}}", shorten <=6pt, shorten >=6pt, Rightarrow, from=1-3, to=3-3]
	\arrow["{i_{TA}}"', curve={height=18pt}, from=3-3, to=1-5]
	\arrow["{\lambda_{TA}}", shorten <=3pt, shorten >=3pt, Rightarrow, from=1, to=0]
	\arrow["{T\beta^{-1}_A}"{pos=0.7}, shorten <=3pt, Rightarrow, from=2, to=1-3]
\end{tikzcd}\]
\end{proof}

\begin{defi}
A pseudomonad $(T,m,i,\mu,\eta,\beta)$ on a 2-category $\ck$ satisfying the above will be referred to as \textit{lax-idempotent}.
\end{defi}

For short, we will denote a lax-idempotent pseudomonad by just the quadruple \[ (T,m,i,\lambda), \] with $\lambda$ being a 2-cell satisfying \eqref{EQ_KZ_psmonad_lambdaAiA}, \eqref{EQ_KZ_psmonad_mAlambdaA}.

\begin{vari}
There are the following variations:
\begin{itemize}
\item in case we have $iT \dashv m$, $T$ is called \textit{colax-idempotent},
\item in case $m \dashv iT$ is an adjoint equivalence, $T$ is called \textit{pseudo-idempotent},
\end{itemize}
\end{vari}

\begin{pozn}[Duals]\label{POZN_duals}
Recall Definition \ref{DEFI_lali_rali_rari_lari} and note that a 1-cell $f$ is a reflector (a lali) in a 2-category $\ck$ if and only if:
\begin{itemize}
\item it is a reflection-inclusion (a rari) in $\ck^{op}$,
\item it is a coreflector (a rali) in $\ck^{co}$,
\item it is a coreflection-inclusion (a lari) in $\ck^{coop}$,
\end{itemize}
we can in particular see that a lax-idempotent pseudomonad $T$ on a 2-category $\ck$ is equivalently:
\begin{itemize}
\item a colax-idempotent pseudomonad $T^{co}$ on $\ck^{co}$,
\item a colax-idempotent pseudo-comonad $T^{op}$ on $\ck^{op}$,
\item a lax-idempotent pseudo-comonad $T^{coop}$ on $\ck^{coop}$.
\end{itemize}
\end{pozn}

\begin{pr}
Examples of lax/colax-idempotent pseudomonads include:
\begin{itemize}
\item Any idempotent 2-monad (where $m: T^2 \Rightarrow m$ is invertible) is lax-idempotent. For instance the identity 2-monad or the terminal 2-monad from Examples \ref{PR_identity2monad} and \ref{PR_terminal2monad}.

\item The initial object 2-monad from Example \ref{PR_freeinitialobject} is lax-idempotent. More generally, pseudomonads that freely adjoin colimits to a category are lax-idempotent. We will meet pseudomonads of this kind in Examples \ref{PR_presheafpsmonad} and \ref{PR_other_cocompletion_psmonads}.

\item The fibration 2-monad from Example \ref{PR_fibration_2monad} is colax-idempotent, see \cite[Proposition 9]{streetfib}.

\item The lax morphism classifier 2-comonad $Q_l$ on the 2-category $\talgs$ for a 2-monad $T$ is usually lax-idempotent. It will be studied in Section \ref{SEC_codescent_objects} and later in chapters \ref{CH_COMPUTATION} and \ref{CH_colax_adjs_lax_idemp_psmons}.
\end{itemize}
\end{pr}

\begin{pr}\label{PR_ps_idemp_psmonads}
Pseudo-idempotent pseudomonads are relatively rare. The motivating examples are the following:

\begin{itemize}
\item There is a pseudomonad on $\cat$ sending a small category $\cc$ to its \textit{Cauchy completion} $\overline{\cc}$ (see \cite[Definition 6.5.8]{handbook}) -- it can be given as the full subcategory of the presheaf category $[\cc^{op},\set]$ spanned by retracts of representables. The unit is given by the co-restricted Yoneda embedding $y: \cc \to \overline{\cc}$. A category is a pseudo algebra if and only if it admits \textit{splittings of idempotents}, i.e. equalizers of pairs consisting of an endomorphism and the identity 1-cell:
\[\begin{tikzcd}
	A && A
	\arrow["{1_A}"', shift right=2, from=1-1, to=1-3]
	\arrow["e", shift left=2, from=1-1, to=1-3]
\end{tikzcd}\]
Since by \cite[Proposition 6.5.9]{handbook}, a category is Cauchy complete if and only if the functor $y_\cc : \cc \to \overline{\cc}$ is an equivalence, this pseudomonad is pseudo-idempotent.
\item The pseudo morphism classifier 2-comonad $Q_p$ on the 2-category $\talgs$ for a 2-monad $T$ is typically pseudo-idempotent -- for instance if $\ck$ admits pseudo limits of arrows, it follows from \cite[Theorem 4.2]{twodim}.
\end{itemize}

\end{pr}

\subsection{Properties}\label{SUBS_lax_idemp_props}

\begin{lemma}\label{THML_lambda_unique_2cells_equal}
Given a lax-idempotent pseudomonad $(T,m,i,\lambda)$ on a 2-category $\ck$, we have the following:
\begin{itemize}
\item the modification $\lambda: Ti \to iT$ is uniquely determined by the axiom \eqref{EQ_KZ_psmonad_mAlambdaA},
\item the following 2-cells are equal:
\begin{equation}\label{EQ_lambdas_motylek}
\begin{tikzcd}[column sep=small]
	&&&&&&&& {T^3X} \\
	&& {T^3X} \\
	{T^2X} &&&& {T^2X} & {T^2X} & TX && {T^2X} \\
	&& {T^3X} \\
	&&&&&&&& {T^3X}
	\arrow["{m_{TX}}", curve={height=-18pt}, from=1-9, to=3-9]
	\arrow["{m_{TX}}", from=2-3, to=3-5]
	\arrow[""{name=0, anchor=center, inner sep=0}, "{Ti_{TX}}"{description}, from=3-1, to=2-3]
	\arrow[""{name=1, anchor=center, inner sep=0}, "{T^2i_X}", curve={height=-30pt}, from=3-1, to=2-3]
	\arrow[""{name=2, anchor=center, inner sep=0}, Rightarrow, no head, from=3-1, to=3-5]
	\arrow[""{name=3, anchor=center, inner sep=0}, "{i_{T^2X}}"', curve={height=30pt}, from=3-1, to=4-3]
	\arrow[""{name=4, anchor=center, inner sep=0}, "{Ti_{TX}}"{description}, from=3-1, to=4-3]
	\arrow["{T^2i_X}", curve={height=-12pt}, from=3-6, to=1-9]
	\arrow["{m_X}"', from=3-6, to=3-7]
	\arrow["{i_{T^2X}}"', curve={height=12pt}, from=3-6, to=5-9]
	\arrow[""{name=5, anchor=center, inner sep=0}, "{i_{TX}}"{description}, curve={height=18pt}, from=3-7, to=3-9]
	\arrow[""{name=5p, anchor=center, inner sep=0}, phantom, from=3-7, to=3-9, start anchor=center, end anchor=center, curve={height=18pt}]
	\arrow[""{name=5p, anchor=center, inner sep=0}, phantom, from=3-7, to=3-9, start anchor=center, end anchor=center, curve={height=18pt}]
	\arrow[""{name=6, anchor=center, inner sep=0}, "{Ti_X}"{description}, curve={height=-18pt}, from=3-7, to=3-9]
	\arrow[""{name=6p, anchor=center, inner sep=0}, phantom, from=3-7, to=3-9, start anchor=center, end anchor=center, curve={height=-18pt}]
	\arrow[""{name=6p, anchor=center, inner sep=0}, phantom, from=3-7, to=3-9, start anchor=center, end anchor=center, curve={height=-18pt}]
	\arrow["{Tm_X}"', from=4-3, to=3-5]
	\arrow["{Tm_X}"', curve={height=18pt}, from=5-9, to=3-9]
	\arrow["{m_{i_X}}", shorten <=7pt, shorten >=10pt, Rightarrow, from=1-9, to=6p]
	\arrow["{\eta_{TX}^{-1}}", shorten <=3pt, shorten >=3pt, Rightarrow, from=2-3, to=2]
	\arrow["{T\lambda_X}", shorten <=4pt, shorten >=4pt, Rightarrow, from=1, to=0]
	\arrow["{T\beta_X^{-1}}", shorten <=3pt, shorten >=3pt, Rightarrow, from=2, to=4-3]
	\arrow["{\lambda_{TX}}", shorten <=4pt, shorten >=4pt, Rightarrow, from=4, to=3]
	\arrow["{\lambda_X}", shorten <=5pt, shorten >=5pt, Rightarrow, from=6p, to=5p]
	\arrow["{i_{m_X}}", shorten <=10pt, shorten >=7pt, Rightarrow, from=5p, to=5-9]
\end{tikzcd}
\end{equation}
\end{itemize}
\end{lemma}
\begin{proof}
Let $\sigma : Ti \to iT$ be a different modification satisfying \eqref{EQ_KZ_psmonad_mAlambdaA}. By taking mates, the 2-cells $\sigma_A, \lambda_A$ are equal if and only if the 2-cells below are equal, which is the case by \eqref{EQ_KZ_psmonad_mAlambdaA}:
\[\begin{tikzcd}[column sep=scriptsize]
	TA && {T^2A} && TA & TA && {T^2A} && TA
	\arrow[""{name=0, anchor=center, inner sep=0}, "{i_{TA}}"', from=1-1, to=1-3]
	\arrow[""{name=0p, anchor=center, inner sep=0}, phantom, from=1-1, to=1-3, start anchor=center, end anchor=center]
	\arrow[""{name=1, anchor=center, inner sep=0}, "{Ti_A}", curve={height=-24pt}, from=1-1, to=1-3]
	\arrow[""{name=1p, anchor=center, inner sep=0}, phantom, from=1-1, to=1-3, start anchor=center, end anchor=center, curve={height=-24pt}]
	\arrow[""{name=2, anchor=center, inner sep=0}, curve={height=24pt}, Rightarrow, no head, from=1-1, to=1-5]
	\arrow[""{name=2p, anchor=center, inner sep=0}, phantom, from=1-1, to=1-5, start anchor=center, end anchor=center, curve={height=24pt}]
	\arrow["{m_A}", from=1-3, to=1-5]
	\arrow[""{name=3, anchor=center, inner sep=0}, "{Ti_A}", curve={height=-24pt}, from=1-6, to=1-8]
	\arrow[""{name=3p, anchor=center, inner sep=0}, phantom, from=1-6, to=1-8, start anchor=center, end anchor=center, curve={height=-24pt}]
	\arrow[""{name=4, anchor=center, inner sep=0}, "{i_{TA}}"', from=1-6, to=1-8]
	\arrow[""{name=4p, anchor=center, inner sep=0}, phantom, from=1-6, to=1-8, start anchor=center, end anchor=center]
	\arrow[""{name=5, anchor=center, inner sep=0}, curve={height=24pt}, Rightarrow, no head, from=1-6, to=1-10]
	\arrow[""{name=5p, anchor=center, inner sep=0}, phantom, from=1-6, to=1-10, start anchor=center, end anchor=center, curve={height=24pt}]
	\arrow["{m_A}", from=1-8, to=1-10]
	\arrow["{\lambda_A}", shorten <=3pt, shorten >=3pt, Rightarrow, from=1p, to=0p]
	\arrow["{\beta_A}", shorten >=2pt, Rightarrow, from=1-3, to=2p]
	\arrow["{\sigma_A}", shorten <=3pt, shorten >=3pt, Rightarrow, from=3p, to=4p]
	\arrow["{\beta_A}", shorten >=2pt, Rightarrow, from=1-8, to=5p]
\end{tikzcd}\]
To prove the \textbf{next claim}, denote the left 2-cell in \eqref{EQ_lambdas_motylek} by ``\text{LHS}”. The 2-cells in \eqref{EQ_lambdas_motylek} are equal if and only if:
\[
i_{m_X}^{-1} \circ \text{LHS} \circ m_{i_X}^{-1} = \lambda_X m_X.
\]
It suffices to prove this after taking mates under an adjunction $m_X \dashv i_{TX}$, i.e. to prove the equality:
\[
i_{TX} \beta_X \circ (i_{m_X}^{-1} \circ \text{LHS} \circ m_{i_X}^{-1})i_{TX} = i_{TX} \beta_X \circ \lambda_X m_X i_{TX}.
\]

\noindent The ``lower part” of the new left-hand side is the following diagram:
\[\begin{tikzcd}
	TX && {T^2X} && {T^3X} && {T^2X} \\
	\\
	&&&& TX
	\arrow["{i_{TX}}", from=1-1, to=1-3]
	\arrow[""{name=0, anchor=center, inner sep=0}, curve={height=30pt}, Rightarrow, no head, from=1-1, to=3-5]
	\arrow[""{name=1, anchor=center, inner sep=0}, "{i_{T^2X}}"', curve={height=24pt}, from=1-3, to=1-5]
	\arrow[""{name=2, anchor=center, inner sep=0}, "{Ti_{TX}}", from=1-3, to=1-5]
	\arrow[""{name=3, anchor=center, inner sep=0}, shift left, curve={height=-30pt}, Rightarrow, no head, from=1-3, to=1-7]
	\arrow["{m_X}"{description}, curve={height=18pt}, from=1-3, to=3-5]
	\arrow["{Tm_X}"', from=1-5, to=1-7]
	\arrow["{i_{m_X}^{-1}}", shorten <=6pt, shorten >=6pt, Rightarrow, from=1-5, to=3-5]
	\arrow["{i_{TX}}"', from=3-5, to=1-7]
	\arrow["{T\beta_X^{-1}}", shorten <=3pt, shorten >=3pt, Rightarrow, from=3, to=1-5]
	\arrow["{\lambda_{TX}}", shorten <=3pt, shorten >=3pt, Rightarrow, from=2, to=1]
	\arrow["{\beta_X}"', shorten <=7pt, shorten >=7pt, Rightarrow, from=1-3, to=0]
\end{tikzcd}\]
It reduces, thanks to the axiom \eqref{EQ_KZ_psmonad_lambdaAiA}, to the following:
\[\begin{tikzcd}
	TX && {T^2X} && {T^3X} && {T^2X} \\
	&& {T^2X} \\
	&&&& TX
	\arrow["{i_{TX}}", from=1-1, to=1-3]
	\arrow["{i_{TX}}"', from=1-1, to=2-3]
	\arrow[""{name=0, anchor=center, inner sep=0}, shift right=2, curve={height=30pt}, Rightarrow, no head, from=1-1, to=3-5]
	\arrow["{Ti_{TX}}", from=1-3, to=1-5]
	\arrow[""{name=1, anchor=center, inner sep=0}, shift left, curve={height=-30pt}, Rightarrow, no head, from=1-3, to=1-7]
	\arrow["{i^{-1}_{i_{TX}}}", Rightarrow, from=1-3, to=2-3]
	\arrow["{Tm_X}"', from=1-5, to=1-7]
	\arrow["{i_{m_X}^{-1}}", shorten <=6pt, shorten >=6pt, Rightarrow, from=1-5, to=3-5]
	\arrow["{i_{T^2X}}"', from=2-3, to=1-5]
	\arrow["{m_X}"{description}, from=2-3, to=3-5]
	\arrow["{i_{TX}}"', from=3-5, to=1-7]
	\arrow["{T\beta_X^{-1}}", shorten <=3pt, shorten >=3pt, Rightarrow, from=1, to=1-5]
	\arrow["{\beta_X}"', shorten <=4pt, shorten >=4pt, Rightarrow, from=2-3, to=0]
\end{tikzcd}\]
Using the local naturality of $i: 1_{\ck} \Rightarrow T$, this becomes the identity on $i_{TX}$.

Consider now the ``upper part” of the new left-hand side:
\[\begin{tikzcd}
	&&&& TX \\
	\\
	TX && {T^2X} && {T^3X} && {T^2X}
	\arrow["{m_{i_X}^{-1}}", shorten <=6pt, shorten >=6pt, Rightarrow, from=1-5, to=3-5]
	\arrow["{Ti_X}", curve={height=-18pt}, from=1-5, to=3-7]
	\arrow["{i_{TX}}", from=3-1, to=3-3]
	\arrow["{m_X}", curve={height=-18pt}, from=3-3, to=1-5]
	\arrow[""{name=0, anchor=center, inner sep=0}, "{Ti_{TX}}"{description}, from=3-3, to=3-5]
	\arrow[""{name=1, anchor=center, inner sep=0}, "{T^2i_X}", curve={height=-30pt}, from=3-3, to=3-5]
	\arrow[""{name=2, anchor=center, inner sep=0}, curve={height=30pt}, Rightarrow, no head, from=3-3, to=3-7]
	\arrow["{m_{TX}}", from=3-5, to=3-7]
	\arrow["{T\lambda_X}", shorten <=4pt, shorten >=4pt, Rightarrow, from=1, to=0]
	\arrow["{\eta_{TX}^{-1}}"{pos=0.3}, shorten >=3pt, Rightarrow, from=3-5, to=2]
\end{tikzcd}\]
Using the local naturality of $i: 1_{\ck} \Rightarrow T$ at the 2-cell $\lambda_X$, we obtain:
\[\begin{tikzcd}
	&& TX & {} && TX \\
	TX && {T^2X} && {T^3X} && {T^2X} \\
	&& {T^2X}
	\arrow["{m_X}", from=1-3, to=1-6]
	\arrow["{i^{-1}_{Ti_X}}", shorten >=2pt, Rightarrow, from=1-3, to=2-3]
	\arrow["{T^2i_X}", from=1-3, to=2-5]
	\arrow[""{name=0, anchor=center, inner sep=0}, draw=none, from=1-4, to=1-6]
	\arrow["{Ti_X}", from=1-6, to=2-7]
	\arrow["{i_{TX}}", curve={height=-30pt}, from=2-1, to=1-3]
	\arrow[""{name=1, anchor=center, inner sep=0}, "{i_{TX}}"{description}, curve={height=24pt}, from=2-1, to=2-3]
	\arrow[""{name=2, anchor=center, inner sep=0}, "{Ti_X}"{description}, curve={height=-24pt}, from=2-1, to=2-3]
	\arrow["{i_{TX}}"', curve={height=30pt}, from=2-1, to=3-3]
	\arrow["{i_{T^2X}}"{description}, from=2-3, to=2-5]
	\arrow["{i_{i_{TX}}}", Rightarrow, from=2-3, to=3-3]
	\arrow["{m_{TX}}", from=2-5, to=2-7]
	\arrow["{Ti_{TX}}"', from=3-3, to=2-5]
	\arrow[""{name=3, anchor=center, inner sep=0}, curve={height=30pt}, Rightarrow, no head, from=3-3, to=2-7]
	\arrow["{m_{i_{X}}^{-1}}", shorten <=3pt, Rightarrow, from=0, to=2-5]
	\arrow["{\lambda_X}", shorten <=6pt, shorten >=6pt, Rightarrow, from=2, to=1]
	\arrow["{\eta_{TX}^{-1}}", shorten <=5pt, shorten >=5pt, Rightarrow, from=2-5, to=3]
\end{tikzcd}\]
Using Lemma \ref{THML_in_every_psmonad}, this reduces to:
\[\begin{tikzcd}
	&& TX && {} & TX \\
	TX && {T^2X} && {T^3X} && {T^2X}
	\arrow["{m_X}", from=1-3, to=1-6]
	\arrow["{i^{-1}_{Ti_X}}", shorten <=1pt, shorten >=1pt, Rightarrow, from=1-3, to=2-3]
	\arrow["{T^2i_X}"{description}, from=1-3, to=2-5]
	\arrow["{m_{i_{X}}^{-1}}", Rightarrow, from=1-5, to=2-5]
	\arrow["{Ti_X}", from=1-6, to=2-7]
	\arrow["{i_{TX}}", curve={height=-30pt}, from=2-1, to=1-3]
	\arrow[""{name=0, anchor=center, inner sep=0}, "{i_{TX}}"{description}, curve={height=24pt}, from=2-1, to=2-3]
	\arrow[""{name=1, anchor=center, inner sep=0}, "{Ti_X}"{description}, curve={height=-24pt}, from=2-1, to=2-3]
	\arrow["{i_{T^2X}}"{description}, from=2-3, to=2-5]
	\arrow[""{name=2, anchor=center, inner sep=0}, curve={height=24pt}, Rightarrow, no head, from=2-3, to=2-7]
	\arrow[""{name=3, anchor=center, inner sep=0}, Rightarrow, draw=none, from=2-3, to=2-7]
	\arrow["{m_{TX}}", from=2-5, to=2-7]
	\arrow["{\lambda_X}", shorten <=6pt, shorten >=6pt, Rightarrow, from=1, to=0]
	\arrow["{\beta_{TX}}", shorten <=3pt, shorten >=3pt, Rightarrow, from=3, to=2]
\end{tikzcd}\]
Finally, using the modification axiom of $\beta$ for the 1-cell $i_X$, this reduces to the following, which is exactly the new right-hand side, so the proof is complete:
\[\begin{tikzcd}
	& {T^2X} \\
	TX && TX && {T^2X}
	\arrow["{m_X}", from=1-2, to=2-3]
	\arrow["{i_{TX}}", from=2-1, to=1-2]
	\arrow[""{name=0, anchor=center, inner sep=0}, curve={height=24pt}, Rightarrow, no head, from=2-1, to=2-3]
	\arrow[""{name=1, anchor=center, inner sep=0}, "{i_{TX}}"{description}, curve={height=24pt}, from=2-3, to=2-5]
	\arrow[""{name=2, anchor=center, inner sep=0}, "{Ti_X}"{description}, curve={height=-24pt}, from=2-3, to=2-5]
	\arrow["{\beta_X}", shorten >=6pt, Rightarrow, from=1-2, to=0]
	\arrow["{\lambda_X}", shorten <=6pt, shorten >=6pt, Rightarrow, from=2, to=1]
\end{tikzcd}\]
\end{proof}

\begin{lemma}\label{THML_Phi_boo}
Let $(T,m,i)$ be a lax-idempotent pseudomonad on a 2-category $\ck$ and let $(A,a,\gamma,\iota)$ be a normal colax $T$-algebra ($\iota: ai_A \Rightarrow 1_A$ is invertible). Then there exists a unique 2-cell $\sigma: 1_{TA} \Rightarrow i_A a$ such that:
\[
(\iota, \sigma): a \dashv i_A : A \to TA.
\]
\end{lemma}
\begin{proof}
The \textbf{uniqueness part} follows from the general fact that the unit in an adjunction is determined by all the other data. To prove the \textbf{existence}, define:
\[\begin{tikzcd}
	{\sigma :=} & TA && {T^2A} && TA \\
	&&& A
	\arrow[""{name=0, anchor=center, inner sep=0}, "{Ti_A}", curve={height=-12pt}, from=1-2, to=1-4]
	\arrow[""{name=1, anchor=center, inner sep=0}, "{i_{TA}}"', curve={height=12pt}, from=1-2, to=1-4]
	\arrow[""{name=2, anchor=center, inner sep=0}, "{1_{TA}}", shift left=3, curve={height=-30pt}, from=1-2, to=1-6]
	\arrow["a"', curve={height=12pt}, from=1-2, to=2-4]
	\arrow["Ta", from=1-4, to=1-6]
	\arrow["{i_a^{-1}}", Rightarrow, from=1-4, to=2-4]
	\arrow["{i_A}"', curve={height=12pt}, from=2-4, to=1-6]
	\arrow["{T\iota^{-1}}", shorten <=4pt, Rightarrow, from=2, to=1-4]
	\arrow["{\lambda_A}", shorten <=3pt, shorten >=3pt, Rightarrow, from=0, to=1]
\end{tikzcd}\]
To prove the \textbf{first triangle identity}, $\iota^{-1} a = a \sigma$, note that by the colax $T$-algebra unit axiom and \eqref{EQ_KZ_psmonad_mAlambdaA}, the left-hand side equals the following:
\[\begin{tikzcd}
	&&&& TA \\
	TA && {T^2A} && TA && A \\
	&& A
	\arrow["\gamma", Rightarrow, from=1-5, to=2-5]
	\arrow["a", from=1-5, to=2-7]
	\arrow[""{name=0, anchor=center, inner sep=0}, curve={height=-30pt}, Rightarrow, no head, from=2-1, to=1-5]
	\arrow[""{name=1, anchor=center, inner sep=0}, "{Ti_A}", curve={height=-12pt}, from=2-1, to=2-3]
	\arrow[""{name=2, anchor=center, inner sep=0}, "{i_{TA}}"', curve={height=12pt}, from=2-1, to=2-3]
	\arrow["a"', curve={height=12pt}, from=2-1, to=3-3]
	\arrow["{m_A}", from=2-3, to=1-5]
	\arrow["Ta"', from=2-3, to=2-5]
	\arrow["{i_a^{-1}}", Rightarrow, from=2-3, to=3-3]
	\arrow["a"', from=2-5, to=2-7]
	\arrow["{i_A}"', curve={height=12pt}, from=3-3, to=2-5]
	\arrow["{\lambda_A}", shorten <=3pt, shorten >=3pt, Rightarrow, from=1, to=2]
	\arrow["{\eta_A}", shorten <=5pt, Rightarrow, from=0, to=2-3]
\end{tikzcd}\]
Using the other colax algebra unit axiom, we see that this is exactly $a \sigma$.
To prove the \textbf{second triangle identity}, we must show that $i_A \iota^{-1} = \sigma i_A$. By the local naturality of $i: 1_{\ck} \Rightarrow T$, the left-hand side equals the following:
\[\begin{tikzcd}
	& TA \\
	A && {T^2A} && TA \\
	& TA && A
	\arrow["{Ti_A}", from=1-2, to=2-3]
	\arrow[""{name=0, anchor=center, inner sep=0}, curve={height=-24pt}, Rightarrow, no head, from=1-2, to=2-5]
	\arrow["{i_{i_A}^{-1}}", shorten <=6pt, shorten >=6pt, Rightarrow, from=1-2, to=3-2]
	\arrow["{i_A}", from=2-1, to=1-2]
	\arrow["{i_A}"', from=2-1, to=3-2]
	\arrow[""{name=1, anchor=center, inner sep=0}, "Ta"{description}, from=2-3, to=2-5]
	\arrow["{i_{TA}}"', from=3-2, to=2-3]
	\arrow["a"', from=3-2, to=3-4]
	\arrow["{i_A}"', from=3-4, to=2-5]
	\arrow["{(T\iota)^{-1}}", shorten <=5pt, shorten >=5pt, Rightarrow, from=0, to=2-3]
	\arrow["{i_a^{-1}}"', shift right=5, shorten <=3pt, Rightarrow, from=1, to=3-4]
\end{tikzcd}\]
Using the axiom \eqref{EQ_KZ_psmonad_lambdaAiA}, this reduces to $\sigma i_A$.
\end{proof}

\noindent Recall the 2-functor $\Phi: \text{Adj}(i) \to \colaxtalgl$ from Theorem \ref{THM_2fun_phi_adj_colaxtalgl}. Denote by $\text{Refl}(i)$ its full sub-2-category spanned by ``reflectors of the unit” -- pairs $(A,(\iota,\sigma): a \dashv i_A)$ for which $\iota$ is invertible. The 2-functor clearly restricts to a 2-fully faithful 2-functor between $\text{Refl}(i)$ and $\ncolaxtalgl$.

The following result captures the fact that ``normal colax algebras are pseudo” for lax-idempotent pseudomonads:

\begin{prop}\label{THMP_normal_colax_algebra_for_lax_idemp_is_pseudo}
Let $(T,m,i)$ be a lax-idempotent pseudomonad and let $(A,a,\gamma,\iota)$ be a normal colax $T$-algebra ($\iota: ai_A \Rightarrow 1_A$ is invertible). Then it is a pseudo $T$-algebra ($\gamma$ is invertible).
\end{prop}
\begin{proof}
By Lemma \ref{THML_Phi_boo}, we know that $(A,a,\gamma,\iota) = \Phi((A,(\iota,\sigma): a \dashv i_A))$ for a unique 2-cell $\sigma$. By Remark \ref{POZN_Phi_sends_LAs_to_psmors}, $\Phi: \text{Refl}(i) \to \colaxtalgl$ sends the 1-cell $a : TA \to A$, regarded as a morphism:
\[
(TA, (\beta_A, \Phi_A) : m_A \dashv i_{TA}) \to (A,(\iota,\sigma): a \dashv i_A),
\]
to a pseudo morphism $(a,\overline{a}): (TA,m_A, \mu^{-1}_A, \beta_A) \to (A,a,\gamma,\iota)$. We will show that the composite $\overline{a} \circ \gamma$ is the identity on $a m_A$. Since $\overline{a}$ is invertible, we obtain that $\gamma = \overline{a}^{-1}$ is also invertible. By taking the mates, it suffices to show that the following equals $a \beta_A$:
\[\begin{tikzcd}
	&&&& TA \\
	TA && {T^2A} && TA && A \\
	&&&& TA
	\arrow["\gamma", Rightarrow, from=1-5, to=2-5]
	\arrow["a", from=1-5, to=2-7]
	\arrow["{i_{TA}}", from=2-1, to=2-3]
	\arrow[""{name=0, anchor=center, inner sep=0}, curve={height=18pt}, Rightarrow, no head, from=2-1, to=3-5]
	\arrow["{m_A}", from=2-3, to=1-5]
	\arrow["Ta"{description}, from=2-3, to=2-5]
	\arrow["{m_A}"', from=2-3, to=3-5]
	\arrow["a"{description}, from=2-5, to=2-7]
	\arrow["{\overline{a}}", Rightarrow, from=2-5, to=3-5]
	\arrow["a"', from=3-5, to=2-7]
	\arrow["{\beta_A}"', shorten >=3pt, Rightarrow, from=2-3, to=0]
\end{tikzcd}\]
Using the colax algebra unit axiom, from this we obtain the 2-cells portrayed above the dotted line in the diagram below. Using the pseudo morphism axiom for $(a,\overline{a})$ we obtain the 2-cells below the dotted line. We see that it then composes to $a \beta_A$:
\[\begin{tikzcd}
	& {T^2A} &&& A \\
	TA && TA && {T^2A} && TA && A \\
	&&&& A
	\arrow["{m_A}", from=1-2, to=2-3]
	\arrow["{i_a}", Rightarrow, from=1-5, to=2-5]
	\arrow["{i_A}"{description}, from=1-5, to=2-7]
	\arrow[""{name=0, anchor=center, inner sep=0}, curve={height=-24pt}, Rightarrow, no head, from=1-5, to=2-9]
	\arrow["{i_{TA}}", from=2-1, to=1-2]
	\arrow[""{name=1, anchor=center, inner sep=0}, Rightarrow, dashed, no head, from=2-1, to=2-3]
	\arrow["a", from=2-3, to=1-5]
	\arrow["{i_{TA}}"{description}, dashed, from=2-3, to=2-5]
	\arrow["a"', from=2-3, to=3-5]
	\arrow["Ta"{description}, dashed, from=2-5, to=2-7]
	\arrow["{i_a^{-1}}", Rightarrow, from=2-5, to=3-5]
	\arrow["a"{description}, dashed, from=2-7, to=2-9]
	\arrow["{i_A}"{description}, from=3-5, to=2-7]
	\arrow[""{name=2, anchor=center, inner sep=0}, curve={height=24pt}, Rightarrow, no head, from=3-5, to=2-9]
	\arrow["{\beta_A}", shorten >=3pt, Rightarrow, from=1-2, to=1]
	\arrow["{\iota^{-1}}"', shorten <=4pt, Rightarrow, from=0, to=2-7]
	\arrow["\iota", shorten >=4pt, Rightarrow, from=2-7, to=2]
\end{tikzcd}\]
\end{proof}

\noindent We moreover have:

\begin{theo}\label{THM_lax_idemp_psmonad_Phi}
Let $(T,m,i, \gamma, \eta, \beta)$ be a pseudomonad on a 2-category $\ck$. The following are equivalent:
\begin{itemize}
\item $T$ is lax-idempotent,
\item the restricted 2-functor:
\[
\Phi : \text{Refl}(i) \to \ncolaxtalgl,
\]
is an isomorphism of 2-categories.
\end{itemize}
Moreover, if this is the case, we have $\ncolaxtalgl = \pstalgl$.
\end{theo}
\begin{proof}
``$\Leftarrow$”: A free algebra $(TA,m_A, \mu^{-1}_A, \beta_A)$ lives in $\ncolaxtalgl$, so there is a unique invertible 2-cell $\Gamma_A: 1_{T^2A} \Rightarrow i_{TA} m_A$ for which we have:
\[\begin{tikzcd}
	{(\beta_A,\Gamma_A):} & TA && {T^2A} & {\text{ in } \ck}
	\arrow[""{name=0, anchor=center, inner sep=0}, "{i_{TA}}"', curve={height=24pt}, from=1-2, to=1-4]
	\arrow[""{name=1, anchor=center, inner sep=0}, "{m_A}"', curve={height=24pt}, from=1-4, to=1-2]
	\arrow["\dashv"{anchor=center, rotate=-90}, draw=none, from=1, to=0]
\end{tikzcd}\]
Since $m, iT$ are pseudonatural and $\beta$ is a modification, the collection $\Gamma_A$ for $A \in \ck$ can also be seen to be a modification due to formal arguments. Thus $T$ is lax-idempotent.

``$\Rightarrow$”: We already know that $\Phi$ is fully faithful and bijective on objects in its image. By Lemma \ref{THML_Phi_boo}, normal pseudo algebras are precisely its image.

The ``moreover” part has been proven in Proposition \ref{THMP_normal_colax_algebra_for_lax_idemp_is_pseudo}.
\end{proof}

\begin{cor}\label{THMC_lax_idemp_1cell_is_uniquely_lax}
Let $(T,m,i)$ be a lax-idempotent pseudomonad on a 2-category $\ck$ and let $(A,a,\gamma,\iota)$, $(B,b,\gamma',\iota')$ be two normal colax $T$-algebras. For any 1-cell $f:A \to B$ there exists a unique 2-cell $\overline{f}$ giving a lax morphism:
\[
(f,\overline{f}): (A,a,\gamma,\iota) \to (B,b,\gamma',\iota').
\]
\noindent Moreover, given two lax morphisms $(f,\overline{f}), (g,\overline{g})$ between normal colax $T$-algebras, any 2-cell $\alpha: f \Rightarrow g$ is an algebra 2-cell $(f,\overline{f}) \Rightarrow (g,\overline{g})$.
\end{cor}
\begin{proof}
This follows from Lemma \ref{THML_Tpsmonad_1cell_mezi_adjs_laxmor}.
\end{proof}

\noindent The structure of a normal colax algebra on an object $A$ is ``unique up to an invertible 2-cell”:

\begin{cor}
Let $(T,m,i)$ be a lax-idempotent pseudomonad on a 2-category $\ck$ and let $(A,a,\gamma,\iota)$, $(A,b,\gamma',\iota')$ be two normal colax $T$-algebras that share their underlying object. Then there exists a unique invertible 2-cell $\overline{\theta}: a \Rightarrow b$ for which the pair $(1_A, \overline{\theta})$ is an algebra isomorphism:
\[
(1_A, \overline{\theta}): (A,a,\gamma,\iota) \to (A,b,\gamma',\iota').
\]
\end{cor}
\begin{proof}
The 2-functor $\Phi : \text{Refl}(i) \to \ncolaxtalgl$ sends isomorphisms to isomorphisms (as does any 2-functor).
\end{proof}

We will also need a ``relative” version of lax-idempotency:

\begin{prop}\label{THMP_KZ_adjunkce}
Let $(\Psi, \Phi): (\epsilon, \eta): F \dashv U: \cd \to \cc$ be a biadjunction. The following are equivalent:
\begin{itemize}
\item the induced pseudomonad $UF$ on $\cc$ is lax-idempotent,
\item there is a modification $\Gamma: F\eta \circ \epsilon F \to 1_{FUF}$ such that $(\Gamma, \Psi^{-1}): F\eta \dashv \epsilon F$ is an adjunction in $\Hom{\cc}{\cd}$,
\item there is a modification $\Theta: 1_{UFU} \to \eta U \circ U \epsilon$ such that $(\Phi^{-1}, \Theta): U\epsilon \dashv \eta U$ is an adjunction in $\Hom{\cd}{\cc}$,
\item the induced pseudocomonad $FU$ on $\cd$ is colax-idempotent.
\end{itemize}
In this case we say that the biadjunction is \textit{lax-idempotent}.
\end{prop}
\begin{proof}
The strict version (involving 2-functors, 2-adjunctions) is proven in \cite{nlab_laxidemp_adj} -- the more general case is analogous. \textbf{A word of caution}: The terminology differs from ours -- the pseudocomonad $FU$ is referred to as ``colax-idempotent” in \cite{nlab_laxidemp_adj}.
\end{proof}

\subsection{Left Kan pseudomonads}\label{SECS_left_kan_psmons}

As we have seen, the notion of a lax-idempotent pseudomonad  contains a large amount of data and axioms. A major simplification can be achieved if one works with their \textit{left Kan pseudomonad} presentation of Marmolejo, Wood \cite{kanext}\footnote{The authors work with the dual version -- \textit{right Kan pseudomonads} -- that involves right Kan extensions.}. In this section we recall all the basic definitions and mention the equivalence of left Kan pseudomonads and lax-idempotent pseudomonads. We also define a special class of left Kan pseudomonads that we call \textit{left Kan 2-monads} -- this is the obvious strict version of the notion.

\begin{defi}
A \textit{left Kan pseudomonad} (\cite{kanext}) $(D,y)$ on a 2-category $\ck$ consists of:

\begin{itemize}
\item A function $D: \ob \ck \to \ob \ck$,
\item For every $A \in \ck$ a 1-cell $y_A: A \to DA$ called its \textit{unit},
\item For every 1-cell $f: A \to DB$ a left Kan extension of $f$ along $y_A$ such that the accompanying 2-cell is invertible:
\begin{equation}\label{EQ_left_kan_ext_mbbd}
\begin{tikzcd}
	A && DA \\
	\\
	&& DB
	\arrow["{f^\mbbd}", from=1-3, to=3-3]
	\arrow["{y_A}", from=1-1, to=1-3]
	\arrow[""{name=0, anchor=center, inner sep=0}, "f"', from=1-1, to=3-3]
	\arrow["{\mbbd_f}", shorten <=6pt, Rightarrow, from=0, to=1-3]
\end{tikzcd}
\end{equation}
\end{itemize}
These are subject to the axioms:
\begin{itemize}
\item For every $A \in \ck$, the identity 2-cell $1_{y_A}$ on $y_A$ exhibits $1_{DA}$ as the left Kan extension of $y_A$ along $y_A$:
\[
\begin{tikzcd}
	A && DA \\
	\\
	&& DA
	\arrow[Rightarrow, no head, from=1-3, to=3-3]
	\arrow["{y_A}", from=1-1, to=1-3]
	\arrow["{y_A}"', from=1-1, to=3-3]
\end{tikzcd}
\]

\item For every $g: B \to DC$ and $f: A \to DB$, the 1-cell $g^\mbbd$ \textit{preserves} the left Kan extension \eqref{EQ_left_kan_ext_mbbd}. This is to say that the following 2-cell exhibits $g^\mbbd f^\mbbd$ as the left Kan extension of $g^\mbbd f$ along $y_A$:
\[\begin{tikzcd}
	A && DA \\
	&& DB \\
	&& DC
	\arrow["{y_A}", from=1-1, to=1-3]
	\arrow[""{name=0, anchor=center, inner sep=0}, "f"', from=1-1, to=2-3]
	\arrow["{f^\mbbd}", from=1-3, to=2-3]
	\arrow["{g^\mbbd}", from=2-3, to=3-3]
	\arrow["{\mbbd_f}"', shorten <=5pt, Rightarrow, from=0, to=1-3]
\end{tikzcd}\]
\end{itemize}
\end{defi}

\begin{defi}\label{DEFI_D_algebra}
A \textit{pseudo $D$-algebra} consists of an object $C \in \ck$ together with a mapping that sends every 1-cell $g: B \to C$ to the left Kan extension of $g$ along $y_B$ such that the accompanying 2-cell is invertible:

\begin{equation}\label{EQ_left_kan_ext_mbbb}
\begin{tikzcd}
	B && DB \\
	\\
	&& C
	\arrow["{y_B}", from=1-1, to=1-3]
	\arrow[""{name=0, anchor=center, inner sep=0}, "g"', from=1-1, to=3-3]
	\arrow["{g^\mbbc}", from=1-3, to=3-3]
	\arrow["{\mbbc_g}", shorten <=6pt, Rightarrow, from=0, to=1-3]
\end{tikzcd}
\end{equation}

\noindent and such that for every $f: A \to DB$, $g^\mbbc$ preserves the left Kan extension \eqref{EQ_left_kan_ext_mbbd}.

A \textit{$D$-pseudomorphism} $h: C \to X$ between pseudo $D$-algebras $C,X$ is a 1-cell $h: C \to X$ that preserves the left Kan extension \eqref{EQ_left_kan_ext_mbbb}. An \textit{algebra 2-cell}: \[ \alpha: h \Rightarrow h' : C \to X, \] is just a 2-cell in $\ck$. It is easy to see that this data assembles into a 2-category that we will call the \textit{2-category of algebras} for the left Kan pseudomonad $D$.
\end{defi}

\begin{theo}\label{THM_korespondence_leftkanpsmon_laxidemppsmon}
There is a correspondence between:
\begin{itemize}
\item left Kan pseudomonads $(D,y)$ on $\ck$,
\item lax-idempotent pseudomonads $(D,m,y)$ on $\ck$.
\end{itemize}
Moreover, the left Kan pseudomonad and lax-idempotent pseudomonad corresponding to one another have biequivalent 2-categories of pseudo algebras and pseudo morphisms, and this biequivalence commutes with the forgetful 2-functors to $\ck$.

\end{theo}
\begin{proof}
For the full proof see \cite[4.1, 4.2]{kanext}, here we sketch only the bits relevant for this thesis. Given a left Kan pseudomonad $(D, y)$, the lax-idempotent pseudomonad is given by a normal pseudofunctor $D: \ck \to \ck$ with action on 1-cells and 2-cells given by $U_D \circ J_D$ from Proposition \ref{THMP_biadjunkce_ck_ck_D}. The components of the unit $y$ become pseudonatural with the pseudonaturality square given by the Kan extension 2-cell:
\[\begin{tikzcd}
	A && DA \\
	B && DB
	\arrow[""{name=0, anchor=center, inner sep=0}, "{y_A}", from=1-1, to=1-3]
	\arrow["{(y_Bf)^\mbbd =:Df}", from=1-3, to=2-3]
	\arrow["f"', from=1-1, to=2-1]
	\arrow[""{name=1, anchor=center, inner sep=0}, "{y_B}"', from=2-1, to=2-3]
	\arrow["\mbbd"', shorten <=4pt, shorten >=4pt, Rightarrow, from=1, to=0]
\end{tikzcd}\]
The multiplication at $A \in \ck$ is given by the morphism $p_{DA}$ again as in  Proposition \ref{THMP_biadjunkce_ck_ck_D}. On the other hand, given a lax-idempotent pseudomonad $(D,m,y)$, the left Kan extension of $f: A \to DB$ along $y_A : A \to DA$ is given by the composite of the pseudonaturality 2-cell $y_f$ and the pseudomonad unitor 2-cell:
\[\begin{tikzcd}
	A && DA \\
	DB && {D^2B} \\
	&& DB
	\arrow[""{name=0, anchor=center, inner sep=0}, "{y_A}", from=1-1, to=1-3]
	\arrow[""{name=0p, anchor=center, inner sep=0}, phantom, from=1-1, to=1-3, start anchor=center, end anchor=center]
	\arrow["f"', from=1-1, to=2-1]
	\arrow[""{name=1, anchor=center, inner sep=0}, "{y_{DB}}"{description}, from=2-1, to=2-3]
	\arrow[""{name=1p, anchor=center, inner sep=0}, phantom, from=2-1, to=2-3, start anchor=center, end anchor=center]
	\arrow[""{name=1p, anchor=center, inner sep=0}, phantom, from=2-1, to=2-3, start anchor=center, end anchor=center]
	\arrow["Df", from=1-3, to=2-3]
	\arrow["{m_B}", from=2-3, to=3-3]
	\arrow[""{name=2, anchor=center, inner sep=0}, Rightarrow, no head, from=2-1, to=3-3]
	\arrow[""{name=2p, anchor=center, inner sep=0}, phantom, from=2-1, to=3-3, start anchor=center, end anchor=center]
	\arrow["{y_f}", shorten <=4pt, shorten >=4pt, Rightarrow, from=0p, to=1p]
	\arrow["\cong"{pos=0.6}, shift left=5, shorten <=3pt, Rightarrow, from=1p, to=2p]
\end{tikzcd}\]
\end{proof}

\noindent In the rest of the thesis we will use the terms ``left Kan pseudomonad” and ``lax-idempotent pseudomonad” interchangeably. We will also use the same notation for the 2-category of algebras -- $\psdalg$ -- regardless whether $D$ is a left Kan pseudomonad or its corresponding lax-idempotent pseudomonad. This will pose no threat, since by Theorem \ref{THM_korespondence_leftkanpsmon_laxidemppsmon} they are biequivalent.

The two statements below are the ``left Kan” versions of well-known statements about lax-idempotent pseudomonads -- we will need these versions in Chapter \ref{CH_colax_adjs_lax_idemp_psmons}.

\begin{defi}\label{DEFI_Kleisli2cat}
By the \textit{Kleisli 2-category} $\ck_D$ associated to the left Kan pseudomonad $(D,y)$ we mean the full sub-2-category of $\psdalg$ spanned by \textit{free $D$-algebras}, that is, algebras whose underlying object is of form $DA$ for some object $A \in \ck$ and the extension operation is given by $(-)^\mbbd$.
\end{defi}

\begin{prop}\label{THMP_biadjunkce_ck_ck_D}
There is a ``free-forgetful” biadjunction given as follows:
\[\begin{tikzcd}
	{(\Psi, \Phi):(p,q):} & \ck && {\ck_D}
	\arrow[""{name=0, anchor=center, inner sep=0}, "{J_D}"', curve={height=24pt}, from=1-2, to=1-4]
	\arrow[""{name=1, anchor=center, inner sep=0}, "{U_D}"', curve={height=24pt}, from=1-4, to=1-2]
	\arrow["\dashv"{anchor=center, rotate=90}, draw=none, from=0, to=1]
\end{tikzcd}\]

\begin{itemize}
\item The right biadjoint $U_D$ is the forgetful 2-functor sending an algebra to its underlying object,
\item the left biadjoint is a normal pseudofunctor sending:
\[
(f: A \to B) \mapsto ((y_B f)^\mbbd: DA \to DB),
\]
\item the counit $p: J_D U_D \Rightarrow 1$ evaluated at the object $DA$ is the following pseudo morphism:
\[
p_{DA} := (1_{DA})^\mbbd: D^2A \to DA,
\]
With its pseudonaturality square at a pseudo morphism $h$ being the canonical isomorphism between $1^\mbbd_{DB} Dh$ and $h 1_{DA}^\mbbd$, as both are the left Kan extensions of $h$ along $y_{DA}$.
\item the unit is given by the unit of the left Kan pseudomonad $y: 1 \Rightarrow U_D J_D$, with the pseudonaturality square at a morphism $h: A \to B$ being given by the canonical isomorphism:
\[\begin{tikzcd}[cramped]
	A && DA \\
	B && DB
	\arrow["Dh", from=1-3, to=2-3]
	\arrow[""{name=0, anchor=center, inner sep=0}, "{y_A}", from=1-1, to=1-3]
	\arrow[""{name=1, anchor=center, inner sep=0}, "{y_B}"', from=2-1, to=2-3]
	\arrow["h"', from=1-1, to=2-1]
	\arrow["{\mbbd_{y_B h}}"', shorten <=4pt, shorten >=4pt, Rightarrow, from=1, to=0]
\end{tikzcd}\]

\item the components of the modifications are given by the canonical isomorphisms:
\begin{align*}
\Psi &: pJ_D \circ J_Dy \cong 1_{J_D},\\
\Phi &: 1_{U_D} \cong U_D p \circ y U_D.
\end{align*}
\end{itemize}
\end{prop}

\noindent This biadjunction is also lax-idempotent:
\begin{prop}\label{THMP_prop_kz_adjunkce}
There exist (non-invertible) modifications $\Gamma, \Theta$ that serve as the unit and the counit of the following adjunctions:
\begin{align*}
(\Phi^{-1}, \Gamma) &: U_D p \dashv y U_D, \\
(\Theta, \Psi^{-1}) &: J_Dy \dashv pJ_D.
\end{align*}
\end{prop}

\begin{pr}\label{PR_fibration_2monad_Kleisli2cat}
Let us describe the Kleisli 2-category for the fibration 2-monad on the 2-category $\ck / C$ from Example \ref{PR_fibration_2monad}. The Kleisli 2-category can be presented as having the objects functors with codomain $C$, while a morphism $f \rightsquigarrow g$ is a 1-cell $\theta: A \to Pg$ making the triangle below left commute: 
\[\begin{tikzcd}
	A && Pg && B && A && B \\
	&& C && C && C && C
	\arrow["\theta", from=1-1, to=1-3]
	\arrow["f"', from=1-1, to=2-3]
	\arrow[""{name=0, anchor=center, inner sep=0}, "{\rho_g}", from=1-3, to=1-5]
	\arrow[""{name=0p, anchor=center, inner sep=0}, phantom, from=1-3, to=1-5, start anchor=center, end anchor=center]
	\arrow["{\pi_g}"', from=1-3, to=2-3]
	\arrow["g", from=1-5, to=2-5]
	\arrow[""{name=1, anchor=center, inner sep=0}, "u", from=1-7, to=1-9]
	\arrow[""{name=1p, anchor=center, inner sep=0}, phantom, from=1-7, to=1-9, start anchor=center, end anchor=center]
	\arrow["f"', from=1-7, to=2-7]
	\arrow["g", from=1-9, to=2-9]
	\arrow[""{name=2, anchor=center, inner sep=0}, Rightarrow, no head, from=2-3, to=2-5]
	\arrow[""{name=2p, anchor=center, inner sep=0}, phantom, from=2-3, to=2-5, start anchor=center, end anchor=center]
	\arrow[""{name=3, anchor=center, inner sep=0}, Rightarrow, no head, from=2-7, to=2-9]
	\arrow[""{name=3p, anchor=center, inner sep=0}, phantom, from=2-7, to=2-9, start anchor=center, end anchor=center]
	\arrow["\chi", shorten <=4pt, shorten >=4pt, Rightarrow, from=0p, to=2p]
	\arrow["\alpha", shorten <=4pt, shorten >=4pt, Rightarrow, from=1p, to=3p]
\end{tikzcd}\]
From the definition of the comma object, this corresponds to pairs $(u,\alpha)$ consisting of a 1-cell $u: A \to B$ and a 2-cell $\alpha: gu \Rightarrow f$ as portrayed above right.

\noindent In other words, the Kleisli 2-category for this 2-monad is isomorphic to the \textit{colax slice 2-category} $\ck // C$\footnote{This 2-category can also be presented as the 2-category of strict coalgebras and lax morphisms for the 2-comonad $(-)\times C$, see \cite[Chapter 5]{laxcomma2cats}.}. Under this identification, we have a 2-adjunction:
\[\begin{tikzcd}
	{\ck/C} && {\ck//C}
	\arrow[""{name=0, anchor=center, inner sep=0}, "J"', curve={height=24pt}, from=1-1, to=1-3]
	\arrow[""{name=1, anchor=center, inner sep=0}, "U"', curve={height=24pt}, from=1-3, to=1-1]
	\arrow["\dashv"{anchor=center, rotate=90}, draw=none, from=0, to=1]
\end{tikzcd}\]

The left 2-adjoint is the canonical inclusion, the right 2-adjoint sends an object -- a 1-cell $f: A \to C$ -- to the comma object projection $\pi_f: Pf \to C$. The counit $p: JU \Rightarrow 1_{\ck //C}$ evaluated at an object $f: A \to C$ is the colax commutative triangle $(\rho_f, \chi): \pi_f \to f$ from \eqref{EQ_comma_object_lonely}.
\end{pr}

\begin{pozn}\label{POZN_Kleisli_bicat}
We may also define the \textit{Kleisli bicategory} associated to a left Kan pseudomonad $(D,y)$, where objects are the objects in $\ck$ and a morphism $A \rightsquigarrow B$ in $\ck_D$ corresponds to a morphism $A \to DB$ in $\ck$. The unit is given by the unit of the pseudomonad, while the composition is defined using the extension operation:
\[\begin{tikzcd}
	A & B & C & \mapsto & A && B
	\arrow["f", squiggly, from=1-1, to=1-2]
	\arrow["g", squiggly, from=1-2, to=1-3]
	\arrow["{g^\mbbd \circ f}", squiggly, from=1-5, to=1-7]
\end{tikzcd}\]
For the full definition of the Kleisli bicategory associated to a pseudomonad, see \cite[Theorem 4.1]{relativepsmons}. We will denote this bicategory by $\text{Kl}(D)$. It is routine to verify that there is a pseudofunctor $N: \text{Kl}(D) \to \ck_D$ sending the Kleisli morphism $f: A \rightsquigarrow B$ to $f^\mbbd : DA \to DB$ and that it is a biequivalence of bicategories. In this thesis we will for the most part use the 2-category presentation since it is easier to work with.
\end{pozn}

\begin{pr}\label{PR_presheafpsmonad}
Given a locally small category $\ca$, denote by $P \ca$ the full subcategory of $[\ca^{op},\set]$ spanned by \textit{small} presheaves, that is, presheaves that are small colimits of representables. The assignment $\ca \mapsto P \ca$ defines a left Kan pseudomonad on the (large) 2-category $\CAT$ of locally small categories, with the unit $y_\ca: \ca \to P\ca$ being given by the Yoneda embedding and the extension operation being given by ordinary left Kan extension along $y_\ca$. These are guaranteed to exist because of the cocompleteness of $P \cb$; and since $y_\ca$ is fully faithful, the accompanying 2-cell is invertible:
\[\begin{tikzcd}
	\ca && P\ca \\
	\\
	&& P\cb
	\arrow["{y_A}", from=1-1, to=1-3]
	\arrow[""{name=0, anchor=center, inner sep=0}, "{\text{Lan}_{y_A}F}", from=1-3, to=3-3]
	\arrow[""{name=1, anchor=center, inner sep=0}, "F"', from=1-1, to=3-3]
	\arrow[shorten <=7pt, shorten >=7pt, Rightarrow, from=1, to=0]
\end{tikzcd}\]

Pseudo $P$-algebras are precisely the cocomplete categories, and pseudo $P$-morphisms are precisely the cocontinuous functors. The Kleisli 2-category $\CAT_P$ thus has presheaf categories as objects and cocontinuous functors as morphisms. In fact, it can be seen to be biequivalent to the bicategory $\prof$ whose objects are locally small categories and whose morphisms $\ca \rightsquigarrow \cb$ are \textit{small profunctors} $H: \cb^{op} \times \ca \to \set$. Here we call a profunctor $H: \cb^{op} \times \ca \to \set$ small if for every $a \in \ca$, the presheaf $H(-,a): \cb^{op} \to \set$ is small (belongs to $P \cb$).

Under this identification, the left biadjoint from the Kleisli biadjunction in Proposition \ref{THMP_biadjunkce_ck_ck_D} $P: \CAT \to \prof$ sends functor $f: \ca \to \cb$ to the profunctor:
\[
\cb(-,f-): \cb^{op} \times \ca \to \set.
\]
We remark that alternatively there is also a 2-monad presentation for this pseudomonad that uses inaccessible cardinals, see \cite[Chapter 7]{onthemonadicity}.
\end{pr}

\begin{prc}\label{PR_other_cocompletion_psmonads}
Similarly, fixing a class $\Phi$ of small categories, for any locally small category $\ca$ there is a full subcategory $P_\Phi \ca \subseteq P\ca$ of the category of small presheaves spanned by (possibly iterated) closure of representables under $\Phi$-colimits. This is the free $\Phi$-cocompletion of $\ca$ (see \cite[Theorem 2.2]{onthemonadicity}) and the assignment $\ca \mapsto P_{\Phi}\ca$ is a Left Kan pseudomonad whose algebras are those categories that admit $\Phi$-shaped colimits.

Also, various cocompletions of categories in which the colimits are required to interact with finite limits (such as in regular, Barr-exact, lextensive categories) can also be captured using left Kan pseudomonads -- this time not on $\CAT$, but on $\text{LEX}$ -- the 2-category of finitely complete categories and left exact functors. See \cite{lex_colimits}.
\end{prc}

\subsection{Left Kan 2-monads}\label{SUBS_leftKan2monads}

As we mentioned, the majority of our examples in the thesis will be 2-monads. In this short section we develop a ``Left Kan” analogue of ``lax-idempotent 2-monads”.

\begin{defi}\label{DEFI_leftkan2monad}
A left Kan pseudomonad $(D,y)$ is a \textit{left Kan 2-monad} if:
\begin{itemize}
\item $\mbbd_f$ is the identity 2-cell for every 1-cell $f:B \to DA$, meaning that: \[ f^\mbbd \circ y_B = f, \]
\item $g^\mbbd f^\mbbd = (g^\mbbd f)^\mbbd$,
\item $y_A^\mbbd = 1_{DA}$.
\end{itemize}
\end{defi}

\noindent Notice that in case of left Kan 2-monads, the biadjunction from Proposition \ref{THMP_biadjunkce_ck_ck_D} becomes a 2-adjunction. Let us also note the following:

\begin{prop}
The correspondence from Theorem \ref{THM_korespondence_leftkanpsmon_laxidemppsmon} restricts to the correspondence between left Kan 2-monads $(D, y)$ and lax-idempotent 2-monads $(D,m,i)$.
\end{prop}
\begin{proof}
``$\Rightarrow$”: Let $(D, y)$ be a left Kan 2-monad. As we outlined in the proof of Theorem \ref{THM_korespondence_leftkanpsmon_laxidemppsmon}, the pseudofunctor $D$ is defined as this left Kan extension:
\[\begin{tikzcd}
	A && DA \\
	B && DB
	\arrow["{y_A}", from=1-1, to=1-3]
	\arrow["{(y_Bf)^\mbbd =:Df}", from=1-3, to=2-3]
	\arrow["f"', from=1-1, to=2-1]
	\arrow["{y_B}"', from=2-1, to=2-3]
\end{tikzcd}\]
If $(f:A \to B,g: B \to C)$ is a composable pair of morphisms, we have:
\[
D(g f) = (y_C g f)^\mbbd = ((y_C g)^\mbbd y_B f)^\mbbd = (y_C g)^\mbbd (y_B f)^\mbbd = Dg Df
\]
Also, $D1_A = y_A^\mbbd = 1_{DA}$ so $D$ is a 2-functor. This also makes $y$ a 2-natural transformation since the pseudonaturality square is the identity. Next, the pseudonaturality square for the multiplication $m: D^2 \Rightarrow D$ is also the identity since both of the triangles below commute:
\[\begin{tikzcd}
	{D^2A} && DA \\
	{D^2B} && DB
	\arrow["{D^2f}"', from=1-1, to=2-1]
	\arrow["{1_{DB}^\mbbd}"', from=2-1, to=2-3]
	\arrow["{1_{DA}^\mbbd}", from=1-1, to=1-3]
	\arrow["Df", from=1-3, to=2-3]
	\arrow["{(Df)^\mbbd}"{description}, from=1-1, to=2-3]
\end{tikzcd}\]

``$\Leftarrow$”: If $(D,m,i)$ is a lax-idempotent 2-monad, the corresponding left Kan extension in the proof of Theorem \ref{THM_korespondence_leftkanpsmon_laxidemppsmon} has the 2-cell component equal to the identity. The other identities in Definition \ref{DEFI_leftkan2monad} are shown by a straightforward manipulation using the 2-monad identities.
\end{proof}

\section{Internal category theory}\label{SEC_int_cat_theory}

Early references to internal category theory include \cite{ehresmann_categories_structures_HIST}, \cite{grothendieck_HIST}.

\medskip
\noindent \textbf{The section is organized as follows}:
\begin{itemize}
\item In Section \ref{SUBS_defi_pr} we recall the 2-category $\cat(\ce)$ of categories internal to a category $\ce$ with pullbacks.
\item In Section \ref{SUBS_cats_and_pshs} we describe the relationship between internal categories and diagrams in $\ce$ and recall results on limits in the 2-category $\cat(\ce)$.
\item Sections \ref{SUBS_2monads_of_form_CatT} and \ref{SUBS_examples_CatT} focus on 2-functors $\cat(\ce) \to \cat(\ce')$ that are induced by pullback-preserving functors $\ce \to \ce'$.
\end{itemize}

\subsection{Definition and examples}\label{SUBS_defi_pr}

\begin{defi}\label{DEFI_internal_category}
Let $\cc$ be a category with pullbacks. An \textit{internal category} in $\ce$ is a diagram in $\ce$ that is subject to axioms we will state below:
\[\begin{tikzcd}
	{X_2} && {X_1} && {X_0}
	\arrow["{\pi_1}"{description}, shift left=4, from=1-1, to=1-3]
	\arrow["\comp"{description}, from=1-1, to=1-3]
	\arrow["{\pi_2}"{description}, shift right=4, from=1-1, to=1-3]
	\arrow["t"{description}, shift right=4, from=1-3, to=1-5]
	\arrow["s"{description}, shift left=4, from=1-3, to=1-5]
	\arrow["\ids"{description}, from=1-5, to=1-3]
\end{tikzcd}\]
We call $s,t$ the \textit{domain}/\textit{codomain} (or also \textit{source}/\textit{target}), $\ids$ the \textit{identities} map and $\comp$ the \textit{composition} map. The axioms include the following equalities:
\begin{align*}
s \circ \ids &= 1_{X_0},\\
t \circ \ids &= 1_{X_0},\\
t \circ \comp &= t \circ \pi_2,\\
s \circ \comp &= s \circ \pi_1.
\end{align*}
Next, the diagram below left must commute and be a pullback in $\ce$:
\begin{equation}\label{EQ_DEFI_internalcat_pb}
\begin{tikzcd}
	{X_2} && {X_1} && {X_3} && {X_2} \\
	\\
	{X_1} && {X_0} && {X_2} && {X_1}
	\arrow["{\pi_2}", from=1-1, to=1-3]
	\arrow["{\pi_1}"', from=1-1, to=3-1]
	\arrow["\lrcorner"{anchor=center, pos=0.125}, draw=none, from=1-1, to=3-3]
	\arrow["s", from=1-3, to=3-3]
	\arrow["{\rho_2}", from=1-5, to=1-7]
	\arrow["{\rho_1}"', from=1-5, to=3-5]
	\arrow["\lrcorner"{anchor=center, pos=0.125}, draw=none, from=1-5, to=3-7]
	\arrow["{\pi_1}", from=1-7, to=3-7]
	\arrow["t"', from=3-1, to=3-3]
	\arrow["{\pi_2}"', from=3-5, to=3-7]
\end{tikzcd}
\end{equation}
It must satisfy the \textit{unit axioms} (here by round brackets we denote the canonical map into the pullback square):
\begin{align*}
\comp \circ (\ids \circ s,1_{X_1}) &= 1_{X_1},\\
\comp \circ (1_{X_1}, \ids \circ t) &= 1_{X_1}.
\end{align*}
Also, denoting by $X_3$ the pullback above right in \eqref{EQ_DEFI_internalcat_pb}, the \textit{associativity axiom} must also be satisfied:
\begin{equation}
\comp \circ \comp_1 =  \comp \circ \comp_2,
\end{equation}
where the $\comp_i: X_3 \to X_2$ are the unique maps into pullbacks:
\[\begin{tikzcd}[column sep=scriptsize]
	{X_3} && {X_2} &&& {X_3} && {X_2} \\
	& {X_2} && {X_1} &&& {X_2} && {X_1} \\
	{X_2} &&&&& {X_2} \\
	& {X_1} && {X_0} &&& {X_1} && {X_0}
	\arrow["{\rho_2}", from=1-1, to=1-3]
	\arrow["{\comp_1}", dashed, from=1-1, to=2-2]
	\arrow["{\rho_1}"', from=1-1, to=3-1]
	\arrow["{\pi_2}", from=1-3, to=2-4]
	\arrow["{\rho_2}", from=1-6, to=1-8]
	\arrow["{\comp_2}", dashed, from=1-6, to=2-7]
	\arrow["{\rho_1}"', from=1-6, to=3-6]
	\arrow["\comp", from=1-8, to=2-9]
	\arrow["{\pi_2}", from=2-2, to=2-4]
	\arrow["{\pi_1}"', from=2-2, to=4-2]
	\arrow["\lrcorner"{anchor=center, pos=0.125}, draw=none, from=2-2, to=4-4]
	\arrow["s", from=2-4, to=4-4]
	\arrow["{\pi_2}", from=2-7, to=2-9]
	\arrow["{\pi_1}"', from=2-7, to=4-7]
	\arrow["\lrcorner"{anchor=center, pos=0.125}, draw=none, from=2-7, to=4-9]
	\arrow["s", from=2-9, to=4-9]
	\arrow["\comp"', from=3-1, to=4-2]
	\arrow["{\pi_1}"', from=3-6, to=4-7]
	\arrow["t"', from=4-2, to=4-4]
	\arrow["t"', from=4-7, to=4-9]
\end{tikzcd}\]
\end{defi}

\begin{pozn}
The above definition is independent on the choice of the pullback on the right side of \eqref{EQ_DEFI_internalcat_pb}.
\end{pozn}

\begin{defi}\label{DEFI_internal_fun}
Let $A,B$ be internal categories in $\ce$. An \textit{internal functor} $F: A \to B$ consists of a pair of 1-cells in $\ce$:
\begin{align*}
F_0 &: A_0 \to B_0,\\
F_1 &: A_1 \to B_1,
\end{align*}
such that:
\begin{align}
s \circ F_1 &= F_0 \circ s,\\
t \circ F_1 &= F_0 \circ t,\label{eqintfun01}\\
F_1 \circ \ids &= \ids \circ F_0,\label{eqintfun_ids} \\
F_1 \circ \comp &= \comp \circ F_2,\label{eqintfun_comp}
\end{align}
here $F_2$ is the unique morphism into the pullback:
\[\begin{tikzcd}
	{A_2} && {A_1} \\
	&& {B_2} && {B_1} \\
	{A_1} \\
	&& {B_1} && {B_0}
	\arrow["{\pi_2}", from=1-1, to=1-3]
	\arrow["{F_2}"', dashed, from=1-1, to=2-3]
	\arrow["{\pi_1}"', from=1-1, to=3-1]
	\arrow["{F_1}", from=1-3, to=2-5]
	\arrow["{\pi_2}", from=2-3, to=2-5]
	\arrow["{\pi_1}"', from=2-3, to=4-3]
	\arrow["\lrcorner"{anchor=center, pos=0.125}, draw=none, from=2-3, to=4-5]
	\arrow["s", from=2-5, to=4-5]
	\arrow["{F_1}"', from=3-1, to=4-3]
	\arrow["t"', from=4-3, to=4-5]
\end{tikzcd}\]
\end{defi}

\begin{defi}
Let $F,G: A \to B$ be internal functors. An \textit{internal natural transformation} $\alpha: F \Rightarrow G$ is a morphism $\alpha : A_0 \to B_1$ such that:
\begin{align}
s \circ \alpha &= F_0,\\
t \circ \alpha &= G_0,\label{eqintnattr1}\\
\comp \circ (\alpha \circ t, F_1) &= \comp \circ (G_1, \alpha \circ s)\label{eqintnattr_naturalita}
\end{align}
\end{defi}

\begin{nota}
Given a category $\ce$ with pullbacks, internal categories, functors, natural transformations in $\ce$ form a 2-category that we will denote by $\cat(\ce)$.
\end{nota}

\begin{pr}
For $\ce = \set$ we obtain small categories, functors and natural transformations. For $\ce = \cat$ we obtain (strict) \textit{double categories} -- they will be studied in Section \ref{SEC_double_cats} and later in Chapter \ref{CH_FS_DBLCATS}.
\end{pr}

\begin{defi}\label{DEFI_ob_D}
Let $\ce$ be a category with pullbacks. An internal category $A$ is \textit{discrete} on an object $X \in \ob \ce$ if $A_i = X$ for $i \in \{ 0,1 \}$ and all the 1-cells in the definition are the identities on $X$. There is a discrete category functor $D: \ce \to \cat(\ce)$ defined by sending an object to the discrete category on it.

We may also define the \textit{underlying object functor} $\ob: \cat(\ce) \to \ce$, which is the assignment:
\[
\cc \mapsto \cc_0.
\]
If $\ce$ admits reflexive coequalizers, there is an additional 2-functor $\pi_0: \cat(\ce) \to \ce$ sending a category $\cc$ to the coequalizer of its source and target maps:
\[\begin{tikzcd}
	{\cc_1} && {\cc_0} && {\pi_0 \cc}
	\arrow["s", shift left=2, from=1-1, to=1-3]
	\arrow["t"', shift right=2, from=1-1, to=1-3]
	\arrow[dashed, from=1-3, to=1-5]
\end{tikzcd}\]
\end{defi}

\begin{prop}\label{THMP_Ddashvobdashvpi0}
Let $\ce$ be a category with pullbacks. Then $D \dashv \ob$ and this adjunction is a coreflection. If in addition $\ce$ admits reflexive coequalizers, the functor $\pi_0$ is left adjoint to $D$, forming an adjoint triple:
\[\begin{tikzcd}
	\ce && {\cat(\ce)}
	\arrow[""{name=0, anchor=center, inner sep=0}, "D"{description}, from=1-1, to=1-3]
	\arrow[""{name=0p, anchor=center, inner sep=0}, phantom, from=1-1, to=1-3, start anchor=center, end anchor=center]
	\arrow[""{name=0p, anchor=center, inner sep=0}, phantom, from=1-1, to=1-3, start anchor=center, end anchor=center]
	\arrow[""{name=1, anchor=center, inner sep=0}, "\ob"', curve={height=24pt}, from=1-3, to=1-1]
	\arrow[""{name=1p, anchor=center, inner sep=0}, phantom, from=1-3, to=1-1, start anchor=center, end anchor=center, curve={height=24pt}]
	\arrow[""{name=2, anchor=center, inner sep=0}, "{\pi_0}", curve={height=-24pt}, from=1-3, to=1-1]
	\arrow[""{name=2p, anchor=center, inner sep=0}, phantom, from=1-3, to=1-1, start anchor=center, end anchor=center, curve={height=-24pt}]
	\arrow["\top"{description}, draw=none, from=0p, to=1p]
	\arrow["\top"{description}, draw=none, from=2p, to=0p]
\end{tikzcd}\]
\end{prop}
\begin{proof}
We clearly have $Z = \ob (DZ)$ for all $Z \in \ce$ -- this will be the unit of the adjunction. Given $Z \in \ce$ and $\cc \in \cat(\ce)$, it is also clear that any 1-cell $f: Z \to \cc_0$ in $\ce$ is of form $\ob F$ for a unique internal functor $F: DZ \to \cc$, proving the universal property of the unit. Next, assume that $\ce$ admits reflexive coequalizers. We will prove that there is a chain of bijections:
\[
\ce(\pi_0 \cc, A) \cong \text{Coeq}((s,t),A) \cong \cat(\ce)(\cc, DA).
\]
The middle set denotes the set of morphisms $f: \cc_0 \to A$ coequalizing the domain and codomain maps $s,t: \cc_1 \to \cc_0$. The first bijection is the definition of the coequalizer. Next, note that any functor $F: \cc \to DA$ is uniquely determined by its objects component $F_0: \cc \to A$. This is because we have the equalities:
\[\begin{tikzcd}
	{\cc_1} && {\cc_0} \\
	\\
	{\cc_0} && A
	\arrow["s", from=1-1, to=1-3]
	\arrow["t"', from=1-1, to=3-1]
	\arrow["{F_1}"{description}, from=1-1, to=3-3]
	\arrow["{F_0}", from=1-3, to=3-3]
	\arrow["{F_0}"', from=3-1, to=3-3]
\end{tikzcd}\]
This also proves that the object component of $F$ equalizes $s,t$.
\end{proof}

\subsection{Categories and presheaves}\label{SUBS_cats_and_pshs}

\begin{nota}\label{NOTA_simplnota}
Denote by $\Delta$ the (skeletal) category of nonempty finite ordinals and order-preserving maps between them. Its objects will be denoted as follows:
\[
[ n ] := \{ 0 \to 1 \to \dots \to n \}.
\]
We will use the usual notation for the generating 1-cells of this category: we denote by $\delta_i: [n] \to [n+1]$ the injective order-preserving maps uniquely determined by the property $\delta^{-1}_i(i) = \emptyset$, and denote by $\sigma_j: [m] \to [m-1]$ those surjective order-preserving maps determined by the property $| \sigma^{-1}_j (j) | = 2$.

Denote by $\Delta_2$ the subcategory of $\Delta$ spanned by the diagram:
\[\begin{tikzcd}
	{[2]} && {[1]} && {[0]}
	\arrow["{\sigma_0}"{description}, from=1-3, to=1-5]
	\arrow["{\delta _0}"{description}, shift left=4, from=1-5, to=1-3]
	\arrow["{\delta_1}"{description}, shift right=4, from=1-5, to=1-3]
	\arrow["{\delta_2}"{description}, shift right=4, from=1-3, to=1-1]
	\arrow["{\delta_1}"{description}, from=1-3, to=1-1]
	\arrow["{\delta_0}"{description}, shift left=4, from=1-3, to=1-1]
\end{tikzcd}\]
\end{nota}

\begin{prop}\label{THMP_nerve}
Let $\ce$ be a category with pullbacks. There is a fully faithful embedding of 1-categories:
\[
N_\ce : \cat(\ce) \to [\Delta^{op}_2,\ce].
\]
\end{prop}
\begin{proof}
The functor $N_\ce$ sends an internal functor $F: A \to B$ to the following:
\[\begin{tikzcd}
	{A_2} && {A_1} && {A_0} \\
	\\
	{B_2} && {B_1} && {B_0}
	\arrow["\comp"{description}, from=1-1, to=1-3]
	\arrow["{\pi_2}"{description}, shift right=5, from=1-1, to=1-3]
	\arrow["{\pi_1}"{description}, shift left=5, from=1-1, to=1-3]
	\arrow["{F_1\times_{B_0}F_1}"{description}, from=1-1, to=3-1]
	\arrow["t"{description}, shift right=5, from=1-3, to=1-5]
	\arrow["s"{description}, shift left=5, from=1-3, to=1-5]
	\arrow["{F_1}"{description}, from=1-3, to=3-3]
	\arrow["\ids"{description}, from=1-5, to=1-3]
	\arrow["{F_0}"{description}, from=1-5, to=3-5]
	\arrow["\comp"{description}, from=3-1, to=3-3]
	\arrow["{\pi_1}"{description}, shift left=4, from=3-1, to=3-3]
	\arrow["{\pi_2}"{description}, shift right=4, from=3-1, to=3-3]
	\arrow["s"{description}, shift left=4, from=3-3, to=3-5]
	\arrow["t"{description}, shift right=4, from=3-3, to=3-5]
	\arrow["\ids"{description}, from=3-5, to=3-3]
\end{tikzcd}\]
By construction, $N_\ce$ is a faithful embedding. Let us now consider a natural transformation $G: N_\ce A \Rightarrow N_\ce B$ -- the naturality forces that $G_2 \sstackrel{!}{=} G_1 \times_{B_0} G_1$ and so $G$ induces a unique functor $A \to B$ with components $G_0, G_1$, proving fullness.
\end{proof}

\begin{defi}\label{DEFI_nerve}
Given an internal category $\cc$ in $\ce$, the diagram below will be called its \textit{nerve}: \[ N_\ce(\cc): \Delta^{op}_2 \to \ce. \]
\end{defi}

\begin{pozn}
A diagram $X: \Delta_2^{op} \to \ce$ lies in the image of $N_\ce$ if an only if an appropriate diagram is a pullback in $\ce$:
\[\begin{tikzcd}
	{X[2]} && {X[1]} \\
	\\
	{X[1]} && {X[0]}
	\arrow["{X\delta_2}", from=1-1, to=1-3]
	\arrow["{X\delta_0}"', from=1-1, to=3-1]
	\arrow["\lrcorner"{anchor=center, pos=0.125}, draw=none, from=1-1, to=3-3]
	\arrow["{X\delta_0}", from=1-3, to=3-3]
	\arrow["{X\delta_1}"', from=3-1, to=3-3]
\end{tikzcd}\]
and also, the associativity and unit laws of an internal category hold for $X$ if we define the pullback $X_3$ as in Definition \ref{DEFI_internal_category}.
\end{pozn}

\begin{lemma}[Preservation and reflection of internal categories]\label{THML_preservation_reflection_int_cats}
Let $\ce, \ce'$ be categories with pullbacks and $T: \ce \to \ce'$ a pullback-preserving functor. The post-composition with $T$ induces a functor between the categories of internal categories that we denote by $\cat(T)$:
\begin{equation}\label{EQ_CatE_CatEprime_Tstar}
\begin{tikzcd}
	{\cat(\ce)} && {\cat(\ce')} \\
	\\
	{[\Delta^{op}_2,\ce]} && {[\Delta^{op}_2,\ce']}
	\arrow["\exists", dashed, from=1-1, to=1-3]
	\arrow["{N_\ce}"', hook, from=1-1, to=3-1]
	\arrow["{N_{\ce'}}", hook, from=1-3, to=3-3]
	\arrow["{T_*}"', from=3-1, to=3-3]
\end{tikzcd}
\end{equation}

\noindent Assume in addition that $T: \ce \to \ce'$ is fully faithful. Then:
\begin{itemize}
\item $\cat(T)$ is a fully faithful functor,
\item a diagram $X:\Delta^{op}_2 \to \ce$ corresponds to an internal category in $\ce$ if and only if the diagram below corresponds to an internal category in $\ce'$:
\[\begin{tikzcd}
	{\Delta^{op}_2} && \ce && {\ce'}
	\arrow["X", from=1-1, to=1-3]
	\arrow["T", from=1-3, to=1-5]
\end{tikzcd}\]
\end{itemize}
\end{lemma}
\begin{proof}
The first part of the statement follows from the fact that $T$ preserves pullbacks.

Assume now that $T$ is fully faithful. This implies that $T_*$ is fully faithful. In particular so is $\cat(T)$ because of the commutativity of \eqref{EQ_CatE_CatEprime_Tstar}.

The last part follows from the fact that $T$ preserves and reflects pullbacks and reflects equalities between 1-cells.

\end{proof}

\begin{pozn}\label{POZN_CatT_2fun_2ff}
More is true than what we have seen in Lemma \ref{THML_preservation_reflection_int_cats}. Given a pullback-preserving functor $T: \ce \to \ce'$, $\cat(T): \cat(\ce) \to \cat(\ce')$ is a 2-functor. It sends an internal natural transformation:
\[
\alpha: G \Rightarrow H: \ca \to \cb,
\]
to the transformation whose corresponding 1-cell in $\ce'$ is given by $T\alpha : T\ca_0 \to T\cb_1$. Moreover, if $T$ is fully faithful, $\cat(T)$ is 2-fully faithful.
\end{pozn}

\begin{lemma}\label{THML_simplicialobject_iff_cat}
Let $\ce$ be a category with pullbacks. The following are equivalent:
\begin{itemize}
\item $X: \Delta_2^{op} \to \ce$ is the nerve of an internal category $\cc$ in $\ce$,
\item $\ce(Z,-) \circ X: \Delta_2^{op} \to \set$ is the nerve of a small category for all $Z \in \ce$,

\end{itemize}
\end{lemma}
\begin{proof}
The diagram $X$ is a category if and only if the composite $y \circ X$ with the Yoneda embedding $y: \ce \to [\ce^{op}, \set]$ is so (this is because of Lemma \ref{THML_preservation_reflection_int_cats} and the fact that $y$ is fully faithful). Since limits in functor categories are pointwise, this latter condition is equivalent to the requirement that the composite $\text{ev}_Z \circ y \circ X$ with the evaluation functor at $Z$ is a small category for each $Z \in \ob \ce$ -- but these diagrams are precisely the diagrams in the second bullet point of this statement.

\end{proof}

In the lemma below and there only, let us denote by $\ce(Z,\cc)$ the small category corresponding to the nerve $\ce(Z,-) \circ N_\ce \cc$.

\begin{lemma}\label{THML_graph_cell_is_category_cell}
Let $\ce$ be a category with pullbacks. We have the following:
\begin{itemize}
\item If $\cc, \cd \in \cat(\ce)$ are internal categories and $F_i: \cc_i \to \cd_i$ 1-cells in $\ce$ for $i \in \{ 0,1 \}$, the following are equivalent:
\begin{itemize}
\item the 1-cells $F_i$ give an internal functor $\cc \to \cd$
\item the functions $\ce(Z,F_i)$ give a functor $\ce(Z,\cc) \to \ce(Z,\cd)$ for all $Z \in \ob \ce$.
\end{itemize}
If this is the case, let us denote the functor by $\ce(Z,F)$ in this lemma.

\item If $F,G: \cc \to \cd$ are internal functors between internal categories and $\alpha: \cc_0 \to \cd_1$ is a 1-cell in $\ce$. The following are equivalent:
\begin{itemize}
\item $\alpha: F \Rightarrow G$ is an internal natural transformation,
\item the function $\ce(Z,\alpha)$ is a natural transformation $\ce(Z,F) \Rightarrow \ce(Z,G)$.
\end{itemize}
\end{itemize}
\end{lemma}
\begin{proof}
The directions ``$\Rightarrow$” for both statements follow from the fact that each $\ce(Z,-)$ preserves pullbacks. The directions ``$\Leftarrow$” follow from the fact that the collection $\ce(Z,-), Z \in \ob \ce$ is jointly faithful and jointly reflects pullbacks.
\end{proof}

\noindent The following has first been observed in \cite[Page 108]{streetfib}:

\begin{theo}\label{THM_Dtilde_2ff}
Let $\ce$ be a category with pullbacks and $D: \ce \to \cat(\ce)$ the discrete category functor (Definition \ref{DEFI_ob_D}). The 2-functor:
\begin{align*}
\widetilde{D} : \cat(\ce) &\to [\ce^{op},\cat],\\
\cc &\mapsto \cat(\ce)(D-,\cc),
\end{align*}
is 2-fully faithful.
\end{theo}
\begin{proof}
Denoting by $y: \ce \to [\ce^{op},\set]$ the Yoneda embedding, notice that we have the following isomorphisms in $\twoCAT$:
\[\begin{tikzcd}
	{\cat(\ce)} && {[\ce^{op},\cat]} \\
	\\
	&& {\cat([\ce^{op},\set])}
	\arrow["{\widetilde{D}}", from=1-1, to=1-3]
	\arrow["{\cat(y)}"', from=1-1, to=3-3]
	\arrow["\cong", from=1-3, to=3-3]
\end{tikzcd}\]
The Yoneda embedding is fully faithful and so $\cat(y)$ is 2-fully faithful by Remark \ref{POZN_CatT_2fun_2ff}. This implies that $\widetilde{D}$ is so as well.
\end{proof}

\begin{cor}\label{THMC_colimit_in_CatE_iff_in_Cat}
Let $\ce$ be a category with pullbacks and denote $\ck := \cat(\ce)$. Consider the 2-functors $W: \cj^{op} \to \cat$, $H: \cj \to \ck$ and let $\lambda: W \Rightarrow \ck(H-,\cc)$ be a cocone. Then the second condition implies the first one:
\begin{itemize}
\item $\lambda$ exhibits the internal category $\cc$ as the $W$-weighted colimit of $H$,
\item for every $Z \in \ob \ce$, the cocone below exhibits the small category $\ck(DZ,\cc)$ as the $W$-weighted colimit of the 2-functor $\ck(DZ, H-):\cj \to \cat$:
\[\begin{tikzcd}
	W && {\ck(H-,\cc)} && {\cat(\ck(DZ,H-),\ck(DZ,\cc))}
	\arrow["\lambda", from=1-1, to=1-3]
	\arrow["{\ck(DZ,-)_{H-, \cc}}", from=1-3, to=1-5]
\end{tikzcd}\]
\end{itemize}
\noindent The same statement holds for limits as well.
\end{cor}
\begin{proof}
The 2-functor $\widetilde{D}$ is 2-fully faithful and thus reflects 2-colimits. This means that $\cc$ is the 2-colimit if $\widetilde{D}\cc = \cat(\ce)(D-,\cc)$ is the 2-colimit in $[\ce^{op},\cat]$. But since colimits in functor categories are pointwise, this is equivalent to the second point in the statement.
\end{proof}

\noindent Finally, let us study 2-limits in the 2-category $\cat(\ce)$. We begin with ordinary limits in its underlying category:

\begin{lemma}\label{THML_CatE_has_1limits}
Let $\ce$ be a category with pullbacks that moreover admits $\cj$-shaped limits (for a small category $\cj$). The category $\cat(\ce)$ then also admits $\cj$-shaped limits and they are ``pointwise”: the following composite creates them:
\[\begin{tikzcd}
	{\cat(\ce)} && {[\Delta^{op}, \ce]} && {[\ob \Delta^{op}, \ce]}
	\arrow["{N_\ce}", from=1-1, to=1-3]
	\arrow["{\text{res}}", from=1-3, to=1-5]
\end{tikzcd}\]

\end{lemma}
\begin{proof}
Let $K: \cj \to \cat(\ce)$ be a functor. Since $\cat(\ce) \xhookrightarrow{} [\Delta^{op}_2, \ce]$ is fully faithful, the result will follow if we show that the limit object $\lim K$ lies in $\cat(\ce)$. But this fact follows immediately from the fact that pullbacks commute with limits.

\end{proof}

\begin{theo}\label{THM_catE_limits}
The 2-category $\cat(\ce)$ for a category $\ce$ with pullbacks admits:
\begin{itemize}
\item powers with $\boldsymbol{2}$,
\item pullbacks,
\item comma objects,
\item lax and oplax limits of arrows.
\end{itemize}
\end{theo}
\begin{proof}

\textbf{Powers with} $\boldsymbol{2}$: The reader may find the full proof in \cite[Proposition 3.19]{johnphd}. What we will offer here is a sketch of a proof demonstrating the methods that will be used in Section \ref{SEC_thecaseckCatE}.

Assume first that $\ce = \set$ -- it is well known that the power of a small category $\cd$ with $\boldsymbol{2}$ is given by the functor category $[\boldsymbol{2}, \cd]$.

Let now $\ce$ be an arbitrary category with pullbacks and let $\cc$ be an internal category in $\ce$. Construct an internal category $\cp$ in $\ce$ such that $\cp_0 = \cc_1$ and $\cp_1$ is the object of ``commutative squares” in $\cc$, in other words, a pullback:
\[\begin{tikzcd}
	{\cp_1} && {\cc_2} \\
	&&&& {\text{ in } \ce} \\
	{\cc_2} && {\cc_1}
	\arrow[from=1-1, to=1-3]
	\arrow[from=1-1, to=3-1]
	\arrow["\lrcorner"{anchor=center, pos=0.125}, draw=none, from=1-1, to=3-3]
	\arrow["\comp", from=1-3, to=3-3]
	\arrow["\comp"', from=3-1, to=3-3]
\end{tikzcd}\]
The domains, codomains, units, compositions can all be defined using pullbacks in $\ce$. Regarded as a diagram $\cp^{op}: \Delta_2 \to \ce$, note that the composition with each functor $\ce(Z,-): \ce \to \set$ yields the nerve of the following category (in $\set$):
\[
[\boldsymbol{2}, \cat(\ce)(DZ,\cc)]
\]
and since this is a category, $\cp$ is an internal category in $\ce$ by Lemma \ref{THML_simplicialobject_iff_cat}. Likewise, one can constructs a pair of functors and a natural transformation in $\cat(\ce)$:
\[\begin{tikzcd}
	\cc && \cp
	\arrow[""{name=0, anchor=center, inner sep=0}, "{P_2}"', curve={height=18pt}, from=1-1, to=1-3]
	\arrow[""{name=1, anchor=center, inner sep=0}, "{P_1}", curve={height=-18pt}, from=1-1, to=1-3]
	\arrow["\lambda", shorten <=5pt, shorten >=5pt, Rightarrow, from=1, to=0]
\end{tikzcd}\]
such that each $\cat(\ce)(DZ,-): \cat(\ce) \to \cat$ sends it to the canonical limit cone with apex $[\boldsymbol{2}, \cat(\ce)(DZ,\cc)]$. And since it is a limit cone for each $Z \in \ce$, it is the limit in $\cat(\ce)$ by Corollary \ref{THMC_colimit_in_CatE_iff_in_Cat}.

\textbf{Pullbacks}: We have seen in Lemma \ref{THML_CatE_has_1limits} that the 1-category $\cat(\ce)$ has pullbacks. It is straightforward to see that any pullback is a 2-pullback -- this again follows from them being pointwise.

\textbf{Comma objects}: The existence of powers with $\boldsymbol{2}$ and pullbacks in $\cat(\ce)$ imply that this 2-category is \textit{representable} \cite{grayreport_HIST}. Any such 2-category admits comma objects by \cite[Proposition 1]{streetfib}.

\textbf{Lax and oplax limits of arrows}: They exist since they are a special case of comma objects, as pointed out in Example \ref{PR_lax_limit_of_arrow}.
\end{proof}

\begin{pozn}
Note that alternatively, one could characterize internal categories to $\ce$ as functors\footnote{$\Delta_3$ would be an appropriate subcategory of $\Delta$ on objects $[0], \dots, [3]$} $\Delta^{op}_3 \to \ce$ or functors $\Delta^{op} \to \ce$ (called \textit{simplicial objects} in $\ce$) that preserve certain pullbacks (analogous to \cite[Proposition 3.3.5]{higheroperads}). We chose our approach since we believe it is better suited to our needs in Chapter \ref{CH_COMPUTATION}, in particular Section \ref{SEC_thecaseckCatE}. 
\end{pozn}

\subsection{2-monads of form $\cat(T)$}\label{SUBS_2monads_of_form_CatT}

\begin{defi}
Recall the following:
\begin{itemize}
\item a natural transformation $\alpha: F \Rightarrow G: \ce \to \cd$ between 1-functors is said to be \textit{cartesian} if all its naturality squares are pullbacks,
\item a monad $(T,\mu,\eta)$ is said to be \textit{cartesian} if $T$ preserves pullbacks and $\mu, \eta$ are cartesian natural transformations.
\end{itemize}
\end{defi}

\begin{prop}\label{THMP_F_gives_CatF_alfa_gives_Catalfa}
Let $\ce, \cd$ be categories with pullbacks.
\begin{itemize}
\item Any pullback-preserving functor $F: \ce \to \cd$ induces a (2-)pullback-preserving 2-functor $\cat(F): \cat(\ce) \to \cat(\cd)$,
\item Any natural transformation $\alpha: F \Rightarrow G: \ce \to \cd$ induces a 2-natural transformation $\cat(\alpha): \cat(F) \Rightarrow \cat(G)$. If $\alpha$ is cartesian, so is $\cat(\alpha)$.
\end{itemize}
\end{prop}
\begin{proof}
We know that there is an induced 2-functor $\cat(F)$ from Lemma \ref{THML_preservation_reflection_int_cats} and Remark \ref{POZN_CatT_2fun_2ff}. The preservation of pullbacks by the 1-functor $\cat(F)$ follows from the fact that they are pointwise in $\cat(\ce)$ and that $F$ preserves them. The fact that it preserves them as a 2-functor follows from the fact that any pullback in $\cat(\ce)$ is a 2-pullback.

The 2-naturality of $\cat(\alpha)$ directly follows the naturality of $\alpha$. The fact that $\cat(\alpha)$ is cartesian follows from the fact that pullbacks in $\cat(\ce)$ are pointwise.
\end{proof}

\begin{prc}
For any category $\ce$ with pullbacks, the representable functors $\ce(Z,-): \ce \to \set$ preserve pullbacks (since they preserve all limits). We may thus consider the induced 2-functors $\cat(\ce(Z,-)): \cat(\ce) \to \cat$. Observe:
\end{prc}

\begin{prop}\label{THMP_catEDZ_is_catEZ}
The following 2-functors are isomorphic:
\[
\cat(\ce)(DZ, -) \cong \cat(\ce(Z,?))(-): \cat(\ce) \to \cat.
\]
\end{prop}
\begin{proof}
It can be seen that in both cases, an internal category $\cc$ in $\ce$ is sent to a category described as follows: The objects are 1-cells $x: Z \to \cc_0$, a morphism $f: x \to y$ is a 1-cell $f: Z \to \cc_1$ with the property that:
\begin{align*}
s \circ f &= x,\\
t \circ f &= y.
\end{align*}
The identity morphism on $x: Z \to \cc_0$ is the 1-cell:
\[\begin{tikzcd}
	Z && {\cc_0} && {\cc_1} & {\text{ in } \ce}
	\arrow["x", from=1-1, to=1-3]
	\arrow["\ids", from=1-3, to=1-5]
\end{tikzcd}\]
The composite of $f: x \to y$, $g: y \to z$ is given by the 1-cell:
\[\begin{tikzcd}
	Z && {\cc_2} && {\cc_1}
	\arrow["{(f,g)}", from=1-1, to=1-3]
	\arrow["\comp", from=1-3, to=1-5]
\end{tikzcd}\]
\end{proof}

\begin{nota}
Let $\ce$ be a category with pullbacks. Proposition \ref{THMP_F_gives_CatF_alfa_gives_Catalfa} in particular implies that any monad $(T,m,i)$ such that $T$ preserves pullbacks gives rise to a 2-monad $(\cat(T), \cat(m), \cat(i))$ on $\cat(\ce)$. We will typically denote this 2-monad also by $T$ and say \textit{$T$ is of form $\cat(T)$} to emphasize that there is an underlying 1-monad. In particular we will use the notation ``$\talgs$” for the 2-category of strict algebras and strict morphisms for the 2-monad $\cat(T)$. This will not cause confusion as the 1-category of algebras for the 1-monad $T$ on $\ce$ will be denoted by $\ce^T$.
\end{nota}

The assignment $\ce \mapsto \cat(\ce)$ for a category $\ce$ with pullbacks preserves EM-objects in the following sense:

\begin{prop}
Let $(T,m,i)$ be a monad on a category $\ce$ with pullbacks and assume that $T$ preserves pullbacks. The following 2-categories are isomorphic:
\[
\talgs \cong \cat(\ce^T).
\]
\end{prop}
\begin{proof}
See \cite[Remark 7.1]{johnphd}
\end{proof}

\begin{defi}\label{DEFI_discropfib}
An internal functor $F: A \to B$ is said to be a \textit{discrete opfibration} if the following is a pullback in $\ce$:
\[\begin{tikzcd}
	{A_1} && {B_1} \\
	\\
	{A_0} && {B_0}
	\arrow["{F_1}", from=1-1, to=1-3]
	\arrow["s"', from=1-1, to=3-1]
	\arrow["\lrcorner"{anchor=center, pos=0.125}, draw=none, from=1-1, to=3-3]
	\arrow["s", from=1-3, to=3-3]
	\arrow["{F_0}"', from=3-1, to=3-3]
\end{tikzcd}\]
If the diagram involving targets is a pullback, $F$ is said to be a \textit{discrete fibration}.
\end{defi}

\begin{pr}[$\ce = \set$]
A functor $F: \ca \to \cb$ between two small categories is a discrete opfibration if for every pair $(a' \in \ob \ca, f: Fa' \to b \in \mor \cb)$ there exists a unique morphism $f' : a' \to b'$ such that $Ff' = f$.
\end{pr}

\begin{prop}\label{THMP_catF_preserves_discropfibs}
Let $\ce, \ce'$ be categories with pullbacks.
\begin{itemize}
\item For any pullback-preserving functor $H: \ce \to \ce'$, the 2-functor $\cat(H)$ preserves discrete fibrations and opfibrations.
\item For any cartesian natural transformation $\alpha: G \Rightarrow H: \ce \to \ce'$ between pullback-preserving functors, the 2-natural transformation $\cat(\alpha)$ is a discrete fibration and opfibration in each component.
\end{itemize}
\end{prop}
\begin{proof}
To prove the \textbf{first claim}, assume that an internal functor $F: A \to B$ is a discrete opfibration. $\cat(H)(F)$ is a discrete opfibration since the relevant square is the one below, and it is a pullback since $H$ preserves them: 
\[\begin{tikzcd}
	{HA_1} && {HB_1} \\
	\\
	{HA_0} && {HB_0}
	\arrow["{HF_1}", from=1-1, to=1-3]
	\arrow["Hs"', from=1-1, to=3-1]
	\arrow["\lrcorner"{anchor=center, pos=0.125}, draw=none, from=1-1, to=3-3]
	\arrow["Hs", from=1-3, to=3-3]
	\arrow["{HF_0}"', from=3-1, to=3-3]
\end{tikzcd}\]

\noindent The \textbf{second claim} follows from the fact that the following naturality square for $\alpha$ is a pullback in $\ce$:
\[\begin{tikzcd}
	{FA_1} && {GA_1} \\
	\\
	{FA_0} && {GA_0}
	\arrow["{\alpha_{A_1}}", from=1-1, to=1-3]
	\arrow["Fs"', from=1-1, to=3-1]
	\arrow["\lrcorner"{anchor=center, pos=0.125}, draw=none, from=1-1, to=3-3]
	\arrow["Gs", from=1-3, to=3-3]
	\arrow["{\alpha_{A_0}}"', from=3-1, to=3-3]
\end{tikzcd}\]
\end{proof}

\noindent The above lemma pinpoints the key properties of 2-monads of form $\cat(T)$ that will be relevant to us later:
\begin{cor}\label{THMC_catT2monad_multiplication_is_discrete_opfib}
Let $(T,m,i)$ be a 2-monad on a 2-category $\ce$ that is of form $\cat(T')$ for a cartesian monad on $\ce$. Then:
\begin{itemize}
\item $T$ preserves discrete (op)fibrations,
\item the multiplication $m: T^2 \Rightarrow T$ is a componentwise discrete (op)fibration.
\end{itemize}
\end{cor}

\subsection{Examples}\label{SUBS_examples_CatT}

\begin{pr}\label{PR_reslan_catT}

The identity 2-monad is of this form. More generally, for a 1-category $\cj$, the reslan 2-monad $T$ on $[ \ob \cj, \cat]$ from Example \ref{PR_reslan2monad} can easily be seen to be of form $\cat(T')$, with the monad $T'$ on $[\ob \cj, \set]$ being given by the 1-categorical left Kan extension followed by restriction.
\end{pr}

\begin{pr}\label{PR_freestrictmoncat2monad_CaTT}
If $T: \cat \to \cat$ is the free monoid monad, $\cat(T)$ is the free strict monoidal category 2-monad from Example \ref{PR_freestrictmoncat2monad}.
\end{pr}

\begin{pr}\label{PR_doublecat_2monad}
There is a cartesian monad $\text{fc}: \gph \to \gph$ on the category $\gph$ of directed multigraphs sending a graph to the free category on it (see \cite[5.1]{higheroperads}). The category $\gph$ is complete so we may consider the 2-category of categories internal to it. An object of the 2-category $\cat(\gph)$ can be identified with a ``$\cat$-graph” -- an example of which we have seen in Example \ref{PR_coinserter_in_Cat}. The 2-monad $\cat(\fc)$ sends a $\cat$-graph $X$ to the one with the same objects and vertical morphisms, and whose horizontal morphisms and squares are ``horizontal paths” of squares in $X$:
\[\begin{tikzcd}
	{a_0} & {a_1} & \dots & {a_{n-1}} & {a_n} \\
	{b_0} & {b_1} & \dots & {b_{n-1}} & {b_n}
	\arrow[""{name=0, anchor=center, inner sep=0}, "{g_1}", from=1-1, to=1-2]
	\arrow["{u_0}"', from=1-1, to=2-1]
	\arrow[from=1-2, to=1-3]
	\arrow[from=1-2, to=2-2]
	\arrow[from=1-3, to=1-4]
	\arrow[""{name=1, anchor=center, inner sep=0}, "{g_n}", from=1-4, to=1-5]
	\arrow[from=1-4, to=2-4]
	\arrow["{u_n}", from=1-5, to=2-5]
	\arrow[""{name=2, anchor=center, inner sep=0}, "{h_1}"', from=2-1, to=2-2]
	\arrow[from=2-2, to=2-3]
	\arrow[from=2-3, to=2-4]
	\arrow[""{name=3, anchor=center, inner sep=0}, "{h_n}"', from=2-4, to=2-5]
	\arrow["{\alpha_1}"', shorten <=4pt, shorten >=4pt, Rightarrow, from=0, to=2]
	\arrow["{\alpha_n}"', shorten <=4pt, shorten >=4pt, Rightarrow, from=1, to=3]
\end{tikzcd}\]

\noindent One can verify that:
\begin{itemize}
\item a strict $T$-algebra is a \text{double category} -- we will encounter them in Section \ref{SEC_double_cats},
\item a strict morphism is a double functor,
\item a lax morphism is a \textit{lax double functor}. A lax double functor $F: X \to Y$ between double categories is an assignment on objects, vertical morphisms, horizontal morphisms, squares, that strictly preserves vertical compositions, and is equipped with the associator and unitor squares in $Y$ that are subject to equations:
\[\begin{tikzcd}
	Fa && Fb && Fc & Fa && Fa \\
	\\
	Fa &&&& Fc & Fa && Fa
	\arrow["Fg", from=1-1, to=1-3]
	\arrow[Rightarrow, no head, from=1-1, to=3-1]
	\arrow["Fh", from=1-3, to=1-5]
	\arrow[Rightarrow, no head, from=1-5, to=3-5]
	\arrow[""{name=0, anchor=center, inner sep=0}, "{1_{Fa}}", from=1-6, to=1-8]
	\arrow[""{name=0p, anchor=center, inner sep=0}, phantom, from=1-6, to=1-8, start anchor=center, end anchor=center]
	\arrow[Rightarrow, no head, from=1-6, to=3-6]
	\arrow[Rightarrow, no head, from=1-8, to=3-8]
	\arrow[""{name=1, anchor=center, inner sep=0}, "{F(hg)}"', from=3-1, to=3-5]
	\arrow[""{name=1p, anchor=center, inner sep=0}, phantom, from=3-1, to=3-5, start anchor=center, end anchor=center]
	\arrow[""{name=2, anchor=center, inner sep=0}, "{F1_a}"', from=3-6, to=3-8]
	\arrow[""{name=2p, anchor=center, inner sep=0}, phantom, from=3-6, to=3-8, start anchor=center, end anchor=center]
	\arrow["{\gamma_{g,h}}", shorten >=7pt, Rightarrow, from=1-3, to=1p]
	\arrow["{\iota_a}", shorten <=9pt, shorten >=9pt, Rightarrow, from=0p, to=2p]
\end{tikzcd}\]
\end{itemize}
\end{pr}

\begin{pr}\label{PR_2category_2monad_CatT}

Fix a set $X$. Denote by $\fc_X$ the restriction of the monad $\fc$ on $\gph$ to the category $\gph_X$ spanned by graphs whose set of object is $X$, and whose morphisms fix the set of objects -- equivalently this is the category $\set^{X \times X}$. $T$ sends an $(X \times X)$-indexed collection of sets $(\cc(A,B))_{(A,B) \in X \times X}$ to the collection:
\[
(T\cc)_{A,B} = \text{Path}_\cc(A,B) = \{ (f_1,\dots, f_m) \mid m \in \mathbb{N}, \text{cod}(f_{i})= \text{dom}(f_{i+1}) \forall i < m \}
\]
$\gph_X$ has pullbacks and $\fc_X$ is cartesian. $\cat(\fc_X)$ is the free 2-category 2-monad from Example \ref{PR_2category_2monad}.
\end{pr}

\begin{nonpr}\label{NOPR_multipl_not_discr_opfib_freesymstrictmoncat}
The free symmetric strict monoidal category 2-monad $S$ on $\cat$ (Example \ref{PR_freesymstrictmoncat2monad}) is cartesian, but not of form $\cat(T)$ -- this is because 2-monads of this form preserve discrete categories (as is easily verified). $S$, on the other hand, does not, because for instance the terminal category $*$ is sent to the permutation grupoid $S*$, which is not discrete. 

This 2-monad also does \textbf{not} satisfy that the multiplication $m: S^2 \Rightarrow S$ is a discrete (op)fibration -- or any (op)fibration for that matter. Consider the multiplication functor $m_*: S^2* \to S*$ for the terminal category. The following pair of an object in $S^2*$ and a morphism in $S*$ can be seen to admit no lift to a morphism in $S^2*$:
\[\begin{tikzcd}
	& {(\bullet,} & {\bullet,} & {\bullet,} & {\bullet)} \\
	{((\bullet),(\bullet,\bullet),(\bullet))} \\
	& {(\bullet,} & {\bullet,} & {\bullet,} & {\bullet)}
	\arrow[from=1-2, to=3-3]
	\arrow[from=1-3, to=3-2]
	\arrow[from=1-4, to=3-5]
	\arrow[from=1-5, to=3-4]
\end{tikzcd}\]
\end{nonpr}

\section{Double categories}\label{SEC_double_cats}

Putting $\ce := \cat$ in the previous section, we obtain the notions of a double category and a double functor. We will write down these definitions explicitly since they will be needed in Chapter \ref{CH_FS_DBLCATS}. We also recall some definitions from \cite{limitsindbl}.

\begin{defi}
A \textit{double category} $X$ is an internal category in $\cat$. In particular it consists of the following diagram in $\cat$:
\[\begin{tikzcd}
	{X_2} && {X_1} && {X_0}
	\arrow["s"{description}, from=1-5, to=1-3]
	\arrow["{d_1}"{description}, shift left=4, from=1-3, to=1-5]
	\arrow["{d_0}"{description}, shift right=4, from=1-3, to=1-5]
	\arrow["{d_1}"{description}, from=1-1, to=1-3]
	\arrow["{d_0}"{description}, shift right=4, from=1-1, to=1-3]
	\arrow["{d_2}"{description}, shift left=4, from=1-1, to=1-3]
\end{tikzcd}\]
We call objects of $X_0$ \textit{the objects} of $X$, morphisms of $X_0$ the \textit{vertical morphisms}, objects of $X_1$ the \textit{horizontal morphisms}, the morphisms of $X_1$ \textit{squares}. The fact that a morphism $\alpha: g \to h$ in $X_1$ has $d_1(\alpha) = u$, $d_0(\alpha) = v$ is portrayed as follows:
\[\begin{tikzcd}
	a && b \\
	\\
	c && d
	\arrow["u"', from=1-1, to=3-1]
	\arrow[""{name=0, anchor=center, inner sep=0}, "g", from=1-1, to=1-3]
	\arrow["v", from=1-3, to=3-3]
	\arrow[""{name=1, anchor=center, inner sep=0}, "h"', from=3-1, to=3-3]
	\arrow["\alpha", shorten <=9pt, shorten >=9pt, Rightarrow, from=0, to=1]
\end{tikzcd}\]
We refer to the composition in categories $X_0, X_1$ as the \textit{vertical composition}. The \textit{horizontal} composition of horizontal morphisms and of squares is given by the functor $d_1: X_2 \to X_1$. An identity morphism in $X_1$ will be called a \textit{vertical identity square}, and a \textit{horizontal identity square} on a vertical morphism $u$ will be the one given by $s(u)$. By an \textit{identity square} in $X$ we mean a square that is either a vertical or a horizontal identity. Objects and horizontal morphisms form a category that we will denote by $h(X)$.
\end{defi}

\begin{pozn}
The category $h(X)$ associated to a double category $X$ can be equivalently defined as $\cat(\ob)(X)$, where $\ob: \cat \to \set$ is the objects functor (Definition \ref{DEFI_ob_D}).
\end{pozn}

\begin{defi}
A \textit{double functor} $F: X \to Y$ between double categories is given by a pair of functors agreeing on objects:
\begin{align*}
F_0 &: X_0 \to Y_0,\\
h(F) &: h(X) \to h(Y),
\end{align*}
together with an assignment on squares:
\[\begin{tikzcd}
	a & b & {} & Fa & Fb & {} \\
	c & d & {} & Fc & Fd & {}
	\arrow[""{name=0, anchor=center, inner sep=0}, "g", from=1-1, to=1-2]
	\arrow["u"', from=1-1, to=2-1]
	\arrow["v", from=1-2, to=2-2]
	\arrow["\mapsto"{description}, draw=none, from=1-3, to=2-3]
	\arrow[""{name=1, anchor=center, inner sep=0}, "Fg", from=1-4, to=1-5]
	\arrow["Fu"', from=1-4, to=2-4]
	\arrow["Fv", from=1-5, to=2-5]
	\arrow["{\in Y}"{description}, draw=none, from=1-6, to=2-6]
	\arrow[""{name=2, anchor=center, inner sep=0}, "h"', from=2-1, to=2-2]
	\arrow[""{name=3, anchor=center, inner sep=0}, "Fh"', from=2-4, to=2-5]
	\arrow["\alpha", shorten <=4pt, shorten >=4pt, Rightarrow, from=0, to=2]
	\arrow["{F\alpha}", shorten <=4pt, shorten >=4pt, Rightarrow, from=1, to=3]
\end{tikzcd}\]
such that all compositions and identities are preserved.
\end{defi}

\begin{nota}
We will use the following notation for the 2-category of double categories and double functors:
\[
\dbl := \cat(\cat).
\]
\end{nota}

\begin{defi}[Duals and transposes]\label{DEFI_dblcat_duals}
Any double category $X$ has its \textit{transpose} $X^T$ obtained by switching vertical and horizontal morphisms and compositions. Any double category $X$ has its \textit{vertical dual} which we denote by $X^{v}$. It is  defined by putting:
\[
(X^{v})_i = X_i^{op} \text{ for } i \in \{ 0,1,2 \}.
\]
Similarly, domain and codomain functors for $X^{v}$ are obtained by applying $(-)^{op}$ on those for $X$.

There is also a notion of a \textit{horizontal dual} $X^{h}$, which is a double category obtained from $X$ by switching $d_0$'s and $d_1$'s.
\end{defi}

\begin{defi}
A double category $X$ is \textit{flat} if any square $\alpha$ is uniquely determined by its boundary.
\end{defi}

\begin{defi}
A double category $X$ is (strictly) \textit{horizontally invariant} if for any two invertible horizontal morphisms $g,h$ and every vertical morphism $u$ there exists a unique square filling the picture\footnote{This definition differs from the one in \cite{limitsindbl} in that we require the filler square to be unique.}:
\[\begin{tikzcd}
	a && b \\
	\\
	c && d
	\arrow[""{name=0, anchor=center, inner sep=0}, "{g \cong}", from=1-1, to=1-3]
	\arrow[""{name=0p, anchor=center, inner sep=0}, phantom, from=1-1, to=1-3, start anchor=center, end anchor=center]
	\arrow[dashed, from=1-1, to=3-1]
	\arrow["u", from=1-3, to=3-3]
	\arrow[""{name=1, anchor=center, inner sep=0}, "{h\cong}"', from=3-1, to=3-3]
	\arrow[""{name=1p, anchor=center, inner sep=0}, phantom, from=3-1, to=3-3, start anchor=center, end anchor=center]
	\arrow["{\exists ! \lambda}", shorten <=9pt, shorten >=9pt, Rightarrow, from=0p, to=1p]
\end{tikzcd}\]

\noindent Such a square is necessarily horizontally invertible. We say that $X$ is (strictly) \textit{vertically invariant} if its transpose $X^T$ is horizontally invariant. A double category that is both horizontally and vertically invariant will be called (strictly) \textit{invariant}.
\end{defi}

\begin{pr}\label{PR_sq_cc}
Let $\cc$ be a category. There is a double category $\sq{\cc}$ such that:
\begin{itemize}
\item objects are the objects of $\cc$,
\item vertical and horizontal morphisms are the morphisms of $\cc$,
\item squares are commutative squares in $\cc$.
\end{itemize}
\end{pr}

\begin{pr}\label{PR_PbSq_cc}
There is a sub-double category $\text{PbSq}(\cc) \subseteq \sq{\cc}$ with the same objects and morphisms, whose squares are the pullback squares in $\cc$.

Another example is a sub-double category $\text{MPbSq}(\cc) \subseteq \text{PbSq}(\cc)$ with the same objects and horizontal morphisms whose vertical morphisms are monomorphisms in the category $\cc$.
\end{pr}

\begin{pr}\label{PR_BOFIB}
Let $\ce$ be a category with pullbacks. There is a double category $\text{BOFib}(\ce)$ such that:
\begin{itemize}
\item objects are categories internal to $\ce$,
\item vertical morphisms are internal functors that are discrete opfibrations (Definition \ref{DEFI_discropfib}),
\item horizontal morphisms are internal functors $F: A \to B$ that are \textit{bijective on objects}, i.e. the object part 1-cell $F_0 : A_0 \to B_0$ is an isomorphism in $\ce$,
\item a square in $\text{BOFib}(\ce)$ is a commutative square in $\cat(\ce)$ (which is automatically a pullback square).
\end{itemize}
\end{pr}

\noindent Note that all of the above examples are flat and invariant.

\section{$T$-multicategories}\label{SEC_Tmulticats}

The idea of multicategories is to have a categorical structure in which the morphisms are allowed to have multiple inputs and one output. The idea of $T$-multicategories is to replace the ``list of inputs” -- which would be an element of the free monoid monad on a set $X$ -- by an element of $TX$, where $T$ is \textbf{any} cartesian monad. $T$-multicategories have been introduced in \cite{burroni_HIST}, for a textbook reference we refer to \cite[Chapter 4]{higheroperads}.

\subsection{Definition and properties}

\begin{ass}
Throughout this section, fix a category $\ce$ with pullbacks and a cartesian monad $(T,\mu, \eta)$ on $\ce$.
\end{ass}

\begin{defi}\label{DEFI_T-graph}
A \textit{$T$-graph} $\cm$ is a span in $\ce$ of this form:
\[\begin{tikzcd}
	{T\cm_0} && {\cm_1} && {\cm_0}
	\arrow["s"', from=1-3, to=1-1]
	\arrow["t", from=1-3, to=1-5]
\end{tikzcd}\]
\begin{itemize}
\item $\cm_0$ will be referred to as the \textit{object of objects} of $\cm$,
\item $\cm_1$ is the \textit{object of morphisms of $\cm$},
\item $s,t$ are the \textit{source, target maps}.
\end{itemize}
\end{defi}

Given a $T$-graph $\cm$, denote by $\cm_2$, $\cm_3$ the following pullbacks:
\[\begin{tikzcd}
	{\cm_2} && {\cm_1} && {\cm_3} && {\cm_1} \\
	\\
	{T\cm_1} && {T\cm_0} && {T\cm_2} & {T\cm_1} & {T\cm_0}
	\arrow["{\pi_2}", from=1-1, to=1-3]
	\arrow["{\pi_1}"', from=1-1, to=3-1]
	\arrow["\lrcorner"{anchor=center, pos=0.125}, draw=none, from=1-1, to=3-3]
	\arrow["s", from=1-3, to=3-3]
	\arrow["{\rho_2}", from=1-5, to=1-7]
	\arrow["{\rho_1}"', from=1-5, to=3-5]
	\arrow["\lrcorner"{anchor=center, pos=0.125}, draw=none, from=1-5, to=3-7]
	\arrow["s", from=1-7, to=3-7]
	\arrow["Tt"', from=3-1, to=3-3]
	\arrow["{T\pi_2}"', from=3-5, to=3-6]
	\arrow["Tt"', from=3-6, to=3-7]
\end{tikzcd}\]

\begin{defi}
A \textit{$T$-multicategory} $\cm$ is a $T$-graph equipped with the \textit{identities} and \textit{composition} 1-cells $\ids: \cm_0 \to \cm_1$, $\comp: \cm_2 \to \cm_1$ that make the following commute:
\[\begin{tikzcd}
	&&&& {\cm_2} & {\cm_1} & {\cm_0} \\
	{\cm_0} && {\cm_0} && {T\cm_1} \\
	&&&& {T^2\cm_0} \\
	{T\cm_0} && {\cm_1} && {T\cm_0} && {\cm_1}
	\arrow["{\pi_2}", from=1-5, to=1-6]
	\arrow["{\pi_1}"', from=1-5, to=2-5]
	\arrow["\comp"{description}, from=1-5, to=4-7]
	\arrow["t", from=1-6, to=1-7]
	\arrow[Rightarrow, no head, from=2-1, to=2-3]
	\arrow["{\eta_{\cm_0}}"', from=2-1, to=4-1]
	\arrow["\ids"{description}, from=2-1, to=4-3]
	\arrow["Ts"', from=2-5, to=3-5]
	\arrow["{\mu_{\cm_0}}"', from=3-5, to=4-5]
	\arrow["t"', from=4-3, to=2-3]
	\arrow["s", from=4-3, to=4-1]
	\arrow["t"', from=4-7, to=1-7]
	\arrow["s", from=4-7, to=4-5]
\end{tikzcd}\]
This data is subject to the \textit{associativity} and \textit{unit laws} -- they can be unpacked for instance from \cite[Definition 4.2.2]{higheroperads}.
\end{defi}

\begin{defi}
A \textit{$T$-graph morphism} $F: \cm \to \cn$ between $T$-multicategories $\cm, \cn$ is given by two 1-cells:
\begin{align*}
F_0 &: \cm_0 \to \cn_0,\\
F_1 &: \cm_1 \to \cn_1,
\end{align*}
satisfying the equation:
\[\begin{tikzcd}
	&& {\cm_1} \\
	{T\cm_0} && {\cn_1} && {\cm_0} \\
	{T\cn_0} &&&& {\cn_0}
	\arrow["s"', from=1-3, to=2-1]
	\arrow["{F_1}", from=1-3, to=2-3]
	\arrow["t", from=1-3, to=2-5]
	\arrow["{TF_0}"', from=2-1, to=3-1]
	\arrow["s", from=2-3, to=3-1]
	\arrow["t"', from=2-3, to=3-5]
	\arrow["{F_0}", from=2-5, to=3-5]
\end{tikzcd}\]
\end{defi}

\begin{defi}
A \textit{$T$-multifunctor} $F: \cm \to \cn$ between $T$-multicategories is a $T$-graph morphism between the underlying $T$-graphs of $\cm$, $\cn$ that makes the following diagrams commute:
\[\begin{tikzcd}
	{\cm_0} && {\cm_1} && {\cm_2} && {\cm_1} \\
	\\
	{\cn_0} && {\cn_1} && {\cn_2} && {\cn_1}
	\arrow["\ids", from=1-1, to=1-3]
	\arrow["{F_0}"', from=1-1, to=3-1]
	\arrow["{F_1}", from=1-3, to=3-3]
	\arrow["\comp", from=1-5, to=1-7]
	\arrow["{F_2}"', from=1-5, to=3-5]
	\arrow["{F_1}", from=1-7, to=3-7]
	\arrow["\ids"', from=3-1, to=3-3]
	\arrow["\comp"', from=3-5, to=3-7]
\end{tikzcd}\]
here $F_2$ is the unique 1-cell into the pullback square:
\[\begin{tikzcd}
	{\cm_2} && {\cm_1} \\
	& {\cn_2} && {\cn_1} \\
	{T\cm_1} \\
	& {T\cn_1} && {T\cn_0}
	\arrow["{\pi_2}", from=1-1, to=1-3]
	\arrow["{F_2}", dashed, from=1-1, to=2-2]
	\arrow["{\pi_1}"', from=1-1, to=3-1]
	\arrow["{F_1}", from=1-3, to=2-4]
	\arrow["{\pi_2}", from=2-2, to=2-4]
	\arrow["{\pi_1}"', from=2-2, to=4-2]
	\arrow["\lrcorner"{anchor=center, pos=0.125}, draw=none, from=2-2, to=4-4]
	\arrow["s", from=2-4, to=4-4]
	\arrow["{TF_1}"', from=3-1, to=4-2]
	\arrow["Tt"', from=4-2, to=4-4]
\end{tikzcd}\]
\end{defi}

\begin{defi}
A \textit{$T$-graph 2-cell} $\alpha: F \Rightarrow G: \cm \to \cn$ is a 1-cell $\alpha: \cm_0 \to \cn_1$ satisfying:
\[\begin{tikzcd}
	{\cm_0} && {\cn_1} && {T\cn_0} & {\cm_0} & {\cn_1} & {\cn_0} \\
	&& {\cn_0}
	\arrow["\alpha", from=1-1, to=1-3]
	\arrow["{F_0}"', from=1-1, to=2-3]
	\arrow["s", from=1-3, to=1-5]
	\arrow["\alpha", from=1-6, to=1-7]
	\arrow["{G_0}"', curve={height=18pt}, from=1-6, to=1-8]
	\arrow["t", from=1-7, to=1-8]
	\arrow["{\eta_{\cn_0}}"', from=2-3, to=1-5]
\end{tikzcd}\]
\end{defi}
\begin{defi}
A \textit{$T$-multinatural transformation} $\alpha: F \Rightarrow G: \cm \to \cn$ is a $T$-graph 2-cell of the underlying $T$-graph morphisms satisfying the further property that:
\[
\comp \circ (T\alpha s, G_1) = \comp \circ (\eta_{\cm_1} F_1, \alpha t),
\]
where the maps above are the unique maps into pullbacks:
\[\begin{tikzcd}[column sep=scriptsize]
	{\cm_1} &&&& {\cm_1} && {\cm_0} \\
	& {\cn_2} && {\cn_1} && {\cn_2} && {\cn_1} \\
	{T\cm_0} &&&& {\cn_1} \\
	& {T\cn_1} && {T\cn_0} && {T\cn_1} && {T\cn_0}
	\arrow["{(T\alpha s, G_1)}", dashed, from=1-1, to=2-2]
	\arrow["{G_1}", curve={height=-24pt}, from=1-1, to=2-4]
	\arrow["s"', from=1-1, to=3-1]
	\arrow["t", from=1-5, to=1-7]
	\arrow["{(\eta_{\cm_1} F_1, \alpha t)}", dashed, from=1-5, to=2-6]
	\arrow["{F_1}"', from=1-5, to=3-5]
	\arrow["\alpha", from=1-7, to=2-8]
	\arrow["{\pi_2}", from=2-2, to=2-4]
	\arrow["{\pi_1}"', from=2-2, to=4-2]
	\arrow["\lrcorner"{anchor=center, pos=0.125}, draw=none, from=2-2, to=4-4]
	\arrow["s", from=2-4, to=4-4]
	\arrow["{\pi_2}", from=2-6, to=2-8]
	\arrow["{\pi_1}"', from=2-6, to=4-6]
	\arrow["\lrcorner"{anchor=center, pos=0.125}, draw=none, from=2-6, to=4-8]
	\arrow["s", from=2-8, to=4-8]
	\arrow["{T\alpha}"', from=3-1, to=4-2]
	\arrow["{\eta_{\cm_1}}"', from=3-5, to=4-6]
	\arrow["Tt"', from=4-2, to=4-4]
	\arrow["Tt"', from=4-6, to=4-8]
\end{tikzcd}\]
\end{defi}

\begin{nota}
$T$-multicategories, $T$-multifunctors and $T$-multinatural transformations form a 2-category that we will denote by $\tmult$.
\end{nota}

\subsection{Examples}

\begin{pr}
For the identity monad on a category $\ce$ with pullbacks, we will refer to $T$-graphs as just \textit{graphs} in $\ce$. Equivalently, the underlying 1-category of $\tgph$ can be defined as the presheaf category $[J^{op},\ce]$, where $J$ is the \textit{walking parallel pair of arrows}:
\[\begin{tikzcd}
	\bullet && \bullet
	\arrow[shift left=2, from=1-1, to=1-3]
	\arrow[shift right=2, from=1-1, to=1-3]
\end{tikzcd}\]

$T$-multicategories are exactly categories internal to $\ce$.
\end{pr}

\begin{pr}
In case $T$ is the free monoid monad on $\set$, a $T$-multicategory is what is referred to as a \textit{(plain) multicategory}. A plain multicategory $\cm$ consists of a set of objects $\ob \cm$ and a set of morphisms $\mor \cm$, each having a source (consisting of a list of objects) and a target (consisting of a single object). We will denote $f \in \mor \cm$ with source $(a_1, \dots, a_n)$ and target $b$ as:
\[
f: a_1, \dots, a_n \to b.
\]
For the exposition and references on plain multicategories, see \cite[Section 2.1]{higheroperads}.
\end{pr}

\begin{pr}
Given the pointed set monad $T := (-)+1$ on $\set$, a $T$-multicategory is a plain multicategory which only contains unary and nullary morphisms. See \cite[Example 4.2.11]{higheroperads}.
\end{pr}

\begin{pr}\label{PR_fc_multicategories}
There is a cartesian monad $\fc: \gph \to \gph$ that sends a directed multigraph to the free category on it -- we have already mentioned it in Example \ref{PR_doublecat_2monad}. $\fc$-multicategories are commonly called \textit{virtual double categories} -- they have been studied in \cite[Chapter 5]{higheroperads}. See also \cite{unified} for their usage in unifying various ``multicategory”-like structures. A virtual double category consists of a set of objects, vertical, horizontal morphisms, and squares, each of which has a path of horizontal morphisms in its domain:
\[\begin{tikzcd}
	{a_1} & {a_2} & \dots & {a_{n-1}} & {a_n} \\
	\\
	{b_1} && {} && {b_2}
	\arrow["{g_1}", from=1-1, to=1-2]
	\arrow["u"', from=1-1, to=3-1]
	\arrow["{g_2}", from=1-2, to=1-3]
	\arrow["{g_{n-2}}", from=1-3, to=1-4]
	\arrow["\alpha", shorten >=6pt, Rightarrow, from=1-3, to=3-3]
	\arrow["{g_{n-1}}", from=1-4, to=1-5]
	\arrow["v", from=1-5, to=3-5]
	\arrow["h"', from=3-1, to=3-5]
\end{tikzcd}\]
There is no composition operation for horizontal morphisms, only for vertical morphisms and squares.
\end{pr}

\begin{pr}\label{PR_opetopes}
Fix a set $X$ and recall the cartesian monad $\fc_X$ on $\set^{X \times X}$ from Example \ref{PR_2category_2monad_CatT}. A $T$-multicategory is what could be called a \textit{many-to-one 2-computad}. It consists of a set of objects, set of morphisms, and set of 2-cells, each having a path of morphisms as its source and a single morphism as its target, for instance:
\[\begin{tikzcd}
	& {a_2} & \dots & {a_{n-1}} \\
	\\
	{b_1} && {} && {b_2}
	\arrow["{g_2}", from=1-2, to=1-3]
	\arrow["{g_{n-2}}", from=1-3, to=1-4]
	\arrow["\alpha", shorten >=6pt, Rightarrow, from=1-3, to=3-3]
	\arrow["{g_{n-1}}", from=1-4, to=3-5]
	\arrow["{g_1}", from=3-1, to=1-2]
	\arrow["h"', from=3-1, to=3-5]
\end{tikzcd}\]
Morphisms can not be composed, but the 2-cells can, in an appropriate sense. Each opetope is in particular a virtual double category (of Example \ref{PR_fc_multicategories}).
\end{pr}

\begin{pr}\label{PR_reslan_multicats}
Fix a category $\cj$ and consider the monad $T$ on the category $[ \ob \cj, \set]$ given by the left Kan extension along $\iota : \ob \cj \to \cj$ followed by the restriction (this is the $\set$-enriched variant of Example \ref{PR_reslan2monad}). The following is well-known in the literature (see for instance \cite{unified} or \cite[Example 4.2.13]{higheroperads}):
\begin{prop}
For the monad $T$ on $[\ob \cj, \set]$ defined above, we have the following equivalence:
\[ \tmult \simeq \cat / \cj. \]
\end{prop}
\end{pr}

\section{Weak limits}\label{SEC_weak_limits}

\begin{defi}\label{DEFI_rali_colimit}
Let $\cx$ be a class of functors between small categories. Let $\ck$ be a 2-category and let $F: \cj \to \ck$, $W: \cj \to \cat$ be two 2-functors. A ($W$-weighted) cone $\lambda: W \Rightarrow \ck(L,F-)$ is said to be an \textit{$\cx$-limit of $F$ weighted by $W$} if for every $A \in \ck$, the whiskering functor
\begin{align}\label{EQ_kappaA}
\begin{split}
\kappa_A : \ck(A,L) &\to [\cj, \cat](W, \ck(A, F?)),\\
\kappa_A : (\theta: A \to L) &\mapsto (\ck(\theta,F?) \circ \lambda),
\end{split}
\end{align}
belongs to the class $\cx$. By an \textit{$\cx$-colimit of $F: \cj \to \ck$ weighted by $W:\cj^{op} \to \cat$} we will mean an $\cx$-limit of the 2-functor $F^{op}: \cj^{op} \to \ck^{op}$ weighted by $W$.
\end{defi}

\begin{vari}
We will encounter a version of this definition where $F,W$ are pseudofunctors and $\lambda$ is pseudonatural (and $\psd{\cj}{\cat}$ is used instead of $[\cj, \cat]$). As not to burden ourselves with further terminology, we will also call this an \textit{$\cx$-limit of $F$ weighted by $W$}.
\end{vari}

\begin{pr}
If $\cx = \{ \text{isomorphisms} \}$, an $\cx$-limit is the same thing as a 2-limit. In case $\cx = \{ \text{equivalences} \}$, we obtain the notion of a bilimit. 
\end{pr}

\begin{pozn}
Conical (left adjoint)-limits of 2-functors have first been introduced \cite[I,7.9.1]{formalcat} under the name \textit{quasi-limits}\footnote{In fact, the definition in \cite{formalcat} is stronger than ours because it requires the existence of a 2-functor picking the limits that is right lax adjoint to the constant embedding 2-functor $\ck \to \ck^\cj$.}.
\end{pozn}

In Chapter \ref{CH_colax_adjs_lax_idemp_psmons} the primary focus will be on $\cx$-(co)limits in which the class $\cx$ consists of adjoint functors. For a better illustration we include this example below:

\begin{prc}
In elementary terms, a conical (right adjoint)-limit of a 2-functor $F: \cj \to \ck$ consists of a cone $\lambda: \Delta L \Rightarrow F$, and, for every $A \in \ck$:
\begin{itemize}
\item a functor $L: \text{Cone}(A,F) \to \ck(A,L)$ (picking the comparison 1-cell into the weak limit object),
\item for every cone $\mu: \Delta A \Rightarrow F$ a modification $\eta: \mu \to \lambda (L\mu)$.
\end{itemize}
This data has the following \textbf{universal property}: given a 1-cell $f: A \to L$ and a modification $\sigma: \mu \to \lambda g$, there exists a unique 2-cell $\overline{\sigma}: L\mu \Rightarrow g$ such that for every $i \in \ob \cj$ the following 2-cells are equal:
\[\begin{tikzcd}
	A && L && Fi & A && L && Fi
	\arrow["g"', from=1-1, to=1-3]
	\arrow[""{name=0, anchor=center, inner sep=0}, "{\mu_i}", curve={height=-30pt}, from=1-1, to=1-5]
	\arrow["{\lambda_i}"', from=1-3, to=1-5]
	\arrow[""{name=1, anchor=center, inner sep=0}, "g"', curve={height=24pt}, from=1-6, to=1-8]
	\arrow[""{name=2, anchor=center, inner sep=0}, "{L\mu}", from=1-6, to=1-8]
	\arrow[""{name=3, anchor=center, inner sep=0}, "{\mu_i}", curve={height=-30pt}, from=1-6, to=1-10]
	\arrow["{\lambda_i}"', from=1-8, to=1-10]
	\arrow["{\sigma_i}"{pos=0.8}, shorten <=3pt, Rightarrow, from=0, to=1-3]
	\arrow["{\exists ! \overline{\sigma}}"', shorten <=3pt, shorten >=3pt, Rightarrow, dashed, from=2, to=1]
	\arrow["{\eta_i}", shorten <=3pt, Rightarrow, from=3, to=1-8]
\end{tikzcd}\]
Perhaps the intuition to have here is that the modification $\eta$ tells us that $\mu$ weakly factorizes through $\lambda$ using the 1-cell $L\mu$. This 1-cell is not unique but it is the ``smallest” 1-cell with this property.
\end{prc}

\begin{pozn}
Because the maps $\kappa_A: \ck(A,L) \to [\cj, \cat](W, \ck(A,F?))$ together form a 2-natural transformation $\kappa: \ck(-,L) \Rightarrow [\cj, \cat](W,\ck(-,F?))$, by the doctrinal adjunction (Theorem \ref{THM_doctrinal_adj} for the 2-monad $T$ on $[ \ob \cj, \cat]$ from Example \ref{PR_reslan2monad}) the above definition is equivalent to requiring that $\kappa$ is a right adjoint in the 2-category $\colax{\ck}{\cat}$.
\end{pozn}

\begin{pozn}[Ordinary weakness]\label{POZN_ordinary_weakness}
In ordinary category theory, a \textit{weak limit} of a 1-functor $F: \cj \to \cc$ is a cone $\lambda: \Delta L \Rightarrow F$ with the property that for every $A \in \ob \cc$, the \textbf{function} $\kappa_A$ (defined analogous to \eqref{EQ_kappaA}) is surjective (see \cite[Definition 3.3]{weakadjfuns_HIST}).

Notice that given a 2-category $\ck$, every conical rali- and lali-limit of a 2-functor $F$ into $\ck$ whose domain is a locally discrete 2-category is automatically a weak limit of the underlying 1-functor in the underlying 1-category of $\ck$. It is not the case that every weak limit of an underlying 1-functor of a 2-functor can be ``enhanced” to be a rali/lali-limit. This is because if the 2-category $\ck$ is locally discrete, the notion of rali-limit coincides with an ordinary limit and not a weak limit.
\end{pozn}

\begin{pozn}[Enriched weakness]
Definition \ref{DEFI_rali_colimit} is a special case of an $(\cv, \ce)$-\textit{enriched weak limit} in the sense of \cite[Section 4]{enrichedweakness}. The enriching category $\cv$ is in our case equal to $\cat$, with the class $\ce$ being the class $\cx$.

In \cite{enrichedweakness}, the authors paid a special attention to the case $\cv = \cat$, with $\cx$ being the class of functors that are surjective equivalences.
\end{pozn}

\begin{pr}
An object $I$ in a 2-category $\ck$ is (right adjoint)-initial if the unique functor into the terminal category is a right adjoint for every object $A \in \ck$:
\[\begin{tikzcd}
	{\ck(I,A)} && {*}
	\arrow[""{name=0, anchor=center, inner sep=0}, "{!}"', curve={height=24pt}, from=1-1, to=1-3]
	\arrow[""{name=1, anchor=center, inner sep=0}, curve={height=24pt}, from=1-3, to=1-1]
	\arrow["\dashv"{anchor=center, rotate=-90}, draw=none, from=1, to=0]
\end{tikzcd}\]
Clearly, this happens if and only if the hom-category $\ck(I,A)$ has an initial object for every $A \in \ck$. A (right adjoint)-initial object is the same thing as a coreflector-initial object, which is the same as a rali-initial object.

For a particular example, consider the 2-category $\text{MonCat}_l$ of monoidal categories and lax monoidal functors. The terminal monoidal category $*$ is rali-initial because for every monoidal category $\mbba$ we have an isomorphism between the category $\text{MonCat}_l(*,\mbba)$ and the category $\text{Mon}(\mbba)$ of monoids in $\mbba$, and this category has an initial object given by the monoidal unit of $\mbba$.
\end{pr}

\begin{pr}
Consider a 2-category $\ck$ with a zero object $0 \in \ck$ and with a further property that the zero morphism $0_{A,B}: A \to B$ is the initial object in $\ck(A,B)$ for every pair of objects $A,B$. The object $0$ is then a conical lali-colimit of any 2-functor $F: \cj \to \ck$. This is because in the definition of a lali-colimit:
\[\begin{tikzcd}
	{\ck(0,A)} && {\text{Cocone}(A,F)}
	\arrow[""{name=0, anchor=center, inner sep=0}, curve={height=24pt}, from=1-1, to=1-3]
	\arrow[""{name=1, anchor=center, inner sep=0}, curve={height=24pt}, from=1-3, to=1-1]
	\arrow["\dashv"{anchor=center, rotate=90}, draw=none, from=0, to=1]
\end{tikzcd}\]
we have $\ck(0,A) \cong *$, and so the question becomes whether the category of cocones of $F$ with apex $A$ has an initial object. In this case, it does: it is given by the cocone whose components are the zero morphisms.
This for instance applies to the poset-enriched categories $\rel$ of sets and relations and $\parf$ of sets and partial functions.
\end{pr}

\begin{pr}\label{PR_colax_slice_2cat_weaklimit_1}
Let $\ck$ be a 2-category with pullbacks and comma objects and consider the slice 2-category $\ck / C$ and the colax slice 2-category $\ck // C$ from Example \ref{PR_fibration_2monad_Kleisli2cat}. It is known that the product of two objects $f_1, f_2$ in $\ck / C$ is the diagonal in the pullback square of $f_1, f_2$ in $\ck$:
\[\begin{tikzcd}
	& L \\
	{A_1} && {A_2} \\
	& C
	\arrow["{f_1}"', from=2-1, to=3-2]
	\arrow["{f_2}", from=2-3, to=3-2]
	\arrow[from=1-2, to=2-1]
	\arrow[from=1-2, to=2-3]
	\arrow["\lrcorner"{anchor=center, pos=0.125, rotate=-45}, draw=none, from=1-2, to=3-2]
\end{tikzcd}\]

One guess would be that this becomes a weak product in $\ck//C$ after applying the inclusion 2-functor $\ck/C \xhookrightarrow{} \ck//C$, but that would be a wrong guess. To calculate the weak product of $f_1, f_2$ in $\ck // C$, we first calculate the product of the comma object projections for $f_1, f_2$ (using the notation from Example \ref{PR_fibration_2monad_Kleisli2cat}) in $\ck/C$ as pictured below left:
\[\begin{tikzcd}
	& L &&& L & {Pf_i} && {A_1} \\
	{Pf_1} && {Pf_2} \\
	& C &&&& C && C
	\arrow["{\tau_1}"', from=1-2, to=2-1]
	\arrow["{\tau_2}", from=1-2, to=2-3]
	\arrow["\lrcorner"{anchor=center, pos=0.125, rotate=-45}, draw=none, from=1-2, to=3-2]
	\arrow["{\tau_i}", from=1-5, to=1-6]
	\arrow["l"', from=1-5, to=3-6]
	\arrow[""{name=0, anchor=center, inner sep=0}, "{\rho_i}", from=1-6, to=1-8]
	\arrow["{\pi_{f_i}}"{description}, from=1-6, to=3-6]
	\arrow["{f_i}", from=1-8, to=3-8]
	\arrow["{\pi_{f_1}}"', from=2-1, to=3-2]
	\arrow["{\pi_{f_2}}", from=2-3, to=3-2]
	\arrow[""{name=1, anchor=center, inner sep=0}, Rightarrow, no head, from=3-6, to=3-8]
	\arrow["{\chi_i}", shorten <=9pt, shorten >=9pt, Rightarrow, from=0, to=1]
\end{tikzcd}\]
Denote $l := \pi_{f_1} \circ \tau_1$. The claim now is that the object $l \in \ck//C$ together with the colax triangles $(\rho_i \tau_i, \chi_i \tau_i): l \to f_i$ (here $\chi_i$ is the comma object square pictured above right) is a coreflector-product of $f_1, f_2$ in $\ck // C$ (in particular a rali-product). We will establish why this is the case in Chapter \ref{CH_colax_adjs_lax_idemp_psmons}, in particular in Theorem \ref{THM_Kleisli2cat_is_corefl_complete}.
\end{pr}

\begin{pr}\label{PR_lax_adjoint_cooplimit}
A variant of a (left adjoint)-colimit has been studied in \cite{colimitsinrel} (with $F,W$ 2-functors and $\lax{\cj}{\cat}$ instead of $[\cj, \cat]$) under the name \textit{lax adjoint cooplimit}. In particular it has been shown in \cite[Theorem 4.2]{colimitsinrel} that under suitable conditions on a category $\cc$, the 2-category $\text{Rel}(\cc)$ of relations in $\cc$ admits (left adjoint)-colimits of $\omega$-chains.
\end{pr}

\begin{defi}\label{DEFI_preservation_of_X_limits}
We say that a 2-functor $H: \ck \to \cl$ \textit{preserves} $\cx$-limits if, whenever $\lambda: W \Rightarrow \ck(L,F-)$ exhibits $L$ as a $\cx$-limit of $F: \cj \to \ck$ weighted by $W: \cj \to \cat$, the composite pictured below  is an $\cx$-limit of $HF$ weighted by $W$:
\[\begin{tikzcd}
	W && {\ck(L, F-)} && {\cl(HL,HF-)}
	\arrow["\lambda", Rightarrow, from=1-1, to=1-3]
	\arrow["H", Rightarrow, from=1-3, to=1-5]
\end{tikzcd}\]
\end{defi}

\begin{pr}
In case $\ck, \cl$ admit comma objects, their preservation as rari-limits has been studied in \cite[Definition 7.1]{yoneda2toposes} where it has been called the \textit{preservation of lax pullbacks up to a right adjoint section}. For instance, given a finitely complete 2-category $\ck$, the 2-functor $(-) \times Z$ has this property for any $Z \in \ck$ (see \cite[Example 7.3]{yoneda2toposes}). In Weber's later work \cite[Theorem 6.1]{familial}, the class of  \textit{familial functors} has been shown to preserve comma objects as lari-limits.
\end{pr}

\begin{pozn}[Uniqueness of rali-limits]\label{POZN_rali_limits_not_unique_up_to_equiv}
Rali-limits are \textbf{not} unique up to an equivalence. It is not even the case that given two rali-limit objects $L_1, L_2$, there exists a left (or right) adjoint 1-cell $L_1 \to L_2$. For a particular example, consider the poset-enriched category $\parf$ of sets and partial functions. 
The empty set $\emptyset$ is the terminal object in $\parf$, in particular it is lali-terminal. The singleton set $*$ is lali-terminal: the ordered set $\parf(A,*)$ has a greatest element given by the total function $!_A: A \to *$. Both are also ``normalized” in the sense that $1_\emptyset$ and $1_*$ are terminal objects in the hom ordered sets they belong to.

Since this 2-category is poset-enriched, equivalent objects would be isomorphic, and an isomorphism in $\parf$ has to be a total function. Thus $*$ and $\emptyset$ can not be equivalent in $\parf$. Moreover, it can be seen that there is no left adjoint 1-cell $* \to \emptyset$. 

We may now also give an example of two non-equivalent left colax adjoints that we have promised in Remark \ref{POZN_lax_adjs_not_unique_up_to_equiv}. It can be seen that given a 2-category $\ck$, the 2-functor $* \to \ck$ picking an object $L$ is a colax left adjoint to the unique 2-functor $\ck \to *$ if and only if $L$ is rali-initial. In this colax adjunction, the modification $\Psi$ is invertible if and only if the 1-cell $1_L$ is the initial object of $\ck(L,L)$. Based on the above paragraphs, the unique 2-functor $\parf^{coop} \to *$ has two left colax adjoints that are not equivalent.
\end{pozn}

\section{Codescent objects}\label{SEC_codescent_objects}

Codescent objects, together with their dual, \textit{descent objects}, play an analogous role in two-dimensional category theory as coequalizers and equalizers do in ordinary category theory (for instance, consult \cite{biadjointtriangles_HIST}, \cite{becksthm}, \cite{axel_HIST}).

An important application of codescent objects with seemingly no one-dimensional analogue is the \textit{coherence theorem for pseudo algebras} \cite{codobjcoh} -- it is a framework that allows to prove that under certain conditions, pseudo algebras are equivalent to strict ones in a canonical way. Statements of this kind have been relevant in category theory since its beginning, the most famous example being the coherence theorem for monoidal categories \cite[3. Theorem 1]{categoriesworking_HIST}. The primary focus in our thesis will be on the \textbf{lax} version of this theorem, a statement about an existence of an adjunction between colax algebra and its strictification.

\subsection{Definition}\label{SUBS_DEFI_CODOBJ}

\begin{defi}\label{DEFI_delta_l}
Denote by $\Delta_l$ the free 2-category generated by the graph:
\[\begin{tikzcd}
	{[2]} && {[1]} && {[0]}
	\arrow["{\sigma_0}"{description}, from=1-3, to=1-5]
	\arrow["{\delta _0}"{description}, shift left=4, from=1-5, to=1-3]
	\arrow["{\delta_1}"{description}, shift right=4, from=1-5, to=1-3]
	\arrow["{\delta_2}"{description}, shift right=4, from=1-3, to=1-1]
	\arrow["{\delta_1}"{description}, from=1-3, to=1-1]
	\arrow["{\delta_0}"{description}, shift left=4, from=1-3, to=1-1]
\end{tikzcd}\]
together with 2-cells:
\begin{align*}
\iota : \sigma_0\delta_0 &\Rightarrow 1 : [0] \to [0], \\
\gamma : \delta_1\delta_0 &\Rightarrow \delta_0^2  : [2] \to [0],
\end{align*}
satisfying the following identities:
\begin{align*}
\delta_0 \delta_1 &= \delta_2 \delta_0, \\ 
\delta_2 \delta_1 &= \delta_1^2, \\ 
1_{[0]} &= \sigma_0 \delta_1. \\
\end{align*}

Give a 2-category $\ck$, we will refer to 2-functors $X: \Delta_l^{op} \to \ck$ as \textit{colax coherence data} in $\ck$.

\begin{vari}
The following variations will be relevant for us:
\begin{itemize}
\item The data is said to be \textit{pseudo} if the images of 2-cells $\gamma, \iota$ are invertible,
\item In case $\iota, \gamma$ go in the other direction, we will call this a \textit{lax coherence data}.
\item The data is said to be \textit{strict} if the images of 2-cells $\gamma, \iota$ are the identities,
\item A strict coherence data $X: \Delta^{op} \to \ck$ is further said to be \textit{reflexive} if it comes equipped with additional morphisms $j_0, j_1$ in $\ck$ (portrayed with the dotted lines below):
\[\begin{tikzcd}
	{X_2} && {X_1} && {X_0}
	\arrow["{d_2}"{description}, shift left=5, from=1-1, to=1-3]
	\arrow["{d_1}"{description}, from=1-1, to=1-3]
	\arrow["{d_0}"{description}, shift right=5, from=1-1, to=1-3]
	\arrow["{j_0}"{description, pos=0.3}, shift left=3, dashed, from=1-3, to=1-1]
	\arrow["{j_1}"{description, pos=0.3}, shift right=3, dashed, from=1-3, to=1-1]
	\arrow["{d_0}"{description}, shift right=4, from=1-3, to=1-5]
	\arrow["{d_1}"{description}, shift left=4, from=1-3, to=1-5]
	\arrow["s"{description}, from=1-5, to=1-3]
\end{tikzcd}\]
that are subject to the additional equations:
\begin{align*}
j_1 s &= j_0 s,\\
d_2 j_1 &= s d_1,\\
d_0 j_0 &= s d_0,\\
d_1 j_1 &= d_0 j_1 = 1_{X_1},\\
d_1 j_0 &= d_2 j_0 = 1_{X_1}.
\end{align*}
\end{itemize}
\end{vari}
\end{defi}

\begin{pr}\label{PR_internal_category_coherence_data}
Given a category $\ce$ with pullbacks, any category $\cc$ internal to $\ce$ gives strict coherence data in $\cat(\ce)$:
\[\begin{tikzcd}
	{\Delta_l^{op}} && {\Delta^{op}_2} && \ce && {\cat(\ce)}
	\arrow[from=1-1, to=1-3]
	\arrow["{N_\ce(\cc)}", from=1-3, to=1-5]
	\arrow["D", from=1-5, to=1-7]
\end{tikzcd}\]
Here the unlabeled functor is the identity on objects and 1-cells and sends the 2-cells to identities as well, and $N_\ce(\cc)$ is the nerve from Definition \ref{DEFI_nerve}. We will refer to this as the \textit{strict coherence data associated to $\cc$} and will (abuse notation) and denote it by $N_\ce(\cc)$.

\end{pr}

\begin{pr}\label{PR_res_of_alg}
Let $(T,m,i)$ be a 2-monad on a 2-category $\ck$. Let $(A,a,\gamma, \iota)$ be a colax $T$-algebra. There is an associated colax coherence data in $\talgs$ given by the diagram:
\[\begin{tikzcd}
	{\res{A,a,\gamma,\iota}:=} & {T^3A} && {T^2A} && TA
	\arrow["{m_{T^2A}}"{description}, shift left=4, from=1-2, to=1-4]
	\arrow["{Tm_{TA}}"{description}, from=1-2, to=1-4]
	\arrow["{T^2a}"{description}, shift right=4, from=1-2, to=1-4]
	\arrow["Ta"{description}, shift right=4, from=1-4, to=1-6]
	\arrow["{m_A}"{description}, shift left=4, from=1-4, to=1-6]
	\arrow["{Ti_A}"{description}, from=1-6, to=1-4]
\end{tikzcd}\]
together with 2-cells:
\begin{align*}
T\iota &: ai_A \Rightarrow 1_a,\\
T\gamma &: am_A \Rightarrow aTa.
\end{align*}
\noindent It is called the \textit{resolution} of the colax $T$-algebra $(A,a,\gamma,\iota)$. It is a colax (or pseudo, strict) coherence data if the algebra is colax (or pseudo, strict).

For a strict $T$-algebra, the underlying data of $\res{A,a}$ in $\ck$ is moreover reflexive -- the additional sections are given by the 1-cells $i_{TA}, Ti_A: TA \to T^2A$.
\end{pr}

\begin{pozn}\label{POZN_Tcartesian_resAa_internal_cat}
Note that in case the 2-monad $T$ is cartesian, the resolution of any \textbf{strict} $T$-algebra $(A,a)$ in $\talgs$ has an additional algebraic structure -- it is the nerve of an internal category in $\ck$. In case $\ck = \cat$, this is a double category. 
\end{pozn}

\begin{defi}\label{DEFI_codescent_object}
Denote by $W : \Delta_l \to \cat$ the 2-functor that regards the objects of $\Delta_l$ as ordinals $[n]= \{ 0 < \dots < n \}$ and sends the arrows to the usual face/degeneracy maps as the notation suggests.

Given a 2-category $\ck$ and a colax coherence data $X$ in $\ck$, a $W$-weighted cocone for $X$ (a 2-natural transformation $W \Rightarrow \ck(X-,A)$) will be referred to as a \textit{coherence cocone} or just a \textit{cocone} for $X$.

By the \textit{codescent object} of $X$ we mean the 2-colimit of $X$ weighted by $W$.
\end{defi}

\begin{nota}\label{NOTA_coherence_data}
Given colax coherence data $X: \Delta_l^{op} \to \ck$, we will denote:
\[
\coh{\ck}{X}{A} := [\Delta_l, \cat](W,\ck(A,X-)).
\]
\end{nota}

\noindent It will be useful to express the data of a cocone and the codescent object for $X$ explicitly. The colax coherence data in $\ck$ is the following diagram in $\ck$:
\[\begin{tikzcd}
	{X_2} && {X_1} && {X_0}
	\arrow["s"{description}, from=1-5, to=1-3]
	\arrow["{d_0}"{description}, shift right=4, from=1-3, to=1-5]
	\arrow["{d_1}"{description}, shift left=4, from=1-3, to=1-5]
	\arrow["{d_2}"{description}, shift left=4, from=1-1, to=1-3]
	\arrow["{d_1}"{description}, from=1-1, to=1-3]
	\arrow["{d_0}"{description}, shift right=4, from=1-1, to=1-3]
\end{tikzcd}\]
together with 2-cells and equations between 1-cells:
\begin{align*}
\iota: d_0 s &\Rightarrow 1_{X_0},\\
\gamma: d_0 d_1 &\Rightarrow d_0^2,\\
d_1 d_0 &= d_0 d_2,\\
d_1 d_2 &= d_1^2,\\
1_{X_0} &= d_1 s.
\end{align*}

\noindent A $W$-weighted cocone with apex $Y$ consists of a tuple $(F: X_0 \to Y, \xi: Fd_1 \Rightarrow Fd_0)$ of a 1-cell and a 2-cell satisfying:
\[\begin{tikzcd}[column sep=scriptsize]
	{X_2} && {X_1} &&& {X_2} && {X_1} \\
	& {X_1} && {X_0} &&&&& {X_0} \\
	{X_1} &&&& {=} & {X_1} && {X_0} \\
	& {X_0} && Y &&& {X_0} && Y
	\arrow["{d_0}", from=1-1, to=1-3]
	\arrow["{d_1}"', from=1-1, to=2-2]
	\arrow["{d_2}"', from=1-1, to=3-1]
	\arrow["{d_0}", from=1-3, to=2-4]
	\arrow["{d_0}", from=1-6, to=1-8]
	\arrow["{d_2}"', from=1-6, to=3-6]
	\arrow[""{name=0, anchor=center, inner sep=0}, "{d_0}", from=1-8, to=2-9]
	\arrow["{d_1}"', from=1-8, to=3-8]
	\arrow[""{name=1, anchor=center, inner sep=0}, "{d_0}"{description}, from=2-2, to=2-4]
	\arrow["{d_1}"', from=2-2, to=4-2]
	\arrow["F", from=2-4, to=4-4]
	\arrow["F", from=2-9, to=4-9]
	\arrow["{d_1}"', from=3-1, to=4-2]
	\arrow["{d_0}", from=3-6, to=3-8]
	\arrow["{d_1}"', from=3-6, to=4-7]
	\arrow[""{name=2, anchor=center, inner sep=0}, "F"{description}, from=3-8, to=4-9]
	\arrow[""{name=3, anchor=center, inner sep=0}, "F"', from=4-2, to=4-4]
	\arrow[""{name=4, anchor=center, inner sep=0}, "F"', from=4-7, to=4-9]
	\arrow["\gamma", shorten <=3pt, Rightarrow, from=1, to=1-3]
	\arrow["\xi"', shorten <=13pt, shorten >=13pt, Rightarrow, from=2, to=0]
	\arrow["\xi"', shorten <=9pt, shorten >=9pt, Rightarrow, from=3, to=1]
	\arrow["\xi", shorten <=3pt, Rightarrow, from=4, to=3-8]
\end{tikzcd}\]
and:
\[\begin{tikzcd}
	{X_0} \\
	\\
	{X_1} && {X_0} & {=} & {1_F} \\
	\\
	{X_0} && Y
	\arrow["s"', from=1-1, to=3-1]
	\arrow[""{name=0, anchor=center, inner sep=0}, Rightarrow, no head, from=1-1, to=3-3]
	\arrow[curve={height=24pt}, Rightarrow, no head, from=1-1, to=5-1]
	\arrow[""{name=1, anchor=center, inner sep=0}, "{d_0}"{description}, from=3-1, to=3-3]
	\arrow["{d_1}"{description}, from=3-1, to=5-1]
	\arrow["F", from=3-3, to=5-3]
	\arrow[""{name=2, anchor=center, inner sep=0}, "F"', from=5-1, to=5-3]
	\arrow["\iota", shorten <=4pt, shorten >=4pt, Rightarrow, from=1, to=0]
	\arrow["\xi", shorten <=9pt, shorten >=9pt, Rightarrow, from=2, to=1]
\end{tikzcd}\]

\noindent Such a cocone $(F,\xi)$ with apex $Y$ is the \textit{codescent object} for the above coherence data if it satisfies the one-dimensional and two-dimensional universal property respectively:
\begin{itemize}
\item Given any other cocone $(G,\psi)$ with apex $Z$, there exists a unique map $\theta: Y \to Z$ such that:
\begin{align}
\theta F &= G,\\
\theta \xi &= \psi.
\end{align}
\item Given any \textit{morphism of cocones} $\rho: (G_1, \psi_1) \to (G_2, \psi_2)$ between cocones with the same apex $Z$, that is, a 2-cell $\rho: G_1 \Rightarrow G_2$ satisfying:
\[\begin{tikzcd}[column sep=scriptsize]
	&& {X_0} &&&&&& {X_0} \\
	\\
	{X_1} &&&& Z & {=} & {X_1} &&&& Z \\
	\\
	&& {X_0} &&&&&& {X_0}
	\arrow["{G_2}", from=1-3, to=3-5]
	\arrow[""{name=0, anchor=center, inner sep=0}, "{G_1}"', from=1-9, to=3-11]
	\arrow[""{name=1, anchor=center, inner sep=0}, "{G_2}", curve={height=-30pt}, from=1-9, to=3-11]
	\arrow["{d_0}", from=3-1, to=1-3]
	\arrow["{d_1}"', from=3-1, to=5-3]
	\arrow["{d_0}", from=3-7, to=1-9]
	\arrow["{d_1}"', from=3-7, to=5-9]
	\arrow["{\psi_2}", shorten <=15pt, shorten >=15pt, Rightarrow, from=5-3, to=1-3]
	\arrow[""{name=2, anchor=center, inner sep=0}, "{G_2}", from=5-3, to=3-5]
	\arrow[""{name=3, anchor=center, inner sep=0}, "{G_1}"', curve={height=30pt}, from=5-3, to=3-5]
	\arrow["{\psi_1}", shorten <=15pt, shorten >=15pt, Rightarrow, from=5-9, to=1-9]
	\arrow["{G_1}"', from=5-9, to=3-11]
	\arrow["\rho", shorten <=5pt, shorten >=5pt, Rightarrow, from=0, to=1]
	\arrow["\rho"', shorten <=5pt, shorten >=5pt, Rightarrow, from=3, to=2]
\end{tikzcd}\]

there exists a unique 2-cell $\rho' : \theta_1 \Rightarrow \theta_2 : Y \to Z$ between the canonical maps from the previous point such that:
\[
\rho' F = \rho.
\]
\end{itemize}

\begin{vari}
We will encounter the following variations:
\begin{itemize}
\item The definition of the \textit{iso-codescent object} is the same as above, except all the cocones involved have the corresponding 2-cell invertible.
\item By the \textit{codescent object} of \textbf{lax} coherence data in $\ck$ we mean the codescent object of colax coherence data in $\ck^{co}$.
\end{itemize}
\end{vari}

\begin{pr}
Any \textit{pointed magma}, that is, a set $M$ together with a binary operation $*: M \times M \to M$ and a distinguished element $e \in M$, gives coherence data in $\cat$ if we regard the sets $*, M, M^2$ as discrete categories:
\[\begin{tikzcd}
	{M \times M} && M && {*}
	\arrow["{\pi_1}"{description}, shift left=4, from=1-1, to=1-3]
	\arrow["{*}"{description}, from=1-1, to=1-3]
	\arrow["{\pi_2}"{description}, shift right=4, from=1-1, to=1-3]
	\arrow["{!}"{description}, shift right=4, from=1-3, to=1-5]
	\arrow["{!}"{description}, shift left=4, from=1-3, to=1-5]
	\arrow["e"{description}, from=1-5, to=1-3]
\end{tikzcd}\]
It can be seen that the codescent object of this diagram is a monoid (regarded as a one-object category) that is the reflection of the pointed magma $M$ along the inclusion functor $\text{Magma}_* \xhookrightarrow{} \text{Mon}$ of pointed magmas into monoids.
\end{pr}

\noindent An example of a codescent object essential for this thesis is given by the following proposition:

\begin{prop}\label{PR_category_is_codobj_of_discretes}
Let $\ce$ be a category with pullbacks and $\cc$ a category internal to $\ce$. Then $\cc$ is the codescent object in $\cat(\ce)$ of its strict coherence data $N_\ce(\cc)$ (of Example \ref{PR_internal_category_coherence_data}).
\end{prop}
\begin{proof}
By Theorem \ref{THM_catE_limits} and Proposition \ref{THMP_powers_with2_implies_1dim_impl_2dim}, it is enough to verify the one-dimensional universal property -- in particular, establish a bijection between the following sets (natural in $\cz$):
\[
\cat(\ce)(\cc,\cz) \cong \coh{\cat(\ce)}{N_{\ce}(\cc)}{\cz}.
\]
Notice first that an element of the right hand side is a pair:
\[(G: \cc_0 \to \cz, \phi: Gs \Rightarrow Gt: \cc_1 \to \cz),\]
of a functor and a natural transformation. Since the $\cc_i$'s are discrete categories, this data amounts to two 1-cells in $\ce$:
\begin{align*}
G: \cc_0 \to \cz_0,\\
\psi: \cc_1 \to \cz_1.
\end{align*}
It is immediate that the cocone axioms for $(G,\psi)$ are satisfied:
\begin{align*}
\comp \circ \psi &= \pi_2 \psi \circ \pi_1 \psi,\\
\psi s &= 1_{X_0},
\end{align*}
if and only if this pair of 1-cells gives a functor $\cc \to \cz$.
\end{proof}

\begin{pr}\label{PR_Kleisliobject_is_codescentobject}
Given a 2-monad $(T,m,i)$ on a 2-category $\ck$ and a colax $T$-algebra $(A,a,\gamma, \iota)$, notice that the pair $(a: TA \to A, \gamma)$ is always a cocone for its resolution $\res{A,a,\gamma, \iota}$ in\footnote{This cocone only lives in $\ck$ and does \textbf{not} lift to $\talgs$.} $\ck$. It is typically not its codescent object though\footnote{It is always a weak codescent object in an appropriate sense however, see Lemma \ref{POZNBAZANT}.}.

For example, consider the identity 2-monad $T := 1_{\ck}$ and a colax algebra $(A,q,\delta,\epsilon)$, that is, a comonad in $\ck$. The codescent object is the Kleisli object for the comonad, which will be typically different from the object $A$.
\end{pr}

\begin{pozn}\label{POZN_reflexive_coherence_data}
Let us also mention that the codescent object of reflexive strict coherence data is a \textit{sifted colimit} -- it commutes with finite products in $\cat$, see \cite[Proposition 8.41]{johnphd}. Such codescent objects play an analogous role to that of coequalizers of reflexive pairs as these are also sifted -- see \cite[Remark 3.2]{algtheories_HIST}.
\end{pozn}

More examples of codescent objects will be given in Chapter \ref{CH_COMPUTATION}.

\begin{warning}
Let us warn the reader about terminology regarding codescent objects in the literature:

\begin{itemize}
\item In this thesis, \cite{johnphd}, \cite{internalalg}, \cite{codobjcoh} the diagram is always called the ``coherence data” and the 2-colimit is the ``codescent object”.

\item In the last paper mentioned, what we call the ``codescent object” and the ``iso-codescent object” is called the ``lax codescent object” and the ``codescent object”.

\item The papers \cite{becksthm}, \cite{axel_HIST} follow a different terminology altogether and use the term ``codescent object” for the diagram (the diagrams they use are what we would call pseudo coherence data), and ``pseudo-coequalizer” for the bicolimit (that we would call a bi-iso-codescent object).
\end{itemize}
\end{warning}

We end the section with a technical remark comparing the codescent objects used in this thesis and \cite{codobjcoh} -- the reader not interested in subtleties may feel free to skip ahead.

\begin{pozn}\label{POZN_technical_lacks_coh_data}
Our \textit{lax coherence data} is a stricter version of Steve Lack's lax coherence data from \cite[Page 228]{codobjcoh}. In particular, we consider his 2-cells $\delta, \kappa, \rho$ to be the identities.

The big difference here is that our codescent object of lax coherence data has the 2-cell component going the \textbf{other way} than in \cite{codobjcoh}. The codescent object of lax coherence data $X$ in the sense of Lack is a pair:
\[
(F: X_0 \to \cy, \xi: Fd_1 \Rightarrow Fd_0),
\]
that satisfies the equations:
\begin{align*}
\xi d_1 &= F \gamma \circ \xi d_0 \circ \xi d_2,\\
F \iota &= \xi s.
\end{align*}

The reason we use the modified codescent object in this thesis is that we will be interested in the left 2-adjoint to the inclusion $\talgs \xhookrightarrow{} \colaxtalgl$ for a 2-monad $T$, as opposed to $\talgs \xhookrightarrow{} \laxtalgl$ as in \cite[Lemma 2.3]{codobjcoh}.

Finally, we remark that according to \cite[Proposition 2.1]{codobjcoh}, Lack's codescent object can be built using simpler colimits -- a coinserter and two coequifiers. By analogous arguments, the same is true for our version of codescent object.
\end{pozn}

\subsection{Strictifying lax structures}

\begin{termi}\label{TERMI_lax_coherence_holds}
Fix a 2-monad $(T,m,i)$ on a 2-category $\ck$. The left 2-adjoints to the following inclusion 2-functors, whenever they exist, will be referred to as the \textit{strictification 2-functors} for colax algebras and lax morphisms respectively:
\begin{equation}\label{EQ_talgs_colaxtalgl_talgl}
\begin{tikzcd}
	\talgs && \colaxtalgl && \talgs && \talgl
	\arrow[""{name=0, anchor=center, inner sep=0}, "J"', curve={height=24pt}, hook, from=1-1, to=1-3]
	\arrow[""{name=1, anchor=center, inner sep=0}, "{(-)'}"', curve={height=24pt}, dashed, from=1-3, to=1-1]
	\arrow[""{name=2, anchor=center, inner sep=0}, "{J_{\bullet}}"', curve={height=24pt}, hook, from=1-5, to=1-7]
	\arrow[""{name=3, anchor=center, inner sep=0}, "{(-)'}"', curve={height=24pt}, dashed, from=1-7, to=1-5]
	\arrow["\dashv"{anchor=center, rotate=-90}, draw=none, from=1, to=0]
	\arrow["\dashv"{anchor=center, rotate=-90}, draw=none, from=3, to=2]
\end{tikzcd}
\end{equation}

\medskip

We say that \textit{coherence for colax algebras holds} for a 2-monad $T$ if the left 2-adjoint to $\talgs \xhookrightarrow{} \colaxtalgl$ exists, and moreover the components of the unit of the corresponding 2-adjunction admit left adjoints in $\ck$.

We say that \textit{coherence for lax morphisms holds} for a 2-monad $T$ if the left 2-adjoint to $\talgs \xhookrightarrow{} \talgl$ exists, and moreover the components of the unit of the corresponding 2-adjunction admit left adjoints in $\ck$.
\end{termi}

\begin{termi}\label{DEFI_Ql_Qp}
In case the 2-adjoint above right exists, there is an induced 2-comonad $Q_l$ on $\talgs$ that we will refer to as the \textit{lax morphism classifier 2-comonad}. Using ``colax” or ``pseudo” in place of ``lax”, we denote the resulting 2-comonads by $Q_c$ and $Q_p$ and call them the \textit{colax} and \textit{pseudo morphism classifier 2-comonads}.
\end{termi}

\noindent Let us capture the relationship that $Q_l$ has with $Q_p$:

\begin{prop}\label{THMP_mor_of_comonads_QlQp}
Assume that $T$ is a 2-monad for which the inclusions:
\begin{align*}
\talgs \xhookrightarrow{} \talgl, \hspace{2cm} \talgs \xhookrightarrow{} \talg,
\end{align*}
admit left 2-adjoints so that we have the induced 2-comonads $Q_l, Q_p$ on $\talgs$. Then there is a morphism of 2-comonads $Q_l \to Q_p$.
\end{prop}
\begin{proof}
Denote the units of these adjunctions by $p_A: A \rightsquigarrow A'$ and $p_A^\dagger: A \rightsquigarrow A^\dagger$ respectively. Since $p^\dagger_A$ is a pseudo morphism, it is in particular a lax morphism and thus there exists a unique strict morphism $\theta_A : A' \to A^\dagger$ making the diagram commute:
\[\begin{tikzcd}
	&& {A'} \\
	A && {A^\dagger} && {\text{ in }\talgs}
	\arrow["{\exists ! \theta_A}", dashed, from=1-3, to=2-3]
	\arrow["{p_A}", squiggly, from=2-1, to=1-3]
	\arrow["{p^\dagger_A }"', squiggly, from=2-1, to=2-3]
\end{tikzcd}\]

\noindent Using the universal property of $A'$, it is readily seen that the maps $\theta_A$'s assemble into a morphism of 2-comonads.
\end{proof}

\begin{pozn}\label{POZN_Talg_is_kleisli}
It is easy to see that $\talgl$ is isomorphic to $(\talgs)_{Q_l}$, the Kleisli 2-category for the 2-comonad $Q_l$. Similarly, $\talg \cong (\talgs)_{Q_p}$.
\end{pozn}

\noindent The theorem below gives necessary and sufficient conditions for the existence of the strictification 2-functors. It is a significant generalization of Theorem \ref{THM_strictification_of_comonads} that has been proven by Steve Lack in \cite{codobjcoh}. The proof we present here is entirely analogous to \cite[Lemma 2.3]{codobjcoh}\footnote{It is not \textit{identical} to it since we use a different version of the codescent object, see Remark \ref{POZN_technical_lacks_coh_data}.}.

\begin{theo}\label{THM_adjunkce_talgs_talgl_colaxtalgl}
Let $T$ be a 2-monad on a 2-category $\ck$.
\begin{itemize}
\item The 2-category $\talgs$ admits codescent objects of resolutions of strict $T$-algebras if and only if the inclusion 2-functor $\talgs \to \talgl$ admits a left 2-adjoint.
\item The 2-category $\talgs$ admits codescent objects of resolutions of colax $T$-algebras if and only if the inclusion 2-functor $\talgs \to \colaxtalgl$ admits a left 2-adjoint.
\end{itemize}
\[\begin{tikzcd}
	\talgs && \colaxtalgl && \talgs && \talgl
	\arrow[""{name=0, anchor=center, inner sep=0}, curve={height=24pt}, hook, from=1-1, to=1-3]
	\arrow[""{name=1, anchor=center, inner sep=0}, "{(-)'}"', curve={height=24pt}, dashed, from=1-3, to=1-1]
	\arrow[""{name=2, anchor=center, inner sep=0}, curve={height=24pt}, hook, from=1-5, to=1-7]
	\arrow[""{name=3, anchor=center, inner sep=0}, "{(-)'}"', curve={height=24pt}, dashed, from=1-7, to=1-5]
	\arrow["\dashv"{anchor=center, rotate=-90}, draw=none, from=1, to=0]
	\arrow["\dashv"{anchor=center, rotate=-90}, draw=none, from=3, to=2]
\end{tikzcd}\]
\end{theo}
\begin{proof}
It is readily verified that there is an isomorphism, 2-natural in strict algebras $(B,b)$ and strict morphisms:
\begin{align*}
\colaxtalgl((A,a,\gamma,\iota),(B,b)) &\cong \coh{\ck}{\text{Res}(A,a,\gamma,\iota)}{(B,b)},\\
(f: A \to B, \overline{f}) &\mapsto (b Tf: TA \to B, b T\overline{f}).
\end{align*}

\noindent A strict $T$-algebra $(A',a')$ is the codescent object of $\text{Res}(A,a,\gamma,\iota)$ if and only if there is an isomorphism, 2-natural in $(B,b)$:
\[
\talgs((A',a'),(B,b)) \cong \coh{\ck}{\text{Res}(A,a,\gamma,\iota)}{(B,b)}.
\]
Composing these two isomorphisms, we obtain the statement.
\end{proof}

\noindent Analogous results hold for pseudo algebras and pseudo morphisms -- except instead of codescent objects, iso-codescent objects are used (see \textbf{Variations} after \ref{NOTA_coherence_data}) -- this has been done in \cite[Theorem 2.6]{codobjcoh}.

\begin{pozn}\label{POZN_komponenta_jednotky_talgs_colaxtalgl}
Assume that codescent objects of resolutions of colax algebras exist in $\talgs$. For a colax $T$-algebra $\mathbb{A} := (A,a,\gamma,\iota)$, denote by $(e, \overline{e})$ the codescent object  of its resolution in $\talgs$. The unit of the above 2-adjunction is then the lax morphism
\[
(ei_A,\overline{e}i_{TA}) : \mathbb{A} \to \mathbb{A}'.
\]

To obtain the component of the counit at a strict algebra $(B,b)$, note that the pair $(b,1_{b m_B})$ is a cocone for $\res{B,b}$ in $\talgs$, so there exists a unique strict algebra morphism $q_{(B,b)}: (B,b)' \to (B,b)$ such that:
\begin{align*}
q_{(B,b)}e_{(B,b)} &= b,\\
q_{(B,b)}\overline{e}_{(B,b)} &= 1.
\end{align*}

\noindent In particular, every strict $T$-algebra $(B,b)$ is a retract in $\talgl$ of an algebra of form $\mathbb{C}'$ for a strict $T$-algebra $\mathbb{C}$.
\end{pozn}

\begin{pozn}
Let us remark that in certain rare cases, the strictification 2-functors may exist due to formal reasons only (without needing to introduce any colimits). Consider a lax-idempotent 2-monad $(T,m,i,\lambda)$ on a 2-category $\ck$ and a pseudo $T$-algebra $(A,a,\gamma,\iota)$. Notice that we have the following 2-natural isomorphisms (this is because the latter hom is isomorphic to $\ck(A,B)$ by Corollary \ref{THMC_lax_idemp_1cell_is_uniquely_lax}):
\[
\talgs((TA,m_A), (B,b)) \cong \ck(A,B) \cong \pstalgl((A,a,\gamma,\iota),(B,b)).
\]

\noindent This means that for lax-idempotent 2-monads, the inclusion $\talgs \xhookrightarrow{} \pstalgl$ always admits a left 2-adjoint -- it sends a pseudo $T$-algebra to the free algebra on its underlying object.
\end{pozn}

\subsection{Lax coherence theorems}\label{SUBS_lax_coh_THMS}

In the remainder of this section, we will concern ourselves with coherence for colax algebras for a 2-monad. This is a straightforward analogue of the coherence for pseudo algebras of \cite[Theorem 3.2]{codobjcoh} and has been observed by Steve Lack in \cite[Theorem 8.5]{moritactxs}. This has also been the subject of the author's Master's thesis \cite[Theorem 3.16]{milomgr}. Another person to treat this topic is Kengo Hirata, see \cite{kengohirata}.

Let $(T,m,i)$ be a 2-monad on a 2-category $\ck$. The following is a lax analogue of Lemma 2.3 of \cite{becksthm} and has been proven in \cite[Theorem 3.14]{milomgr}. Recall the weak limit terminology from Section \ref{SEC_weak_limits}.

\begin{lemma}\label{POZNBAZANT}
Given a colax $T$-algebra $(A,a,\gamma,\iota)$, the cocone $(a,\gamma)$ is a (right adjoint)-codescent object of the diagram $\res{A,a,\gamma,\iota}$ in $\ck$.
\end{lemma}
\begin{proof}
We have to show that for every $B \in \ck$, the following functor admits a left adjoint:
\begin{align*}
\kappa: \ck(A,B) &\to \coh{\ck}{\res{A,a,\gamma,\iota}}{B},\\
\kappa: \theta &\mapsto (\theta a, \theta \gamma).
\end{align*}
This left adjoint will be given by the following functor:
\begin{align*}
(i_A)^* &: \coh{\ck}{\res{A,a,\gamma,\iota}}{B} \to \ck(A,B),\\
(i_A)^* &: (f,\overline{f}) \mapsto fi_A.
\end{align*}

\noindent The counit $\epsilon : (i_A)^* \kappa \Rightarrow 1_{\ck(A,B)}$ of the adjunction evaluated at $\theta : A \to B$ is given by:
\[
\epsilon_\theta := \theta \iota : \theta a i_A \Rightarrow \theta.
\]
The unit $\eta: 1_{\coh{\ck}{\res{A,a,\gamma,\iota}}{B}} \Rightarrow \kappa (i_A)^*$ of the adjunction evaluated at a cocone $(e,\overline{e})$ is given by:
\[
\eta_{(e,\overline{e})} := \overline{e}i_{TA}: (e,\overline{e}) \to (ei_A a, e i_A \gamma).
\]

\noindent The first triangle identity follows from the cocone axiom:
\[
\epsilon_{i^*_A (e,\overline{e})} \circ i_A^* \eta_{(e,\overline{e})} = e T\iota i_A \circ \overline{e} i_{TA} i_A = 1.
\]

\noindent The second triangle identity follows from the $T$-algebra unit axiom:
\[
\kappa \epsilon_{\theta} \circ \eta_{\kappa \theta} = \theta \iota a \circ \theta \gamma i_{TA} = 1.
\]

\end{proof}

\begin{theo}[Coherence for colax algebras]\label{THM_lax_coherence_colax_algebras_EXPOSITION}
Assume $T$ is a 2-monad on a 2-category $\ck$ for which $\talgs$ admits codescent objects of resolutions of colax algebras. Assume that $T$ (or equivalently, the forgetful 2-functor $U: \talgs \to \ck$) preserves these codescent objects. Then coherence for colax algebras holds for $T$.
\end{theo}
\begin{proof}
Denote by $(A',a')$ the codescent object in $\talgs$ of the resolution of a colax $T$-algebra $(A,a,\gamma,\iota)$ and denote by $(e:TA \to A', \overline{e})$ the colimit cocone. If $U$ preserves it, there is an isomorphism as displayed below left:

\begin{equation}\label{EQ_PF_kAprime_cohK_kA}
\begin{tikzcd}
	{\ck(A',-)} && {\coh{\ck}{\res{A,a,\gamma,\iota}}{-}} && {\ck(A,-)}
	\arrow["\cong"{description}, draw=none, from=1-1, to=1-3]
	\arrow[""{name=0, anchor=center, inner sep=0}, "{(i_A)^*}"', curve={height=24pt}, from=1-3, to=1-5]
	\arrow[""{name=1, anchor=center, inner sep=0}, "\kappa"', curve={height=24pt}, from=1-5, to=1-3]
	\arrow["\dashv"{anchor=center, rotate=-90}, draw=none, from=1, to=0]
\end{tikzcd}
\end{equation}
By Lemma \ref{POZNBAZANT}, there is an adjunction that is displayed above right. By the Yoneda lemma for 2-categories, we obtain the following adjunction in $\ck$.
\[\begin{tikzcd}
	A && {A'}
	\arrow[""{name=0, anchor=center, inner sep=0}, "{ei_A}"', curve={height=18pt}, from=1-1, to=1-3]
	\arrow[""{name=1, anchor=center, inner sep=0}, curve={height=18pt}, dashed, from=1-3, to=1-1]
	\arrow["\dashv"{anchor=center, rotate=-90}, draw=none, from=1, to=0]
\end{tikzcd}\]

\end{proof}

\begin{pozn}\label{POZN_explicit_description}
Explicitly, given a colax $T$-algebra $\mathbb{A} = (A,a,\gamma,\iota)$, the left adjoint to the 1-cell $ei_A: A \to A'$ is the unique 1-cell $\widetilde{q}_{\mathbb{A}}: A' \to A$ in $\ck$ with the property that:
\begin{align*}
\widetilde{q}_{\mathbb{A}} e &= a,\\
\widetilde{q}_{\mathbb{A}} \overline{e} &= \gamma.
\end{align*}
The counit of the adjunction is the coassociator for $\mathbb{A}$:
\[
\iota : ai_A = \widetilde{q} ei_A \Rightarrow 1_A.
\]
The unit $\widetilde{\eta}: 1_{A'} \Rightarrow ei_A \widetilde{q}_{\mathbb{A}}$ is the unique 2-cell with the property that:
\[
\widetilde{\eta} e = \overline{e} i_{TA}.
\]

\noindent In total, we have:
\begin{equation}\label{EQ_A_Aprime}
\begin{tikzcd}
	{(\iota, \widetilde{\eta}):} & A && {A'} & {\text{ in } \ck}
	\arrow[""{name=0, anchor=center, inner sep=0}, "{ei_A}"', curve={height=18pt}, from=1-2, to=1-4]
	\arrow[""{name=1, anchor=center, inner sep=0}, "{\widetilde{q}_{\mathbb{A}}}"', curve={height=18pt}, from=1-4, to=1-2]
	\arrow["\dashv"{anchor=center, rotate=-90}, draw=none, from=1, to=0]
\end{tikzcd}
\end{equation}
\end{pozn}

\begin{pozn}\label{POZN_adjunkce_AAprime_cant_lift}
Note that in contrast to the coherence theorem for pseudo algebras (\cite[Theorem 3.2]{codobjcoh}), we can not in general lift \eqref{EQ_A_Aprime} up to $\colaxtalgl$. It can always be done for strict algebras and lax morphisms though, see Theorem \ref{THM_coh_for_lax_mors}. In the colax algebra case, the best we can do is give $\widetilde{q}$ the structure of a \textbf{colax} morphism (see Theorem \ref{THM_doctrinal_adj}).
\end{pozn}

From the explicit description in Remark \ref{POZN_explicit_description}, we may immediately deduce:

\begin{cor}[Coherence for normal colax algebras]\label{THMC_coherence_for_normal_colax}
Assume $T$ is a 2-monad on a 2-category $\ck$. Assume that $\talgs$ admits codescent objects of resolutions of colax algebras and that $T$ preserves these codescent objects. For every \textbf{normal} colax $T$-algebra $\mathbb{A} = (A,a,\gamma,\iota)$, the object $A$ is reflectively embedded in the object $A'$ underlying the strictification $\mathbb{A}'$.
\end{cor}

Let us record here conditions that guarantee lax-idempotency of the lax morphism classifier 2-comonad $Q_l$:

\begin{prop}\label{THMP_Ql_lax_idemp_Qp_ps_idemp}
Let $T$ be a 2-monad on a 2-category $\ck$ such that the left 2-adjoint to the inclusion $\talgs \xhookrightarrow{} \talgl$ exist.  Then the lax morphism classifier 2-comonad $Q_l$ is lax-idempotent if either:
\begin{itemize}
\item $\ck$ admits oplax limits of arrows,
\item $T$ preserves codescent objects of resolutions of strict algebras.
\end{itemize}
\end{prop}
\begin{proof}
For the first part, see \cite[Lemma 2.5]{enhanced2cats}. Assume now that $T$ preserves codescent objects of resolutions of strict algebras. Let $(B,b)$ be a strict algebra. Theorem \ref{THM_lax_coherence_colax_algebras_EXPOSITION} guarantees the existence of the following reflection in $\ck$:
\[\begin{tikzcd}
	{(1, \widetilde{\eta}):} & B && {B'} & {\text{ in } \ck}
	\arrow[""{name=0, anchor=center, inner sep=0}, "{ei_B}"', curve={height=18pt}, from=1-2, to=1-4]
	\arrow[""{name=1, anchor=center, inner sep=0}, "{\widetilde{q}_{\mathbb{B}}}"', curve={height=18pt}, from=1-4, to=1-2]
	\arrow["\dashv"{anchor=center, rotate=-90}, draw=none, from=1, to=0]
\end{tikzcd}\]

By the explicit description of this adjunction in Remark \ref{POZN_explicit_description} and by the explicit description of the counit of the 2-adjunction between strict and lax morphisms in  Remark \ref{POZN_komponenta_jednotky_talgs_colaxtalgl} we may conclude that in fact $q_{\mathbb{B}} \sstackrel{!}{=} \widetilde{q}_{\mathbb{B}}$. Thus the left adjoint lifts to $\talgs$ and in particular to $\talgl$. To lift the whole adjunction to $\talgl$, it must be shown that:
\begin{itemize}
\item $\widetilde{\eta}: 1_{\mathbb{B}'} \Rightarrow (ei_B, \overline{e}i_{TB})q_{\mathbb{B}}: \mathbb{B}' \to \mathbb{B}'$ is an algebra 2-cell,
\item the triangle identities $q_{\mathbb{B}} \widetilde{\eta} = 1$, $\widetilde{\eta} (ei_B, \overline{e}i_{TB}) = 1$ hold.
\end{itemize}
Each of these two statements is proven in a straightforward way using the universal property of the codescent object $B'$ so we omit the proof. We have now proven that the 2-adjunction between $\talgs$ and $\talgl$ is colax-idempotent. By Proposition \ref{THMP_KZ_adjunkce}, $Q_l$ is lax-idempotent.
\end{proof}

Analogous results hold for the pseudo morphism classifier 2-comonad $Q_p$, except that pseudo limits of arrows need to be used -- this is also captured by \cite[Lemma 2.5]{enhanced2cats}. A re-formulation of the above proposition gives us a result about coherence for lax morphisms:

\begin{theo}[Coherence for lax morphisms]\label{THM_coh_for_lax_mors}
Let $(T,m,i)$ be a 2-monad on a 2-category $\ck$ for which the lax morphism strictification 2-functor exists:
\begin{equation}\label{EQ_talgs_talgl}
\begin{tikzcd}
	{(q,p):} & \talgs && \talgl
	\arrow[""{name=0, anchor=center, inner sep=0}, "J"', curve={height=24pt}, hook, from=1-2, to=1-4]
	\arrow[""{name=1, anchor=center, inner sep=0}, "{(-)'}"', curve={height=24pt}, dashed, from=1-4, to=1-2]
	\arrow["\dashv"{anchor=center, rotate=-90}, draw=none, from=1, to=0]
\end{tikzcd}
\end{equation}
Then coherence for lax morphisms holds for $T$ if either of the following conditions is true:
\begin{itemize}
\item $\ck$ admits oplax limits of arrows,
\item $T$ preserves codescent objects of resolutions of strict algebras.
\end{itemize}
In fact, in both cases the left adjoint is given by the counit of \eqref{EQ_talgs_talgl} and the whole adjunction lives in $\talgl$:
\[\begin{tikzcd}
	{(A,a)} && {(A',a')} & {}
	\arrow[""{name=0, anchor=center, inner sep=0}, "{p_{\mathbb{A}}}"', curve={height=24pt}, squiggly, from=1-1, to=1-3]
	\arrow[""{name=1, anchor=center, inner sep=0}, "{q^{\bullet}_{(A,a)}}"', curve={height=24pt}, two heads, from=1-3, to=1-1]
	\arrow["\dashv"{anchor=center, rotate=-90}, draw=none, from=1, to=0]
\end{tikzcd}\]
\end{theo}

\noindent Let us also note that there is a possible generalization of Theorem \ref{THM_lax_coherence_colax_algebras_EXPOSITION}, where $T$ preserving the relevant codescent objects may be weakened to $T$ preserving them only as a (right adjoint)-colimit (of Section \ref{SEC_weak_limits}), see Chapter \ref{CH_generalized_colax_coh}.

\subsection{Examples}

We will now consult examples that we consider well known. More examples will be encountered in Section \ref{SEC_examples_lax_coh}.

\begin{pr}\label{PR_identity2monad_laxcoh}
Coherence for colax algebras holds for the identity 2-monad on a 2-category $\ck$ provided it admits Kleisli objects of comonads. Theorem \ref{THM_lax_coherence_colax_algebras_EXPOSITION} specializes to Theorem \ref{THM_laxcoherencethm_for_identity2monad} which we have seen in Section \ref{SEKCE_monads_in_2cat}.

\end{pr}

\begin{pozn}\label{POZN_strictification_colax_algs_AND_Kleisli_objects}
Given a 2-monad $(T,m,i)$ on a 2-category $\ck$ that admits all codescent objects, we have the two 2-adjunctions pictured below (here $\mathbb{U}^T$ is the 2-functor that sends a colax algebra to its underlying comonad, see Remark \ref{POZN_laxalg_underlying_monad}):
\[\begin{tikzcd}
	\talgs && \colaxtalgl \\
	\\
	\ck && {\comnd{\ck}}
	\arrow[""{name=0, anchor=center, inner sep=0}, curve={height=24pt}, from=1-1, to=1-3]
	\arrow["{U^T}"', from=1-1, to=3-1]
	\arrow[""{name=1, anchor=center, inner sep=0}, "{(-)'}"', curve={height=24pt}, dashed, from=1-3, to=1-1]
	\arrow["{\mathbb{U}^T}", from=1-3, to=3-3]
	\arrow[""{name=2, anchor=center, inner sep=0}, curve={height=24pt}, from=3-1, to=3-3]
	\arrow[""{name=3, anchor=center, inner sep=0}, "{\text{Kl}(-)}"{description}, curve={height=24pt}, dashed, from=3-3, to=3-1]
	\arrow["\dashv"{anchor=center, rotate=-90}, draw=none, from=1, to=0]
	\arrow["\dashv"{anchor=center, rotate=-90}, draw=none, from=3, to=2]
\end{tikzcd}\]
Notice that the non-dotted square commutes. By taking mates, we get a  natural transformation between the left 2-adjoints:
\[\begin{tikzcd}
	\talgs && \colaxtalgl \\
	\\
	\ck && {\comnd{\ck}}
	\arrow["{U^T}"', from=1-1, to=3-1]
	\arrow[""{name=0, anchor=center, inner sep=0}, "{(-)'}"', from=1-3, to=1-1]
	\arrow[""{name=0p, anchor=center, inner sep=0}, phantom, from=1-3, to=1-1, start anchor=center, end anchor=center]
	\arrow["{\mathbb{U}^T}", from=1-3, to=3-3]
	\arrow[""{name=1, anchor=center, inner sep=0}, "{\text{Kl}(-)}", from=3-3, to=3-1]
	\arrow[""{name=1p, anchor=center, inner sep=0}, phantom, from=3-3, to=3-1, start anchor=center, end anchor=center]
	\arrow["\chi"', shorten <=9pt, shorten >=9pt, Rightarrow, from=1p, to=0p]
\end{tikzcd}\]
For each colax $T$-algebra $\mathbb{A} := (A,a,\gamma, \iota)$, this gives a canonical 1-cell $A_{ai_A} \to A'$ in $\ck$ between the Kleisli object of the underlying comonad of $\mathbb{A}$ and the underlying object of the strictification of $\mathbb{A}$. Explicitly, the 2-monad unit induces a \textit{morphism of colax coherence data}, i.e. a 2-natural transformation:
\[
\res{A,ai_A, \gamma i_{TA} i_A, \iota} \Rightarrow \res{A,a,\gamma,\iota}: \Delta^{op}_2 \to \ck.
\]
The components of the 2-natural transformation are as follows:
\[\begin{tikzcd}
	A && A && A \\
	\\
	{T^3A} && {T^2A} && TA
	\arrow["{1_A}"{description}, no head, from=1-1, to=1-3]
	\arrow["{ai_A}"{description}, shift right=4, from=1-1, to=1-3]
	\arrow["{1_A}"{description}, shift left=4, no head, from=1-1, to=1-3]
	\arrow["{i_{T^2A} i_{TA} i_A}"{description}, from=1-1, to=3-1]
	\arrow["{ai_A}"{description}, shift right=4, from=1-3, to=1-5]
	\arrow["{1_A}"{description}, shift left=4, no head, from=1-3, to=1-5]
	\arrow["{i_{TA} i_A}"{description}, from=1-3, to=3-3]
	\arrow["{1_A}"{description}, no head, from=1-5, to=1-3]
	\arrow["{i_A}"{description}, from=1-5, to=3-5]
	\arrow["{m_{T^2A}}"{description}, shift left=4, from=3-1, to=3-3]
	\arrow["{Tm_{TA}}"{description}, from=3-1, to=3-3]
	\arrow["{T^2a}"{description}, shift right=4, from=3-1, to=3-3]
	\arrow["Ta"{description}, shift right=4, from=3-3, to=3-5]
	\arrow["{m_A}"{description}, shift left=4, from=3-3, to=3-5]
	\arrow["{Ti_A}"{description}, from=3-5, to=3-3]
\end{tikzcd}\]
This implies that if $(f: TA \to A', \xi: f m_A \Rightarrow f Ta)$ is the codescent object of $\res{A,a,\gamma,\iota}$ and $(g: A \to A_t, \psi: g \Rightarrow g ai_A)$ is the Kleisli object of the underlying comonad, there exists a unique comparison map $\theta: A_t \to A'$ such that $\theta g = fi_A$ and such that the following 2-cells are equal:
\[\begin{tikzcd}[column sep=scriptsize]
	&&& TA \\
	A & TA & {T^2A} && {A'} & A && {A_t} & {A'} \\
	&& A & TA &&& A & TA
	\arrow["f", from=1-4, to=2-5]
	\arrow["\xi", shorten <=6pt, shorten >=6pt, Rightarrow, from=1-4, to=3-4]
	\arrow["{i_A}", from=2-1, to=2-2]
	\arrow[curve={height=-12pt}, Rightarrow, no head, from=2-2, to=1-4]
	\arrow["{i_{TA}}", from=2-2, to=2-3]
	\arrow["a"', from=2-2, to=3-3]
	\arrow["{m_A}", from=2-3, to=1-4]
	\arrow["Ta"{description}, from=2-3, to=3-4]
	\arrow[""{name=0, anchor=center, inner sep=0}, "g", from=2-6, to=2-8]
	\arrow["{fi_A}", curve={height=-24pt}, from=2-6, to=2-9]
	\arrow["{ai_A}"', from=2-6, to=3-7]
	\arrow["\theta"', from=2-8, to=2-9]
	\arrow["{i_A}"', from=3-3, to=3-4]
	\arrow["f"', from=3-4, to=2-5]
	\arrow["g"', from=3-7, to=2-8]
	\arrow["{i_A}"', from=3-7, to=3-8]
	\arrow["f"', from=3-8, to=2-9]
	\arrow["\psi", shorten <=3pt, Rightarrow, from=0, to=3-7]
\end{tikzcd}\]
\end{pozn}

\begin{pr}[Coherence for colax functors]\label{PR_reslan_laxcoh}
Let $\cj$ be a small 2-category and consider the reslan 2-monad $T$ on $[\ob \cj, \cat]$ from Example \ref{PR_reslan2monad}. This 2-monad preserves all colimits, so since $[ \ob \cj, \cat]$ is 2-cocomplete, so is $\talgs = [\cj, \cat]$. By Theorem \ref{THM_adjunkce_talgs_talgl_colaxtalgl}, the inclusions into lax morphisms and colax algebras both admit 2-adjoints. Again, because the 2-monad preserves all colimits, by Theorem \ref{THM_lax_coherence_colax_algebras_EXPOSITION}, coherence for colax algebras holds for this 2-monad.
\end{pr}

\begin{prc}[Coherence for lax morphisms]

We mentioned in Remark \ref{POZN_reflexive_coherence_data} that the codescent object of a resolution of a strict $T$-algebra is a sifted colimit. 

Coherence for lax morphisms is thus true for any 2-monad $T$ on a 2-cocomplete 2-category $\ck$ that preserves sifted colimits -- the preservation of resolutions of strict algebras by $T$ and $T^2$ guarantees that those 2-colimits lift to $\talgs$ and so the strictification 2-functors exist by Theorem \ref{THM_adjunkce_talgs_talgl_colaxtalgl}. The rest now follows from Theorem \ref{THM_coh_for_lax_mors}.

For a concrete example, take the monoidal category 2-monads on $\cat$ from Examples \ref{PR_freestrictmoncat2monad} and \ref{PR_freesymstrictmoncat2monad} -- as observed in \cite[Examples 4.3.7]{internalalg}, they preserve sifted colimits.

A whole class of examples satisfying coherence for lax morphisms is the one consisting of \textit{strongly finitary 2-monads} on $\cat$ (see \cite[Definition 8.14]{johnphd}) -- they preserve sifted colimits by \cite[Corollary 8.21]{johnphd}.

\end{prc}

%% file: PhD_CH_3_FS_DBL.tex
\chapter{Factorization systems and double categories}\label{CH_FS_DBLCATS}

\section{Introduction}\label{SEC_intro_dbl_cats_fs}

In this chapter, we will study the relationship between factorization systems and double categories. A double category consists of objects, vertical morphisms, horizontal morphisms and squares just like the one pictured below:
\[\begin{tikzcd}
	a && b \\
	\\
	c && d
	\arrow[""{name=0, anchor=center, inner sep=0}, "g", from=1-1, to=1-3]
	\arrow[""{name=0p, anchor=center, inner sep=0}, phantom, from=1-1, to=1-3, start anchor=center, end anchor=center]
	\arrow["u"', from=1-1, to=3-1]
	\arrow["v", from=1-3, to=3-3]
	\arrow[""{name=1, anchor=center, inner sep=0}, "h"', from=3-1, to=3-3]
	\arrow[""{name=1p, anchor=center, inner sep=0}, phantom, from=3-1, to=3-3, start anchor=center, end anchor=center]
	\arrow["\alpha", shorten <=9pt, shorten >=9pt, Rightarrow, from=0p, to=1p]
\end{tikzcd}\]
In a general double category you cannot compose vertical morphisms with horizontal ones, but if you could, you might interpret the above square $\alpha$ as telling us that the morphism $v \circ g$ (horizontal followed by vertical) can be factored as $h \circ u$ (vertical followed by horizontal) -- this is reminiscent of ordinary factorization systems on a category.

Taking this philosophy to heart, we assign to a double category $X$ a certain \textit{category of corners} $\cnr{X}$ (a concept introduced by Mark Weber in \cite{internalalg}), in which composition of vertical and horizontal morphisms is possible, and for which squares in $X$ turn into commutative squares in $\cnr{X}$. Regarding the double category as a diagram $X: \Delta^{op}_2 \to \cat$, producing $\cnr{X}$ amounts to taking the codescent object of $X$ (Definition \ref{DEFI_codescent_object}) -- we will focus on this aspect in Section \ref{SEC_codobjs_doublecats} of the next chapter.

\medskip

\noindent \textbf{The chapter is organized as follows}:

\begin{itemize}
\item In Section \ref{SEC_dbl_cats_corners} we introduce a slight generalization of crossed double categories of \cite{internalalg}. We also describe the category of corners construction for this class of double categories and mention some examples.

\item In Section \ref{SEC_fs_dbl_cats} we establish two equivalences: the first is the equivalence between strict factorization systems and double categories for which every top right corner can be uniquely filled into a square:
\[\begin{tikzcd}
	a & b \\
	& c
	\arrow["g", from=1-1, to=1-2]
	\arrow["u", from=1-2, to=2-2]
\end{tikzcd}\]
The second is the equivalence between orthogonal factorization systems and a special kind of a crossed double category.
\end{itemize}

\section{Double categories and corners}\label{SEC_dbl_cats_corners}

\subsection{Crossed double categories}\label{SECS_pscr}

Crossed double categories (a generalization of \textit{crossed simplicial groups} of \cite{crossedsimpgps}) were introduced by Mark Weber in \cite{internalalg} to calculate various internal algebra classifiers. For instance, if $S$ is the free symmetric strict monoidal category 2-monad on $\cat$, the bar construction (also called a \textit{resolution}, see Definition \ref{PR_internal_category_coherence_data}) $\res{*}$ of a terminal $S$-algebra $*$ has the structure of a crossed double category.

In that paper, any crossed double category can be turned into a category ``in the best possible way” -- this is called the \textit{category of corners} construction. In the case of $\res{*}$ the category of corners construction produces the \textit{free symmetric strict monoidal category containing a commutative monoid}, which happens to be the category $\text{FinSet}$ of finite ordinals and all functions between them.

In this chapter we consider a slight generalization of crossed double categories, obtained by dropping the ``splitness” assumption on the opfibration that appears in the definition. This allows us to consider a bigger class of examples -- the ones for which there is no canonical choice of ``opcartesian lifts”. We then present an analogue of the category of corners construction for this wider class of double categories and prove some of its key properties. All of this is in preparation for Section \ref{SEC_fs_dbl_cats} where we show that under some conditions the category of corners admits a strict or an orthogonal factorization system.

\subsubsection{Definition and examples}

\begin{defi}\label{DEFI_crossed_dblcat}
A double category $X$ is said to be \textit{crossed} if $d_0 : X_1 \to X_0$ is an opfibration and $d_1: X_2 \to X_1, s : X_0 \to X_1$ are morphisms of opfibrations\footnote{Note that the map $d_0^2 = d_0 \circ d_0 : X_2 \to X_0$ is an opfibration since it is a composite of an opfibration and a pullback of an opfibration.}:
\[\begin{tikzcd}
	{X_0} && {X_1} && {X_2} \\
	\\
	&& {X_0}
	\arrow["{s_0}", from=1-1, to=1-3]
	\arrow["{d_1}"', from=1-5, to=1-3]
	\arrow["{d_0}", from=1-3, to=3-3]
	\arrow["{d_0^2}", from=1-5, to=3-3]
	\arrow[Rightarrow, no head, from=1-1, to=3-3]
\end{tikzcd}\]
\end{defi}

\noindent In elementary terms, this is to say the following:

A square $\kappa$ (as below left) is said to be \textit{opcartesian} (by which we mean $d_0$-opcartesian if regarded as a morphism in $X_1$) when given any square $\alpha$ (below right):
\[\begin{tikzcd}
	a && b && a && b \\
	{\widehat{b}} && c && d && e
	\arrow[""{name=0, anchor=center, inner sep=0}, "g", from=1-1, to=1-3]
	\arrow["{\widehat{u}}"', from=1-1, to=2-1]
	\arrow["u", from=1-3, to=2-3]
	\arrow[""{name=1, anchor=center, inner sep=0}, "g", from=1-5, to=1-7]
	\arrow[""{name=1p, anchor=center, inner sep=0}, phantom, from=1-5, to=1-7, start anchor=center, end anchor=center]
	\arrow["w"', from=1-5, to=2-5]
	\arrow["x", from=1-7, to=2-7]
	\arrow[""{name=2, anchor=center, inner sep=0}, "{\widehat{g}}"', from=2-1, to=2-3]
	\arrow[""{name=3, anchor=center, inner sep=0}, "h"', from=2-5, to=2-7]
	\arrow[""{name=3p, anchor=center, inner sep=0}, phantom, from=2-5, to=2-7, start anchor=center, end anchor=center]
	\arrow["\kappa", shorten <=4pt, shorten >=4pt, Rightarrow, from=0, to=2]
	\arrow["\alpha", shorten <=4pt, shorten >=4pt, Rightarrow, from=1p, to=3p]
\end{tikzcd}\]
and a factorization of $x$ as $vu: b \to c \to e$, there exists a unique square $\beta$ so that the following equality of squares holds:
\begin{equation}\label{eqopcartesiansquare}
\begin{tikzcd}
	a && b && a && b \\
	{\widehat{b}} && c & {=} &&& c \\
	d && e && d && e
	\arrow[""{name=0, anchor=center, inner sep=0}, "g", from=1-1, to=1-3]
	\arrow["{\widehat{u}}"', from=1-1, to=2-1]
	\arrow["w"', curve={height=24pt}, from=1-1, to=3-1]
	\arrow["u", from=1-3, to=2-3]
	\arrow[""{name=1, anchor=center, inner sep=0}, "g", from=1-5, to=1-7]
	\arrow["w"', from=1-5, to=3-5]
	\arrow["u", from=1-7, to=2-7]
	\arrow[""{name=2, anchor=center, inner sep=0}, "{\widehat{g}}"{description}, from=2-1, to=2-3]
	\arrow["\theta"', from=2-1, to=3-1]
	\arrow["v", from=2-3, to=3-3]
	\arrow["v", from=2-7, to=3-7]
	\arrow[""{name=3, anchor=center, inner sep=0}, "h"', from=3-1, to=3-3]
	\arrow[""{name=4, anchor=center, inner sep=0}, "h"', from=3-5, to=3-7]
	\arrow["\kappa", shorten <=4pt, shorten >=4pt, Rightarrow, from=0, to=2]
	\arrow["\alpha", shorten <=9pt, shorten >=9pt, Rightarrow, from=1, to=4]
	\arrow["{\exists ! \beta}", shorten <=4pt, shorten >=4pt, Rightarrow, from=2, to=3]
\end{tikzcd}
\end{equation}

\noindent To say that $d_0: X_1 \to X_0$ is an opfibration is to say that any tuple $(g,f)$ of a ``composable” pair of a  horizontal and a vertical morphism (as pictured below) can be filled to an opcartesian square:
\begin{equation}\label{eqd0opfib}
\begin{tikzcd}
	a && b \\
	\\
	&& c
	\arrow["g", from=1-1, to=1-3]
	\arrow["f", from=1-3, to=3-3]
\end{tikzcd}
\end{equation}
Such tuples will be referred to as (top-right) \textit{corners}.

Finally, to say that $d_1: X_2 \to X_1$ and $s:X_0 \to X_1$ are morphisms of opfibrations is to say that the opcartesian squares are closed under horizontal composition and that every horizontal identity square is opcartesian.

We will denote by $\pscr$ the full subcategory of $\dbl$ spanned by crossed double categories.

\begin{pozn}[Split version]
Crossed double categories studied by Mark Weber (\cite{internalalg}) are defined as in Definition \ref{DEFI_crossed_dblcat}, except it is required that $d_0: X_1 \to X_0$ is a \textbf{split} opfibration and the maps $d_1, s$ are (strict) morphisms of split opfibrations.

This amounts to, for every top-right corner $(g,f)$, having \textbf{a choice} of an opcartesian square $\kappa_{g,f}$ filling the corner, requiring that identity squares are chosen opcartesian, and moreover vertical and horizontal composition of chosen opcartesian squares is chosen opcartesian. In this thesis we will call them \textit{split crossed} to emphasize the presence of chosen filler squares.
\end{pozn}

\begin{pozn}[Dual version]\label{POZN_pscocr}
There is a dual version of a crossed double category that we will call \textit{co-crossed}; it is obtained by replacing ``opfibration” by ``fibration” everywhere. Note then that the double category $X$ is co-crossed if and only if $X^{v}$ is crossed.
\end{pozn}

\begin{pozn}
The requirement that every corner can be filled into an opcartesian square is equivalent to saying that every corner can be filled into a pre-opcartesian square and vertical composition of pre-opcartesian squares is pre-opcartesian. By a \textit{pre-opcartesian} square we mean a square satisfying \eqref{eqopcartesiansquare} only for squares for which $v$ is the identity. This equivalence is proven in \cite[Proposition 8.1.7]{handbook2} for a general fibration.
\end{pozn}

\begin{pr}\label{PR_sqcc_cocrossed}
If $\cc$ is a category with pullbacks, the double category $\sq{\cc}$ is co-crossed: a square is cartesian (with respect to the codomain functor $d_0: \cc^2 \to \cc$) if and only if it is a pullback square in $\cc$. Clearly, vertical and horizontal composition of pullback squares yields a pullback square, and identity squares are pullbacks.

For similar reasons, the following double categories are all  co-crossed. In each of these, every square is cartesian:
\begin{itemize}
\item $\text{PbSq}(\cc)$, $\text{MPbSq}(\cc)$ (Example \ref{PR_PbSq_cc}),
\item $\text{BOFib}(\ce)$ (Example \ref{PR_BOFIB}).
\end{itemize}
\end{pr}

\begin{pr}\label{pr_awfs_dblcat}
Let $(L, R)$ be an \textit{algebraic weak factorization system} (\cite[2.2]{awfs}) on a category $\cc$ with pullbacks. There is an associated double category $R\text{-}\mathbb{A}lg$ of $R$-algebras. It can be shown that the codomain functor of this double category is a fibration (\cite[Proposition 8]{awfs}) and moreover, $d_1, s$ are morphisms of fibrations. Thus $R\text{-}\mathbb{A}lg$ is a co-crossed double category.
\end{pr}

Of note is that it is not these examples that will play a  role later on, but rather their vertical duals (that are crossed).

\begin{pr}\label{pr_crosseddblcats}
Let now $S$ be the free symmetric strict monoidal category 2-monad on $\cat$ (Example \ref{PR_freesymstrictmoncat2monad}). Denote by $*$ the terminal $S$-algebra. Since $S$ is cartesian, the resolution $\res{*}$ is an internal category in $\cat$ and so is a double category. Moreover, it is (split) crossed by \cite[Example 4.4.5]{internalalg}. This double category has finite ordinals as objects, order-preserving maps as horizontal morphisms, permutations as vertical morphisms and squares being commutative squares in $\set$. Analogous results hold if we instead consider the free braided strict monoidal category 2-monad on $\cat$.
\end{pr}
\color{black}

\noindent A special class of the crossed double category will be of interest to us:

\subsubsection{Codomain-discrete double categories}

\begin{defi}\label{deficoddiscr}
A double category $X$ will be called \textit{codomain-discrete} if every top-right corner can be uniquely filled into a square: 
\[\begin{tikzcd}
	a && b \\
	\\
	\bullet && c
	\arrow[""{name=0, anchor=center, inner sep=0}, "g", from=1-1, to=1-3]
	\arrow["u", from=1-3, to=3-3]
	\arrow[dashed, from=1-1, to=3-1]
	\arrow[""{name=1, anchor=center, inner sep=0}, dashed, from=3-1, to=3-3]
	\arrow["{\exists !}", shorten <=9pt, shorten >=9pt, Rightarrow, from=0, to=1]
\end{tikzcd}\]
\end{defi}

This property amounts to the codomain functor $d_0: X_1 \to X_0$ being a discrete opfibration. In that case, $d_1, s$ are automatically morphisms of opfibrations and thus any codomain-discrete double category is crossed.

We denote by $\coddiscr \subseteq \dbl$ the full subcategory spanned by codomain-discrete double categories.

\begin{pozn}
Codomain-discrete double categories first appeared in  \cite[2.3]{crossedsimpgps} as double categories satisfying the \textit{star condition}.
\end{pozn}

Of note is the fact that every codomain-discrete double category is flat but not necessarily invariant as the following example demonstrates:

\begin{pr}\label{pr_atimesbdblcat}
Let $\ca, \cb$ be categories. There is a double category $X_{\ca,\cb}$ such that:
\begin{itemize}
\item Objects are the objects of $\ca \times \cb$,
\item vertical morphisms are morphisms in $\ca \times \cb$ of form $(f,1_b)$.
\item horizontal morphisms are morphisms in $\ca \times \cb$ of form $(1_a,g)$,
\item a square is a commutative square in $\ca \times \cb$.
\end{itemize}
This double category is clearly codomain-discrete and thus flat. It is not invariant because for example if we have $\theta, \psi$ distinct isomorphisms in $\ca$, the following can not be filled into a square:
\[\begin{tikzcd}
	{(a,b)} && {(a',b)} \\
	\\
	{(a,b')} && {(a',b')}
	\arrow["{(\theta,1_b)}", from=1-1, to=1-3]
	\arrow["{(1_{a'},f)}", from=1-3, to=3-3]
	\arrow["{(\psi,1_{b'})}"', from=3-1, to=3-3]
\end{tikzcd}\]
\end{pr}

\begin{pr}
Given a 2-monad $T$ on $\cat$ of form $\cat(T')$ for $T'$ a cartesian monad on $\set$, the transpose of the resolution $\res{A,a}$ of a strict $T$-algebra is codomain-discrete. We will encounter this class of examples in Chapter \ref{CH_COMPUTATION}.
\end{pr}

\subsection{The category of corners}\label{subsekce_catofcnrs}

In \cite{internalalg}, given a crossed double category $X$, the category of corners $\cnr{X}$ is constructed in two steps: first, a 2-category of corners $\cb$ is constructed, and then the category of corners is obtained by taking connected components in each hom category of $\cb$, i.e. $\cnr{X} = (\pi_0)_* \cb$. In order to avoid very long proofs, we define $\cnr{X}$ for our notion of a crossed double category $X$ straight away without ever introducing $\cb$ (which in the absence of the splitness assumption would only be a bicategory). The universal property of this construction will be studied in Subsection \ref{SEC_codobjs_doublecats} of the next chapter.

\subsubsection{Definition and examples}

\begin{defi}\label{defi2cell}
Let $X$ be a double category and assume we are given two bottom-left corners $(e,m), (e',m')$ for which the domains of $e,e'$ and codomains of $m,m'$ agree (as pictured below). A \textit{2-cell} $\beta$ between them (that we denote by $\beta: (e,m) \Rightarrow (e',m')$) is a square as below for which $d_1(\beta) \circ e = e'$:
\[\begin{tikzcd}
	a \\
	{a'} && b \\
	{a''} && b
	\arrow["e"', from=1-1, to=2-1]
	\arrow["\theta"', from=2-1, to=3-1]
	\arrow[Rightarrow, no head, from=2-3, to=3-3]
	\arrow["{e'}"', curve={height=18pt}, from=1-1, to=3-1]
	\arrow[""{name=0, anchor=center, inner sep=0}, "{m'}"', from=3-1, to=3-3]
	\arrow[""{name=1, anchor=center, inner sep=0}, "m", from=2-1, to=2-3]
	\arrow["\beta", shorten <=4pt, shorten >=4pt, Rightarrow, from=1, to=0]
\end{tikzcd}\]
\end{defi}

\begin{con}\label{concnrpbsqs}
Let $X$ be a crossed double category. Define the  \textit{category of corners} $\cnr{X}$ as follows. Its objects are the objects of $X$, while a morphism $a \to b$ is an equivalence class of corners, denoted $[u, g]:a \to b$:
\[\begin{tikzcd}
	a \\
	{a'} & b
	\arrow["u"', from=1-1, to=2-1]
	\arrow["g"', from=2-1, to=2-2]
\end{tikzcd}\]
Here two corners $(u,g), (v,h)$ are \textit{equivalent} if and only if there exists a zigzag of 2-cells between them:
\[\begin{tikzcd}
	& {(u_1,g_1)} && {(u_n,g_n)} \\
	{(u,g)} && \dots && {(v,h)}
	\arrow["{\alpha_1}"', Rightarrow, from=1-2, to=2-1]
	\arrow["{\beta_1}", Rightarrow, from=1-2, to=2-3]
	\arrow["{\alpha_n}"', Rightarrow, from=1-4, to=2-3]
	\arrow["{\beta_n}", Rightarrow, from=1-4, to=2-5]
\end{tikzcd}\]

The identity on $a \in X$ is the equivalence class $[1_a, 1_a]$, while the composite $[v,h] \circ [u,g]$ is defined to be the equivalence class of corners obtained by filling the middle corner \textbf{with a choice} of an opcartesian square:
\begin{equation}\label{EQ_composite_of_corners}
\begin{tikzcd}
	a \\
	{a'} && b \\
	{\widehat{b}} && {b'} && c
	\arrow["u"', from=1-1, to=2-1]
	\arrow["{\widehat{v}\circ u}"', curve={height=24pt}, from=1-1, to=3-1]
	\arrow[""{name=0, anchor=center, inner sep=0}, "g", from=2-1, to=2-3]
	\arrow[""{name=0p, anchor=center, inner sep=0}, phantom, from=2-1, to=2-3, start anchor=center, end anchor=center]
	\arrow["{\widehat{v}}"', from=2-1, to=3-1]
	\arrow["v", from=2-3, to=3-3]
	\arrow[""{name=1, anchor=center, inner sep=0}, "{\widehat{g}}"{description}, from=3-1, to=3-3]
	\arrow[""{name=1p, anchor=center, inner sep=0}, phantom, from=3-1, to=3-3, start anchor=center, end anchor=center]
	\arrow["{h\circ \widehat{g}}"', curve={height=24pt}, from=3-1, to=3-5]
	\arrow["h", from=3-3, to=3-5]
	\arrow["\kappa", shorten <=4pt, shorten >=4pt, Rightarrow, from=0p, to=1p]
\end{tikzcd}
\end{equation}
The composite is well-defined on the equivalence classes. To see this, consider two 2-cells $\alpha: (u,g) \Rightarrow (u',g')$, $\beta: (v,h) \Rightarrow (v',h')$ as pictured below:
\[\begin{tikzcd}
	a & a &&& b & b \\
	& {a'} && b && {b'} && c \\
	{a''} & {a''} && b & {b''} & {b''} && c
	\arrow[Rightarrow, no head, from=1-1, to=1-2]
	\arrow["{u'}"', from=1-1, to=3-1]
	\arrow["u", from=1-2, to=2-2]
	\arrow[Rightarrow, no head, from=1-5, to=1-6]
	\arrow["{v'}"', from=1-5, to=3-5]
	\arrow["v", from=1-6, to=2-6]
	\arrow[""{name=0, anchor=center, inner sep=0}, "g", from=2-2, to=2-4]
	\arrow[""{name=0p, anchor=center, inner sep=0}, phantom, from=2-2, to=2-4, start anchor=center, end anchor=center]
	\arrow["\theta"', from=2-2, to=3-2]
	\arrow[Rightarrow, no head, from=2-4, to=3-4]
	\arrow[""{name=1, anchor=center, inner sep=0}, "h", from=2-6, to=2-8]
	\arrow[""{name=1p, anchor=center, inner sep=0}, phantom, from=2-6, to=2-8, start anchor=center, end anchor=center]
	\arrow["{\theta'}"', from=2-6, to=3-6]
	\arrow[Rightarrow, no head, from=2-8, to=3-8]
	\arrow[Rightarrow, no head, from=3-1, to=3-2]
	\arrow[""{name=2, anchor=center, inner sep=0}, "{g'}"', from=3-2, to=3-4]
	\arrow[""{name=2p, anchor=center, inner sep=0}, phantom, from=3-2, to=3-4, start anchor=center, end anchor=center]
	\arrow[Rightarrow, no head, from=3-5, to=3-6]
	\arrow[""{name=3, anchor=center, inner sep=0}, "{h'}"', from=3-6, to=3-8]
	\arrow[""{name=3p, anchor=center, inner sep=0}, phantom, from=3-6, to=3-8, start anchor=center, end anchor=center]
	\arrow["\alpha", shorten <=4pt, shorten >=4pt, Rightarrow, from=0p, to=2p]
	\arrow["\beta", shorten <=4pt, shorten >=4pt, Rightarrow, from=1p, to=3p]
\end{tikzcd}\]
and consider two composites:
\[\begin{tikzcd}
	a &&&& a \\
	{a'} && b && {a'} && b \\
	{\widehat{b}} && {b'} & c & {\widetilde{b}} && {b''} & c
	\arrow["u", from=1-1, to=2-1]
	\arrow["{u'}", from=1-5, to=2-5]
	\arrow[""{name=0, anchor=center, inner sep=0}, "g", from=2-1, to=2-3]
	\arrow[""{name=0p, anchor=center, inner sep=0}, phantom, from=2-1, to=2-3, start anchor=center, end anchor=center]
	\arrow["{\widehat{v}}"', from=2-1, to=3-1]
	\arrow["v", from=2-3, to=3-3]
	\arrow[""{name=1, anchor=center, inner sep=0}, "{g'}", from=2-5, to=2-7]
	\arrow[""{name=1p, anchor=center, inner sep=0}, phantom, from=2-5, to=2-7, start anchor=center, end anchor=center]
	\arrow["{\widehat{v'}}"', from=2-5, to=3-5]
	\arrow["{v'}", from=2-7, to=3-7]
	\arrow[""{name=2, anchor=center, inner sep=0}, "{\widehat{g}}"', from=3-1, to=3-3]
	\arrow[""{name=2p, anchor=center, inner sep=0}, phantom, from=3-1, to=3-3, start anchor=center, end anchor=center]
	\arrow["h", from=3-3, to=3-4]
	\arrow[""{name=3, anchor=center, inner sep=0}, "{\widehat{g'}}"', from=3-5, to=3-7]
	\arrow[""{name=3p, anchor=center, inner sep=0}, phantom, from=3-5, to=3-7, start anchor=center, end anchor=center]
	\arrow["{h'}", from=3-7, to=3-8]
	\arrow["\kappa", shorten <=4pt, shorten >=4pt, Rightarrow, from=0p, to=2p]
	\arrow["{\kappa'}", shorten <=4pt, shorten >=4pt, Rightarrow, from=1p, to=3p]
\end{tikzcd}\]
Since $\kappa$ is opcartesian, there exists a unique square $\tau$ satisfying the equation below:
\[\begin{tikzcd}
	{a'} && b && {a'} && b \\
	{\widehat{b}} && {b'} & {=} & {a''} && {b'} \\
	{\widetilde{b}} && {b''} && {\widetilde{b}} && {b''}
	\arrow[""{name=0, anchor=center, inner sep=0}, "g", from=1-1, to=1-3]
	\arrow[""{name=0p, anchor=center, inner sep=0}, phantom, from=1-1, to=1-3, start anchor=center, end anchor=center]
	\arrow["{\widehat{v}}"', from=1-1, to=2-1]
	\arrow["v", from=1-3, to=2-3]
	\arrow[""{name=1, anchor=center, inner sep=0}, "g", from=1-5, to=1-7]
	\arrow[""{name=1p, anchor=center, inner sep=0}, phantom, from=1-5, to=1-7, start anchor=center, end anchor=center]
	\arrow["\theta"', from=1-5, to=2-5]
	\arrow[Rightarrow, no head, from=1-7, to=2-7]
	\arrow[""{name=2, anchor=center, inner sep=0}, "{\widehat{g}}"{description}, from=2-1, to=2-3]
	\arrow[""{name=2p, anchor=center, inner sep=0}, phantom, from=2-1, to=2-3, start anchor=center, end anchor=center]
	\arrow[""{name=2p, anchor=center, inner sep=0}, phantom, from=2-1, to=2-3, start anchor=center, end anchor=center]
	\arrow["\sigma"', from=2-1, to=3-1]
	\arrow["{\theta'}", from=2-3, to=3-3]
	\arrow[""{name=3, anchor=center, inner sep=0}, "{g'}"{description}, from=2-5, to=2-7]
	\arrow[""{name=3p, anchor=center, inner sep=0}, phantom, from=2-5, to=2-7, start anchor=center, end anchor=center]
	\arrow[""{name=3p, anchor=center, inner sep=0}, phantom, from=2-5, to=2-7, start anchor=center, end anchor=center]
	\arrow["{\widehat{v'}}"', from=2-5, to=3-5]
	\arrow["{v'}", from=2-7, to=3-7]
	\arrow[""{name=4, anchor=center, inner sep=0}, "{\widehat{g'}}"', from=3-1, to=3-3]
	\arrow[""{name=4p, anchor=center, inner sep=0}, phantom, from=3-1, to=3-3, start anchor=center, end anchor=center]
	\arrow[""{name=5, anchor=center, inner sep=0}, "{\widehat{g'}}"', from=3-5, to=3-7]
	\arrow[""{name=5p, anchor=center, inner sep=0}, phantom, from=3-5, to=3-7, start anchor=center, end anchor=center]
	\arrow["\kappa", shorten <=4pt, shorten >=4pt, Rightarrow, from=0p, to=2p]
	\arrow["\alpha", shorten <=4pt, shorten >=4pt, Rightarrow, from=1p, to=3p]
	\arrow["\tau", shorten <=4pt, shorten >=4pt, Rightarrow, from=2p, to=4p]
	\arrow["{\kappa'}", shorten <=4pt, shorten >=4pt, Rightarrow, from=3p, to=5p]
\end{tikzcd}\]

The following now provides a 2-cell between $(v,h) \circ (u,g)$ and $(v',h') \circ (u',g')$:
\[\begin{tikzcd}
	a & a \\
	& {a'} \\
	{a''} & {\widehat{b}} && {b'} && c \\
	{\widetilde{b}} & {\widetilde{b}} && {b''} && c
	\arrow[Rightarrow, no head, from=1-1, to=1-2]
	\arrow["{u'}"', from=1-1, to=3-1]
	\arrow["u", from=1-2, to=2-2]
	\arrow["{\widehat{v}}", from=2-2, to=3-2]
	\arrow["{\widehat{v'}}"', from=3-1, to=4-1]
	\arrow[""{name=0, anchor=center, inner sep=0}, "{\widehat{g}}", from=3-2, to=3-4]
	\arrow["\sigma"', from=3-2, to=4-2]
	\arrow[""{name=1, anchor=center, inner sep=0}, "h", from=3-4, to=3-6]
	\arrow["{\theta'}", from=3-4, to=4-4]
	\arrow[Rightarrow, no head, from=3-6, to=4-6]
	\arrow[Rightarrow, no head, from=4-1, to=4-2]
	\arrow[""{name=2, anchor=center, inner sep=0}, "{\widehat{g'}}"', from=4-2, to=4-4]
	\arrow[""{name=3, anchor=center, inner sep=0}, "{h'}"', from=4-4, to=4-6]
	\arrow["\tau", shorten <=4pt, shorten >=4pt, Rightarrow, from=0, to=2]
	\arrow["\beta", shorten <=4pt, shorten >=4pt, Rightarrow, from=1, to=3]
\end{tikzcd}\]

\noindent Note that the above implies that the composition of two corners is independent on the choice of the opcartesian square $\kappa$: just take the 2-cells $\alpha, \beta$ to be the identities.
\end{con}

\begin{prop}\label{lemmacnrXiscat}
$\cnr{X}$ is a category.
\end{prop}
\begin{proof}
Let $[f,g]: a \to b$ be a morphism in $\cnr{X}$. To show that $[f,g] \circ [1_a,1_a] = [f,g]$, note that the horizontal identity square on $f$ is by definition opcartesian so we might as well use it for the composite (the composite is independent of the choice) and the result follows.
\[\begin{tikzcd}
	a \\
	a && a \\
	c && c && b
	\arrow[Rightarrow, no head, from=1-1, to=2-1]
	\arrow[""{name=0, anchor=center, inner sep=0}, Rightarrow, no head, from=2-1, to=2-3]
	\arrow["f", from=2-3, to=3-3]
	\arrow["g", from=3-3, to=3-5]
	\arrow["f"', from=2-1, to=3-1]
	\arrow[""{name=1, anchor=center, inner sep=0}, Rightarrow, no head, from=3-1, to=3-3]
	\arrow["f"', curve={height=24pt}, from=1-1, to=3-1]
	\arrow["g"', curve={height=24pt}, from=3-1, to=3-5]
	\arrow[shorten <=4pt, shorten >=4pt, Rightarrow, from=0, to=1]
\end{tikzcd}\]
Analogously $[1_b,1_b] \circ [f,g] = [f,g]$. Consider now a composable triple $[f,g]$, $[f',g']$, $[f'',g'']$ and fill it to form a single corner as depicted below:
\[\begin{tikzcd}
	a \\
	{a'} && b \\
	\bullet && {b'} && c \\
	\bullet && \bullet && {c'} && d
	\arrow["f", from=1-1, to=2-1]
	\arrow[""{name=0, anchor=center, inner sep=0}, "g", from=2-1, to=2-3]
	\arrow[""{name=0p, anchor=center, inner sep=0}, phantom, from=2-1, to=2-3, start anchor=center, end anchor=center]
	\arrow[dashed, from=2-1, to=3-1]
	\arrow["{f'}", from=2-3, to=3-3]
	\arrow[""{name=1, anchor=center, inner sep=0}, dashed, from=3-1, to=3-3]
	\arrow[""{name=1p, anchor=center, inner sep=0}, phantom, from=3-1, to=3-3, start anchor=center, end anchor=center]
	\arrow[""{name=1p, anchor=center, inner sep=0}, phantom, from=3-1, to=3-3, start anchor=center, end anchor=center]
	\arrow[dashed, from=3-1, to=4-1]
	\arrow[""{name=2, anchor=center, inner sep=0}, "{g'}", from=3-3, to=3-5]
	\arrow[""{name=2p, anchor=center, inner sep=0}, phantom, from=3-3, to=3-5, start anchor=center, end anchor=center]
	\arrow[from=3-3, to=4-3]
	\arrow["{f''}", from=3-5, to=4-5]
	\arrow[""{name=3, anchor=center, inner sep=0}, dashed, from=4-1, to=4-3]
	\arrow[""{name=3p, anchor=center, inner sep=0}, phantom, from=4-1, to=4-3, start anchor=center, end anchor=center]
	\arrow[""{name=4, anchor=center, inner sep=0}, dashed, from=4-3, to=4-5]
	\arrow[""{name=4p, anchor=center, inner sep=0}, phantom, from=4-3, to=4-5, start anchor=center, end anchor=center]
	\arrow["{g''}", from=4-5, to=4-7]
	\arrow["{\kappa_1}", shorten <=4pt, shorten >=4pt, Rightarrow, from=0p, to=1p]
	\arrow["{\kappa_3}", shorten <=4pt, shorten >=4pt, Rightarrow, from=1p, to=3p]
	\arrow["{\kappa_2}", shorten <=4pt, shorten >=4pt, Rightarrow, from=2p, to=4p]
\end{tikzcd}\]
Denote this composite corner by $[f'',g''] \circ [f',g'] \circ [f,g]$ and call it the \textit{ternary composite}. Now to define the composition $[f'',g''] \circ [f',g']$, choose the square $\kappa_2$ as above. To define $([f'',g''] \circ [f',g']) \circ [f,g]$, choose the square $\kappa_3 \circ \kappa_1$ (as a vertical composite of opcartesian squares, it is opcartesian). We see that $([f'',g''] \circ [f',g']) \circ [f,g]$ is equal to the ternary composite. By an analogous argument (and using that opcartesian squares are closed under horizontal composites), $[f'',g''] \circ ([f',g'] \circ [f,g])$ also equals this ternary composite and thus composition is associative.

\end{proof}

If $F: X \to Y$ is any double functor between crossed double categories, there is an induced functor $\cnr{F}: \cnr{X} \to \cnr{Y}$ sending:
\begin{equation}\label{mapsto_cnrF}
\begin{tikzcd}
	a &&& Fa \\
	{a'} & b & \mapsto & {Fa'} & Fb
	\arrow["u"', from=1-1, to=2-1]
	\arrow["g", from=2-1, to=2-2]
	\arrow["Fu"', from=1-4, to=2-4]
	\arrow["Fg", from=2-4, to=2-5]
\end{tikzcd}
\end{equation}
This gives us a functor $\cnr{-}: \pscr \to \cat$.

\begin{pr}\label{prspans}
Let $\cc$ be a category with pullbacks and consider the double category $\text{PbSq}(\cc)^{v}$. The category $\cnr{\text{PbSq}(\cc)^{v}}$ has as objects the objects of $\cc$, while a morphism $a \to b$ is an equivalence class of corners (usually called \textit{spans}) like this:
\[\begin{tikzcd}
	a \\
	{a'} & b
	\arrow[from=2-1, to=1-1]
	\arrow[from=2-1, to=2-2]
\end{tikzcd}\]
Note that there is a span isomorphism between the two spans $(u,g), (v,h)$ if and only if there is a 2-cell between them if we regard them as corners.
\[\begin{tikzcd}
	a &&&&& a \\
	& {a'} && \iff && {a'} && b \\
	{a''} && b &&& {a'} && b
	\arrow["u"', from=2-2, to=1-1]
	\arrow["g", from=2-2, to=3-3]
	\arrow["v", from=3-1, to=1-1]
	\arrow["{\exists \cong}"', from=3-1, to=2-2]
	\arrow["h"', from=3-1, to=3-3]
	\arrow["u"', from=2-6, to=1-6]
	\arrow[""{name=0, anchor=center, inner sep=0}, "g", from=2-6, to=2-8]
	\arrow[Rightarrow, no head, from=2-8, to=3-8]
	\arrow[""{name=1, anchor=center, inner sep=0}, "g"', from=3-6, to=3-8]
	\arrow[from=3-6, to=2-6]
	\arrow["v", curve={height=-24pt}, from=3-6, to=1-6]
	\arrow["\exists", shorten <=4pt, shorten >=4pt, Rightarrow, from=0, to=1]
\end{tikzcd}\]
The composition of corners is defined using pullbacks. In other words, we have \\$\cnr{\text{PbSq}(\cc)^{v}} = \text{Span}(\cc)$, the category of isomorphism classes of spans in $\cc$.
\end{pr}

\begin{pr}\label{prparc}
Let $\cc$ be a category with pullbacks and consider the double category $\text{MPbSq}(\cc)^{v}$. By a similar reasoning as above we obtain that the category of corners corresponding to this double category is isomorphic to the category $\text{Par}(\cc)$ of \textit{partial maps} in $\cc$, as defined in \cite[p. 246]{restrcats}.
\end{pr}

\begin{pr}\label{prcofun}
Let $\ce$ be a category with pullbacks and consider now the double category $\text{BOFib}(\ce)^{v}$. Morphisms in $\cnr{\text{BOFib}(\ce)^{v}}$ are equivalence classes of spans $(F,G)$ where $F$ is a bijection on objects and $G$ is a discrete opfibration as pictured below:
\[\begin{tikzcd}
	\ca \\
	{\ca'} && \cb
	\arrow["F", from=2-1, to=1-1]
	\arrow["G", from=2-1, to=2-3]
\end{tikzcd}\]
This category of corners is isomorphic to $\text{Cof}(\ce)$, the category of internal categories and \textit{cofunctors}, see for instance \cite[Theorem 18]{bryce_cof}.
\end{pr}

\begin{pr}
If $(L,R)$ is an algebraic weak factorization system on a category $\cc$ with pullbacks, consider the co-crossed double category $R\text{-}\mathbb{A}lg$ of $R$-algebras (Example \ref{pr_awfs_dblcat}). Its vertical dual $R\text{-}\mathbb{A}lg^{v}$ is thus crossed. Its category of corners construction now gives the \textit{category of weak maps} $\boldsymbol{Wk}_l(L,R)$ associated to the system $(L,R)$, see \cite[Section 3.4 and  Remark 13]{awfs2}.
\end{pr}

\begin{pr}
If $X$ is split crossed, our category of corners construction agrees with that of \cite[Corollary 5.4.5]{internalalg}, as is easily verified.

Recall from Example \ref{pr_crosseddblcats} the free symmetric strict monoidal category 2-monad $S$ and the crossed double category $\res{*}$.

The objects of $\cnr{\res{*}}$ are finite ordinals, while a morphism $m \to n$ is an equivalence class of corners consisting of a permutation followed by an order-preserving map:
\[\begin{tikzcd}
	m \\
	m & n
	\arrow["\phi"', from=1-1, to=2-1]
	\arrow["f"', from=2-1, to=2-2]
\end{tikzcd}\]
It can be proven that $\cnr{\res{*}} = \text{FinSet}$, the category of finite ordinals and all functions (see \cite[Theorem 6.3.1]{internalalg}).

If $B$ is the free braided strict monoidal category 2-monad, $\cnr{\res{*}} = \text{Vine}$, the category with objects being natural numbers and morphisms being \textit{vines}, that is, ``braids for which the strings can merge” (see \cite[Theorem 6.3.2]{internalalg}).
\end{pr}
\color{black}

\begin{pozn}\label{pozn_coddiscrequivrel}
If $X$ is codomain-discrete, note that there is a 2-cell $\beta$ between corners $(u,g)$, $(v,h)$ if and only if $u=v$, $g=h$ and $\beta = 1_g$ is the identity square. Thus the category $\cnr{X}$ has corners as morphisms with no equivalence relation involved. This has also been observed in \cite[Corollary 5.4.7]{internalalg}.
\end{pozn}

\begin{prc}
For instance, to any 2-category $\ck$ with a distinguished object $C$ one can associate a double category $K$ such that:
\begin{itemize}
\item the objects are 1-cells $f: A \to C$ in $\ck$,
\item a horizontal morphism $u: f \to g$ is a commutative triangle:
\[\begin{tikzcd}
	A && B \\
	\\
	& C
	\arrow["u", from=1-1, to=1-3]
	\arrow["f"', from=1-1, to=3-2]
	\arrow["g", from=1-3, to=3-2]
\end{tikzcd}\]
\item a vertical morphism $f \to f'$ is a 2-cell $\alpha: f \Rightarrow f'$ (and thus exists only if the domains of $f,f'$ agree),
\item Given:
\begin{align*}
f,f' &: A \to C,\\
g,g' &: B \to C,\\
u,v &: A \to B,\\
\alpha &: f \Rightarrow f',\\
\beta &: g \Rightarrow g',
\end{align*}
this data is a boundary of a square:
\[\begin{tikzcd}
	f && g \\
	\\
	{f'} && {g'}
	\arrow[""{name=0, anchor=center, inner sep=0}, "u", from=1-1, to=1-3]
	\arrow["\alpha"', from=1-1, to=3-1]
	\arrow["\beta", from=1-3, to=3-3]
	\arrow[""{name=1, anchor=center, inner sep=0}, "v"', from=3-1, to=3-3]
	\arrow["\exists", shorten <=9pt, shorten >=9pt, Rightarrow, from=0, to=1]
\end{tikzcd}\]
if and only if $u = v$ and:
\[\begin{tikzcd}
	& B &&&& B \\
	A && C & {=} & A && C
	\arrow[""{name=0, anchor=center, inner sep=0}, "g"', curve={height=24pt}, from=1-2, to=2-3]
	\arrow[""{name=1, anchor=center, inner sep=0}, "{g'}", from=1-2, to=2-3]
	\arrow["{g'}", from=1-6, to=2-7]
	\arrow["u", from=2-1, to=1-2]
	\arrow["f"', curve={height=30pt}, from=2-1, to=2-3]
	\arrow["v", from=2-5, to=1-6]
	\arrow[""{name=2, anchor=center, inner sep=0}, "f"', curve={height=24pt}, from=2-5, to=2-7]
	\arrow[""{name=3, anchor=center, inner sep=0}, "{f'}", from=2-5, to=2-7]
	\arrow["\beta"', shorten <=4pt, shorten >=4pt, Rightarrow, from=0, to=1]
	\arrow["\alpha"', shorten <=3pt, shorten >=3pt, Rightarrow, from=2, to=3]
\end{tikzcd}\]
\end{itemize}
Notice however that given just $f,u,g,\beta,g'$, it is forced that:
\begin{align*}
f' &\sstackrel{!}{=} g'u,\\
\alpha &\sstackrel{!}{=} \beta u.
\end{align*}
This means that the double category is codomain-discrete. The category $\cnr{\ck_K}$ can be seen to be the \textit{lax slice category over} $C \in \ob \ck$. We have seen its dual version in Example \ref{PR_fibration_2monad_Kleisli2cat}.
\end{prc}

\subsubsection{Some properties of $\cnr{X}$}

The following proposition captures the idea that ``a square in $X$ turns into a commutative square in $\cnr{X}$”:

\begin{prop}\label{prop_naturalityinsquares}
Let $X$ be a crossed double category. We have:
\[\begin{tikzcd}[column sep=scriptsize]
	a && b &&& a && b \\
	&&& {\text{in } X} & \Rightarrow &&&& {\text{commutes in } \cnr{X}} \\
	c && d &&& c && d
	\arrow[""{name=0, anchor=center, inner sep=0}, "m", from=1-1, to=1-3]
	\arrow["e", from=1-3, to=3-3]
	\arrow["{e'}"', from=1-1, to=3-1]
	\arrow[""{name=1, anchor=center, inner sep=0}, "{m'}"', from=3-1, to=3-3]
	\arrow["{[e',1]}"', from=1-6, to=3-6]
	\arrow["{[1,m]}", from=1-6, to=1-8]
	\arrow["{[e,1]}", from=1-8, to=3-8]
	\arrow["{[1,m']}"', from=3-6, to=3-8]
	\arrow["\exists"', shorten <=9pt, shorten >=9pt, Rightarrow, from=0, to=1]
\end{tikzcd}\]
\end{prop}
\begin{proof}
Denote $[u,g] := [e,1] \circ [1,m]$ and denote by $\kappa$ the opcartesian square we used for this composition. From opcartesianness there is a unique square $\beta$ such that:
\[\begin{tikzcd}
	a && b && a &&& b \\
	{\widehat{b}'} && d & {=} \\
	c && d && c &&& d
	\arrow[""{name=0, anchor=center, inner sep=0}, "m", from=1-1, to=1-3]
	\arrow["e", from=1-3, to=2-3]
	\arrow["u"', from=1-1, to=2-1]
	\arrow[""{name=1, anchor=center, inner sep=0}, "g"{description}, from=2-1, to=2-3]
	\arrow[""{name=2, anchor=center, inner sep=0}, "{m'}"', from=3-1, to=3-3]
	\arrow[Rightarrow, no head, from=2-3, to=3-3]
	\arrow[""{name=3, anchor=center, inner sep=0}, "m", from=1-5, to=1-8]
	\arrow["e", from=1-8, to=3-8]
	\arrow["{e'}"', from=1-5, to=3-5]
	\arrow[""{name=4, anchor=center, inner sep=0}, "{m'}"', from=3-5, to=3-8]
	\arrow["\theta"', from=2-1, to=3-1]
	\arrow["{e'}"', curve={height=24pt}, from=1-1, to=3-1]
	\arrow["\alpha"', shorten <=9pt, shorten >=9pt, Rightarrow, from=3, to=4]
	\arrow["\kappa", shorten <=4pt, shorten >=4pt, Rightarrow, from=0, to=1]
	\arrow["{\exists ! \beta}", shorten <=4pt, shorten >=4pt, Rightarrow, from=1, to=2]
\end{tikzcd}\]
This square $\beta$ now exhibits the equality $[1,m'] \circ [e',1] = [e',m'] = [u,g] = [e,1] \circ [1,m]$.

\end{proof}

The additional assumption requiring that every square is opcartesian simplifies the description of $\cnr{X}$ for a crossed double category $X$:

\begin{lemma}\label{lemmaXpscr_everysquareopcart}
Let $X$ be a crossed double category in which every square is opcartesian. Then any 2-cell $\beta: [e,m] \Rightarrow [e',m']$ between corners is vertically invertible. In particular two corners in $\cnr{X}$ are equivalent if and only if there exists a single (vertically invertible) 2-cell between them.
\end{lemma}
\begin{proof}
Consider a square $\beta$ as pictured below. Since the vertical identity square $1_g$ on the morphism $g$ is opcartesian, there exists a unique square $\gamma$ such that:
\[\begin{tikzcd}
	a && b && a && b \\
	c && b & {=} \\
	a && b && a && b
	\arrow["u"', from=1-1, to=2-1]
	\arrow[""{name=0, anchor=center, inner sep=0}, "g", from=1-1, to=1-3]
	\arrow[Rightarrow, no head, from=1-3, to=2-3]
	\arrow[""{name=1, anchor=center, inner sep=0}, "h"{description}, from=2-1, to=2-3]
	\arrow["v"', from=2-1, to=3-1]
	\arrow[Rightarrow, no head, from=2-3, to=3-3]
	\arrow[""{name=2, anchor=center, inner sep=0}, "g"', from=3-1, to=3-3]
	\arrow["g", from=1-5, to=1-7]
	\arrow[Rightarrow, no head, from=1-7, to=3-7]
	\arrow[Rightarrow, no head, from=1-5, to=3-5]
	\arrow["g"', from=3-5, to=3-7]
	\arrow[curve={height=24pt}, Rightarrow, no head, from=1-1, to=3-1]
	\arrow["{\exists ! \gamma}", shorten <=4pt, shorten >=4pt, Rightarrow, from=1, to=2]
	\arrow["\beta", shorten <=4pt, shorten >=4pt, Rightarrow, from=0, to=1]
\end{tikzcd}\]
Thus $\gamma \circ \beta = 1_g$. Post-composing with $\beta$, we get:
\[
\beta \circ \gamma \circ \beta = 1_h \circ \beta
\]
Since $\beta$ is opcartesian, we get $\beta \circ \gamma = 1_h$ as well.
\end{proof}

\begin{nota}\label{notecnrX_verticalhorizontalcnrs}
Given a crossed double category $X$, denote by $\ce_X$ the class of corners in $\cnr{X}$ of the form $[f,1_b]$ and by $\cm_X$ the class of corners of the form $[1_a, g]$. We will call these \textit{vertical} and \textit{horizontal} corners.
\end{nota}

\begin{prop}\label{THMcancprops_weakorthog}
Let $X$ be a crossed double category. Then the class $\ce_X$ has the right cancellation property. Both classes $\ce_X, \cm_X$ contain all identities and are closed under composition. We also have:
\[
\cnr{X} = \cm_X \circ \ce_X.
\]
Moreover, if every square is opcartesian in $X$, we have that $\ce_X$ is weakly orthogonal to $\cm_X$:
\[
\ce_X \boxslash \cm_X.
\]
\end{prop}
\begin{proof}
To show the cancellation property, assume $[s,t] \circ [e,1] = [u,1] \in \ce_X$. We then have a square as pictured below left:
\[\begin{tikzcd}
	a \\
	b &&&&& b \\
	{b'} && c &&& {a'} && c \\
	c && c &&& c && c
	\arrow["s", from=2-1, to=3-1]
	\arrow[""{name=0, anchor=center, inner sep=0}, "t", from=3-1, to=3-3]
	\arrow[""{name=1, anchor=center, inner sep=0}, Rightarrow, no head, from=4-1, to=4-3]
	\arrow[Rightarrow, no head, from=3-3, to=4-3]
	\arrow["e", from=1-1, to=2-1]
	\arrow["\theta"', from=3-1, to=4-1]
	\arrow["u"', curve={height=24pt}, from=1-1, to=4-1]
	\arrow["\theta"', from=3-6, to=4-6]
	\arrow[""{name=2, anchor=center, inner sep=0}, Rightarrow, no head, from=4-6, to=4-8]
	\arrow[""{name=3, anchor=center, inner sep=0}, "t", from=3-6, to=3-8]
	\arrow[Rightarrow, no head, from=3-8, to=4-8]
	\arrow["s", from=2-6, to=3-6]
	\arrow["{\theta s}"', curve={height=24pt}, from=2-6, to=4-6]
	\arrow["\beta", shorten <=4pt, shorten >=4pt, Rightarrow, from=0, to=1]
	\arrow["\beta", shorten <=4pt, shorten >=4pt, Rightarrow, from=3, to=2]
\end{tikzcd}\]
The same square $\beta$ now exhibits the equality $[s,t] = [\theta s, 1]$ (pictured above right) and thus $[s,t] \in \ce_X$. The fact that $\ce_X, \cm_X$ contain identities and are closed under composition, as well as the fact that $\cnr{X} = \cm_X \circ \ce_X$ are obvious.

Assume now that every square is opcartesian in $X$. Since $\cnr{X} = \cm_X \circ \ce_X$, to prove weak orthogonality we first show that given two factorizations $[e,m] =  [e',m']$ of the same morphism, there exists a morphism of factorizations between them, i.e.:
\begin{equation}\label{EQ_morfismus_faktorizaci}
\begin{tikzcd}
	a && {a'} && b \\
	a && {a''} && b
	\arrow["{[e,1]}", from=1-1, to=1-3]
	\arrow["{[1,m]}", from=1-3, to=1-5]
	\arrow["{[e',1]}"', from=2-1, to=2-3]
	\arrow["{[1,m']}"', from=2-3, to=2-5]
	\arrow[Rightarrow, no head, from=1-1, to=2-1]
	\arrow[Rightarrow, no head, from=1-5, to=2-5]
	\arrow[dashed, from=1-3, to=2-3]
\end{tikzcd}
\end{equation}
Since $[e,m] = [e',m']$ and thanks to Lemma \ref{lemmaXpscr_everysquareopcart}, there exists a single invertible square like this:
\[\begin{tikzcd}
	a \\
	{a'} && b \\
	{a''} && b
	\arrow["{e'}"', curve={height=24pt}, from=1-1, to=3-1]
	\arrow["e"', from=1-1, to=2-1]
	\arrow[""{name=0, anchor=center, inner sep=0}, "m", from=2-1, to=2-3]
	\arrow["\theta"', from=2-1, to=3-1]
	\arrow[""{name=1, anchor=center, inner sep=0}, "{m'}"', from=3-1, to=3-3]
	\arrow[Rightarrow, no head, from=2-3, to=3-3]
	\arrow["\beta", shorten <=4pt, shorten >=4pt, Rightarrow, from=0, to=1]
\end{tikzcd}\]
It is now easy to verify that the corner $[\theta, 1]$ makes both squares in the above diagram commute. Consider now the general lifting problem for the classes $\ce_X$, $\cm_X$ (as pictured below left). By a standard argument, the filler can be defined to be the corner $[\theta_1 v, g \theta_2]$, where the corner $[\theta_1, \theta_2]$ is obtained from the diagram below right:

\[\begin{tikzcd}
	a && b && a & y & b \\
	{a'} && {b'} && a & x & b
	\arrow["{[u,g]}", from=1-1, to=1-3]
	\arrow["{[e,1]}"', from=1-1, to=2-1]
	\arrow["{[1,m]}", from=1-3, to=2-3]
	\arrow["{[u,1]}", from=1-5, to=1-6]
	\arrow[equals, from=1-5, to=2-5]
	\arrow["{[1,mg]}", from=1-6, to=1-7]
	\arrow[equals, from=1-7, to=2-7]
	\arrow[dashed, from=2-1, to=1-3]
	\arrow["{[v,h]}"', from=2-1, to=2-3]
	\arrow["{[ve,1]}"', from=2-5, to=2-6]
	\arrow["{[\theta_1, \theta_2]}"', dashed, from=2-6, to=1-6]
	\arrow["{[1,h]}"', from=2-6, to=2-7]
\end{tikzcd}\]

\end{proof}

\noindent The classes $(\ce_X,\cm_X)$ give rise to an ordinary weak factorization system $(\widetilde{\ce_X},\widetilde{\cm_X})$ for which the first class is obtained by closing $\ce_X$ under codomain-retracts and the second is obtained from $\cm_X$ by closing it under domain retracts. This is a well-known result so we omit its proof.

\begin{pr}\label{pr_spansWFS}
Recall the example $\cnr{\pbsq{\cc}^{v}} = \spn{\cc}$. Since every square is opcartesian (a pullback), we obtain that the class $\ce_{\pbsq{\cc}^{v}}$ is weakly orthogonal to $\cm_{\pbsq{\cc}^{v}}$. Note that in this case, both classes are already closed under the required retracts. We obtain:

\begin{prop}\label{prop_wfsonspans}
The two canonical classes of morphisms in the category $\spn(\cc)$ form a weak factorization system.
\end{prop}
\end{pr}

\color{black}

\section{Factorization systems and double categories}\label{SEC_fs_dbl_cats}

In this section we will be putting additional hypotheses on the crossed double category $X$ to ensure that the classes $(\ce_X,\cm_X)$ of morphisms in $\cnr{X}$ have more desirable properties (namely, form a strict or an orthogonal factorization system). This gives us the direction:
\[
\text{double categories } \rightsquigarrow \text{ factorization systems.}
\]
For the opposite direction, we introduce a construction that sends two classes $(\ce,\cm)$ of morphisms in a category $\cc$ to a certain double category $D_{\ce,\cm}$ of commutative squares.

We begin in Subsection \ref{subsekce_sfs} by doing this to the case of strict factorization systems, showing that the mappings $X \mapsto (\ce_X, \cm_X)$ and $(\ce,\cm) \mapsto D_{\ce,\cm}$ induce an equivalence between the categories of strict factorization systems and codomain-discrete double categories.

In Subsection \ref{subsekce_ofs} we prove analogous results for the categories of orthogonal factorization systems and  \textit{factorization double categories}: a symmetric variant of crossed double categories whose bottom-left corners satisfy a certain joint monicity property.

\subsection{Strict factorization systems}\label{subsekce_sfs}

In \cite{distrlaws} it has been shown that distributive laws in $\spanc$ can equivalently be described as strict factorization systems. Given a codomain-discrete double category $X$, the category of corners $\cnr{X}$ can be constructed using a distributive law in $\spanc(\cat)$\footnote{This is the original construction of $\cnr{X}$ for a codomain-discrete double category $X$ in \cite{internalalg}} -- this gives a first hint that there is a relationship between double categories and factorization systems.

\begin{defi}
A \textit{strict factorization system} $(\ce,\cm)$ on a category $\cc$ consists of two wide sub-categories $\ce, \cm \subseteq \cc$ such that for every morphism $f \in \cc$ there exist unique $e \in \ce, m \in \cm$ such that: $f = m \circ e$.
\end{defi}

\begin{defi}\label{defi_SFS}
Denote by $\cs \cf \cs$ the category whose:
\begin{itemize}
\item objects are strict factorization systems $\ce \subseteq \cc \supseteq \cm$,
\item morphisms $(\ce \subseteq \cc \supseteq \cm) \to (\ce' \subseteq \cc' \supseteq \cm')$ are functors $F: \cc \to \cc'$ satisfying:
\begin{align*}
F(\ce) &\subseteq \ce',\\
F(\cm) &\subseteq \cm'.
\end{align*}
\end{itemize}
\end{defi}

\begin{lemma}\label{lemma_coddiscrimpliessfs}
Let $X$ be codomain-discrete. Then the classes $(\ce_X,\cm_X)$ form a strict factorization system on $\cnr{X}$. There is a functor $\coddiscr \to \sfs$ given by the assignment $X \mapsto (\ce_X,\cm_X)$.
\end{lemma}
\begin{proof}
Recall from Remark \ref{pozn_coddiscrequivrel} that for such $X$, two morphisms $[u,g], [v,h]$ in $\cnr{X}$ are equal if and only if $u = v, g = h$. From this it follows that the factorization $[u,g] = [1,g] \circ [u,1]$ is unique.

If $H: X \to Y$ is a double functor, the induced functor $\cnr{H}$ from $\cnr{X}$ to $\cnr{Y}$ (see \eqref{mapsto_cnrF}) satisfies $\cnr{H}(\ce_X) \subseteq \ce_Y$, $\cnr{H}(\cm_X) \subseteq \cm_Y$ and thus is a morphism in $\sfs$.
\end{proof}

Denote the above functor simply by $\cnr{-}: \coddiscr \to \sfs$.

\begin{pr}
Let $\ca, \cb$ be categories and $X_{\ca,\cb}$ the codomain-discrete double category from Example \ref{pr_atimesbdblcat}. The category of corners $\cnr{X_{\ca,\cb}}$ is isomorphic to just $\ca \times \cb$ and this category admits a strict factorization system $(\ce,\cm)$, where:
\begin{align*}
\ce &:= \{ (1_a, f) \mid a \in \ca, f \in \cb \},\\
\cm &:= \{ (g, 1_b) \mid g \in \ca, b \in \cb \}.
\end{align*}
\end{pr}
\color{black}

\begin{con}\label{conDcecm}
Let $(\ce,\cm)$ be two classes of morphisms in a category $\cc$, both closed under composition and containing all identities. We define a double category $D_{\ce,\cm}$ as follows:
\begin{itemize}
\item The objects are the objects of $\cc$,
\item the category of objects and vertical morphisms is $\ce$,
\item the category of objects and horizontal morphisms is $\cm$,
\item the squares are commutative squares in $\cc$.
\end{itemize}

If we have two classes $(\ce,\cm), (\ce',\cm')$ on categories $\cc, \cc'$ and $F: \cc \to \cc'$ a functor satisfying $F(\ce) \subseteq \ce'$ and $F(\cm) \subseteq \cm'$, there is an induced double functor: \[ D_F : D_{\ce,\cm} \to D_{\ce',\cm'},\] defined in an obvious way.
\end{con}

\begin{lemma}\label{lemma_sfs_to_coddiscr}
Let $(\ce,\cm)$ be a strict factorization system on a category $\cc$. Then $D_{\ce,\cm}$ is codomain-discrete. The assignment $(\ce,\cm) \mapsto D_{\ce,\cm}$ induces a functor $\sfs \to \coddiscr$.
\end{lemma}
\begin{proof}
Every morphism in $\cc$ of form $e \circ m$ can be uniquely factored as $m' \circ e'$ with $e' \in \ce$ and $m' \in \cm$. But this precisely means that $D_{\ce,\cm}$ is codomain-discrete.

In the above construction we have seen that $(\ce,\cm) \mapsto D_{\ce,\cm}$ is functorial, the rest is now obvious.
\end{proof}

\begin{theo}\label{thmekvivalence_sfs_coddiscr}
The functor $\cnr{-}: \coddiscr \to \sfs$ is an equivalence of categories, with the equivalence inverse being the functor $D: \sfs \to \coddiscr$.
\end{theo}
\begin{proof}
We will show that there are natural isomorphisms $1 \cong D \circ \cnr{-}$ and\\$\cnr{-} \circ D \cong 1$:

To see that $1 \cong D \circ \cnr{-}$, let $X$ be a codomain-discrete double category. First note that by Remark \ref{pozn_coddiscrequivrel}, the identity-on-objects functor $X_0 \to \ce_X$ sending a morphism $e \mapsto (e,1)$ is an isomorphism. Similarly we have $h(X) \cong \cm_X$ via the identity-on-objects functor $m \mapsto (1,m)$.

There is now a double functor $X \to D_{\ce_X, \cm_X}$ that is identity on objects and whose vertical morphism and horizontal morphism components are given by the functors described above. To see that it is well-defined on squares, we would need to prove the direction ``$\Rightarrow$” in the picture below:
\begin{equation}\label{eq_iffsquares}
\begin{tikzcd}[column sep=scriptsize]
	a && b &&& a && b \\
	&&& {\text{in } X} & \iff &&&& {\text{in }D_{\ce_X,\cm_X}} \\
	c && d &&& c && d
	\arrow[""{name=0, anchor=center, inner sep=0}, "m", from=1-1, to=1-3]
	\arrow["e", from=1-3, to=3-3]
	\arrow["{e'}"', from=1-1, to=3-1]
	\arrow[""{name=1, anchor=center, inner sep=0}, "{m'}"', from=3-1, to=3-3]
	\arrow["{(e',1)}"', from=1-6, to=3-6]
	\arrow[""{name=2, anchor=center, inner sep=0}, "{(m,1)}", from=1-6, to=1-8]
	\arrow["{(e,1)}", from=1-8, to=3-8]
	\arrow[""{name=3, anchor=center, inner sep=0}, "{(m',1)}"', from=3-6, to=3-8]
	\arrow["\exists"', shorten <=9pt, shorten >=9pt, Rightarrow, from=0, to=1]
	\arrow["\exists"', shorten <=9pt, shorten >=9pt, Rightarrow, from=2, to=3]
\end{tikzcd}
\end{equation}
But this direction follows from from Proposition \ref{prop_naturalityinsquares}. Since both double categories are flat, to show that this double functor is an isomorphism it suffices to show the direction ``$\Leftarrow$” in Diagram \eqref{eq_iffsquares}. Consider the square used for the composition $(e,1) \circ (m,1)$:
\[\begin{tikzcd}
	a && b \\
	\\
	\bullet && d
	\arrow[""{name=0, anchor=center, inner sep=0}, "m", from=1-1, to=1-3]
	\arrow["e", from=1-3, to=3-3]
	\arrow["{\widehat{e}}"', from=1-1, to=3-1]
	\arrow[""{name=1, anchor=center, inner sep=0}, "{\widehat{m}}"', from=3-1, to=3-3]
	\arrow[shorten <=9pt, shorten >=9pt, Rightarrow, from=0, to=1]
\end{tikzcd}\]
From the commutativity of the square on the above right we get $e' = \widehat{e}$, $m' = \widehat{m}$ and thus we obtain the left square in Diagram \eqref{eq_iffsquares}.

To see that $\cnr{-} \circ D \cong 1$, consider now a strict factorization system $(\ce,\cm)$ on $\cc$, we then have a functor $\cc \to \cnr{D_{\ce,\cm}}$ that is identity on objects and sends $f \mapsto (e,m)$, where $f = me$ is the unique factorization with $e \in \ce, m \in \cm$. The uniqueness of factorizations also guarantees that this functor is fully faithful, so it is an isomorphism. It also clearly preserves the classes $\ce, \cm$. 
\end{proof}

\subsection{Orthogonal factorization systems}\label{subsekce_ofs}

\begin{defi}\label{DEFI_ofs}
An \textit{orthogonal factorization system} $(\ce,\cm)$ on a category $\cc$ consists of two wide sub-categories $\ce, \cm \subseteq \cc$ satisfying\footnote{Note that this definition is equivalent to the more standard one in which the orthogonality $\ce \perp \cm$ appears. See \cite[Theorem 3.7]{joyalfs}}:
\begin{itemize}
\item For every morphism $f \in \cc$ there exist $e \in \ce, m \in \cm$ such that $f = m \circ e$, and if $f = m'e'$ is a second factorization with $e' \in \ce, m' \in \cm$, there exists a unique morphism $\theta$ so that the following diagram commutes:
\[\begin{tikzcd}
	a && {a'} && b \\
	a && {a''} && b
	\arrow["e", from=1-1, to=1-3]
	\arrow["m", from=1-3, to=1-5]
	\arrow["{e'}"', from=2-1, to=2-3]
	\arrow["{m'}"', from=2-3, to=2-5]
	\arrow[Rightarrow, no head, from=1-1, to=2-1]
	\arrow[Rightarrow, no head, from=1-5, to=2-5]
	\arrow["{\exists ! \theta}", dashed, from=1-3, to=2-3]
\end{tikzcd}\]
\item we have that $\ce \cap \cm = \{ \text{isomorphisms in }\cc \}$.
\end{itemize}
\end{defi}

In the same way as in Definition \ref{defi_SFS} define the category $\ofs$ with objects orthogonal factorization systems and morphisms being functors preserving both classes.

To describe orthogonal factorization systems as certain double categories, we will introduce a more symmetric version of a crossed double category. Given a double category $X$, denote:
\begin{equation}\label{eq_Xhvezdnicka}
X^* := ((X^v)^h)^T.
\end{equation}
This is the double category obtained from $X$ by taking both the vertical and horizontal opposites as well as the transpose.

\begin{defi}\label{DEFI_bicartesian_square}
A square $\lambda$ in a double category $X$ will be called \textit{bicartesian} if it is opcartesian in both $X$ and $X^*$. In elementary terms, this means that given a square $\alpha$ with the same top-right corner as $\lambda$, there exist unique squares $\epsilon, \delta$ so that both the bottom left composite and the bottom right composite are equal to the square $\alpha$:
\[\begin{tikzcd}
	a && b \\
	c && d & a & b & a && a && b \\
	e && d & e & d & e && c && d
	\arrow[""{name=0, anchor=center, inner sep=0}, "g", from=1-1, to=1-3]
	\arrow["u"', from=1-1, to=2-1]
	\arrow["v", from=1-3, to=2-3]
	\arrow[""{name=1, anchor=center, inner sep=0}, "h"{description}, from=2-1, to=2-3]
	\arrow["\theta"', from=2-1, to=3-1]
	\arrow[Rightarrow, no head, from=2-3, to=3-3]
	\arrow[""{name=2, anchor=center, inner sep=0}, "g", from=2-4, to=2-5]
	\arrow["l"', from=2-4, to=3-4]
	\arrow["v", from=2-5, to=3-5]
	\arrow[""{name=3, anchor=center, inner sep=0}, Rightarrow, no head, from=2-6, to=2-8]
	\arrow["l"', from=2-6, to=3-6]
	\arrow[""{name=4, anchor=center, inner sep=0}, "g", from=2-8, to=2-10]
	\arrow["u"', from=2-8, to=3-8]
	\arrow["v", from=2-10, to=3-10]
	\arrow[""{name=5, anchor=center, inner sep=0}, "k"', from=3-1, to=3-3]
	\arrow[""{name=6, anchor=center, inner sep=0}, "k"', from=3-4, to=3-5]
	\arrow[""{name=7, anchor=center, inner sep=0}, "\psi"', from=3-6, to=3-8]
	\arrow[""{name=8, anchor=center, inner sep=0}, "h"', from=3-8, to=3-10]
	\arrow["\lambda", shorten <=4pt, shorten >=4pt, Rightarrow, from=0, to=1]
	\arrow["{\exists !}", shorten <=4pt, shorten >=4pt, Rightarrow, from=1, to=5]
	\arrow["\alpha"', shorten <=4pt, shorten >=4pt, Rightarrow, from=2, to=6]
	\arrow["{\exists !}", shorten <=4pt, shorten >=4pt, Rightarrow, from=3, to=7]
	\arrow["\lambda", shorten <=4pt, shorten >=4pt, Rightarrow, from=4, to=8]
\end{tikzcd}\]
\end{defi}

\begin{defi}
A double category $X$ is \textit{top-right bicrossed} if every top-right corner can be filled into a bicartesian square, and moreover bicartesian squares are closed under horizontal and vertical compositions and contain vertical and horizontal identities.
\end{defi}

\noindent In a top-right bicrossed double category $X$, the two conditions of being opcartesian in $X$ and in $X^*$ can be expressed as a single condition as follows:

\begin{lemma}\label{LEMMA_bicartdblopcart}
Let $\lambda$ be a bicartesian square in a top-right bicrossed double category $X$ and let $\alpha$ be any square with the same top-right corner. Then there exists a unique square $\beta$ such that this equation holds:
\begin{equation}\label{EQ_lambda_bicrossed_alternative}
\begin{tikzcd}
	a && a && b && a &&&& b \\
	c && c && d & {=} \\
	e && c && d && e &&&& d
	\arrow[Rightarrow, no head, from=1-1, to=1-3]
	\arrow["u", from=1-1, to=2-1]
	\arrow["l"', curve={height=18pt}, from=1-1, to=3-1]
	\arrow[""{name=0, anchor=center, inner sep=0}, "g", from=1-3, to=1-5]
	\arrow["u"', from=1-3, to=2-3]
	\arrow["v", from=1-5, to=2-5]
	\arrow[""{name=1, anchor=center, inner sep=0}, "g", from=1-7, to=1-11]
	\arrow["l"', from=1-7, to=3-7]
	\arrow["v", from=1-11, to=3-11]
	\arrow[""{name=2, anchor=center, inner sep=0}, Rightarrow, no head, from=2-1, to=2-3]
	\arrow["\theta", from=2-1, to=3-1]
	\arrow[""{name=3, anchor=center, inner sep=0}, "h"', from=2-3, to=2-5]
	\arrow[Rightarrow, no head, from=2-3, to=3-3]
	\arrow[Rightarrow, no head, from=2-5, to=3-5]
	\arrow[""{name=4, anchor=center, inner sep=0}, "\psi"', from=3-1, to=3-3]
	\arrow["k"', curve={height=18pt}, from=3-1, to=3-5]
	\arrow["h"', from=3-3, to=3-5]
	\arrow[""{name=5, anchor=center, inner sep=0}, "k"', from=3-7, to=3-11]
	\arrow["\lambda", shorten <=4pt, shorten >=4pt, Rightarrow, from=0, to=3]
	\arrow["\alpha", shorten <=9pt, shorten >=9pt, Rightarrow, from=1, to=5]
	\arrow["{\exists ! \beta}", shorten <=4pt, shorten >=4pt, Rightarrow, from=2, to=4]
\end{tikzcd}
\end{equation}
\end{lemma}
\begin{proof}
From the definition of opcartesianness in $X^*$ there is a unique square $\gamma$ such that:
\begin{equation}\label{eq_gammalambdahoriz}
\begin{tikzcd}
	a && a && b && a && b \\
	e && c && d && e && d
	\arrow[""{name=0, anchor=center, inner sep=0}, Rightarrow, no head, from=1-1, to=1-3]
	\arrow["l"', from=1-1, to=2-1]
	\arrow[""{name=1, anchor=center, inner sep=0}, "g", from=1-3, to=1-5]
	\arrow["u"', from=1-3, to=2-3]
	\arrow[""{name=2, anchor=center, inner sep=0}, "v", from=1-5, to=2-5]
	\arrow[""{name=3, anchor=center, inner sep=0}, "g", from=1-7, to=1-9]
	\arrow[""{name=4, anchor=center, inner sep=0}, "l"', from=1-7, to=2-7]
	\arrow["v", from=1-9, to=2-9]
	\arrow[""{name=5, anchor=center, inner sep=0}, "\psi"', from=2-1, to=2-3]
	\arrow[""{name=6, anchor=center, inner sep=0}, "h"', from=2-3, to=2-5]
	\arrow[""{name=7, anchor=center, inner sep=0}, "k"', from=2-7, to=2-9]
	\arrow["\gamma", shorten <=4pt, shorten >=4pt, Rightarrow, from=0, to=5]
	\arrow["\lambda", shorten <=4pt, shorten >=4pt, Rightarrow, from=1, to=6]
	\arrow["{=}"{description}, draw=none, from=2, to=4]
	\arrow["\alpha", shorten <=4pt, shorten >=4pt, Rightarrow, from=3, to=7]
\end{tikzcd}
\end{equation}

Because the horizontal identity square on a vertical morphism $u$ is opcartesian in $X$ ($X$ is bicrossed), there exists a unique square $\gamma'$ such that:
\[\begin{tikzcd}
	a && a && a && a \\
	c && c & {=} \\
	e && c && e && c
	\arrow[Rightarrow, no head, from=1-1, to=1-3]
	\arrow["u"', from=1-1, to=2-1]
	\arrow["u", from=1-3, to=2-3]
	\arrow[""{name=0, anchor=center, inner sep=0}, Rightarrow, no head, from=1-5, to=1-7]
	\arrow["l"', from=1-5, to=3-5]
	\arrow["u", from=1-7, to=3-7]
	\arrow[""{name=1, anchor=center, inner sep=0}, Rightarrow, no head, from=2-1, to=2-3]
	\arrow["\theta"', from=2-1, to=3-1]
	\arrow[Rightarrow, no head, from=2-3, to=3-3]
	\arrow[""{name=2, anchor=center, inner sep=0}, "\psi"', from=3-1, to=3-3]
	\arrow[""{name=3, anchor=center, inner sep=0}, "\psi"', from=3-5, to=3-7]
	\arrow["\gamma", shorten <=9pt, shorten >=9pt, Rightarrow, from=0, to=3]
	\arrow["{\gamma'}", shorten <=4pt, shorten >=4pt, Rightarrow, from=1, to=2]
\end{tikzcd}\]
This gives us the \textbf{existence}. To prove the  \textbf{uniqueness}, let $\beta$ be a different square satisfying the equation in the statement. Then the composite $\beta \circ 1^\bullet_u$ (vertical composite of $\beta$ and the horizontal identity square on $u$) satisfies the equation \eqref{eq_gammalambdahoriz}. Because $\lambda$ is opcartesian in $X^*$, this forces $\beta \circ  1^\bullet_{u} = \gamma = \gamma' \circ 1^\bullet_{u}$. Because $1^\bullet_{u}$ is opcartesian in $X$, this in turn forces $\beta = \gamma'$.
\end{proof}

\begin{pozn}\label{POZN_bicrossed_is_bicartesian_NEW}
On the other hand, given a top-right bicrossed double category $X$ and a square $\lambda$ with the property described in Lemma \ref{LEMMA_bicartdblopcart}, it is easily seen that $\lambda$ is bicartesian. For example for the opcartesianness in $X^*$, the \textbf{existence} of a square $\delta$ such that the rightmost picture in Definition \ref{DEFI_bicartesian_square} equals $\alpha$, is immediate. For the \textbf{uniqueness}, note that any square $\delta$ as in this definition factors uniquely like this (because the horizontal identity square on $u$ is opcartesian in $X$):
\[\begin{tikzcd}
	a && a && a && a \\
	&&& {=} & c && c \\
	e && c && e && c
	\arrow[""{name=0, anchor=center, inner sep=0}, Rightarrow, no head, from=1-1, to=1-3]
	\arrow["l"', from=1-1, to=3-1]
	\arrow["u", from=1-3, to=3-3]
	\arrow[Rightarrow, no head, from=1-5, to=1-7]
	\arrow["u"', from=1-5, to=2-5]
	\arrow["u", from=1-7, to=2-7]
	\arrow[""{name=1, anchor=center, inner sep=0}, Rightarrow, no head, from=2-5, to=2-7]
	\arrow["{\widehat{u}}"', from=2-5, to=3-5]
	\arrow[Rightarrow, no head, from=2-7, to=3-7]
	\arrow[""{name=2, anchor=center, inner sep=0}, "\psi"', from=3-1, to=3-3]
	\arrow[""{name=3, anchor=center, inner sep=0}, "\psi"', from=3-5, to=3-7]
	\arrow["\delta", shorten <=9pt, shorten >=9pt, Rightarrow, from=0, to=2]
	\arrow["{\exists ! \delta'}", shorten <=4pt, shorten >=4pt, Rightarrow, from=1, to=3]
\end{tikzcd}\]
\end{pozn}

\begin{pozn}[Representable definition]\label{POZN_repr_defi_bicrossed}
For any double category $X$, there is a category $\text{BotLeft}(X)$ whose objects are the objects of $X$, and for which a morphism $a \to b$ is a square of this form:
\[\begin{tikzcd}
	a & a \\
	b & a
	\arrow["l"', from=1-1, to=2-1]
	\arrow[""{name=0, anchor=center, inner sep=0}, Rightarrow, no head, from=1-2, to=1-1]
	\arrow[""{name=1, anchor=center, inner sep=0}, "k"', from=2-1, to=2-2]
	\arrow[Rightarrow, no head, from=2-2, to=1-2]
	\arrow["\alpha", shorten <=4pt, shorten >=4pt, Rightarrow, from=0, to=1]
\end{tikzcd}\]
Given a top-right corner $(g,v)$ in the double category:
\[\begin{tikzcd}
	x & y \\
	& z
	\arrow["g", from=1-1, to=1-2]
	\arrow["v", from=1-2, to=2-2]
\end{tikzcd}\]
There is a functor $\sq{(g,v),-}:\text{BotLeft}(X) \to \set$ sending an object $a \in X$ to the set of squares with top-right corner being the pair $(g,v)$ that have the object $a$ in the bottom-left spot. Given a morphism $\alpha: a \to b$ as above, the function $\sq{(g,v),a} \to \sq{(g,v),b}$ is given by the assignment:
\[\begin{tikzcd}
	z && y && x & x & y \\
	&&& \mapsto & a & a & z \\
	a && x && b & a & z
	\arrow[""{name=0, anchor=center, inner sep=0}, "g", from=1-1, to=1-3]
	\arrow["u"', from=1-1, to=3-1]
	\arrow["v", from=1-3, to=3-3]
	\arrow["u"', from=1-5, to=2-5]
	\arrow[Rightarrow, no head, from=1-6, to=1-5]
	\arrow[""{name=1, anchor=center, inner sep=0}, "g", from=1-6, to=1-7]
	\arrow["u"', from=1-6, to=2-6]
	\arrow["v", from=1-7, to=2-7]
	\arrow["l"', from=2-5, to=3-5]
	\arrow[""{name=2, anchor=center, inner sep=0}, Rightarrow, no head, from=2-6, to=2-5]
	\arrow[""{name=3, anchor=center, inner sep=0}, "h"', from=2-6, to=2-7]
	\arrow[""{name=4, anchor=center, inner sep=0}, "h"', from=3-1, to=3-3]
	\arrow[""{name=5, anchor=center, inner sep=0}, "k"', from=3-5, to=3-6]
	\arrow[Rightarrow, no head, from=3-6, to=2-6]
	\arrow["h"', from=3-6, to=3-7]
	\arrow[Rightarrow, no head, from=3-7, to=2-7]
	\arrow["\beta"', shorten <=9pt, shorten >=9pt, Rightarrow, from=0, to=4]
	\arrow["\beta"', shorten <=4pt, shorten >=4pt, Rightarrow, from=1, to=3]
	\arrow["\alpha"', shorten <=4pt, shorten >=4pt, Rightarrow, from=2, to=5]
\end{tikzcd}\]

\noindent We have:
\begin{prop}
Given a top-right corner $(g,v)$ in a double category $X$, the following are equivalent:
\begin{itemize}
\item the functor $\sq{(g,v),-}: \text{BotLeft}(X) \to \set$ is representable,
\item the top-right corner $(g,v)$ can be filled into a square $\lambda$ satisfying the property \eqref{EQ_lambda_bicrossed_alternative}.
\end{itemize}
\end{prop}
\end{pozn}

\begin{pr}
In the double category $\sq{\cc}^{v}$, a square is bicartesian if and only if it is a pullback square. If $\cc$ has pullbacks, the double category $\sq{\cc}^{v}$ is top-right bicrossed.

Both $\text{MPbSq}(\cc)^{v}$ (for a category $\cc$ with pullbacks) and $\text{BOFib}(\ce)^{v}$ are top-right bicrossed with every square being bicartesian.
\end{pr}
\color{black}

We will now focus on double categories in which every top-right corner can be filled into a square and every square is bicartesian. Such double categories are automatically crossed so we may again use the category of corners construction. Notice that in this case, two corners $(e,m), (e',m')$ are equivalent if and only if there is a single 2-cell between them (Lemma \ref{lemmaXpscr_everysquareopcart}), and moreover every 2-cell has the following form (since the vertical identities on morphisms are bicartesian):
\begin{equation}\label{PIC_twocellforbicrosseddblcat}
\begin{tikzcd}
	a && a \\
	{a'} && {a'} && b \\
	{a''} && {a'} && b
	\arrow[Rightarrow, no head, from=2-5, to=3-5]
	\arrow["m", from=2-3, to=2-5]
	\arrow["e", from=1-3, to=2-3]
	\arrow[Rightarrow, no head, from=1-1, to=1-3]
	\arrow["e"', from=1-1, to=2-1]
	\arrow[""{name=0, anchor=center, inner sep=0}, Rightarrow, no head, from=2-1, to=2-3]
	\arrow["\theta"', from=2-1, to=3-1]
	\arrow[Rightarrow, no head, from=2-3, to=3-3]
	\arrow[""{name=1, anchor=center, inner sep=0}, "\psi"', from=3-1, to=3-3]
	\arrow["m"', from=3-3, to=3-5]
	\arrow["{e'}"', curve={height=24pt}, from=1-1, to=3-1]
	\arrow["{m'}"', curve={height=24pt}, from=3-1, to=3-5]
	\arrow["\beta", shorten <=4pt, shorten >=4pt, Rightarrow, from=0, to=1]
\end{tikzcd}
\end{equation}

\noindent Recall Notation \ref{notecnrX_verticalhorizontalcnrs}.

\begin{lemma}\label{lemma_isossubcxcm}
Let $X$ be a double category in which every top-right corner can be filled into a square and every square is bicartesian. Then:
\[
\ce_X \cap \cm_X \subseteq \{ \text{isomorphisms in } \cnr{X} \}.
\]
\end{lemma}
\begin{proof}
Let $[u,1] = [1,h]$ be a morphism in the intersection. There is then a 2-cell as follows (see the remark above diagram \eqref{PIC_twocellforbicrosseddblcat}):
\[\begin{tikzcd}
	a && a \\
	{a'} && {a'} && {a'} \\
	a && {a'} && {a'}
	\arrow["u", from=1-3, to=2-3]
	\arrow[Rightarrow, no head, from=2-3, to=2-5]
	\arrow["u"', from=1-1, to=2-1]
	\arrow["\theta"', from=2-1, to=3-1]
	\arrow[""{name=0, anchor=center, inner sep=0}, "h"', from=3-1, to=3-3]
	\arrow[Rightarrow, no head, from=2-3, to=3-3]
	\arrow[""{name=1, anchor=center, inner sep=0}, Rightarrow, no head, from=2-1, to=2-3]
	\arrow[Rightarrow, no head, from=3-3, to=3-5]
	\arrow[Rightarrow, no head, from=2-5, to=3-5]
	\arrow[Rightarrow, no head, from=1-1, to=1-3]
	\arrow[curve={height=18pt}, Rightarrow, no head, from=1-1, to=3-1]
	\arrow["h"', curve={height=18pt}, from=3-1, to=3-5]
	\arrow["\beta", shorten <=4pt, shorten >=4pt, Rightarrow, from=1, to=0]
\end{tikzcd}\]
Since every square is bicartesian, $\theta, h$ are isomorphisms, and so is $\theta^{-1} = u$. Hence $[u,1]$ is an isomorphism in $\cnr{X}$ with the inverse being $[u^{-1},1]$.

\end{proof}

\begin{lemma}
Let $X$ be a double category in which every top-right corner can be filled into a square and every square is bicartesian. An equivalence class of corners $[v,h]$ is invertible in $\cnr{X}$ if and only if both $v$ and $h$ are isomorphisms.
\end{lemma}
\begin{proof}
Let $[v,h]$ be an isomorphism in $\cnr{X}$ with inverse $[u,g]$ as pictured together with the inverse laws below:
\begin{equation}\label{eq_ugvhinverse}
\begin{tikzcd}[column sep=small]
	a & a &&&&&& b & b \\
	{a'} & {a'} && b &&&& {b'} & {b'} && a \\
	a & a && {b'} && a && b & b && {a'} && b \\
	a & a && {b'} && a && b & b && {a'} && b
	\arrow[Rightarrow, no head, from=1-1, to=1-2]
	\arrow["u"', from=1-1, to=2-1]
	\arrow[curve={height=24pt}, Rightarrow, no head, from=1-1, to=4-1]
	\arrow["u", from=1-2, to=2-2]
	\arrow[Rightarrow, no head, from=1-8, to=1-9]
	\arrow["v"', from=1-8, to=2-8]
	\arrow[curve={height=24pt}, Rightarrow, no head, from=1-8, to=4-8]
	\arrow["v", from=1-9, to=2-9]
	\arrow["{\widehat{v}}"', from=2-1, to=3-1]
	\arrow[""{name=0, anchor=center, inner sep=0}, "g", from=2-2, to=2-4]
	\arrow[""{name=0p, anchor=center, inner sep=0}, phantom, from=2-2, to=2-4, start anchor=center, end anchor=center]
	\arrow["{\widehat{v}}"', from=2-2, to=3-2]
	\arrow["v", from=2-4, to=3-4]
	\arrow["{\widehat{u}}"', from=2-8, to=3-8]
	\arrow[""{name=1, anchor=center, inner sep=0}, "h", from=2-9, to=2-11]
	\arrow[""{name=1p, anchor=center, inner sep=0}, phantom, from=2-9, to=2-11, start anchor=center, end anchor=center]
	\arrow["{\widehat{u}}"', from=2-9, to=3-9]
	\arrow["u", from=2-11, to=3-11]
	\arrow[""{name=2, anchor=center, inner sep=0}, Rightarrow, no head, from=3-1, to=3-2]
	\arrow[""{name=2p, anchor=center, inner sep=0}, phantom, from=3-1, to=3-2, start anchor=center, end anchor=center]
	\arrow["\theta"', from=3-1, to=4-1]
	\arrow[""{name=3, anchor=center, inner sep=0}, "{\widehat{g}}"', from=3-2, to=3-4]
	\arrow[""{name=3p, anchor=center, inner sep=0}, phantom, from=3-2, to=3-4, start anchor=center, end anchor=center]
	\arrow[Rightarrow, no head, from=3-2, to=4-2]
	\arrow["h", from=3-4, to=3-6]
	\arrow[Rightarrow, no head, from=3-6, to=4-6]
	\arrow[""{name=4, anchor=center, inner sep=0}, Rightarrow, no head, from=3-8, to=3-9]
	\arrow[""{name=4p, anchor=center, inner sep=0}, phantom, from=3-8, to=3-9, start anchor=center, end anchor=center]
	\arrow["{\theta'}"{description}, from=3-8, to=4-8]
	\arrow[""{name=5, anchor=center, inner sep=0}, "{\widehat{h}}"', from=3-9, to=3-11]
	\arrow[""{name=5p, anchor=center, inner sep=0}, phantom, from=3-9, to=3-11, start anchor=center, end anchor=center]
	\arrow[Rightarrow, no head, from=3-9, to=4-9]
	\arrow["g", from=3-11, to=3-13]
	\arrow[Rightarrow, no head, from=3-13, to=4-13]
	\arrow[""{name=6, anchor=center, inner sep=0}, "\psi"{description}, from=4-1, to=4-2]
	\arrow[""{name=6p, anchor=center, inner sep=0}, phantom, from=4-1, to=4-2, start anchor=center, end anchor=center]
	\arrow[curve={height=24pt}, Rightarrow, no head, from=4-1, to=4-6]
	\arrow["{\widehat{g}}"', from=4-2, to=4-4]
	\arrow["h"', from=4-4, to=4-6]
	\arrow[""{name=7, anchor=center, inner sep=0}, "{\psi'}"{description}, from=4-8, to=4-9]
	\arrow[""{name=7p, anchor=center, inner sep=0}, phantom, from=4-8, to=4-9, start anchor=center, end anchor=center]
	\arrow[curve={height=24pt}, Rightarrow, no head, from=4-8, to=4-13]
	\arrow["{\widehat{h}}"', from=4-9, to=4-11]
	\arrow["g"', from=4-11, to=4-13]
	\arrow[shorten <=4pt, shorten >=4pt, Rightarrow, from=0p, to=3p]
	\arrow[shorten <=4pt, shorten >=4pt, Rightarrow, from=1p, to=5p]
	\arrow[shorten <=4pt, shorten >=4pt, Rightarrow, from=2p, to=6p]
	\arrow[shorten <=4pt, shorten >=4pt, Rightarrow, from=4p, to=7p]
\end{tikzcd}
\end{equation}

From the pictures below we obtain the following equalities:
\begin{align*}
[v,1] \circ [u,g] &= [1,\widehat{g}\psi], \\
[u,g] \circ [1,h] &= [\theta' \widehat{u},1].
\end{align*}
\[\begin{tikzcd}[column sep=scriptsize]
	a & a \\
	{a'} & {a'} && b &&& {b'} & {b'} && a \\
	a & a && {b'} &&& b & b && {a'} && b \\
	a & a && {b'} &&& b & b && {a'} && b
	\arrow[Rightarrow, no head, from=1-1, to=1-2]
	\arrow["u"', from=1-1, to=2-1]
	\arrow[curve={height=24pt}, Rightarrow, no head, from=1-1, to=4-1]
	\arrow["u", from=1-2, to=2-2]
	\arrow["{\widehat{v}}"', from=2-1, to=3-1]
	\arrow[""{name=0, anchor=center, inner sep=0}, "g", from=2-2, to=2-4]
	\arrow[""{name=0p, anchor=center, inner sep=0}, phantom, from=2-2, to=2-4, start anchor=center, end anchor=center]
	\arrow["{\widehat{v}}"', from=2-2, to=3-2]
	\arrow["v", from=2-4, to=3-4]
	\arrow[Rightarrow, no head, from=2-7, to=2-8]
	\arrow["{\widehat{u}}"', from=2-7, to=3-7]
	\arrow["{\theta'\widehat{u}}"', curve={height=24pt}, from=2-7, to=4-7]
	\arrow[""{name=1, anchor=center, inner sep=0}, "h", from=2-8, to=2-10]
	\arrow[""{name=1p, anchor=center, inner sep=0}, phantom, from=2-8, to=2-10, start anchor=center, end anchor=center]
	\arrow["{\widehat{u}}"', from=2-8, to=3-8]
	\arrow["u", from=2-10, to=3-10]
	\arrow[""{name=2, anchor=center, inner sep=0}, Rightarrow, no head, from=3-1, to=3-2]
	\arrow[""{name=2p, anchor=center, inner sep=0}, phantom, from=3-1, to=3-2, start anchor=center, end anchor=center]
	\arrow["\theta"', from=3-1, to=4-1]
	\arrow[""{name=3, anchor=center, inner sep=0}, "{\widehat{g}}"', from=3-2, to=3-4]
	\arrow[""{name=3p, anchor=center, inner sep=0}, phantom, from=3-2, to=3-4, start anchor=center, end anchor=center]
	\arrow[Rightarrow, no head, from=3-2, to=4-2]
	\arrow[Rightarrow, no head, from=3-4, to=4-4]
	\arrow[""{name=4, anchor=center, inner sep=0}, Rightarrow, no head, from=3-7, to=3-8]
	\arrow[""{name=4p, anchor=center, inner sep=0}, phantom, from=3-7, to=3-8, start anchor=center, end anchor=center]
	\arrow["{\theta'}"', from=3-7, to=4-7]
	\arrow[""{name=5, anchor=center, inner sep=0}, "{\widehat{h}}"', from=3-8, to=3-10]
	\arrow[""{name=5p, anchor=center, inner sep=0}, phantom, from=3-8, to=3-10, start anchor=center, end anchor=center]
	\arrow[Rightarrow, no head, from=3-8, to=4-8]
	\arrow["g", from=3-10, to=3-12]
	\arrow[Rightarrow, no head, from=3-12, to=4-12]
	\arrow[""{name=6, anchor=center, inner sep=0}, "\psi"', from=4-1, to=4-2]
	\arrow[""{name=6p, anchor=center, inner sep=0}, phantom, from=4-1, to=4-2, start anchor=center, end anchor=center]
	\arrow["{\widehat{g}\psi}"', curve={height=24pt}, from=4-1, to=4-4]
	\arrow["{\widehat{g}}"', from=4-2, to=4-4]
	\arrow[""{name=7, anchor=center, inner sep=0}, "{\psi'}"{description}, from=4-7, to=4-8]
	\arrow[""{name=7p, anchor=center, inner sep=0}, phantom, from=4-7, to=4-8, start anchor=center, end anchor=center]
	\arrow[curve={height=24pt}, Rightarrow, no head, from=4-7, to=4-12]
	\arrow["{\widehat{h}}"', from=4-8, to=4-10]
	\arrow["g"', from=4-10, to=4-12]
	\arrow[shorten <=4pt, shorten >=4pt, Rightarrow, from=0p, to=3p]
	\arrow[shorten <=4pt, shorten >=4pt, Rightarrow, from=1p, to=5p]
	\arrow[shorten <=4pt, shorten >=4pt, Rightarrow, from=2p, to=6p]
	\arrow[shorten <=4pt, shorten >=4pt, Rightarrow, from=4p, to=7p]
\end{tikzcd}\]

Now consider the following composite:
\[
[v \theta' \widehat{u},1] = [v,1] \circ [\theta' \widehat{u},1] = [v,1] \circ [u,g] \circ [1,h] = [1,\widehat{g}\psi] \circ [1,h] = [1,\widehat{g}\psi h]
\]

This composite belongs both in $\ce_X$ and $\cm_X$, so as in the proof of Lemma \ref{lemma_isossubcxcm} we obtain that $(\widehat{g}\psi) h$ is an isomorphism that we denote by $\Theta$. This implies that $\Theta^{-1}(\widehat{g}\psi) h = 1$ and so $h$ is a split monomorphism. Since $h (\widehat{g}\psi) = 1$ by Equation \eqref{eq_ugvhinverse}, $h$ is a split epimorphism. Thus $h$ is an isomorphism and by similar reasoning, $v$ is also an isomorphism. Conversely, if $[v,h]$ is a corner for which $v,h$ are isomorphisms, it is an isomorphism in $\cnr{X}$ with the inverse given by $[v^{-1},1] \circ [1,h^{-1}]$.
\end{proof}

Note that for $X = \pbsq{\cc}^{v}$ the above lemma gives the usual folklore characterization of isomorphisms in the category $\spn{\cc}$ of spans.

\begin{lemma}\label{lemma_isossupcxcm}
Let $X$ be a double category in which every top-right corner can be filled into a square and every square is bicartesian. Assume in addition that $X$ is invariant. Then:
\[
\ce_X \cap \cm_X \supseteq \{ \text{isomorphisms in } \cnr{X} \}.
\]
\end{lemma}
\begin{proof}

We will show that when $X$ is horizontally invariant, we have:
\[
[u,g] \in \cnr{X}\text{, }g\text{ is an isomorphism}\Rightarrow [u,g] \in \ce_X.
\]

Let $[u,g]$ be such a corner. From horizontal invariance we get the square (pictured below), that exhibits the equality $[u,g] = [\theta u, 1]$:
\[\begin{tikzcd}
	a \\
	{a'} & b \\
	b & b
	\arrow["u"', from=1-1, to=2-1]
	\arrow["{\theta u}"', curve={height=24pt}, from=1-1, to=3-1]
	\arrow[""{name=0, anchor=center, inner sep=0}, "g", from=2-1, to=2-2]
	\arrow[""{name=0p, anchor=center, inner sep=0}, phantom, from=2-1, to=2-2, start anchor=center, end anchor=center]
	\arrow["\theta"', from=2-1, to=3-1]
	\arrow[Rightarrow, no head, from=2-2, to=3-2]
	\arrow[""{name=1, anchor=center, inner sep=0}, Rightarrow, no head, from=3-1, to=3-2]
	\arrow[""{name=1p, anchor=center, inner sep=0}, phantom, from=3-1, to=3-2, start anchor=center, end anchor=center]
	\arrow[shorten <=4pt, shorten >=4pt, Rightarrow, from=0p, to=1p]
\end{tikzcd}\]
Thus, $[u,g] \in \ce_X$. Dually, if $X$ is vertically invariant, we have:
\[
[u,g] \in \cnr{X}\text{, }u\text{ is an isomorphism}\Rightarrow [u,g] \in \cm_X.
\]

Now if $[u,g]$ is an isomorphism, by previous lemma both $u,g$ are isomorphisms, and by the above implications $[u,g]$ belongs to both $\ce_X$ and $\cm_X$.
\end{proof}

The following somewhat technical notion is introduced in this chapter to guarantee the uniqueness of factorizations up to a unique morphism:

\begin{defi}\label{defi_jointly_monic}
A bottom-left corner $(\pi_1, \pi_2)$ is said to be \textit{jointly monic} if, given squares $\kappa_1$, $\kappa_2$ as pictured below:
\[\begin{tikzcd}
	a &&&& {a'} && {a'} && {a'} && {a'} \\
	{a'} && b && {a''} && {a'} && {a''} && {a'}
	\arrow["{\pi_1}", from=1-1, to=2-1]
	\arrow[""{name=0, anchor=center, inner sep=0}, Rightarrow, no head, from=1-5, to=1-7]
	\arrow["\theta"', from=1-5, to=2-5]
	\arrow[""{name=1, anchor=center, inner sep=0}, Rightarrow, no head, from=1-9, to=1-11]
	\arrow["{\theta'}"', from=1-9, to=2-9]
	\arrow["{\pi_2}", from=2-1, to=2-3]
	\arrow[""{name=2, anchor=center, inner sep=0}, "\psi"', from=2-5, to=2-7]
	\arrow[Rightarrow, no head, from=2-7, to=1-7]
	\arrow[""{name=3, anchor=center, inner sep=0}, "{\psi'}"', from=2-9, to=2-11]
	\arrow[Rightarrow, no head, from=2-11, to=1-11]
	\arrow["{\kappa_1}"', shorten <=4pt, shorten >=4pt, Rightarrow, from=0, to=2]
	\arrow["{\kappa_2}"', shorten <=4pt, shorten >=4pt, Rightarrow, from=1, to=3]
\end{tikzcd}\]
we have the following implication:
\[
(\theta \pi_1 = \theta' \pi_1 \land \pi_2 \psi = \pi_2 \psi' )\hspace{1cm} \Rightarrow \hspace{1cm} \theta = \theta', \psi = \psi'.
\]

In the examples below, we will consider bottom-left corners in double categories of form $X^v$. Those are the same as top-left corners $(\pi_1, \pi_2)$ in $X$ having the property that, when given squares $\kappa_1$, $\kappa_2$ pictured below:
\[\begin{tikzcd}
	{a'} && b && {a''} && {a'} && {a''} && {a'} \\
	a &&&& {a'} && {a'} && {a'} && {a'}
	\arrow["{\pi_2}", from=1-1, to=1-3]
	\arrow["{\pi_1}"', from=1-1, to=2-1]
	\arrow[""{name=0, anchor=center, inner sep=0}, "\psi", from=1-5, to=1-7]
	\arrow[""{name=0p, anchor=center, inner sep=0}, phantom, from=1-5, to=1-7, start anchor=center, end anchor=center]
	\arrow["\theta"', from=1-5, to=2-5]
	\arrow[Rightarrow, no head, from=1-7, to=2-7]
	\arrow[""{name=1, anchor=center, inner sep=0}, "{\psi'}", from=1-9, to=1-11]
	\arrow[""{name=1p, anchor=center, inner sep=0}, phantom, from=1-9, to=1-11, start anchor=center, end anchor=center]
	\arrow["{\theta'}"', from=1-9, to=2-9]
	\arrow[Rightarrow, no head, from=1-11, to=2-11]
	\arrow[""{name=2, anchor=center, inner sep=0}, Rightarrow, no head, from=2-5, to=2-7]
	\arrow[""{name=2p, anchor=center, inner sep=0}, phantom, from=2-5, to=2-7, start anchor=center, end anchor=center]
	\arrow[""{name=3, anchor=center, inner sep=0}, Rightarrow, no head, from=2-9, to=2-11]
	\arrow[""{name=3p, anchor=center, inner sep=0}, phantom, from=2-9, to=2-11, start anchor=center, end anchor=center]
	\arrow["{\kappa_1}", shorten <=4pt, shorten >=4pt, Rightarrow, from=0p, to=2p]
	\arrow["{\kappa_2}", shorten <=4pt, shorten >=4pt, Rightarrow, from=1p, to=3p]
\end{tikzcd}\]
We have the following implication:
\[
(\pi_1 \theta = \pi_1 \theta' \land \pi_2 \psi = \pi_2 \psi' ) \Rightarrow \theta = \theta', \psi = \psi'.
\]
We will say the top-left corner in $X$ is \textit{jointly monic} if the corresponding bottom-left corner is so in $X^v$.
\end{defi}

\begin{pr}\label{prsq_jointmonicity}
In $\sq{\cc}$ a top-left corner $(\pi_1,\pi_2)$ is jointly monic if and only if the pair of morphisms $(\pi_1,\pi_2)$ is jointly monic in the category $\cc$.
In $\pbsq{\cc}$ a top-left corner is jointly monic if and only if the pair is jointly monic in $\cc$ \textbf{with respect to all isomorphisms}.
\end{pr}

\begin{pr}\label{prmpbsq_pbsqsjointmonicity}
In the double category $\text{MPbSq}(\cc)$ every top-left corner is jointly monic as we now show: Let there be squares $\kappa_1, \kappa_2$ as in the definition, then $\theta = \psi$, $\theta' = \psi'$ and the equality $\pi_1 \theta = \pi_1 \theta'$ forces $\theta = \theta'$ because $\pi_1$ is a monomorphism.
\end{pr}

\begin{pr}\label{PR_BOFIB_pbsqsjointmonicity}
Every top-left corner is jointly monic in the double category $\text{BOFib}(\ce)$ as well. To see this, let there be a top-left corner and squares as pictured below:
\[\begin{tikzcd}
	C && B && D && C && D && C \\
	\\
	A &&&& C && C && C && C
	\arrow["G", from=1-1, to=1-3]
	\arrow["F"', from=1-1, to=3-1]
	\arrow[""{name=0, anchor=center, inner sep=0}, "{P'}", from=1-9, to=1-11]
	\arrow["H"', from=1-5, to=3-5]
	\arrow[""{name=1, anchor=center, inner sep=0}, "P", from=1-5, to=1-7]
	\arrow[Rightarrow, no head, from=1-7, to=3-7]
	\arrow[""{name=2, anchor=center, inner sep=0}, Rightarrow, no head, from=3-5, to=3-7]
	\arrow["{H'}"', from=1-9, to=3-9]
	\arrow[""{name=3, anchor=center, inner sep=0}, Rightarrow, no head, from=3-9, to=3-11]
	\arrow[Rightarrow, no head, from=1-11, to=3-11]
	\arrow["{\kappa_2}", shorten <=9pt, shorten >=9pt, Rightarrow, from=0, to=3]
	\arrow["{\kappa_1}", shorten <=9pt, shorten >=9pt, Rightarrow, from=1, to=2]
\end{tikzcd}\]

\noindent Assume that $FH = FH'$, $GP = GP'$. We again get $H = P, H' = P'$. Now since $F$ is a bijection on objects, we have $H_0 = H'_0$ (the object parts of the functors agree). Because $G$ is a discrete opfibration, the square below is a pullback and we obtain $H_1 = H'_1$ as well:
\[\begin{tikzcd}
	{D_1} \\
	& {C_1} && {B_1} \\
	\\
	& {C_0} && {B_0}
	\arrow["s", from=2-4, to=4-4]
	\arrow["s"', from=2-2, to=4-2]
	\arrow["{G_1}", from=2-2, to=2-4]
	\arrow["{G_0}"', from=4-2, to=4-4]
	\arrow["\lrcorner"{anchor=center, pos=0.125}, draw=none, from=2-2, to=4-4]
	\arrow["{H_1'}", shift left=2, from=1-1, to=2-2]
	\arrow["{H_1}"', shift right=2, from=1-1, to=2-2]
	\arrow["{sH_0 = sH'_0}"', curve={height=24pt}, from=1-1, to=4-2]
	\arrow["{G_1H_1' = G_1H_1}", curve={height=-24pt}, from=1-1, to=2-4]
\end{tikzcd}\]
\end{pr}

\begin{lemma}\label{lemma_orthog}
Let $X$ be a double category in which every top-right corner can be filled into a square and every square is bicartesian. Assume further that every bottom-left corner in $X$ is jointly monic. Then the $(\ce_X,\cm_X)$-factorization of a morphism in $\cnr{X}$ is unique up to a unique morphism.
\end{lemma}
\begin{proof}
Assume that $[e,m] = [e',m']$ are two $(\ce_X,\cm_X)$-factorizations of a morphism in $\cnr{X}$. We wish to show that there is a unique morphism between them:
\begin{equation}\label{eqmorfismusfaktorizaci}
\begin{tikzcd}
	a && {a'} && b \\
	a && {a''} && b
	\arrow["{[e,1]}", from=1-1, to=1-3]
	\arrow["{[1,m]}", from=1-3, to=1-5]
	\arrow["{[e',1]}"', from=2-1, to=2-3]
	\arrow["{[1,m']}"', from=2-3, to=2-5]
	\arrow[Rightarrow, no head, from=1-1, to=2-1]
	\arrow[Rightarrow, no head, from=1-5, to=2-5]
	\arrow[dashed, from=1-3, to=2-3]
\end{tikzcd}
\end{equation}

As in the proof of Theorem \ref{THMcancprops_weakorthog}, one such morphism is given by the corner $[\theta,1]$, where $\theta$ is the domain of the 2-cell square between $(e,m)$ and $(e',m')$:
\[\begin{tikzcd}
	a && a \\
	{a'} && {a'} && b \\
	{a''} && {a'} && b
	\arrow["e", from=1-3, to=2-3]
	\arrow["m", from=2-3, to=2-5]
	\arrow["e"', from=1-1, to=2-1]
	\arrow["\theta"', from=2-1, to=3-1]
	\arrow[""{name=0, anchor=center, inner sep=0}, "\psi"', from=3-1, to=3-3]
	\arrow[Rightarrow, no head, from=2-3, to=3-3]
	\arrow[""{name=1, anchor=center, inner sep=0}, Rightarrow, no head, from=2-1, to=2-3]
	\arrow["m"', from=3-3, to=3-5]
	\arrow[Rightarrow, no head, from=2-5, to=3-5]
	\arrow[Rightarrow, no head, from=1-1, to=1-3]
	\arrow["{e'}"', curve={height=18pt}, from=1-1, to=3-1]
	\arrow["{m'}"', curve={height=18pt}, from=3-1, to=3-5]
	\arrow["\beta", shorten <=4pt, shorten >=4pt, Rightarrow, from=1, to=0]
\end{tikzcd}\]

Assume that there is a different morphism $[s,t]: a' \to a''$ making both squares in \eqref{eqmorfismusfaktorizaci} commute. The commutativity of these two squares gives the following 2-cells:

\[\begin{tikzcd}[column sep=small]
	a && a \\
	{a'} && {a'} &&&& {a'} && {a'} \\
	x && x && {a''} && x && x && {a''} && b \\
	{a''} && x && {a''} && {a'} && x && {a''} && b
	\arrow[Rightarrow, no head, from=1-1, to=1-3]
	\arrow["e", from=1-1, to=2-1]
	\arrow["{e'}"', curve={height=24pt}, from=1-1, to=4-1]
	\arrow["e", from=1-3, to=2-3]
	\arrow["s"', from=2-1, to=3-1]
	\arrow["s", from=2-3, to=3-3]
	\arrow[Rightarrow, no head, from=2-7, to=2-9]
	\arrow["s"', from=2-7, to=3-7]
	\arrow[curve={height=18pt}, Rightarrow, no head, from=2-7, to=4-7]
	\arrow["s", from=2-9, to=3-9]
	\arrow[""{name=0, anchor=center, inner sep=0}, Rightarrow, no head, from=3-1, to=3-3]
	\arrow[""{name=0p, anchor=center, inner sep=0}, phantom, from=3-1, to=3-3, start anchor=center, end anchor=center]
	\arrow["{\theta'}"{description}, from=3-1, to=4-1]
	\arrow["t", from=3-3, to=3-5]
	\arrow[Rightarrow, no head, from=3-3, to=4-3]
	\arrow[Rightarrow, no head, from=3-5, to=4-5]
	\arrow[""{name=1, anchor=center, inner sep=0}, Rightarrow, no head, from=3-7, to=3-9]
	\arrow[""{name=1p, anchor=center, inner sep=0}, phantom, from=3-7, to=3-9, start anchor=center, end anchor=center]
	\arrow["{\widetilde{\theta}}"', from=3-7, to=4-7]
	\arrow["t", from=3-9, to=3-11]
	\arrow[Rightarrow, no head, from=3-9, to=4-9]
	\arrow["{m'}", from=3-11, to=3-13]
	\arrow[Rightarrow, no head, from=3-13, to=4-13]
	\arrow[""{name=2, anchor=center, inner sep=0}, "{\psi'}"{description}, from=4-1, to=4-3]
	\arrow[""{name=2p, anchor=center, inner sep=0}, phantom, from=4-1, to=4-3, start anchor=center, end anchor=center]
	\arrow[curve={height=18pt}, Rightarrow, no head, from=4-1, to=4-5]
	\arrow["t"', from=4-3, to=4-5]
	\arrow[""{name=3, anchor=center, inner sep=0}, "{\widetilde{\psi}}"{description}, from=4-7, to=4-9]
	\arrow[""{name=3p, anchor=center, inner sep=0}, phantom, from=4-7, to=4-9, start anchor=center, end anchor=center]
	\arrow["m"', curve={height=24pt}, from=4-7, to=4-13]
	\arrow["t"', from=4-9, to=4-11]
	\arrow["{m'}"{description}, from=4-11, to=4-13]
	\arrow["{\beta'}", shorten <=4pt, shorten >=4pt, Rightarrow, from=0p, to=2p]
	\arrow["{\widetilde{\beta}}", shorten <=4pt, shorten >=4pt, Rightarrow, from=1p, to=3p]
\end{tikzcd}\]

We claim that we have the following:
\begin{align}\label{eqthetathetatilda}
\begin{split}
\theta \widetilde{\theta} &= \theta', \\
\widetilde{\psi} \psi &= \psi'.
\end{split}
\end{align}

The square $\beta'$ will then exhibit the equality $[s,t] = [\theta, 1]$:
\[\begin{tikzcd}
	{a'} && {a'} \\
	x && x && {a''} \\
	{a''} && x && {a''}
	\arrow[Rightarrow, no head, from=1-1, to=1-3]
	\arrow["s"', from=1-1, to=2-1]
	\arrow["\theta"', curve={height=24pt}, from=1-1, to=3-1]
	\arrow["s", from=1-3, to=2-3]
	\arrow[""{name=0, anchor=center, inner sep=0}, Rightarrow, no head, from=2-1, to=2-3]
	\arrow[""{name=0p, anchor=center, inner sep=0}, phantom, from=2-1, to=2-3, start anchor=center, end anchor=center]
	\arrow["{\theta'}"', from=2-1, to=3-1]
	\arrow["t", from=2-3, to=2-5]
	\arrow[Rightarrow, no head, from=2-3, to=3-3]
	\arrow[Rightarrow, no head, from=2-5, to=3-5]
	\arrow[""{name=1, anchor=center, inner sep=0}, "{\psi'}"', from=3-1, to=3-3]
	\arrow[""{name=1p, anchor=center, inner sep=0}, phantom, from=3-1, to=3-3, start anchor=center, end anchor=center]
	\arrow[curve={height=24pt}, Rightarrow, no head, from=3-1, to=3-5]
	\arrow["t"', from=3-3, to=3-5]
	\arrow["{\beta'}", shorten <=4pt, shorten >=4pt, Rightarrow, from=0p, to=1p]
\end{tikzcd}\]

\noindent To prove the claim that \eqref{eqthetathetatilda} holds, note that the corner $(se, m't)$ is jointly monic in $X$ and we have the following squares that relate the 1-cells in \eqref{eqthetathetatilda}:
\[\begin{tikzcd}
	x & x & x && x && x \\
	{a'} & {a'} & x \\
	{a''} & {a'} & x && {a''} && x
	\arrow[equals, from=1-1, to=1-2]
	\arrow["{\widetilde{\theta}}"', from=1-1, to=2-1]
	\arrow[""{name=0, anchor=center, inner sep=0}, equals, from=1-2, to=1-3]
	\arrow["{\widetilde{\theta}}"', from=1-2, to=2-2]
	\arrow[equals, from=1-3, to=2-3]
	\arrow[""{name=1, anchor=center, inner sep=0}, equals, from=1-5, to=1-7]
	\arrow["{\theta'}"', from=1-5, to=3-5]
	\arrow[equals, from=1-7, to=3-7]
	\arrow[""{name=2, anchor=center, inner sep=0}, equals, from=2-1, to=2-2]
	\arrow["\theta"', from=2-1, to=3-1]
	\arrow[""{name=3, anchor=center, inner sep=0}, "{\widetilde{\psi}}"', from=2-2, to=2-3]
	\arrow[equals, from=2-2, to=3-2]
	\arrow[equals, from=2-3, to=3-3]
	\arrow[""{name=4, anchor=center, inner sep=0}, "\psi"', from=3-1, to=3-2]
	\arrow["{\widetilde{\psi}}"', from=3-2, to=3-3]
	\arrow[""{name=5, anchor=center, inner sep=0}, "{\psi'}"', from=3-5, to=3-7]
	\arrow["{\widetilde{\beta}}"', shorten <=4pt, shorten >=4pt, Rightarrow, from=0, to=3]
	\arrow["{\beta'}", shorten <=9pt, shorten >=9pt, Rightarrow, from=1, to=5]
	\arrow["\beta", shorten <=4pt, shorten >=4pt, Rightarrow, from=2, to=4]
\end{tikzcd}\]

Thus, to show \eqref{eqthetathetatilda} it suffices to show:
\begin{align*}
\theta \widetilde{\theta} se &= \theta' se, \\
m't \widetilde{\psi} \psi &= m't \psi'.
\end{align*}
The first equality holds because:
\[
\theta \widetilde{\theta} s e = \theta e = e' = \theta' s e,
\]
while the second equality holds because:
\[
m't\widetilde{\psi} \psi = m \psi = m' = m' t \psi'.
\]
\end{proof}

\noindent We therefore propose the following terminology:

\begin{defi}
By an (orthogonal) \textit{factorization double category} we mean a double category $X$ with the following properties:
\begin{itemize}
\item $X$ is invariant,
\item every top-right corner in $X$ can be filled into a square and every square is bicartesian,
\item every bottom-left corner in $X$ is jointly monic.
\end{itemize}
Denote by $\factdbl$ the full subcategory of $\text{Dbl}$ consisting of factorization double categories.
\end{defi}

\begin{pozn}
Any factorization double category is automatically flat: Given two squares $\alpha, \lambda$ with the same boundary, by Lemma \ref{LEMMA_bicartdblopcart} there exists a unique square $\beta$ as follows:
\[\begin{tikzcd}
	a && a && b && a &&&& b \\
	d && d && d & {=} \\
	d && d && d && d &&&& c
	\arrow[""{name=0, anchor=center, inner sep=0}, "f", from=1-3, to=1-5]
	\arrow["g", from=1-5, to=2-5]
	\arrow["u"', from=1-3, to=2-3]
	\arrow[""{name=1, anchor=center, inner sep=0}, "v"', from=2-3, to=2-5]
	\arrow[""{name=2, anchor=center, inner sep=0}, "\psi"', from=3-1, to=3-3]
	\arrow["\theta"', from=2-1, to=3-1]
	\arrow[""{name=3, anchor=center, inner sep=0}, "f", from=1-7, to=1-11]
	\arrow["g", from=1-11, to=3-11]
	\arrow["u"', from=1-7, to=3-7]
	\arrow[""{name=4, anchor=center, inner sep=0}, "v"', from=3-7, to=3-11]
	\arrow[""{name=5, anchor=center, inner sep=0}, Rightarrow, no head, from=2-1, to=2-3]
	\arrow[Rightarrow, no head, from=2-3, to=3-3]
	\arrow[Rightarrow, no head, from=1-1, to=1-3]
	\arrow[Rightarrow, no head, from=2-5, to=3-5]
	\arrow["u"', from=1-1, to=2-1]
	\arrow["v"', from=3-3, to=3-5]
	\arrow["v"', curve={height=18pt}, from=3-1, to=3-5]
	\arrow["u"', curve={height=18pt}, from=1-1, to=3-1]
	\arrow["\lambda", shorten <=4pt, shorten >=4pt, Rightarrow, from=0, to=1]
	\arrow["\alpha", shorten <=9pt, shorten >=9pt, Rightarrow, from=3, to=4]
	\arrow["{\exists ! \beta}", shorten <=4pt, shorten >=4pt, Rightarrow, from=5, to=2]
\end{tikzcd}\]
Now in $X^{v}$ the corner $(u,v)$ is jointly monic and we have $\theta \circ u = 1_d \circ u$, $v \circ \psi = v \circ 1_d$. Thus $\theta = 1_d$ and $\psi = 1_d$. Invariance now forces $\beta = \Box_d$, the identity square on the morphism $1_d$ and so $\alpha = \lambda$.
\end{pozn}

\noindent Combining Lemmas \ref{lemma_isossubcxcm}, \ref{lemma_isossupcxcm}, \ref{lemma_orthog} we obtain:

\begin{prop}\label{THM_dbl_to_ofs}
Let $X$ be a factorization double category. Then the two classes of vertical and horizontal corners, $(\ce_X,\cm_X)$, form an orthogonal factorization system on the category $\cnr{X}$.
\end{prop}

\begin{pr}[Partial maps]

Let $\cc$ be a category with pullbacks and consider the double category $\text{MPbSq}(\cc)^{v}$. It is obviously flat and invariant with every square bicartesian. We have seen that top-left corners in its vertical dual are jointly monic in Example \ref{prmpbsq_pbsqsjointmonicity}. $\text{MPbSq}(\cc)^{v}$ is thus a factorization double category.

Combined with the description of the category of corners from Example \ref{prparc} and Proposition \ref{THM_dbl_to_ofs} we obtain that the category $\text{Par}(\cc)$ of objects and partial maps admits an orthogonal factorization system given by vertical corners followed by horizontal ones:
\[\begin{tikzcd}
	A &&&& A \\
	{A'} && {A'} && A && B
	\arrow[Rightarrow, no head, from=2-1, to=2-3]
	\arrow["\iota", hook, from=2-1, to=1-1]
	\arrow["f"', from=2-5, to=2-7]
	\arrow[Rightarrow, no head, from=1-5, to=2-5]
\end{tikzcd}\]
In \cite{restrcats} these are called the \textit{domains} and \textit{total maps} in $\cc$.
\end{pr}

\begin{pr}[Categories and cofunctors]
If $\ce$ is a category with pullbacks, the double category $\text{BOFib}(\ce)^{v}$ is a factorization double category. In the category\\ $\cnr{\text{BOFib}(\ce)^{v}} = \text{Cof}(\ce)$ every morphism can be factored as (the opposite of) a bijection on objects functor followed by a discrete opfibration (as mentioned in \cite[Theorem 18]{bryce_cof}). By the results in this section, these classes form an orthogonal factorization system on $\text{Cof}(\ce)$.
\end{pr}

\begin{pr}\label{pr_xp_fib}
Given a fibration $P: \ce \to \cb$, there is an associated sub-double category $X_P \subseteq \sq{\ce}$ whose vertical morphisms are $P$-vertical morphisms (those that are sent to isomorphisms by $P$), horizontal morphisms are cartesian lifts of morphisms in $\cb$, and squares are commutative squares in $\ce$.
\begin{prop}
$X_P$ is a factorization double category.
\end{prop}
\begin{proof}
Invariance is straightforward. To show joint monicity, assume we are given the data as in Definition \ref{defi_jointly_monic}. Note then that from the existence of squares $\kappa_1, \kappa_2$ in $(X^P)^v$ it follows that $P\theta = (P\psi)^{-1}$, and from $\theta \pi_1 = \theta' \pi_1$ we have $P\theta = P\theta'$. Since $\pi_2$ is a cartesian lift, the following picture forces $\psi \circ \theta' = 1$ and thus $\psi = \psi'$:
\[\begin{tikzcd}
	& {a'} && P && {Pa'} \\
	{a'} && b & \mapsto & {Pa'} && Pb
	\arrow["{\pi_2}", from=1-2, to=2-3]
	\arrow["{\pi_2}"', from=2-1, to=2-3]
	\arrow["{P\pi_2}", from=1-6, to=2-7]
	\arrow["{P\pi_2}"', from=2-5, to=2-7]
	\arrow[Rightarrow, no head, from=2-5, to=1-6]
	\arrow["{\psi \circ \theta'}", from=2-1, to=1-2]
\end{tikzcd}\]
Given a top-right corner $\lambda, u$, the bicartesian filler square is given by the cartesian lift of the pair $(P\lambda, b')$ and the unique canonical comparison morphism:
\[\begin{tikzcd}
	a && b && a & b & {b'} \\
	{a'} && {b'} && Pa && Pb
	\arrow[""{name=0, anchor=center, inner sep=0}, "\lambda", from=1-1, to=1-3]
	\arrow["u", from=1-3, to=2-3]
	\arrow[""{name=1, anchor=center, inner sep=0}, "{\lambda_{P\lambda, b}}"', from=2-1, to=2-3]
	\arrow["{\exists !}"', from=1-1, to=2-1]
	\arrow["\lambda", from=1-5, to=1-6]
	\arrow["u", from=1-6, to=1-7]
	\arrow[""{name=2, anchor=center, inner sep=0}, "P\lambda"', from=2-5, to=2-7]
	\arrow[shorten <=4pt, shorten >=4pt, Rightarrow, from=0, to=1]
	\arrow["\mapsto"{marking}, draw=none, from=1-6, to=2]
\end{tikzcd}\]
\end{proof}
The category of corners $\cnr{X_P}$ is isomorphic to the category $\ce$ via the functor sending an equivalence class $[u,\lambda]$ to the composite $\lambda \circ u$.

From the results in this section we obtain that the category $\ce$ admits an orthogonal factorization system given by the class of $P$-vertical morphisms followed by cartesian morphisms. Factorization systems associated to fibrations are a special case of \textit{simple reflective factorization systems} associated to prefibrations and have been studied in \cite{factorizationOFSrosicky}.
\end{pr}

\begin{nonpr}[Spans]

Let $\cc$ be a category with pullbacks and consider the double category $\text{PbSq}(\cc)^{v}$ of (opposite) pullback squares. This is not a factorization double category because not every top-left corner in $\text{PbSq}(\cc)$ is jointly monic (see Example \ref{prsq_jointmonicity}). We can not thus use Theorem \ref{THM_dbl_to_ofs} to obtain an orthogonal factorization system on $\cnr{\text{PbSq}(\cc)^{v}} = \text{Span}(\cc)$ and in fact, the two canonical classes of spans in $\spn(\cc)$ do not form one\footnote{They do form a weak factorization system as we have seen in Example \ref{pr_spansWFS}}.

To see this, let $\cc = \set$, denote by $\text{sw}: 2 \to 2$ the non-identity automorphism of the two-element set. Note that the class $[sw,1]: 2 \to 2$ in $\text{Span}$ is not the identity morphism. Consider now the span $(!,!): * \to *$:
\[\begin{tikzcd}
	{*} \\
	2 && {*}
	\arrow["{!}", from=2-1, to=1-1]
	\arrow["{!}", from=2-1, to=2-3]
\end{tikzcd}\]

\noindent As both $[sw,1]$ and $[1,1]$ in the place of the dotted line make the diagram below commute, we obtain that the factorization is not unique up to a unique isomorphism and thus the classes are not orthogonal.
\[\begin{tikzcd}
	{*} && 2 && {*} \\
	{*} && 2 && {*}
	\arrow[dashed, from=1-3, to=2-3]
	\arrow["{[!,1]}"', from=1-3, to=1-1]
	\arrow["{[!,1]}", from=2-3, to=2-1]
	\arrow["{[1,!]}"', from=2-3, to=2-5]
	\arrow["{[1,!]}", from=1-3, to=1-5]
	\arrow[Rightarrow, no head, from=1-1, to=2-1]
	\arrow[Rightarrow, no head, from=1-5, to=2-5]
\end{tikzcd}\]

\end{nonpr}
\color{black}
\bigskip

Recall now the assignment $(\ce,\cm) \mapsto D_{\ce,\cm}$ from Construction \ref{conDcecm}. We have:

\begin{prop}\label{lemma_ofs_to_dbl}
Let $(\ce,\cm)$ be an orthogonal factorization system on a category $\cc$. Then $D_{\ce,\cm}$ is a factorization double category. The assignment $(\ce,\cm) \to D_{\ce,\cm}$ induces a functor $\ofs \to \factdbl$.
\end{prop}
\begin{proof}
Let us verify each point:
\begin{itemize}

\item \textbf{Invariance}: We show the horizontal invariance, the vertical invariance is done similarly. Because the classes $\ce,\cm$ are closed under composition, given $u \in \ce$ and two isomorphisms $\theta, \psi$, the composite $\psi^{-1} u \theta$ gives the unique square with the given boundary:
\[\begin{tikzcd}
	a && b \\
	\\
	c && d
	\arrow[""{name=0, anchor=center, inner sep=0}, "{\theta \cong}", from=1-1, to=1-3]
	\arrow["u", from=1-3, to=3-3]
	\arrow[""{name=1, anchor=center, inner sep=0}, "\psi\cong"', from=3-1, to=3-3]
	\arrow[dashed, from=1-1, to=3-1]
	\arrow[shorten <=9pt, shorten >=9pt, Rightarrow, from=0, to=1]
\end{tikzcd}\]

\item \textbf{Filling corners into squares, every square bicartesian}: Consider a top-right corner as pictured below:
\[\begin{tikzcd}
	a && b \\
	\\
	{a'} && c
	\arrow[""{name=0, anchor=center, inner sep=0}, "{m'}", from=1-1, to=1-3]
	\arrow["{e'}", from=1-3, to=3-3]
	\arrow["e"', dashed, from=1-1, to=3-1]
	\arrow[""{name=1, anchor=center, inner sep=0}, "m"', dashed, from=3-1, to=3-3]
	\arrow[shorten <=9pt, shorten >=9pt, Rightarrow, from=0, to=1]
\end{tikzcd}\]
The filler square is given by the $(\ce,\cm)$-factorization of the morphism $e'm'$ in $\cc$. Next, let there be a square as pictured below left:
\[\begin{tikzcd}
	a && b && a && b \\
	{a'} && c && {a''} && c
	\arrow[""{name=0, anchor=center, inner sep=0}, "g", from=1-1, to=1-3]
	\arrow[""{name=0p, anchor=center, inner sep=0}, phantom, from=1-1, to=1-3, start anchor=center, end anchor=center]
	\arrow["e"', from=1-1, to=2-1]
	\arrow["u", from=1-3, to=2-3]
	\arrow[""{name=1, anchor=center, inner sep=0}, "g", from=1-5, to=1-7]
	\arrow[""{name=1p, anchor=center, inner sep=0}, phantom, from=1-5, to=1-7, start anchor=center, end anchor=center]
	\arrow["{e'}"', from=1-5, to=2-5]
	\arrow["u", from=1-7, to=2-7]
	\arrow[""{name=2, anchor=center, inner sep=0}, "m"', from=2-1, to=2-3]
	\arrow[""{name=2p, anchor=center, inner sep=0}, phantom, from=2-1, to=2-3, start anchor=center, end anchor=center]
	\arrow[""{name=3, anchor=center, inner sep=0}, "{m'}"', from=2-5, to=2-7]
	\arrow[""{name=3p, anchor=center, inner sep=0}, phantom, from=2-5, to=2-7, start anchor=center, end anchor=center]
	\arrow[shorten <=4pt, shorten >=4pt, Rightarrow, from=0p, to=2p]
	\arrow[shorten <=4pt, shorten >=4pt, Rightarrow, from=1p, to=3p]
\end{tikzcd}\]
We wish to show that it is bicartesian. Assume there is another square (pictured above right) with the same top-right corner.

Because $me = m'e'$ are two factorizations of the same morphism, there is a unique isomorphism $\theta \in \ce \cap \cm$ between them (pictured below left). It then gives a comparison square between the first square and the second square, as pictured below right:
\[\begin{tikzcd}[column sep=scriptsize]
	& a &&& a && b && a && b \\
	{a''} && {a'} && {a'} && c & {=} \\
	& c &&& {a''} && c && {a''} && c
	\arrow["{e'}"', from=1-2, to=2-1]
	\arrow["e", from=1-2, to=2-3]
	\arrow[""{name=0, anchor=center, inner sep=0}, "g", from=1-5, to=1-7]
	\arrow[""{name=0p, anchor=center, inner sep=0}, phantom, from=1-5, to=1-7, start anchor=center, end anchor=center]
	\arrow["e"', from=1-5, to=2-5]
	\arrow["u", from=1-7, to=2-7]
	\arrow[""{name=1, anchor=center, inner sep=0}, "g", from=1-9, to=1-11]
	\arrow[""{name=1p, anchor=center, inner sep=0}, phantom, from=1-9, to=1-11, start anchor=center, end anchor=center]
	\arrow["{e'}"', from=1-9, to=3-9]
	\arrow["u", from=1-11, to=3-11]
	\arrow["{m'}"', from=2-1, to=3-2]
	\arrow["{\cong \theta}"{description}, from=2-3, to=2-1]
	\arrow["m", from=2-3, to=3-2]
	\arrow[""{name=2, anchor=center, inner sep=0}, "m"{description}, from=2-5, to=2-7]
	\arrow[""{name=2p, anchor=center, inner sep=0}, phantom, from=2-5, to=2-7, start anchor=center, end anchor=center]
	\arrow[""{name=2p, anchor=center, inner sep=0}, phantom, from=2-5, to=2-7, start anchor=center, end anchor=center]
	\arrow["\theta"', from=2-5, to=3-5]
	\arrow[Rightarrow, no head, from=2-7, to=3-7]
	\arrow[""{name=3, anchor=center, inner sep=0}, "{m'}"', from=3-5, to=3-7]
	\arrow[""{name=3p, anchor=center, inner sep=0}, phantom, from=3-5, to=3-7, start anchor=center, end anchor=center]
	\arrow[""{name=4, anchor=center, inner sep=0}, "{m'}"', from=3-9, to=3-11]
	\arrow[""{name=4p, anchor=center, inner sep=0}, phantom, from=3-9, to=3-11, start anchor=center, end anchor=center]
	\arrow[shorten <=4pt, shorten >=4pt, Rightarrow, from=0p, to=2p]
	\arrow[shorten <=9pt, shorten >=9pt, Rightarrow, from=1p, to=4p]
	\arrow[shorten <=4pt, shorten >=4pt, Rightarrow, from=2p, to=3p]
\end{tikzcd}\]
This gives the \textbf{existence}. To prove the \textbf{uniqueness}, assume that there is a different comparison square. Its vertical domain map then gives the morphism of factorizations $(e,m),(e',m')$ and is thus forced to be equal to $\theta$. Thus the square is opcartesian in $D_{\ce,\cm}$. The proof that it is also opcartesian in ${(D_{\ce,\cm})}^*$ is done the same way.

\item \textbf{Joint monicity}: Let $(\pi_1, \pi_2)$, $\kappa_1$, $\kappa_2$ be the data in $D_{\ce,\cm}$ as pictured below:
\[\begin{tikzcd}[column sep=scriptsize]
	a &&&& {a'} && {a'} && {a'} && {a'} \\
	{a'} && b && {a''} && {a'} && {a''} && {a'}
	\arrow["{\pi_2}", from=2-1, to=2-3]
	\arrow["{\pi_1}", from=1-1, to=2-1]
	\arrow[""{name=0, anchor=center, inner sep=0}, "{\psi'}"', from=2-9, to=2-11]
	\arrow["\theta"', from=1-5, to=2-5]
	\arrow[""{name=1, anchor=center, inner sep=0}, "\psi"', from=2-5, to=2-7]
	\arrow[Rightarrow, no head, from=2-7, to=1-7]
	\arrow[""{name=2, anchor=center, inner sep=0}, Rightarrow, no head, from=1-5, to=1-7]
	\arrow["{\theta'}"', from=1-9, to=2-9]
	\arrow[""{name=3, anchor=center, inner sep=0}, Rightarrow, no head, from=1-9, to=1-11]
	\arrow[Rightarrow, no head, from=2-11, to=1-11]
	\arrow["{\kappa_2}"', shorten <=4pt, shorten >=4pt, Rightarrow, from=3, to=0]
	\arrow["{\kappa_1}"', shorten <=4pt, shorten >=4pt, Rightarrow, from=2, to=1]
\end{tikzcd}\]

Assume $\theta \pi_1 = \theta' \pi_1 $ and $\pi_2 \psi = \pi_2 \psi'$. We have the following:
\begin{align*}
\psi' \theta \pi_1 &= \pi_1,\\
\pi_2 &= \pi_2 \psi' \theta.
\end{align*}

In other words, $\psi' \theta$ is an endomorphism of the factorization:
\[\begin{tikzcd}
	a && {a'} && b \\
	a && {a'} && b
	\arrow["{\pi_1}", from=1-1, to=1-3]
	\arrow["{\pi_2}", from=1-3, to=1-5]
	\arrow["{\pi_1}"', from=2-1, to=2-3]
	\arrow["{\pi_2}"', from=2-3, to=2-5]
	\arrow[Rightarrow, no head, from=1-1, to=2-1]
	\arrow[Rightarrow, no head, from=1-5, to=2-5]
	\arrow["{\psi'\circ \theta}", from=1-3, to=2-3]
\end{tikzcd}\]
From orthogonality we get $\psi' \theta = 1$. Since $\psi \circ \theta$ is an $(\ce,\cm)$-factorization of the identity, both of $\psi$ and $\theta$ are isomorphisms together with $\psi' \theta = 1$ we get:
\[
\psi = \theta^{-1} = \psi'.
\]
Applying inverses, we get $\theta = \theta'$.
\end{itemize}
Thus $D_{\ce,\cm}$ is a factorization double category. The functoriality is straightforward.
\end{proof}

Analogous to Theorem \ref{thmekvivalence_sfs_coddiscr}, we have:
\begin{theo}\label{thm_equiv_ofs_factdbl}
The functor $\cnr{-}: \factdbl \to \ofs$ is an equivalence of categories, with the equivalence inverse being the functor $D: \ofs \to \factdbl$.
\end{theo}
\begin{proof}
We will again show that there are natural isomorphisms $1 \cong D \circ \cnr{-}$ and $\cnr{-} \circ D \cong 1$:

To see that $1 \cong D \circ \cnr{-}$, let $X$ be a factorization double category.
Consider again the identity-on-objects functor $X_0 \to \ce_X$ that sends the morphism $e \mapsto [e,1]$, we would like to show that it is an isomorphism. To see the fullness, consider a corner $[u,g]$ for which $g$ is invertible. By invariance, we obtain a unique square as below left that exhibits the equality $[\tau u, 1 ] = [u,g]$. To see faithfulness, assume we have $[e,1] = [e',1]$, then there is a 2-cell as pictured below right:
\[\begin{tikzcd}
	a &&& a \\
	b & c && {a'} & {a'} \\
	c & c && {a'} & {a'}
	\arrow["u"', from=1-1, to=2-1]
	\arrow["{\tau \circ u}"', curve={height=24pt}, from=1-1, to=3-1]
	\arrow["e", from=1-4, to=2-4]
	\arrow["{e'}"', curve={height=18pt}, from=1-4, to=3-4]
	\arrow[""{name=0, anchor=center, inner sep=0}, "g", from=2-1, to=2-2]
	\arrow["\tau"', from=2-1, to=3-1]
	\arrow[equals, from=2-2, to=3-2]
	\arrow[""{name=1, anchor=center, inner sep=0}, equals, from=2-4, to=2-5]
	\arrow["\theta"', from=2-4, to=3-4]
	\arrow[equals, from=2-5, to=3-5]
	\arrow[""{name=2, anchor=center, inner sep=0}, equals, from=3-1, to=3-2]
	\arrow[""{name=3, anchor=center, inner sep=0}, equals, from=3-4, to=3-5]
	\arrow[shorten <=4pt, shorten >=4pt, Rightarrow, from=0, to=2]
	\arrow[shorten <=4pt, shorten >=4pt, Rightarrow, from=1, to=3]
\end{tikzcd}\]
By invariance, the square is forced to be the identity and thus $e = e'$. Analogously there is an isomorphism $h(X) \cong \cm_X$.  Define a double functor $X \to D_{\ce_X, \cm_X}$ so that it is identity on objects, on vertical morphisms sends $e \mapsto [e,1]$ and on horizontal morphisms sends $m \mapsto [1,m]$. Because both double categories are flat, to show that it is well-defined on cells and an isomorphism, it is enough to prove the following:
\[\begin{tikzcd}[column sep=scriptsize]
	a && b &&& a && b \\
	&&& {\text{in } X} & \iff &&&& {\text{in }D_{\ce_X,\cm_X}} \\
	c && d &&& c && d
	\arrow[""{name=0, anchor=center, inner sep=0}, "m", from=1-1, to=1-3]
	\arrow["{e'}"', from=1-1, to=3-1]
	\arrow["e", from=1-3, to=3-3]
	\arrow[""{name=1, anchor=center, inner sep=0}, "{[1,m]}", from=1-6, to=1-8]
	\arrow["{[e',1]}"', from=1-6, to=3-6]
	\arrow["{[e,1]}", from=1-8, to=3-8]
	\arrow[""{name=2, anchor=center, inner sep=0}, "{m'}"', from=3-1, to=3-3]
	\arrow[""{name=3, anchor=center, inner sep=0}, "{[1,m']}"', from=3-6, to=3-8]
	\arrow["\exists"', shorten <=9pt, shorten >=9pt, Rightarrow, from=0, to=2]
	\arrow["\exists"', shorten <=9pt, shorten >=9pt, Rightarrow, from=1, to=3]
\end{tikzcd}\]

The direction ``$\Rightarrow$” follows already from Proposition \ref{prop_naturalityinsquares}. For the direction ``$\Leftarrow$” assume the right square commutes. If we denote $[e,1] \circ [1,m] = [u,g]$, we obtain the square required above left as the following composite, where the upper square is the square used for the composition of the corners, and the lower square exists from the equality $[u,g] = [e',m']$:
\[\begin{tikzcd}
	a && b \\
	c && d \\
	c && d
	\arrow[""{name=0, anchor=center, inner sep=0}, "m", from=1-1, to=1-3]
	\arrow[""{name=0p, anchor=center, inner sep=0}, phantom, from=1-1, to=1-3, start anchor=center, end anchor=center]
	\arrow["u"', from=1-1, to=2-1]
	\arrow["{e'}"', curve={height=24pt}, from=1-1, to=3-1]
	\arrow["e", from=1-3, to=2-3]
	\arrow[""{name=1, anchor=center, inner sep=0}, "g"{description}, from=2-1, to=2-3]
	\arrow[""{name=1p, anchor=center, inner sep=0}, phantom, from=2-1, to=2-3, start anchor=center, end anchor=center]
	\arrow[""{name=1p, anchor=center, inner sep=0}, phantom, from=2-1, to=2-3, start anchor=center, end anchor=center]
	\arrow["\theta"', from=2-1, to=3-1]
	\arrow[Rightarrow, no head, from=2-3, to=3-3]
	\arrow[""{name=2, anchor=center, inner sep=0}, "{m'}"', from=3-1, to=3-3]
	\arrow[""{name=2p, anchor=center, inner sep=0}, phantom, from=3-1, to=3-3, start anchor=center, end anchor=center]
	\arrow[shorten <=4pt, shorten >=4pt, Rightarrow, from=0p, to=1p]
	\arrow[shorten <=4pt, shorten >=4pt, Rightarrow, from=1p, to=2p]
\end{tikzcd}\]

To see that $\cnr{-} \circ D \cong 1$, let $(\ce,\cm)$ be an orthogonal factorization system on a category $\cc$ and define an identity-on-objects functor $F: \cnr{D_{\ce,\cm}} \to \cc$ so that it sends the class of corners $[e,m]$ to $m \circ e$. This is well defined because if $[e,m] = [e',m']$, there is a 2-cell between $(e,m), (e',m')$ in $D_{\ce,\cm}$:
\[\begin{tikzcd}
	a && a \\
	{a'} && {a'} && b \\
	{a''} && {a'} && b
	\arrow[Rightarrow, no head, from=1-1, to=1-3]
	\arrow["e"', from=1-1, to=2-1]
	\arrow["{e'}"', curve={height=24pt}, from=1-1, to=3-1]
	\arrow["e", from=1-3, to=2-3]
	\arrow[""{name=0, anchor=center, inner sep=0}, Rightarrow, no head, from=2-1, to=2-3]
	\arrow[""{name=0p, anchor=center, inner sep=0}, phantom, from=2-1, to=2-3, start anchor=center, end anchor=center]
	\arrow["\theta"', from=2-1, to=3-1]
	\arrow["m", from=2-3, to=2-5]
	\arrow[Rightarrow, no head, from=2-3, to=3-3]
	\arrow[Rightarrow, no head, from=2-5, to=3-5]
	\arrow[""{name=1, anchor=center, inner sep=0}, "\psi"', from=3-1, to=3-3]
	\arrow[""{name=1p, anchor=center, inner sep=0}, phantom, from=3-1, to=3-3, start anchor=center, end anchor=center]
	\arrow["{m'}"', curve={height=24pt}, from=3-1, to=3-5]
	\arrow["m"', from=3-3, to=3-5]
	\arrow[shorten <=4pt, shorten >=4pt, Rightarrow, from=0p, to=1p]
\end{tikzcd}\]
But since squares in this double category are commutative squares in $\cc$, we get: $me = m \psi \theta e = m'e'$. Functoriality of $F$ is straightforward, faithfulness follows from the fact that $(\ce,\cm)$-factorizations are unique up to a unique isomorphism, and fullness follows because every morphism in $\cc$ has an $(\ce,\cm)$-factorization. Thus $F$ is an isomorphism.
\end{proof}

\begin{pozn}
Let us point out that a more conceptual point of view for why factorization systems are equivalent to certain double categories is offered in the following chapter, in particular in Theorem \ref{THM_fundamentaladj} of Section \ref{SEC_codobjs_doublecats}.  Therein we will also show that there is an equivalence between \textbf{weak} factorization systems and certain double categories, although we have not identified which double-categorical properties such double categories must satisfy.
\end{pozn}

%% file: PhD_CH_4_COMPUTATION.tex
\chapter[Turning lax structures into strict ones]{Turning lax structures into strict ones using codescent objects}\label{CH_COMPUTATION}
\section{Introduction}\label{SEC_intro_computation}

In this chapter we will compute codescent objects of colax coherence data.

The \textbf{first goal} is to put the results from Chapter \ref{CH_FS_DBLCATS} into the broader perspective of 2-category theory -- the category of corners associated to a crossed double category $X$ is in fact its codescent object if $X$ is regarded as a strict coherence data in $\cat$ (Proposition \ref{THMP_cnr_is_cod_obj}). This approach allows us to generalize the equivalences from Chapter \ref{CH_FS_DBLCATS} to wider classes of factorization systems.

The \textbf{second goal} is to show that for a class of 2-monads of form $\cat(T)$ on a 2-category $\cat(\ce)$ (for $\ce$ with pullbacks), a coherence theorem for colax algebras (\ref{TERMI_lax_coherence_holds}) holds: for every colax $T$-algebra $\mathbb{A}$, there is a strict $T$-algebra $\mathbb{A}'$ and an adjunction in $\cat(\ce)$ of their underlying categories:
\[\begin{tikzcd}
	A && {A'}
	\arrow[""{name=0, anchor=center, inner sep=0}, curve={height=24pt}, from=1-1, to=1-3]
	\arrow[""{name=1, anchor=center, inner sep=0}, curve={height=24pt}, from=1-3, to=1-1]
	\arrow["\dashv"{anchor=center, rotate=-90}, draw=none, from=1, to=0]
\end{tikzcd}\]
At the same time we will provide an explicit description for $\mathbb{A}'$ for various examples including colax monoidal categories -- it will be given by a more general variant of the \textit{category of corners} construction (a concept introduced by Mark Weber \cite{internalalg}, an instance of which we have already seen in Chapter \ref{CH_FS_DBLCATS}).

The algebra $\mathbb{A}'$ is given by the codescent object of the resolution $\res{\mathbb{A}}$ of $\mathbb{A}$. Since $\cat(T)$ is cartesian, in case $\mathbb{A}$ is strict, $\res{\mathbb{A}}$ is an internal category in $\cat(\ce)$ (Remark \ref{POZN_Tcartesian_resAa_internal_cat}) and so the codescent object can be computed by essentially using the methods of \cite{internalalg} (and we have done so in \cite{milo1}).

The problem is that when one considers colax algebras, the resolution is no longer an internal category since the required simplicial identities between 1-cells do not hold (there are non-identity 2-cells present instead). $\res{\mathbb{A}}$ does however still admit a lot of structure resembling an internal category -- in this thesis we will call them \textit{codomain-colax categories}. For $\ce = \set$ this is a ``double category”-like structure consisting of objects, vertical, horizontal morphisms, and squares, that can be composed both horizontally and vertically. The difference from double categories is that the composition of horizontal arrows does not preserve the codomain: in the picture below, the codomain of $g_2 \circ g_1$ will not in general agree with that of $g_2$. There will however always be a comparison vertical arrow between them:
\[\begin{tikzcd}
	a &&&& {\widehat{c}} \\
	\\
	a && b && c
	\arrow["{g_2 \circ g_1}", from=1-1, to=1-5]
	\arrow["\gamma", from=1-5, to=3-5]
	\arrow["{g_1}"', from=3-1, to=3-3]
	\arrow["{g_2}"', from=3-3, to=3-5]
\end{tikzcd}\]

Assuming $\ce$ is a category with pullbacks, the 2-category $\cat(\ce)$ automatically admits certain 2-limits (Corollary \ref{THM_catE_limits}). What the author finds curious is that it automatically admits certain 2-colimits -- codescent objects of special kinds of coherence data, as they can be built using pullbacks in $\ce$. Bourke has observed that $\cat(\ce)$ always admits codescent objects of \textit{cateads}\footnote{A special kind of internal category in $\cat(\ce)$. We briefly mention the case for $\ce = \set$ in Definition \ref{DEFI_catead} and afterwards.}, see \cite[Theorem 3.65]{johnphd}. Weber has observed in  \cite[5.4.7 Corollary]{internalalg} that $\cat(\ce)$ admits codescent objects of internal categories in which the domain functor is a discrete fibration. In this chapter we will extend this last result to codomain-colax categories with the same property.

\medskip

\noindent \textbf{The chapter is structured as follows}:

\begin{itemize}
\item In Section \ref{SEC_codobjs_doublecats} we notice that when the coherence data $X$ is actually a category in $\cat$ (a double category of \ref{SEC_double_cats}), each cocone can be equivalently described as a pair of functors from the horizontal and vertical categories of $X$. In the case of crossed double categories, we show that the codescent object is precisely the category of corners $\cnr{X}$ studied in the previous chapter. We end the section with a generalization of the results from sections \ref{subsekce_sfs} and \ref{subsekce_ofs} -- we show that both of these are a restriction of a more general adjunction between double categories and certain pairs of functors with a common codomain.
\item In Section \ref{SEC_codomain_colax_cats} we introduce the aforementioned \textit{codomain-colax categories} -- a special kind of coherence data that resembles a double category, except the target of a composition of two horizontal arrows does not have to be equal to the target of the latter of the two. In Subsection \ref{SEC_thecaseckCatE} we give a formula for the codescent object of a special kind of a codomain-colax category in $\cat(\ce)$ using pullbacks in $\ce$. In Subsection \ref{SEC_thecaseckCat} we give the formula for the codescent object of a general codomain-colax category for the case $\ce = \set$ in terms of generators and relations.

\item In Section \ref{SEC_application_lax_coh} we use the results of the previous subsection to immediately deduce that the assumptions of Theorem \ref{THM_lax_coherence_colax_algebras_EXPOSITION} are satisfied for 2-monads of form $\cat(T)$, thus proving that coherence for colax algebras holds for this class of 2-monads. We demonstrate the constructions on colax monoidal categories, lax functors between 2-categories, lax double functors between double categories and colax functors $\cj \to \cat$ (for $\cj$ a 1-category).
\end{itemize}

\section{Codescent objects of double categories}\label{SEC_codobjs_doublecats}

\noindent \textbf{The section is organized as follows}:
\begin{itemize}
\item Subsection \ref{SUBS_cocone_of_doublecat}, we show that in case the coherence data $X$ is an internal category to $\cat(\ce)$, a cocone can equivalently be described as a pair of functors agreeing on objects satisfying a ``naturality” condition.

\item Subsection \ref{SUBS_COMPUTATION_Examples} gives examples of codescent objects of crossed double categories and cateads.

\item In Subsection \ref{SECS_fundamental_adj} we prove that there is an adjunction between double categories and certain pairs of functors that restricts to the equivalences proven in Theorems \ref{thmekvivalence_sfs_coddiscr} and \ref{thm_equiv_ofs_factdbl}.
\end{itemize}

\subsection{Cocone of a double category}\label{SUBS_cocone_of_doublecat}

\begin{defi}
Let $\ce$ be a category with pullbacks. An \textit{internal double category} in $\ce$ is a category internal to $\cat(\ce)$.
\end{defi}

\begin{con}\label{CON_pair}
Let $X$ be an internal double category in a category $\ce$ with pullbacks. Define a 1-functor $\pair(X,-): \cat(\ce) \to \set$ as follows: the set $\pair(X,\cc)$ consists of pairs of (internal) functors $(F: X_0 \to \cc, \xi: h(X) \to \cc)$ that agree on objects and satisfy the following \textit{naturality condition}:
\[
\comp^\cc \circ (\xi_1 \circ s, F_1 \circ d^1_0) = \comp^\cc \circ (F_1 \circ d^1_1, \xi_1 \circ t),
\]
where $s,t$ are the source/target 1-cells for the internal category $X_1$ and $d_0, d_1$ are the domain/codomain internal functors $X_1 \to X_0$.

Given an internal functor $H: \cc \to \cc'$, the function:
\[
\pair(X,H): \pair(X,\cc) \to \pair(X,\cc'),
\]
is the assignment $(F,\xi) \mapsto (HF, H\xi)$.
\end{con}

\begin{pozn}[$\ce = \set$]
For a double category $X$, the naturality condition for a pair $(F: X_0 \to \cc,\xi: h(X) \to \cc)$ says that:
\[\begin{tikzcd}
	a & b \\
	c & d & {\text{ exists in }X} & \Rightarrow & {Fv \circ \xi(g) = \xi(h) \circ Fu \hspace{0.5cm}\text{ in } \cc}
	\arrow[""{name=0, anchor=center, inner sep=0}, "g", from=1-1, to=1-2]
	\arrow["u"', from=1-1, to=2-1]
	\arrow["v", from=1-2, to=2-2]
	\arrow[""{name=1, anchor=center, inner sep=0}, "h"', from=2-1, to=2-2]
	\arrow["{\exists \alpha}", shorten <=4pt, shorten >=4pt, Rightarrow, from=0, to=1]
\end{tikzcd}\]
\end{pozn}

\begin{prop}\label{THMP_cohX_cong_pairX}
Let $X$ be an internal double category in a category $\ce$ with pullbacks. There is a natural isomorphism of functors:
\[
\coh{\cat(\ce)}{X}{-} \cong \pair(X,-): \cat(\ce) \to \set.
\]
\end{prop}
\begin{proof}
A cocone $(F,\xi)$ for $X$ with apex $\cy$ consists of a functor $F:X_0 \to \cy$ and a natural transformation $\xi: Fd_1 \Rightarrow Fd_0: X_1 \to \cy$. In other words, it is a data consisting of three 1-cells in $\ce$: $(F_0, F_1, \xi)$. The set of such triples of 1-cells is in bijection with tuples $(F,\widetilde{\xi})$ of graph morphisms $X_0 \to \cc$, $h(X) \to \cc$ that agree on objects. Now, it is easily seen that the cocone axioms for $(F,\xi)$ hold if and only if $\widetilde{\xi}$ satisfies the functor axioms. Also, the 1-cell $\xi$ satisfies the axioms for a natural transformation if and only if the pair $(F,\widetilde{\xi})$ satisfies the naturality condition in Construction \ref{CON_pair}.
\end{proof}

From now on we will thus not distinguish between cocones for (an internal) double category $X$ and pairs of functors $(F: X_0 \to \cc, \xi: h(X) \to \cc)$ agreeing on objects satisfying the naturality condition. Likewise, this pair being the codescent object of $X$ is equivalent to requiring that for any other such pair $(G: X_0 \to \cd,\psi: h(X) \to \cd)$, there exists a unique functor making the following commute:
\[\begin{tikzcd}
	&& \cc \\
	{X_0} &&&& {h(X)} \\
	&& \cd
	\arrow[dashed, from=1-3, to=3-3]
	\arrow["F", from=2-1, to=1-3]
	\arrow["G"', from=2-1, to=3-3]
	\arrow["\xi"', from=2-5, to=1-3]
	\arrow["\psi", from=2-5, to=3-3]
\end{tikzcd}\]
This one condition suffices, since the two-dimensional universal property follows automatically from the one-dimensional one in $\cat(\ce)$ by Proposition \ref{THMP_powers_with2_implies_1dim_impl_2dim} and Corollary \ref{THM_catE_limits}.

We will denote the codescent object of a double category $X$ by $\cod{X}$. Analogous to Definition \ref{DEFI_dblcat_duals}, to every double category $X$ internal to $\ce$ we may also associate its \textit{transpose} $X^T$. We have:

\begin{prop}[Invariance under transposition]\label{THMP_invariance_under_transposition}
Let $X$ be a double category in a category $\ce$ with pullbacks. We have: $\cod{X} \cong \cod{X^T}$.
\end{prop}
\begin{proof}
It is enough to show that there is a natural bijection between the sets of $W$-weighted codescent cocones:
\[
\text{Cocone}(X,\cc) \cong \text{Cocone}(X^T, \cc)
\]
The bijection is given by $(F,\xi) \mapsto (\xi,F)$.
\end{proof}

\subsection{Examples}\label{SUBS_COMPUTATION_Examples}

Recall the category of corners $\cnr{X}$ associated to a crossed double category $X$ from Construction \ref{concnrpbsqs}. We have:

\begin{prop}\label{THMP_cnr_is_cod_obj}
Let $X$ be a crossed double category. Then the pair of functors $(F: X_0 \to \cnr{X},\xi: h(X) \to \cnr{X})$ sending $u \mapsto [u,1]$, $g \mapsto [1,g]$ is the codescent object of $X$.
\end{prop}
\begin{proof}
The naturality condition has already been verified in Proposition \ref{prop_naturalityinsquares}. Let now $(G: X_0 \to Y, \psi: h(X) \to Y)$ be a different cocone. We have to show that there is a unique functor $\theta: \cnr{X} \to Y$ commuting with the functors:
\[\begin{tikzcd}
	&& {\cnr{X}} \\
	{X_0} &&&& {h(X)} \\
	&& Y
	\arrow["F", from=2-1, to=1-3]
	\arrow["\xi"', from=2-5, to=1-3]
	\arrow["G"', from=2-1, to=3-3]
	\arrow["\psi", from=2-5, to=3-3]
	\arrow[dashed, from=1-3, to=3-3]
\end{tikzcd}\]
The above commutativity forces $\theta(a) = Fa = \xi(a)$ for an object $a \in \cnr{X}$, and on morphisms forces $\theta([u,g]) = \psi(g) \circ Gu$. It is routine to verify that this mapping is well defined and a functor.
\end{proof}

\begin{defi}\label{DEFI_catead}
A \textit{catead} $X$ in $\cat$ is a double category with the property that:
\[\begin{tikzcd}
	& x & y && x & z \\
	{\forall(u,h):} & z & y & {\text{and }\forall(g,v):} & x & y
	\arrow[""{name=0, anchor=center, inner sep=0}, dashed, from=1-2, to=1-3]
	\arrow["u"', from=1-2, to=2-2]
	\arrow[Rightarrow, no head, from=1-3, to=2-3]
	\arrow[""{name=1, anchor=center, inner sep=0}, "g", from=1-5, to=1-6]
	\arrow[Rightarrow, no head, from=1-5, to=2-5]
	\arrow["v", from=1-6, to=2-6]
	\arrow[""{name=2, anchor=center, inner sep=0}, "h"', from=2-2, to=2-3]
	\arrow[""{name=3, anchor=center, inner sep=0}, dashed, from=2-5, to=2-6]
	\arrow["{\exists!}", shorten <=4pt, shorten >=4pt, Rightarrow, dashed, from=0, to=2]
	\arrow["{\exists!}", shorten <=4pt, shorten >=4pt, Rightarrow, dashed, from=1, to=3]
\end{tikzcd}\]
and moreover for every square $\alpha$ there exists a unique pair $(\kappa_1, \kappa_2)$ of squares of the following form whose composite equals $\alpha$:
\[\begin{tikzcd}
	a && b && a && b \\
	&&& {=} & a && d \\
	c && d && c && d
	\arrow[""{name=0, anchor=center, inner sep=0}, "g", from=1-1, to=1-3]
	\arrow["u"', from=1-1, to=3-1]
	\arrow["v", from=1-3, to=3-3]
	\arrow[""{name=1, anchor=center, inner sep=0}, "g", from=1-5, to=1-7]
	\arrow[Rightarrow, no head, from=1-5, to=2-5]
	\arrow["v", from=1-7, to=2-7]
	\arrow[""{name=2, anchor=center, inner sep=0}, dashed, from=2-5, to=2-7]
	\arrow["u"', from=2-5, to=3-5]
	\arrow[Rightarrow, no head, from=2-7, to=3-7]
	\arrow[""{name=3, anchor=center, inner sep=0}, "h"', from=3-1, to=3-3]
	\arrow[""{name=4, anchor=center, inner sep=0}, "h"', from=3-5, to=3-7]
	\arrow["\alpha", shorten <=9pt, shorten >=9pt, Rightarrow, from=0, to=3]
	\arrow["{\kappa_1}", shorten <=4pt, shorten >=4pt, Rightarrow, dashed, from=1, to=2]
	\arrow["{\kappa_2}", shorten <=4pt, shorten >=4pt, Rightarrow, dashed, from=2, to=4]
\end{tikzcd}\]
\end{defi}

\noindent The following has been observed in \cite[Remark 2.87]{johnphd}. Let us offer an alternative proof exploiting our description of cocones for a double category:

\begin{prop}\label{THMP_codobj_of_catead_in_CAT}
Let $X$ be a catead in $\cat$. Then: $\cod{X} = h(X)$.
\end{prop}
\begin{proof}
All that needs to be proven is that there is a natural bijection between cocones $(F,\xi)$ for $X$ with apex $\cc \in \cat$ and functors $h(X) \to \cc$.

Let $(F,\xi)$ be a cocone. For any vertical $u: x \to z$, there exists a unique square like this:
\[\begin{tikzcd}
	x & z \\
	z & z
	\arrow[""{name=0, anchor=center, inner sep=0}, "{\rho_u}", dashed, from=1-1, to=1-2]
	\arrow["u"', from=1-1, to=2-1]
	\arrow[Rightarrow, no head, from=1-2, to=2-2]
	\arrow[""{name=1, anchor=center, inner sep=0}, Rightarrow, no head, from=2-1, to=2-2]
	\arrow["{\exists !}", shorten <=4pt, shorten >=4pt, Rightarrow, dashed, from=0, to=1]
\end{tikzcd}\]
The naturality in squares of $X$ now forces $Fu \sstackrel{!}{=} \xi(\rho_u)$, so $F$ is uniquely determined by $\xi: h(X) \to \cc$. Conversely, given a functor $\psi: h(X) \to \cc$, define the assignment $G: X_0 \to \cc$ by agreeing with $\psi$ on objects and sending $Gu := \psi(\rho_u)$. That this assignment preserves the identities is obvious. To see that $G(vu) = G(v)G(u)$, we will show that $\rho_{vu} = \rho_v \circ \rho_u$. This fact follows from the equality of the squares below:
\[\begin{tikzcd}
	x && a && x & z & a \\
	z &&& {=} & z & z & a \\
	a && a && a & a & a
	\arrow[""{name=0, anchor=center, inner sep=0}, "{\rho_{vu}}", from=1-1, to=1-3]
	\arrow["u"', from=1-1, to=2-1]
	\arrow[Rightarrow, no head, from=1-3, to=3-3]
	\arrow[""{name=1, anchor=center, inner sep=0}, "{\rho_u}", from=1-5, to=1-6]
	\arrow["u"', from=1-5, to=2-5]
	\arrow["{\rho_v}", from=1-6, to=1-7]
	\arrow[Rightarrow, no head, from=1-6, to=2-6]
	\arrow[Rightarrow, no head, from=1-7, to=2-7]
	\arrow["v"', from=2-1, to=3-1]
	\arrow[""{name=2, anchor=center, inner sep=0}, Rightarrow, no head, from=2-5, to=2-6]
	\arrow["v"', from=2-5, to=3-5]
	\arrow[""{name=3, anchor=center, inner sep=0}, "{\rho_v}", from=2-6, to=2-7]
	\arrow["v"', from=2-6, to=3-6]
	\arrow[Rightarrow, no head, from=2-7, to=3-7]
	\arrow[""{name=4, anchor=center, inner sep=0}, Rightarrow, no head, from=3-1, to=3-3]
	\arrow[Rightarrow, no head, from=3-5, to=3-6]
	\arrow[""{name=5, anchor=center, inner sep=0}, Rightarrow, no head, from=3-6, to=3-7]
	\arrow[shorten <=9pt, shorten >=9pt, Rightarrow, from=0, to=4]
	\arrow[shorten <=4pt, shorten >=4pt, Rightarrow, from=1, to=2]
	\arrow[shorten <=4pt, shorten >=4pt, Rightarrow, from=3, to=5]
\end{tikzcd}\]

To prove the naturality in squares, we first prove it for the squares of the following forms:
\[\begin{tikzcd}
	{a)} & a & d & {b)} & a & b \\
	& c & d && a & d
	\arrow[""{name=0, anchor=center, inner sep=0}, "g", from=1-2, to=1-3]
	\arrow["u"', from=1-2, to=2-2]
	\arrow[equals, from=1-3, to=2-3]
	\arrow[""{name=1, anchor=center, inner sep=0}, "g", from=1-5, to=1-6]
	\arrow[equals, from=1-5, to=2-5]
	\arrow["v", from=1-6, to=2-6]
	\arrow[""{name=2, anchor=center, inner sep=0}, "h"', from=2-2, to=2-3]
	\arrow[""{name=3, anchor=center, inner sep=0}, "h"', from=2-5, to=2-6]
	\arrow[shorten <=4pt, shorten >=4pt, Rightarrow, from=0, to=2]
	\arrow[shorten <=4pt, shorten >=4pt, Rightarrow, from=1, to=3]
\end{tikzcd}\]
In \textbf{a)}, we have $\psi(g) = \psi(h) \circ Gu$ since $g = h \circ \rho_u$, which follows from the equality of squares below:
\[\begin{tikzcd}
	a & d && a & c & d \\
	c & d && c & c & d
	\arrow[""{name=0, anchor=center, inner sep=0}, "g", from=1-1, to=1-2]
	\arrow["u"', from=1-1, to=2-1]
	\arrow[Rightarrow, no head, from=1-2, to=2-2]
	\arrow[""{name=1, anchor=center, inner sep=0}, "{\rho_u}", from=1-4, to=1-5]
	\arrow["u"', from=1-4, to=2-4]
	\arrow["h", from=1-5, to=1-6]
	\arrow[Rightarrow, no head, from=1-5, to=2-5]
	\arrow[Rightarrow, no head, from=1-6, to=2-6]
	\arrow[""{name=2, anchor=center, inner sep=0}, "h"', from=2-1, to=2-2]
	\arrow[""{name=3, anchor=center, inner sep=0}, Rightarrow, no head, from=2-4, to=2-5]
	\arrow["h"', from=2-5, to=2-6]
	\arrow[shorten <=4pt, shorten >=4pt, Rightarrow, from=0, to=2]
	\arrow[shorten <=4pt, shorten >=4pt, Rightarrow, from=1, to=3]
\end{tikzcd}\]
Now, \textbf{b)} is proven in a similar way so we omit the proof. Finally, the general case follows since every square factors as a square \textbf{b)} followed by a square \textbf{a)}.
\end{proof}

\begin{pr}\label{PR_sqcc_cod_obj}
Consider the double category of commutative squares $\sq{\cc}$ associated to a category $\cc$. This is easily seen to be a catead, so we obtain $\cod{\sq{\cc}} = \cc$.
\end{pr}

\begin{pozn}\label{POZN_codobj_general_doublecat}
The codescent object of a general double category $X$ can be described in terms of generators and relations -- this will be done in Construction \ref{CON_Cnr_pro_codomaincolaxcat}.
\end{pozn}

\subsection{Fundamental adjunction}\label{SECS_fundamental_adj}

With the framework of codescent objects we can generalize the correspondences from Section \ref{SEC_fs_dbl_cats}.

\begin{defi}
Denote by $\cocone$ the subcategory of $\cat^{\bullet \rightarrow \bullet \leftarrow \bullet}$ that consists of cospans of functors:
\[\begin{tikzcd}
	\cx && \cc && \cy
	\arrow["F", from=1-1, to=1-3]
	\arrow["\xi"', from=1-5, to=1-3]
\end{tikzcd}\]
such that $\ob \cx = \ob \cc = \ob \cy$ and $F, \xi$ are identity-on-object functors.

\noindent A morphism $(F,\xi) \to (G,\psi)$ in this category is a triple of functors $(V,\theta,H)$ so that the following commutes in $\cat$:
\[\begin{tikzcd}
	\cx && \cc && \cy \\
	\\
	{\cx'} && \cd && {\cy'}
	\arrow["F", from=1-1, to=1-3]
	\arrow["V"', from=1-1, to=3-1]
	\arrow["\theta", from=1-3, to=3-3]
	\arrow["\xi"', from=1-5, to=1-3]
	\arrow["H", from=1-5, to=3-5]
	\arrow["G"', from=3-1, to=3-3]
	\arrow["\psi", from=3-5, to=3-3]
\end{tikzcd}\]
We refer to the objects of this category as \textit{cocones} in $\cat$. Given such a cocone $(F,\xi)$, we will call the category $\cc$ its \textit{apex}.
\end{defi}

\begin{con}
Similar to Construction \ref{conDcecm}, assume that we are given a cocone $(F: \cx \to \cc,\xi: \cy \to \cc)$ and define a double category $D_{F,\xi}$ such that:
\begin{itemize}
\item Objects are the objects of $\cc$,
\item vertical morphisms are the morphisms of $\cx$,
\item horizontal morphisms are the morphisms of $\cy$,
\item the bottom left square is a square in $D_{F,\xi}$ if and only if the bottom right diagram commutes in $\cc$:
\[\begin{tikzcd}
	a && b && a && b \\
	&&& {} &&&& {} \\
	c && d && c && d
	\arrow[""{name=0, anchor=center, inner sep=0}, "h", from=1-1, to=1-3]
	\arrow["{v'}", from=1-3, to=3-3]
	\arrow["v"', from=1-1, to=3-1]
	\arrow[""{name=1, anchor=center, inner sep=0}, "{h'}"', from=3-1, to=3-3]
	\arrow["Fv"', from=1-5, to=3-5]
	\arrow["{\xi_h}", from=1-5, to=1-7]
	\arrow["{Fv'}", from=1-7, to=3-7]
	\arrow["{\xi_{h'}}"', from=3-5, to=3-7]
	\arrow[shorten <=9pt, shorten >=9pt, Rightarrow, from=0, to=1]
\end{tikzcd}\]
\end{itemize}
Given a morphism $(V,\theta,H): (F,\xi) \to (G,\psi)$ between cocones with apexes $\cc, \cd$, there is an induced double functor $D_{V,\theta,H}$ whose vertical part is given by $V$ and horizontal part given by $H$. To see that it is well-defined on squares, we need to show that:
\[
Fv' \circ \xi_h = \xi_{h'} \circ Fv\hspace{1.5cm}\text{in } \cc,
\]
implies:
\[
GVv \circ \psi H g = \psi H h \circ GVu \hspace{1cm}\text{in } \cd.
\]
But this is the case since by the definition of the morphism $(V,\theta,H)$, both sides are equal to the following:
\[
\theta F v \circ \theta \xi g = \theta (F v \circ \xi g) = \theta (\xi h \circ F u) = \theta \xi h \circ \theta F u.
\]
\end{con}

\noindent Denote by $\cod{-}: \dbl \to \cocone$ the functor sending the double category $X$ to the pair of functors corresponding to its codescent cocone (see Lemma \ref{THMP_cohX_cong_pairX}).

\begin{theo}\label{THM_fundamentaladj}
There is an adjunction:
\[\begin{tikzcd}
	{\text{Cocone}(\cat)\phantom{....}} &&& \dbl
	\arrow[""{name=0, anchor=center, inner sep=0}, "D"', curve={height=30pt}, from=1-1, to=1-4]
	\arrow[""{name=1, anchor=center, inner sep=0}, "{\cod{-}}"', curve={height=30pt}, from=1-4, to=1-1]
	\arrow["\dashv"{anchor=center, rotate=-90}, draw=none, from=1, to=0]
\end{tikzcd}\]
\end{theo}
\begin{proof}
Let us define the unit $\eta_X: X \to D_{\cod{X}}$ evaluated at a double category $X$ to be the unique double functor that is the identity on objects, vertical and horizontal morphisms. To prove the universal property of the unit, we need to show that given a double functor $H: X \to D_{G,\psi}$, where $(G: \cy \to \cc, \psi: \cz \to \cc)$ is a cocone, there exists a unique morphism $(V,\theta, S): (F,\xi) \to (G,\psi)$ from the colimit cocone that makes the following commute:
\[\begin{tikzcd}
	&& {D_{G,\psi}} \\
	\\
	X && {D_{\cod{X}}}
	\arrow["H", from=3-1, to=1-3]
	\arrow["{\eta_X}"', from=3-1, to=3-3]
	\arrow["{D(V,\theta,S)}"', dashed, from=3-3, to=1-3]
\end{tikzcd}\]
This equation forces $V = H_0$ and $S = h(X)$. Moreover, the requirement that $(V,\theta, S)$ is a morphism in $\cocone$ makes the pair $(GH_0, \psi h(H))$ into a cocone for $X$. This implies that there exists a unique functor $\theta: \cod{X} \to \cc$ making the diagram below commute, proving the universal property:
\[\begin{tikzcd}
	{X_0} && {\cod{X}} && {h(X)} \\
	\\
	\cy && \cc && \cz
	\arrow["F", from=1-1, to=1-3]
	\arrow["{H_0}"', from=1-1, to=3-1]
	\arrow["{\exists ! \theta}", dashed, from=1-3, to=3-3]
	\arrow["\xi"', from=1-5, to=1-3]
	\arrow["{h(H)}", from=1-5, to=3-5]
	\arrow["G"', from=3-1, to=3-3]
	\arrow["\psi", from=3-5, to=3-3]
\end{tikzcd}\]
\end{proof}

\begin{pozn}
Appropriately restricting the adjunction above, one obtains the equivalences from Theorem \ref{thmekvivalence_sfs_coddiscr} and Theorem \ref{thm_equiv_ofs_factdbl}.
\end{pozn}

\noindent The adjunction in Theorem \ref{THM_fundamentaladj} can also be seen to restrict to the equivalence of weak factorization systems and certain double categories as we will now show. Denote by $\cw \cf \cs$ the full sub-category of $\cocone$ spanned by triples $(F: \ce \to \cc, \xi: \cy \to \cm)$, where $F, \xi$ are category inclusions and the pair $(\ce, \cm)$ is a \textit{weak factorization systems} on $\cc$. We have:

\begin{prop}
Let $(\ce,\cm)$ be a weak factorization system on a category $\cc$. Then $(\ce,\cm)$ is the codescent object of the double category $D_{\ce,\cm}$:
\[
\cod{D_{\ce,\cm}} \cong (\ce,\cm).
\]
\end{prop}
\begin{proof}
Let $(G,\psi)$ be a cocone with apex $\cd$. We want to show that there is a unique comparison functor $\theta: \cc \to \cd$ making the following commute:
\[\begin{tikzcd}
	\ce && \cc && \cm \\
	\\
	&& \cd
	\arrow[hook, from=1-1, to=1-3]
	\arrow["G"', from=1-1, to=3-3]
	\arrow["{\exists ! \theta}"{description}, dashed, from=1-3, to=3-3]
	\arrow[hook', from=1-5, to=1-3]
	\arrow["\psi", from=1-5, to=3-3]
\end{tikzcd}\]
The commutativity forces $\theta$ to send:
\begin{alignat*}{2}
e &\mapsto Ge &\hspace{0.5cm}& \text{if }e \in \mor \ce,\\
m &\mapsto \psi m && \text{if }m \in \mor \cm.
\end{alignat*}
As any morphism $f$ in $\cc$ factors as $m_f \circ e_f$ with $e_f \in \ce$, $m_f \in \cm$, we are forced to put:
\[
\theta f \sstackrel{!}{=} \psi m_f \circ G e_f.
\]
To see that this is well-defined, take two $(\ce,\cm)$-factorizations of a morphism in $\cc$ as pictured in the square below:
\[\begin{tikzcd}
	A && {A''} \\
	\\
	{A'} && B
	\arrow["{e'}", from=1-1, to=1-3]
	\arrow["e"', from=1-1, to=3-1]
	\arrow["{m'}", from=1-3, to=3-3]
	\arrow["{\exists \tau}", dashed, from=3-1, to=1-3]
	\arrow["m"', from=3-1, to=3-3]
\end{tikzcd}\]
By the weak orthogonality, there exists a diagonal filler $\tau$. Consider the $(\ce,\cm)$-factoriza\-tion of $\tau$ and notice that the following are squares in $D_{\ce,\cm}$:
\[\begin{tikzcd}
	A && A && {A'} && B \\
	{A'} \\
	X && {A''} && X & {A''} & B
	\arrow[""{name=0, anchor=center, inner sep=0}, Rightarrow, no head, from=1-1, to=1-3]
	\arrow[""{name=0p, anchor=center, inner sep=0}, phantom, from=1-1, to=1-3, start anchor=center, end anchor=center]
	\arrow["e"', from=1-1, to=2-1]
	\arrow["{e'}", from=1-3, to=3-3]
	\arrow[""{name=1, anchor=center, inner sep=0}, "m", from=1-5, to=1-7]
	\arrow[""{name=1p, anchor=center, inner sep=0}, phantom, from=1-5, to=1-7, start anchor=center, end anchor=center]
	\arrow["{e_\tau}"', from=1-5, to=3-5]
	\arrow[Rightarrow, no head, from=1-7, to=3-7]
	\arrow["{e_\tau}"', from=2-1, to=3-1]
	\arrow[""{name=2, anchor=center, inner sep=0}, "{m_\tau}"', from=3-1, to=3-3]
	\arrow[""{name=2p, anchor=center, inner sep=0}, phantom, from=3-1, to=3-3, start anchor=center, end anchor=center]
	\arrow["{m_\tau}"', from=3-5, to=3-6]
	\arrow["{m'}"', from=3-6, to=3-7]
	\arrow[shorten <=9pt, shorten >=9pt, Rightarrow, from=0p, to=2p]
	\arrow[shorten <=7pt, shorten >=4pt, Rightarrow, from=1p, to=3-6]
\end{tikzcd}\]
In particular:
\[
\theta (me) = \psi m \circ Ge \sstackrel{(A)}{=} \psi m' \circ \psi m_\tau \circ Ge_\tau \circ Ge \sstackrel{(B)}{=} \psi m' \circ Ge' = \theta (m'e'),
\]
where $(A)$ is the naturality of the cocone $(G,\psi)$ in the second portrayed square and $(B)$ is the naturality in the first square -- thus the assignment $\theta$ is well-defined. It also clearly preserves the identities. To see that it preserves compositions, consider a pair $(f,g)$ as below:
\[\begin{tikzcd}
	A && B && C \\
	& {X_f} && {X_g} \\
	&& Y \\
	&& {X_{gf}}
	\arrow["f", from=1-1, to=1-3]
	\arrow["{e_f}"', from=1-1, to=2-2]
	\arrow["{e_{gf}}"', curve={height=30pt}, from=1-1, to=4-3]
	\arrow["g", from=1-3, to=1-5]
	\arrow["{e_g}"', from=1-3, to=2-4]
	\arrow["{m_f}"', from=2-2, to=1-3]
	\arrow["e"', from=2-2, to=3-3]
	\arrow["{m_g}"', from=2-4, to=1-5]
	\arrow["m"', from=3-3, to=2-4]
	\arrow["\tau", from=3-3, to=4-3]
	\arrow["{m_{gf}}"', curve={height=30pt}, from=4-3, to=1-5]
\end{tikzcd}\]
We have:
\begin{align*}
\theta(g) \circ \theta (f) &= \psi m_g  \circ Ge_g  \circ \psi m_f  \circ Ge_f \\
&= \psi m_g  \circ \psi m  \circ Ge  \circ Ge_f \\
&\sstackrel{(*)}{=} \psi m_{gf}  \circ \psi m_\tau  \circ Ge_\tau  \circ Ge  \circ Ge_f \\
&\sstackrel{(*)}{=} \psi m_{gf}  \circ Ge_{gf}\\
&= \theta(gf),
\end{align*}
where $(*)$'s again follow from the existence of appropriate squares in $D_{\ce,\cm}$.
\end{proof}

\begin{lemma}
Consider an adjunction:
\[\begin{tikzcd}
	{(\epsilon,\eta):} & \cc && \cd
	\arrow[""{name=0, anchor=center, inner sep=0}, "F"', curve={height=24pt}, from=1-2, to=1-4]
	\arrow[""{name=1, anchor=center, inner sep=0}, "U"', curve={height=24pt}, from=1-4, to=1-2]
	\arrow["\dashv"{anchor=center, rotate=90}, draw=none, from=0, to=1]
\end{tikzcd}\]
And assume that $\cd' \subseteq \cd$ is a full subcategory spanned by objects $D \in \ob \cd$ for which the component $\epsilon_D$ is an isomorphism. Then there is a full subcategory $\cc' \subseteq \cc$ for which the adjunction restricts to an equivalence:
\[
\cc' \simeq \cd'.
\]
\end{lemma}
\begin{proof}
Notice first that $\cd'$ is \textit{replete} -- if $D' \in \ob \cd'$ and $\tau: D' \cong K$ is an isomorphism in $\cd$, then $K \in \ob \cd'$. This follows from the fact that in the naturality square for $\epsilon$, the morphisms $\epsilon_{D'}, \tau, FU\tau$ are isomorphisms -- this implies that $\epsilon_K$ is as well:
\[\begin{tikzcd}
	{FUD'} && {D'} \\
	\\
	FUK && K
	\arrow["{\epsilon_{D'}}", from=1-1, to=1-3]
	\arrow["{FU\tau}"', from=1-1, to=3-1]
	\arrow["\tau", from=1-3, to=3-3]
	\arrow["{\epsilon_K}"', from=3-1, to=3-3]
\end{tikzcd}\]

Define $\cc'$ to be the full subcategory spanned by the objects of form $UD'$, $D' \in \ob \cd'$. Notice that $F: \cc \to \cd$ restricts to a functor $F: \cc' \to \cd'$ since if $UD' \in \ob \cc'$, we have $FUD' \in \ob \cd'$ because there is an isomorphism $FUD' \cong D'$. Next, if $UD' \in \ob \cc'$, the morphism $\eta_{UD'}$ is an isomorphism since by the second triangle identity for $F \dashv U$ we have $\eta_{UD'} = U(\epsilon_{D'})^{-1}$.
\end{proof}

\begin{cor}\label{THMC_wfs_equiv_to_dbl}
The adjunction in Theorem \ref{THM_fundamentaladj} restricts to an equivalence between weak factorization systems and certain double categories.
\end{cor}

\noindent As of the time of writing this thesis, we do not have an elementary characterization available for such double categories.

\begin{pozn}
The adjunction in Theorem \ref{THM_fundamentaladj} is moreover \textbf{idempotent}, as we will now show. Notice first that the counit $\epsilon: \cod{D} \Rightarrow 1_{\text{Cocone}(\cat)}$ evaluated at a cocone $(G,\psi)$ is the unique comparison 1-cell from the colimit:
\[\begin{tikzcd}
	{X_0} && {\cod{X}} && {h(X)} \\
	\\
	{X_0} && \cc && {h(X)}
	\arrow["F", from=1-1, to=1-3]
	\arrow[Rightarrow, no head, from=1-1, to=3-1]
	\arrow["{\exists !}", dashed, from=1-3, to=3-3]
	\arrow["\xi"', from=1-5, to=1-3]
	\arrow[Rightarrow, no head, from=1-5, to=3-5]
	\arrow["G"', from=3-1, to=3-3]
	\arrow["\psi", from=3-5, to=3-3]
\end{tikzcd}\]
$D$ sends the morphism $\epsilon_{(G,\psi)}$ to a double functor that is the identity on objects, vertical and horizontal morphisms. For the functors below:
\[\begin{tikzcd}
	{D_{(G,\psi)}} && {D_{\cod{D_{(G,\psi)}}}}
	\arrow["{\eta_{D_{(G,\psi)}}}"', shift right=2, from=1-1, to=1-3]
	\arrow["{D\epsilon_{(G,\psi)}}"', shift right=2, from=1-3, to=1-1]
\end{tikzcd}\]
not only is $D\epsilon_{(G,\psi)}$ retract of $\eta_{D_{(G,\psi)}}$, they are actually \textbf{mutually inverse}. This is simply because the unit $\eta_Y$ is also a double functor that is the identity on objects, vertical and horizontal morphisms (and also because of the fact that both double categories in question are flat). This proves the adjunction is idempotent.

It thus generates an idempotent monad $T$ on $\dbl$. A $T$-algebra has to be necessarily flat, but it is not the case that \textbf{any} flat double category is a $T$-algebra: consider the double category $X$ freely generated by the pinwheel double graph:
\[\begin{tikzcd}
	0 & 1 & {} & 2 \\
	& 3 & 4 & 5 \\
	6 & 7 & 8 \\
	9 & {} & 10 & 11
	\arrow[""{name=0, anchor=center, inner sep=0}, from=1-1, to=1-2]
	\arrow[""{name=0p, anchor=center, inner sep=0}, phantom, from=1-1, to=1-2, start anchor=center, end anchor=center]
	\arrow[from=1-1, to=3-1]
	\arrow[from=1-2, to=1-4]
	\arrow[from=1-2, to=2-2]
	\arrow[shorten <=2pt, Rightarrow, from=1-3, to=2-3]
	\arrow[from=1-4, to=2-4]
	\arrow[""{name=1, anchor=center, inner sep=0}, from=2-2, to=2-3]
	\arrow[""{name=1p, anchor=center, inner sep=0}, phantom, from=2-2, to=2-3, start anchor=center, end anchor=center]
	\arrow[from=2-2, to=3-2]
	\arrow[""{name=2, anchor=center, inner sep=0}, from=2-3, to=2-4]
	\arrow[""{name=2p, anchor=center, inner sep=0}, phantom, from=2-3, to=2-4, start anchor=center, end anchor=center]
	\arrow[from=2-3, to=3-3]
	\arrow[from=2-4, to=4-4]
	\arrow[""{name=3, anchor=center, inner sep=0}, from=3-1, to=3-2]
	\arrow[""{name=3p, anchor=center, inner sep=0}, phantom, from=3-1, to=3-2, start anchor=center, end anchor=center]
	\arrow[from=3-1, to=4-1]
	\arrow[""{name=4, anchor=center, inner sep=0}, from=3-2, to=3-3]
	\arrow[""{name=4p, anchor=center, inner sep=0}, phantom, from=3-2, to=3-3, start anchor=center, end anchor=center]
	\arrow[shorten >=3pt, Rightarrow, from=3-2, to=4-2]
	\arrow[from=3-3, to=4-3]
	\arrow[from=4-1, to=4-3]
	\arrow[""{name=5, anchor=center, inner sep=0}, from=4-3, to=4-4]
	\arrow[""{name=5p, anchor=center, inner sep=0}, phantom, from=4-3, to=4-4, start anchor=center, end anchor=center]
	\arrow[shorten <=9pt, shorten >=9pt, Rightarrow, from=0p, to=3p]
	\arrow[shorten <=4pt, shorten >=4pt, Rightarrow, from=1p, to=4p]
	\arrow[shorten <=9pt, shorten >=9pt, Rightarrow, from=2p, to=5p]
\end{tikzcd}\]
It can not be isomorphic to $D_{\cnr{X}}$ since the latter contains an additional square -- the big "composite" square with corners $0, 2, 9, 11$. Presumably, given a flat double category $X$, $TX$ is the universal flat double category associated to it in which such arrangements as the one above can be evaluated.
\end{pozn}

\begin{cor}\label{THMC_adj_cat_dbl}
There is a reflection:
\[\begin{tikzcd}
	\cat &&& \dbl
	\arrow[""{name=0, anchor=center, inner sep=0}, "{\cod{-}}"', curve={height=30pt}, from=1-4, to=1-1]
	\arrow[""{name=1, anchor=center, inner sep=0}, "{\text{Sq}(-)}"', curve={height=30pt}, from=1-1, to=1-4]
	\arrow["\dashv"{anchor=center, rotate=-90}, draw=none, from=0, to=1]
\end{tikzcd}\]
\end{cor}
\begin{proof}
Note that there is an adjunction:
\[\begin{tikzcd}
	\cat &&& \coconecat
	\arrow[""{name=0, anchor=center, inner sep=0}, curve={height=30pt}, from=1-4, to=1-1]
	\arrow[""{name=1, anchor=center, inner sep=0}, curve={height=30pt}, from=1-1, to=1-4]
	\arrow["\dashv"{anchor=center, rotate=-90}, draw=none, from=0, to=1]
\end{tikzcd}\]
where the right adjoint sends a category $\ca$ to the identity cospan $(1_\ca, 1_\ca)$ and the left adjoint sends a cospan to the common codomain. Combining this adjunction with the one from Theorem \ref{THM_fundamentaladj} gives the result. Moreover, the counit of the adjunction is an isomorphism, see Example \ref{PR_sqcc_cod_obj}.
\end{proof}

\begin{pozn}
The above corollary is an instance of a more general phenomenon involving diagonal 2-functors for weighted (co)limits, see \cite[Section 3.1]{johnphd}.
\end{pozn}

\begin{pozn}
The above corollary is also a part of a higher-dimensional picture. Recall the adjunction from Proposition \ref{THMP_Ddashvobdashvpi0} for the case $\ce = \set$:
\[\begin{tikzcd}
	\set &&& \cat
	\arrow[""{name=0, anchor=center, inner sep=0}, "{\pi_{0}}"', curve={height=30pt}, from=1-4, to=1-1]
	\arrow[""{name=1, anchor=center, inner sep=0}, "D"', curve={height=30pt}, from=1-1, to=1-4]
	\arrow["\dashv"{anchor=center, rotate=-90}, draw=none, from=0, to=1]
\end{tikzcd}\]
This induces a change-of-base adjunction between $\set\text{-}\cat$ and $\cat\text{-}\cat$, portrayed as the adjunction below left:
\[\begin{tikzcd}
	\cat &&& {2\text{-}\cat} &&& \dbl
	\arrow[""{name=0, anchor=center, inner sep=0}, "{\pi_{0*}}"', curve={height=30pt}, from=1-4, to=1-1]
	\arrow[""{name=1, anchor=center, inner sep=0}, "\iota"', curve={height=30pt}, from=1-1, to=1-4]
	\arrow[""{name=2, anchor=center, inner sep=0}, "{2\text{-}\cnr{X}}"', curve={height=30pt}, from=1-7, to=1-4]
	\arrow[""{name=3, anchor=center, inner sep=0}, "{\text{LaxSq(-)}}"', curve={height=30pt}, from=1-4, to=1-7]
	\arrow["\dashv"{anchor=center, rotate=-90}, draw=none, from=0, to=1]
	\arrow["\dashv"{anchor=center, rotate=-90}, draw=none, from=2, to=3]
\end{tikzcd}\]
Here $\iota$ regards a category as a locally discrete 2-category and $\pi_{0*}$ applies $\pi_0$ on each hom of a 2-category to produce a 1-category. For a 2-category $\ck$, $\text{LaxSq}(\ck)$ is a double category with vertical and horizontal morphisms those of $\ck$ and with a square:
\[\begin{tikzcd}
	a && b \\
	\\
	c && d
	\arrow[""{name=0, anchor=center, inner sep=0}, "h", from=1-1, to=1-3]
	\arrow["v"', from=1-1, to=3-1]
	\arrow["{v'}", from=1-3, to=3-3]
	\arrow[""{name=1, anchor=center, inner sep=0}, "{h'}"', from=3-1, to=3-3]
	\arrow["\alpha", shorten <=9pt, shorten >=9pt, Rightarrow, from=0, to=1]
\end{tikzcd}\]
being a 2-cell $\alpha : v'h \Rightarrow h'v$ in $\ck$. For a double category $X$, there exists a 2-category $2\text{-}\cnr{X}$ whose underlying 1-category consists of paths of vertical and horizontal morphisms of $X$, and whose 2-cells are generated by the squares of $X$:
\[\begin{tikzcd}
	a && b \\
	&&& {:(h,v') \Rightarrow (v,h')} \\
	c && d
	\arrow[""{name=0, anchor=center, inner sep=0}, "h", from=1-1, to=1-3]
	\arrow["v"', from=1-1, to=3-1]
	\arrow["{v'}", from=1-3, to=3-3]
	\arrow[""{name=1, anchor=center, inner sep=0}, "{h'}"', from=3-1, to=3-3]
	\arrow["\alpha", shorten <=9pt, shorten >=9pt, Rightarrow, from=0, to=1]
\end{tikzcd}\]
modulo relations given by the composition of squares and identity squares in $X$. These two functors together form an adjunction displayed above right. Composing these adjunctions, we obtain the adjunction from Corollary \ref{THMC_adj_cat_dbl}.
\end{pozn}

\begin{pozn}
Let us also compare the results in this thesis with the results of David Forsman in \cite{forsman}. In that paper, a strict factorization system $(\ce, \cm)$ on a category $\ce$ is called a \textit{variance} on $\cc$ if $(\cm, \ce)$ is also a strict factorization system on $\cc$.

It is easily observed that the adjunction Theorem \ref{THM_fundamentaladj} restricts to an equivalence between variances and double categories in which not just every top-right corner admits a unique filler, but every bottom-left corner admits a unique filler as well. \cite[Theorem 1.4]{forsman} is then a special case of Proposition \ref{THMP_cnr_is_cod_obj} applied to the double category $D_{\ce,\cm}$ associated to a variance $(\ce,\cm)$ on $\cc$.
\end{pozn}

\section{Codescent objects of colax coherence data}\label{SEC_codomain_colax_cats}

\noindent \textbf{The section is organized as follows}:
\begin{itemize}
\item In Subsection \ref{SUBS_codomain_colax_cats} we introduce codomain-colax categories, prove some of their properties, and mention the important classes of examples we will come across.
\item In Subsection \ref{SEC_thecaseckCatE} we prove that to any codomain-colax category $X$ in $\cat(\ce)$ (with an additional property) we may assign an internal category of corners $\cnr{X}$ (Construction \ref{CON_cnrX_for_X_domdiscr_codcolaxcat}) built using pullbacks of $\ce$, and that this is the codescent object of $X$ (Theorem \ref{THM_catE_admits_codobj_of_domdiscr_codlax}). Moreover, any 2-functor of form $\cat(T)$ automatically preserves these kinds of codescent objects (Theorem \ref{THM_CatT_preserves_special_codobjs})
\item In Subsection \ref{SUBS_grothendieck} we outline the relationship between codomain-colax categories and the Grothendieck construction.
\item In Subsection \ref{SEC_thecaseckCat} we give the formula for the codescent object of a codomain-colax category for the case $\ce = \cat$ in terms of generators and relations.
\end{itemize}

\subsection{Codomain-colax categories}\label{SUBS_codomain_colax_cats}

\begin{defi}
Let $\ck$ be a 2-category. By a \textbf{codomain-colax category} in $\ck$ we mean a diagram in $\ck$ of the following form:
\begin{center}
\begin{tikzpicture}
    \node (A) at (-1,0) {$X_2$};
    \node (B) at (2.5,0) {$X_1$};
    \node (C) at (6,0) {$X_0$};
    \node (D) at (-4.5,0) {$X_3$};
    \draw[transform canvas={yshift=3ex},->] (B) -- node[fill=white]{$j_1$} (A);
    \draw[transform canvas={yshift=6ex},->] (A) -- node[fill=white]{$d^\bullet_2$} (B);
    \draw[transform canvas={yshift=0ex},->] (A) --node[fill=white]{ $d^\bullet_1$}  (B);
    \draw[transform canvas={yshift=-3ex},->] (B) -- node[fill=white]{ $j_0$} (A);
    \draw[transform canvas={yshift=-6ex},->] (A) -- node[fill=white]{ $d^\bullet_0$} (B);

    \draw[transform canvas={yshift=3ex},->] (B) -- node[fill=white]{ $d_1$} (C);
    \draw[transform canvas={yshift=0ex},->] (C) --node[fill=white]{ $s$}  (B);
    \draw[transform canvas={yshift=-3ex},->] (B) -- node[fill=white]{ $d_0$} (C);
    
    \draw[transform canvas={yshift=4.5ex},->] (D) --node[fill=white]{ $d_3^{\bullet \bullet}$}  (A);
    \draw[transform canvas={yshift=1.5ex},->] (D) --node[fill=white]{ $d_2^{\bullet \bullet}$}  (A);
    \draw[transform canvas={yshift=-1.5ex},->] (D) --node[fill=white]{ $d_1^{\bullet \bullet}$}  (A);
    \draw[transform canvas={yshift=-4.5ex},->] (D) --node[fill=white]{ $d_0^{\bullet \bullet}$}  (A);
    
\node at (-1,1.5) {};
\node at (-1,1.5) {};

\node at (-1,-1.5) {};
\node at (-4.5,0) {};
\end{tikzpicture}
\end{center}
together with 2-cells:
\begin{align*}
\gamma :  d_0 d_1^\bullet &\Rightarrow d_0 d_0^\bullet, \\
\widehat{\gamma} : d_0^\bullet d_1^{\bullet \bullet} &\Rightarrow d_0^\bullet d_0^{\bullet \bullet}, \\
\iota : d_0 s &\Rightarrow 1_{X_0}, \\
\widehat{\iota} : d_0^\bullet j_0 &\Rightarrow 1_{X_1}.
\end{align*}

\noindent The first two  will be called the \textit{coassociators}, the second two the \textit{counitors}. The 1-cells in this data must satisfy the following ``simplicial identities”\footnote{These are all the usual simplicial identities except for the identities that would occur if $\gamma, \iota, \widehat{\gamma}, \widehat{\iota}$ were identity 2-cells.}:
\begin{align}
d_1 d_1^\bullet &= d_1 d_2^\bullet, \label{eqclmonadassoc}\\
d_1^\bullet d_1^{\bullet \bullet} &= d_1^\bullet d_2^{\bullet \bullet}, \label{eqclTmonadassoc}\\
d_2^\bullet d_2^{\bullet \bullet} &= d_2^\bullet d_3^{\bullet \bullet}, \label{eqclmonadassocT}
\end{align}
\begin{align}
d_1 s &= 1_{X_0}, \label{eqclmonadunit1LOWER}\\
d_1^\bullet j_0 &= 1_{X_1},\label{eqclmonadunit1}\\
d_2^\bullet j_1 &= 1_{X_1}, \label{eqclmonadunit2}\\
d_2 j_1 &= 1_{X_1},\label{eqclmonadunit3}
\end{align}
\begin{align}
\begin{split}\label{eqclunit_naturality}
s d_0 &= d_0^\bullet j_1,\\
j_0 s &= j_1 s.
\end{split}
\end{align}
\begin{align}
\begin{split}\label{eqclmultipl_naturality}
d_0 d_2^\bullet &= d_1 d_0^\bullet,\\
s d_1 &= d_2^\bullet j_0,\\
d_1^\bullet d_3^{\bullet \bullet} &= d_2^\bullet d_1^{\bullet \bullet},\\
d_0^\bullet d_2^{\bullet \bullet} &= d_1^\bullet d_0^{\bullet \bullet},\\
d_0^\bullet d_3^{\bullet \bullet} &= d_2^\bullet d_0^{\bullet \bullet}.
\end{split}
\end{align}

Next, the following three pairs of 2-cells are required to be equal (to ease on the notation, we omit the bullet symbols $\bullet$ in the upper index):
\begin{equation}\label{eqcllaxalgassoc}
\begin{tikzcd}[column sep=small]
	& {X_2} && {X_1} &&& {X_2} && {X_1} \\
	{X_3} && {X_2} && {X_0} & {X_3} && {X_1} && {X_0} \\
	& {X_2} && {X_1} &&& {X_2} && {X_1}
	\arrow["{d_1}", from=1-2, to=1-4]
	\arrow["{d_0}", from=1-4, to=2-5]
	\arrow["\gamma", shorten <=6pt, shorten >=6pt, Rightarrow, from=1-4, to=3-4]
	\arrow[""{name=0, anchor=center, inner sep=0}, "{d_1}", from=1-7, to=1-9]
	\arrow["{d_0}"', from=1-7, to=2-8]
	\arrow["{d_0}", from=1-9, to=2-10]
	\arrow["{d_2}", from=2-1, to=1-2]
	\arrow[""{name=1, anchor=center, inner sep=0}, "{d_1}", from=2-1, to=2-3]
	\arrow["{d_0}"', from=2-1, to=3-2]
	\arrow["{d_1}", from=2-3, to=1-4]
	\arrow["{d_0}"', from=2-3, to=3-4]
	\arrow["{d_2}", from=2-6, to=1-7]
	\arrow["{d_0}"', from=2-6, to=3-7]
	\arrow["{d_0}"{description}, from=2-8, to=2-10]
	\arrow["{d_0}"', from=3-2, to=3-4]
	\arrow["{d_0}"', from=3-4, to=2-5]
	\arrow["{d_1}", from=3-7, to=2-8]
	\arrow[""{name=2, anchor=center, inner sep=0}, "{d_0}"', from=3-7, to=3-9]
	\arrow["{d_0}"', from=3-9, to=2-10]
	\arrow["\gamma", shift left=5, shorten <=3pt, Rightarrow, from=0, to=2-8]
	\arrow["{\widehat{\gamma}}", shift left=5, shorten <=3pt, Rightarrow, from=1, to=3-2]
	\arrow["\gamma", shorten >=3pt, Rightarrow, from=2-8, to=2]
\end{tikzcd}
\end{equation}
\begin{equation}\label{eqcl2-natur1}
\begin{tikzcd}[column sep=small]
	& {X_1} && {X_2} &&& {X_1} && {X_2} \\
	{X_0} && {X_1} && {X_1} & {X_0} && {X_0} && {X_1}
	\arrow["{j_1}", from=1-2, to=1-4]
	\arrow["{d_0}", from=1-4, to=2-5]
	\arrow["{j_1}", from=1-7, to=1-9]
	\arrow["{d_0}", from=1-7, to=2-8]
	\arrow["{d_0}", from=1-9, to=2-10]
	\arrow["s", from=2-1, to=1-2]
	\arrow["s"', from=2-1, to=2-3]
	\arrow["{j_0}", from=2-3, to=1-4]
	\arrow[""{name=0, anchor=center, inner sep=0}, Rightarrow, no head, from=2-3, to=2-5]
	\arrow["s", from=2-6, to=1-7]
	\arrow[""{name=1, anchor=center, inner sep=0}, Rightarrow, no head, from=2-6, to=2-8]
	\arrow["s"', from=2-8, to=2-10]
	\arrow["{\widehat{\iota}}"', shorten >=3pt, Rightarrow, from=1-4, to=0]
	\arrow["\iota"', shorten >=3pt, Rightarrow, from=1-7, to=1]
\end{tikzcd}
\end{equation}
\begin{equation}\label{eqcl2-natur2}
\begin{tikzcd}
	& {X_1} &&& {X_1} \\
	{X_2} && {X_0} & {X_2} && {X_0} \\
	& {X_1} &&& {X_2} \\
	{X_3} && {X_1} & {X_3} && {X_1} \\
	& {X_2} &&& {X_2}
	\arrow["{d_0}", from=1-2, to=2-3]
	\arrow["\gamma", shorten <=6pt, shorten >=6pt, Rightarrow, from=1-2, to=3-2]
	\arrow["{d_0}", from=1-5, to=2-6]
	\arrow["{d_1}", from=2-1, to=1-2]
	\arrow["{d_0}"', from=2-1, to=3-2]
	\arrow["{d_1}", from=2-4, to=1-5]
	\arrow["{d_0}"', from=3-2, to=2-3]
	\arrow["{d_2}"{description}, from=3-5, to=1-5]
	\arrow["{d_0}", from=3-5, to=4-6]
	\arrow["{\widehat{\gamma}}", shorten <=6pt, shorten >=6pt, Rightarrow, from=3-5, to=5-5]
	\arrow["{d_3}", from=4-1, to=2-1]
	\arrow["{d_0}"', from=4-1, to=5-2]
	\arrow["{d_1}"', from=4-3, to=2-3]
	\arrow["{d_3}", from=4-4, to=2-4]
	\arrow["{d_1}", from=4-4, to=3-5]
	\arrow["{d_0}"', from=4-4, to=5-5]
	\arrow["{d_1}"', from=4-6, to=2-6]
	\arrow["{d_2}"{description}, from=5-2, to=3-2]
	\arrow["{d_0}"', from=5-2, to=4-3]
	\arrow["{d_0}"', from=5-5, to=4-6]
\end{tikzcd}
\end{equation}

and both of these 2-cells are required to be equal to the identity on $d_0: X_1 \to X_0$:
\begin{equation}\label{eqcllaxalgunit}
\begin{tikzcd}[column sep=scriptsize]
	&&&&& {X_1} &&&& {X_1} \\
	&&& {X_1} &&&& {X_2} \\
	{X_1} && {X_2} && {X_0} &&& {X_1} \\
	&&& {X_1} && {X_0} &&&& {X_0}
	\arrow[Rightarrow, no head, from=1-6, to=1-10]
	\arrow["{j_1}"', from=1-6, to=2-8]
	\arrow["{d_0}"', from=1-6, to=4-6]
	\arrow["{d_0}", from=1-10, to=4-10]
	\arrow["{d_0}", from=2-4, to=3-5]
	\arrow["\gamma", shorten <=6pt, shorten >=6pt, Rightarrow, from=2-4, to=4-4]
	\arrow[""{name=0, anchor=center, inner sep=0}, "{d_1}"{description}, from=2-8, to=1-10]
	\arrow["{d_0}"', from=2-8, to=3-8]
	\arrow[curve={height=-18pt}, Rightarrow, no head, from=3-1, to=2-4]
	\arrow["{j_0}", from=3-1, to=3-3]
	\arrow[""{name=1, anchor=center, inner sep=0}, curve={height=18pt}, Rightarrow, no head, from=3-1, to=4-4]
	\arrow["{d_1}", from=3-3, to=2-4]
	\arrow["{d_0}"{description}, from=3-3, to=4-4]
	\arrow[""{name=2, anchor=center, inner sep=0}, "{d_0}"{description}, from=3-8, to=4-10]
	\arrow["{d_0}"', from=4-4, to=3-5]
	\arrow["s", from=4-6, to=3-8]
	\arrow[""{name=3, anchor=center, inner sep=0}, Rightarrow, no head, from=4-6, to=4-10]
	\arrow["\gamma", shorten <=9pt, shorten >=9pt, Rightarrow, from=0, to=2]
	\arrow["{\widehat{\iota}}"', shift right=4, shorten <=4pt, shorten >=4pt, Rightarrow, from=3-3, to=1]
	\arrow["\iota", shorten >=3pt, Rightarrow, from=3-8, to=3]
\end{tikzcd}
\end{equation}

\noindent Finally, we also require the following squares to be pullbacks:
\[\begin{tikzcd}
	{X_3} && {X_2} && {X_2} && {X_1} \\
	\\
	{X_2} && {X_1} && {X_1} && {X_0}
	\arrow["{d_3}", from=1-1, to=1-3]
	\arrow["{d_0}"', from=1-1, to=3-1]
	\arrow["\lrcorner"{anchor=center, pos=0.125}, draw=none, from=1-1, to=3-3]
	\arrow["{d_0}", from=1-3, to=3-3]
	\arrow["{d_2}", from=1-5, to=1-7]
	\arrow["{d_0}"', from=1-5, to=3-5]
	\arrow["\lrcorner"{anchor=center, pos=0.125}, draw=none, from=1-5, to=3-7]
	\arrow["{d_0}", from=1-7, to=3-7]
	\arrow["{d_2}"', from=3-1, to=3-3]
	\arrow["{d_1}"', from=3-5, to=3-7]
\end{tikzcd}\]
\end{defi}

\begin{vari}
A \textit{codomain-lax category} in $\ck$ is a codomain-colax category in $\ck^{co}$. There is also a \textit{domain-colax category} which is defined as a codomain-colax category except the order of 1-cells in the defining diagram from top to bottom is reversed (this means for instance that each $d^{\bullet \bullet}_i$, $i \in \{0,1,2,3 \}$ is swapped with $d^{\bullet \bullet}_{3 - i}$ and so on). Each of these has a \textbf{pseudo} variant in which we require the invertibility of all present 2-cells.
\end{vari}

\begin{pozn}
Notice that the 1-cells $d_2: X_3 \to X_2$, $j_1: X_1 \to X_2$ are uniquely determined by their equations as the 1-cells going into pullback squares:
\begin{equation}\label{diagramd2j1}
\begin{tikzcd}
	{X_3} && {X_2} &&& {X_1} \\
	& {X_2} && {X_1} &&& {X_2} && {X_1} \\
	{X_2} &&&&& {X_0} \\
	& {X_1} && {X_0} &&& {X_1} && {X_0}
	\arrow["{d_1}"', from=4-2, to=4-4]
	\arrow["{d_0}", from=2-4, to=4-4]
	\arrow["{d_0}"', from=2-2, to=4-2]
	\arrow["{d_2}", from=2-2, to=2-4]
	\arrow["{d_2}", from=1-3, to=2-4]
	\arrow["{d_2}"{pos=0.6}, dashed, from=1-1, to=2-2]
	\arrow["{d_3}", from=1-1, to=1-3]
	\arrow["{d_0}"', from=1-1, to=3-1]
	\arrow["{d_1}"', from=3-1, to=4-2]
	\arrow["\lrcorner"{anchor=center, pos=0.125}, draw=none, from=2-2, to=4-4]
	\arrow["{d_2}", from=2-7, to=2-9]
	\arrow[curve={height=-12pt}, Rightarrow, no head, from=1-6, to=2-9]
	\arrow["{j_1}", dashed, from=1-6, to=2-7]
	\arrow["{d_0}"', from=1-6, to=3-6]
	\arrow["{d_0}"', from=2-7, to=4-7]
	\arrow["{d_0}", from=2-9, to=4-9]
	\arrow["s"', from=3-6, to=4-7]
	\arrow["{d_1}"', from=4-7, to=4-9]
	\arrow["\lrcorner"{anchor=center, pos=0.125}, draw=none, from=2-7, to=4-9]
\end{tikzcd}
\end{equation}

The 1-cells $d_1: X_3 \to X_2$, $j_0 : X_1 \to X_2$ are the 1-cells into the following pullbacks (note that this can not be taken as their definition because $d_1, j_0$ are already used for the legs of those cones):
\begin{equation}\label{diagramd1j0}
\begin{tikzcd}[column sep=scriptsize]
	{X_3} && {X_2} &&& {X_1} && {X_0} \\
	& {X_2} && {X_1} &&& {X_2} && {X_1} \\
	\\
	& {X_1} && {X_0} &&& {X_1} && {X_0}
	\arrow["{d_0}"', from=2-2, to=4-2]
	\arrow["{d_1}"', from=4-2, to=4-4]
	\arrow["{d_2}", from=2-2, to=2-4]
	\arrow["{d_0}", from=2-4, to=4-4]
	\arrow["{d_1}"{pos=0.6}, dashed, from=1-1, to=2-2]
	\arrow["\lrcorner"{anchor=center, pos=0.125}, draw=none, from=2-2, to=4-4]
	\arrow["{d_3}", from=1-1, to=1-3]
	\arrow["{d_1}", from=1-3, to=2-4]
	\arrow["{d_0}"', from=2-7, to=4-7]
	\arrow["s", from=1-8, to=2-9]
	\arrow["{d_0}", from=2-9, to=4-9]
	\arrow["{d_1}"{description}, from=4-7, to=4-9]
	\arrow["{d_2}", from=2-7, to=2-9]
	\arrow["{d_1}", from=1-6, to=1-8]
	\arrow["{j_0}", dashed, from=1-6, to=2-7]
	\arrow["\lrcorner"{anchor=center, pos=0.125}, draw=none, from=2-7, to=4-9]
	\arrow["{d_0\circ j_0}"', curve={height=12pt}, from=1-6, to=4-7]
	\arrow["{d_0 \circ d_1}"', curve={height=12pt}, from=1-1, to=4-2]
\end{tikzcd}
\end{equation}
\end{pozn}

\begin{prc}
Any category internal to $\cat(\ce)$ (a double category in $\ce$) is a codomain-colax category where all the 2-cells $\iota, \widehat{\iota}$, $\gamma, \widehat{\gamma}$ are the identities.
\end{prc}

\begin{prc}[$\ck = \cat$]\label{prclcat}
Let us describe a codomain-colax category $X$ in $\cat$ in elementary terms.
\begin{itemize}
\item $X_0$ is the category of \textit{objects} and \textit{vertical morphisms},
\item $X_1$ is the category whose objects are called the \textit{horizontal morphisms} and whose morphisms are called the \textit{squares}, with the composition operation being called the \textit{vertical composition}. As for double categories, the fact that a morphism $\alpha: g \to h$ in $X_1$ has $d_1(\alpha) = u$, $d_0(\alpha) = v$ is portrayed as follows:
\[\begin{tikzcd}
	a && b \\
	\\
	c && d
	\arrow["u"', from=1-1, to=3-1]
	\arrow[""{name=0, anchor=center, inner sep=0}, "g", from=1-1, to=1-3]
	\arrow["v", from=1-3, to=3-3]
	\arrow[""{name=1, anchor=center, inner sep=0}, "h"', from=3-1, to=3-3]
	\arrow["\alpha", shorten <=9pt, shorten >=9pt, Rightarrow, from=0, to=1]
\end{tikzcd}\]
\item $X_2$ is the category of \textit{composable pairs} of horizontal morphisms and squares. So an object of this category is a tuple $(g_1,g_2)$ of horizontal morphisms with $d_0(g_1) = d_1(g_2)$,
\item the maps $d_0, d_2: X_2 \to X_1$ are the pullback projections, i.e. on a composable pair they are defined as follows:
\begin{align*}
d_0 &: (g_1,g_2) \mapsto g_2, \\
d_2 &: (g_1,g_2) \mapsto g_1.
\end{align*}
\item $d_1:X_2 \to X_1$ is the \textit{horizontal composition operation} on horizontal morphisms and squares. We denote $d_1(g_1,g_2)$ by $g_2 \circ g_1$ if $g_i$ are horizontal morphisms. The composition preserves the domains (by axiom \eqref{eqclmonadassoc}) but in general not codomains, so a general horizontal composition of squares looks like this (where $\widehat{f'} = d_0(\beta * \alpha)$):
\[\begin{tikzcd}
	a & b & c && a && {\widehat{c}} \\
	{a'} & {b'} & {c'} & \mapsto & {a'} && {\widehat{c'}}
	\arrow[""{name=0, anchor=center, inner sep=0}, "{h_1}", from=1-1, to=1-2]
	\arrow[""{name=1, anchor=center, inner sep=0}, "{h_2}", from=1-2, to=1-3]
	\arrow["{f'}", from=1-3, to=2-3]
	\arrow["f"', from=1-1, to=2-1]
	\arrow[""{name=2, anchor=center, inner sep=0}, "{h_1'}"', from=2-1, to=2-2]
	\arrow[""{name=3, anchor=center, inner sep=0}, "{h_2'}"', from=2-2, to=2-3]
	\arrow["{\widehat{f'}}", from=1-7, to=2-7]
	\arrow["f"', from=1-5, to=2-5]
	\arrow[""{name=4, anchor=center, inner sep=0}, "{h_2 \circ h_1}", from=1-5, to=1-7]
	\arrow[""{name=5, anchor=center, inner sep=0}, "{h'_2 \circ h'_1}"', from=2-5, to=2-7]
	\arrow["v", from=1-2, to=2-2]
	\arrow["\beta", shorten <=4pt, shorten >=4pt, Rightarrow, from=1, to=3]
	\arrow["\alpha", shorten <=4pt, shorten >=4pt, Rightarrow, from=0, to=2]
	\arrow["{\beta * \alpha}"{description}, shorten <=4pt, shorten >=4pt, Rightarrow, from=4, to=5]
\end{tikzcd}\]
\item $s : X_0 \to X_1$ is the \textit{horizontal identity} functor. We denote its action on morphisms as follows:
\[\begin{tikzcd}
	a && a && {[a]} \\
	& \mapsto \\
	b && b && {[b]}
	\arrow["f"', from=1-1, to=3-1]
	\arrow[""{name=0, anchor=center, inner sep=0}, "{s(a)}", from=1-3, to=1-5]
	\arrow["f"', from=1-3, to=3-3]
	\arrow[""{name=1, anchor=center, inner sep=0}, "{s(b)}"', from=3-3, to=3-5]
	\arrow["{[f]}", from=1-5, to=3-5]
	\arrow["{s(f)}", shorten <=9pt, shorten >=9pt, Rightarrow, from=0, to=1]
\end{tikzcd}\]
By \eqref{eqclmonadunit1LOWER} the square $s(f)$ has the morphism $f$ as its domain.

\item The functor $j_1: X_1 \to X_2$ is given on objects (because of \eqref{diagramd2j1}) by the assignment $g \mapsto (g, s(b))$. The functor $j_0: X_1 \to X_2$ is given (because of \eqref{diagramd1j0}) by the assignment $j_0: g \mapsto (s(a),\widehat{g})$, where $\widehat{g} = d_0 j_0 (g)$.

\item $X_3$ is the category of horizontally composable triples of horizontal morphisms and squares. The maps $d_3, d_0: X_3 \to X_2$ are pullback projections, so they send:
\begin{align*}
d_3 &: (g_1,g_2,g_3) \mapsto (g_1, g_2), \\
d_0 &: (g_1,g_2,g_3) \mapsto (g_2, g_3).
\end{align*}
$d_2 : X_3 \to X_2$ is the unique map into the pullback \eqref{diagramd2j1} -- as such, it sends $(g_1,g_2,g_3)$ to $(g_1,g_3 \circ g_2)$. For the map $d_1: X_3 \to X_2$, the equation \eqref{diagramd1j0} tells us that:
\[
d_1: (g_1,g_2,g_3) \mapsto (g_2 \circ g_1,\widehat{g_3}),
\]
where $\widehat{g_3} = d_0 d_1 (g_1,g_2,g_3)$.

\item The 2-cells in the definition of $X$ measure to what extent does "taking horizontal composition" and "taking horizontal identities" preserve the co\-domains.  The components of $\iota$ and $\widehat{\iota}$ are given in the picture below:
\[\begin{tikzcd}
	a && {[a]} && {\widehat{b}} \\
	\\
	&& a && b
	\arrow["{s(a)}", from=1-1, to=1-3]
	\arrow[""{name=0, anchor=center, inner sep=0}, "{\widehat{g}}", from=1-3, to=1-5]
	\arrow[""{name=0p, anchor=center, inner sep=0}, phantom, from=1-3, to=1-5, start anchor=center, end anchor=center]
	\arrow["{\iota_a}"', from=1-3, to=3-3]
	\arrow[from=1-5, to=3-5]
	\arrow[""{name=1, anchor=center, inner sep=0}, "g"', from=3-3, to=3-5]
	\arrow[""{name=1p, anchor=center, inner sep=0}, phantom, from=3-3, to=3-5, start anchor=center, end anchor=center]
	\arrow["{\widehat{\iota}_{g}}", shorten <=9pt, shorten >=9pt, Rightarrow, from=0p, to=1p]
\end{tikzcd}\]
In the picture above, $\widehat{g} = d_0 j_0 (g)$. The domain of the square $\widehat{\iota}_g$ is $\iota_a$ by \eqref{eqcl2-natur1}.

\item Similarly, the components of $\gamma$ and $\widehat{\gamma}$ are given as follows:
\[\begin{tikzcd}
	a &&&& {\widehat{c}} && \bullet \\
	\\
	a && b && c && d
	\arrow["{g_2 \circ g_1}", from=1-1, to=1-5]
	\arrow[""{name=0, anchor=center, inner sep=0}, "{\widehat{g_3}}", from=1-5, to=1-7]
	\arrow["{\gamma_{g_1,g_2}}"', from=1-5, to=3-5]
	\arrow["{d_0 (\widehat{\gamma}_{g_1,g_2,g_3})}", from=1-7, to=3-7]
	\arrow["{g_1}"', from=3-1, to=3-3]
	\arrow["{g_2}"', from=3-3, to=3-5]
	\arrow[""{name=1, anchor=center, inner sep=0}, "{g_3}"', from=3-5, to=3-7]
	\arrow["{\widehat{\gamma}_{g_1,g_2,g_3}}"{description}, shorten <=9pt, shorten >=9pt, Rightarrow, from=0, to=1]
\end{tikzcd}\]
In the picture above, $\widehat{g_3} = d_0d_1(g_1,g_2,g_3)$. The domain of the square $\widehat{\gamma}_{g_1,g_2,g_3}$ is $\gamma_{g_1,g_2}$ by \eqref{eqcl2-natur2}. 
\item The axioms \eqref{eqclmonadassoc}, \eqref{eqclmonadunit2}, \eqref{eqclmonadunit1} assert that the composition of horizontal morphisms is in a certain sense associative and unital -- this will be described in detail in Lemma \ref{THML_codcolax_are_assoc_unital_ABSTRAKTNI}.
\item The axiom \eqref{eqcllaxalgassoc} is the assertion that for any composable triple of horizontal morphisms $(g_1,g_2,g_3)$, we have:
\begin{equation}\label{eqcllaxalgassocALG}
d_0 (\widehat{\gamma}_{g_1,g_2,g_3}) \circ \gamma_{g_2 g_1, g_3} = \gamma_{g_2,g_3} \circ \gamma_{g_1, g_3 g_2}
\end{equation}
\item The axiom \eqref{eqcllaxalgunit} says that for any horizontal $g: a \to b$, we have (denoting $\widehat{g} := d_0 j_0 (g)$):
\begin{align}
d_0 \widehat{\iota_g} \circ \gamma_{s(a), \widehat{g}} &= 1_b \label{eqcllaxalgunitALG1} \\
\iota_b \circ \gamma_{g,s(b)}  &= 1_b. \label{eqcllaxalgunitALG2}
\end{align}

\end{itemize}

\medskip

\noindent It will be useful to unwrap what it means for $\gamma$ to be a natural transformation. Consider the following pair of composable squares and their composite:
\[\begin{tikzcd}
	a && b && c && a &&& {\widehat{c}} \\
	\\
	x && y && z && x &&& {\widehat{z}}
	\arrow[""{name=0, anchor=center, inner sep=0}, "{g_1}", from=1-1, to=1-3]
	\arrow[""{name=0p, anchor=center, inner sep=0}, phantom, from=1-1, to=1-3, start anchor=center, end anchor=center]
	\arrow["u"', from=1-1, to=3-1]
	\arrow[""{name=1, anchor=center, inner sep=0}, "{g_2}", from=1-3, to=1-5]
	\arrow[""{name=1p, anchor=center, inner sep=0}, phantom, from=1-3, to=1-5, start anchor=center, end anchor=center]
	\arrow["v"{description}, from=1-3, to=3-3]
	\arrow["w", from=1-5, to=3-5]
	\arrow[""{name=2, anchor=center, inner sep=0}, "{g_2 \circ g_1}", from=1-7, to=1-10]
	\arrow[""{name=2p, anchor=center, inner sep=0}, phantom, from=1-7, to=1-10, start anchor=center, end anchor=center]
	\arrow["u"', from=1-7, to=3-7]
	\arrow["{\widehat{w}}", from=1-10, to=3-10]
	\arrow[""{name=3, anchor=center, inner sep=0}, "{h_1}"', from=3-1, to=3-3]
	\arrow[""{name=3p, anchor=center, inner sep=0}, phantom, from=3-1, to=3-3, start anchor=center, end anchor=center]
	\arrow[""{name=4, anchor=center, inner sep=0}, "{h_2}"', from=3-3, to=3-5]
	\arrow[""{name=4p, anchor=center, inner sep=0}, phantom, from=3-3, to=3-5, start anchor=center, end anchor=center]
	\arrow[""{name=5, anchor=center, inner sep=0}, "{h_2 \circ h_1}"', from=3-7, to=3-10]
	\arrow[""{name=5p, anchor=center, inner sep=0}, phantom, from=3-7, to=3-10, start anchor=center, end anchor=center]
	\arrow["{\kappa_1}", shorten <=9pt, shorten >=9pt, Rightarrow, from=0p, to=3p]
	\arrow["{\kappa_2}", shorten <=9pt, shorten >=9pt, Rightarrow, from=1p, to=4p]
	\arrow["{\kappa_2 *\kappa_1}", shorten <=9pt, shorten >=9pt, Rightarrow, from=2p, to=5p]
\end{tikzcd}\]
This is a morphism $(\kappa_1, \kappa_2): (g_1, g_2) \to (g'_1, g'_2)$ in $X_2$. The naturality square of $\gamma$ for this morphism says that the following square commutes in $X_0$:
\begin{equation}\label{gammanaturalita}
\begin{tikzcd}
	{\widehat{c}} && c \\
	&&& {\text{ in } X_0} \\
	{\widehat{z}} && z
	\arrow["{\gamma_{g_1,g_2}}", from=1-1, to=1-3]
	\arrow["{\widehat{w}}"', from=1-1, to=3-1]
	\arrow["w", from=1-3, to=3-3]
	\arrow["{\gamma_{h_1,h_2}}"', from=3-1, to=3-3]
\end{tikzcd}
\end{equation}
\end{prc}

A simple example of a codomain-colax category that is \textbf{not} a double category is given as follows:

\begin{pr}
Any comonad $(A,q,\delta,\epsilon)$ on a small category $A$ gives a codomain-colax double category as follows.
\begin{itemize}
\item the objects are the objects of $A$,
\item the vertical morphisms are the morphisms of $A$,
\item the horizontal morphisms are in bijection with $\ob A$ -- we denote by $\widetilde{a}$ the horizontal morphism corresponding to $a \in \ob A$. Its domain and codomain are the following:
\[
\widetilde{a}: a \to qa.
\]
\item the squares are in bijection with $\mor A$ and we again denote them using tilde. The horizontal composition of horizontal morphisms and squares is then given as follows:
\[\begin{tikzcd}
	a && qa && {q^2a} && a && qa \\
	&&&&& \mapsto \\
	b && qb && {q^2b} && b && qb
	\arrow[""{name=0, anchor=center, inner sep=0}, "{\widetilde{a}}", from=1-1, to=1-3]
	\arrow["u"', from=1-1, to=3-1]
	\arrow[""{name=1, anchor=center, inner sep=0}, "{\widetilde{qa}}", from=1-3, to=1-5]
	\arrow["qu", from=1-3, to=3-3]
	\arrow["{q^2u}", from=1-5, to=3-5]
	\arrow[""{name=2, anchor=center, inner sep=0}, "{\widetilde{a}}", from=1-7, to=1-9]
	\arrow["u"', from=1-7, to=3-7]
	\arrow["qu", from=1-9, to=3-9]
	\arrow[""{name=3, anchor=center, inner sep=0}, "{\widetilde{b}}"', from=3-1, to=3-3]
	\arrow[""{name=4, anchor=center, inner sep=0}, "{\widetilde{qb}}"', from=3-3, to=3-5]
	\arrow[""{name=5, anchor=center, inner sep=0}, "{\widetilde{b}}"', from=3-7, to=3-9]
	\arrow["{\widetilde{u}}", shorten <=9pt, shorten >=9pt, Rightarrow, from=0, to=3]
	\arrow["{\widetilde{qu}}", shorten <=9pt, shorten >=9pt, Rightarrow, from=1, to=4]
	\arrow["{\widetilde{u}}", shorten <=9pt, shorten >=9pt, Rightarrow, from=2, to=5]
\end{tikzcd}\]
\item every horizontal morphism is the ``horizontal identity” on its domain. The counitors use the comonad counit:
\[\begin{tikzcd}
	a && qa && {q^2a} \\
	\\
	&& a && qa
	\arrow["{\widetilde{a}}", from=1-1, to=1-3]
	\arrow[""{name=0, anchor=center, inner sep=0}, "{\widetilde{qa}}", from=1-3, to=1-5]
	\arrow["{\epsilon_a}"', from=1-3, to=3-3]
	\arrow["{q \epsilon_a}", from=1-5, to=3-5]
	\arrow[""{name=1, anchor=center, inner sep=0}, "{\widetilde{a}}"', from=3-3, to=3-5]
	\arrow["{\widetilde{\epsilon_a}}", shorten <=9pt, shorten >=9pt, Rightarrow, from=0, to=1]
\end{tikzcd}\]
\item the coassociators use the comonad comultiplication:
\[\begin{tikzcd}
	a &&&& qa && {q^2a} \\
	\\
	a && qa && {q^2a} && {q^3a}
	\arrow["{\widetilde{a}}", from=1-1, to=1-5]
	\arrow[""{name=0, anchor=center, inner sep=0}, "{\widetilde{qa}}", from=1-5, to=1-7]
	\arrow["{\delta_a}"', from=1-5, to=3-5]
	\arrow["{q \delta_a}", from=1-7, to=3-7]
	\arrow["{\widetilde{a}}"', from=3-1, to=3-3]
	\arrow["{\widetilde{qa}}"', from=3-3, to=3-5]
	\arrow[""{name=1, anchor=center, inner sep=0}, "{\widetilde{q^2a}}"', from=3-5, to=3-7]
	\arrow["{\widetilde{\delta_a}}", shorten <=9pt, shorten >=9pt, Rightarrow, from=0, to=1]
\end{tikzcd}\]
\end{itemize}

\end{pr}

\noindent The above example is a special case of the class of examples that motivates the introduction of codomain-colax categories:

\begin{prc}\label{PR_cartesian2monad_res_codcolaxcat}
Let $T$ be a cartesian 2-monad on a 2-category $\ck$. Any colax $T$-algebra $(A, a, \gamma, \iota)$ gives rise to a codomain-colax category:
\begin{center}
\begin{tikzpicture}
    \node (A) at (-1,0) {$T^3\ca$};
    \node (B) at (2.5,0) {$T^2\ca$};
    \node (C) at (6,0) {$T\ca$};
    \node (D) at (-4.5,0) {$T^4\ca$};
    \draw[transform canvas={yshift=3ex},->] (B) -- node[fill=white]{$Ti_{T\ca}$} (A);
    \draw[transform canvas={yshift=6ex},->] (A) -- node[fill=white]{$ m_{T\ca}$} (B);
    \draw[transform canvas={yshift=0ex},->] (A) --node[fill=white]{ $ Tm_{\ca}$}  (B);
    \draw[transform canvas={yshift=-3ex},->] (B) -- node[fill=white]{ $T^2i_{\ca}$} (A);
    \draw[transform canvas={yshift=-6ex},->] (A) -- node[fill=white]{ $ T^2a$} (B);

    \draw[transform canvas={yshift=3ex},->] (B) -- node[fill=white]{ $m_\ca$} (C);
    \draw[transform canvas={yshift=0ex},->] (C) --node[fill=white]{ $Ti_{\ca}$}  (B);
    \draw[transform canvas={yshift=-3ex},->] (B) -- node[fill=white]{ $Ta$} (C);

    \draw[transform canvas={yshift=4.5ex},->] (D) --node[fill=white]{ $ m_{T^2\ca}$}  (A);
    \draw[transform canvas={yshift=1.5ex},->] (D) --node[fill=white]{ $ Tm_{T\ca}$}  (A);
    \draw[transform canvas={yshift=-1.5ex},->] (D) --node[fill=white]{ $ T^2m_{\ca}$}  (A);
    \draw[transform canvas={yshift=-4.5ex},->] (D) --node[fill=white]{ $ T^3a$}  (A);
    
\node at (-1,1.5) {};
\node at (-1,1.5) {};

\node at (-1,-1.5) {};
\node at (-4.5,0) {};
\end{tikzpicture}
\end{center}
\noindent with 2-cells:
\begin{align*}
T\gamma &: TaTm_\ca \Rightarrow TaT^2a, \\
T\iota &: TaTi_\ca \Rightarrow 1_{T\ca},\\
T^2\gamma &: T^2aT^2m_\ca \Rightarrow T^2aT^3a,\\
T^2\iota &: T^2aT^2i_\ca \Rightarrow 1_{T^2A}.
\end{align*}
The axioms \eqref{eqcllaxalgassoc} and \eqref{eqcllaxalgunit} are the colax algebra axioms. The other axioms are as follows:
\begin{itemize}
\item \eqref{eqclmonadassoc} - \eqref{eqclmonadassocT}: 2-monad associativity laws,
\item \eqref{eqclmonadunit1LOWER} - \eqref{eqclmonadunit3}: 2-monad unit laws,
\item \eqref{eqclunit_naturality}, \eqref{eqcl2-natur1}: the 2-naturality of $i: 1 \Rightarrow T$,
\item \eqref{eqclmultipl_naturality}, \eqref{eqcl2-natur2}: the 2-naturality of $m: T^2 \Rightarrow T$.
\end{itemize}

\end{prc}

\begin{nota}
In a codomain-colax double category $X$ the convention we have is that the picture below represents a composable pair $(g_1,g_2)$ of horizontal morphisms, but \textbf{not} their composite (as the codomain of the composite will not necessarily be $c$):
\[\begin{tikzcd}
	a & b & c
	\arrow["{g_1}"', from=1-1, to=1-2]
	\arrow["{g_2}"', from=1-2, to=1-3]
\end{tikzcd}\]

\noindent Next, consider the diagrams:
\[\begin{tikzcd}
	a & a & c & d \\
	b & b & c & d
	\arrow[Rightarrow, no head, from=1-1, to=1-2]
	\arrow["{f_1}"', from=1-1, to=2-1]
	\arrow["{f_2}", from=1-2, to=2-2]
	\arrow[""{name=0, anchor=center, inner sep=0}, "{g_1}", from=1-3, to=1-4]
	\arrow[Rightarrow, no head, from=1-3, to=2-3]
	\arrow[Rightarrow, no head, from=1-4, to=2-4]
	\arrow[Rightarrow, no head, from=2-1, to=2-2]
	\arrow[""{name=1, anchor=center, inner sep=0}, "{g_2}"', from=2-3, to=2-4]
	\arrow[shorten <=4pt, shorten >=4pt, Rightarrow, from=0, to=1]
\end{tikzcd}\]
by the top right picture we mean a square in $X$, or in other words a morphism $g_1 \to g_2$ in $X_1$. By the top-left picture we mean \textbf{the assertion} that morphisms $f_1,f_2$ in $X_1$ are equal. Note that in a general $X$, there are no ``horizontal identities” which would have the same domain and codomain, so the top left square would make no sense otherwise.
\end{nota}

\begin{lemma}\label{THML_codcolax_are_assoc_unital_ABSTRAKTNI}
The horizontal composition in a codomain-colax double category $X$ in $\cat$ is \textbf{associative} in the sense that given a composable triple of morphisms:
\[\begin{tikzcd}
	a & b & c & d
	\arrow["{g_1}"', from=1-1, to=1-2]
	\arrow["{g_2}"', from=1-2, to=1-3]
	\arrow["{g_3}"', from=1-3, to=1-4]
\end{tikzcd}\]
We have:
\[
(g_3 \circ g_2) \circ g_1 = \widehat{g_3} \circ (g_2 \circ g_1),
\]
where we denote $\widehat{g_3} = d_0 d_1 (g_1, g_2, g_3)$.

\vspace{0.5cm}

\noindent The composition is also \textbf{unital} in the sense that given $g: a \to b$, both pairs:
\[\begin{tikzcd}
	a && {[a]} && {d_0 (\widehat{g})} \\
	a && b && {[b]}
	\arrow["{s(a)}", from=1-1, to=1-3]
	\arrow["{\widehat{g}}", from=1-3, to=1-5]
	\arrow["g", from=2-1, to=2-3]
	\arrow["{s(b)}", from=2-3, to=2-5]
\end{tikzcd}\]

(where we denote $\widehat{g}= d_0j_0(g)$) compose to $g$.
\end{lemma}
\begin{proof}
The associativity follows since:
\begin{align*}
\widehat{g_3} \circ (g_2 \circ g_1) &= d_1 (g_2 \circ g_1, \widehat{g_3})\\
&= d_1 d_1 (g_1, g_2, g_3)\\
&\sstackrel{\eqref{eqclmonadassoc}}{=} d_1 d_2 (g_1, g_2, g_3)\\
&= d_1 (g_1, g_3 \circ g_2)\\
&= (g_3 \circ g_2) \circ g_1.
\end{align*}

The unitality follows because:
\begin{align*}
s(b) \circ g &= d_1 (g, s(b)) = d_1 j_1(g) \sstackrel{\eqref{eqclmonadunit2}}{=} g, \\
\widehat{g} \circ s(a) &= d_1 (\widehat{g},s (a)) = d_1 j_0(g) \sstackrel{\eqref{eqclmonadunit1}}{=} g.
\end{align*}
\end{proof}

\begin{termi}
Clearly, every codomain-colax category $X$ in $\ck$ has an underlying colax coherence data in $\ck$. By the \textit{codescent object} of $X$ or the \textit{(coherence) cocone} for $X$ we mean such for its underlying colax coherence data.
\end{termi}

\begin{pozn}
Given a codomain-colax double category $X$, the cocone for its underlying coherence data consists of a pair $(F: X_0 \to \cy, \xi: Fd_1 \Rightarrow Fd_0)$ of a functor and a natural transformation satisfying equations from \ref{NOTA_coherence_data}: for all $(g_1,g_2) \in X_2$, $a \in \ob X_0$:
\begin{align}
F\gamma_{(g_1,g_2)} \circ \xi_{g_2 \circ g_1} &= \xi_{g_2}\circ \xi_{g_1},\label{colaxcohcocone1} \\
F\iota_a \circ \xi_{s(a)} &= 1_{Fa} \label{colaxcohcocone2} 
\end{align}
\end{pozn}

\subsection{Codescent for internal categories}\label{SEC_thecaseckCatE}

\noindent The following is a generalization of codomain-discrete double categories from Definition \ref{deficoddiscr}:

\begin{defi}
Let $\ce$ be a category with pullbacks and let $X$ be a codomain-colax category in $\cat(\ce)$. We say that $X$ is \textit{domain-codiscrete} if the map $d_1 : X_1 \to X_0$ is a discrete fibration.
\end{defi}

The \textbf{goal} of this section is to show that 2-categories of form $\cat(\ce)$ always admit codescent objects of this kind of codomain-colax categories (Theorem \ref{THM_catE_admits_codobj_of_domdiscr_codlax}). Also that 2-functors of form $\cat(F)$ automatically preserve these diagrams and their colimits (Theorem \ref{THM_CatT_preserves_special_codobjs}).

The motivating example of a domain-codiscrete, codomain-colax category in $\cat(\ce)$ is given by the following:

\begin{prop}\label{prdomdiscrcodlax}
Let $\ce$ be a category with pullbacks and let $T$ be a 2-monad on $\cat(\ce)$ of form $\cat(T)$. For any colax $T$-algebra $(A,a,\gamma,\iota)$, its resolution $\res{A,a,\gamma,\iota}$ is a domain-codiscrete, codomain-colax category in $\talgs$.
\end{prop}
\begin{proof}
We know it is a codomain-colax category since $T$ is cartesian (Example \ref{PR_cartesian2monad_res_codcolaxcat}). Next, the domain functor $m_A : T^2A \to TA$ is given by the 2-monad multiplication, which is a discrete opfibration by Proposition \ref{THMC_catT2monad_multiplication_is_discrete_opfib}.

\end{proof}

\begin{nonpr}
The resolution of the terminal $S$-algebra for the free symmetric strict monoidal category 2-monad is \textbf{not} domain-codiscrete since the multiplication is not a discrete opfibration, see Non-example \ref{NOPR_multipl_not_discr_opfib_freesymstrictmoncat}.
\end{nonpr}

\begin{con}\label{CON_cnrX_for_X_domdiscr_codcolaxcat}
Assume now $\ck = \cat(\ce)$ for a category $\ce$ with pullbacks and let $X$ be a \textbf{domain-codiscrete}, codomain-colax category in $\ck$. We will proceed to construct the \textit{internal category of corners} $\cnr{X}$ associated to $X$. Define the object of objects as $\cnr{X}_0 := X_{00}$. The object of morphisms will be the pullback below left:
\[\begin{tikzcd}
	{\cnr{X}_1} && {X_{01}} && {X_{11}} && {X_{10}} \\
	\\
	{X_{10}} && {X_{00}} && {X_{01}} && {X_{00}}
	\arrow["{\pi_2}", from=1-1, to=1-3]
	\arrow["{\pi_1}"', from=1-1, to=3-1]
	\arrow["\lrcorner"{anchor=center, pos=0.125}, draw=none, from=1-1, to=3-3]
	\arrow["s", from=1-3, to=3-3]
	\arrow["t", from=1-5, to=1-7]
	\arrow["{d_1}"', from=1-5, to=3-5]
	\arrow["\lrcorner"{anchor=center, pos=0.125}, draw=none, from=1-5, to=3-7]
	\arrow["{d_1}", from=1-7, to=3-7]
	\arrow["{d_0}"', from=3-1, to=3-3]
	\arrow["t"', from=3-5, to=3-7]
\end{tikzcd}\]
The square above right is a pullback because $X$ is domain-codiscrete. In what follows we will employ the philosophy of generalized elements: for instance a generalized element $Z \to X_{11}$ of $X_{11}$ can be presented as a ``bottom-left corner”, i.e. a pair:
\[
(u: Z \to X_{01}, h: Z \to X_{10}),
\]
of generalized elements of $X_{01}, X_{10}$ with the property that $t \circ u = d_1 \circ h$. Define the identities map $\ids : X_{00} \to \cnr{X}_1$ to be the generalized element pictured below (here $\iota$ is the counitor for the codomain-lax category $X$):
\[
(s_0:X_{00} \to X_{10}, \iota: X_{00} \to X_{01}).
\]
The underlying reflexive graph of $\cnr{X}$ is given by the following:
\[\begin{tikzcd}
	{\cnr{X}_1} && {X_{00}}
	\arrow["{d_1 \circ \pi_1}", shift left=3, from=1-1, to=1-3]
	\arrow["{t \circ \pi_2}"', shift right=3, from=1-1, to=1-3]
	\arrow["\ids"{description}, from=1-3, to=1-1]
\end{tikzcd}\]

\noindent To define the composition arrow:
\[
\comp: \cnr{X}_2 \to \cnr{X}_1,
\]
it suffices to specify its action on a generalized element $((g,u),(h,v)): Z \to \cnr{X}_2$. The composite $\comp \circ ((g,u),(h,v)): Z \to \cnr{X}_1$ is then of form $(c_1,c_2)$, where the first component is the element (here $(u,h): Z \to X_{11}$):
\[\begin{tikzcd}
	Z && {X_{20}} && {X_{10}}
	\arrow["{(g,s \circ (u,h))}", from=1-1, to=1-3]
	\arrow["{d_1}", from=1-3, to=1-5]
\end{tikzcd}\]
For the second component, note that the following triple gives a generalized element of the ``object of composable triples” of $X_0$, i.e. a 1-cell $Z \to X_{03}$  (here by $\gamma$ we denote the coassociator for $X$):
\[\begin{tikzcd}[row sep=small]
	Z && {X_{20}} && {X_{01}} \\
	Z && {X_{11}} && {X_{01}} \\
	Z &&&& {X_{01}}
	\arrow["{(g,s \circ (u,h))}", from=1-1, to=1-3]
	\arrow["\gamma", from=1-3, to=1-5]
	\arrow["{(u,h)}", from=2-1, to=2-3]
	\arrow["{d_0}", from=2-3, to=2-5]
	\arrow["v", from=3-1, to=3-5]
\end{tikzcd}\]
Denoting the ternary composition arrow for $X_0$ by $\comp_3: X_{03} \to X_{01}$, define $c_2$ to be the arrow:
\[\begin{tikzcd}[row sep=small]
	Z &&& {X_{03}} && {X_{01}}
	\arrow["{((g,s \circ (u,h), d_0 \circ (u,h), v)}", from=1-1, to=1-4]
	\arrow["{\comp_3}", from=1-4, to=1-6]
\end{tikzcd}\]
\end{con}

\begin{pr}[$\ce = \set$]
The objects of $\cnr{X}$ are the objects of $X_0$. A morphism $a \to b$ is a ``top-right corner” in $X$, i.e. a pair of a horizontal morphism followed by a vertical one:
\[\begin{tikzcd}
	a & x \\
	& b
	\arrow["g", from=1-1, to=1-2]
	\arrow["u", from=1-2, to=2-2]
\end{tikzcd}\]
Because the domain functor $d_1 : X_1 \to X_0$ is a discrete fibration, any bottom-left corner can be uniquely filled into a square:
\[\begin{tikzcd}
	y && \bullet \\
	\\
	b && x
	\arrow[""{name=0, anchor=center, inner sep=0}, dashed, from=1-1, to=1-3]
	\arrow["u"', from=1-1, to=3-1]
	\arrow[dashed, from=1-3, to=3-3]
	\arrow[""{name=1, anchor=center, inner sep=0}, "g"', from=3-1, to=3-3]
	\arrow["{\exists !}", shorten <=9pt, shorten >=9pt, Rightarrow, from=0, to=1]
\end{tikzcd}\]

\noindent The composite of two composable corners $(g,u)$, $(h,v)$ is given by the corner:
\[
(\widehat{h} \circ g, v \circ \widehat{u} \circ \gamma_{g,\widehat{h}}),
\]
as pictured below:
\[\begin{tikzcd}
	a &&&& {d_0(\widehat{h} \circ g)} \\
	a && {a'} && {\widehat{b'}} \\
	&& b && {b'} \\
	&&&& c
	\arrow["{\widehat{h} \circ g}", from=1-1, to=1-5]
	\arrow["{\gamma_{g,\widehat{h}}}", from=1-5, to=2-5]
	\arrow["g", from=2-1, to=2-3]
	\arrow[""{name=0, anchor=center, inner sep=0}, "{\widehat{h}}", from=2-3, to=2-5]
	\arrow[""{name=0p, anchor=center, inner sep=0}, phantom, from=2-3, to=2-5, start anchor=center, end anchor=center]
	\arrow["u"', from=2-3, to=3-3]
	\arrow["{\widehat{u}}", from=2-5, to=3-5]
	\arrow[""{name=1, anchor=center, inner sep=0}, "h"', from=3-3, to=3-5]
	\arrow[""{name=1p, anchor=center, inner sep=0}, phantom, from=3-3, to=3-5, start anchor=center, end anchor=center]
	\arrow["v", from=3-5, to=4-5]
	\arrow["{\exists !}", shorten <=4pt, shorten >=4pt, Rightarrow, from=0p, to=1p]
\end{tikzcd}\]

\noindent The identity morphism on $a$ is the corner:
\[\begin{tikzcd}
	a && {[a]} \\
	&& a
	\arrow["{s(a)}", from=1-1, to=1-3]
	\arrow["{\iota_a}", from=1-3, to=2-3]
\end{tikzcd}\]
\end{pr}

\begin{pozn}[Meditation on preservation]\label{POZN_catF_preserves_Cnr}
A codomain-colax category in a 2-catego\-ry $\ck$, being a diagram in $\ck$ satisfying some equations between 1-cells and 2-cells and also having some pullback properties, is preserved by any pullback-preserving 2-functor $H: \ck \to \cl$. Let us denote by $H_* X$ the resulting codomain-colax category in $\cl$.

If $X$ is a codomain-colax category in the 2-category $\cat(\ce)$ that is moreover domain-codiscrete, and $H: \cat(\ce) \to \cat(\ce')$ is a 2-functor of form $\cat(H')$, the resulting codomain-colax category $H_*X$ is again domain-codiscrete since 2-functors of this class preserve discrete (op)fibrations by Proposition \ref{THMP_catF_preserves_discropfibs}. These 2-functors moreover preserve the category of corners:

\begin{lemma}\label{THML_CnrFX_is_FCnrX}
Let $X$ be a codomain-colax, domain-codiscrete category in $\cat(\ce)$ and $H: \cat(\ce) \to \cat(\ce')$ a 2-functor of form $\cat(H')$. The following constructions are isomorphic if regarded as diagrams $\Delta^{op}_2 \to \ce'$:
\[
\cnr{H_* X} \cong H' \circ \cnr{X}.
\]
\end{lemma}
\begin{proof}
This follows from the fact that that both $X$ and $\cnr{X}$ are built from the 1-cells of $\ce$ on which $H'$ is applied.
\end{proof}
\end{pozn}

\begin{lemma}\label{THML_codlaxX_cnrX_is_internal_cat}
Given a codomain-colax category $X$ in $\cat(\ce)$, the structure $\cnr{X}$ thus described is an internal category in $\cat(\ce)$.
\end{lemma}
\begin{proof}
Assume first that $\ce = \set$. We verify the left and right unit laws and the associativity.
\begin{itemize}
\item \textbf{The left unit law}: The composition $(s(b), \iota_b) \circ (g,u)$ equals $(g,u)$ due to Lemma \ref{THML_codcolax_are_assoc_unital_ABSTRAKTNI}, the axiom \eqref{eqcllaxalgunitALG2} (square (b)) and the naturality of $\iota$ (square (a)), as displayed in this picture:
\[\begin{tikzcd}
	a &&&& x && x && x \\
	a && x && {[x]} && {[x]} & {(b)} \\
	&& b && {[b]} & {(a)} & x && x \\
	&&&& b && b && b
	\arrow["g", from=2-1, to=2-3]
	\arrow["u"', from=2-3, to=3-3]
	\arrow[""{name=0, anchor=center, inner sep=0}, "{s(b)}"', from=3-3, to=3-5]
	\arrow["{\iota_b}", from=3-5, to=4-5]
	\arrow[""{name=1, anchor=center, inner sep=0}, "{s(x)}", from=2-3, to=2-5]
	\arrow["{[u]}", from=2-5, to=3-5]
	\arrow["g", from=1-1, to=1-5]
	\arrow["{\gamma_{g,s(x)}}", from=1-5, to=2-5]
	\arrow[Rightarrow, no head, from=2-5, to=2-7]
	\arrow[Rightarrow, no head, from=4-5, to=4-7]
	\arrow[Rightarrow, no head, from=1-5, to=1-7]
	\arrow["{\gamma_{g,s(x)}}", from=1-7, to=2-7]
	\arrow[Rightarrow, no head, from=1-7, to=1-9]
	\arrow[Rightarrow, no head, from=4-7, to=4-9]
	\arrow["{\iota_x}", from=2-7, to=3-7]
	\arrow["u", from=3-7, to=4-7]
	\arrow["u", from=3-9, to=4-9]
	\arrow[Rightarrow, no head, from=3-7, to=3-9]
	\arrow[Rightarrow, no head, from=1-9, to=3-9]
	\arrow["{s(u)}", shorten <=4pt, shorten >=4pt, Rightarrow, from=1, to=0]
\end{tikzcd}\]

\item \textbf{The right unit law}: The composition $(g,u) \circ (s(a), \iota_a)$ is $(g,u)$ due to Lemma \ref{THML_codcolax_are_assoc_unital_ABSTRAKTNI} and the axiom  \eqref{eqcllaxalgunitALG1} (square $(*)$) as displayed in this picture:
\[\begin{tikzcd}
	a &&&& x && x \\
	a && {[a]} && {\widehat{x}} & {(*)} \\
	&& a && x && x \\
	&&&& b && b
	\arrow["g", from=1-1, to=1-5]
	\arrow[equals, from=1-5, to=1-7]
	\arrow["{\gamma_{s(a),d_0j_0(g)}}", from=1-5, to=2-5]
	\arrow[equals, from=1-7, to=3-7]
	\arrow["{s(a)}", from=2-1, to=2-3]
	\arrow[""{name=0, anchor=center, inner sep=0}, "{\widehat{g}}", from=2-3, to=2-5]
	\arrow[""{name=0p, anchor=center, inner sep=0}, phantom, from=2-3, to=2-5, start anchor=center, end anchor=center]
	\arrow["{\iota_a}"', from=2-3, to=3-3]
	\arrow["{d_0(\widehat{\iota_g})}", from=2-5, to=3-5]
	\arrow[""{name=1, anchor=center, inner sep=0}, "g"', from=3-3, to=3-5]
	\arrow[""{name=1p, anchor=center, inner sep=0}, phantom, from=3-3, to=3-5, start anchor=center, end anchor=center]
	\arrow[equals, from=3-5, to=3-7]
	\arrow["u", from=3-5, to=4-5]
	\arrow["u", from=3-7, to=4-7]
	\arrow[equals, from=4-5, to=4-7]
	\arrow["{\widehat{\iota_g}}", shorten <=4pt, shorten >=4pt, Rightarrow, from=0p, to=1p]
\end{tikzcd}\]

\item \textbf{Associativity}: Consider a triple of corners like this and formally fill them into a single corner (note that we are not composing anything yet):
\[\begin{tikzcd}
	a && x && {b'} && {y'} \\
	&& b && y && {c'} \\
	&&&& c && z \\
	&&&&&& d
	\arrow["{g_1}", from=1-1, to=1-3]
	\arrow[""{name=0, anchor=center, inner sep=0}, "{g'_2}", from=1-3, to=1-5]
	\arrow["{u_1}"', from=1-3, to=2-3]
	\arrow[""{name=1, anchor=center, inner sep=0}, "{g_3''}", from=1-5, to=1-7]
	\arrow["{u_1'}", from=1-5, to=2-5]
	\arrow["{u_1''}", from=1-7, to=2-7]
	\arrow[""{name=2, anchor=center, inner sep=0}, "{g_2}"', from=2-3, to=2-5]
	\arrow[""{name=3, anchor=center, inner sep=0}, "{g_3'}"{description}, from=2-5, to=2-7]
	\arrow["{u_2}", from=2-5, to=3-5]
	\arrow["{u_2'}", from=2-7, to=3-7]
	\arrow[""{name=4, anchor=center, inner sep=0}, "{g_3}"', from=3-5, to=3-7]
	\arrow["{u_3}", from=3-7, to=4-7]
	\arrow["{\kappa_1}", shorten <=4pt, shorten >=4pt, Rightarrow, from=0, to=2]
	\arrow["{\kappa_3}", shorten <=4pt, shorten >=4pt, Rightarrow, from=1, to=3]
	\arrow["{\kappa_2}", shorten <=4pt, shorten >=4pt, Rightarrow, from=3, to=4]
\end{tikzcd}\]

The composite $((g_3, u_3) \circ (g_2, u_2)) \circ (g_1, u_1)$ equals:
\[\begin{tikzcd}
	a &&&&& {\widehat{\widehat{y''}}} \\
	a & x &&&& {\widehat{y'}} \\
	& b &&&& {\widehat{c'}} \\
	& b && y && {c'} \\
	&&& c && z \\
	&&&&& d
	\arrow["{(g_3'' \circ g_2') \circ g_1}", from=1-1, to=1-6]
	\arrow["{\gamma_{g_1, g_3'' \circ g_2'}}", from=1-6, to=2-6]
	\arrow["{g_1}", from=2-1, to=2-2]
	\arrow[""{name=0, anchor=center, inner sep=0}, "{g_3'' \circ g_2'}", from=2-2, to=2-6]
	\arrow[""{name=0p, anchor=center, inner sep=0}, phantom, from=2-2, to=2-6, start anchor=center, end anchor=center]
	\arrow["{u_1}"', from=2-2, to=3-2]
	\arrow["{\widehat{u_1''}}", from=2-6, to=3-6]
	\arrow[""{name=1, anchor=center, inner sep=0}, "{g_3' \circ g_2}"', from=3-2, to=3-6]
	\arrow[""{name=1p, anchor=center, inner sep=0}, phantom, from=3-2, to=3-6, start anchor=center, end anchor=center]
	\arrow["{\gamma_{g_2, g'_3}}", from=3-6, to=4-6]
	\arrow["{g_2}", from=4-2, to=4-4]
	\arrow[""{name=2, anchor=center, inner sep=0}, "{g'_3}", from=4-4, to=4-6]
	\arrow[""{name=2p, anchor=center, inner sep=0}, phantom, from=4-4, to=4-6, start anchor=center, end anchor=center]
	\arrow["{u_2}"', from=4-4, to=5-4]
	\arrow["{u_2'}", from=4-6, to=5-6]
	\arrow[""{name=3, anchor=center, inner sep=0}, "{g_3}"', from=5-4, to=5-6]
	\arrow[""{name=3p, anchor=center, inner sep=0}, phantom, from=5-4, to=5-6, start anchor=center, end anchor=center]
	\arrow["{u_3}", from=5-6, to=6-6]
	\arrow["{\kappa_3 \circ_H \kappa_1}", shorten <=4pt, shorten >=4pt, Rightarrow, from=0p, to=1p]
	\arrow["{\kappa_2}", shorten <=4pt, shorten >=4pt, Rightarrow, from=2p, to=3p]
\end{tikzcd}\]

The composite $(g_3, u_3) \circ ((g_2,u_2) \circ (g_1,u_1))$ equals:
\[\begin{tikzcd}
	a &&&&&& {\widehat{\widehat{y'}}} \\
	a &&&& {\widehat{b'}} && {\widehat{y'}} \\
	a && x && {b'} && {y'} \\
	&& b && y && {c'} \\
	&&&& c && z \\
	&&&&&& d
	\arrow["{\widehat{g''_3} \circ (g_2' \circ g_1)}", from=1-1, to=1-7]
	\arrow["{\gamma_{g_2' \circ g_1, g_3'''}}", from=1-7, to=2-7]
	\arrow["{g_2' \circ g_1}", from=2-1, to=2-5]
	\arrow[""{name=0, anchor=center, inner sep=0}, "{\widehat{g''_3}}", from=2-5, to=2-7]
	\arrow["{\gamma_{g_1, g'_2}}"', from=2-5, to=3-5]
	\arrow["{d_0(\widehat{\gamma}_{g_1, g'_2, g''_3})}", from=2-7, to=3-7]
	\arrow["{g_1}", from=3-1, to=3-3]
	\arrow[""{name=1, anchor=center, inner sep=0}, "{g'_2}", from=3-3, to=3-5]
	\arrow["{u_1}"', from=3-3, to=4-3]
	\arrow[""{name=2, anchor=center, inner sep=0}, "{g_3''}"{description}, from=3-5, to=3-7]
	\arrow["{u_1'}", from=3-5, to=4-5]
	\arrow["{u_1''}", from=3-7, to=4-7]
	\arrow[""{name=3, anchor=center, inner sep=0}, "{g_2}"', from=4-3, to=4-5]
	\arrow["{u_2}", from=4-5, to=5-5]
	\arrow["{u_2'}", from=4-7, to=5-7]
	\arrow[""{name=4, anchor=center, inner sep=0}, "{g_3}"', from=5-5, to=5-7]
	\arrow["{u_3}", from=5-7, to=6-7]
	\arrow["{\widehat{\gamma}_{g_1, g'_2, g''_3}}", shorten <=4pt, shorten >=4pt, Rightarrow, from=0, to=2]
	\arrow["{\kappa_1}", shorten <=4pt, shorten >=4pt, Rightarrow, from=1, to=3]
	\arrow["{\kappa_2 \circ \kappa_3}", shorten <=9pt, shorten >=9pt, Rightarrow, from=2, to=4]
\end{tikzcd}\]

The vertical components of both composites are equal because of the following equalities in $X_0$:
\[\begin{tikzcd}
	{\widehat{\widehat{y'}}} && {\widehat{y'}} && {\widehat{c'}} \\
	\\
	{\widehat{y'}} && {y'} && {c'}
	\arrow[""{name=0, anchor=center, inner sep=0}, "{\gamma_{g_1,g_3'' \circ g_2'}}", from=1-1, to=1-3]
	\arrow["{\gamma_{g_2' \circ g_1,g_3''}}"', from=1-1, to=3-1]
	\arrow["{\widehat{u_1''}}", from=1-3, to=1-5]
	\arrow["{\gamma_{g_2',g_3''}}"{description}, from=1-3, to=3-3]
	\arrow["{(A)}"{description}, draw=none, from=1-5, to=3-3]
	\arrow["{\gamma_{g_2,g_3'}}", from=1-5, to=3-5]
	\arrow[""{name=1, anchor=center, inner sep=0}, "{d_0(\widehat{\gamma}_{g_1, g_2', g_3''})}"', from=3-1, to=3-3]
	\arrow["{u_1''}"', from=3-3, to=3-5]
	\arrow["{(B)}"{description}, draw=none, from=0, to=1]
\end{tikzcd}\]

\end{itemize}
In the diagram above, $(A)$ due to the naturality of $\gamma$ (see \eqref{gammanaturalita} in Example \ref{prclcat}). The equality $(B)$ is the "associativity law" of a codomain-colax category (see \eqref{eqcllaxalgassocALG} in Example \ref{prclcat} for the composable triple $(g_1,g_2',g_3'')$). The horizontal components are equal by Lemma \ref{THML_codcolax_are_assoc_unital_ABSTRAKTNI}.

Assume now that $\ce$ is arbitrary with pullbacks. Observe that by Lemma \ref{THML_simplicialobject_iff_cat}, $\cnr{X}$ is an internal category in $\ce$ if and only if for each $Z \in \ob \ce$, the diagram:
\[
\ce(Z,-) \circ \cnr{X} : \Delta^{op}_2 \to \set,
\]
corresponds to a small category. But this is the case since by Lemma \ref{THML_CnrFX_is_FCnrX}, we have:
\[
\ce(Z,-) \circ \cnr{X} \cong \cnr{\cat(\ce(Z,-))_* X},
\]
and the latter is a category by the first part of the proof.
\end{proof}

\noindent Lemma \ref{THML_CnrFX_is_FCnrX} can be thus expressed as follows: Given a domain-codiscrete codomain-colax category $X$ in $\cat(\ce)$, for any 2-functor $H: \cat(\ce) \to \cat(\ce')$ of form $\cat(H')$, we have:
\begin{equation}\label{EQ_CnrHX_cong_HCnrX}
\cnr{H_* X} \cong H\cnr{X}.
\end{equation}

\begin{con}\label{constrcodlaxcodcocone}
Let $X$ be a codomain-colax, domain-codiscrete category in $\cat(\ce)$. Define a graph morphism $F: X_0 \to \cnr{X}$ as follows:
\begin{align*}
F_0 &:= 1_{X_{00}},\\
F_1 &:= (s \circ s^{X_0}, \comp^{X_0} \circ (\iota \circ s^{X_0}, 1_{X_{01}})): X_{01} \to \cnr{X}_1.
\end{align*}
and define a 1-cell $\xi$ in $\ce$ as follows:
\[
\xi := (1_{X_{10}}, \ids^{X_0} \circ d_0): X_{00} \to \cnr{X}_1.
\]
We will denote this pair by $(F_X, \xi_X)$ to emphasize it is associated to the codomain-colax category $X$.
\end{con}

\begin{lemma}\label{THML_FHX_xiHX_is_HFX_HxiX}
Let $H: \cat(\ce) \to \cat(\ce')$ be a 2-functor of form $\cat(H')$. Under the identification \eqref{EQ_CnrHX_cong_HCnrX}, we have:
\begin{align}
F_{H_*X} &= H' F_X,\\
\xi_{H_* X} &= H'(\xi_X).
\end{align}
\end{lemma}

\begin{lemma}
Let $X$ be a codomain-colax, domain-codiscrete category in $\cat(\ce)$. The graph morphism $F: X_0 \to \cnr{X}$ is an internal functor and $\xi$ is an internal natural transformation $\xi: Fd_1 \Rightarrow Fd_0$.
\end{lemma}
\begin{proof}
Let us first prove the case $\ce = \set$. In $\cat$, $F$ is a mapping:
\[\begin{tikzcd}
	& a && a && {[a]} \\
	{F:} && \mapsto &&& a \\
	& b &&&& b
	\arrow["f"', from=1-2, to=3-2]
	\arrow["{s(a)}", from=1-4, to=1-6]
	\arrow["{\iota_a}", from=1-6, to=2-6]
	\arrow["f", from=2-6, to=3-6]
\end{tikzcd}\]
and the component of $\xi$ at a horizontal morphism $g: a \to b$ is the following corner in $\cnr{X}$:
\[\begin{tikzcd}
	& a && b \\
	{\xi_{g}:=} \\
	&&& b
	\arrow["g", from=1-2, to=1-4]
	\arrow[Rightarrow, no head, from=3-4, to=1-4]
\end{tikzcd}\]

\noindent Let us prove the following:
\begin{itemize}
\item \textbf{$F$ is a functor}: It clearly sends identities to identities. Next, consider the composition $Fv \circ Fu$ as in this picture:
\[\begin{tikzcd}
	a &&& {[a]} && {[a]} & {[a]} \\
	a & {[a]} && {[[a]]} & {[[a]]} & {[[a]]} \\
	& a && {[a]} & {[a]} & {[a]} & {[a]} \\
	& b && {[b]} & a & a & a \\
	&&& b & b & b & b \\
	&&& c & c & c & c
	\arrow["{s(a)}", from=1-1, to=1-4]
	\arrow[Rightarrow, no head, from=1-4, to=1-6]
	\arrow["{\gamma_{s(a), s([a])}}", from=1-4, to=2-4]
	\arrow[Rightarrow, no head, from=1-6, to=1-7]
	\arrow["{\gamma_{s(a), s([a])}}", from=1-6, to=2-6]
	\arrow[""{name=0, anchor=center, inner sep=0}, Rightarrow, no head, from=1-7, to=3-7]
	\arrow[""{name=0p, anchor=center, inner sep=0}, phantom, from=1-7, to=3-7, start anchor=center, end anchor=center]
	\arrow["{s(a)}", from=2-1, to=2-2]
	\arrow[""{name=1, anchor=center, inner sep=0}, "{s([a])}", from=2-2, to=2-4]
	\arrow["{\iota_a}"', from=2-2, to=3-2]
	\arrow[Rightarrow, no head, from=2-4, to=2-5]
	\arrow["{[\iota_a]}", from=2-4, to=3-4]
	\arrow[Rightarrow, no head, from=2-5, to=2-6]
	\arrow["{[\iota_a]}", from=2-5, to=3-5]
	\arrow["{\iota_{[a]}}", from=2-6, to=3-6]
	\arrow[""{name=2, anchor=center, inner sep=0}, "{s(a)}"{description}, from=3-2, to=3-4]
	\arrow["u"', from=3-2, to=4-2]
	\arrow[Rightarrow, no head, from=3-4, to=3-5]
	\arrow["{[u]}", from=3-4, to=4-4]
	\arrow["{(A)}"{description}, draw=none, from=3-5, to=3-6]
	\arrow["{\iota_a}", from=3-5, to=4-5]
	\arrow[Rightarrow, no head, from=3-6, to=3-7]
	\arrow["{\iota_a}", from=3-6, to=4-6]
	\arrow["{\iota_a}", from=3-7, to=4-7]
	\arrow[""{name=3, anchor=center, inner sep=0}, "{s(b)}"', from=4-2, to=4-4]
	\arrow["{(A)}"{description}, draw=none, from=4-4, to=4-5]
	\arrow["{\iota_b}"', from=4-4, to=5-4]
	\arrow[Rightarrow, no head, from=4-5, to=4-6]
	\arrow["u", from=4-5, to=5-5]
	\arrow["u", from=4-6, to=5-6]
	\arrow["u", from=4-7, to=5-7]
	\arrow[Rightarrow, no head, from=5-4, to=5-5]
	\arrow["v"', from=5-4, to=6-4]
	\arrow["v", from=5-5, to=6-5]
	\arrow["v", from=5-6, to=6-6]
	\arrow["v", from=5-7, to=6-7]
	\arrow[Rightarrow, no head, from=6-4, to=6-5]
	\arrow[Rightarrow, no head, from=6-5, to=6-6]
	\arrow[Rightarrow, no head, from=6-6, to=6-7]
	\arrow["{s(\iota_a)}", shorten <=4pt, shorten >=4pt, Rightarrow, from=1, to=2]
	\arrow["{(B)}"{description}, draw=none, from=2-6, to=0p]
	\arrow["{s(u)}", shorten <=4pt, shorten >=4pt, Rightarrow, from=2, to=3]
\end{tikzcd}\]

The squares labelled $(A)$ commute because of the naturality of $\iota$. The square $(B)$ commutes because of \eqref{eqcllaxalgunitALG2}. The horizontal pair composes to $s(c)$ by Lemma \ref{THML_codcolax_are_assoc_unital_ABSTRAKTNI}. In total we obtain that this corner equals $F(vu)$.

\item \textbf{$\xi$ is a natural transformation}: Consider a square $\alpha$ in $X$:
\[\begin{tikzcd}
	a && b \\
	\\
	c && d
	\arrow[""{name=0, anchor=center, inner sep=0}, "g", from=1-1, to=1-3]
	\arrow["u"', from=1-1, to=3-1]
	\arrow["v", from=1-3, to=3-3]
	\arrow[""{name=1, anchor=center, inner sep=0}, "h"', from=3-1, to=3-3]
	\arrow["\alpha", shorten <=9pt, shorten >=9pt, Rightarrow, from=0, to=1]
\end{tikzcd}\]

\noindent we wish to show that $Fv \circ \xi_g = \xi_{h} \circ Fu$ in $\cnr{X}$. We will show that both sides of this equation simplify to the corner $(g,v)$. The composites on both sides become:
\[\begin{tikzcd}
	a && b && a && b \\
	a & b & {[b]} && a & {[a]} & {\widehat{b}} \\
	& b & {[b]} &&& a & b \\
	&& b &&& c & d \\
	&& d &&&& d
	\arrow["g", from=1-1, to=1-3]
	\arrow["{\gamma_{g,s(b)}}", from=1-3, to=2-3]
	\arrow["g", from=1-5, to=1-7]
	\arrow["{\gamma_{s(a),\widehat{g}}}", from=1-7, to=2-7]
	\arrow["g", from=2-1, to=2-2]
	\arrow["{s(b)}", from=2-2, to=2-3]
	\arrow[Rightarrow, no head, from=2-2, to=3-2]
	\arrow[Rightarrow, no head, from=2-3, to=3-3]
	\arrow["{s(a)}", from=2-5, to=2-6]
	\arrow[""{name=0, anchor=center, inner sep=0}, "{\widehat{g}}", from=2-6, to=2-7]
	\arrow[""{name=0p, anchor=center, inner sep=0}, phantom, from=2-6, to=2-7, start anchor=center, end anchor=center]
	\arrow["{\iota_a}"', from=2-6, to=3-6]
	\arrow["{d_0(\widehat{\iota_g})}", from=2-7, to=3-7]
	\arrow["{s(b)}"', from=3-2, to=3-3]
	\arrow["{\iota_b}", from=3-3, to=4-3]
	\arrow[""{name=1, anchor=center, inner sep=0}, "g"{description}, from=3-6, to=3-7]
	\arrow[""{name=1p, anchor=center, inner sep=0}, phantom, from=3-6, to=3-7, start anchor=center, end anchor=center]
	\arrow[""{name=1p, anchor=center, inner sep=0}, phantom, from=3-6, to=3-7, start anchor=center, end anchor=center]
	\arrow["u"', from=3-6, to=4-6]
	\arrow["v", from=3-7, to=4-7]
	\arrow["v", from=4-3, to=5-3]
	\arrow[""{name=2, anchor=center, inner sep=0}, "h"', from=4-6, to=4-7]
	\arrow[""{name=2p, anchor=center, inner sep=0}, phantom, from=4-6, to=4-7, start anchor=center, end anchor=center]
	\arrow[Rightarrow, no head, from=4-7, to=5-7]
	\arrow["{\widehat{\iota_g}}", shorten <=4pt, shorten >=4pt, Rightarrow, from=0p, to=1p]
	\arrow["\alpha", shorten <=4pt, shorten >=4pt, Rightarrow, from=1p, to=2p]
\end{tikzcd}\]
The horizontal parts both compose to $g$ because of Lemma \ref{THML_codcolax_are_assoc_unital_ABSTRAKTNI}. The vertical part on the left side reduces to $v$ because of the axiom \eqref{eqcllaxalgunitALG2} and on the right side because of the axiom \eqref{eqcllaxalgunitALG1}.
\end{itemize}
The case for a general $\ce$ with pullbacks now follows from Lemma \ref{THML_FHX_xiHX_is_HFX_HxiX} and Lemma \ref{THML_graph_cell_is_category_cell}.
\end{proof}

\begin{lemma}
The pair $(F,\xi)$ from Construction \ref{constrcodlaxcodcocone} constitutes a cocone for the codomain-colax category $X$.
\end{lemma}
\begin{proof}
Let us first prove the case $\ce = \set$. We have to verify the cocone axioms \eqref{colaxcohcocone1}, \eqref{colaxcohcocone2}:
For all $(g,h) \in X_2$, $a \in X_0$:

\begin{align*}
F(\gamma_{(g,h)}) \circ \xi_{hg} &\sstackrel{?}{=} \xi_{h}\circ \xi_{g} \\
F(\iota_a) \circ \xi_{s(a)} &\sstackrel{?}{=} 1_{Fa}
\end{align*}
\noindent The first one holds since both sides are equal to $(h \circ g, \gamma_{g,h})$:
\[\begin{tikzcd}
	a && {\widehat{c}} & {\widehat{c}} \\
	a & {\widehat{c}} & {[\widehat{c}]} && a && {\widehat{c}} \\
	& {\widehat{c}} & {[\widehat{c}]} && a & b & c \\
	&& {\widehat{c}} & {\widehat{c}} && b & c \\
	&& c & c &&& c
	\arrow["{h \circ g}", from=1-1, to=1-3]
	\arrow[Rightarrow, no head, from=1-3, to=1-4]
	\arrow["{\gamma_{hg, s(\widehat{c})}}", from=1-3, to=2-3]
	\arrow[""{name=0, anchor=center, inner sep=0}, Rightarrow, no head, from=1-4, to=4-4]
	\arrow["{h \circ g}", from=2-1, to=2-2]
	\arrow["{s(\widehat{c})}", from=2-2, to=2-3]
	\arrow[Rightarrow, no head, from=2-2, to=3-2]
	\arrow[""{name=1, anchor=center, inner sep=0}, Rightarrow, no head, from=2-3, to=3-3]
	\arrow["{h \circ g}", from=2-5, to=2-7]
	\arrow["{\gamma_{g,h}}", from=2-7, to=3-7]
	\arrow["{s(\widehat{c})}"', from=3-2, to=3-3]
	\arrow["{\iota_{\widehat{c}}}", from=3-3, to=4-3]
	\arrow["g", from=3-5, to=3-6]
	\arrow["h", from=3-6, to=3-7]
	\arrow[Rightarrow, no head, from=3-6, to=4-6]
	\arrow[Rightarrow, no head, from=3-7, to=4-7]
	\arrow[Rightarrow, no head, from=4-3, to=4-4]
	\arrow["{\gamma_{g,h}}"', from=4-3, to=5-3]
	\arrow["{\gamma_{g,h}}"', from=4-4, to=5-4]
	\arrow["h"', from=4-6, to=4-7]
	\arrow[Rightarrow, no head, from=4-7, to=5-7]
	\arrow[Rightarrow, no head, from=5-3, to=5-4]
	\arrow["{\eqref{eqcllaxalgunitALG2}}"{description}, draw=none, from=1, to=0]
\end{tikzcd}\]

\noindent The second one holds because the composite on the left-hand side equals:
\[\begin{tikzcd}
	a && {[a]} && {[a]} \\
	a & {[a]} & {[[a]]} \\
	& {[a]} & {[[a]]} \\
	&& {[a]} && {[a]} \\
	&& a && a
	\arrow["{s(a)}", from=1-1, to=1-3]
	\arrow[Rightarrow, no head, from=1-3, to=1-5]
	\arrow["{\gamma_{s(a),s([a])}}", from=1-3, to=2-3]
	\arrow[""{name=0, anchor=center, inner sep=0}, Rightarrow, no head, from=1-5, to=4-5]
	\arrow["{s(a)}", from=2-1, to=2-2]
	\arrow["{s([a])}", from=2-2, to=2-3]
	\arrow[Rightarrow, no head, from=2-2, to=3-2]
	\arrow[""{name=1, anchor=center, inner sep=0}, Rightarrow, no head, from=2-3, to=3-3]
	\arrow["{s([a])}"', from=3-2, to=3-3]
	\arrow["{\iota_{[a]}}", from=3-3, to=4-3]
	\arrow[Rightarrow, no head, from=4-3, to=4-5]
	\arrow["{\iota_a}"', from=4-3, to=5-3]
	\arrow["{\iota_a}", from=4-5, to=5-5]
	\arrow[Rightarrow, no head, from=5-3, to=5-5]
	\arrow["{\eqref{eqcllaxalgunitALG2}}"{description}, draw=none, from=1, to=0]
\end{tikzcd}\]

The case of a general $\ce$ with pullbacks now follows from the fact the cocone axioms are equations in $\ce$, and those are jointly reflected by the collection $\ce(Z,-)$, $Z \in \ce$ (this is because the Yoneda embedding $y: \ce \to [\ce^{op},\set]$ is fully faithful).
\end{proof}

\begin{theo}\label{THM_catE_admits_codobj_of_domdiscr_codlax}
Let $\ck = \cat(\ce)$ for a category $\ce$ with pullbacks and let $X$ be a domain-codiscrete, codomain-colax category in $\ck$. Then the pair $(F,\xi)$ from Construction \ref{constrcodlaxcodcocone} is the codescent object of  $X$.
\end{theo}
\begin{proof}
Let us first prove this in the case of $\ce = \set$. Let $(G,\psi)$ be a cocone with apex $\cy$. This means that we have:
\begin{align}
G(\gamma_{g,g'}) \circ \psi_{g'g} &= \psi_{g'}\circ \psi_{g} \label{eqGcocone1ABSTRACT} \\
G(\iota_a) \circ \psi_{s(a)} &= 1_{Ga} \label{eqGcocone2ABSTRACT}
\end{align}

\noindent If there is to be a morphism of cocones $\theta: \cnr{X} \to \cy$, it must be true that:
\begin{align*}
\theta F(f) &= G(f) \\
\theta \xi_g &= \psi_g
\end{align*}
It is easy to see that every corner $(g,u)$ (as pictured below left) factorizes as the middle corner followed by the rightmost one:
\[\begin{tikzcd}
	a && b & a && b & b && {[b]} \\
	&&&&&&&& b \\
	&& c &&& b &&& c
	\arrow["g", from=1-1, to=1-3]
	\arrow["u", from=1-3, to=3-3]
	\arrow["g", from=1-4, to=1-6]
	\arrow[Rightarrow, no head, from=1-6, to=3-6]
	\arrow["{s(b)}", from=1-7, to=1-9]
	\arrow["{\iota_b}", from=1-9, to=2-9]
	\arrow["u", from=2-9, to=3-9]
\end{tikzcd}\]
In other words, $(g,u) = Fu \circ \xi_g$. Note that this factorization is unique in the sense that if $(g,u) = Fv \circ \xi_h$, we have $v = u$ and $h = g$. We are thus forced to put:
\[
\theta(g,u) = Gu \circ \psi_g,
\]
and this assignment is well-defined. It also sends identities to identities because:
\[
\theta(s(a),\iota_a) = G\iota_a \circ \psi_{s(a)} \hspace{0.2cm}\sstackrel{(\ref{eqGcocone2ABSTRACT})}{=}\hspace{0.2cm} 1_{Ga}
\]
To show that it preserves compositions, take a composable pair of corners $(g,u), (v,h)$ and their composite as pictured below:
\[\begin{tikzcd}
	a &&&& {\widehat{x}} \\
	a && b && x \\
	&& c && d \\
	&&&& e
	\arrow["{\widehat{h} \circ g}", from=1-1, to=1-5]
	\arrow["{\gamma_{g,\widehat{h}}}", from=1-5, to=2-5]
	\arrow["g", from=2-1, to=2-3]
	\arrow[""{name=0, anchor=center, inner sep=0}, "{\widehat{h}}", from=2-3, to=2-5]
	\arrow[""{name=0p, anchor=center, inner sep=0}, phantom, from=2-3, to=2-5, start anchor=center, end anchor=center]
	\arrow["u"', from=2-3, to=3-3]
	\arrow["{\widehat{u}}", from=2-5, to=3-5]
	\arrow[""{name=1, anchor=center, inner sep=0}, "h"', from=3-3, to=3-5]
	\arrow[""{name=1p, anchor=center, inner sep=0}, phantom, from=3-3, to=3-5, start anchor=center, end anchor=center]
	\arrow["v"', from=3-5, to=4-5]
	\arrow["{\exists !}", shorten <=4pt, shorten >=4pt, Rightarrow, from=0p, to=1p]
\end{tikzcd}\]

We compute:
\begin{align*}
\theta(h,v) \circ \theta(g,u) &= Gv \circ \psi_h \circ Gu \circ \psi_g \\
                              &\sstackrel{(*)}{=} Gv \circ G\widehat{u} \circ \psi_{\widehat{h}} \circ \psi_g \\
							  &\sstackrel{(\ref{eqGcocone1ABSTRACT})}{=}\hspace{0.2cm} Gv \circ G\widehat{u} \circ G(\gamma_{\widehat{h},g}) \circ \psi_{\widehat{h} g} \\
							  &= G(v \widehat{u} \gamma_{\widehat{h},g}) \psi_{\widehat{h} g} \\
							  &= \theta (\widehat{h} g, v \widehat{u} \gamma_{g,\widehat{h}})\\
							  &= \theta ((h,v) \circ (g,u)).
\end{align*}
In the above, $(*)$ is the naturality $\psi$.

\medskip

\noindent Let us now assume that $\ce$ is a general category with pullbacks. According to Corollary \ref{THMC_colimit_in_CatE_iff_in_Cat}, it is enough to prove that for every $Z \in \ob \ce$, the following pair:
\[
(\cat(\ce)(DZ,F_X),\cat(\ce)(DZ,\xi_X)),
\]
is the codescent object of the codomain-colax double category corresponding to the diagram $\cat(\ce)(DZ,-)_*X$. By Lemma \ref{THML_CnrFX_is_FCnrX} and Lemma \ref{THML_FHX_xiHX_is_HFX_HxiX}, we know that there is an isomorphism of cocones:
\[
(\cat(\ce)(DZ,F_X),\cat(\ce)(DZ,\xi_X)) \cong (F_{\cat(\ce(Z,-))_* X}, \xi_{\cat(\ce(Z,-))_* X}),
\]
and the latter is the codescent object by what we have proven in the first part of the proof.
\end{proof}

\begin{theo}\label{THM_CatT_preserves_special_codobjs}
Let $\ce, \ce'$ be categories with pullbacks and $H: \ce \to \ce'$ a pullback-preserving functor. Then the 2-functor $\cat(H)$ preserves:
\begin{itemize}
\item domain-codiscrete codomain-colax categories,
\item their codescent objects.
\end{itemize}
\end{theo}
\begin{proof}
This follows from Remark \ref{POZN_catF_preserves_Cnr}, \eqref{EQ_CnrHX_cong_HCnrX}, Lemma \ref{THML_FHX_xiHX_is_HFX_HxiX} and Theorem \ref{THM_catE_admits_codobj_of_domdiscr_codlax}.
\end{proof}

\begin{cor}\label{THMC_CatT_preserves_codobj_of_resAa}
Let $\ce$ be a category with pullbacks. Any 2-monad $T$ on $\cat(\ce)$ that is of form $\cat(T)$ preserves codescent objects of resolutions of colax $T$-algebras.
\end{cor}
\begin{proof}
Resolution of a colax $T$-algebra is domain-codiscrete, codomain-colax category by Example \ref{prdomdiscrcodlax}.
\end{proof}

\subsection{Grothendieck construction and codescent}\label{SUBS_grothendieck}

The main source of examples of codomain-colax categories are the resolutions of colax algebras for a cartesian 2-monad $T$ (Example \ref{PR_cartesian2monad_res_codcolaxcat}). Let us give a nontrivial example that is not of this form -- a codomain-colax double category associated to a lax functor.

Let $\ce$ be a 1-category and $P : \ce \to \cat$ be a colax functor. There is a codomain-colax double category associated to $P$ which we denote by $X_P$. It has:
\begin{itemize}
\item objects the pairs $(I,x)$ with $I \in \ob \ce, x \in PI$,
\item a vertical morphism $v: (I,x) \to (I,y)$ is a morphism $v: x \to y$ in $P(I)$,
\item a horizontal morphism $\alpha : (I,x) \to (J, y)$ is a morphism $\alpha : I \to J$ in $\ce$ with the property that $P(\alpha)(x) = y$,
\item a square below left exists if and only if the condition below right is satisfied:
\[\begin{tikzcd}[column sep=scriptsize]
	{(I,x)} && {(J,P(\alpha)(x))} \\
	&&& \iff & {\alpha = \beta \phantom{X}\land \phantom{X} v = P(\alpha)(u)} \\
	{(I,y)} && {(J,P(\beta)(y))}
	\arrow[""{name=0, anchor=center, inner sep=0}, "\alpha", from=1-1, to=1-3]
	\arrow["u"', from=1-1, to=3-1]
	\arrow["v", from=1-3, to=3-3]
	\arrow[""{name=1, anchor=center, inner sep=0}, "\beta"', from=3-1, to=3-3]
	\arrow["\exists", shorten <=13pt, shorten >=13pt, Rightarrow, from=0, to=1]
\end{tikzcd}\]
We can see that such a square is uniquely determined by its bottom-left corner, so the this codomain-colax category is domain-codiscrete.
\item the horizontal composition and the coassociators are portrayed below (as always, $\gamma$ is the coassociator for the colax functor):
\[\begin{tikzcd}[column sep=tiny]
	{(I,x)} &&&& {(K,P(\beta \circ \alpha)(x))} && {(K,P(\delta)P(\beta \circ \alpha)(x))} \\
	\\
	{(I,x)} && {(J,P(\alpha)(x))} && {(K,P(\beta)P(\alpha)(x))} && {(K,P(\delta)P(\beta)P(\alpha)(x))}
	\arrow["{\beta \circ \alpha}", from=1-1, to=1-5]
	\arrow[""{name=0, anchor=center, inner sep=0}, "\delta", from=1-5, to=1-7]
	\arrow[""{name=0p, anchor=center, inner sep=0}, phantom, from=1-5, to=1-7, start anchor=center, end anchor=center]
	\arrow["{\gamma_{\alpha,\beta,x}}"', from=1-5, to=3-5]
	\arrow["{P(\delta)(\gamma_{\alpha,\beta,x})}", from=1-7, to=3-7]
	\arrow["\alpha"', from=3-1, to=3-3]
	\arrow["\beta"', from=3-3, to=3-5]
	\arrow[""{name=1, anchor=center, inner sep=0}, "\delta"', from=3-5, to=3-7]
	\arrow[""{name=1p, anchor=center, inner sep=0}, phantom, from=3-5, to=3-7, start anchor=center, end anchor=center]
	\arrow[shorten <=9pt, shorten >=9pt, Rightarrow, from=0p, to=1p]
\end{tikzcd}\]
\item The counitors and the horizontal identity are portrayed below (they use the counitor of $P$):
\[\begin{tikzcd}
	{(I,x)} && {(I,P(1_I)(x))} && {(I,P(\alpha)P(1_I)(x))} \\
	\\
	&& {(I,x)} && {(J,P(\alpha)(x))}
	\arrow["{1_I}", from=1-1, to=1-3]
	\arrow[""{name=0, anchor=center, inner sep=0}, "\alpha", from=1-3, to=1-5]
	\arrow[""{name=0p, anchor=center, inner sep=0}, phantom, from=1-3, to=1-5, start anchor=center, end anchor=center]
	\arrow["{\iota_{I,x}}"', from=1-3, to=3-3]
	\arrow["{P(\alpha)(\iota_{I,x})}", from=1-5, to=3-5]
	\arrow[""{name=1, anchor=center, inner sep=0}, "\alpha"', from=3-3, to=3-5]
	\arrow[""{name=1p, anchor=center, inner sep=0}, phantom, from=3-3, to=3-5, start anchor=center, end anchor=center]
	\arrow[shorten <=9pt, shorten >=9pt, Rightarrow, from=0p, to=1p]
\end{tikzcd}\]
\end{itemize}

\noindent The codescent object has the objects pairs $(I, x \in PI)$, with a morphism $(I,x) \to (J,y)$ being the corner below left, or equivalently, the lax triangle in $\cat$ portrayed below right:
\[\begin{tikzcd}
	{(I,x)} && {(J,P(\alpha)(x))} &&& {*} \\
	\\
	&& {(J,y)} && {P(I)} && {P(J)}
	\arrow["\alpha", from=1-1, to=1-3]
	\arrow["u", from=1-3, to=3-3]
	\arrow["x"', curve={height=12pt}, from=1-6, to=3-5]
	\arrow[""{name=0, anchor=center, inner sep=0}, "y", curve={height=-12pt}, from=1-6, to=3-7]
	\arrow["{P(\alpha)}"', from=3-5, to=3-7]
	\arrow["u", shorten <=10pt, shorten >=10pt, Rightarrow, from=3-5, to=0]
\end{tikzcd}\]
This is the \textit{Grothendieck construction} associated to the colax functor, see e.g. \cite[Definition 1.3.1]{elephant}. We will denote it by $\int F$. It is well known that $\int F$ is the lax colimit of $F$ -- as we have just seen, it is also an ordinary 2-colimit -- the codescent object.

Further, notice that the category $\ce$ is an apex for a cocone $(G,\psi)$ for $X_P$. The functor $G: (X_P)_0 \to \ce$ is \textit{componentwise constant}, sending every morphism to the identity morphism:
\[\begin{tikzcd}
	& {(I,x)} \\
	{G:} && \mapsto & {1_I} \\
	& {(I,y)}
	\arrow["u"', from=1-2, to=3-2]
\end{tikzcd}\]
The natural transformation $\psi$, evaluated at $\alpha: (I,x) \to (J,y)$ (which is a horizontal morphism of $X_P$) is just the morphism $\alpha: I \to J$. The cocone axioms are an easy exercise. By the universal property of the codescent object, we obtain the existence of a unique functor $\pi: \int F \to \ce$, which is of course the usual projection functor that can be later shown to be an opfibration.

Note that this example is not the same as the one to be given in Example \ref{PR_reslan_laxcoh_explicit}, where we also associate a domain-codiscrete, codomain-colax category to a colax functor $F: \ce \to \cat$: the former is a codomain-colax category in $\cat$, while the latter lives in $[\ce, \cat]$.

Also, compare this example to Example \ref{pr_xp_fib} where we associated a crossed double category to a fibration -- presumably codomain-colax double categories are to crossed double category as lax functors are to fibrations.

\subsection{Colimit formula in $\cat$}\label{SEC_thecaseckCat}

Any codescent object can be always built using coinserters and coequifiers (see Remark \ref{POZN_technical_lacks_coh_data}). Using the explicit description of those 2-colimits in $\cat$ (Subsection \ref{SUBS_colimits_in_Cat}), one can deduce the formula for codescent objects of codomain-colax double categories in $\cat$ using generators and relations -- this is what we will do in this subsection.

\begin{con}\label{CON_Cnr_pro_codomaincolaxcat}
Let $X$ be a codomain-colax category in $\cat$. Define the category $F\cnr{X}$ as follows. The objects are the objects of $X_0$, while the morphisms are ``composable paths” $(f_1, \dots, f_n)$ where each $f_i$ is either a vertical or a horizontal morphism in $X$. The identity morphism on each object is given by the empty path, and the composition is given by concatenation.

There is a congruence on $F\cnr{X}$ generated by the following data:
\begin{itemize}
\item $(g_1,g_2) \sim (g_2 \circ g_1, \gamma_{g_1,g_2})$ if $g_1$, $g_2$ are both horizontal,
\item $(s(a),\iota_a) \sim ()_a$ for each object $a \in X$,
\item $(g,v) \sim (u,h)$ if there is a square $\alpha$ as pictured below:
\[\begin{tikzcd}
	a & b \\
	c & d
	\arrow[""{name=0, anchor=center, inner sep=0}, "g", from=1-1, to=1-2]
	\arrow[""{name=0p, anchor=center, inner sep=0}, phantom, from=1-1, to=1-2, start anchor=center, end anchor=center]
	\arrow["u"', from=1-1, to=2-1]
	\arrow["v", from=1-2, to=2-2]
	\arrow[""{name=1, anchor=center, inner sep=0}, "h"', from=2-1, to=2-2]
	\arrow[""{name=1p, anchor=center, inner sep=0}, phantom, from=2-1, to=2-2, start anchor=center, end anchor=center]
	\arrow["\alpha", shorten <=4pt, shorten >=4pt, Rightarrow, from=0p, to=1p]
\end{tikzcd}\]
\item $(f_1, f_2) \sim (f_2 \circ f_1)$ if both $f_1, f_2$ are vertical,
\item $(1_a) \sim ()_a$ for the vertical identity on each $a \in X$.
\end{itemize}
Denote $\cnr{X} := F\cnr{X}/\sim$.
\end{con}

\begin{theo}
The above category is the codescent object of $X$:
\[
\cnr{X} = \cod{X}.
\]
\end{theo}
\begin{proof}
Define $F: X_0 \to \cnr{X}$ as the functor that is the identity on objects and sends a morphism $f$ to the equivalence class $[f]$. Next, define the natural transformation \[ \xi: Fd_1 \Rightarrow Fd_0: X_1 \to \cnr{X}, \] by: $\xi_g := [g]$. The naturality of $\xi$ as well as the cocone axioms $F\gamma_{g_1,g_2} \circ \xi_{g_2 g_1} = \xi_{g_2} \circ \xi_{g_1}$, $F\iota_a \circ \xi_{s(a)} = 1_{Fa}$ follow from the defining congruence for $\cnr{X}$.

To prove the \textbf{one-dimensional universal property}, let $(G,\psi)$ be a cocone with apex $\cc$. The equations for the comparison functor $\theta: \cnr{X} \to \cc$ force that:
\begin{alignat*}{2}
\theta (a) &= Ga &\phantom{...}& \text{ if }a \in \ob X_0,\\
\theta (f) &= Gf && \text{ if }f\text{ is vertical,}\\
\theta (g) &= \psi_g && \text{ if }g\text{ is horizontal.}
\end{alignat*}
Let us show that this is well-defined on equivalence classes. The facts that $\theta [g_1, g_2] = \theta[g_2 \circ g_1, \gamma_{g_1, g_2}]$, $\theta [s(a), \iota_a] = 1_{Ga}$ follow from the cocone axioms for $(G,\psi)$. The fact that $\theta [g,v] = \theta [u,h]$ whenever there is a square $\alpha$ present (as in \ref{CON_Cnr_pro_codomaincolaxcat}) follows from the naturality of $\psi$. The assignment $\theta$ also preserves composites since they are done by concatenation. Preservation of identities is clear.

Again, since $\cat$ admits powers with $\boldsymbol{2}$, the \textbf{two-dimensional universal property} follows automatically (see Proposition \ref{THMP_powers_with2_implies_1dim_impl_2dim} and Theorem \ref{THM_catE_limits}).
\end{proof}

\begin{pozn}
For completeness, we also provide the formula for Lack's codescent object of lax coherence data (see Remark \ref{POZN_technical_lacks_coh_data}). Given a codomain-colax category $X$, the codescent object in the sense of Lack can be described as follows. Take again the category $F\cnr{X}$ from the previous construction, but this time build the congruence as follows:
\begin{itemize}
\item $(g_2 g_1) \sim (g_1, g_2, \gamma_{g_1,g_2})$ if $g_1, g_2$ are both horizontal,
\item $(s(a)) \sim (\iota_a)$,
\item $(g,v) \sim (u,h)$ if there is a square $\alpha$ as pictured below:
\[\begin{tikzcd}
	a & b \\
	c & d
	\arrow[""{name=0, anchor=center, inner sep=0}, "g", from=1-1, to=1-2]
	\arrow[""{name=0p, anchor=center, inner sep=0}, phantom, from=1-1, to=1-2, start anchor=center, end anchor=center]
	\arrow["u"', from=1-1, to=2-1]
	\arrow["v", from=1-2, to=2-2]
	\arrow[""{name=1, anchor=center, inner sep=0}, "h"', from=2-1, to=2-2]
	\arrow[""{name=1p, anchor=center, inner sep=0}, phantom, from=2-1, to=2-2, start anchor=center, end anchor=center]
	\arrow["\alpha", shorten <=4pt, shorten >=4pt, Rightarrow, from=0p, to=1p]
\end{tikzcd}\]
\item $(f_1, f_2) \sim (f_2 \circ f_1)$ if both $f_1, f_2$ are vertical,
\item $(1_a) \sim ()_a$ for the vertial identity on each $a \in X$.
\end{itemize}

\begin{prop}
The above category is the codescent object (in the sense of Steve Lack) of $X$:
\[
\cnr{X} = \cod{X}.
\]
\end{prop}
\end{pozn}

Let us end with a small observation that tells us that given certain categories with structure (for instance strict monoidal categories), any functor $F$ between them that preserves the structure in a lax way can equivalently be encoded as a pair of ordinary functors, with the latter one in the pair measuring to what extent $F$ preserves the structure:

\begin{pozn}
Given any cartesian 2-monad $T$ on $\cat$ and a lax morphism between two strict $T$-algebras:
\[(f,\overline{f}): (A,a) \to (B,b),\]
by the proof of Theorem \ref{THM_adjunkce_talgs_talgl_colaxtalgl}, this is equivalently a cocone $(bTf, bT\overline{f})$ for the coherence data $\res{A,a}$ in $\talgs$. Since $T$ is cartesian, $\res{A,a}$ is a double category (see Example \ref{PR_res_of_alg}). In particular, by Proposition \ref{THMP_cohX_cong_pairX} this cocone is equivalently a pair of functors:
\begin{align*}
F &: TA \to B,\\
\xi &: h(\res{A,a}) \to B,
\end{align*}
satisfying the naturality condition (Construction \ref{CON_pair}). In particular, the lax morphism $(f,\overline{f})$ is strict if and only if the functor $\xi$ is \textit{componentwise constant}, i.e. sends every morphism to an identity morphism.
\end{pozn}

\section{Lax coherence for 2-monads of form $\cat(T)$}\label{SEC_application_lax_coh}

\begin{prop}\label{THMP_catT_strictification_exists}
Let $(T,m,i)$ be a 2-monad of form $\cat(T')$ for a cartesian monad $T'$ on a category $\ce$ with pullbacks. The inclusion 2-functors admit left 2-adjoints:
\[\begin{tikzcd}
	\talgs && \talgl && \talgs && \colaxtalgl
	\arrow[""{name=0, anchor=center, inner sep=0}, curve={height=24pt}, hook, from=1-1, to=1-3]
	\arrow[""{name=1, anchor=center, inner sep=0}, curve={height=24pt}, dashed, from=1-3, to=1-1]
	\arrow[""{name=2, anchor=center, inner sep=0}, curve={height=24pt}, hook, from=1-5, to=1-7]
	\arrow[""{name=3, anchor=center, inner sep=0}, curve={height=24pt}, dashed, from=1-7, to=1-5]
	\arrow["\dashv"{anchor=center, rotate=-90}, draw=none, from=1, to=0]
	\arrow["\dashv"{anchor=center, rotate=-90}, draw=none, from=3, to=2]
\end{tikzcd}\]
\end{prop}
\begin{proof}

By Theorem \ref{THM_adjunkce_talgs_talgl_colaxtalgl}, these 2-adjoints are given by computing the codescent objects of resolutions in $\talgs$ of colax $T$-algebras. By Proposition \ref{prdomdiscrcodlax}, each resolution is a domain-codiscrete, codomain-colax category. $\cat(\ce)$ admits codescent objects of these by Theorem \ref{THM_catE_admits_codobj_of_domdiscr_codlax}.
 Finally, since $T = \cat(T')$ preserves these codescent objects by Theorem \ref{THM_CatT_preserves_special_codobjs}, they uniquely lift\footnote{This is a $\cat$-enriched version of \cite[Theorem 5.6.5.]{context}.} along the forgetful 2-functor $U: \talgs \to \cat(\ce)$ so they exist in $\talgs$.
\end{proof}

\begin{prop}\label{THMP_CatT_Ql_lax_idemp}
Let $(T,m,i)$ be a 2-monad of form $\cat(T')$ for a cartesian monad $T'$ on a category $\ce$ with pullbacks. Then the lax morphism classifier 2-comonad $Q_l$ on $\talgs$ is lax-idempotent.
\end{prop}
\begin{proof}
This follows from Proposition \ref{THMP_Ql_lax_idemp_Qp_ps_idemp}. In fact, \textbf{both} of the conditions are true in this case -- $\cat(\ce)$ admits oplax limits of arrows by Corollary \ref{THM_catE_limits}, and $T$ preserves the relevant codescent objects by Theorem \ref{THM_CatT_preserves_special_codobjs}.
\end{proof}

\begin{theo}\label{THM_coherence_for_colax_algs_CatT}
Let $(T,m,i)$ be a 2-monad of form $\cat(T')$ for a cartesian monad $T'$ on a category $\ce$ with pullbacks. Then coherence for colax algebras holds for $T$.
\end{theo}
\begin{proof}
By Proposition \ref{THMP_catT_strictification_exists}, colax algebras can be strictified. By Corollary \ref{THMC_CatT_preserves_codobj_of_resAa}, $T$ preserves codescent objects of resolutions of colax algebras and so the result follows from Theorem \ref{THM_lax_coherence_colax_algebras_EXPOSITION}.
\end{proof}

The reader may want to compare this result with the coherence result for pseudo algebras for this class of 2-monads that has been established in the literature -- we will do so in the remainder of this section. Recall that it has been proven in \cite[Theorem 3.4]{ageneralcohresult} that for any 2-monad $T$ on $\cat$ that preserves functors that are bijective on objects, every pseudo algebra is equivalent to a strict one. In \cite[Theorem 7.1]{eckmann_hilton}, it has been observed that this readily generalizes to 2-monads on $\cat(\ce)$ that preserve internal functors that are ``bijective on objects”, meaning the objects-component 1-cell is an isomorphism in $\ce$. From \cite[Theorem 4.6]{codobjcoh} it follows that this situation implies that the \textit{coherence for pseudo algebras} holds for this class of 2-monads, meaning the following:

\begin{theo}
Let $\ce$ be a category with pullbacks and let $T$ be a 2-monad on $\cat(\ce)$ which preserves functors that are bijective on objects. Then the inclusion of strict algebras to pseudo algebras admits a left 2-adjoint:
\[\begin{tikzcd}
	\talgs && \pstalg
	\arrow[""{name=0, anchor=center, inner sep=0}, curve={height=24pt}, hook, from=1-1, to=1-3]
	\arrow[""{name=1, anchor=center, inner sep=0}, "{(-)^{\dagger}}"', curve={height=24pt}, dashed, from=1-3, to=1-1]
	\arrow["\dashv"{anchor=center, rotate=-90}, draw=none, from=1, to=0]
\end{tikzcd}\]
Moreover, each component of the unit of this 2-adjunction is an equivalence in $\cat(\ce)$ (and in particular in $\pstalg$ by Corollary \ref{THMC_doct_adj_corollary}).
\end{theo}

Any 2-monad of form $\cat(T)$ preserves bijections on objects. In particular, looking at the examples from Subsection \ref{SUBS_examples_CatT}, the following structures are all canonically equivalent to their strict version: monoidal categories, bicategories, pseudo double categories, pseudo weights (pseudofunctors $\cj \to \cat$, with $\cj$ a 1-category).

\subsection{Examples}\label{SEC_examples_lax_coh}

\begin{nota}
In the examples below, given a colax $T$-algebra $(A,a,\gamma,\iota)$ for a 2-monad of form $\cat(T)$, we will denote:
\[
\cnr{A} := \cnr{\res{A,a,\gamma,\iota}}.
\]
\end{nota}

Each of the 2-monads below is of form $\cat(T)$ so coherence for colax algebras and coherence for lax morphisms holds for them:

\begin{pr}[Coherence for colax monoidal categories]\label{PR_freestrictmoncat2monad_coherence}
Consider the free strict monoi\-dal category 2-monad $T$ on $\cat$ (Example \ref{PR_freestrictmoncat2monad}). and a colax monoidal category $(\ca,a,\gamma,\iota)$. Objects and vertical morphisms of $\res{\ca,a,\gamma,\iota}$ are the objects and morphisms of $T\ca$. A horizontal morphism is a \textit{partial evaluation} (\cite{parteval}), that is, an object $\alpha \in \ob T^2\ca$ whose domain is given by $m_{\ca}(\alpha)$ and whose codomain is given by $(Ta)(\alpha)$. For instance:
\[\begin{tikzcd}
	{(a_1,a_2,a_3)} && {(a_1\otimes a_2, I, [a_3])}
	\arrow["{((a_1,a_2),(),(a_3))}", from=1-1, to=1-3]
\end{tikzcd}\]
The strict monoidal category $\cnr{\ca}$ has tuples of objects of $\ca$ as \textbf{objects}, with \textbf{morphisms} consisting of two components: a partial evaluation and a tuple of morphisms of $\ca$:
\[\begin{tikzcd}
	{(a_1, a_2, a_3)} &&& {(I, a_1 \otimes a_2, I, [a_3])} \\
	&&& {(b_1, b_2, b_3, b_4)}
	\arrow["{((),(a_1, a_2),(),(a_3))}", from=1-1, to=1-4]
	\arrow["{(f_1,f_2,f_3,f_4)}", from=1-4, to=2-4]
\end{tikzcd}\]
The \textbf{identity morphism} is given by the following corner:
\[\begin{tikzcd}
	{(a_1, \dots, a_n)} && {([a_1], \dots, [a_n])} \\
	&& {(a_1, \dots, a_n)}
	\arrow["{((a_1), \dots, (a_n))}", from=1-1, to=1-3]
	\arrow["{(\iota_{a_1}, \dots, \iota_{a_n})}", from=1-3, to=2-3]
\end{tikzcd}\]

For an example of a \textbf{composition} of two corners, consider the pair:
\[\begin{tikzcd}
	{(a_1, a_2, a_3)} && {(a_1 \otimes a_2,I, [a_3])} \\
	&& {(b_1, b_2, b_3)} && {([b_1], b_2 \otimes b_3)} \\
	&&&& {(c_1, c_2)}
	\arrow["{((a_1, a_2),(),(a_3))}", from=1-1, to=1-3]
	\arrow["{(f_1, f_2, f_3)}"', from=1-3, to=2-3]
	\arrow["{((b_1),(b_2,b_3))}"', from=2-3, to=2-5]
	\arrow["{(g_1, g_2)}", from=2-5, to=3-5]
\end{tikzcd}\]

\noindent Using the fact that the multiplication is a discrete opfibration, the bottom-left corner admits a unique filler into a square (a morphism in $T^2\ca$):
\[\begin{tikzcd}[cramped]
	{(a_1, a_2, a_3)} && {(a_1 \otimes a_2,I, [a_3])} && {([a_1 \otimes a_2],I \otimes [a_3])} \\
	&& {(b_1, b_2, b_3)} && {([b_1], b_2 \otimes b_3)} \\
	&&&& {(c_1, c_2)}
	\arrow["{((a_1, a_2),(),(a_3))}", from=1-1, to=1-3]
	\arrow["{(f_1, f_2, f_3)}"', from=1-3, to=2-3]
	\arrow["{((b_1),(b_2,b_3))}"', from=2-3, to=2-5]
	\arrow["{(g_1, g_2)}", from=2-5, to=3-5]
	\arrow["{((a_1 \otimes a_2),(I, [a_3]))}", dashed, from=1-3, to=1-5]
	\arrow["{([f_1], f_2 \otimes f_3)}"', dashed, from=1-5, to=2-5]
\end{tikzcd}\]

\noindent The final step consists of composing the partial evaluations and appropriately compensating for the fact that the horizontal codomain changed by pre-composing the vertical component with the components of the coassociator for $\ca$:
\[\begin{tikzcd}
	{(a_1, a_2, a_3)} &&&& {(a_1\otimes a_2, [a_3])} \\
	&&&& {([a_1 \otimes a_2],I \otimes [a_3])} \\
	\\
	&&&& {(c_1, c_2)}
	\arrow["{((a_1,a_2),(a_3))}", from=1-1, to=1-5]
	\arrow["{(\gamma_{((a_1,a_2))}, \gamma_{((),(a_3))})}"', from=1-5, to=2-5]
	\arrow["{(g_1 \circ [f_1], g_2 \circ f_2 \otimes f_3)}"', from=2-5, to=4-5]
\end{tikzcd}\]

The \textbf{strict monoidal structure} of $\cnr{\ca}$ is given by concatenation, we will denote it by $\boxplus: T\cnr{\ca} \to \cnr{\ca}$. The \textbf{unit} $(P,\overline{P}): (\ca, a, \gamma, \iota) \to (\cnr{\ca},\boxplus)$ is a lax monoidal functor whose underlying functor sends:
\[\begin{tikzcd}
	a && {(a)} && {([a])} \\
	&&&& {(a)} \\
	& \mapsto \\
	b &&&& {(b)}
	\arrow["f"', from=1-1, to=4-1]
	\arrow["{s(a)}", from=1-3, to=1-5]
	\arrow["{(\iota_a)}", from=1-5, to=2-5]
	\arrow["{(f)}", from=2-5, to=4-5]
\end{tikzcd}\]
The associator $\overline{P}$, evaluated for instance at a tuple $(a_1,a_2,a_3)$, is a morphism:
\[
Pa_1 \boxplus Pa_2 \boxplus Pa_3 = (a_1,a_2,a_3) \to P(a_1 \otimes a_2 \otimes a_3) \text{ in }\cnr{\ca}.
\]
It is given by this corner:
\[\begin{tikzcd}
	{(a_1,a_2,a_3)} && {(a_1 \otimes a_2 \otimes a_3)} \\
	&& {(a_1 \otimes a_2 \otimes a_3)}
	\arrow["{((a_1,a_2,a_3))}", from=1-1, to=1-3]
	\arrow[Rightarrow, no head, from=1-3, to=2-3]
\end{tikzcd}\]

\noindent In the adjunction guaranteed by the coherence for colax algebras:
\[\begin{tikzcd}
	\ca && {\cnr{\ca}} & {\text{ in }\cat}
	\arrow[""{name=0, anchor=center, inner sep=0}, "P"', curve={height=24pt}, from=1-1, to=1-3]
	\arrow[""{name=1, anchor=center, inner sep=0}, "Q"', curve={height=24pt}, from=1-3, to=1-1]
	\arrow["\dashv"{anchor=center, rotate=-90}, draw=none, from=1, to=0]
\end{tikzcd}\]
The \textbf{left adjoint} $Q$ uses the coassociator of $\ca$, and for instance sends:
\[\begin{tikzcd}[column sep=small]
	{(a_1, a_2, a_3)} &&& {(I, a_1 \otimes a_2, I, [a_3])} && {a_1 \otimes a_2 \otimes a_3} \\
	&&&& \mapsto & {I \otimes (a_1 \otimes a_2) \otimes I \otimes [a_3]} \\
	&&& {(b_1, b_2, b_3, b_4)} && {b_1 \otimes b_2 \otimes b_3 \otimes b_4}
	\arrow["{((),(a_1, a_2),(),(a_3))}", shift left=3, from=1-1, to=1-4]
	\arrow["{(f_1,f_2,f_3,f_4)}"', from=1-4, to=3-4]
	\arrow["{\gamma_{((),(a_1, a_2),(),(a_3))}}", from=1-6, to=2-6]
	\arrow["{f_1 \otimes f_2 \otimes f_3 \otimes f_4}", from=2-6, to=3-6]
\end{tikzcd}\]
It can be either seen explicitly or deduced from Corollary \ref{THMC_coherence_for_normal_colax} that:
\begin{cor}\label{THMC_coherence_normal_colax_moncat}
Every category with normal colax monoidal structure can be reflectively embedded into a category with strict monoidal structure such that the embedding is lax monoidal and the reflector is colax monoidal.
\end{cor}
\end{pr}

To make the exposition easier to digest, in the next example we will consider only the strictification of lax morphisms, and not colax algebras:

\begin{pr}[Coherence for lax functors between 2-categories]\label{prcodlax-colaxfunclass}
Fix a set $X$ and consider the free 2-category 2-monad $T$ on $\set^{X \times X}$ from Example \ref{PR_2category_2monad}. Any 2-category $\cc$ with object-set $X$ (regarded as a $T$-algebra) gives rise to its resolution $U\res{\cc}$, which is a diagram in $\cat^{X \times X}$. Denote its codescent object by $\cnr{\cc}$ -- this is a $\cat$-graph with the set of objects being $X$ that moreover has the structure of a 2-category.

As colimits in $\cat^{X \times X}$ are computed pointwise, $\cnr{\cc}(x,y)$ is given by the codescent object of $\res{\cc}(x,y)$. Because the 2-monad multiplication $m: T^2 \Rightarrow T$ is a pointwise discrete opfibration, each $\res{\cc}(x,y)$ is a domain-codiscrete double category and so $\cnr{\cc}(x,y)$ can be computed using the category of corners construction as follows. 

Objects in $\cnr{\cc}(x,y)$ are the objects of $T\cc(x,y)$, that is, paths of morphisms in the 2-category $\cc$. Morphisms are corners whose first component is given by a \textit{partial evaluation 2-cell} (an object of $T^2\cc(x,y)$) and the second component is given by a tuple of 2-cells in $\cc$ (a morphism in $T\cc(x,y)$). For instance this morphism $(f_1, f_2, f_3,f_4) \to (g_1, g_2, g_3)$:

\begin{equation}\label{EQ_morphism_in_cnrCC}
\begin{tikzcd}
	{a_1} && {a_2} && {a_3} && {a_4} && {a_5} \\
	&&& {a_3} && {a_5}
	\arrow["{f_1}", from=1-1, to=1-3]
	\arrow[""{name=0, anchor=center, inner sep=0}, "{f_2 \circ f_1}"{pos=0.7}, from=1-1, to=2-4]
	\arrow[""{name=1, anchor=center, inner sep=0}, "{g_1}"', curve={height=30pt}, from=1-1, to=2-4]
	\arrow["{f_2}", from=1-3, to=1-5]
	\arrow["{f_3}", from=1-5, to=1-7]
	\arrow["{f_4}", from=1-7, to=1-9]
	\arrow[""{name=2, anchor=center, inner sep=0}, "{g_2}"', curve={height=30pt}, from=2-4, to=2-6]
	\arrow[""{name=3, anchor=center, inner sep=0}, "{f_4 \circ f_3}"{description}, from=2-4, to=2-6]
	\arrow[""{name=4, anchor=center, inner sep=0}, Rightarrow, no head, from=2-6, to=1-9]
	\arrow[""{name=5, anchor=center, inner sep=0}, "{g_3}"', curve={height=30pt}, from=2-6, to=1-9]
	\arrow["{\alpha_1}", shorten <=4pt, shorten >=4pt, Rightarrow, from=0, to=1]
	\arrow["{((f_1,f_2),(f_3,f_4),())}"{pos=0.4}, shorten <=3pt, shorten >=3pt, Rightarrow, from=1-5, to=3]
	\arrow["{\alpha_2}", shorten <=4pt, shorten >=4pt, Rightarrow, from=3, to=2]
	\arrow["{\alpha_3}", shorten <=4pt, shorten >=4pt, Rightarrow, from=4, to=5]
\end{tikzcd}
\end{equation}
The 2-category structure of $\cnr{\cc}$ is given by concatenation of paths and tuples of 2-cells.

The canonical identity-on-objects lax functor $P: \cc \to \cnr{\cc}$ is given on morphisms by the assignment $f \mapsto (f)$. Its action on 2-cells is given by $\alpha \mapsto (\alpha)$. It is a lax functor -- its associator and unitor again use the additional ``partial evaluation” 2-cells present in $\cnr{\cc}$. In the colax coherence adjunction:
\[\begin{tikzcd}
	\cc && {\cnr{\cc}} & {\text{ in }\twocat_l}
	\arrow[""{name=0, anchor=center, inner sep=0}, "P"', curve={height=24pt}, from=1-1, to=1-3]
	\arrow[""{name=1, anchor=center, inner sep=0}, "Q"', curve={height=24pt}, from=1-3, to=1-1]
	\arrow["\dashv"{anchor=center, rotate=-90}, draw=none, from=1, to=0]
\end{tikzcd}\]
$Q: \cnr{\cc} \to \cc$ is a 2-functor and for instance sends the 2-cell \eqref{EQ_morphism_in_cnrCC} in $\cnr{\cc}$ to the horizontal composite 2-cell $\alpha_3 * \alpha_2 * \alpha_1$ in $\cc$.

We finish by remarking that $\cnr{\cc}$ actually classifies \textbf{all} lax functors, not just the ones that are identity-on-objects:

\begin{prop}
\textbf{Any} lax functor $F: \cc \to \cl$ admits a unique factorization through a strict 2-functor:
\[\begin{tikzcd}
	\cc && \cl \\
	\\
	{\cnr{\cc}}
	\arrow["F", squiggly, from=1-1, to=1-3]
	\arrow["P"', squiggly, from=1-1, to=3-1]
	\arrow["{\exists ! F'}"', from=3-1, to=1-3]
\end{tikzcd}\]
\end{prop}
\begin{proof}
We will prove \textbf{the existence} first. First use Proposition \ref{THMP_2Catl_OFS} to factorize $F$ as $F_\circ \circ F_\bullet$, with $F_\circ$ fully faithful 2-functor and $F_\bullet$ an identity-on-objects lax functor. By the universal property of $\cnr{\cc}$, there is a unique 2-functor making the triangle commute:
\[\begin{tikzcd}
	\cc && {\widetilde{\cc}} && \cl \\
	\\
	{\cnr{\cc}}
	\arrow["{F_\bullet}", squiggly, from=1-1, to=1-3]
	\arrow["P"', squiggly, from=1-1, to=3-1]
	\arrow["{F_\circ}", from=1-3, to=1-5]
	\arrow["{\exists ! F'_\bullet}"', from=3-1, to=1-3]
\end{tikzcd}\]
We see that the 2-functor $F_\circ \circ F'_\bullet$ gives us a solution. To prove \textbf{the uniqueness}, let $G: \cnr{\cc} \to \cl$ be a different 2-functor making the diagram in this proposition commute. Factorize $G$ as an identity-on-objects 2-functor $G_\circ$ followed by a fully faithful one -- $G_\bullet$. Again by Proposition \ref{THMP_2Catl_OFS}, there is a unique isomorphism making the following commute:
\[\begin{tikzcd}
	&&& {\widetilde{\cc}} \\
	\cc &&&&&& \cl \\
	&& {\cnr{\cc}} && {\widetilde{\cnr{\cc}}}
	\arrow["{F_\circ}", from=1-4, to=2-7]
	\arrow["{F_\bullet}", squiggly, from=2-1, to=1-4]
	\arrow["P"', squiggly, from=2-1, to=3-3]
	\arrow["{G_\bullet}"', squiggly, from=3-3, to=3-5]
	\arrow["{\exists ! \Phi}"', dashed, from=3-5, to=1-4]
	\arrow["{G_\circ}"', from=3-5, to=2-7]
\end{tikzcd}\]
The left-hand square makes $\Phi$ an identity-on-objects 2-functor, so since we have:
\[
F'_\bullet \circ P = F_\bullet = \Phi \circ G_\bullet \circ P.
\]
The universal property of $\cnr{\cc}$ now forces:
\[
F'_\bullet = \Phi \circ G_\bullet.
\]
In particular, this proves the uniqueness, since:
\[
F_\circ \circ F'_\bullet = F_\circ \circ \Phi \circ G_\bullet = G_\circ \circ G_\bullet = G.
\]
\end{proof}
\noindent Note that this description of the lax functor classifier 2-category has been sketched in \cite[Page 246]{elephant}. 
\end{pr}

\begin{pr}[Coherence for lax double functors between double categories]\label{prcodlax-dblcats}
Consider the double category 2-monad $T$ from Example \ref{PR_doublecat_2monad}. Given a strict $T$-algebra $X$, which is a double category, the double category $\cnr{X}$ has the same objects, vertical morphisms as $X$, horizontal morphisms are paths of horizontal morphisms of $X$, and each square in $\cnr{X}$ has two components: a \textit{partial evaluation square} followed by a tuple of squares in $X$, for instance:
\[\begin{tikzcd}
	{a_1} && {a_2} && {a_3} && {a_4} && {a_5} \\
	{a_1} &&& {a_3} && {a_5} &&& {a_5} \\
	{b_1} &&& {b_2} && {b_3} &&& {b_4}
	\arrow["{f_1}", from=1-1, to=1-3]
	\arrow[Rightarrow, no head, from=1-1, to=2-1]
	\arrow["{f_2}", from=1-3, to=1-5]
	\arrow["{f_3}", from=1-5, to=1-7]
	\arrow["{f_4}", from=1-7, to=1-9]
	\arrow[Rightarrow, no head, from=1-9, to=2-9]
	\arrow[""{name=0, anchor=center, inner sep=0}, "{f_2 \circ f_1}", from=2-1, to=2-4]
	\arrow["{u_1}"', from=2-1, to=3-1]
	\arrow[""{name=1, anchor=center, inner sep=0}, "{f_4 \circ f_3}"{description}, from=2-4, to=2-6]
	\arrow["{u_2}"', from=2-4, to=3-4]
	\arrow[""{name=2, anchor=center, inner sep=0}, "{1_{a_5}}", from=2-6, to=2-9]
	\arrow["{u_3}"', from=2-6, to=3-6]
	\arrow["{u_4}"', from=2-9, to=3-9]
	\arrow[""{name=3, anchor=center, inner sep=0}, "{g_1}"', from=3-1, to=3-4]
	\arrow[""{name=4, anchor=center, inner sep=0}, "{g_2}"', from=3-4, to=3-6]
	\arrow[""{name=5, anchor=center, inner sep=0}, "{g_3}"', from=3-6, to=3-9]
	\arrow["{((f_1,f_2),(f_3,f_4),())}"{pos=0.4}, shorten <=3pt, shorten >=3pt, Rightarrow, from=1-5, to=1]
	\arrow["{\alpha_1}", shorten <=4pt, shorten >=4pt, Rightarrow, from=0, to=3]
	\arrow["{\alpha_2}", shorten <=4pt, shorten >=4pt, Rightarrow, from=1, to=4]
	\arrow["{\alpha_3}", shorten <=4pt, shorten >=4pt, Rightarrow, from=2, to=5]
\end{tikzcd}\]
More succinctly, such a square might be portrayed as follows (and has been done so in the paper we cite below):
\[\begin{tikzcd}
	{a_1} && {a_2} && {a_3} & {a_4} & {a_5} \\
	{b_1} &&&& {b_2} && {b_3} &&& {b_4}
	\arrow["{f_1}"{pos=0.7}, from=1-1, to=1-3]
	\arrow["{u_1}"', from=1-1, to=2-1]
	\arrow["{f_2}", from=1-3, to=1-5]
	\arrow["{f_3}", from=1-5, to=1-6]
	\arrow["{u_2}"', from=1-5, to=2-5]
	\arrow["{f_4}", from=1-6, to=1-7]
	\arrow["{u_3}"', from=1-7, to=2-7]
	\arrow[""{name=0, anchor=center, inner sep=0}, "{u_4}", from=1-7, to=2-10]
	\arrow[""{name=1, anchor=center, inner sep=0}, "{g_1}"', from=2-1, to=2-5]
	\arrow[""{name=2, anchor=center, inner sep=0}, "{g_2}"', from=2-5, to=2-7]
	\arrow[""{name=3, anchor=center, inner sep=0}, "{g_3}"', from=2-7, to=2-10]
	\arrow["{\alpha_1}", shorten <=3pt, shorten >=3pt, Rightarrow, from=1-3, to=1]
	\arrow["{\alpha_2}", shorten <=3pt, shorten >=3pt, Rightarrow, from=1-6, to=2]
	\arrow["{\alpha_3}", shift right=5, shorten <=2pt, shorten >=2pt, Rightarrow, from=0, to=3]
\end{tikzcd}\]

The dual of the construction $\cnr{X}$ (for the \textbf{colax} double functor classifier) agrees with the construction $\mathbb{P}\text{ath } X$ of \cite[The construction 1.1, Proposition 1.19]{pathsindblcats}. In~addition:
\begin{itemize}
\item the fact that the colax morphism classifier $Q_c$ is colax-idempotent has been proven in \cite[Theorem 1.21]{pathsindblcats},
\item the fact that the colax morphism classifier admits an internal strict factorization system was proven in \cite[Proposition 1.5]{pathsindblcats} -- this is an internal version of Lemma \ref{lemma_coddiscrimpliessfs}.
\end{itemize}

\end{pr}

\begin{pr}\label{PR_reslan_laxcoh_explicit}
Let $\cj$ be a 1-category and consider the reslan 2-monad $T$ on the 2-category $[\ob \cj, \cat]$ from Example \ref{PR_reslan2monad}. We already know from Example \ref{PR_reslan_laxcoh} that coherence for colax algebras holds for this 2-monad -- in this example we will provide explicit descriptions. The 2-monad is cartesian and of form $\cat(T)$ (see Example \ref{PR_reslan_catT}).

Given a colax functor $W: \cj \to \cat$, the resolution $\res{W}$ is a diagram in $[ \ob \cj, \cat ]$ -- at each component this is a codomain-colax, domain-codiscrete category. Unpacking the data, we  obtain that for any $b \in \cj$, $\cnr{W}(b)$ is a category such that the objects are pairs $(f: j \to b, x \in Wj)$, and a morphism:
\[
(f: j \to b, x \in Wj) \to (g: j' \to b, y \in Wj'),
\]
is a pair $(h: j \to j', \alpha: W(h)(x) \to y)$ of a morphism in $\cj$ and a morphism in $Wj'$ such that:
\[\begin{tikzcd}
	j && {j'} \\
	& b
	\arrow["h", from=1-1, to=1-3]
	\arrow["f"', from=1-1, to=2-2]
	\arrow["g", from=1-3, to=2-2]
\end{tikzcd}\]
Given any 1-cell $k: b \to b'$, the functor $\cnr{W}b \to \cnr{W}b'$ is defined on objects and morphisms as follows:
\begin{align*}
(f: j \to b, x \in Wj) &\mapsto (kf : j \to b', x \in Wj),\\
(h:j \to j'',\alpha: W(h)(x) \to y) &\mapsto (h,\alpha).
\end{align*}

Theorem \ref{THM_lax_coherence_colax_algebras_EXPOSITION} in this case says that for every colax functor $W: \cj \to \cat$ and every $b \in \cj$, there is an adjunction:
\[\begin{tikzcd}
	Wb && {\cnr{W}b} & {\text{ in } \cat}
	\arrow[""{name=0, anchor=center, inner sep=0}, "{p_{W,b}}"', curve={height=24pt}, from=1-1, to=1-3]
	\arrow[""{name=1, anchor=center, inner sep=0}, "{q_{W,b}}"', curve={height=24pt}, from=1-3, to=1-1]
	\arrow["\dashv"{anchor=center, rotate=-90}, draw=none, from=1, to=0]
\end{tikzcd}\]
The right adjoint is defined as follows -- it uses $W$'s counitor:
\[\begin{tikzcd}
	x && {(1_b,x \in Wb)} \\
	& \mapsto \\
	y && {(1_b,x \in Wb)}
	\arrow["\alpha"', from=1-1, to=3-1]
	\arrow["{(1_b, \alpha \circ \iota_{b,x}: W(1_b)x \to y)}", from=1-3, to=3-3]
\end{tikzcd}\]
The left adjoint is defined like this -- it uses $W$'s coassociator:
\[\begin{tikzcd}
	{(f: j \to b, x \in Wj)} && {W(f)(x)} \\
	& \mapsto & {W(gh)(x)} \\
	&& {W(g)(W(h)(x))} \\
	{(g: j' \to b, y \in Wj')} && {W(g)(y)}
	\arrow["{(h,\alpha)}"', from=1-1, to=4-1]
	\arrow[Rightarrow, no head, from=1-3, to=2-3]
	\arrow["{\gamma_{b,x}}", from=2-3, to=3-3]
	\arrow["{W(g)(\alpha)}", from=3-3, to=4-3]
\end{tikzcd}\]

The dual of this construction -- the computation of the left 2-adjoint to the inclusion $[\cj, \cat] \to \laxhomc{\cj}{\cat}$, has been given in \cite[Theorem 3]{twoconstr} (colax natural transformations have been called \textit{right lax natural transformations}). Also, \cite[Theorem 2]{twoconstr} describes a canonical adjunction between a lax functor and its strictification.

Finally, let us keep the promise we gave in Remark \ref{POZN_2limit_is_not_laxlimit} and give an example of a weighted limit that is not a lax limit. Let $J := (r: a \to b)$ be the free walking arrow category. Using the above description, we see that when $W: J \to \cat$ is such that $Wa = \emptyset$, $\cnr{W}(a) = \emptyset$. If $Wa \not = \emptyset$, pick an object $x_0 \in Wa$. Then $\cnr{W}b$ is not discrete, since it contains the following non-identity morphism:
\[\begin{tikzcd}
	{(r,x_0)} && {(1_b,Wr(x_0)).}
	\arrow["r", from=1-1, to=1-3]
\end{tikzcd}\]
We see that $\widetilde{W} := *$, having $\widetilde{W}a$ nonempty and $\widetilde{W}b$ discrete, is not of form $\cnr{W}$ for any $W: \cj \to \cat$.
\end{pr}

\begin{pozn}[Colax morphism classifiers]\label{POZN_colax_mor_classifier}
We can apply dualities to compute \textbf{colax} morphism classifier for a 2-monad $T$ on $\cat$ of form $\cat(T')$ as follows. First note that the opposite category 2-functor $(-)^{op}: \cat^{co} \to \cat$ induces a 2-isomorphism
\begin{align*}
\talgc &\cong T^{co}\text{-Alg}_l,\\
(A,a) &\mapsto (A^{op},a^{op}).
\end{align*}
This implies that a $T$-algebra $(B,b)$ is the colax $T$-morphism classifier for $(A,a)$ if and only if $(B^{op},b^{op})$ is the lax $T^{co}$-morphism classifier for $(A^{op},a^{op})$.

Now let $(A,a)$ be a strict $T$-algebra. The lax-$T^{co}$-morphism classifier for $(A^{op},a^{op})$ is a $T^{co}$-algebra $\cnr{A^{op}}$, and thus the colax morphism classifier is given by the formula:
\[
(A,a)' = \cnr{A^{op}}^{op}.
\]

For instance, the colax monoidal functor classifier for a strict monoidal category $(\ca,\otimes)$ again has tuples of objects in $\ca$ as object, and a morphism $(a_1,a_2,a_3,a_4) \to$\\ $\to (b_1,b_2,b_3,b_4)$ is a corner (or rather, a \textit{cospan}) like this:
\[\begin{tikzcd}
	{(a_1,a_2,a_3,a_4)} \\
	{(b_1, b_2\otimes b_3, I, b_4)} &&& {(b_1,b_2,b_3,b_4)}
	\arrow["{(f_1,f_2,f_3,f_4)}"', from=1-1, to=2-1]
	\arrow["{((b_1),(b_2,b_3),(),(b_4))}", from=2-4, to=2-1]
\end{tikzcd}\]
\end{pozn}

%% file: PhD_CH_5_COLAX_ADJ_LAX_IDEMP_PSMON.tex
\chapter[Colax adjunctions and KZ pseudomonads]{Colax adjunctions and lax-idempotent pseudomonads}\label{CH_colax_adjs_lax_idemp_psmons}
\section{Introduction}\label{SEC_intro_colax_adj_lax_idemp_psm}

The primary motivation for this chapter is to develop lax analogues of results in two-dimensional monad theory, as proven in the papers \cite{twodim} and \cite{onsemiflexible}. Under mild assumptions, various results have been proven in \cite{twodim} about $T$-algebras and pseudo morphisms, for instance their bicocompleteness or the existence of change-of-base biadjunctions between 2-categories of algebras and pseudo morphisms for two different 2-monads $S,T$.

On the other hand, there are fewer known results about 2-categories $\talgl$ of algebras and lax morphisms between them. Limits in $\talgl$ have been well-understood (\cite{limitsforlax}, \cite{enhanced2cats}), but not much has been proven about colimits. This was for a good reason: 2-colimits or even bicolimits often do not exist in those 2-categories. Our task in this chapter is to consider a weakened version of a bicolimit (rali-colimits of Section \ref{SEC_weak_limits}) and show that, under mild assumptions, 2-categories of lax morphisms are in fact cocomplete in this weak sense. Another task we have is to establish change-of-base theorems for algebras and lax morphisms. Again, the notion that works for pseudo morphisms -- biadjunctions -- will have to be replaced by a weaker one -- colax adjunctions.

The 2-category $\talg$ of algebras and pseudo morphisms can often be described as the Kleisli 2-category for a certain pseudo-idempotent 2-comonad (Remark \ref{POZN_Talg_is_kleisli}). A key observation to be made is that many statements and proofs about $\talg$ in papers \cite{twodim}, \cite{onsemiflexible} are very formal and are in fact true for any pseudo-idempotent 2-comonad on a 2-category. They also easily dualize to pseudo-idempotent 2-monads. Since we are interested in the lax world, we are naturally led to the study of Kleisli 2-categories for lax-idempotent pseudomonads, using the formalism of left Kan pseudomonads (recalled in Section \ref{SECS_left_kan_psmons}). The usage of pseudomonads instead of 2-monads will allow us to consider a wider array of examples such as the small presheaf pseudomonad, and lets us prove that the bicategory $\text{Prof}$ of locally small categories and small profunctors is weakly complete in the sense of the previous paragraph.

As mentioned, colax adjunctions are inevitable when working with lax morphisms. The definition of a (co)lax adjunction is hard to work with because it contains a large amount of data. The first result in this chapter, Theorem \ref{THM_Bunge}, shows that a left colax adjoint $F$ to a pseudofunctor $U$ can be more conveniently given by a collection of 1-cells $y_A: A \to UFA$ satisfying certain ``relative $U$-left Kan extension” conditions. This is an extension of the work of Bunge and Gray (\cite{bunge}, \cite{formalcat}) where this has been proven for the case when $U$ is a 2-functor. A result of this kind is similar to how left Kan pseudomonads provide a more convenient description of lax-idempotent pseudomonads. We will use this theorem to obtain results on colax adjunctions involving the Kleisli 2-category for a lax-idempotent pseudomonad (Theorem \ref{THM_BIG_lax_adj_thm}), and the dual of \textit{this} result will be used to obtain results on colax adjunctions involving $T$-algebras and lax $T$-morphisms for a 2-monad (Theorem \ref{THM_big_colax_adj_THM}).

\medskip

\noindent \textbf{The chapter is organized as follows}:

\begin{itemize}
\item In Section \ref{SEKCE_Bunge} we prove the generalization of Bunge's and Gray's results on colax adjunctions to the setting of pseudofunctors: we show that there is a correspondence between left colax adjoint pseudofunctors to a pseudofunctor $U$ and collections of 1-cells $y_A: A \to UFA$ satisfying the aforementioned relative $U$-left Kan extension conditions (Theorem \ref{THM_korespondence_lax_adj_u_coh_ext}).

\item In Section \ref{SEKCE_Kleisli2cat} we first give (an essentially folklore) characterization of algebras for a lax-idempotent pseudomonad in terms of the existence of certain adjoints (Proposition \ref{THM_characterization_D-algs_D-corefl-incl}). We then use this characterization and the generalized Bunge's and Gray's result to prove that when given a lax-idempotent pseudomonad $D$ on $\ck$, any left biadjoint $\ck \to \cl$ that factorizes through the Kleisli 2-category $\ck_D$ gives rise to a colax left adjoint $\ck_D \to \cl$ (Theorem \ref{THM_BIG_lax_adj_thm}). We list various applications, for instance the weak completeness of $\ck_D$ (Theorem \ref{THM_Kleisli2cat_is_corefl_complete}) provided that $\ck$ is bicomplete, or that there is a canonical colax adjunction between $\ck_D$ and the 2-category of pseudo $D$-algebras (Corollary \ref{THM_lax_adjunkce_EM_Kleisli}).

\item In Section \ref{SEKCE_applications_to_twodim} we spell out what these results in particular say about the 2-category $\talgl$ of strict algebras and lax morphisms for a 2-monad $T$. This includes the aforementioned colax change-of-base theorem (Corollary \ref{THM_COR_change_of_base_SalgTalg}) as well as the weak cocompleteness result for $\talgl$ (Theorem \ref{THM_talgl_reflector_cocomplete}).
\end{itemize}

\section{Relative Kan extensions and colax adjunctions}\label{SEKCE_Bunge}

In \cite{bunge}, Bunge introduced the notion of a relative Kan extension with respect to a 2-functor $U$ and showed that for a collection $y_A: A \to U\overline{A}$ of 1-cells that admit these extensions (and satisfy certain coherence conditions), there is an induced left colax adjoint $F$ to $U$, where $F$ is a colax functor (\cite[Theorem 4.1]{bunge}). She also proves a partial converse to this result (\cite[Theorem 4.3]{bunge}). Note that at the same time these results also appeared in Gray's work (\cite[I,7.8.]{formalcat}).

In this section, we generalize these results to the case where $U$ is a pseudofunctor and, on the other hand, refine it by identifying conditions under which the colax left adjoint $F$ is actually a pseudofunctor. This enables us to describe, in Theorem \ref{THM_korespondence_lax_adj_u_coh_ext}, a symmetric relationship between $U$-extensions and colax adjunctions. We will see an application of these results to the settings of lax-idempotent pseudomonads in Section \ref{SEKCE_Kleisli2cat}.

\begin{defi}
Let $U: \cd \to \cc$ be a pseudofunctor and let\footnote{In what follows, the codomains of 1-cells $y_A$ will be of form $UB$, where $B \in \cd$. Abusing notation, we will denote $B := FA$, since this will be the action of soon-to-be defined colax functor $F$ on the object $A$.} $y_A: A \to UFA$, $f: A \to UB$ be 1-cells of $\cc$. The \textit{left $U$-extension} of $f$ along $y_A$ is a pair $(f',\psi_f)$ of a 1-cell in $\cd$ and a 2-cell in $\cc$ (as pictured below) with the property that for any other pair $(g,\alpha)$, there is a unique 2-cell $\theta: f' \Rightarrow g$ in $\cd$ such that the following 2-cells are equal:
\[\begin{tikzcd}
	A && UFA &&& A && UFA \\
	&&&& {=} \\
	&& UB &&&&& UB
	\arrow["{y_A}", from=1-1, to=1-3]
	\arrow[""{name=0, anchor=center, inner sep=0}, "f"', from=1-1, to=3-3]
	\arrow[""{name=1, anchor=center, inner sep=0}, "{Uf'}"{description, pos=0.7}, from=1-3, to=3-3]
	\arrow[""{name=2, anchor=center, inner sep=0}, "Ug", curve={height=-24pt}, from=1-3, to=3-3]
	\arrow["{y_A}", from=1-6, to=1-8]
	\arrow[""{name=3, anchor=center, inner sep=0}, "f"', from=1-6, to=3-8]
	\arrow["Ug", from=1-8, to=3-8]
	\arrow["{\psi_f}"', shorten <=5pt, Rightarrow, from=0, to=1-3]
	\arrow["\alpha", shorten <=5pt, Rightarrow, from=3, to=1-8]
	\arrow["U\theta", shorten <=5pt, shorten >=5pt, Rightarrow, from=1, to=2]
\end{tikzcd}\]
\end{defi}

\begin{defi}\label{DEFI_coherently_closed_for_Uext}
Let $U: \cd \to \cc$ be a pseudofunctor. We say that a collection of 1-cells $y_A: A \to UFA$ for each object $A \in \cc$ is \textit{coherently closed for $U$-extensions} if:
\begin{itemize}
\item for every $f: A \to UB$ we have a choice of a $U$-extension $(f^\mbbd, \mbbd_f)$,
\item the following composite 2-cell exhibits $1_{UY}^\mbbd \circ (y_{UY} f)^\mbbd$ as the left $U$-extension of $f$ along $y_X$:
\begin{equation}\label{EQ_DEFI_coherently_closed_for_Uext}
\begin{tikzcd}
	X &&& UFX \\
	\\
	UY &&& UFUY \\
	\\
	&&& UY
	\arrow[""{name=0, anchor=center, inner sep=0}, "{y_X}", from=1-1, to=1-4]
	\arrow[""{name=0p, anchor=center, inner sep=0}, phantom, from=1-1, to=1-4, start anchor=center, end anchor=center]
	\arrow["f"', from=1-1, to=3-1]
	\arrow[""{name=1, anchor=center, inner sep=0}, "{y_{UY}}"{description}, from=3-1, to=3-4]
	\arrow[""{name=1p, anchor=center, inner sep=0}, phantom, from=3-1, to=3-4, start anchor=center, end anchor=center]
	\arrow[""{name=1p, anchor=center, inner sep=0}, phantom, from=3-1, to=3-4, start anchor=center, end anchor=center]
	\arrow["{U(y_{UY} f)^\mbbd}"{description}, from=1-4, to=3-4]
	\arrow[""{name=2, anchor=center, inner sep=0}, Rightarrow, no head, from=3-1, to=5-4]
	\arrow[""{name=2p, anchor=center, inner sep=0}, phantom, from=3-1, to=5-4, start anchor=center, end anchor=center]
	\arrow["{U1_{UY}^\mbbd}"{description}, from=3-4, to=5-4]
	\arrow[""{name=3, anchor=center, inner sep=0}, "{U(1_{UY}^\mbbd \circ (y_{UY} f)^\mbbd)}", shift left=5, curve={height=-30pt}, from=1-4, to=5-4]
	\arrow[""{name=3p, anchor=center, inner sep=0}, phantom, from=1-4, to=5-4, start anchor=center, end anchor=center, shift left=5, curve={height=-30pt}]
	\arrow[""{name=4, anchor=center, inner sep=0}, curve={height=-6pt}, draw=none, from=1-4, to=5-4]
	\arrow["{\mbbd_{y_{UY} f}}"', shorten <=9pt, shorten >=9pt, Rightarrow, from=1p, to=0p]
	\arrow["{\mbbd_{1_{UY}}}"', shorten <=4pt, shorten >=4pt, Rightarrow, from=2p, to=1p]
	\arrow["{\gamma^{-1}}"{pos=0.7}, shorten <=14pt, Rightarrow, from=4, to=3p]
\end{tikzcd}
\end{equation}
Equivalently, this is to say that the canonical 2-cell $f^\mbbd \Rightarrow 1^{\mbbd}_Y \circ (y_Y f)^\mbbd$ is invertible.
\end{itemize}
\end{defi}

\noindent In what follows, we will sometimes omit the lower index for the 2-cell components $\mathbb{D}_f$ of the left $U$-extensions.

\begin{theo}\label{THM_Bunge}
Let $U: \cd \to \cc$ be a pseudofunctor and $y_A : A \to UFA$ a collection of 1-cells coherently closed for $U$-extensions. Then:
\begin{itemize}
\item the mapping $A \mapsto FA$ can be extended to a colax functor $F: \cc \to \cd$,
\item $y$ can be extended to a colax natural transformation $1_\cc \Rightarrow UF$,
\item there exists a colax natural transformation $\epsilon: FU \Rightarrow 1_\cd$ and a modification\\ $\Phi: 1_U \to U\epsilon \circ yU$.

\medskip

\noindent Assume moreover the \textit{composition} and \textit{unit} axioms for $U$-extensions: the diagram below left is a $U$-extension of $y_A$ along $y_A$, and the diagram below right is the $U$-extension of $y_C gf$ along $y_A$:
\begin{equation}\label{EQ_U_extensions_unit_comp}
\begin{tikzcd}
	A && UFA && A && UFA \\
	&&&& B && UFB \\
	&& UFA && C && UFC
	\arrow[""{name=0, anchor=center, inner sep=0}, "{y_A}", from=1-5, to=1-7]
	\arrow["U(y_Bf)^\mbbd"', from=1-7, to=2-7]
	\arrow["f"', from=1-5, to=2-5]
	\arrow["g"', from=2-5, to=3-5]
	\arrow[""{name=1, anchor=center, inner sep=0}, "{y_C}"', from=3-5, to=3-7]
	\arrow["U(y_Cg)^\mbbd"', from=2-7, to=3-7]
	\arrow["{y_A}", from=1-1, to=1-3]
	\arrow["{y_A}"', shift right=3, from=1-1, to=3-3]
	\arrow[""{name=2, anchor=center, inner sep=0}, "{U1_{FA}}", curve={height=-24pt}, from=1-3, to=3-3]
	\arrow[""{name=3, anchor=center, inner sep=0}, curve={height=24pt}, Rightarrow, no head, from=1-3, to=3-3]
	\arrow[""{name=4, anchor=center, inner sep=0}, "{U((y_Cg)^\mbbd (y_Bf)^\mbbd)}", shift left=5, curve={height=-30pt}, from=1-7, to=3-7]
	\arrow[""{name=4p, anchor=center, inner sep=0}, phantom, from=1-7, to=3-7, start anchor=center, end anchor=center, shift left=5, curve={height=-30pt}]
	\arrow[""{name=5, anchor=center, inner sep=0}, "{y_B}"{description}, from=2-5, to=2-7]
	\arrow["{\iota^{-1}}", shorten <=10pt, shorten >=10pt, Rightarrow, from=3, to=2]
	\arrow["\mbbd", shorten <=4pt, shorten >=4pt, Rightarrow, from=5, to=0]
	\arrow["\mbbd", shorten <=4pt, shorten >=4pt, Rightarrow, from=1, to=5]
	\arrow["{\gamma^{-1}}", shorten >=6pt, Rightarrow, from=2-7, to=4p]
\end{tikzcd}
\end{equation}

\noindent Then:
\item $F$ is a pseudofunctor,
\item there is an invertible modification $\Psi: \epsilon F \circ Fy \to 1_F$ and all this data gives a colax adjunction:
\[
(\Psi,\Phi): (\epsilon, y): F \ddashv U: \cd \to \cc.
\]
\end{itemize}
\end{theo}
\begin{proof}
Denote by $(f^\mbbd,\mbbd)$ the choice of a $U$-extension of $f: A \to UB$. Define the colax functor $F: \cc \to \cd$ on a morphism $f: A \to B$ as:
\[
Ff := (y_B \circ f)^\mbbd.
\]
We then have the following $U$-extension:
\[\begin{tikzcd}
	A && UFA \\
	\\
	B && UFB
	\arrow[""{name=0, anchor=center, inner sep=0}, "{y_A}", from=1-1, to=1-3]
	\arrow["UFf", from=1-3, to=3-3]
	\arrow[""{name=1, anchor=center, inner sep=0}, "{y_B}"', from=3-1, to=3-3]
	\arrow["f"', from=1-1, to=3-1]
	\arrow["\mbbd"', shorten <=9pt, shorten >=9pt, Rightarrow, from=1, to=0]
\end{tikzcd}\]
Define the action of $F$ on a 2-cell $\alpha$ as the unique 2-cell making the following equal:
\begin{equation}\label{EQP_Falpha}
\begin{tikzcd}
	A && UFA &&&& A && UFA \\
	&&&& {=} \\
	B && UFB &&&& B && UFB
	\arrow[""{name=0, anchor=center, inner sep=0}, "{y_A}", from=1-1, to=1-3]
	\arrow[""{name=0p, anchor=center, inner sep=0}, phantom, from=1-1, to=1-3, start anchor=center, end anchor=center]
	\arrow[""{name=1, anchor=center, inner sep=0}, "UFf"', from=1-3, to=3-3]
	\arrow[""{name=1p, anchor=center, inner sep=0}, phantom, from=1-3, to=3-3, start anchor=center, end anchor=center]
	\arrow[""{name=2, anchor=center, inner sep=0}, "{y_B}"', from=3-1, to=3-3]
	\arrow[""{name=2p, anchor=center, inner sep=0}, phantom, from=3-1, to=3-3, start anchor=center, end anchor=center]
	\arrow["f"', from=1-1, to=3-1]
	\arrow[""{name=3, anchor=center, inner sep=0}, "UFg", curve={height=-30pt}, from=1-3, to=3-3]
	\arrow[""{name=3p, anchor=center, inner sep=0}, phantom, from=1-3, to=3-3, start anchor=center, end anchor=center, curve={height=-30pt}]
	\arrow[""{name=4, anchor=center, inner sep=0}, "f"', curve={height=30pt}, from=1-7, to=3-7]
	\arrow[""{name=4p, anchor=center, inner sep=0}, phantom, from=1-7, to=3-7, start anchor=center, end anchor=center, curve={height=30pt}]
	\arrow[""{name=5, anchor=center, inner sep=0}, "g", from=1-7, to=3-7]
	\arrow[""{name=5p, anchor=center, inner sep=0}, phantom, from=1-7, to=3-7, start anchor=center, end anchor=center]
	\arrow[""{name=6, anchor=center, inner sep=0}, "{y_A}", from=1-7, to=1-9]
	\arrow[""{name=6p, anchor=center, inner sep=0}, phantom, from=1-7, to=1-9, start anchor=center, end anchor=center]
	\arrow[""{name=7, anchor=center, inner sep=0}, "{y_B}"', from=3-7, to=3-9]
	\arrow[""{name=7p, anchor=center, inner sep=0}, phantom, from=3-7, to=3-9, start anchor=center, end anchor=center]
	\arrow["UFg", from=1-9, to=3-9]
	\arrow["\mbbd", shorten <=9pt, shorten >=9pt, Rightarrow, from=2p, to=0p]
	\arrow["{U(\exists !)}", shorten <=3pt, shorten >=3pt, Rightarrow, from=1p, to=3p]
	\arrow["\alpha"{pos=0.4}, shorten <=3pt, shorten >=3pt, Rightarrow, from=4p, to=5p]
	\arrow["\mbbd"', shorten <=9pt, shorten >=9pt, Rightarrow, from=7p, to=6p]
\end{tikzcd}
\end{equation}

The above equation makes $y$ locally natural. The coassociator and the counitor 2-cells $\gamma': F(gf) \Rightarrow Fg \circ Ff$, $\iota': F1_A \Rightarrow 1_{FA}$ for $F$ are given as the unique 2-cells satisfying these equations:
\begin{equation}\label{EQP_gamma_prime}
\begin{tikzcd}
	A && UFA &&& A && UFA \\
	&&&& {=} & B && UFB \\
	C && UFC &&& C && UFC
	\arrow[""{name=0, anchor=center, inner sep=0}, "{y_A}", from=1-6, to=1-8]
	\arrow[""{name=0p, anchor=center, inner sep=0}, phantom, from=1-6, to=1-8, start anchor=center, end anchor=center]
	\arrow["UFf"', from=1-8, to=2-8]
	\arrow[""{name=1, anchor=center, inner sep=0}, "{y_B}"{description}, from=2-6, to=2-8]
	\arrow[""{name=1p, anchor=center, inner sep=0}, phantom, from=2-6, to=2-8, start anchor=center, end anchor=center]
	\arrow[""{name=1p, anchor=center, inner sep=0}, phantom, from=2-6, to=2-8, start anchor=center, end anchor=center]
	\arrow["f"', from=1-6, to=2-6]
	\arrow["g"', from=2-6, to=3-6]
	\arrow[""{name=2, anchor=center, inner sep=0}, "{y_C}"', from=3-6, to=3-8]
	\arrow[""{name=2p, anchor=center, inner sep=0}, phantom, from=3-6, to=3-8, start anchor=center, end anchor=center]
	\arrow["UFg"', from=2-8, to=3-8]
	\arrow[""{name=3, anchor=center, inner sep=0}, "{U(Fg \circ Ff)}", shift left=3, curve={height=-30pt}, from=1-8, to=3-8]
	\arrow[""{name=3p, anchor=center, inner sep=0}, phantom, from=1-8, to=3-8, start anchor=center, end anchor=center, shift left=3, curve={height=-30pt}]
	\arrow["gf"', from=1-1, to=3-1]
	\arrow[""{name=4, anchor=center, inner sep=0}, "{y_C}"', from=3-1, to=3-3]
	\arrow[""{name=4p, anchor=center, inner sep=0}, phantom, from=3-1, to=3-3, start anchor=center, end anchor=center]
	\arrow[""{name=5, anchor=center, inner sep=0}, "{y_A}", from=1-1, to=1-3]
	\arrow[""{name=5p, anchor=center, inner sep=0}, phantom, from=1-1, to=1-3, start anchor=center, end anchor=center]
	\arrow[""{name=6, anchor=center, inner sep=0}, "{UF(gf)}"', from=1-3, to=3-3]
	\arrow[""{name=6p, anchor=center, inner sep=0}, phantom, from=1-3, to=3-3, start anchor=center, end anchor=center]
	\arrow[""{name=7, anchor=center, inner sep=0}, "{U(Fg \circ Ff)}"{pos=0.3}, curve={height=-30pt}, from=1-3, to=3-3]
	\arrow[""{name=7p, anchor=center, inner sep=0}, phantom, from=1-3, to=3-3, start anchor=center, end anchor=center, curve={height=-30pt}]
	\arrow[""{name=8, anchor=center, inner sep=0}, shift left=5, draw=none, from=1-8, to=3-8]
	\arrow[""{name=8p, anchor=center, inner sep=0}, phantom, from=1-8, to=3-8, start anchor=center, end anchor=center, shift left=5]
	\arrow["\mbbd", shorten <=4pt, shorten >=4pt, Rightarrow, from=1p, to=0p]
	\arrow["\mbbd", shorten <=9pt, shorten >=9pt, Rightarrow, from=4p, to=5p]
	\arrow["{U(\exists !)}", shorten <=6pt, shorten >=6pt, Rightarrow, from=6p, to=7p]
	\arrow["\mbbd", shorten <=4pt, shorten >=4pt, Rightarrow, from=2p, to=1p]
	\arrow["{\gamma^{-1}}", shorten <=3pt, Rightarrow, from=8p, to=3p]
\end{tikzcd}
\end{equation}

\begin{equation}\label{EQP_iota_prime}
\begin{tikzcd}
	A && UFA &&& A && UFA \\
	&&&& {=} \\
	A && UFA &&& A && UFA
	\arrow["{y_A}", from=1-6, to=1-8]
	\arrow[""{name=0, anchor=center, inner sep=0}, Rightarrow, no head, from=1-8, to=3-8]
	\arrow[""{name=0p, anchor=center, inner sep=0}, phantom, from=1-8, to=3-8, start anchor=center, end anchor=center]
	\arrow[Rightarrow, no head, from=1-6, to=3-6]
	\arrow["{y_A}"', from=3-6, to=3-8]
	\arrow[""{name=1, anchor=center, inner sep=0}, "{U1_{FA}}", curve={height=-30pt}, from=1-8, to=3-8]
	\arrow[""{name=1p, anchor=center, inner sep=0}, phantom, from=1-8, to=3-8, start anchor=center, end anchor=center, curve={height=-30pt}]
	\arrow[""{name=2, anchor=center, inner sep=0}, "{UF1_A}"', from=1-3, to=3-3]
	\arrow[""{name=2p, anchor=center, inner sep=0}, phantom, from=1-3, to=3-3, start anchor=center, end anchor=center]
	\arrow[""{name=3, anchor=center, inner sep=0}, "{U1_{FA}}", curve={height=-30pt}, from=1-3, to=3-3]
	\arrow[""{name=3p, anchor=center, inner sep=0}, phantom, from=1-3, to=3-3, start anchor=center, end anchor=center, curve={height=-30pt}]
	\arrow[""{name=4, anchor=center, inner sep=0}, "{y_A}", from=1-1, to=1-3]
	\arrow[Rightarrow, no head, from=1-1, to=3-1]
	\arrow[""{name=5, anchor=center, inner sep=0}, "{y_A}"', from=3-1, to=3-3]
	\arrow["{U(\exists !)}", shorten <=6pt, shorten >=6pt, Rightarrow, from=2p, to=3p]
	\arrow["\mbbd", shorten <=9pt, shorten >=9pt, Rightarrow, from=5, to=4]
	\arrow["{\iota^{-1}}", shorten <=6pt, shorten >=6pt, Rightarrow, from=0p, to=1p]
\end{tikzcd}
\end{equation}
The colax functor axioms for $F$ follow from those of $U$ and can be readily proven using the universal property of $U$-extensions. The above equations also make $y$ into a colax natural transformation $y: 1_\cc \Rightarrow UF$. Next, define $\epsilon_B: FUB \to B$ and $\Phi_B$ as the $U$-extension of the identity on $UB$ along $y_{UB}$:
\[\begin{tikzcd}
	UB && UFUB \\
	\\
	&& UB
	\arrow[""{name=0, anchor=center, inner sep=0}, "{y_{UB}}", from=1-1, to=1-3]
	\arrow[""{name=0p, anchor=center, inner sep=0}, phantom, from=1-1, to=1-3, start anchor=center, end anchor=center]
	\arrow["{U\epsilon_B}", from=1-3, to=3-3]
	\arrow[""{name=1, anchor=center, inner sep=0}, Rightarrow, no head, from=1-1, to=3-3]
	\arrow[""{name=1p, anchor=center, inner sep=0}, phantom, from=1-1, to=3-3, start anchor=center, end anchor=center]
	\arrow["{\Phi_B}"', shorten <=4pt, shorten >=4pt, Rightarrow, from=1p, to=0p]
\end{tikzcd}\]
The colax naturality square for $\epsilon$ at a 1-cell $h:B \to C$ is the unique 2-cell $\epsilon_h$ making the 2-cells below equal (it is guaranteed to uniquely exist because part of the diagram below left is a $U$-extensions \eqref{EQ_DEFI_coherently_closed_for_Uext}):
\begin{equation}\label{EQP_Phi_epsilon}
\begin{minipage}{0.8\textwidth}
\adjustbox{scale=0.9,center}{
\begin{tikzcd}
	UB && UFUB && UB && UFUB \\
	\\
	UC & UFUC &&& UC && UB \\
	\\
	&& UC &&&& UC
	\arrow[""{name=0, anchor=center, inner sep=0}, "{y_{UC}}", from=3-1, to=3-2]
	\arrow["{U\epsilon_C}"{description}, from=3-2, to=5-3]
	\arrow[""{name=1, anchor=center, inner sep=0}, curve={height=24pt}, Rightarrow, no head, from=3-1, to=5-3]
	\arrow[""{name=1p, anchor=center, inner sep=0}, phantom, from=3-1, to=5-3, start anchor=center, end anchor=center, curve={height=24pt}]
	\arrow["Uh"', from=1-1, to=3-1]
	\arrow[""{name=2, anchor=center, inner sep=0}, "{y_{UB}}", from=1-1, to=1-3]
	\arrow["UFUh"{description}, from=1-3, to=3-2]
	\arrow["{y_{UB}}", from=1-5, to=1-7]
	\arrow["{U\epsilon_B}"{description}, from=1-7, to=3-7]
	\arrow["Uh"', from=1-5, to=3-5]
	\arrow[curve={height=24pt}, Rightarrow, no head, from=3-5, to=5-7]
	\arrow["Uh"', from=3-7, to=5-7]
	\arrow[""{name=3, anchor=center, inner sep=0}, curve={height=24pt}, Rightarrow, no head, from=1-5, to=3-7]
	\arrow[""{name=4, anchor=center, inner sep=0}, "{U(h \circ \epsilon_B)}", curve={height=-30pt}, from=1-3, to=5-3]
	\arrow[""{name=4p, anchor=center, inner sep=0}, phantom, from=1-3, to=5-3, start anchor=center, end anchor=center, curve={height=-30pt}]
	\arrow[""{name=5, anchor=center, inner sep=0}, "{U(\epsilon_C \circ FUh)}"{description, pos=0.7}, curve={height=18pt}, from=1-3, to=5-3]
	\arrow[""{name=5p, anchor=center, inner sep=0}, phantom, from=1-3, to=5-3, start anchor=center, end anchor=center, curve={height=18pt}]
	\arrow[""{name=5p, anchor=center, inner sep=0}, phantom, from=1-3, to=5-3, start anchor=center, end anchor=center, curve={height=18pt}]
	\arrow[""{name=6, anchor=center, inner sep=0}, "{U(h \circ \epsilon_B)}", shift left=5, curve={height=-30pt}, from=1-7, to=5-7]
	\arrow[""{name=6p, anchor=center, inner sep=0}, phantom, from=1-7, to=5-7, start anchor=center, end anchor=center, shift left=5, curve={height=-30pt}]
	\arrow["\mbbd", shorten <=14pt, shorten >=14pt, Rightarrow, from=0, to=2]
	\arrow["{\Phi_B}", shorten <=9pt, Rightarrow, from=3, to=1-7]
	\arrow["{\Phi_C}", shorten <=6pt, shorten >=6pt, Rightarrow, from=1p, to=3-2]
	\arrow["{U\epsilon_h}", shorten <=10pt, shorten >=10pt, Rightarrow, from=5p, to=4p]
	\arrow["{\gamma^{-1}}", shorten >=2pt, Rightarrow, from=3-2, to=5p]
	\arrow["{\gamma^{-1}}", shorten >=6pt, Rightarrow, from=3-7, to=6p]
\end{tikzcd}
}
\end{minipage}
\end{equation}
This also makes $\Phi$ into a modification $1_U \to U\epsilon \circ yU$.
Let us now consider the additional assumptions. It is clear that $F$ will be a pseudofunctor. Define $\Psi_A: \epsilon_{FA} \circ Fy_A \Rightarrow 1_{FA}$ as the unique 2-cell making the two 2-cells below equal:
\begin{equation}\label{EQP_Psi}
\begin{minipage}{0.8\textwidth}
\adjustbox{scale=0.9,center}{
\begin{tikzcd}
	A &&& UFA && A && UFA \\
	\\
	UFA && UFUFA &&& UFA \\
	\\
	&&& UFA &&&& UFA
	\arrow[""{name=0, anchor=center, inner sep=0}, "{y_{UFA}}", from=3-1, to=3-3]
	\arrow[""{name=0p, anchor=center, inner sep=0}, phantom, from=3-1, to=3-3, start anchor=center, end anchor=center]
	\arrow["{U\epsilon_{FA}}"{description}, from=3-3, to=5-4]
	\arrow[""{name=1, anchor=center, inner sep=0}, curve={height=24pt}, Rightarrow, no head, from=3-1, to=5-4]
	\arrow[""{name=1p, anchor=center, inner sep=0}, phantom, from=3-1, to=5-4, start anchor=center, end anchor=center, curve={height=24pt}]
	\arrow["{y_A}"', from=1-1, to=3-1]
	\arrow[""{name=2, anchor=center, inner sep=0}, "{y_A}", from=1-1, to=1-4]
	\arrow[""{name=2p, anchor=center, inner sep=0}, phantom, from=1-1, to=1-4, start anchor=center, end anchor=center]
	\arrow["{UFy_A}"{description}, from=1-4, to=3-3]
	\arrow[""{name=3, anchor=center, inner sep=0}, "{U1_{FA}}", curve={height=-30pt}, from=1-4, to=5-4]
	\arrow[""{name=3p, anchor=center, inner sep=0}, phantom, from=1-4, to=5-4, start anchor=center, end anchor=center, curve={height=-30pt}]
	\arrow[""{name=4, anchor=center, inner sep=0}, "{U1_{FA}}", curve={height=-30pt}, from=1-8, to=5-8]
	\arrow[""{name=4p, anchor=center, inner sep=0}, phantom, from=1-8, to=5-8, start anchor=center, end anchor=center, curve={height=-30pt}]
	\arrow[""{name=5, anchor=center, inner sep=0}, Rightarrow, no head, from=1-8, to=5-8]
	\arrow[""{name=5p, anchor=center, inner sep=0}, phantom, from=1-8, to=5-8, start anchor=center, end anchor=center]
	\arrow["{y_A}"', from=1-6, to=3-6]
	\arrow[curve={height=24pt}, Rightarrow, no head, from=3-6, to=5-8]
	\arrow[""{name=6, anchor=center, inner sep=0}, "{U(\epsilon_{FA} \circ Fy_A)}"{description, pos=0.6}, from=1-4, to=5-4]
	\arrow[""{name=6p, anchor=center, inner sep=0}, phantom, from=1-4, to=5-4, start anchor=center, end anchor=center]
	\arrow[""{name=6p, anchor=center, inner sep=0}, phantom, from=1-4, to=5-4, start anchor=center, end anchor=center]
	\arrow["{y_A}", from=1-6, to=1-8]
	\arrow["\mbbd", shorten <=14pt, shorten >=14pt, Rightarrow, from=0p, to=2p]
	\arrow["{\iota^{-1}}"{pos=0.4}, shorten <=6pt, shorten >=6pt, Rightarrow, from=5p, to=4p]
	\arrow["{U(\exists !)}", shorten <=6pt, shorten >=6pt, Rightarrow, from=6p, to=3p]
	\arrow["{\gamma^{-1}}", shorten >=5pt, Rightarrow, from=3-3, to=6p]
	\arrow["{\Phi_{FA}}", shorten <=7pt, shorten >=7pt, Rightarrow, from=1p, to=3-3]
\end{tikzcd}
}
\end{minipage}
\end{equation}
By the assumption, $\Psi_A$ is invertible. This equality also proves the first swallowtail identity. What remains to be proved is the following:
\begin{itemize}
\item $\epsilon$ is colax natural,
\item $\Phi$ is a modification,
\item the second swallowtail identity.
\end{itemize}
These are all straightforward computations and we will prove them in the Appendix as Lemma \ref{THML_epsilon_phi_second_swallowtail}.
\end{proof}

\begin{pozn}[A special case]\label{POZN_almost_a_biadj}
Observe from the above proof that if each $U$-extension $(f^\mbbd, \mbbd_f)$ has the 2-cell $\mbbd_f$ invertible, in the colax adjunction obtained:
\[
(\Psi,\Phi): (\epsilon, y): F \ddashv U: \cd \to \cc,
\]
the modification $\Phi$ is invertible (in addition to $\Psi$ being invertible) and $y$ is pseudonatural. It is still not a biadjunction since $\epsilon$ has no reason to be pseudonatural.
\end{pozn}

\noindent We may use the hard work from the previous proof to give a full proof of the folklore result characterizing biadjunctions:

\begin{cor}\label{THM_folklore_biadjunctions}
The following are equivalent for a pseudofunctor $U: \cd \to \cc$:
\begin{enumerate}
\item there exists a biadjunction $(\Psi, \Phi): (\epsilon, \eta): F \dashv U$,
\item there exists a pseudofunctor $F: \cc \to \cd$, two pseudonatural transformations\\$\epsilon: FU \Rightarrow 1_\cd$, $\eta: 1_\cc \Rightarrow UF$ and two  invertible modifications $\Psi: \epsilon F \circ F\eta \to 1_F$ and $\Phi: 1_U \to U\epsilon \circ \eta U$ (not necessarily satisfying the swallowtail identities),
\item there is a pseudofunctor $F: \cc \to \cd$ and a collection of equivalences:
\[
\lambda_{A,B}: \cd(FA, B) \simeq \cc(A,UB),
\]
for each pair of objects $A \in \ob \cc$, $B \in \ob \cd$, that is pseudonatural in $A,B$,
\item the same as above, but each $\lambda_{A,B}$ is an adjoint equivalence,
\item there is a function $F: \ob \cc \to \ob \cd$ and 1-cells $\eta_A: A \to UFA$ for each object $A \in \ob \cc$ such that the induced functor:
\[
U(-)\eta_A: \cd(FA,B) \to \cc(A,UB),
\]
is an equivalence for each pair of objects $A \in \ob \cc$, $B \in \ob \cd$.
\end{enumerate}
\end{cor}
\begin{proof}
``$1) \Rightarrow 2)$” is obvious. For ``$2) \Rightarrow 3)$” note that the following functors are equivalence inverses to each other:
\[\begin{tikzcd}
	{\cd(FA, B)} && {\cc(A,UB)}
	\arrow[""{name=0, anchor=center, inner sep=0}, "{(\eta_A)^* \circ U:}", curve={height=-12pt}, from=1-1, to=1-3]
	\arrow[""{name=1, anchor=center, inner sep=0}, "{(\epsilon_B)_* \circ F}", curve={height=-12pt}, from=1-3, to=1-1]
	\arrow["\simeq"{description}, draw=none, from=1, to=0]
\end{tikzcd}\]
as witnessed by these invertible 2-cells:
\[\begin{tikzcd}
	FA &&& A && UFA \\
	FUFA && FA & UB && UFUB \\
	FUB && B &&& UB
	\arrow["{F\eta_A}", from=1-1, to=2-1]
	\arrow[""{name=0, anchor=center, inner sep=0}, equals, from=1-1, to=2-3]
	\arrow[""{name=0p, anchor=center, inner sep=0}, phantom, from=1-1, to=2-3, start anchor=center, end anchor=center]
	\arrow[""{name=1, anchor=center, inner sep=0}, "{F(Uh\eta_A)}"', shift right=5, curve={height=30pt}, from=1-1, to=3-1]
	\arrow[""{name=1p, anchor=center, inner sep=0}, phantom, from=1-1, to=3-1, start anchor=center, end anchor=center, shift right=5, curve={height=30pt}]
	\arrow[""{name=2, anchor=center, inner sep=0}, curve={height=6pt}, draw=none, from=1-1, to=3-1]
	\arrow[""{name=2p, anchor=center, inner sep=0}, phantom, from=1-1, to=3-1, start anchor=center, end anchor=center, curve={height=6pt}]
	\arrow[""{name=3, anchor=center, inner sep=0}, "{\eta_A}", from=1-4, to=1-6]
	\arrow[""{name=3p, anchor=center, inner sep=0}, phantom, from=1-4, to=1-6, start anchor=center, end anchor=center]
	\arrow["g"', from=1-4, to=2-4]
	\arrow["UFg"', from=1-6, to=2-6]
	\arrow[""{name=4, anchor=center, inner sep=0}, "{U(\epsilon_B Fg)}", shift left=5, curve={height=-30pt}, from=1-6, to=3-6]
	\arrow[""{name=4p, anchor=center, inner sep=0}, phantom, from=1-6, to=3-6, start anchor=center, end anchor=center, shift left=5, curve={height=-30pt}]
	\arrow[""{name=5, anchor=center, inner sep=0}, curve={height=-6pt}, draw=none, from=1-6, to=3-6]
	\arrow[""{name=5p, anchor=center, inner sep=0}, phantom, from=1-6, to=3-6, start anchor=center, end anchor=center, curve={height=-6pt}]
	\arrow[""{name=6, anchor=center, inner sep=0}, "{\epsilon_{FA}}"{description, pos=0.7}, from=2-1, to=2-3]
	\arrow[""{name=6p, anchor=center, inner sep=0}, phantom, from=2-1, to=2-3, start anchor=center, end anchor=center]
	\arrow[""{name=6p, anchor=center, inner sep=0}, phantom, from=2-1, to=2-3, start anchor=center, end anchor=center]
	\arrow["FUh", from=2-1, to=3-1]
	\arrow["h", from=2-3, to=3-3]
	\arrow[""{name=7, anchor=center, inner sep=0}, "{\eta_{UB}}"{description, pos=0.3}, from=2-4, to=2-6]
	\arrow[""{name=7p, anchor=center, inner sep=0}, phantom, from=2-4, to=2-6, start anchor=center, end anchor=center]
	\arrow[""{name=7p, anchor=center, inner sep=0}, phantom, from=2-4, to=2-6, start anchor=center, end anchor=center]
	\arrow[""{name=8, anchor=center, inner sep=0}, equals, from=2-4, to=3-6]
	\arrow[""{name=8p, anchor=center, inner sep=0}, phantom, from=2-4, to=3-6, start anchor=center, end anchor=center]
	\arrow["{U\epsilon_B}"', from=2-6, to=3-6]
	\arrow[""{name=9, anchor=center, inner sep=0}, "{\epsilon_B}"', from=3-1, to=3-3]
	\arrow[""{name=9p, anchor=center, inner sep=0}, phantom, from=3-1, to=3-3, start anchor=center, end anchor=center]
	\arrow["\gamma"{pos=0.2}, shorten >=17pt, Rightarrow, from=1p, to=2p]
	\arrow["{\gamma^{-1}}"{pos=0.7}, shorten <=17pt, shorten >=3pt, Rightarrow, from=5p, to=4p]
	\arrow["{\Psi_A}", shorten <=2pt, shorten >=2pt, Rightarrow, from=6p, to=0p]
	\arrow["{\eta_g}", shorten <=4pt, shorten >=4pt, Rightarrow, from=7p, to=3p]
	\arrow["{\Phi_B}"', shorten <=2pt, shorten >=2pt, Rightarrow, from=8p, to=7p]
	\arrow["{\epsilon_h}"', shorten <=4pt, shorten >=4pt, Rightarrow, from=9p, to=6p]
\end{tikzcd}\]
For the pseudonaturality in $A,B$, we may argue as in Remark \ref{POZN_colax_adj_adjoint_homcats_corefl_refl}.
``$3) \Rightarrow 4)$”:
For any pair $A,B$, the equivalence $\lambda_{A,B}$ can be promoted to an adjoint equivalence. By the Doctrinal adjunction (Theorem \ref{THM_doctrinal_adj}), the resulting adjoint equivalence lives in $\Hom{\cc \times \cd}{\cat}$.
``$4) \Rightarrow 5)$”:
By the Yoneda Lemma for pseudofunctors, the pseudonatural equivalence:
\[
\lambda_{A,-}: \cd(FA, - ) \to \cc(A,U-),
\]
is isomorphic to a pseudonatural equivalence whose component at $X$ sends $g: FA \to X$ to $Ug \circ \eta_A$ for some 1-cell $\eta_A : A \to UFA$.

``$5) \Rightarrow 1)$”:
We will prove that $\eta_A$'s are coherently closed for $U$-extensions. Notice first that since $U(-)\eta_A$ is fully faithful, given any 1-cell $f: A \to UB$, any invertible 2-cell $\kappa: f \cong Uf' \circ \eta_A$ is a $U$-extension. Since $U(-)\eta_A$ is essentially surjective, \textbf{choose} a $U$-extension $(f',\mathbb{D}_f)$ for each 1-cell $f: A \to UB$ -- because of what we remarked in the previous sentence, this is a collection coherently closed for $U$-extensions. For the same reason, the additional conditions in Theorem \ref{THM_Bunge} are satisfied and so there is a colax adjunction $(\Psi, \Phi): (\epsilon, \eta): F \ddashv U$ with $\Psi$ invertible. Going through the proof of this theorem, we immediately see that $\Phi$ is invertible and $\eta, \epsilon$ are pseudonatural -- we thus obtain a biadjunction.

\end{proof}

\begin{theo}\label{THM_Bunge_opak}
Let $(\Psi, \Phi): (\epsilon, y): F \ddashv U: \cd \to \cc$ be a colax adjunction between pseudofunctors in which $\Psi$ is invertible. Then:
\begin{itemize}
\item the components of the unit $y_A : A \to UFA$ are coherently closed for $U$-extensions,
\item the unit and composition axioms \eqref{EQ_U_extensions_unit_comp} for $U$-extensions hold.
\end{itemize}
\end{theo}
\begin{proof}

Recall that by Proposition \ref{THMP_colax_adj_adjoint_homcats}  we have the following adjunction:
\[\begin{tikzcd}
	{\cd(FA, B)} && {\cc(A,UB)}
	\arrow[""{name=0, anchor=center, inner sep=0}, "{(y_A)^* \circ U:}", curve={height=-24pt}, from=1-1, to=1-3]
	\arrow[""{name=1, anchor=center, inner sep=0}, "{(\epsilon_B)_* \circ F}", curve={height=-24pt}, from=1-3, to=1-1]
	\arrow["\dashv"{anchor=center, rotate=90}, draw=none, from=1, to=0]
\end{tikzcd}\]

Denote by $\mbbd_g$ the unit of this adjunction evaluated at $g: A \to UB$ and denote $g^\mbbd := \epsilon_B Fg$. By the definition of the unit, the pair $(g^\mbbd, \mbbd_g)$ is the left $U$-extension of $g$ along $y_A$. Next, notice that for $f: A \to B$, the invertible 2-cell that we will denote by $\beth_f$:
\[\begin{tikzcd}
	FA && FB && FUFB \\
	\\
	&&&& FB
	\arrow["Ff", from=1-1, to=1-3]
	\arrow[""{name=0, anchor=center, inner sep=0}, "{F(y_B \circ f)}", shift left, curve={height=-30pt}, from=1-1, to=1-5]
	\arrow[""{name=0p, anchor=center, inner sep=0}, phantom, from=1-1, to=1-5, start anchor=center, end anchor=center, shift left, curve={height=-30pt}]
	\arrow[""{name=1, anchor=center, inner sep=0}, draw=none, from=1-1, to=1-5]
	\arrow[""{name=1p, anchor=center, inner sep=0}, phantom, from=1-1, to=1-5, start anchor=center, end anchor=center]
	\arrow["{Fy_B}", from=1-3, to=1-5]
	\arrow[""{name=2, anchor=center, inner sep=0}, Rightarrow, no head, from=1-3, to=3-5]
	\arrow["{\epsilon_{FB}}", from=1-5, to=3-5]
	\arrow["{\gamma^{-1}}"{pos=0.3}, shorten >=9pt, Rightarrow, from=0p, to=1p]
	\arrow["{\Psi_B}", shorten >=6pt, Rightarrow, from=1-5, to=2]
\end{tikzcd}\]
satisfies the following equality (this again follows from a swallowtail identity):
\[\begin{tikzcd}
	A && UFA &&& A && UFA \\
	&&&& {=} \\
	B && UFB &&& B && UFB
	\arrow[""{name=0, anchor=center, inner sep=0}, "{y_A}", from=1-1, to=1-3]
	\arrow["f"', from=1-1, to=3-1]
	\arrow[""{name=1, anchor=center, inner sep=0}, "{U(y_B f)^\mbbd}"'{pos=0.7}, from=1-3, to=3-3]
	\arrow[""{name=1p, anchor=center, inner sep=0}, phantom, from=1-3, to=3-3, start anchor=center, end anchor=center]
	\arrow[""{name=2, anchor=center, inner sep=0}, "UFf", curve={height=-30pt}, from=1-3, to=3-3]
	\arrow[""{name=2p, anchor=center, inner sep=0}, phantom, from=1-3, to=3-3, start anchor=center, end anchor=center, curve={height=-30pt}]
	\arrow[""{name=3, anchor=center, inner sep=0}, "{y_A}", from=1-6, to=1-8]
	\arrow["f"', from=1-6, to=3-6]
	\arrow["UFf", from=1-8, to=3-8]
	\arrow[""{name=4, anchor=center, inner sep=0}, "{y_B}"', from=3-1, to=3-3]
	\arrow[""{name=5, anchor=center, inner sep=0}, "{y_B}"', from=3-6, to=3-8]
	\arrow["{U\beth_f}", shorten <=6pt, shorten >=6pt, Rightarrow, from=1p, to=2p]
	\arrow["{\mbbd_{y_Bf}}", shorten <=9pt, shorten >=9pt, Rightarrow, from=4, to=0]
	\arrow["{y_f}"', shorten <=9pt, shorten >=9pt, Rightarrow, from=5, to=3]
\end{tikzcd}\]
This proves that $(UFf, y_f)$ is also a $U$-extension of $y_B f$ along $y_A$. Next, notice that for an object $B$, the invertible 2-cell that we will denote by $\Xi_B$:
\[\begin{tikzcd}
	FUB && FUB && B
	\arrow[""{name=0, anchor=center, inner sep=0}, Rightarrow, no head, from=1-1, to=1-3]
	\arrow[""{name=1, anchor=center, inner sep=0}, "{F1_{UB}}", curve={height=-24pt}, from=1-1, to=1-3]
	\arrow["{\epsilon_B}", from=1-3, to=1-5]
	\arrow["\iota", shorten <=3pt, shorten >=3pt, Rightarrow, from=1, to=0]
\end{tikzcd}\]

satisfies this equality:
\[\begin{tikzcd}
	UB && UFUB &&& UB && UFUB \\
	&&&& {=} \\
	&& UB &&&&& UB
	\arrow[""{name=0, anchor=center, inner sep=0}, Rightarrow, no head, from=1-1, to=3-3]
	\arrow[""{name=0p, anchor=center, inner sep=0}, phantom, from=1-1, to=3-3, start anchor=center, end anchor=center]
	\arrow["{y_{UB}}", from=1-1, to=1-3]
	\arrow[""{name=1, anchor=center, inner sep=0}, "{U1^\mbbd_{UB}}"{pos=0.8}, from=1-3, to=3-3]
	\arrow[""{name=1p, anchor=center, inner sep=0}, phantom, from=1-3, to=3-3, start anchor=center, end anchor=center]
	\arrow[""{name=1p, anchor=center, inner sep=0}, phantom, from=1-3, to=3-3, start anchor=center, end anchor=center]
	\arrow[""{name=2, anchor=center, inner sep=0}, "{U\epsilon_{B}}", shift left=4, curve={height=-30pt}, from=1-3, to=3-3]
	\arrow[""{name=2p, anchor=center, inner sep=0}, phantom, from=1-3, to=3-3, start anchor=center, end anchor=center, shift left=4, curve={height=-30pt}]
	\arrow["{y_{UB}}", from=1-6, to=1-8]
	\arrow[""{name=3, anchor=center, inner sep=0}, "{U\epsilon_B}", from=1-8, to=3-8]
	\arrow[""{name=3p, anchor=center, inner sep=0}, phantom, from=1-8, to=3-8, start anchor=center, end anchor=center]
	\arrow[""{name=4, anchor=center, inner sep=0}, Rightarrow, no head, from=1-6, to=3-8]
	\arrow[""{name=4p, anchor=center, inner sep=0}, phantom, from=1-6, to=3-8, start anchor=center, end anchor=center]
	\arrow["{\mbbd_{1_{UB}}}", shorten <=7pt, shorten >=7pt, Rightarrow, from=0p, to=1p]
	\arrow["{U\Xi_B}", shorten <=8pt, shorten >=8pt, Rightarrow, from=1p, to=2p]
	\arrow["{\Phi_{B}}", shorten <=7pt, shorten >=7pt, Rightarrow, from=4p, to=3p]
\end{tikzcd}\]

This proves that $(U\epsilon_B, \Phi_B)$ is also a $U$-extension of $1_{UB}$ along $y_{UB}$. Using these two isomorphisms of $U$-extensions, we see that the pair consisting of the 1-cell $1^\mbbd_{UY} \circ (y_{UY} f)^\mbbd$ and the composite 2-cell \eqref{EQ_DEFI_coherently_closed_for_Uext} in Definition \ref{DEFI_coherently_closed_for_Uext} is a $U$-extension since it is isomorphic to the $U$-extension $(f^\mbbd, \mbbd_f)$. We thus obtain that the collection $y_A: A \to UFA$ is coherently closed for $U$-extensions.

Let us now prove the composition and unit axioms \eqref{EQ_U_extensions_unit_comp}. The proof that the pair $(1_{FA}, \iota^{-1}y_A)$ (portrayed below) is a $U$-extension follows immediately from the fact that it is isomorphic to the $U$-extension $(y_A^\mbbd, \mbbd_{y_A})$ via the modification $\Psi_A$ (this is the first swallowtail identity):
\[\begin{tikzcd}
	A && UFA &&& A && UFA \\
	&&&& {=} \\
	&& UFA &&&&& UFA
	\arrow["{y_A}", from=1-1, to=1-3]
	\arrow["{y_A}"', from=1-1, to=3-3]
	\arrow[""{name=0, anchor=center, inner sep=0}, "{U1_{FA}}", curve={height=-30pt}, from=1-3, to=3-3]
	\arrow[""{name=0p, anchor=center, inner sep=0}, phantom, from=1-3, to=3-3, start anchor=center, end anchor=center, curve={height=-30pt}]
	\arrow[""{name=1, anchor=center, inner sep=0}, Rightarrow, no head, from=1-3, to=3-3]
	\arrow[""{name=1p, anchor=center, inner sep=0}, phantom, from=1-3, to=3-3, start anchor=center, end anchor=center]
	\arrow["{y_A}", from=1-6, to=1-8]
	\arrow[""{name=2, anchor=center, inner sep=0}, "{Uy_A^\mbbd}"{pos=0.8}, from=1-8, to=3-8]
	\arrow[""{name=2p, anchor=center, inner sep=0}, phantom, from=1-8, to=3-8, start anchor=center, end anchor=center]
	\arrow[""{name=2p, anchor=center, inner sep=0}, phantom, from=1-8, to=3-8, start anchor=center, end anchor=center]
	\arrow[""{name=3, anchor=center, inner sep=0}, "{y_A}"', from=1-6, to=3-8]
	\arrow[""{name=3p, anchor=center, inner sep=0}, phantom, from=1-6, to=3-8, start anchor=center, end anchor=center]
	\arrow[""{name=4, anchor=center, inner sep=0}, "{U1_{FA}}", shift left=5, curve={height=-30pt}, from=1-8, to=3-8]
	\arrow[""{name=4p, anchor=center, inner sep=0}, phantom, from=1-8, to=3-8, start anchor=center, end anchor=center, shift left=5, curve={height=-30pt}]
	\arrow["{\iota^{-1}}", shorten <=6pt, shorten >=6pt, Rightarrow, from=1p, to=0p]
	\arrow["{\mbbd_{y_A}}", shorten <=7pt, shorten >=7pt, Rightarrow, from=3p, to=2p]
	\arrow["{U\Psi_A}", shorten <=8pt, shorten >=8pt, Rightarrow, from=2p, to=4p]
\end{tikzcd}\]

Again using the isomorphism $\beth$ from above, the question whether the 2-cell below right is a $U$-extension is equivalent to asking whether the 2-cell below left is a $U$-extension:

\[\begin{tikzcd}
	A && UFA &&& A && UFA \\
	B && UFB &&& B && UFB \\
	C && UFC &&& C && UFC
	\arrow[""{name=0, anchor=center, inner sep=0}, "{y_A}", from=1-1, to=1-3]
	\arrow["f"', from=1-1, to=2-1]
	\arrow["UFf", from=1-3, to=2-3]
	\arrow[""{name=1, anchor=center, inner sep=0}, "{UF(gf)}", shift left=5, curve={height=-30pt}, from=1-3, to=3-3]
	\arrow[""{name=1p, anchor=center, inner sep=0}, phantom, from=1-3, to=3-3, start anchor=center, end anchor=center, shift left=5, curve={height=-30pt}]
	\arrow[""{name=2, anchor=center, inner sep=0}, "{y_A}", from=1-6, to=1-8]
	\arrow["f"', from=1-6, to=2-6]
	\arrow["{U(y_Bf)^\mbbd}"', from=1-8, to=2-8]
	\arrow[""{name=3, anchor=center, inner sep=0}, "{U((y_Cg)^\mbbd (y_Bf)^\mbbd)}", shift left=5, curve={height=-30pt}, from=1-8, to=3-8]
	\arrow[""{name=3p, anchor=center, inner sep=0}, phantom, from=1-8, to=3-8, start anchor=center, end anchor=center, shift left=5, curve={height=-30pt}]
	\arrow[""{name=4, anchor=center, inner sep=0}, "{y_B}"{description}, from=2-1, to=2-3]
	\arrow["g"', from=2-1, to=3-1]
	\arrow["UFg", from=2-3, to=3-3]
	\arrow[""{name=5, anchor=center, inner sep=0}, "{y_B}"{description}, from=2-6, to=2-8]
	\arrow["g"', from=2-6, to=3-6]
	\arrow["{U(y_Cg)^\mbbd}"', from=2-8, to=3-8]
	\arrow[""{name=6, anchor=center, inner sep=0}, "{y_C}"', from=3-1, to=3-3]
	\arrow[""{name=7, anchor=center, inner sep=0}, "{y_C}"{description}, from=3-6, to=3-8]
	\arrow["{y_f}"', shorten <=4pt, shorten >=4pt, Rightarrow, from=4, to=0]
	\arrow["\gamma"{pos=0.3}, Rightarrow, from=2-3, to=1p]
	\arrow["{\mbbd_{y_Bf}}", shorten <=4pt, shorten >=4pt, Rightarrow, from=5, to=2]
	\arrow["\gamma", Rightarrow, from=2-8, to=3p]
	\arrow["{y_g}"', shorten <=4pt, shorten >=4pt, Rightarrow, from=6, to=4]
	\arrow["{\mbbd_{y_Cg}}", shorten <=4pt, shorten >=4pt, Rightarrow, from=7, to=5]
\end{tikzcd}\]
But this 2-cell equals $y_{gf}$ and is thus a $U$-extension by what we have proven above.

\end{proof}

\begin{theo}\label{THM_korespondence_lax_adj_u_coh_ext}
Fix a pseudofunctor $U: \cd \to \cc$ between 2-categories. The following are equivalent for a collection of 1-cells $\{ y_A: A \to UFA \mid A \in \cc \}$:
\begin{itemize}
\item the collection $y_A$ is coherently closed for $U$-extensions and satisfies composition and unit axioms \eqref{EQ_U_extensions_unit_comp},
\item there is a colax adjunction $(\Psi, \Phi): (\epsilon, \eta): F \ddashv U$ for which $\Psi$ is invertible, $F$ is a pseudofunctor and the 1-cell component of the unit at each $A \in \cc$ equals  $y_A$.
\end{itemize}
\end{theo}

\begin{pozn}
In the above theorem, we do not have a one-to-one correspondence; instead, there is a suitable ``equivalence” between these two concepts. Starting with coherent $U$-extensions $(f^\mbbd, \mbbd_f)$ of $f$ along $y_A$, producing a colax adjunction and then going back to $U$-extensions gives the $U$-extension:
\[
(\epsilon_B Ff, \gamma^{-1}y_A \circ U\epsilon_B y_f \circ \Phi_B f),
\]
which in general will not be equal to $(f^\mbbd, \mbbd_f)$ (but will be canonically isomorphic to it). Similarly, starting with left colax adjoint $F$, going to $U$-extensions and back only gives a pseudofunctor isomorphic to $F$.
\end{pozn}

\noindent In our applications to two-dimensional monad theory, we will encounter this very special case of $U$-extensions:

\begin{defi}\label{DEFI_strictly_closed_for_Uexts}
Let $U: \cd \to \cc$ be a 2-functor. We will say that a collection of 1-cells $y_A: A \to UFA$ is \textit{strictly closed for $U$-extensions} if:
\begin{itemize}
\item for every $f: A \to UB$ there is a $U$-extension $(f^\mbbd, 1_f)$ along $y_A$ with the 2-cell component being the identity,
\item $y_A^\mbbd = 1_{FA}$, 
\item for $f: X \to Y$, $g: Y \to Z$, denoting $Ff := (y_Y \circ f)^\mbbd$, we have $Ff \circ Fg = F(fg)$,
\item for $f: A \to UB$, denoting $\epsilon_Y := (1_Y)^\mbbd$, we have $\epsilon_Y \circ Ff = f^\mbbd$.
\end{itemize}
\end{defi}

\begin{pozn}\label{POZN_strict_BungesTHM}
It is clear from the proof of Theorem \ref{THM_Bunge} that a collection strictly closed for $U$-extension gives rise to a colax adjunction $(\epsilon, y): F \ddashv U$ for which:
\begin{itemize}
\item y is a 2-natural transformation,
\item $F$ is a 2-functor,
\item the modifications $\Phi, \Psi$ are the identities.
\end{itemize}
(This will in general not be a 2-adjunction because $\epsilon$ will  only be colax natural.)
\end{pozn}

\section{On the Kleisli 2-category for a left Kan pseudomonad}\label{SEKCE_Kleisli2cat}

This section is devoted to studying the Kleisli 2-category for a general left Kan pseudomonad $(D,y)$ on a 2-category $\ck$. For this section, recall Proposition \ref{THMP_biadjunkce_ck_ck_D} and its notation.

In Subsection \ref{SUBS_3_charakterizace_algeber} we prove a result characterizing the pseudo $D$-algebra structure on an object in terms of the existence of certain adjoints (Theorem \ref{THM_characterization_D-algs_D-corefl-incl}).

In Subsection \ref{SUBS_3_big_lax_adj_thm} we use this result and Theorem \ref{THM_Bunge} to prove that any left biadjoint $\ck \to \cl$ that factorizes through the Kleisli 2-category gives rise to a lax left adjoint $\ck_D \to \cl$. We list several applications, one of which is the assertion that there is a canonical colax adjunction between EM and Kleisli 2-categories for left Kan pseudomonads.

A further application is given in Section \ref{SEC_weak_limits_Kleisli2cat} -- in Theorem \ref{THM_Kleisli2cat_is_corefl_complete} we show that whenever the base 2-category $\ck$ admits $J$-indexed bilimits, the Kleisli 2-category for a left Kan pseudomonad on $\ck$ will admit them as reflector-limits.

\subsection{A characterization of algebras}\label{SUBS_3_charakterizace_algeber}

\begin{defi}\label{DEFI_D_coreflector}
Let $F: \ck \to \cl$ be a pseudofunctor. We will call a morphism $f: A \to B$ in $\ck$ an \textit{$F$-coreflector} if $Ff$ is a coreflector in the 2-category $\cl$. Similarly for the other variants from Definition \ref{DEFI_lali_rali_rari_lari}.
\end{defi}

\begin{pr}
Let $P: \CAT \to \prof$ be the canonical inclusion pseudofunctor. In Example \ref{PR_presheaf_psmonad_charakterizace_algeber} below we will show that a functor $f: \ca \to \cb$ between locally small categories is a $P$-coreflection-inclusion if and only if it is fully faithful and satisfies a certain smallness condition.
\end{pr}

\begin{pr}\label{PR_doctrinal_adj}
Consider the lax morphism classifier 2-comonad $Q_l$ associated to a 2-monad $T$ on a 2-category $\ck$ (see \ref{DEFI_Ql_Qp}). Denote by $J: \talgs \to \talgl$ the canonical inclusion to the Kleisli 2-category and by $U: \talgs \to \ck$ the forgetful 2-functor. By the doctrinal adjunction (Theorem \ref{THM_doctrinal_adj}), a strict morphism is a $J$-reflector if and only if it is a $U$-reflector, that is, the underlying morphism in $\ck$ is a reflector.
\end{pr}

\begin{pr}
Given bicategories $\ck, \cm$ and a \textit{proarrow equipment}, which is a certain pseudofunctor $(-)_*: \ck \to \cm$ (see \cite[Section 1]{abstractproarrows}), a 1-cell in $\ck$ is a $(-)_*$-coreflection-inclusion if and only if it is fully faithful in the sense of Wood \cite[Page 10]{abstractproarrows}.
\end{pr}

\begin{pr}
Given a lax-idempotent pseudomonad $P$ on a 2-category $\ck$, 1-cells in $\ck$ that are $P$-left adjoints have been studied in the literature (\cite{bungebicomma}, \cite{distrlawsadmiss}, \cite{adjfunthms}) under the name \textit{$P$-admissible} 1-cells.
\end{pr}

\noindent The following lemma is the left Kan pseudomonad version of \cite[Theorem 3.4]{threepieces}:

\begin{lemma}\label{THM_lemma_Dcoreflincl_is_JDcoreflincl}
Let $(D,y)$ be a left Kan pseudomonad on $\ck$. Denote by $D: \ck \to \ck$ the corresponding endo-pseudofunctor and by $J_D: \ck \to \ck_D$ the inclusion to the Kleisli 2-category. The following are equivalent for a 1-cell $f: B \to C$:
\begin{itemize}
\item $f$ is a $D$-coreflection-inclusion,
\item $f$ is a $J_D$-coreflection-inclusion.
\end{itemize}
\end{lemma}
\begin{proof}
``$(1) \Rightarrow (2)$” follows from \cite[Proposition 1.3]{bungebicomma}: namely, the right adjoint to $Df$ in $\ck$ is actually a pseudo morphism and thus is an adjoint in $\ck_D$. ``$(2) \Rightarrow (1)$” is obvious because we have the forgetful 2-functor $U_D: \ck_D \to \ck$ that satisfies $D = U_D J_D$.
\end{proof}

\begin{lemma}\label{THM_Lemma_Kan_extensions}
The following holds in a 2-category $\ck$:
\begin{itemize}
\item Let $f \dashv u: B \to A$ be an adjunction with unit $\eta$ and let $(\mbbd,g^\mbbd)$ be the left Kan extension of $g: A' \to C$ along $y: A' \to A$. Then the diagram below left exhibits $g^\mbbd u$ as the left Kan extension of $g$ along $f y$:
\[\begin{tikzcd}
	B &&&&& C \\
	A && A &&& B \\
	{A'} &&&& C & {B'} && D
	\arrow["u", from=1-1, to=2-3]
	\arrow["h", curve={height=-24pt}, from=1-6, to=3-8]
	\arrow["f", from=2-1, to=1-1]
	\arrow["\eta"'{pos=0.1}, shift right=5, shorten >=5pt, Rightarrow, from=2-1, to=1-1]
	\arrow[Rightarrow, no head, from=2-1, to=2-3]
	\arrow["{g^\mathbb{D}}", from=2-3, to=3-5]
	\arrow["k", from=2-6, to=1-6]
	\arrow["hk", curve={height=-18pt}, from=2-6, to=3-8]
	\arrow["y", from=3-1, to=2-1]
	\arrow["{\mathbb{D}}"', shift right=5, shorten <=1pt, shorten >=1pt, Rightarrow, from=3-1, to=2-1]
	\arrow["g"', from=3-1, to=3-5]
	\arrow["\iota", from=3-6, to=2-6]
	\arrow["\alpha"', shift right=5, shorten >=2pt, Rightarrow, from=3-6, to=2-6]
	\arrow["f"', from=3-6, to=3-8]
\end{tikzcd}\]
\item In the diagram above right, suppose that the top and outer diagrams are left Kan extensions of $f$. If left Kan extensions along $k$ exist in $\ck$ and have an invertible unit, then $\alpha$ is the left Kan extension of $f$ along $\iota$.
\end{itemize}
\end{lemma}
\begin{proof}
The first point follows by composing the following bijections. For a 1-cell $h: B \to C$, the first one is given by the adjunction $f \dashv u$, the second one is given by the definition of $g^\mbbd$:
\[
\ck(B,C)(g^\mbbd u, h) \cong \ck(A,C)(g^\mbbd, h f) \cong \ck(A',C)(g, hfy).
\]
In the second point, assume we have a 2-cell $\beta$ as pictured below, and we want to find a unique 2-cell solving this equation:
\begin{equation}\label{EQ_PROOF_alphathetabeta}
\begin{tikzcd}
	B &&&& B \\
	&&& {=} \\
	{B'} && D && {B'} && D
	\arrow["f"', from=3-1, to=3-3]
	\arrow["\iota", from=3-1, to=1-1]
	\arrow[""{name=0, anchor=center, inner sep=0}, "hk"{description}, from=1-1, to=3-3]
	\arrow["\iota", from=3-5, to=1-5]
	\arrow["f"', from=3-5, to=3-7]
	\arrow[""{name=1, anchor=center, inner sep=0}, "l", curve={height=-30pt}, from=1-1, to=3-3]
	\arrow[""{name=2, anchor=center, inner sep=0}, "l", curve={height=-18pt}, from=1-5, to=3-7]
	\arrow["\alpha"', shorten <=6pt, shorten >=6pt, Rightarrow, from=3-1, to=0]
	\arrow["?", shorten <=5pt, shorten >=5pt, Rightarrow, from=0, to=1]
	\arrow["\beta"', shorten <=8pt, shorten >=8pt, Rightarrow, from=3-5, to=2]
\end{tikzcd}
\end{equation}

First note that we have a unique 2-cell $\theta'$ making the following diagram equal (here $(l^\mbba, \mbba)$ is the left Kan extension of $l$ along $k$):
\[\begin{tikzcd}
	C &&&&& C & {} \\
	B & {} &&& {=} & B & {} \\
	{B'} & {} && D && {B'} & {} && D
	\arrow["f"', from=3-1, to=3-4]
	\arrow["\iota", from=3-1, to=2-1]
	\arrow["hk"{description}, from=2-1, to=3-4]
	\arrow["k", from=2-1, to=1-1]
	\arrow[""{name=0, anchor=center, inner sep=0}, "h"', curve={height=-6pt}, from=1-1, to=3-4]
	\arrow[""{name=1, anchor=center, inner sep=0}, "{l^\mbba}", curve={height=-30pt}, from=1-1, to=3-4]
	\arrow["f"', from=3-6, to=3-9]
	\arrow["k", from=2-6, to=1-6]
	\arrow["\iota", from=3-6, to=2-6]
	\arrow["{l^\mbba}", curve={height=-30pt}, from=1-6, to=3-9]
	\arrow["l"{description}, curve={height=-6pt}, from=2-6, to=3-9]
	\arrow["\beta", shorten <=1pt, shorten >=4pt, Rightarrow, from=3-7, to=2-7]
	\arrow["\alpha"'{pos=0.3}, shift left=5, shorten >=9pt, Rightarrow, from=3-2, to=2-2]
	\arrow["\mbba"{pos=0.7}, shorten >=3pt, Rightarrow, from=2-7, to=1-7]
	\arrow["{\theta'}"', shorten <=4pt, shorten >=4pt, Rightarrow, from=0, to=1]
\end{tikzcd}\]
Clearly, $\theta := \mbba^{-1} \circ \theta' k$ solves the equation \eqref{EQ_PROOF_alphathetabeta}, giving us \textbf{the existence} part of the proof. To show \textbf{the uniqueness}, let $\phi$ be a different 2-cell solving \eqref{EQ_PROOF_alphathetabeta}. Note that there exists a unique 2-cell $\phi'$ solving the following:
\[\begin{tikzcd}
	C &&&& C \\
	\\
	B &&& {=} & B \\
	&& D &&&& D
	\arrow["k", from=3-1, to=1-1]
	\arrow[""{name=0, anchor=center, inner sep=0}, "hk"', from=3-1, to=4-3]
	\arrow["hk"', from=3-5, to=4-7]
	\arrow["k", from=3-5, to=1-5]
	\arrow[""{name=1, anchor=center, inner sep=0}, "{l^\mbba}", curve={height=-30pt}, from=1-5, to=4-7]
	\arrow[""{name=2, anchor=center, inner sep=0}, "h"{description}, from=1-5, to=4-7]
	\arrow[""{name=3, anchor=center, inner sep=0}, "{l^\mbba}", curve={height=-30pt}, from=1-1, to=4-3]
	\arrow[""{name=4, anchor=center, inner sep=0}, "l"{description, pos=0.3}, curve={height=-24pt}, from=3-1, to=4-3]
	\arrow["{\phi'}"', shorten <=5pt, shorten >=5pt, Rightarrow, from=2, to=1]
	\arrow["\mbba", shorten <=4pt, shorten >=4pt, Rightarrow, from=4, to=3]
	\arrow["\phi"{pos=0.3}, shorten <=3pt, shorten >=3pt, Rightarrow, from=0, to=4]
\end{tikzcd}\]
Pre-pasting this with $\alpha$ and using the diagram above this one, we see that $\phi' = \theta'$. From this we obtain:
\[
\mbba^{-1} \circ \theta' k = \mbba^{-1} \circ \phi' k = \mbba^{-1} \circ \mbba \circ \phi = \phi.
\]

\end{proof}

\begin{prop}\label{THM_characterization_D-algs_D-corefl-incl}
Let $(D, y)$ be a left Kan pseudomonad on a 2-category $\ck$. Denote by $J_D$ the inclusion to the Kleisli 2-category and by $D$ the endo-pseudofunctor associated to the left Kan pseudomonad. The following are equivalent for an object $A \in \ck$:
\begin{enumerate}
\item $A$ admits the structure of a pseudo $D$-algebra,
\item for every object $B \in \ck$, the left Kan extension of a 1-cell $B \to A$ along $y_B: B \to DB$ exists and has invertible unit. In other words, the 2-functor $\ck(-,A): \ck^{op} \to \cat$ sends each $y_B$ to a coreflector,
\item $y_A$ is a reflection-inclusion (admits a left adjoint with an invertible counit),
\item $\ck(-,A)$ sends $J_D$-coreflection-inclusions in $\ck$ to coreflectors in $\cat$,
\item $\ck(-,A)$ sends $D$-coreflection-inclusions in $\ck$ to coreflectors in $\cat$.
\end{enumerate}
\end{prop}
\begin{proof}
The equivalence ``$(1) \Leftrightarrow (3)$” is well known, the lax-idempotent pseudomonad version has been done in our treatment in Section \ref{SUBS_lax_idemp_props}. The implication ``$(1) \Rightarrow (2)$” is obvious.

For ``$(2) \Rightarrow (3)$”, denote by $(a: DA \to A,\mbba)$ the left Kan extension of $1_A$ along $y_A$. Because the identity 2-cell on $DA$ exibits $1_{DA}$ as the left Kan extension of $y_A$ along $y_A$, there exists a unique 2-cell $\eta$ making these 2-cells equal:
\[\begin{tikzcd}
	& DA &&&&& A \\
	A && A & DA & A & DA && DA
	\arrow[""{name=0, anchor=center, inner sep=0}, Rightarrow, no head, from=2-1, to=2-3]
	\arrow["{y_A}"', from=2-3, to=2-4]
	\arrow["{y_A}", from=2-1, to=1-2]
	\arrow["a", from=1-2, to=2-3]
	\arrow["{y_A}", from=2-5, to=2-6]
	\arrow[""{name=1, anchor=center, inner sep=0}, Rightarrow, no head, from=2-6, to=2-8]
	\arrow["a", from=2-6, to=1-7]
	\arrow["{y_A}", from=1-7, to=2-8]
	\arrow["\eta", shorten <=3pt, Rightarrow, from=1, to=1-7]
	\arrow["\mbba"', shorten <=3pt, Rightarrow, from=0, to=1-2]
\end{tikzcd}\]

We will now show that $(\mbba^{-1}, \eta): a \dashv y_A$ is an adjunction. The triangle identity $y_A \mbba^{-1} \circ \eta y_A = 1_{y_A}$ is guaranteed by the above formula -- let us prove the other one:
\[
\mbba^{-1}a \circ a\eta  = 1_a.
\]
Because $a$ is the left Kan extension along $y_A$, it suffices to prove that both sides of this equation become equal after pre-composing them with $y_A$. It then becomes:
\[
\mbba^{-1}ay_A \circ a\eta y_A = \mbba^{-1}ay_A \circ a y_A \mbba = \mbba^{-1}ay_A \circ \mbba a y_A = 1_{ay_A}.
\]

``$(4) \Leftrightarrow (5)$” follows from Lemma \ref{THM_lemma_Dcoreflincl_is_JDcoreflincl} and ``$(5) \Rightarrow (2)$” is obvious since $y_B$ is a $D$-coreflection-inclusion.

We will now prove ``$(2) \Rightarrow (5)$”. Let $f: B \to C$ such that there is an adjunction in $\ck_D$ where the \textbf{unit} $\eta$ is invertible:
\[\begin{tikzcd}
	{(\epsilon,\eta):} & DB && DC
	\arrow[""{name=0, anchor=center, inner sep=0}, "Df"', curve={height=24pt}, from=1-2, to=1-4]
	\arrow[""{name=1, anchor=center, inner sep=0}, "r"', curve={height=24pt}, from=1-4, to=1-2]
	\arrow["\dashv"{anchor=center, rotate=90}, draw=none, from=0, to=1]
\end{tikzcd}\]
We wish to show that the functor $f^*: \ck(C,A) \to \ck(B,A)$ has a left adjoint with invertible unit. We will define this left adjoint by the following formula:
\[
L: (g: B \to A) \mapsto (g^\mbba \circ r \circ y_C : C \to A)
\]
Define the component of the unit $\widetilde{\eta}$ at $g: B \to A$ as the following composite 2-cell:
\[\begin{tikzcd}
	C && {} && DC \\
	B && DB & {} & {} & DB & A \\
	&&& {}
	\arrow["{y_B}", from=2-1, to=2-3]
	\arrow["f", from=2-1, to=1-1]
	\arrow["{y_C}", from=1-1, to=1-5]
	\arrow["r", from=1-5, to=2-6]
	\arrow[Rightarrow, no head, from=2-3, to=2-6]
	\arrow["{g^\mbba}", from=2-6, to=2-7]
	\arrow["Df", from=2-3, to=1-5]
	\arrow["g"', curve={height=30pt}, from=2-1, to=2-7]
	\arrow["\eta", shorten <=3pt, Rightarrow, from=2-5, to=1-5]
	\arrow["{\mbbd^{-1}}", shorten <=2pt, Rightarrow, from=2-3, to=1-3]
	\arrow["\mbba"{pos=0.7}, shorten <=9pt, Rightarrow, from=3-4, to=2-4]
\end{tikzcd}\]

We wish to show that this has the universal property of the unit, in other words, $g^\mbba \circ r \circ y_C$ is the left Kan extension of $g: B \to A$ along $f: B \to C$.

By Lemma \ref{THM_Lemma_Kan_extensions} \textbf{point 1}, $g^\mbba \circ r$ is the left Kan extension of $g$ along $Df \circ y_B$. Equivalently, it is the left Kan extension of $g$ along $y_C f$ with the accompanying 2-cell being given by the 2-cell above. Since $g^\mbba r$ is a $D$-morphism, $g^\mbba r$ (with the identity 2-cell component) is the left Kan extension of $g^\mbba r y_C$ along $y_C$. By Lemma \ref{THM_Lemma_Kan_extensions} \textbf{point 2}, for $h := g^\mbba r$, $k := y_C$, $\iota := f$, $f := g$ and $\alpha$ the 2-cell above, the result follows.

\end{proof}

\begin{pozn}
Using the terminology of \cite[Definition 1.2]{kaninj}, in Theorem \ref{THM_characterization_D-algs_D-corefl-incl}, the equivalence ``$(1) \Leftrightarrow (4)$” says that an object $A$ is a pseudo $D$-algebra if and only if it is \textit{left Kan injective} with respect to the class of 1-cells given by $J_D$-coreflection-inclusions.

Let us also note that a version of ``$(1) \Rightarrow (5)$” for $D$-left adjoints in Theorem \ref{THM_characterization_D-algs_D-corefl-incl} has already been proven in \cite[Proposition 1.5]{bungebicomma}.
\end{pozn}

\begin{pozn}\label{POZN_strict_charakterizace_D_algeber}
Given a left Kan 2-monad $(D,y)$, a pseudo $D$-algebra $C$ will be said to be \textit{normal} if the left Kan extension 2-cell $\mbbc_f$ in Definition \ref{DEFI_D_algebra} is the identity for all 1-cells $f$. Notice that a variation of Proposition \ref{THM_characterization_D-algs_D-corefl-incl} may be proven for normal pseudo $D$-algebras, where we replace all invertible 2-cells by identities, for instance replace a ``reflector” by a ``lali”.
\end{pozn}

\noindent In the remainder of this section we will demonstrate Proposition \ref{THM_characterization_D-algs_D-corefl-incl} on the case of the small presheaf pseudomonad from Example \ref{PR_presheafpsmonad}. An application to the lax morphism classifier 2-comonads will be described in Section \ref{SEKCE_applications_to_twodim}.

\begin{pr}\label{PR_presheaf_psmonad_charakterizace_algeber}
Consider the small presheaf pseudomonad $P$ on $\CAT$. Note that in the virtual double category $\PROF$ of locally small categories and \textbf{all} profunctors (\cite[Example 2.9]{unified}), for any functor $f: \ca \to \cb$, the small profunctor $Pf = \cb(-,f-): \cb^{op} \times \ca \to \set$ has a right adjoint:
\[\begin{tikzcd}
	\ca && \cb
	\arrow[""{name=0, anchor=center, inner sep=0}, "{\cb(-,f-)}"', curve={height=24pt}, from=1-1, to=1-3]
	\arrow[""{name=1, anchor=center, inner sep=0}, "{\cb(f-,-)}"', curve={height=24pt}, from=1-3, to=1-1]
	\arrow["\dashv"{anchor=center, rotate=90}, draw=none, from=0, to=1]
\end{tikzcd}\]

We will call the functor $f: \ca \to \cb$ \textit{small} if the right adjoint is also a small profunctor. Clearly, this happens if and only if $Pf$ has a right adjoint in the bicategory $\prof$ of locally small categories and small profunctors.

Next, note that the unit of the adjunction is a collection of functions for every pair $(a',a'') \in \ca^{op} \times \ca$ like this:
\begin{align*}
\ca(a',a'') &\to \int^{b \in \cb} \cb(fa',b) \times \cb(b,fa''),\\
(\theta: a' \to a'') &\mapsto [1_{fa'},f(\theta)].
\end{align*}
As is readily seen, the unit is invertible if and only if $f$ is fully faithful. So a functor $f: \ca \to \cb$ is a $P$-coreflection-inclusion if and only if it is fully faithful and small. The precomposition functor $f^* : \CAT(\cb, \cc) \to \CAT(\ca, \cc)$ is a coreflector if and only if left Kan extensions along $f$ exist in $\cc$. Theorem \ref{THM_characterization_D-algs_D-corefl-incl} for the small presheaf pseudomonad now gives a folklore result: a category $\cc$ is cocomplete if and only if left Kan extensions along small fully faithful functors exist in $\cc$.
\end{pr}

\begin{pozn}[One-dimensional situation]\label{POZN_idemp_monad_char_algeber}
For an idempotent monad $(T,\mu, \eta)$ on a category $\cc$ (regarded as a 2-monad on a locally discrete 2-category), Proposition \ref{THM_characterization_D-algs_D-corefl-incl} in particular gives the folklore characterization of its algebras -- $A$ is a $T$-algebra if and only if $\eta_A$ is an isomorphism.
\end{pozn}

\subsection{Colax adjunctions out of the Kleisli 2-category}\label{SUBS_3_big_lax_adj_thm}

\begin{prop}\label{THM_bigadjthm_HL_is_D-algebra}
Let $(D, y)$ be a left Kan pseudomonad on a 2-category $\ck$ and assume there are pseudofunctors $G,H$ and a biadjunction as pictured below:
\[\begin{tikzcd}
	{(\sigma, \tau):(s,c):} & \ck && {\ck_D} && \cl
	\arrow["{J_D}"', from=1-2, to=1-4]
	\arrow["G"', from=1-4, to=1-6]
	\arrow[""{name=0, anchor=center, inner sep=0}, "H"', curve={height=30pt}, from=1-6, to=1-2]
	\arrow["\dashv"{anchor=center, rotate=90}, draw=none, from=1-4, to=0]
\end{tikzcd}\]
Then for every object $L \in \cl$, the object $HL$ admits the structure of a pseudo $D$-algebra.
\end{prop}
\begin{proof}
By Proposition \ref{THM_characterization_D-algs_D-corefl-incl}, it suffices to show that $\ck(-,HL): \ck^{op} \to \cat$ sends $J_D$-coreflection-inclusions to coreflectors. Notice that we have the following pseudonatural equivalence:
\[
\ck(-,HL) \simeq \cl(GJ_D-,L) = \cl(G-,L) \circ J_D.
\]
Now, by definition, $J_D$ sends $J_D$-coreflection-inclusions to coreflection-inclusions. Since $\cl(G-,L)$ is a (contravariant) pseudofunctor, it sends coreflection-inclusions to coreflectors. We thus obtain the result.
\end{proof}

\begin{pozn}\label{POZN_explicit_form_h_L}
Recall the transformation $p: J_D U_D \Rightarrow 1_{\ck_D}$ and the modification \\$\Psi: pJ_D \circ J_D y \cong 1_{J_D}$ from Proposition \ref{THMP_biadjunkce_ck_ck_D}. Going through the proof of Proposition \ref{THM_characterization_D-algs_D-corefl-incl} for the case of $HL$, we see that the algebra multiplication map $h_L : DHL \to HL$ (the reflector of the morphism $y_{HL}: HL \to DHL$) is given by the following composite:
\[\begin{tikzcd}
	DHL && {HGD^2HL} && HGDHL && HL
	\arrow["{c_{DHL}}", from=1-1, to=1-3]
	\arrow["{HGp_{DHL}}", from=1-3, to=1-5]
	\arrow["{Hs_L}", from=1-5, to=1-7]
\end{tikzcd}\]
Also, the counit of the adjunction $h_L \dashv y_{HL}$, an invertible 2-cell $\epsilon_L: h_L y_{HL} \Rightarrow 1_{HL}$, is given by the following:
\[\begin{tikzcd}
	& DHL & {HGD^2HL} && HGDHL \\
	HL & HGDHL & {} &&& HL \\
	&& {}
	\arrow["{y_{HL}}", from=2-1, to=1-2]
	\arrow["{c_{DHL}}", from=1-2, to=1-3]
	\arrow["{HGp_{DHL}}", from=1-3, to=1-5]
	\arrow["{Hs_L}", from=1-5, to=2-6]
	\arrow["{c_{HL}}"', from=2-1, to=2-2]
	\arrow["{HGDy_{HL}}"{description}, from=2-2, to=1-3]
	\arrow[curve={height=18pt}, Rightarrow, no head, from=2-2, to=1-5]
	\arrow[shift right=2, curve={height=30pt}, Rightarrow, no head, from=2-1, to=2-6]
	\arrow["{c_{y_{HL}}}"', shorten <=2pt, Rightarrow, from=1-2, to=2-2]
	\arrow["{(HG\Psi)_{HL}}"{pos=0.4}, shorten >=4pt, Rightarrow, from=1-3, to=2-3]
	\arrow["{\tau^{-1}_L}"{pos=0.4}, shorten >=1pt, Rightarrow, from=2-3, to=3-3]
\end{tikzcd}\]
\end{pozn}

\begin{theo}[\textbf{The main colax adjunction theorem}]\label{THM_BIG_lax_adj_thm}
Let $(D, y)$ be a left Kan pseudomonad on a 2-category $\ck$. Any biadjunction whose left adjoint factorizes through the Kleisli 2-category $\ck_D$ induces a colax adjunction pictured below:
\[\begin{tikzcd}
	\ck && {\ck_D} && \cl & \rightsquigarrow & {\ck_D} && \cl
	\arrow["{J_D}"', from=1-1, to=1-3]
	\arrow["G"', from=1-3, to=1-5]
	\arrow[""{name=0, anchor=center, inner sep=0}, "H"', curve={height=30pt}, from=1-5, to=1-1]
	\arrow[""{name=1, anchor=center, inner sep=0}, "G"', curve={height=24pt}, from=1-7, to=1-9]
	\arrow[""{name=2, anchor=center, inner sep=0}, "{J_DH}"', curve={height=24pt}, from=1-9, to=1-7]
	\arrow["\dashv"{anchor=center, rotate=90}, draw=none, from=1-3, to=0]
	\arrow["\ttop"{description}, draw=none, from=1, to=2]
\end{tikzcd}\]
\end{theo}
\begin{proof}
Denote again the unit, counit and the modifications of the biadjunction as follows:
\begin{align*}
s &: GJ_D H \Rightarrow 1_\cl, &&\sigma: sGJ_D \circ GJ_Dc \cong 1_{GJ_D} ,\\
c &: 1_\ck \Rightarrow HGJ_D,  &&\tau : 1_{H} \cong Hs \circ cH.
\end{align*}
We will show that the components of the counit $s_L: GDHL \to L$ are coherently closed for $G$-lifts. By (the dual of) Theorem \ref{THM_Bunge}, there is a right colax adjoint to $G$. We will then prove that it is isomorphic to $J_D H$.

Let us first prove the following: given a 1-cell $l: GDA \to L$ in $\cl$, \textbf{any} pair $(Dl', \lambda)$ where $l': A \to HL$ is a 1-cell and $\lambda$ is an invertible 2-cell as pictured below exhibits $Dl'$ as the right $G$-lift of $l$ along $s_L$:
\[\begin{tikzcd}
	L && GDHL \\
	\\
	&& GDA
	\arrow["{GDl'}"', from=3-3, to=1-3]
	\arrow["{s_L}"', from=1-3, to=1-1]
	\arrow[""{name=0, anchor=center, inner sep=0}, "l", from=3-3, to=1-1]
	\arrow["\lambda"', shorten >=6pt, Rightarrow, from=1-3, to=0]
\end{tikzcd}\]

By Theorem \ref{THM_bigadjthm_HL_is_D-algebra}, $HL$ has the structure of a $D$-algebra. Denoting its multiplication map by $h_L$ as in Remark \ref{POZN_explicit_form_h_L}, we have the following composite adjunction with an invertible counit:
\begin{equation}\label{EQ_hL_yA_yHL_JD}
\begin{tikzcd}
	{\ck(A,HL)} && {\ck(A,DHL)} && {\ck_D(DA,DHL)}
	\arrow[""{name=0, anchor=center, inner sep=0}, "{(y_{HL})_*}"', from=1-1, to=1-3]
	\arrow[""{name=0p, anchor=center, inner sep=0}, phantom, from=1-1, to=1-3, start anchor=center, end anchor=center]
	\arrow[""{name=1, anchor=center, inner sep=0}, "{J_D}"', curve={height=30pt}, from=1-1, to=1-5]
	\arrow[""{name=2, anchor=center, inner sep=0}, "{(h_L)_*}"', curve={height=30pt}, from=1-3, to=1-1]
	\arrow[""{name=2p, anchor=center, inner sep=0}, phantom, from=1-3, to=1-1, start anchor=center, end anchor=center, curve={height=30pt}]
	\arrow[""{name=3, anchor=center, inner sep=0}, "{p_{DHL}J_D(-)}"', from=1-3, to=1-5]
	\arrow[""{name=4, anchor=center, inner sep=0}, "{U_D(-)y_A}"', curve={height=30pt}, from=1-5, to=1-3]
	\arrow[""{name=5, anchor=center, inner sep=0}, "{(-)^\#}"', curve={height=60pt},shift right=3, from=1-5, to=1-1]
	\arrow["\dashv"{anchor=center, rotate=-93}, draw=none, from=2p, to=0p]
	\arrow["\simeq"{description}, draw=none, from=3, to=4]
	\arrow["{:=}"{marking, allow upside down}, draw=none, from=5, to=1-3]
	\arrow["\cong"{description}, draw=none, from=1-3, to=1]
\end{tikzcd}
\end{equation}

Notice that there is an isomorphism:
\[
\beth: s_L GJ_D(-)^\# \cong s_L G(-): \ck_D(DA,DHL) \to \cl(GDA,L),
\]
with the component at $f: DA \to DHL$ being given by the 2-cell\footnote{We have not shown it in this diagram, but the 2-cell has to be pre-composed with the associators for the pseudofunctor $GD$ so that its source really equals $s_L GDf^\#$.}:
\begin{equation}\label{EQ_sigma_s_s_p_psi}
\begin{minipage}{0.8\textwidth}
\adjustbox{scale=0.9,center}{
\begin{tikzcd}
	{GD^2HL} && {GDHGD^2HL} && GDHGDHL && GDHL \\
	\\
	{GD^2A} && {GD^2HL} && GDHL && L \\
	&& GDA \\
	GDA
	\arrow[""{name=0, anchor=center, inner sep=0}, Rightarrow, no head, from=1-1, to=3-3]
	\arrow[""{name=0p, anchor=center, inner sep=0}, phantom, from=1-1, to=3-3, start anchor=center, end anchor=center]
	\arrow[""{name=1, anchor=center, inner sep=0}, "{GDc_{DHL}}", from=1-1, to=1-3]
	\arrow[""{name=1p, anchor=center, inner sep=0}, phantom, from=1-1, to=1-3, start anchor=center, end anchor=center]
	\arrow["{s_{GD^2HL}}"{description}, from=1-3, to=3-3]
	\arrow[""{name=2, anchor=center, inner sep=0}, "{Gp_{DHL}}"{description}, from=3-3, to=3-5]
	\arrow[""{name=2p, anchor=center, inner sep=0}, phantom, from=3-3, to=3-5, start anchor=center, end anchor=center]
	\arrow[""{name=3, anchor=center, inner sep=0}, "{s_L}"', from=3-5, to=3-7]
	\arrow[""{name=3p, anchor=center, inner sep=0}, phantom, from=3-5, to=3-7, start anchor=center, end anchor=center]
	\arrow[""{name=4, anchor=center, inner sep=0}, "{GDHGp_{DHL}}", from=1-3, to=1-5]
	\arrow[""{name=4p, anchor=center, inner sep=0}, phantom, from=1-3, to=1-5, start anchor=center, end anchor=center]
	\arrow[""{name=4p, anchor=center, inner sep=0}, phantom, from=1-3, to=1-5, start anchor=center, end anchor=center]
	\arrow["{s_{GDHL}}"{description}, from=1-5, to=3-5]
	\arrow[""{name=5, anchor=center, inner sep=0}, "{GDHs_L}", from=1-5, to=1-7]
	\arrow[""{name=5p, anchor=center, inner sep=0}, phantom, from=1-5, to=1-7, start anchor=center, end anchor=center]
	\arrow["{s_L}", from=1-7, to=3-7]
	\arrow["Gf"', from=4-3, to=3-5]
	\arrow[""{name=6, anchor=center, inner sep=0}, "{Gp_{DA}}", from=3-1, to=4-3]
	\arrow["GDf", from=3-1, to=1-1]
	\arrow["{GDy_A}", from=5-1, to=3-1]
	\arrow[""{name=7, anchor=center, inner sep=0}, Rightarrow, no head, from=5-1, to=4-3]
	\arrow["{(Gp)_f}", Rightarrow, from=3-3, to=4-3]
	\arrow[""{name=8, anchor=center, inner sep=0}, "{GDh_L}", shift left=2, curve={height=-30pt}, from=1-1, to=1-7]
	\arrow[""{name=8p, anchor=center, inner sep=0}, phantom, from=1-1, to=1-7, start anchor=center, end anchor=center, shift left=2, curve={height=-30pt}]
	\arrow["{s_{Gp_{DHL}}}", shorten <=9pt, shorten >=9pt, Rightarrow, from=4p, to=2p]
	\arrow["{s_{s_L}}", shorten <=9pt, shorten >=9pt, Rightarrow, from=5p, to=3p]
	\arrow["{\sigma_{DHL}}", shorten <=4pt, shorten >=4pt, Rightarrow, from=1p, to=0p]
	\arrow["{(G\Psi)_A}"', shorten <=4pt, shorten >=4pt, Rightarrow, from=6, to=7]
	\arrow["{}"{description, pos=0.2}, draw=none, from=8p, to=4p]
\end{tikzcd}
}
\end{minipage}
\end{equation}

\noindent We have the following chain of bijections:
\begin{align*}
\ck_D(DA,DHL)(f,Dl') &\stackrel{(A)}{\cong} \ck(A,HL)(f^\#,l') \\
&\stackrel{(B)}{\cong} \cl(GDA,L)(s_L \circ GDf^\#, s_L \circ GDl')\\
&\stackrel{(C)}{\cong} \cl(GDA,L)(s_L \circ Gf, s_L \circ GDl')\\
&\stackrel{(D)}{\cong} \cl(GDA,L)(s_L \circ Gf, l).\\
\end{align*}
The bijection \textbf{(A)} follows from the adjunction \eqref{EQ_hL_yA_yHL_JD} above. \textbf{(B)} is given by the action on morphisms of the following functor:
\begin{equation}\label{EQ_sLcircGD}
s_L \circ GJ_D(-): \ck(A,HL) \to \cl(GDA,L).
\end{equation}
This functor is (by assumption) an equivalence -- in particular it is fully faithful. \textbf{(C)} is given by the pre-composition with $\beth^{-1}_f$ and \textbf{(D)} is given by the post-composition with $\lambda$. To conclude that $(Dl', \lambda)$ is a $G$-lift, it has to be shown that the composite bijection is given by the assignment $\alpha \mapsto \lambda \circ s_LG\alpha$. Equivalently, the composite of the first three bijections is the assignment $\alpha \mapsto s_LG\alpha$. We prove this fact in the appendix as Lemma \ref{THML_big_l_adj_thm}.

\medskip

The pair $(Dl', \lambda)$ is thus a $G$-lift. Because the functor \eqref{EQ_sLcircGD} is essentially surjective, such a $G$-lift is guaranteed to always exist. Make now a \textbf{choice} of a lift for every $l: GDA \to L$ and denote it by $(Dl^\mbbl, \mbbl)$. To prove that our choice is coherently closed for $G$-lifts, the following 2-cell needs to be shown to be a $G$-lift of $l: GDA \to L$ along $s_L$:
\[\begin{tikzcd}
	L && GDHL \\
	\\
	GDA && GDHGDA \\
	\\
	&& GDA
	\arrow["{GD(l s_{GDA})^\mbbl}"{description}, from=3-3, to=1-3]
	\arrow[""{name=0, anchor=center, inner sep=0}, "{s_L}"', from=1-3, to=1-1]
	\arrow[""{name=0p, anchor=center, inner sep=0}, phantom, from=1-3, to=1-1, start anchor=center, end anchor=center]
	\arrow["l", from=3-1, to=1-1]
	\arrow[""{name=1, anchor=center, inner sep=0}, "{s_{GDA}}"', from=3-3, to=3-1]
	\arrow[""{name=1p, anchor=center, inner sep=0}, phantom, from=3-3, to=3-1, start anchor=center, end anchor=center]
	\arrow[""{name=1p, anchor=center, inner sep=0}, phantom, from=3-3, to=3-1, start anchor=center, end anchor=center]
	\arrow["{GD1_A^\mbbl}"{description}, from=5-3, to=3-3]
	\arrow[""{name=2, anchor=center, inner sep=0}, Rightarrow, no head, from=5-3, to=3-1]
	\arrow[""{name=2p, anchor=center, inner sep=0}, phantom, from=5-3, to=3-1, start anchor=center, end anchor=center]
	\arrow[""{name=3, anchor=center, inner sep=0}, "{G(D(l s_{GDA})^\mbbl \circ D1_A^\mbbl)}"', shift right=8, curve={height=30pt}, from=5-3, to=1-3]
	\arrow["\mbbl"', shorten <=9pt, shorten >=9pt, Rightarrow, from=0p, to=1p]
	\arrow["\mbbl", shorten <=4pt, shorten >=4pt, Rightarrow, from=1p, to=2p]
	\arrow["\gamma"{pos=0.3}, shorten >=4pt, Rightarrow, from=3-3, to=3]
\end{tikzcd}\]
Since this composite 2-cell is invertible and the 1-cell above right is (isomorphic) to a 1-cell of form $D h$ for a 1-cell $h$ in $\ck$, this is a $G$-lift by what we have proven above. For the same reasons, the unit and composition axioms in the assumptions of Theorem \ref{THM_Bunge} are satisfied.

We thus have a right colax adjoint to $G: \ck_D \to \cl$, let us denote it by $R: \cl \to \ck_D$. Since the pseudonaturality square of the counit $s$ is a $G$-lift (this again follows from what we have proven at the beginning of the proof), for any 1-cell $l: L \to K$ there exists a unique invertible 2-cell $\delta_l: Rl \Rightarrow J_D Hl$ making the following diagrams equal:
\[\begin{tikzcd}
	L && GDHL &&& L && GDHL \\
	&&&& {=} \\
	K && GDHK &&& K && GDHK
	\arrow[""{name=0, anchor=center, inner sep=0}, "GDHl", from=3-3, to=1-3]
	\arrow[""{name=1, anchor=center, inner sep=0}, "{s_L}"', from=1-3, to=1-1]
	\arrow[""{name=1p, anchor=center, inner sep=0}, phantom, from=1-3, to=1-1, start anchor=center, end anchor=center]
	\arrow["l", from=3-1, to=1-1]
	\arrow[""{name=2, anchor=center, inner sep=0}, "{s_{L}}", from=3-3, to=3-1]
	\arrow[""{name=2p, anchor=center, inner sep=0}, phantom, from=3-3, to=3-1, start anchor=center, end anchor=center]
	\arrow[""{name=3, anchor=center, inner sep=0}, "GRl", curve={height=-30pt}, from=1-3, to=3-3]
	\arrow["l", from=3-6, to=1-6]
	\arrow[""{name=4, anchor=center, inner sep=0}, "{s_K}", from=3-8, to=3-6]
	\arrow[""{name=5, anchor=center, inner sep=0}, "{s_L}"', from=1-8, to=1-6]
	\arrow["GRl"', from=3-8, to=1-8]
	\arrow["{s_l}", shorten <=9pt, shorten >=9pt, Rightarrow, from=1p, to=2p]
	\arrow["{G\delta_l}"', shorten <=6pt, shorten >=6pt, Rightarrow, from=3, to=0]
	\arrow["\mbbl", shorten <=9pt, shorten >=9pt, Rightarrow, from=5, to=4]
\end{tikzcd}\]
It is now routine to verify that this data gives an invertible icon $\delta: R \Rightarrow J_D H$ (which is an isomorphism in $\psd{\cl}{\ck_D}$), proving that the pseudofunctor $J_D H$ is right colax adjoint to $G: \ck_D \to \cl$ as well.
\end{proof}

\begin{pozn}\label{POZN_almost_a_biadj_big_l_adj_thm}
Since for our $G$-lifts, the accompanying 2-cell is always invertible, by (the dual of) Remark \ref{POZN_almost_a_biadj} we see that the colax adjunction we obtained in Theorem \ref{THM_BIG_lax_adj_thm}:
\[\begin{tikzcd}
	{(\Psi, \Phi): (s, \eta):} & {\ck_D} && \cl
	\arrow[""{name=0, anchor=center, inner sep=0}, "G"', curve={height=24pt}, from=1-2, to=1-4]
	\arrow[""{name=1, anchor=center, inner sep=0}, "{J_DH}"', curve={height=24pt}, from=1-4, to=1-2]
	\arrow["\ttop"{description}, draw=none, from=0, to=1]
\end{tikzcd}\]
is almost a biadjunction: $\Psi, \Phi$ are invertible and the counit $s$ is pseudonatural.
\end{pozn}

This next observation will be (only) needed in the proof of Theorem \ref{THM_Kleisli2cat_is_corefl_complete}:

\begin{pozn}\label{POZN_bigladjthm_hom_adjunkce}
By Proposition \ref{THMP_colax_adj_adjoint_homcats}, for each $L \in \cl$, $DA \in \ob \ck_D$ there is an adjunction between hom categories like this:
\[\begin{tikzcd}
	{\cl(GDA, L)} && {\ck_D(DA,J_DHL)}
	\arrow[""{name=0, anchor=center, inner sep=0}, curve={height=-24pt}, from=1-1, to=1-3]
	\arrow[""{name=1, anchor=center, inner sep=0}, "{(s_L)_* \circ G}", curve={height=-24pt}, from=1-3, to=1-1]
	\arrow["\dashv"{anchor=center, rotate=90}, draw=none, from=1, to=0]
\end{tikzcd}\]
Further, by the remark above and Remark \ref{POZN_colax_adj_adjoint_homcats_corefl_refl}, this is a \textbf{reflection} -- the counit is invertible because it is built out of $s$ and $\Psi$. The collection of these adjunctions is also pseudonatural in the variable $L \in \cl$.
\end{pozn}

\medskip

\noindent Our first application will be the following:
\begin{cor}\label{THMC_leftKanpsmon_rightcolaxadj_to_identity}
Given a left Kan pseudomonad $(D,y)$ on a 2-category $\ck$, the biadjunction between the base 2-category and the Kleisli 2-category induces a colax adjunction on the Kleisli 2-category:
\[\begin{tikzcd}
	\ck && {\ck_D} & \rightsquigarrow & {\ck_D} && {\ck_D}
	\arrow[""{name=0, anchor=center, inner sep=0}, "{J_D}"', curve={height=24pt}, from=1-1, to=1-3]
	\arrow[""{name=1, anchor=center, inner sep=0}, "{F_D}"', curve={height=24pt}, from=1-3, to=1-1]
	\arrow[""{name=2, anchor=center, inner sep=0}, curve={height=24pt}, Rightarrow, no head, from=1-5, to=1-7]
	\arrow[""{name=3, anchor=center, inner sep=0}, "{J_DF_D}"', curve={height=24pt}, from=1-7, to=1-5]
	\arrow["\dashv"{anchor=center, rotate=90}, draw=none, from=0, to=1]
	\arrow["\ttop"{description}, draw=none, from=2, to=3]
\end{tikzcd}\]
\end{cor}

In Remark \ref{POZN_rali_limits_not_unique_up_to_equiv} we gave an elementary example demonstrating that a colax adjoint to a pseudofunctor is not unique up to an equivalence. The above corollary gives a non-elementary example -- $J_D F_D$ is clearly not equivalent to the identity 2-functor on $\ck_D$.

The next corollary is a categorification of the fact that for an idempotent monad, the Kleisli and EM-categories are equivalent:

\begin{cor}\label{THM_lax_adjunkce_EM_Kleisli}
Given a left Kan pseudomonad $(D, y)$, the associated free-forgetful biadjunction induces a colax adjunction between the Kleisli 2-category and the 2-category of algebras:
\[\begin{tikzcd}
	\ck && {\text{Ps-}D\text{-Alg}} & \rightsquigarrow & {\ck_D} && {\text{Ps-}D\text{-Alg}}
	\arrow[""{name=0, anchor=center, inner sep=0}, "{F^D}"', curve={height=24pt}, from=1-1, to=1-3]
	\arrow[""{name=1, anchor=center, inner sep=0}, "{U^D}"', curve={height=24pt}, from=1-3, to=1-1]
	\arrow[""{name=2, anchor=center, inner sep=0}, curve={height=24pt}, hook, from=1-5, to=1-7]
	\arrow[""{name=3, anchor=center, inner sep=0}, "{J_D\circ U^D}"', curve={height=24pt}, from=1-7, to=1-5]
	\arrow["\dashv"{anchor=center, rotate=90}, draw=none, from=0, to=1]
	\arrow["\ttop"{description}, draw=none, from=2, to=3]
\end{tikzcd}\]
\end{cor}

\noindent The following is a change-of-base-style theorem:
\begin{cor}\label{THM_change_of_base_kleisli}
Let $D$ be a lax-idempotent pseudomonad on a 2-category $\ck$ and $T$ be a pseudomonad on a 2-category $\cl$. Assume that:

\begin{itemize}
\item there is a biadjunction between the base 2-categories, as portrayed below left:
\[\begin{tikzcd}
	&&&& \ck && \cl \\
	\ck && \cl \\
	&&&& {\ck_D} && {\cl_T}
	\arrow["L", from=1-5, to=1-7]
	\arrow["{J_D}"', from=1-5, to=3-5]
	\arrow["{J_T}", from=1-7, to=3-7]
	\arrow[""{name=0, anchor=center, inner sep=0}, "L"', curve={height=24pt}, from=2-1, to=2-3]
	\arrow[""{name=1, anchor=center, inner sep=0}, "R"', curve={height=24pt}, from=2-3, to=2-1]
	\arrow["{L^\#}"', dashed, from=3-5, to=3-7]
	\arrow["\dashv"{anchor=center, rotate=90}, draw=none, from=0, to=1]
\end{tikzcd}\]
\item the left biadjoint admits an extension to the Kleisli 2-categories, as portrayed above right.
\end{itemize}
Then there is an induced colax adjunction between the Kleisli 2-categories:
\[\begin{tikzcd}
	{\ck_D} && {\cl_T}
	\arrow[""{name=0, anchor=center, inner sep=0}, "{L^\#}"', curve={height=24pt}, from=1-1, to=1-3]
	\arrow[""{name=1, anchor=center, inner sep=0}, curve={height=24pt}, from=1-3, to=1-1]
	\arrow["\ttop"{description}, draw=none, from=0, to=1]
\end{tikzcd}\]
\end{cor}
\begin{proof}
Composing the Kleisli biadjunction with the $L \dashv R$ biadjunction we obtain the following biadjunction on which we can apply the theorem:
\[\begin{tikzcd}
	\ck && \cl && {\cl_T} \\
	\\
	&& {\ck_D}
	\arrow[""{name=0, anchor=center, inner sep=0}, "L"', curve={height=24pt}, from=1-1, to=1-3]
	\arrow[""{name=1, anchor=center, inner sep=0}, "R"', curve={height=24pt}, from=1-3, to=1-1]
	\arrow[""{name=2, anchor=center, inner sep=0}, "T"', curve={height=24pt}, from=1-3, to=1-5]
	\arrow[""{name=3, anchor=center, inner sep=0}, "{F_T}"', curve={height=24pt}, from=1-5, to=1-3]
	\arrow["{J_D}"', curve={height=18pt}, from=1-1, to=3-3]
	\arrow["{L^\#}"', curve={height=18pt}, from=3-3, to=1-5]
	\arrow["\dashv"{anchor=center, rotate=90}, draw=none, from=2, to=3]
	\arrow["\dashv"{anchor=center, rotate=90}, draw=none, from=0, to=1]
\end{tikzcd}\]
\end{proof}

\begin{pr}
Consider a 2-category $\ck$ with comma objects and pullbacks and denote by $D$ the fibration 2-monad on $\ck/C$ and $T$ the fibration 2-monad on $\ck/D$ (Example \ref{PR_fibration_2monad}). Any 1-cell $k: C \to D$ gives a 2-functor $k_*: \ck/C \to \ck/D$ with a right 2-adjoint $k^*$ given by pulling back. The 2-functor $k_*$ clearly extends to the colax slice 2-categories, hence giving rise to a \textbf{lax} adjunction between the colax slices:
\[\begin{tikzcd}
	{\ck//C} && {\ck//D}
	\arrow[""{name=0, anchor=center, inner sep=0}, "{k_*}"', curve={height=24pt}, from=1-1, to=1-3]
	\arrow[""{name=1, anchor=center, inner sep=0}, "{D \circ F_D \circ k^*}"', curve={height=24pt}, from=1-3, to=1-1]
	\arrow["\ttop"{description}, draw=none, from=0, to=1]
\end{tikzcd}\]
\end{pr}

\begin{pr}
In the next section (Corollary \ref{THM_COR_change_of_base_SalgTalg})  we will see how, when given a morphism of 2-monads $\theta: S \to T$, this gives rise to a colax adjunction between $\talgl$ and $S\text{-Alg}_l$.
\end{pr}

\begin{pozn}[Left Kan 2-monads]\label{POZN_leftkan2monads_big_l_adj}

Assume that $(D,y)$ is a left Kan 2-monad and that we have the same starting biadjunction as in Theorem \ref{THM_BIG_lax_adj_thm}, except now the modifications $\sigma, \tau$ are the identities and the counit $s$ is 2-natural. Going through the proof, note that $s_L \circ GJ_D(-): \ck(A,HL) \to \cl(GDA,L)$ is an isomorphism of categories: for each $l: GDA \to L$ there is a \textbf{unique} $l^\mbbl: A \to HL$ such that $s_L \circ GDl^\mbbl = l$. Because $J_D: \ck \to \ck_D$ is now a 2-functor and because of the uniqueness of each $l^\mbbl$, the collection $s_L: GDHL \to L$ is  strictly closed for $G$-lifts (dual of Definition \ref{DEFI_strictly_closed_for_Uexts}). By Remark \ref{POZN_strict_BungesTHM} we obtain a colax adjunction for which the modifications are the identities and the counit $s$ is 
2-natural.
\end{pozn}

\begin{pozn}[One-dimensional situation]
A lax-idempotent 2-monad $(T,m,i)$ on a locally discrete 2-category $\cc$ is precisely an idempotent 1-monad. The one-dimensional version of Theorem \ref{THM_BIG_lax_adj_thm} says that for idempotent monads, any adjunction as below left induces an adjunction as below right:
\[\begin{tikzcd}
	\cc && {\cc_T} && \cd & \rightsquigarrow & {\cc_T} && \cd
	\arrow["{J_D}"', from=1-1, to=1-3]
	\arrow["G"', from=1-3, to=1-5]
	\arrow[""{name=0, anchor=center, inner sep=0}, "H"', curve={height=30pt}, from=1-5, to=1-1]
	\arrow[""{name=1, anchor=center, inner sep=0}, "G"', curve={height=24pt}, from=1-7, to=1-9]
	\arrow[""{name=2, anchor=center, inner sep=0}, "{J_DH}"', curve={height=24pt}, from=1-9, to=1-7]
	\arrow["\dashv"{anchor=center, rotate=90}, draw=none, from=1-3, to=0]
	\arrow["\dashv"{anchor=center, rotate=90}, draw=none, from=1, to=2]
\end{tikzcd}\]
It can be seen that for instance Corollary \ref{THM_lax_adjunkce_EM_Kleisli} implies that for an idempotent monad, the Kleisli and EM-categories are equivalent.
\end{pozn}

\subsection{Weak limits in the Kleisli 2-category}\label{SEC_weak_limits_Kleisli2cat}

\begin{theo}\label{THM_Kleisli2cat_is_corefl_complete}
Assume $\ck$ is a 2-category that admits weighted $J$-indexed bilimits and assume $D$ is a lax-idempotent pseudomonad on $\ck$. Then, the Kleisli 2-category $\ck_D$ admits weighted $J$-indexed reflector-limits.
\end{theo}
\begin{proof}
Let $G: \cp \to \ck_D$ be a 2-functor, write $\widetilde{G}$ for the 2-functor given by the assignment $DA \mapsto \ck_D(DA,G?)$. Denoting again by $U_D$ the canonical pseudofunctor $\ck_D \to \ck$, there is a biadjunction where the right biadjoint sends a weight $W$ to the bilimit of $U_D G: \cp \to \ck$ weighted by $W$:
\[\begin{tikzcd}
	\ck && {\ck_D} && {\psd{\cp}{\cat}^{op}}
	\arrow["{J_D}"', hook, from=1-1, to=1-3]
	\arrow["{\widetilde{G}}"', from=1-3, to=1-5]
	\arrow[""{name=0, anchor=center, inner sep=0}, "{\{-,U_DG \}}"', curve={height=30pt}, from=1-5, to=1-1]
	\arrow["\top"{description}, draw=none, from=1-3, to=0]
\end{tikzcd}\]
This is because of the following equivalences that are pseudonatural with respect to weights $W \in [\cp, \cat]$ and objects $A \in \ck$:
\begin{align*}
\ck(A, \{ W, U_D G \}) &\simeq \text{Psd}[\cp,\cat](W, \ck(A,U_D G-))\\
&\simeq \text{Psd}[\cp,\cat](W, \ck_D(DA,G-))\\
&= \text{Psd}[\cp,\cat]^{op}(\ck_D(DA,G-),W).
\end{align*}

\noindent Transferring the identity across these equivalences, we see that the \textbf{counit} of the biadjunction is the composite pseudonatural transformation:
\[\begin{tikzcd}
	W && {\ck(\{ W, U_DG \}, U_D G-)} && {\ck_D(D \{ W, U_D G \}, G- ),}
	\arrow["{\lambda_W}", Rightarrow, from=1-1, to=1-3]
	\arrow["{(p_{G-})_* \circ J_D}", Rightarrow, from=1-3, to=1-5]
\end{tikzcd}\]

\noindent where $\lambda_W$ is the $W$-weighted bilimit cone for $U_DG : \cp \to \ck$. By Theorem \ref{THM_BIG_lax_adj_thm} this induces a colax adjunction with the same counit:
\[\begin{tikzcd}
	{\ck_D} && {\text{Psd}[\cp,\cat]^{op}}
	\arrow[""{name=0, anchor=center, inner sep=0}, "{\widetilde{G}}"', from=1-1, to=1-3]
	\arrow[""{name=0p, anchor=center, inner sep=0}, phantom, from=1-1, to=1-3, start anchor=center, end anchor=center]
	\arrow[""{name=1, anchor=center, inner sep=0}, "{J_D\{-,U_DG \}}"', curve={height=24pt}, from=1-3, to=1-1]
	\arrow[""{name=1p, anchor=center, inner sep=0}, phantom, from=1-3, to=1-1, start anchor=center, end anchor=center, curve={height=24pt}]
	\arrow["\ttop"{description}, draw=none, from=0p, to=1p]
\end{tikzcd}\]

\noindent By Remark \ref{POZN_bigladjthm_hom_adjunkce} there is a reflection between hom categories as pictured below:

\[\begin{tikzcd}[column sep=scriptsize]
	{\ck_D(B,J_D \{ W, U_DG\})\phantom{...............}} && {\phantom{.................}\text{Psd}[\cp,\cat](W,\ck_D(B,G-))}
	\arrow[""{name=0, anchor=center, inner sep=0}, "{\theta \mapsto \ck_D(\theta,G-) \circ (p_{G-})_* \circ J_D \circ \lambda_W}"', curve={height=24pt}, from=1-1, to=1-3]
	\arrow[""{name=1, anchor=center, inner sep=0}, curve={height=24pt}, from=1-3, to=1-1]
	\arrow["\dashv"{anchor=center, rotate=90}, draw=none, from=0, to=1]
\end{tikzcd}\]
This exhibits the pseudonatural transformation $(p_{G-})_* \circ J_D \circ \lambda_W$ as the reflector-limit of $G$ weighted by $W$.
\end{proof}

\begin{pr}\label{PR_colax_slice_2cat_weaklimit_2}
The proof of Theorem \ref{THM_Kleisli2cat_is_corefl_complete} gives a concrete way to compute limits. Consider the colax slice 2-category from Example \ref{PR_fibration_2monad_Kleisli2cat}, which is the Kleisli 2-category for the (colax-idempotent) fibration 2-monad from Example \ref{PR_fibration_2monad}. We can see that the process of computing \textbf{rali}-product of two objects $(f_1: A_1 \to C, f_2: A_2 \to C)$ in the colax-slice 2-category $\ck//C$ agrees with the process described in Example \ref{PR_colax_slice_2cat_weaklimit_1}.
\end{pr}

\begin{pozn}\label{POZN_strict_weaklimits}
Given a left Kan 2-monad $(D,y)$, going through the proof of Theorem \ref{THM_Kleisli2cat_is_corefl_complete} (and considering Remark \ref{POZN_leftkan2monads_big_l_adj}) we may now replace $\psd{\cp}{\cat}$ by $[\cp, \cat]$ and the result can be changed to the claim that $\ck_D$ admits $J$-indexed lali-limit whenever $\ck$ admits them as 2-limits.
\end{pozn}

\section{Applications to two-dimensional monad theory}\label{SEKCE_applications_to_twodim}

In Subsection \ref{SUBS_4_lax_flexibility} we will prove a theorem characterizing $T$-algebras for a 2-monad that are \textit{lax-flexible} -- a lax analogue of the classical notion of flexibility.

In Section \ref{SUBS_4_on_various_colax_adjs} we prove the existence of various colax adjunctions involving the 2-category $\talgl$ of $T$-algebras and lax morphisms for a 2-monad $T$. The results proven here give the lax versions of results in \cite[Section 5]{twodim}.

\begin{defi}\label{DEFI_property_L}
Let $T$ be a 2-monad on a 2-category $\ck$. We will say that it satisfies \textbf{Property L} if the inclusion $\talgs \xhookrightarrow{} \talgl$ admits a left 2-adjoint and the corresponding lax-morphism classifier 2-comonad $Q_l$ on $\talgs$ (see \ref{DEFI_Ql_Qp}) is lax-idempotent.
\end{defi}

\begin{prc}
Given a category $\ce$ with pullbacks, every 2-monad $T$ of form $\cat(T)$ satisfies \textbf{Property L}. We have proven this in Theorem \ref{THMP_catT_strictification_exists} and Proposition \ref{THMP_CatT_Ql_lax_idemp}.
\end{prc}

By Theorem \ref{THM_adjunkce_talgs_talgl_colaxtalgl}, Proposition \ref{THMP_Ql_lax_idemp_Qp_ps_idemp}, a 2-monad $T$ on $\ck$ will have this property when $\ck$ admits oplax limits of an arrow and $\talgs$ is sufficiently cocomplete (admits codescent objects). This is the case for instance if $\ck$ is complete and cocomplete and $T$ has a \textit{rank} (preserves $\alpha$-filtered colimits), see \cite[Theorem 3.8]{twodim}.

To apply the (appropriate dual of the) results developed in Section \ref{SEKCE_Kleisli2cat} to the lax-idempotent 2-comonad $Q_l$, notice that (with the hint of Remark \ref{POZN_duals}) that this means dualizing all the results from $\ck$ to $\ck^{coop}$.

For instance ``coreflection-inclusion” gets replaced by ``reflector”.

\subsection{Lax flexibility}\label{SUBS_4_lax_flexibility}

For this section, recall the notions of \textit{semiflexible} and \textit{flexible} algebras for a 2-monad $T$ from \cite[Remark 4.5, page 22]{twodim}. By \cite[Proposition 1]{onsemiflexible}, a $T$-algebra $(A,a)$ is semi-flexible if and only if it admits the structure of a pseudo $Q_p$-coalgebra. A \textit{pie} $T$-algebra was then defined to be a $T$-algebra that admits a strict $Q_p$-coalgebra structure. This motivates us to define:

\begin{defi}\label{DEFI_lax_semiflexible_flexible_pie}
Let $T$ be a 2-monad on a 2-category $\ck$ that satisfies \textbf{Property L}. A $T$-algebra $(A,a)$ is said to be:
\begin{itemize}
\item \textit{lax-semiflexible} if it admits a pseudo $Q_l$-coalgebra structure.
\item \textit{lax-flexible} if it admits a normal pseudo $Q_l$-coalgebra structure.
\item \textit{lax-pie} if it admits a strict $Q_l$-coalgebra structure.
\end{itemize}
\end{defi}

\begin{pozn}\label{POZN_lax_semiflexible_is_semiflexible}
Every lax-\textbf{Y} $T$-algebra is \textbf{Y}, where $\textbf{Y} \in \{ \text{flexible, semiflexible, pie} \}$. This is because of the fact that by Proposition \ref{THMP_mor_of_comonads_QlQp} there is an induced 2-functor from pseudo $Q_l$-coalgebras to pseudo $Q_p$-coalgebras that commutes with the 2-functors that forget the coalgebra structure (and thus keeps the $T$-algebra structure intact). 
\end{pozn}

\begin{pr}
In Corollary \ref{THM_free_forgetful_colax_adj} we will see that every free $T$-algebra is lax-flexible; this is a strengthening of the fact that every free $T$-algebra is flexible (\cite[Corollary 5.6]{twodim}).
\end{pr}

\begin{pr}
Fix a 1-category $\cj$ and consider the reslan 2-monad $T$ on $[ \ob \cj, \cat]$ whose algebras are the weights (2-functors) $\cj \to \cat$ (Example \ref{PR_reslan2monad}). Weights that index lax limits are lax-pie: they are precisely the free $Q_l$-coalgebras, i.e. those of the form $Q_l W$ (see also Remark \ref{POZN_2limit_is_not_laxlimit}). Since a lax limit is in general not a pseudo-limit \cite[Remark 5.5]{twodim}, not every pie algebra is lax-pie.
\end{pr}

Following Example \ref{PR_doctrinal_adj} and Remark \ref{POZN_strict_charakterizace_D_algeber}, the application of (the dual of) Proposition \ref{THM_characterization_D-algs_D-corefl-incl} to the lax-idempotent 2-comonad $Q_l$ provides a lax version of \cite[Theorem 20 a)]{onsemiflexible}. It reads as:

\begin{theo}\label{THM_charakterizace_lax_semiflexible_algs}
Let $T$ be a 2-monad on a 2-category $\ck$ satisfying \textbf{Property L} and denote by $U: \talgs \to \ck$ the forgetful 2-functor. A $T$-algebra is lax-semiflexible, semiflexible, if and only if, respectively:
\begin{itemize}
\item $\talgs((A,a),-): \talgs \to \cat$ sends $U$-reflectors to reflectors in $\cat$.
\item $\talgs((A,a),-): \talgs \to \cat$ sends $U$-lalis to lalis in $\cat$.
\end{itemize}
\end{theo}

\noindent Lax-pie $T$-algebras will be briefly studied in Chapter \ref{CH_coKZfication} -- Section \ref{SEC_lax_pie_algebras}.

\subsection{Colax adjunctions and rali-cocompleteness of lax morphisms}\label{SUBS_4_on_various_colax_adjs}

Considering Remark \ref{POZN_leftkan2monads_big_l_adj}, the application of Theorem \ref{THM_bigadjthm_HL_is_D-algebra} and Theorem \ref{THM_BIG_lax_adj_thm} for the 2-comonad $Q_l$ reads as:

\begin{theo}\label{THM_big_colax_adj_THM}
Let $T$ be a 2-monad satisfying \textbf{Property L}. Any 2-adjunction below left induces a colax adjunction below right:
\[\begin{tikzcd}
	\talgs && \talgl && \cl & \rightsquigarrow & \talgl && \cl
	\arrow["J"', from=1-1, to=1-3]
	\arrow["G"', from=1-3, to=1-5]
	\arrow[""{name=0, anchor=center, inner sep=0}, "H"', curve={height=30pt}, from=1-5, to=1-1]
	\arrow[""{name=1, anchor=center, inner sep=0}, "JH"', from=1-7, to=1-9]
	\arrow[""{name=2, anchor=center, inner sep=0}, "G"', curve={height=24pt}, from=1-9, to=1-7]
	\arrow["\dashv"{anchor=center, rotate=-90}, draw=none, from=0, to=1-3]
	\arrow["\tbot"{description}, draw=none, from=1, to=2]
\end{tikzcd}\]
Moreover, for every $L \in \cl$, the $T$-algebra $HL$ is lax-flexible.
\end{theo}

\begin{cor}\label{THM_free_forgetful_colax_adj}
The free-forgetful adjunction for a 2-monad $T$ on a 2-category $\ck$ satisfying \textbf{Property L} induces a colax adjunction between $\talgl$ and $\ck$. In particular, every free $T$-algebra is lax-flexible.
\[\begin{tikzcd}[column sep=scriptsize]
	\talgs && \talgl && \ck & \rightsquigarrow & \talgl && \ck
	\arrow[""{name=0, anchor=center, inner sep=0}, "{F^T}"', curve={height=30pt}, from=1-5, to=1-1]
	\arrow[""{name=1, anchor=center, inner sep=0}, "U"', curve={height=24pt}, from=1-7, to=1-9]
	\arrow[""{name=2, anchor=center, inner sep=0}, "{J F^T}"', curve={height=24pt}, from=1-9, to=1-7]
	\arrow["J"', from=1-1, to=1-3]
	\arrow["U"', from=1-3, to=1-5]
	\arrow["\tbot"{description}, draw=none, from=1, to=2]
	\arrow["\dashv"{anchor=center, rotate=-90}, draw=none, from=0, to=1-3]
\end{tikzcd}\]
\end{cor}

Considering Remark \ref{POZN_lax_semiflexible_is_semiflexible}, the claim that free algebras are lax-flexible is a strengthening of the claim that they are flexible, which has been proven in \cite[Corollary 5.6]{twodim}.

\begin{cor}\label{THM_COR_change_of_base_SalgTalg}
Let $T,S$ be two 2-monads on a 2-category $\ck$ satisfying \textbf{Property L} and let $\theta: S \to T$ be a 2-monad morphism\footnote{A monad morphism (as in \textbf{Variations} for Definition \ref{DEFI_comonadfun}) for the case $\ck = \twoCAT$}. Assume that the induced 2-functor $\theta^*: \talgs \to S\text{-Alg}_s$ admits a left 2-adjoint $\theta_*$ (this is the case when $\ck$ is complete and cocomplete and $T$ is finitary, see \cite[Theorem 3.9]{twodim}). Then there is an induced colax adjunction between $\talgl$ and $\salgl$:
\[\begin{tikzcd}[column sep=scriptsize]
	\talgs && {S\text{-Alg}_s} && {S\text{-Alg}_l} & \rightsquigarrow & \talgl && {S\text{-Alg}_l} \\
	\\
	&& \talgl
	\arrow[""{name=0, anchor=center, inner sep=0}, "{\theta^*}"', curve={height=24pt}, from=1-1, to=1-3]
	\arrow[""{name=1, anchor=center, inner sep=0}, "{\theta_*}"', curve={height=24pt}, from=1-3, to=1-1]
	\arrow[""{name=2, anchor=center, inner sep=0}, "J"', curve={height=24pt}, from=1-3, to=1-5]
	\arrow[""{name=3, anchor=center, inner sep=0}, "{(-)'}"', curve={height=24pt}, from=1-5, to=1-3]
	\arrow["J"', curve={height=18pt}, from=1-1, to=3-3]
	\arrow["{\theta^*}"', curve={height=18pt}, from=3-3, to=1-5]
	\arrow[""{name=4, anchor=center, inner sep=0}, "{\theta^*}"', curve={height=24pt}, from=1-7, to=1-9]
	\arrow[""{name=5, anchor=center, inner sep=0}, curve={height=24pt}, from=1-9, to=1-7]
	\arrow["\dashv"{anchor=center, rotate=90}, draw=none, from=2, to=3]
	\arrow["\dashv"{anchor=center, rotate=90}, draw=none, from=0, to=1]
	\arrow["\ttop"{description}, draw=none, from=4, to=5]
\end{tikzcd}\]
\end{cor}

\begin{theo}\label{THM_talgl_reflector_cocomplete}
Let $T$ be a 2-monad on a 2-category $\ck$ that admits oplax limits of an arrow. Assume that $\talgs$ is cocomplete (in particular $T$ satisfies \textbf{Property L}). Then $\talgl$ is rali-cocomplete.
\end{theo}
\begin{proof}
This follows from (the dual of) Remark \ref{POZN_strict_weaklimits}.
\end{proof}

\begin{pozn}
By Remark \ref{POZN_ordinary_weakness}, this in particular shows that $\talgl$ is weakly cocomplete.
\end{pozn}

\begin{cor}\label{THMC_examples_of_rali_cocomplete_2cats}
The following 2-categories are rali-cocomplete:
\begin{enumerate}
\item for a 2-category $\cj$, the 2-category $\lax{\cj}{\cat}$ of 2-functors $\cj \to \cat$, lax natural transformations and modifications,
\item the 2-category of monoidal categories and lax-monoidal functors and its symmetric/braided variants,
\item the 2-category of small 2-categories, lax functors, and icons,
\item for a set $\Phi$ of small categories, the 2-category $\Phi\text{-Colim}_l$ of small categories that admit a choice of $J$-indexed colimits for $J \in \Phi$ and \textbf{all} functors between them. 
\end{enumerate}
\end{cor}
\begin{proof}
Each of these is a 2-category of form $\talgl$, where $T$ is a 2-monad with a rank on a complete and cocomplete 2-category (and so satisfies \textbf{Property L}). In the list above, $T$ is one of the following:
\begin{enumerate}
\item the reslan 2-monad from Example \ref{PR_reslan2monad},
\item the 2-monad for (nonstrict) monoidal categories, see \cite[5.5]{companion},
\item the 2-category 2-monad from Example \ref{PR_2category_2monad},
\item the 2-monad $T$ described in \cite[Theorem 6.1]{onthemonadicity} whose strict $T$-morphisms are functors that preserve the choices of $\Phi$-colimits. Lax $T$-morphisms are all functors because this 2-monad is lax-idempotent by \cite[Theorem 6.3]{onthemonadicity}.
\end{enumerate}
\end{proof}

\begin{pozn}
There is also a dual version for the 2-category $\talgc$ of $T$-algebras and \textbf{colax} $T$-morphisms. If $\talgs$ is sufficiently cocomplete, there is the colax morphism classifier 2-comonad and if $\ck$ admits lax limits of arrows, $Q_c$ is colax-idempotent (duals of Theorem \ref{THM_adjunkce_talgs_talgl_colaxtalgl} and Proposition \ref{THMP_Ql_lax_idemp_Qp_ps_idemp}). If $\talgs$ is 2-cocomplete, $\talgc$ can be seen to be lali-cocomplete.
\end{pozn}

%% file: PhD_CH_6_COKZFICATION.tex
\chapter{Further topics on lax-idempotent 2-monads}\label{CH_coKZfication}

\noindent \textbf{This final research chapter consists of two independent sections}:
\begin{itemize}
\item In Section \ref{SEC_coKZfication} we describe a process that can (under suitable condition) turn a 2-monad into a colax-idempotent one, categorifying the ``idempotentiation” process of Fakir \cite{fakir_HIST}. We apply this to show that the lax morphism classifier 2-comonad $Q_l$ (of \ref{DEFI_Ql_Qp}) enjoys a certain universal property.
\item In Section \ref{SEC_lax_pie_algebras} we prove that given any cartesian monad $T$ on a category $\ce$ with pullbacks, $T$-multicategories are comonadic over $\cat(T)$-algebras. We then use this result to show that \textit{lax-pie} algebras (lax analogue of pie algebras of \cite{onsemiflexible}) for the 2-monad $\cat(T)$ are equivalent to $T$-multicategories.
\end{itemize}

\section{Turning a 2-monad into a colax-idempotent one}\label{SEC_coKZfication}
\subsection{Introduction}

In \cite{fakir_HIST}, Sabah Fakir observed that any monad $T$ on a complete, well-powered category $\cc$ can be turned into an idempotent monad $T'$ on $\cc$, and this process is in fact a coreflection to the inclusion $\text{IdempMonad}(\cc) \xhookrightarrow{} \text{Monad}(\cc)$ of idempotent monads on $\cc$ to monads on $\cc$. The proof uses a transfinite sequence of monads, each built from the previous one using equalizers, with the monad $T'$ being given by the limit (which exists by well-poweredness).

A higher-dimensional analogue of this kind, namely, a result about turning a 2-monad $T$ into a lax-idempotent one $T'$ on the same 2-category, has appeared in the literature -- in \cite[Theorem 8.7]{onpropertylike} this has been proven under the assumption that $\ck$ is a locally finitely presentable 2-category. In that case, the resulting lax-idempotent 2-monad is moreover finitary.

In this chapter we give a small step towards a more direct categorification of Fakir's result -- we describe a process that takes a 2-monad $T$ on a 2-category $\ck$ with coreflexive descent objects and produces a new 2-monad $T'$ on $\ck$. We then show that if $T$ preserves coreflexive descent objects, this process terminates after 1 step and yields a colax-idempotent 2-monad. Since the upcoming sections are rather technical, let us provide:

\medskip

\noindent \textbf{Overview of \ref{SEC_coKZfication}}:

\begin{itemize}
\item Subsection \ref{SUBS_coKZ_process} describes a process that takes a 2-monad $(T,m,i)$ on a 2-category $\ck$ with coreflexive descent objects and produces a new 2-monad $(D,\mu,\eta)$. In particular:
\begin{itemize}
\item the endo-2-functor $D$ is obtained from Lemma \ref{THML_coKZ_D_is_2fun},
\item the unit $\eta$ is obtained from Lemma \ref{THML_coKZ_cone_for_eta},
\item the multiplication $\mu$ is obtained from Lemma \ref{THML_coKZ_D_is_2fun},
\item the 2-monad axioms are proven in Lemma \ref{THML_e_is_2monad_morphism},
\item there is a 2-monad morphism $e: D \to T$ again by Lemma \ref{THML_e_is_2monad_morphism}.
\end{itemize}

\item Subsection \ref{SUBS_coKZ_one_step_termination} shows that if $T$ preserves coreflexive descent objects, the 2-monad $D$ is colax-idempotent. We proceed as follows:
\begin{itemize}
\item an auxiliary modification $\delta: e D \circ \eta D \to e D \circ D \eta$ is obtained from Lemma \ref{THML_coKZ_conemor_for_delta},
\item a modification $\lambda: \eta S \to S \eta$ is obtained from Lemma \ref{THML_coKZ_conemor_for_lambda},
\item the colax-idempotent 2-monad axiom is proven in Theorem \ref{THM_one_step_kzfication}.
\end{itemize}

\item In Subsection \ref{SUBS_coKZ_coreflection} we show that for $T$ preserving coreflexive descent objects, the 2-monad morphism $e: D \to T$ is the coreflection of $T$ into colax-idempotent 2-monads (Theorem \ref{THM_coKZfication}).

\item Finally, in Subsection \ref{SEC_coKZfication_further_props} we focus on what the construction does to some relevant classes of examples (colax-idempotent 2-monads and 2-comonads on $\talgs$ generated by the free-forgetful adjunction).
\end{itemize}

\subsection{The process}\label{SUBS_coKZ_process}

Let $(T,m,i)$ be a 2-monad on a 2-category $\ck$. Each object $X \in \ck$ determines a strict \textit{op-coherence data} (coherence data in $\ck^{op}$) that we will refer to as its \textit{resolution}:
\[\begin{tikzcd}
	{\opres{X}:=} & TX && {T^2X} && {T^3X}
	\arrow["{i_{TX}}", shift left=5, from=1-2, to=1-4]
	\arrow["{Ti_X}"', shift right=5, from=1-2, to=1-4]
	\arrow["{Ti_{TX}}"{description}, from=1-4, to=1-6]
	\arrow["{i_{T^2X}}", shift left=5, from=1-4, to=1-6]
	\arrow["{T^2i_X}"', shift right=5, from=1-4, to=1-6]
	\arrow["{m_X}"{description}, from=1-4, to=1-2]
\end{tikzcd}\]

\begin{pozn}
$\opres{X}$ is a \textit{coreflexive} op-coherence data (the dual of the \textbf{Variation} in Definition \ref{DEFI_delta_l}) -- it automatically comes equipped with additional retractions:
\[
Tm_X, m_{TX}: T^3X \to T^2X.
\]
As such, it behaves in a similar manner as to how coreflexive equalizers do in ordinary category theory.
\end{pozn}

For the rest of this chapter, assume that $\ck$ admits descent objects of strict coreflexive op-coherence data. For $X \in \ob \ck$, denote by $(e_X: DX \to TX, \xi_X)$ the descent object of $\opres{X}$:
\[\begin{tikzcd}
	&& TX \\
	DX &&&& {T^2X} \\
	&& TX
	\arrow["{i_{TX}}", from=1-3, to=2-5]
	\arrow["{\xi_X}", shorten <=10pt, shorten >=6pt, Rightarrow, from=1-3, to=3-3]
	\arrow["{e_X}", from=2-1, to=1-3]
	\arrow["{e_X}"', from=2-1, to=3-3]
	\arrow["{Ti_X}"', from=3-3, to=2-5]
\end{tikzcd}\]

As such, it satisfies the following cone axioms:
\begin{align}
m_X \xi_X &= 1, \label{eqconeaxiom1}\\
Ti_{TX} \xi_X &= T^2i_X \xi_X \circ i_{T^2X} \xi_X.\label{eqconeaxiom2}
\end{align}

\noindent The statements below follow from the standard fact that limits in functor categories are pointwise, so we omit the proofs:

\begin{lemma}\label{THML_coKZ_D_is_2fun}
The assignment $X \mapsto DX$ gives a 2-functor $D:\ck \to \ck$. Moreover:
\begin{itemize}
\item The components $e_X$ give a 2-natural transformation $e: D \Rightarrow T$,
\item The components $\xi_X$ give a modification $\xi: iT \circ e \to Ti \circ e$.
\end{itemize}
\end{lemma}

\begin{lemma}\label{THML_cone_in_KK}
Given a 2-functor $S: \ck \to \ck$, 2-natural transformation $g: S \Rightarrow T$ and a modification:
\[\begin{tikzcd}
	&& T \\
	S &&&& {T^2} \\
	&& T
	\arrow["iT", Rightarrow, from=1-3, to=2-5]
	\arrow["\psi", shorten <=10pt, shorten >=6pt, from=1-3, to=3-3]
	\arrow["g", Rightarrow, from=2-1, to=1-3]
	\arrow["g"', Rightarrow, from=2-1, to=3-3]
	\arrow["Ti"', Rightarrow, from=3-3, to=2-5]
\end{tikzcd}\]
that has the property that for each $A \in \ck$, $(g_A, \psi_A)$ is a cone for $\opres{A}$, then the collection $\theta_A: SA \to DA$ of comparison maps from the cone $(g_A,\psi_A)$ to $(e_A, \xi_A)$ assembles into a 2-natural transformation $\theta: S \Rightarrow D$.

Moreover, given another such pair $(h: S \Rightarrow T, \phi: iT \circ h \to Ti \circ h)$ and a modification $\overline{l}: h \to g$ such that for each $A$, $\overline{l}_A: h_A \Rightarrow g_A$ is a morphism $(h_A, \phi_A) \to (g_A, \psi_A)$ of cones, the collection of 2-cells $\overline{\theta}_A$, each unique with the property that:
\[\begin{tikzcd}
	SA && TA & {=} & SA && DA && TA
	\arrow[""{name=0, anchor=center, inner sep=0}, "{h_A}", curve={height=-24pt}, from=1-1, to=1-3]
	\arrow[""{name=1, anchor=center, inner sep=0}, "{g_A}"', curve={height=24pt}, from=1-1, to=1-3]
	\arrow[""{name=2, anchor=center, inner sep=0}, "{\theta_A}", curve={height=-24pt}, from=1-5, to=1-7]
	\arrow[""{name=3, anchor=center, inner sep=0}, "{\theta'_A}"', curve={height=24pt}, from=1-5, to=1-7]
	\arrow["{e_A}", from=1-7, to=1-9]
	\arrow["{\overline{l_A}}", shorten <=6pt, shorten >=6pt, Rightarrow, from=0, to=1]
	\arrow["{\overline{\theta}_A}", shorten <=6pt, shorten >=6pt, Rightarrow, from=2, to=3]
\end{tikzcd}\]
Assembles into a modification:
\[
\overline{\theta}: \theta \to \theta' : S \Rightarrow D : \ck \to \ck.
\]
\end{lemma}

\begin{pozn}
The previous lemma says that $(e,\xi)$ is the descent object of the diagram $\opres{-}$ in the $\CAT$-category $[\ck, \ck]$.
\end{pozn}

\noindent The following is obvious:

\begin{lemma}\label{THML_coKZ_cone_for_eta}
The pair $(i_X, 1_{Ti_X i_X})$ is a cone for $\opres{X}$.
\end{lemma}

There thus exists a unique map $\eta_X : X \to DX$ such that:
\begin{align}
e_X \eta_X &= i_X,\label{eqeta} \\
\xi_X \eta_X &= 1_{Ti_X i_X}.\label{eqeta2}
\end{align}

By Lemma \ref{THML_cone_in_KK}, this gives a 2-natural transformation $\eta: 1_{\ck} \Rightarrow D$.

\begin{lemma}\label{THML_coKZ_cone_for_mu}
The pair $(m_X \circ Te_X \circ e_{DX}, \psi_X)$ (with the 2-cell $\psi_X$ portrayed below) is a cone for $\opres{DX}$:
\[\begin{tikzcd}
	&&& {T^2X} && TX \\
	&& DTX && {T^3X} \\
	{D^2X} &&& {T^2X} &&& {T^2X} \\
	&& TDX && {T^3X} \\
	& DTX && {T^2X} && TX
	\arrow["{m_X}", from=1-4, to=1-6]
	\arrow["{i_{T^2X}}", from=1-4, to=2-5]
	\arrow["{\xi_{TX}}"', shorten <=6pt, shorten >=6pt, Rightarrow, from=1-4, to=3-4]
	\arrow["{i_{TX}}", from=1-6, to=3-7]
	\arrow["{e_{TX}}", from=2-3, to=1-4]
	\arrow["{e_{TX}}"', from=2-3, to=3-4]
	\arrow["{Tm_X}"{description}, from=2-5, to=3-7]
	\arrow["{De_X}", from=3-1, to=2-3]
	\arrow["{e_{DX}}"', from=3-1, to=4-3]
	\arrow["{De_X}"', from=3-1, to=5-2]
	\arrow["{Ti_{TX}}"', from=3-4, to=2-5]
	\arrow[Rightarrow, no head, from=3-4, to=3-7]
	\arrow["{Ti_{TX}}", from=3-4, to=4-5]
	\arrow["{T\xi_X}"', shorten <=6pt, shorten >=6pt, Rightarrow, from=3-4, to=5-4]
	\arrow["{Te_X}", from=4-3, to=3-4]
	\arrow["{Te_X}"', from=4-3, to=5-4]
	\arrow["{m_{TX}}"{description}, from=4-5, to=3-7]
	\arrow["{e_{TX}}"', from=5-2, to=5-4]
	\arrow["{T^2i_X}"', from=5-4, to=4-5]
	\arrow["{m_X}"', from=5-4, to=5-6]
	\arrow["{Ti_X}"', from=5-6, to=3-7]
\end{tikzcd}\]

\end{lemma}
\begin{proof}
The cone axiom \eqref{eqconeaxiom1} is immediate from the naturality of $m: T^2 \Rightarrow T$. The cone axiom \eqref{eqconeaxiom2} follows since\footnote{In the equations below, we will \underline{underline} the terms that are about to change in the next step.}:
\begin{align*}
T^2i_X \psi_X \circ i_{T^2X} \psi_X \sstackrel{(A)}{=}\\[10pt]
T^2i_X m_{TX} T\xi_X e_{DX} \circ T^2i_XTm_X T^2e_X \xi_{DX} \circ \\ \circ \underline{i_{T^2X}m_{TX}T\xi_X} e_{DX} \circ i_{T^2X}Tm_X T^2e_X \xi_{DX} =\\[10pt]
T^2i_X m_{TX} T\xi_X e_{DX} \circ \underline{T^2i_X Tm_X T^2e_X \xi_{DX} \circ} \\ \underline{\circ Tm_{TX} T^2\xi_X i_{TDX} e_{DX}} \circ i_{T^2X} Tm_X T^2e_X \xi_{DX} \sstackrel{(*)}{=}\\[10pt]
\underline{T^2i_X m_{TX}} T\xi_X e_{DX} \circ Tm_{TX}\underline{T^2\xi_X Ti_{DX}} e_{DX} \circ \\ \circ \underline{Tm_{TX} T^2i_{TX}} T^2e_X \xi_{DX} \circ \underline{i_{T^2X}Tm_X} T^2e_X \xi_{DX} =\\[10pt]
m_{T^2X}T^3i_X T\xi_X e_{DX} \circ \underline{Tm_{TX} Ti_{T^2X}} T\xi_X e_{DX} \circ \\ \circ T^2m_X \underline{T^2i_{TX} T^2e_X} \xi_{DX} \circ T^2m_X T^3e_X i_{T^2DX} \xi_{DX} =\\[10pt]
\underline{m_{T^2X} T^3i_X T\xi_X e_{DX} \circ m_{T^2X} Ti_{T^2X} T\xi_X e_{DX}} \circ \\ \circ \underline{T^2m_X T^3e_X T^2i_{DX} \xi_{DX} \circ T^2m_X T^3e_X i_{T^2DX} \xi_{DX}} \sstackrel{(C)}{=}\\[10pt]
\underline{m_{T^2X} T^2i_{TX}} T\xi_X e_{DX} \circ T^2m_X T^3e_X Ti_{TDX} \xi_{DX} =\\
Ti_{TX} m_{TX} T\xi_X e_{DX} \circ Ti_{TX} Tm_X T^2e_X \xi_{DX} =\\
Ti_{TX} \psi_X.
\end{align*}
here $(A)$ uses the modification axiom $\xi_{TX} De_X = T^2e_X \xi_{DX}$, $(*)$ is the middle-four interchange rule for 2-cells, $(C)$ is the cone axiom \eqref{eqconeaxiom2} for $(e_X, \xi_X)$.
\end{proof}

\noindent Thus there exists a unique 1-cell $\mu_X: D^2X \to DX$ such that:
\begin{align}
e_X \mu_X &= m_X Te_X e_{DX},\label{eqmu1}\\
\xi_X \mu_X &= m_{TX} T\xi_X e_{DX} \circ Tm_X T^2e_X \xi_{DX}.\label{eqmu2}
\end{align}

By Lemma \ref{THML_cone_in_KK}, $\mu: D^2 \Rightarrow D$ is a 2-natural transformation.

\begin{lemma}\label{THML_e_is_2monad_morphism}
$(D,\mu,\eta)$ is a 2-monad and the 2-natural transformation $e: D \Rightarrow T$ is a 2-monad morphism\footnote{A monad morphism (as in \textbf{Variations} for Definition \ref{DEFI_comonadfun}) for the case $\ck = \twoCAT$} $(D,\mu,\eta) \to (T,m,i)$.
\end{lemma}
\begin{proof}
\textbf{The first monad unit axiom} $\mu_X \eta_{DX} = 1_{DX}$: By the one-dimensional universal property of the descent object, this will follow if we show that:
\begin{align*}
e_X \mu_X \eta_{DX} &= e_X,\\
\xi_X \mu_X \eta_{DX} &= \xi_X.
\end{align*}
The first equation holds because:
\[
e_X \mu_X \eta_{DX} \sstackrel{\eqref{eqmu1}}{=} m_X Te_X e_{DX} \eta_{DX} \sstackrel{\eqref{eqeta}}{=} m_X Te_X i_{DX} = m_X i_{TX} e_X = e_X.
\]
The second equation holds because:
\begin{align*}
\xi_X \mu_X \eta_{DX} &\sstackrel{\eqref{eqmu2}}{=} (m_{TX} T \xi_X e_{DX} \circ  Tm_X T^2e_X \xi_{DX} ) \eta_{DX}\\
&\sstackrel{\eqref{eqeta2}}{=} m_{TX} T\xi_X e_{DX} \eta_{DX} \\
&\sstackrel{\eqref{eqeta}}{=} m_{TX} T\xi_X i_{DX} \\
&= m_{TX} i_{T^2X} \xi_X \\
&= \xi_X.
\end{align*}

\noindent \textbf{The second monad unit axiom} $\mu_X D\eta_X = 1_{DX}$: We proceed in an analogous way. Observe that we have:
\begin{align*}
e_X \mu_X D\eta_X &\sstackrel{\eqref{eqmu1}}{=} m_X Te_X e_{DX} D\eta_X \\
&= m_X Te_X T\eta_X e_X \\
&\sstackrel{\eqref{eqeta}}{=} m_X Ti_X e_X \\
&= e_X,
\end{align*}
and:
\begin{align*}
\xi_X \mu_X D\eta_X &\sstackrel{\eqref{eqmu2}}{=} (m_{TX} T\xi_X e_{DX} \circ  Tm_X T^2e_X \xi_{DX})D\eta_X \\
&= m_{TX} T\xi_X e_{DX} D\eta_X \circ  Tm_X T^2e_X \xi_{DX} D\eta_X \\
&= m_{TX} T\xi_X T\eta_X e_X \circ  Tm_X T^2e_X \xi_{DX} D\eta_X \\
&\sstackrel{\eqref{eqeta2}}{=}  Tm_X T^2e_X \xi_{DX} D\eta_X \\
&=  Tm_X T^2e_X T^2\eta_X \xi_X \\
&\sstackrel{\eqref{eqeta}}{=}  Tm_X T^2i_X \xi_X \\
&= \xi_X.
\end{align*}

\textbf{The associativity axiom} $\mu_X \mu_{DX} = \mu_X D\mu_X$:
The equality holds after post-composing with $e_X$:
\begin{align*}
e_X \mu_X D\mu_X &\sstackrel{\eqref{eqmu1}}{=} m_X Te_X e_{DX} D\mu_X \\
                 &= m_X Te_X T\mu_X e_{D^2X} \\
                 &\sstackrel{\eqref{eqmu1}}{=} m_X Tm_X T^2e_X Te_{DX} e_{D^2X} \\
				 &= m_X m_{TX} T^2e_X Te_{DX} e_{D^2X} \\
                 &= m_X Te_X m_{DX} Te_{DX} e_{D^2X} \\
                 &\sstackrel{\eqref{eqmu1}}{=} m_X Te_X e_{DX} \mu_{DX} \\
                 &\sstackrel{\eqref{eqmu1}}{=} e_X \mu_X \mu_{DX}.
\end{align*}

It also holds after post-whiskering with $\xi_X$:
\begin{align*}
\xi_X \mu_X \mu_{DX} =\\
m_{TX} T\xi_X \underline{e_{DX} \mu_{DX}} \circ Tm_X T^2e_X \underline{\xi_{DX} \mu_{DX}}=\\[10pt]
m_{TX} \underline{T\xi_X m_{DX}} Te_{DX} e_{D^2X} \circ Tm_X \underline{T^2e_X m_{TDX}} T\xi_{DX} e_{D^2X} \circ \\ \circ Tm_X \underline{T^2e_X Tm_{DX}} T^2e_{DX}\xi_{D^2X} =\\[10pt]
\underline{m_{TX} m_{T^2X}} T^2\xi_X Te_{DX} e_{D^2X} \circ \underline{Tm_X m_{T^2X}} T^3e_X T\xi_{DX} e_{D^2X} \circ \\ \circ \underline{Tm_X Tm_{TX}} T^3e_X T^2e_{DX}\xi_{D^2X} =\\[10pt]
\underline{m_{TX} Tm_{TX} T^2\xi_X Te_{DX} e_{D^2X} \circ m_{TX} T^2m_X T^3e_X T\xi_{DX} e_{D^2X}} \circ \\ \circ Tm_X \underline{T^2m_X T^3e_X T^2e_{DX}}\xi_{D^2X} =\\[10pt]
m_{TX} T\xi_X \underline{T\mu_X e_{D^2X}} \circ Tm_X T^2e_X \underline{T^2\mu_X \xi_{D^2X}} =\\
m_{TX} T\xi_X e_{DX} D\mu_X \circ Tm_X T^2e_X \xi_{DX} D\mu_X =\\
\xi_X \mu_X D\mu_X.
\end{align*}

Finally, the 2-monad morphism axioms for $e: D \Rightarrow T$:
\[\begin{tikzcd}
	& X && {D^2X} && DX \\
	&&&&&& TX \\
	DX && TX & TDX && {T^2X}
	\arrow["{\eta_X}"', from=1-2, to=3-1]
	\arrow["{i_X}", from=1-2, to=3-3]
	\arrow["{\mu_X}", from=1-4, to=1-6]
	\arrow["{e_{DX}}"', from=1-4, to=3-4]
	\arrow["{e_X}", from=1-6, to=2-7]
	\arrow["{e_X}"', from=3-1, to=3-3]
	\arrow["{Te_X}"', from=3-4, to=3-6]
	\arrow["{m_X}"', from=3-6, to=2-7]
\end{tikzcd}\]
are exactly the identities \eqref{eqeta}, \eqref{eqmu1} for $\eta_X$, $\mu_X$.
\end{proof}

\subsection{One-step termination}\label{SUBS_coKZ_one_step_termination}

Assume that $(T,m,i)$ is a 2-monad on a 2-category $\ck$ such that both $T$ and $T^2$ preserve the descent objects of resolutions of objects $X \in \ob \ck$. Our goal will be to show that $(D,\mu,\eta)$ is colax-idempotent. For this, we will construct a modification $\lambda : \eta D \to D\eta: D \Rightarrow D^2$ satisfying the appropriate equations (see Theorem \ref{THM_kz_monad_TFAE}).

\begin{nota}
Given a cone $(g: B \to TA,\psi)$ for $\opres{A}$ and a 1-cell $h: C \to B$, denote by $h^*(g,\psi)$ the cone $(gh, \psi h)$.
\end{nota}

\begin{lemma}\label{THML_coKZ_conemor_for_delta}
The 2-cell $\xi_X$ is a cone morphism:
\[
(e_{DX}\eta_{DX})^*(Te_X,T\xi_X) \to (e_{DX}D\eta_X)^*(Te_X,T\xi_X).
\]
\end{lemma}
\begin{proof}
Notice that $\xi_X$ really \textbf{is} a 2-cell $Te_X e_{DX} \eta_{DX} \Rightarrow Te_X e_{DX} D\eta_X$ since its domain and codomain are:
\begin{align*}
i_{TX} e_X = Te_X i_{DX} = Te_X e_{DX} \eta_{DX}, \\
Ti_X e_X = Te_X T\eta_X e_X = Te_X e_{DX} D\eta_X.
\end{align*}
The cone morphism axiom says that we must have:
\[\begin{tikzcd}
	DX && {D^2X} && DX && {D^2X} \\
	{D^2X} && TDX && {D^2X} && TDX \\
	TDX && {T^2X} & {=} & TDX && {T^2X} \\
	\\
	{T^2X} && {T^3X} && {T^2X} && {T^3X}
	\arrow[""{name=0, anchor=center, inner sep=0}, "{\eta_{DX}}", from=1-1, to=1-3]
	\arrow[""{name=0p, anchor=center, inner sep=0}, phantom, from=1-1, to=1-3, start anchor=center, end anchor=center]
	\arrow["{D\eta_X}"', from=1-1, to=2-1]
	\arrow["{e_{DX}}", from=1-3, to=2-3]
	\arrow["{\eta_{DX}}", from=1-5, to=1-7]
	\arrow["{D\eta_X}"', from=1-5, to=2-5]
	\arrow["{e_{DX}}", from=1-7, to=2-7]
	\arrow["{e_{DX}}"', from=2-1, to=3-1]
	\arrow["{Te_X}", from=2-3, to=3-3]
	\arrow["{e_{DX}}"', from=2-5, to=3-5]
	\arrow["{\xi_X}"', shorten <=9pt, shorten >=9pt, Rightarrow, from=2-7, to=2-5]
	\arrow["{Te_X}", from=2-7, to=3-7]
	\arrow["{Te_X}"{description}, from=2-7, to=5-5]
	\arrow[""{name=1, anchor=center, inner sep=0}, "{Te_X}"{description}, from=3-1, to=3-3]
	\arrow[""{name=1p, anchor=center, inner sep=0}, phantom, from=3-1, to=3-3, start anchor=center, end anchor=center]
	\arrow[""{name=1p, anchor=center, inner sep=0}, phantom, from=3-1, to=3-3, start anchor=center, end anchor=center]
	\arrow["{Te_X}"', from=3-1, to=5-1]
	\arrow["{Ti_{TX}}", from=3-3, to=5-3]
	\arrow[""{name=2, anchor=center, inner sep=0}, "{Te_X}"', from=3-5, to=5-5]
	\arrow[""{name=3, anchor=center, inner sep=0}, "{Ti_{TX}}", from=3-7, to=5-7]
	\arrow[""{name=4, anchor=center, inner sep=0}, "{T^2i_X}"', from=5-1, to=5-3]
	\arrow[""{name=4p, anchor=center, inner sep=0}, phantom, from=5-1, to=5-3, start anchor=center, end anchor=center]
	\arrow["{T^2i_X}"', from=5-5, to=5-7]
	\arrow["{\xi_X}", shorten <=9pt, shorten >=9pt, Rightarrow, from=0p, to=1p]
	\arrow["{T\xi_X}", shorten <=9pt, shorten >=9pt, Rightarrow, from=1p, to=4p]
	\arrow["{T\xi_X}"{pos=0.3}, shorten <=14pt, shorten >=41pt, Rightarrow, from=3, to=2]
\end{tikzcd}\]

This holds since we have:
\begin{align*}
T\xi_X e_{DX} D\eta_X \circ Ti_{TX}\xi_X &= T\xi_X T\eta_X e_X \circ Ti_{TX}\xi_X \\
&\sstackrel{\eqref{eqeta2}}{=} Ti_{TX} \xi_X \\
&\sstackrel{\eqref{eqconeaxiom2}}{=} T^2i_X \xi_X \circ i_{T^2X} \xi_X \\
&= T^2i_X \xi_X \circ T\xi_X i_{DX} \\
&\sstackrel{\eqref{eqeta}}{=} T^2i_X \xi_X \circ T\xi_X e_{DX} \eta_{DX}.
\end{align*}
\end{proof}

\noindent Because $(Te_X, T\xi_X)$ is by assumption the descent object of $T\opres{X}$, there is a unique 2-cell:
\[
\delta_X : e_{DX} \eta_{DX} \Rightarrow e_{DX} D\eta_X,
\]
with the property:
\begin{equation}\label{eqdelta}
(Te_X)\delta_X = \xi_X.
\end{equation}
By Lemma \ref{THML_cone_in_KK}, $\delta: eD \circ \eta D \to eD \circ D\eta$ is a modification.

\begin{lemma}\label{THML_coKZ_conemor_for_lambda}
The 2-cell $\delta_X$ is a cone morphism:
\[
(\eta_{DX})^*(e_{DX},\xi_{DX}) \to (D\eta_X)^* (e_{DX},\xi_{DX}).
\]
\end{lemma}
\begin{proof}
We wish to show that:
\[\begin{tikzcd}
	DX && {D^2X} && DX && {D^2X} \\
	\\
	{D^2X} && TDX & {=} & {D^2X} && TDX \\
	\\
	TDX && {T^2DX} && TDX && {T^2DX}
	\arrow[""{name=0, anchor=center, inner sep=0}, "{\eta_{DX}}", from=1-1, to=1-3]
	\arrow["{D\eta_X}"', from=1-1, to=3-1]
	\arrow["{e_{DX}}", from=1-3, to=3-3]
	\arrow["{\eta_{DX}}", from=1-5, to=1-7]
	\arrow["{D\eta_X}"', from=1-5, to=3-5]
	\arrow["{e_{DX}}", from=1-7, to=3-7]
	\arrow[""{name=1, anchor=center, inner sep=0}, "{e_{DX}}"{description, pos=0.7}, from=1-7, to=5-5]
	\arrow[""{name=2, anchor=center, inner sep=0}, "{Te_X}"{description}, from=3-1, to=3-3]
	\arrow["{e_{DX}}"', from=3-1, to=5-1]
	\arrow["{i_{TDX}}", from=3-3, to=5-3]
	\arrow["{e_{DX}}"', from=3-5, to=5-5]
	\arrow["{i_{TDX}}", from=3-7, to=5-7]
	\arrow[""{name=3, anchor=center, inner sep=0}, "{Ti_{DX}}"', from=5-1, to=5-3]
	\arrow["{Ti_{DX}}"', from=5-5, to=5-7]
	\arrow["{\delta_X}", shorten <=9pt, shorten >=9pt, Rightarrow, from=0, to=2]
	\arrow["{\delta_X}"', shift right=5, shorten <=5pt, Rightarrow, from=1, to=3-5]
	\arrow["{\xi_{DX}}", shorten <=9pt, shorten >=9pt, Rightarrow, from=2, to=3]
	\arrow["{\xi_{DX}}", shift left=5, shorten >=5pt, Rightarrow, from=3-7, to=1]
\end{tikzcd}\]

Because by assumption, $(T^2e_X, T^2\xi_X)$ is the descent object of $T^2 \opres{X}$, it suffices to show that they are equal after post-whiskering with $T^2e_X$. We then obtain:
\begin{align*}
T^2e_X(LHS) &= T^2e_X \xi_{DX} D\eta_X \circ T^2e_X i_{TDX} \delta_X \\
&= T^2e_X \xi_{DX} D\eta_X \circ i_{T^2X}Te_X \delta_X \\
&\sstackrel{\eqref{eqdelta}}{=} T^2e_X \xi_{DX} D\eta_X \circ i_{T^2X}\xi_X \\
&= T^2e_X T^2\eta_X \xi_X \circ i_{T^2X} \xi_X \\
&\sstackrel{\eqref{eqeta}}{=} T^2 i_X \xi_X \circ i_{T^2X} \xi_X \\
&\sstackrel{\eqref{eqconeaxiom2}}{=} Ti_{TX}\xi_X \\
&\sstackrel{\eqref{eqdelta}}{=} Ti_{TX} Te_X \delta_X\\
&\sstackrel{\eqref{eqeta2}}{=} Ti_{TX} Te_X \delta_X \circ T^2e_X \xi_{DX} \eta_{DX}\\
&= T^2e_X Ti_{DX}  \delta_X \circ T^2e_X \xi_{DX} \eta_{DX}\\
&= T^2e_X(RHS).
\end{align*}
\end{proof}

\noindent Because $(e_{DX},\xi_{DX})$ is a limit cone, there exists a unique 2-cell:
\[
\lambda_X: \eta_{DX} \Rightarrow D\eta_X,
\]
satisfying
\begin{equation}\label{eqlambda}
e_{DX} \lambda_X = \delta_X.
\end{equation}
By Lemma \ref{THML_cone_in_KK}, this is a modification $\lambda: \eta D \to D \eta$.

\begin{theo}\label{THM_one_step_kzfication}
Let $(T,m,i)$ be a 2-monad that preserves descent objects of coreflexive op-coherence data. Then $(D,\mu,\eta, \lambda)$ is colax-idempotent.
\end{theo}
\begin{proof}
By Theorem \ref{THM_kz_monad_TFAE}, it remains to show that:
\begin{align*}
\mu_X \lambda_X &= 1,\\
\lambda_X \eta_X &= 1.
\end{align*}

It suffices to prove the first equation after post-whiskering with $e_X$:
\[
e_X \mu_X \lambda_X \sstackrel{\eqref{eqmu1}}{=} m_X Te_X e_{DX} \lambda_X \sstackrel{\eqref{eqlambda}}{=} m_X Te_X \delta_X \\ \sstackrel{\eqref{eqdelta}}{=} m_X \xi_X \sstackrel{\eqref{eqconeaxiom1}}{=} 1.
\]
It suffices to prove the second equation after post-composing with $e_{DX}$:
\[
e_{DX}\lambda_X \eta_X \sstackrel{\eqref{eqlambda}}{=} \delta_X \eta_X = 1,
\]
where this last equality can be proven by post-whiskering with $Te_X$, in which case it is the identity by \eqref{eqdelta} and \eqref{eqeta2}.
\end{proof}

\subsection{The universal property of $D$}\label{SUBS_coKZ_coreflection}

In the remainder of this section we will prove that this ``co-KZ-fication” of $(T,m,i)$ enjoys a universal property. Recall first that 2-monad morphisms between lax-idempotent 2-monads ``preserve” the lax-idempotence 2-cell in the following sense:

\begin{lemma}\label{THML_2monadmor_preserves_lambda_sigma}
Let $(T,m,i,\lambda)$, $(S,n,j,\sigma)$ be two lax-idempotent 2-monads on $\ck$ and let $\theta: T \Rightarrow S$ be a 2-monad morphism. Then $\theta$ commutes with the 2-cells $\lambda, \sigma$ in the sense that there is the following equality:
\[\begin{tikzcd}
	&&&&& SX \\
	&&&&& TSX \\
	{\sigma_X \theta_X =} & TX && {T^2X} &&&& {S^2X} \\
	&&&&& STX \\
	&&&&& SX
	\arrow["{Sj_X}", curve={height=-18pt}, from=1-6, to=3-8]
	\arrow["{\theta_{SX}}", from=2-6, to=3-8]
	\arrow["{\theta_X}", curve={height=-30pt}, from=3-2, to=1-6]
	\arrow["{Tj_X}"{description}, curve={height=-30pt}, from=3-2, to=2-6]
	\arrow[""{name=0, anchor=center, inner sep=0}, "{i_{TX}}"{description}, curve={height=12pt}, from=3-2, to=3-4]
	\arrow[""{name=1, anchor=center, inner sep=0}, "{Ti_X}"{description}, curve={height=-12pt}, from=3-2, to=3-4]
	\arrow["{j_{TX}}"{description}, curve={height=18pt}, from=3-2, to=4-6]
	\arrow["{\theta_X}"', curve={height=24pt}, from=3-2, to=5-6]
	\arrow["{T\theta_X}", from=3-4, to=2-6]
	\arrow["{\theta_{TX}}"', from=3-4, to=4-6]
	\arrow["{S\theta_X}"', from=4-6, to=3-8]
	\arrow["{j_{SX}}"', curve={height=18pt}, from=5-6, to=3-8]
	\arrow["{\lambda_X}", shorten <=3pt, shorten >=3pt, Rightarrow, from=1, to=0]
\end{tikzcd}\]
\end{lemma}
\begin{proof}
By transposing under the adjunction $Sj_X \dashv n_X$, we see that these 2-cells are equal if and only if they are equal after post-whiskering with the 2-monad multiplication $n_X:S^2X \to SX$. The left-hand side is then the identity by the axiom \eqref{EQ_KZ_psmonad_mAlambdaA}. Because of the 2-monad morphism axiom, the right-hand side also reduces to the identity:
\[
n_X S\theta_X \theta_{TX} \lambda_X = \theta_X m_X \lambda_X \sstackrel{\eqref{EQ_KZ_psmonad_mAlambdaA}}{=} 1.
\]
\end{proof}

\begin{theo}\label{THM_coKZfication}
Let $(T,m,i)$ be a 2-monad on a 2-category $\ck$ which admits descent objects of coreflexive op-coherence data and assume that $T$ preserves them. Then the 2-monad morphism $e$ from Lemma \ref{THML_e_is_2monad_morphism} exhibits $(D,\mu,\eta)$ as the coreflection of $(T,m,i)$ along the inclusion:
\[
\text{ColaxIdemp}(\ck) \xhookrightarrow{} \mndioo{\twoCAT}{\ck},
\]
of the full subcategory of the category of 2-monads and 2-monads morphisms on $\ck$ spanned by colax-idempotent 2-monads.
\end{theo}
\begin{proof}
We have to prove that for any 2-monad morphism $k: S \to T$ with $S$ colax-idempotent, there is a unique 2-monad morphism $\theta: S \to D$ making the diagram commute:
\[\begin{tikzcd}
	{(D,\mu,\eta,\lambda)} && {(T,m,i)} \\
	\\
	&& {(S,n,j,\sigma)}
	\arrow["e", from=1-1, to=1-3]
	\arrow["{\exists ! \theta}", dashed, from=3-3, to=1-1]
	\arrow["{\forall k}"', from=3-3, to=1-3]
\end{tikzcd}\]
The proof is straightforward and will be done in steps:
\begin{itemize}
\item we show that for each $X \in \ck$, there is a cone $(k_X, \overline{k_X})$ for $\opres{X}$, giving us a unique 1-cell $\theta_X : SX \to DX$ to the descent object, or, a 2-natural transformation $\theta: S \Rightarrow D$,
\item we show that $\theta$ is a 2-monad morphism $\theta: S \to D$ that makes the above diagram commute,
\item we show that any other 2-monad morphism $\tau: S \to D$ making the above diagram commute equals $\theta$.
\end{itemize}

\noindent \textbf{Step 1}:
We claim that the pair $(k_X, \overline{k}_X)$, where $\overline{k}_X$ is the 2-cell bellow, is a cone for $\opres{X}$:
\[\begin{tikzcd}
	&& TX \\
	SX && {S^2X} && STX && {T^2X} \\
	\\
	&&& TX
	\arrow["{j_{TX}}"{description}, from=1-3, to=2-5]
	\arrow["{i_{TX}}", curve={height=-12pt}, from=1-3, to=2-7]
	\arrow["{k_X}", curve={height=-18pt}, from=2-1, to=1-3]
	\arrow[""{name=0, anchor=center, inner sep=0}, "{Sj_X}"', curve={height=12pt}, from=2-1, to=2-3]
	\arrow[""{name=1, anchor=center, inner sep=0}, "{j_{SX}}", curve={height=-12pt}, from=2-1, to=2-3]
	\arrow["{Si_X}"', shift right=3, curve={height=30pt}, from=2-1, to=2-5]
	\arrow["{k_X}"', curve={height=30pt}, from=2-1, to=4-4]
	\arrow["{Sk_X}"', from=2-3, to=2-5]
	\arrow["{k_{TX}}"{description}, from=2-5, to=2-7]
	\arrow["{Ti_X}", from=4-4, to=2-7]
	\arrow["{\sigma_X}"', shorten <=3pt, shorten >=3pt, Rightarrow, from=1, to=0]
\end{tikzcd}\]
The axiom \eqref{eqconeaxiom1} follows from the 2-monad morphism axiom:
\[
m_X k_{TX} Sk_X \sigma_X = k_X n_X \sigma_X \sstackrel{\eqref{EQ_KZ_psmonad_mAlambdaA}}{=} 1.
\]
To see that the axiom \eqref{eqconeaxiom2} holds, we need to show that the 2-cell below equals the 2-cell $Ti_{TX} k_{TX} Sk_X \sigma_X$:
\[\begin{tikzcd}[column sep=scriptsize]
	SX && {S^2X} && STX && {T^2X} & {T^3X} \\
	\\
	&& TX &&& {T^2X} \\
	& {S^2X} &&& STX
	\arrow[""{name=0, anchor=center, inner sep=0}, "{Sj_X}"', curve={height=12pt}, from=1-1, to=1-3]
	\arrow[""{name=1, anchor=center, inner sep=0}, "{j_{SX}}", curve={height=-12pt}, from=1-1, to=1-3]
	\arrow["{Si_X}"{description}, shift right=3, curve={height=30pt}, from=1-1, to=1-5]
	\arrow["{k_X}", curve={height=18pt}, from=1-1, to=3-3]
	\arrow[""{name=2, anchor=center, inner sep=0}, "{j_{SX}}"{description, pos=0.7}, curve={height=12pt}, from=1-1, to=4-2]
	\arrow[""{name=3, anchor=center, inner sep=0}, "{Sj_X}"', shift right=5, curve={height=30pt}, from=1-1, to=4-2]
	\arrow["{Sk_X}", from=1-3, to=1-5]
	\arrow["{k_{TX}}", from=1-5, to=1-7]
	\arrow["{i_{T^2X}}", from=1-7, to=1-8]
	\arrow["{Ti_X}"{description}, from=3-3, to=1-7]
	\arrow["{i_{TX}}"{description}, from=3-3, to=3-6]
	\arrow["{j_{TX}}"{description}, from=3-3, to=4-5]
	\arrow["{T^2i_X}"', from=3-6, to=1-8]
	\arrow["{Sk_X}"', from=4-2, to=4-5]
	\arrow["{k_{TX}}"', from=4-5, to=3-6]
	\arrow["{\sigma_X}"', shorten <=3pt, shorten >=3pt, Rightarrow, from=1, to=0]
	\arrow["{\sigma_X}"', shorten <=5pt, shorten >=5pt, Rightarrow, from=2, to=3]
\end{tikzcd}\]
Notice that we have:
\begin{align*}
\underline{T^2i_X k_{TX} Sk_X} \sigma_X &= k_{T^2X} \underline{STi_X Sk_X} \sigma_X \\
&= k_{T^2X} Sk_{TX} \underline{S^2i_X} \sigma_X \\
&= k_{T^2X} Sk_{TX} \underline{S^2k_X S^2j_X \sigma_X} \\
&= k_{T^2X} Sk_{TX} \sigma_{TX} Sk_X Sj_X.
\end{align*}
The composite thus reduces to:
\[\begin{tikzcd}[column sep=scriptsize]
	&& TX \\
	SX && {S^2X} && STX && {T^2X} && {T^3X} \\
	&&&&&& TSTX \\
	&&&& {S^2TX} &&& {ST^2X}
	\arrow["{j_{TX}}", curve={height=-12pt}, dashed, from=1-3, to=2-5]
	\arrow["{k_X}", curve={height=-18pt}, dashed, from=2-1, to=1-3]
	\arrow[""{name=0, anchor=center, inner sep=0}, "{Sj_X}"', curve={height=12pt}, from=2-1, to=2-3]
	\arrow[""{name=1, anchor=center, inner sep=0}, "{j_{SX}}", curve={height=-12pt}, from=2-1, to=2-3]
	\arrow["{Sk_X}", from=2-3, to=2-5]
	\arrow["{k_{TX}}", from=2-5, to=2-7]
	\arrow["{i_{STX}}"{description}, from=2-5, to=3-7]
	\arrow[""{name=2, anchor=center, inner sep=0}, "{Sj_{TX}}"', curve={height=24pt}, from=2-5, to=4-5]
	\arrow[""{name=3, anchor=center, inner sep=0}, "{j_{STX}}", curve={height=-24pt}, from=2-5, to=4-5]
	\arrow["{i_{T^2X}}", from=2-7, to=2-9]
	\arrow["{Tk_{TX}}"{description}, from=3-7, to=2-9]
	\arrow["{k_{STX}}"{description}, from=4-5, to=3-7]
	\arrow["{Sk_{TX}}"', from=4-5, to=4-8]
	\arrow["{k_{T^2X}}"', from=4-8, to=2-9]
	\arrow["{\sigma_X}"', shorten <=3pt, shorten >=3pt, Rightarrow, from=1, to=0]
	\arrow["{\sigma_{TX}}"', shorten <=10pt, shorten >=10pt, Rightarrow, from=3, to=2]
\end{tikzcd}\]
Using the dotted equality and the axiom \eqref{EQ_KZ_psmonad_lambdaAiA} for the modification $\sigma$, the whiskered 2-cell $\sigma_{TX}$ disappears. Now the fact that the remaining 2-cell $k_{T^2X} Sk_{TX} Sj_{TX} Sk_X \sigma_X$ equals $Ti_{TX} k_{TX} Sk_X \sigma_X$ follows from this manipulation:
\[
k_{T^2X} Sk_{TX} Sj_{TX} = k_{T^2X} Si_{TX} =  Ti_{TX} k_{TX}.
\]

There thus exists a unique 1-cell $\theta_X : SX \to DX$ satisfying:
\begin{align}
e_X \theta_X &= k_X, \label{eqtheta1}\\
\xi_X \theta_X &= k_{TX} Sk_X \sigma_X. \label{eqtheta2}
\end{align}

By Lemma \ref{THML_cone_in_KK}, this is a 2-natural transformation $\theta: S \Rightarrow D$.

\noindent \textbf{Step 2}: To see that $\theta$ is a 2-monad morphism $S \to D$, we have to verify the following:
\[\begin{tikzcd}
	& X && {S^2X} && SX && DX \\
	SX && DX & SDX && {D^2X}
	\arrow["{j_X}"', from=1-2, to=2-1]
	\arrow["{\eta_X}", from=1-2, to=2-3]
	\arrow["{n_X}", from=1-4, to=1-6]
	\arrow["{S\theta_X}"', from=1-4, to=2-4]
	\arrow["{\theta_X}", from=1-6, to=1-8]
	\arrow["{\theta_X}"', from=2-1, to=2-3]
	\arrow["{\theta_{DX}}"', from=2-4, to=2-6]
	\arrow["{\mu_X}"', from=2-6, to=1-8]
\end{tikzcd}\]

The equation above left holds since it holds after post-composing and post-whiskering with $e_X$, $\xi_X$:
\begin{align*}
&e_X \theta_X j_X \sstackrel{\eqref{eqtheta1}}{=} k_X j_X = i_X \sstackrel{\eqref{eqeta}}{=} e_X \eta_X,\\
&\xi_X \theta_X j_X \phantom{.}\sstackrel{\eqref{eqtheta2}}{=}\phantom{.} k_{TX} Sk_X \sigma_X j_X \sstackrel{\eqref{EQ_KZ_psmonad_lambdaAiA}}{=} \phantom{.}1\phantom{.} \sstackrel{\eqref{eqeta2}}{=} \xi_X \eta_X.
\end{align*}

Likewise, the second equation will follow if we show that:
\begin{align*}
e_X \theta_X n_X &= e_X \mu_X \theta_{DX} S\theta_X,\\
\xi_X \theta_X n_X &= \xi_X \mu_X \theta_{DX} S\theta_X.
\end{align*}
The first equation follows since the following diagram commutes:
\[\begin{tikzcd}[column sep=small]
	&&&&& SX \\
	{S^2X} &&&&& TSX &&& DX \\
	&&& DSX \\
	&& STX && DTX && {T^2X} && TX \\
	&& SDX && {D^2X} && DX
	\arrow["{\theta_X}", from=1-6, to=2-9]
	\arrow["{k_X}", from=1-6, to=4-9]
	\arrow["{n_X}", from=2-1, to=1-6]
	\arrow["{k_{SX}}", from=2-1, to=2-6]
	\arrow["{\theta_{SX}}", from=2-1, to=3-4]
	\arrow["{Sk_X}", from=2-1, to=4-3]
	\arrow["{S\theta_X}"', from=2-1, to=5-3]
	\arrow["{Tk_X}", from=2-6, to=4-7]
	\arrow["{e_X}", from=2-9, to=4-9]
	\arrow["{e_{SX}}", from=3-4, to=2-6]
	\arrow["{Dk_X}", from=3-4, to=4-5]
	\arrow["{\theta_{TX}}"', from=4-3, to=4-5]
	\arrow["{e_{TX}}"', from=4-5, to=4-7]
	\arrow["{m_X}", from=4-7, to=4-9]
	\arrow["{Se_X}"', from=5-3, to=4-3]
	\arrow["{\theta_{DX}}"', from=5-3, to=5-5]
	\arrow["{De_X}", from=5-5, to=4-5]
	\arrow["{\mu_X}"', from=5-5, to=5-7]
	\arrow["{e_X}"', from=5-7, to=4-9]
\end{tikzcd}\]

To prove the second equation, note that the right-hand side reduces to:
\begin{align*}
\xi_X \mu_X \theta_{DX} S\theta_X &\sstackrel{\eqref{eqmu2}}{=} m_{TX} T\xi_X \underline{e_{DX} \theta_{DX}} S\theta_X \circ Tm_X \xi_{TX} \underline{De_X \theta_{DX}} S\theta_X \\
&= m_{TX} \underline{T\xi_X k_{DX}} S\theta_X \circ Tm_X \underline{\xi_{TX} \theta_{TX}} Se_X S\theta_X \\
&= m_{TX} k_{T^2X} \underline{S\xi_X S\theta_X} \circ Tm_X k_{T^2X} Sk_{TX} \sigma_{TX} \underline{Se_X S\theta_X} \\
&= \underline{m_{TX} k_{T^2X} Sk_{TX}} S^2k_X S\sigma_X \circ Tm_X k_{T^2X} Sk_{TX} \underline{\sigma_{TX} Sk_X} \\
&=  k_{TX} \underline{n_{TX} S^2k_X} S\sigma_X \circ Tm_X \underline{k_{T^2X} Sk_{TX} S^2k_X} \sigma_{SX} \\
&\sstackrel{(*)}{=} k_{TX} Sk_X n_{SX} S\sigma_X \circ \underline{Tm_X T^2k_X Tk_{SX}} k_{S^2X} \sigma_{SX} \\
&= k_{TX} Sk_X n_{SX} S\sigma_X \circ Tk_X \underline{Tn_X k_{S^2X}} \sigma_{SX} \\
&= k_{TX} Sk_X n_{SX} S\sigma_X \circ \underline{Tk_X k_{SX}} Sn_X \sigma_{SX} \\
&= k_{TX} Sk_X n_{SX} S\sigma_X \circ k_{TX} Sk_X Sn_X \sigma_{SX}.
\end{align*}
Here $(*)$ uses the naturality of $k: S \Rightarrow T$ as follows:
\[\begin{tikzcd}
	{S^3X} && {S^2TX} \\
	&&& {ST^2X} \\
	{TS^2X} && TSTX \\
	& {T^2SX} && {T^3X}
	\arrow["{S^2k_X}", from=1-1, to=1-3]
	\arrow["{k_{S^2X}}"', from=1-1, to=3-1]
	\arrow["{Sk_{TX}}", from=1-3, to=2-4]
	\arrow["{k_{STX}}"{description}, from=1-3, to=3-3]
	\arrow["{k_{T^2X}}", from=2-4, to=4-4]
	\arrow["{TSk_X}", from=3-1, to=3-3]
	\arrow["{Tk_{SX}}"', from=3-1, to=4-2]
	\arrow["{Tk_{TX}}", from=3-3, to=4-4]
	\arrow["{T^2k_X}"', from=4-2, to=4-4]
\end{tikzcd}\]
The left-hand side reduces to $\xi_X \theta_X n_X \sstackrel{\eqref{eqtheta2}}{=} k_{TX} Sk_X \sigma_X n_X$. The left and right side are now equal by Lemma \ref{THML_lambda_unique_2cells_equal}. Thus $\theta: S \to D$ is a 2-monad morphism. The fact that it makes the diagram in the statement commute follows from its definition -- \eqref{eqtheta1}.

\noindent \textbf{Step 3}: Let $\tau: S \Rightarrow D$ be another 2-monad morphism making the triangle in the statement commute. By Lemma \ref{THML_2monadmor_preserves_lambda_sigma}, we have that $\lambda_X \tau_X = \tau_{DX} S\tau_X \sigma_X$. In particular:
\begin{align*}
\xi_X \tau_X &\sstackrel{\eqref{eqdelta}}{=} Te_X \delta_X \tau_X \\
			&\sstackrel{\eqref{eqlambda}}{=} Te_X e_{DX} \lambda_X \tau_X \\
			&= Te_X \underline{e_{DX} \tau_{DX}} S\tau_X \sigma_X \\
			&= Te_X \underline{k_{DX} S\tau_X} \sigma_X \\
			&= \underline{Te_X T\tau_X} k_{SX} \sigma_X \\
			&= Tk_X k_{SX} \sigma_X.
\end{align*}
But $\theta_X: SX \to DX$ is the unique 1-cell satisfying $k_X = e_X \theta_X$ and $\xi_X \theta_X = k_{TX} Sk_X \sigma_X$ -- thus $\theta_X = \tau_X$.
\end{proof}

\subsection{Examples}\label{SEC_coKZfication_further_props}

Observe:

\begin{prop}
A 2-monad $(T,m,i)$ is colax-idempotent if and only if for every $X$ there is a 2-cell $\sigma_X: i_{TX} \Rightarrow Ti_X$ making the pair $(1_{TX}, \sigma_X)$ the descent object of $\opres{X}$.
\end{prop}
\begin{proof}
``$\Rightarrow$”: Let $(T,m,i,\sigma)$ be a colax-idempotent 2-monad. Let us show that $(1_{TX}, \sigma_X)$ is a cone for $\opres{X}$. The cone axiom $m_X \sigma_X = 1_{TX}$ is the axiom for $\sigma$ (see \eqref{EQ_KZ_psmonad_mAlambdaA}). To verify the cone axiom \eqref{eqconeaxiom2}, using the local naturality of $i$, the right-hand side of \eqref{eqconeaxiom2} becomes:
\[\begin{tikzcd}
	TX && {T^2X} && {T^3X}
	\arrow[""{name=0, anchor=center, inner sep=0}, "{Ti_X}"', curve={height=24pt}, from=1-1, to=1-3]
	\arrow[""{name=1, anchor=center, inner sep=0}, "{i_{TX}}", curve={height=-24pt}, from=1-1, to=1-3]
	\arrow[""{name=2, anchor=center, inner sep=0}, "{T^2i_X}"', curve={height=24pt}, from=1-3, to=1-5]
	\arrow[""{name=3, anchor=center, inner sep=0}, "{Ti_{TX}}", curve={height=-24pt}, from=1-3, to=1-5]
	\arrow["{\sigma_X}", shorten <=6pt, shorten >=6pt, Rightarrow, from=1, to=0]
	\arrow["{T\sigma_X}", shorten <=6pt, shorten >=6pt, Rightarrow, from=3, to=2]
\end{tikzcd}\]
Using the middle-four-interchange rule and the axiom \eqref{EQ_KZ_psmonad_lambdaAiA}, this becomes $Ti_{TX} \sigma_X$.

To prove that this cone is the limit, let $(g: A \to TX, \alpha: i_{TX} g \Rightarrow Ti_X g)$ be a different cone for $\opres{X}$. We would like to show that $\alpha = \sigma_X g$. By taking mates, this equality holds if and only if the following composites are equal:
\[\begin{tikzcd}[column sep=scriptsize]
	& TX \\
	A && {T^2X} && A & TX && {T^2X} \\
	& TX && TX &&&&& TX
	\arrow["{i_{TX}}", from=1-2, to=2-3]
	\arrow["\alpha", shorten <=6pt, shorten >=6pt, Rightarrow, from=1-2, to=3-2]
	\arrow["g", from=2-1, to=1-2]
	\arrow["g"', from=2-1, to=3-2]
	\arrow["{m_X}", from=2-3, to=3-4]
	\arrow["g", from=2-5, to=2-6]
	\arrow[""{name=0, anchor=center, inner sep=0}, "{i_{TX}}", curve={height=-24pt}, from=2-6, to=2-8]
	\arrow[""{name=1, anchor=center, inner sep=0}, "{Ti_X}"{description}, curve={height=24pt}, from=2-6, to=2-8]
	\arrow[""{name=2, anchor=center, inner sep=0}, shift right, curve={height=30pt}, Rightarrow, no head, from=2-6, to=3-9]
	\arrow["{m_X}", from=2-8, to=3-9]
	\arrow["{Ti_X}"{description}, from=3-2, to=2-3]
	\arrow[""{name=3, anchor=center, inner sep=0}, curve={height=24pt}, Rightarrow, no head, from=3-2, to=3-4]
	\arrow["{\eta_X^{-1}}", shorten <=6pt, shorten >=6pt, Rightarrow, from=2-3, to=3]
	\arrow["{\sigma_X}", shorten <=6pt, shorten >=6pt, Rightarrow, from=0, to=1]
	\arrow["{\eta_X^{-1}}"{pos=0.3}, shift left=5, shorten <=6pt, shorten >=6pt, Rightarrow, from=2-8, to=2]
\end{tikzcd}\]
This is the case because of the cone axiom \eqref{eqconeaxiom1}. The two-dimensional universal property is immediate because the 1-cell component of the limit cone is the identity.

``$\Leftarrow$”: Let $(1_{TX}, \sigma_X)$ be the descent object of $\opres{X}$. All that needs to be done is to show that $\sigma_X i_X = 1$. This follows from the fact that the pair $(i_X: X \to TX, 1_{i_{TX} i_X})$ is a cone for $\opres{X}$.
\end{proof}

This proposition tells us two things -- first, if $T$ is a colax-idempotent 2-monad on $\ck$, this 2-category automatically admits descent objects of diagrams of form $\opres{X}$. We may thus proceed with the process from Section \ref{SUBS_coKZ_process}. The second thing it tells us is that this process will do nothing: we will have $D = T$, and tracing through the proofs in \ref{SUBS_coKZ_process}, we will also have $\mu = m$, $\eta = i$.

\bigskip

Our next application will concern lax-idempotent 2-comonads. This amounts to replacing $\ck$ by $\ck^{coop}$, and the descent objects involved here will become the codescent objects of Section \ref{SUBS_DEFI_CODOBJ}. Fix now a 2-monad $(T,m,i)$ on a 2-category $\ck$, and consider the 2-comonad $FU$ on $\talgs$ generated by the free-forgetful 2-adjunction:
\[\begin{tikzcd}
	\talgs && \ck
	\arrow[""{name=0, anchor=center, inner sep=0}, "U"', curve={height=24pt}, from=1-1, to=1-3]
	\arrow[""{name=1, anchor=center, inner sep=0}, "F"', curve={height=24pt}, from=1-3, to=1-1]
	\arrow["\dashv"{anchor=center, rotate=-90}, draw=none, from=1, to=0]
\end{tikzcd}\]
The resolution associated to an object $(A,a) \in \talgs$ is precisely the resolution of an algebra (Example \ref{PR_res_of_alg}). Assume that $\talgs$ admits codescent objects of these resolutions. By Theorem \ref{THM_adjunkce_talgs_talgl_colaxtalgl} we can thus observe that this process turns the 2-comonad $FU$ into the lax morphism classifier 2-comonad $Q_l$. We obtain a theorem observed by Richard Garner in \cite{ncatcafe}:

\begin{theo}[Universal property of the lax morphism classifier]\label{THM_univ_prop_of_Ql}
Assume that $T$ is a 2-monad on a 2-category $\ck$ for which the inclusion $\talgs \to \talgl$ admits a left 2-adjoint so that the lax morphism classifier 2-comonad $Q_l$ exists. Assume that either:
\begin{itemize}
\item $T$ preserves reflexive codescent objects,
\item $\ck$ admits oplax limits of arrows.
\end{itemize}
Then $Q_l$ is the reflection of the 2-comonad $FU$ along the inclusion of lax-idempotent 2-comonads to 2-comonads on $\talgs$.
\end{theo}

\noindent Let us finish with the one-dimensional special case. Let $(T,m,i)$ be a 2-monad on a locally discrete 2-category. In other words, a monad on a category. Notice that a coreflexive descent object is the same thing as a coreflexive equalizer -- our process in Section \ref{SUBS_coKZ_process} reduces to that of Fakir \cite{fakir_HIST}. Since a 2-monad $(D,\mu,\eta)$ on such a 2-category is colax-idempotent if and only if it is an idempotent monad, our results in particular imply:

\begin{prop}
Let $(T,m,i)$ be a monad on a category $\cc$ that admits reflexive equalizers and $T$ preserves them. Then $T$ admits a coreflection along the inclusion of idempotent monads to monads on $\cc$:
\[
\text{IdempMonad}(\cc) \xhookrightarrow{} \text{Monad}(\cc).
\]
\end{prop}

\section{Characterization of lax-pie algebras}\label{SEC_lax_pie_algebras}

\subsection{Comonadicity of multicategories}

\begin{ass}
Throughout this section, let $\ce$ be a category with pullbacks and let $T$ be a 2-monad on $\cat(\ce)$ of form $\cat(T)$.
\end{ass}

\noindent Recall the following statements about the relationship between $T$-algebras and $T$-multicate\-gories:

\begin{lemma}\label{THML_ff_functor_talgl_multicat}
There is a fully faithful functor:
\[
I: \colaxtalgl \xhookrightarrow{} \tmult.
\]
\end{lemma}
\begin{proof}
This follows from \cite[Theorem 9.13]{unified}.
\end{proof}

\begin{theo}\label{THM_talgs_tmult_adj}
Denote by $U: \talgs \to \tmult$ the restriction of the functor from Proposition \ref{THML_ff_functor_talgl_multicat}. The functor $U$ admits a left adjoint:
\[\begin{tikzcd}
	\talgs && {T\text{-Multicat}}
	\arrow[""{name=0, anchor=center, inner sep=0}, "U"', curve={height=12pt}, from=1-1, to=1-3]
	\arrow[""{name=1, anchor=center, inner sep=0}, "F"', curve={height=12pt}, from=1-3, to=1-1]
	\arrow["\dashv"{anchor=center, rotate=-90}, draw=none, from=1, to=0]
\end{tikzcd}\]
\end{theo}
\begin{proof}
For the full proof see \cite{generaloperadsandmulticats}, we will only describe the action of the functor $F$. Consider a $T$-multicategory $\cm$:
\[\begin{tikzcd}
	{T\cm_0} && {\cm_1} && {\cm_0}
	\arrow["s"', from=1-3, to=1-1]
	\arrow["t", from=1-3, to=1-5]
\end{tikzcd}\]

The $T$-algebra $F\cm$ given as follows. The underlying category has the following graph:
\[\begin{tikzcd}
	&&&& {T\cm_1} \\
	&& {T^2\cm_1} &&&& {T\cm_0} \\
	{T\cm_0}
	\arrow["Ts"', from=1-5, to=2-3]
	\arrow["Tt", from=1-5, to=2-7]
	\arrow["{m_{\cm_0}}"', from=2-3, to=3-1]
\end{tikzcd}\]

The algebra structure map $a: TF\cm \to F\cm$ has components:
\begin{align*}
a_0 &:= m_{\cm_0}: T^2\cm_0 \to T\cm_0,\\
a_1 &:= m_{\cm_1} : T^2\cm_1 \to T\cm_1.
\end{align*}

Now given a $T$-multifunctor $f: \cm \to \cn$, the $T$-algebra morphism $Ff$ has components:
\begin{align*}
(Ff)_0 &:= Tf_0 ,\\
(Ff)_1 &:= Tf_1 .
\end{align*}
\end{proof}

We will aim to show that the above adjunction is 2-comonadic. We will need the following elementary lemma:

\begin{lemma}\label{THML_pbpresfun_pres_corefleq}
Let $\ce$ be a category and $f,g: A \to B$ a \textit{coreflexive pair} in $\ce$, that is, there is a common retraction $r: B \to A$ for both $f$ and $g$. Then the diagram below left is a pullback in $\ce$ if and only if the diagram below right is an equalizer in $\ce$:
\[\begin{tikzcd}
	E && A \\
	&&& E && A && B \\
	A && B
	\arrow["e", from=1-1, to=1-3]
	\arrow["e"', from=1-1, to=3-1]
	\arrow["g", from=1-3, to=3-3]
	\arrow["e", from=2-4, to=2-6]
	\arrow["f", shift left=2, from=2-6, to=2-8]
	\arrow["g"', shift right=2, from=2-6, to=2-8]
	\arrow["f"', from=3-1, to=3-3]
\end{tikzcd}\]
In particular:
\begin{itemize}
\item any category with pullbacks admits coreflexive equalizers,
\item any pullback-preserving functor preserves coreflexive equalizers.
\end{itemize}
\end{lemma}
\begin{proof}
The direction „$\Rightarrow$“ is obvious (and we don't even need the pair to be coreflexive). The direction „$\Leftarrow$“ follows from the fact that in any commutative square like this:
\[\begin{tikzcd}
	X && A \\
	\\
	A && B
	\arrow["f"', from=3-1, to=3-3]
	\arrow["y", from=1-1, to=1-3]
	\arrow["g", from=1-3, to=3-3]
	\arrow["x"', from=1-1, to=3-1]
\end{tikzcd}\]
It is forced that $x=y$, which implies that $x$ equalizes $f,g$ and that there is a canonical map into the square on the left.
\end{proof}

\begin{theo}\label{THM_Tmulticats_comonadic}
Let $\ce$ be a category with pullbacks and $(T,m,i)$ a cartesian monad on $\ce$. The functor $F: \tmult \to \talgs$  from Theorem \ref{THM_talgs_tmult_adj} is comonadic.
\end{theo}
\begin{proof}
By (the dual of) the \textit{reflexive tripleability theorem} \cite[Proposition 5.5.8]{context}, this will follow if we show that:
\begin{itemize}
\item $\tmult$ admits coreflexive equalizers,
\item $F$ reflects isomorphisms,
\item $F$ preserves coreflexive equalizers.
\end{itemize}

It is straightforward to prove that $\tmult$ admits any limits that $\ce$ admits and $T: \ce \to \ce$ preserves. In particular this is true for pullbacks and so for coreflexive equalizers by Lemma \ref{THML_pbpresfun_pres_corefleq}.

The fact that $F$ reflects isomorphisms is implied by $T: \ce \to \ce$ reflecting isomorphisms. $T$ does so because since the diagram below is a pullback ($i: 1_\ce \Rightarrow T$ is cartesian), the fact that $Tu$ is an isomorphism implies that $u$ is an isomorphism since isomorphisms are stable under pullback:
\[\begin{tikzcd}
	{\cm_i} && {T\cm_i} \\
	\\
	{\cn_i} && {T\cn_i}
	\arrow["{Tf_i}", from=1-3, to=3-3]
	\arrow["{f_i}"', from=1-1, to=3-1]
	\arrow["{i_{\cm_i}}", from=1-1, to=1-3]
	\arrow["{i_{\cn_i}}"', from=3-1, to=3-3]
	\arrow["\lrcorner"{anchor=center, pos=0.125}, draw=none, from=1-1, to=3-3]
\end{tikzcd}\]

To prove that $F$ preserves coreflexive equalizers, let the following be an equalizer of a coreflexive pair in $\tmult$:
\[\begin{tikzcd}
	\co && \cm && \cn & {\text{ in } \tmult}
	\arrow["e", from=1-1, to=1-3]
	\arrow["f", shift left=2, from=1-3, to=1-5]
	\arrow["g"', shift right=2, from=1-3, to=1-5]
	\arrow["r"', curve={height=30pt}, from=1-5, to=1-3]
\end{tikzcd}\]

We wish to show that the following is an equalizer in $\talgs$:
\[\begin{tikzcd}
	{F\co} && {F\cm} && {F\cn}
	\arrow["Fe", from=1-1, to=1-3]
	\arrow["Ff", shift left=2, from=1-3, to=1-5]
	\arrow["Fg"', shift right=2, from=1-3, to=1-5]
	\arrow["Fr"', curve={height=30pt}, from=1-5, to=1-3]
\end{tikzcd}\]
This is equivalent to showing that this is an equalizer after applying the forgetful functor $U: \talgs \to \cat(\ce)$, which in turn is equivalent to saying these two diagrams are equalizers in $\ce$:
\[\begin{tikzcd}
	{(F\co)_i} && {(F\cm)_i} && {(F\cn)_i} & {i \in \{ 0,1 \}}
	\arrow["{(Fe)_i}", from=1-1, to=1-3]
	\arrow["{(Ff)_i}", shift left=2, from=1-3, to=1-5]
	\arrow["{(Fg)_i}"', shift right=2, from=1-3, to=1-5]
\end{tikzcd}\]

By the definition of $F$, these are just the diagrams:
\begin{equation}\label{EQP_TMTN_equalizer}
\begin{tikzcd}
	{T\co_i} && {T\cm_i} && {T\cn_i} & {i \in \{ 0,1 \}}
	\arrow["{Te_i}", from=1-1, to=1-3]
	\arrow["{Tf_i}", shift left=2, from=1-3, to=1-5]
	\arrow["{Tg_i}"', shift right=2, from=1-3, to=1-5]
\end{tikzcd}
\end{equation}
Since limits in $\tmult$ are computed pointwise, the following diagrams are equalizers in $\ce$:
\[\begin{tikzcd}
	{\co_i} && {\cm_i} && {\cn_i} & {i \in \{ 0,1 \}}
	\arrow["{e_i}", from=1-1, to=1-3]
	\arrow["{f_i}", shift left=2, from=1-3, to=1-5]
	\arrow["{g_i}"', shift right=2, from=1-3, to=1-5]
\end{tikzcd}\]
Since $T: \ce \to \ce$ preserves coreflexive equalizers by Lemma \ref{THML_pbpresfun_pres_corefleq}, the diagram \eqref{EQP_TMTN_equalizer} is an equalizer, which completes the proof.
\end{proof}

\begin{pr}\label{PR_doublecat_2monad_comonadicity}
For the free double category 2-monad \ref{PR_doublecat_2monad}, the comonadicity in Theorem \ref{THM_Tmulticats_comonadic} has been observed\footnote{To be more precise, its dual, involving ``$T$-comulticategories”, called \textit{oplax double categories} therein.} in \cite[Theorem 2.13.]{pathsindblcats}.
\end{pr}

\subsection{Multicategories and codescent}\label{SUBS_multicats_codescent}

In this section we sketch an alternative approach to 2-monads of form $\cat(T)$ using $T$-multicategories.

\begin{lemma}\label{THML_adjunction_descent}
Any adjunction as on the left in which the functor $J: \ce \to \cc$ is fully faithful induces an adjunction as on the right:
\[\begin{tikzcd}
	\cd & \ce & \cc & \rightsquigarrow & \cd && \ce
	\arrow["U"', from=1-1, to=1-2]
	\arrow["J"', from=1-2, to=1-3]
	\arrow[""{name=0, anchor=center, inner sep=0}, "F"', curve={height=24pt}, from=1-3, to=1-1]
	\arrow[""{name=1, anchor=center, inner sep=0}, "U"', from=1-5, to=1-7]
	\arrow[""{name=2, anchor=center, inner sep=0}, "FJ"', curve={height=24pt}, from=1-7, to=1-5]
	\arrow["\dashv"{anchor=center, rotate=-90}, draw=none, from=0, to=1-2]
	\arrow["\dashv"{anchor=center, rotate=-90}, draw=none, from=2, to=1]
\end{tikzcd}\]
\end{lemma}
\begin{proof}
This follows from the fact that for each $E \in \ob \ce$, $D \in \ob \cd$ we have:
\[
\cd(FJE, D) \cong \cc(JE, JUD) \cong \ce(E,UD),
\]
where the first bijection follows from the adjunction in the statement, and the other bijection follows from the full and faithfulness of $J$.
\end{proof}

Let $(T,m,i)$ be a cartesian monad on a category $\ce$ with pullbacks. Recall Theorem \ref{THM_talgs_tmult_adj} about the adjunction between algebras and multicategories and also the embedding of lax morphisms to multicategories from Lemma \ref{THML_ff_functor_talgl_multicat}:
\[\begin{tikzcd}
	\talgs &&&& \tmult \\
	&& \talgl
	\arrow[""{name=0, anchor=center, inner sep=0}, "U"', from=1-1, to=1-5]
	\arrow["J"', hook, from=1-1, to=2-3]
	\arrow[""{name=1, anchor=center, inner sep=0}, "F"', curve={height=24pt}, from=1-5, to=1-1]
	\arrow["I"', from=2-3, to=1-5]
	\arrow["\dashv"{anchor=center, rotate=-91}, draw=none, from=1, to=0]
\end{tikzcd}\]
Since $I$ is fully faithful and $F \dashv U$, by the above lemma we obtain $FI \dashv J$. In particular this gives an alternative proof to the strict algebra, lax morphism part of Proposition \ref{THMP_catT_strictification_exists}. It also unravels a connection between $T'$-multicategories and the category of corners construction from Section \ref{SEC_thecaseckCatE}:

\begin{cor}
Let $\ce$ be a category with pullbacks and let $T$ be a 2-monad of form $\cat(T')$ on $\cat(\ce)$. Given a $T$-algebra $(A,a)$, the category of corners $\cnr{A}$ associated to the domain-codiscrete category $\res{A,a}$ is the free $T$-algebra on the underlying $T'$-multicategory of $(A,a)$.
\end{cor}

\noindent Recall the definition of a \textit{lax-pie algebra} for a 2-monad from Definition \ref{DEFI_lax_semiflexible_flexible_pie} -- those are the $T$-algebras that admit the structure of a strict $Q_l$-coalgebras for the lax-morphism classifier 2-comonad on $\talgs$. For 2-monads of form $\cat(T)$, they can be characterized as follows.

\begin{theo}
Let $\ce$ be a category with pullbacks and $T$ a 2-monad on $\cat(\ce)$ of form $\cat(T)$. There is an equivalence between the category of $T$-multicategories and the category of lax-pie $T$-algebras:
\[\begin{tikzcd}
	\tmult && {Q_l\text{-Coalg}_s} \\
	\\
	& \talgs
	\arrow["\simeq", from=1-1, to=1-3]
	\arrow["F"', from=1-1, to=3-2]
	\arrow["{F^Q}", from=1-3, to=3-2]
\end{tikzcd}\]
Here $F$ is the free algebra on a multicategory functor from Theorem \ref{THM_talgs_tmult_adj} and $F^Q$ is the forgetful functor, sending a $Q_l$-coalgebra to its underlying $T$-algebra.
\end{theo}
\begin{proof}
By \textbf{definition}, the functor $F^Q$ is comonadic and the corresponding adjunction generates the comonad $Q_l$. By Theorem \ref{THM_Tmulticats_comonadic} and the observations after Lemma \ref{THML_adjunction_descent}, the same is true for $F$ -- hence the result.
\end{proof}

We will devote the rest of this section to giving an informal description of this equivalence. Denoting the counit and unit of the adjunction on the right side of \eqref{EQ_talgs_colaxtalgl_talgl} by $(q,p)$, the lax morphism classifier 2-comonad on $\talgs$ is the triple $(Q_l, p', q)$. Given a strict $T$-algebra $\mathbb{A} := (A,a)$, a $Q_l$-coalgebra structure consists of an algebra morphism $G: \mathbb{A} \to Q_l \mathbb{A}$ satisfying the $Q_l$-coalgebra axioms (in $\talgs$):
\begin{equation}\label{EQ_Ql_coalgebra_axioms}
\begin{tikzcd}
	{Q_l^2 \mathbb{A}} && {Q_l \mathbb{A}} && {Q_l \mathbb{A}} && {\mathbb{A}} \\
	\\
	{Q_l \mathbb{A}} && {\mathbb{A}} && {\mathbb{A}}
	\arrow["{p'_{\mathbb{A}}}"', from=1-3, to=1-1]
	\arrow["{q_{\mathbb{A}}}", from=1-5, to=1-7]
	\arrow["{G'}", from=3-1, to=1-1]
	\arrow["G"', from=3-3, to=1-3]
	\arrow["G", from=3-3, to=3-1]
	\arrow["G", from=3-5, to=1-5]
	\arrow[Rightarrow, no head, from=3-5, to=1-7]
\end{tikzcd}
\end{equation}

The functor $\tmult \to \qlcoalgs$ sends a multicategory $\cm$ to the $T$-algebra $F\cm$ (which can be observed to have a canonical $Q_l$-coalgebra structure). Its equivalence inverse sends a $Q_l$-coalgebra $(A,a,G)$ to its \textit{prime $T$-multicategory} that we will denote by $A^p$. The object of objects is the following pullback in $\ce$:
\[\begin{tikzcd}
	{A^p_0} && {A_0} \\
	\\
	{A_0} && {\cnr{\mathbb{A}}_0 = TA_0}
	\arrow["{\iota^p_0}", from=1-1, to=1-3]
	\arrow["{\iota^p_0}"', from=1-1, to=3-1]
	\arrow["\lrcorner"{anchor=center, pos=0.125}, draw=none, from=1-1, to=3-3]
	\arrow["{i_{\ca_0}}", from=1-3, to=3-3]
	\arrow["{G_0}"', from=3-1, to=3-3]
\end{tikzcd}\]

Since both 1-cells $G_0, i_{A_0}$ have a common section $a_0: T\ca_0 \to \ca_0$, the pullback projections are equal. In particular, it is an equalizer and so $A_0^p$ is a subobject of $A_0$ -- let us call it the \textit{subobject of prime elements} of the coalgebra $(A,a,G)$. The object of morphisms of $A^p$ will consist of ``morphisms with prime target”:
\[\begin{tikzcd}
	{\ca^p_1} && {\ca_0^p} \\
	\\
	{\ca_1} && {\ca_0}
	\arrow["\lrcorner"{anchor=center, pos=0.125}, draw=none, from=1-1, to=3-3]
	\arrow["{\iota_p}", from=1-3, to=3-3]
	\arrow["{\iota^p_1}"', from=1-1, to=3-1]
	\arrow["t"', from=3-1, to=3-3]
	\arrow["{t^p}", from=1-1, to=1-3]
\end{tikzcd}\]

\noindent The underlying $T$-multigraph is then given by the span:
\[\begin{tikzcd}
	{TA_0^p} && {A^p_1} && {A^p_0}
	\arrow["{\overline{G_0}\circ s \circ \iota_1^p}"', from=1-3, to=1-1]
	\arrow["{t^p}", from=1-3, to=1-5]
\end{tikzcd}\]
where $\overline{G_0}: \ca_0 \to T\ca^p_0$ is the canonical map into pullback:
\[\begin{tikzcd}
	{\ca_0} \\
	& {T\ca_0^p} && {T\ca_0} \\
	\\
	& {T\ca_0} && {T^2\ca_0}
	\arrow["{T\iota_p}"', from=2-2, to=4-2]
	\arrow[from=2-2, to=2-4]
	\arrow["{Ti_{\ca_0}}", from=2-4, to=4-4]
	\arrow["{TG_0}"', from=4-2, to=4-4]
	\arrow["{\exists ! \overline{G_0}}"', from=1-1, to=2-2]
	\arrow["{G_0}"', curve={height=18pt}, from=1-1, to=4-2]
	\arrow["{G_0}", curve={height=-12pt}, from=1-1, to=2-4]
	\arrow["\lrcorner"{anchor=center, pos=0.125}, draw=none, from=2-2, to=4-4]
\end{tikzcd}\]

\subsection{Examples}

\noindent In the example below, we use a similar reasoning to that of \cite[Prologue]{pathsindblcats} -- from this paper we also borrow the \textit{prime object} terminology. We mention the results of \cite{pathsindblcats} in Example \ref{PR_doublecat_2monad_laxpie}.

\begin{pr}\label{PR_freestrictmoncat2monad_laxpie}
Consider the free strict monoidal category 2-monad $T$ on $\cat$ from Example \ref{PR_freestrictmoncat2monad}. Since Example \ref{PR_freestrictmoncat2monad_CaTT}, we know it is of form $\cat(T)$ for the free monoid monad $T$ on $\set$.

We know the explicit description of $p,q,Q_l$ from Example \ref{PR_freestrictmoncat2monad_coherence}. Given a lax monoidal functor $(F,\overline{F}): (\ca,\otimes,I) \to (\cb, \otimes', I')$, applying the strictification to it produces a strict monoidal functor:
\[
\cnr{F,\overline{F}} : \cnr{\ca} \to \cnr{\cb},
\]
that is given by sending the corner below left to the one below right:
\[\begin{tikzcd}[column sep=tiny]
	{(a_1, a_2, a_3)} && {(I, a_1 \otimes a_2, a_3)} & {(Fa_1,Fa_2,Fa_3)} && {(I', Fa_1 \odot Fa_2, Fa_3)} \\
	&&&&& {(FI,F(a_1\otimes a_2), Fa_3)} \\
	&& {(b_1, b_2, b_3)} &&& {(Fb_1, Fb_2, Fb_3)}
	\arrow["{((),(a_1, a_2),(a_3))}", shift left=2, from=1-1, to=1-3]
	\arrow["{(f_1,f_2,f_3)}"', from=1-3, to=3-3]
	\arrow["{((),(Fa_1, Fa_2),(Fa_3))}", shift left=2, from=1-4, to=1-6]
	\arrow["{(\overline{F}_{()},\overline{F}_{a_1,a_2},1_{Fa_3})}"', from=1-6, to=2-6]
	\arrow["{(Ff_1,Ff_2,Ff_3)}"', from=2-6, to=3-6]
\end{tikzcd}\]

Fix a strict monoidal category $(\ca,\otimes,I)$. In particular, the comultiplication for the 2-comonad is a strict monoidal functor $\cnr{P,\overline{P}}: \cnr{\ca} \to \cnr{\cnr{\ca}}$ that would send the corner above left to the following morphism in $\cnr{\cnr{\ca}}$:
\[\begin{tikzcd}
	{[(a_1), (a_2), (a_3)]} &&& {[(),(a_1,a_2),(a_3)]} \\
	&&& {[(I),(a_1\otimes a_2), (a_3)]} \\
	&&& {[(b_1), (b_2), (b_3)]}
	\arrow["{[[\phantom{.}],[(a_1),(a_2)],[(a_3)]]}", from=1-1, to=1-4]
	\arrow["{[\overline{P}_{()},\overline{P}_{a_1,a_2},1_{(a_3)}]}"', from=1-4, to=2-4]
	\arrow["{[(f_1), (f_2), (f_3)]}"', from=2-4, to=3-4]
\end{tikzcd}\]

Assume that the strict monoidal category $(\ca,\otimes,I)$ admits a $Q_l$-coalgebra structure $G: \ca \to \cnr{\ca}$, and call an object $x \in \ob \ca$ \textit{prime} if $Gx = (x)$. The coalgebra counit axiom \eqref{EQ_Ql_coalgebra_axioms} tells us that every object $a \in \ob \ca$ can be written as $a_1 \otimes \dots \otimes a_n$ for a tuple $Ga = (a_1, \dots, a_n)$, and the coassociativity axiom tells us that each of these $a_i$'s is in fact prime. So every object admits a prime factorization. This factorization is unique since if:
\[
b_1 \otimes \dots \otimes b_m = a_1 \otimes \dots \otimes a_n,
\]
for $b_i$, $a_j$ prime, we have:
\begin{align*}
(b_1, \dots, b_m) &= (b_1) \boxplus \dots \boxplus (b_m)\\
&\sstackrel{(*)}{=} G(b_1 \dots b_m) \\
&= G(a_1 \dots a_n) \\
&= (a_1, \dots, a_n).
\end{align*}
Here $(*)$ is true since each $b_i$ is prime and $G$ is a strict monoidal functor. Next, let us focus on morphisms in $\ca$. If $f: a \to b \in \mor \ca$ is a morphism in which $b$ is prime, notice that by the definition of morphisms $\cnr{\ca}$ and the counit axiom \eqref{EQ_Ql_coalgebra_axioms}, $G$ ``does nothing” to $f$:
\[\begin{tikzcd}
	& a && {(a_1, \dots, a_n)} && {(a)} \\
	{G:} && \mapsto &&&& {\text{ if } b \text{ prime}} \\
	& b &&&& {(b)}
	\arrow["f"', from=1-2, to=3-2]
	\arrow["{((a_1,\dots,a_n))}", from=1-4, to=1-6]
	\arrow["{(f)}"', from=1-6, to=3-6]
\end{tikzcd}\]
Now, the \textit{prime multicategory} $\ca^p$ is defined such that the objects are the prime objects of $\ca$, and morphisms $a_1, \dots, a_n \to b$ are precisely the morphisms:
\[
a_1 \otimes \dots \otimes a_n \to b \hspace{0.5cm}\text{in } \ca.
\]
It is also clear that $(\ca, \otimes, I)$ is the free monoidal category on the prime multicategory $\ca^p$.
\end{pr}

\begin{pr}\label{PR_doublecat_2monad_laxpie}
Recall the double category 2-monad from \ref{PR_doublecat_2monad} (see also its coherence result in Example \ref{prcodlax-dblcats} and comonadicity result in Example \ref{PR_doublecat_2monad_comonadicity}).

Using similar arguments as those in Example \ref{PR_freestrictmoncat2monad_laxpie}, one may observe that a double category $X$ is lax-pie if and only if there is an underlying virtual double category $X^p$ and $X$ is a free double category on $X^p$. The dual of this result has been observed in \cite[Theorem 2.13.]{pathsindblcats}.
\end{pr}

\begin{pr}
Given a 1-category $\cj$ and considering the reslan 2-monad $T$ on the 2-category $[ \ob \cj, \cat]$, a 2-functor $P: \cj \to \cat$ can be seen to be lax-pie if and only if there is a functor $F: \ce \to \cj$ and $P$ is isomorphic to the comma category functor (see also Example \ref{PR_reslan_multicats}):
\[
j \mapsto F \downarrow j.
\]
\end{pr}

\begin{pr}[2-categories]
For the 2-category 2-monad (see Example \ref{PR_2category_2monad}, Example \ref{PR_opetopes}) it can be seen that a 2-category $\cc$ is lax-pie if and only if there is an underlying opetope $\cc^p$, and $\cc$ is the free 2-category on it.

\end{pr}

\begin{pozn}[Comparison with pie algebras]
In \cite[Theorem 3]{onsemiflexible}, \textbf{pie} algebras (recalled in Section \ref{SUBS_4_lax_flexibility}) for a suitable class of 2-monads $T$ have been characterized as those $T$-algebras that are ``free on the level of objects”. In particular:
\begin{itemize}
\item a strict monoidal category is pie $\iff$ the underlying monoid on objects is free,
\item a 2-functor $F: \cj \to \cat$ is pie $\iff$ the composite $\ob \circ F: \cj \to \set$ is a free $\text{reslan}$-algebra, i.e. a coproduct of representable presheaves on $\cj$,
\item a 2-category is pie $\iff$ the underlying 1-category is free on a graph.
\end{itemize}
We see from our examples that ``being lax-pie” means in particular ``being pie” together with a condition on the morphisms -- it is a strictly stronger condition. This is always the case since by the general theory (Remark \ref{POZN_lax_semiflexible_is_semiflexible}) we have observed that every lax-pie algebra is pie.
\end{pozn}

%% file: PhD_CH_7_Concluding.tex
\chapter{Open questions}\label{CH_concludingremarks}

Let us finish with ideas for further development. We begin with concrete questions not addressed in the thesis, but for which the author believes there is a straightforward answer:

\begin{itemize}
\item We have seen that weak factorization systems are equivalent to certain double categories (Corollary \ref{THMC_wfs_equiv_to_dbl}). What elementary double categorical properties characterize them?
\item What are the two-dimensional aspects of theorems  \ref{THM_fundamentaladj}, \ref{THM_coKZfication}, \ref{THM_Tmulticats_comonadic}?
\item What is the generalization of Example \ref{PR_reslan_laxcoh_explicit} to the setting where $J$ is a 2-category?
\item What are the examples of lax-semiflexible, lax-flexible algebras for the 2-monads of form $\cat(T)$?
\item What is an example of a 2-monad $T$ for which the coherence for colax algebras does \textbf{NOT} hold? Is it given by a lax analogue of \cite[Example 3.1]{codobjcoh}?
\item Prove that for 2-monads of form $\cat(T')$, lax-pie $T$-algebras that are \textit{continuous}\footnote{The $Q_l$-coalgebra structure map, which is a left adjoint to the counit of $Q_l$, has a further left adjoint} correspond to \textit{representable $T'$-multicategories}, and as such, are equivalent to pseudo $T$-algebras. For the free strict monoidal category 2-monad, this recovers the equivalence between representable multicategories and monoidal categories proven in \cite{reprmulticats}.
\end{itemize}

\noindent Less specific questions that the author intends to pursue in the future include:
\begin{itemize}
\item What can be said about colimits in Kleisli 2-categories for lax-idempotent pseudomonads?
\item What is the relationship between codomain-colax categories and crossed double categories? Some connections have been sketched in Subsection \ref{SUBS_grothendieck}.
\item What are colax algebras for lax-idempotent pseudomonads? Is it always a comonad together with a ``property-like structure”, like it was in Example \ref{PR_freeinitialobject}?
\item What are the lax-pie algebras for more general classes of 2-monads? For the free symmetric strict monoidal category 2-monad (Example \ref{PR_freesymstrictmoncat2monad}), are they equivalent to symmetric multicategories?
\item How to generalize Section \ref{SEC_coKZfication} and turn a wider class of 2-monads into colax-idempotent ones?
\end{itemize}

%% file: PhD_CH_A_Further_remarks.tex
\chapter{A generalization of coherence for colax algebras}\label{CH_generalized_colax_coh}

In this chapter we mention a possible generalization of Theorem \ref{THM_lax_coherence_colax_algebras_EXPOSITION} for which we have no examples. John Bourke observed that the coherence theorem for pseudo algebras holds if and only if $T$ preserves the relevant (iso)-codescent objects as a bicolimit (mentioned in \cite[page 11]{coherencethreedim_HIST}, the details may be found in author's Master's thesis \cite[Proposition 3.1]{milomgr}). A first step in studying the lax analogue of this phenomenon would be to replace preservation as a bicolimit by something weaker. It seems that (right adjoint)-colimits (Definition \ref{DEFI_rali_colimit}) that satisfy an additional property work well for this task:

\begin{defi}
A (right adjoint)-limit of $F$ weighted by $W$ is said to be \textit{coherent} if the collection of right adjoints to each $\kappa_A$ is 2-natural in $A$, thus forming an adjunction in $[\ck^{op},\cat]$:
\[\begin{tikzcd}
	{\ck(-,L)} && {[\cj, \cat](W, \ck(-, F?))}
	\arrow[""{name=0, anchor=center, inner sep=0}, "\kappa"', curve={height=24pt}, from=1-1, to=1-3]
	\arrow[""{name=1, anchor=center, inner sep=0}, curve={height=24pt}, from=1-3, to=1-1]
	\arrow["\dashv"{anchor=center, rotate=-90}, draw=none, from=1, to=0]
\end{tikzcd}\]
We say that a 2-functor $H: \ck \to \cl$ prerves a (right adjoint)-limit \textit{coherently} if the composite $H \circ \lambda$ is a coherent colimit of $HF$ weighted by $W$.
\end{defi}

\noindent The generalization of Theorem \ref{THM_lax_coherence_colax_algebras_EXPOSITION} now proceeds as follows:

\begin{theo}[Generalized coherence for colax algebras]
Assume $T$ is a 2-monad on a 2-category $\ck$. Assume that $\talgs$ admits codescent objects of resolutions of colax algebras.  Assume the forgetful 2-functor $U: \talgs \to \ck$ coherently preserves these codescent objects as a (left adjoint)-colimit. Then for every colax $T$-algebra $(A,a,\gamma,\iota)$, the unit of the 2-adjunction between $\talgs$ and $\colaxtalgl$ has a left adjoint in $\ck$.
\end{theo}
\begin{proof}
Notice first that the (right adjoint)-codescent object from Lemma \ref{POZNBAZANT} is coherent. The proof now proceeds  exactly the same as the proof of Theorem \ref{THM_lax_coherence_colax_algebras_EXPOSITION}, except that in the diagram \eqref{EQ_PF_kAprime_cohK_kA}, there will be an adjunction (guaranteed by the preservation as a (right adjoint)-colimit) instead of an isomorphism.
\end{proof}

%% file: PhD_CH_A_AUXILIARY_LEMMAS.tex
\chapter{Auxiliary lemmas for Chapter 5}\label{CH_APPENDIX_AUXILIARY}

\begin{lemma}\label{THML_epsilon_phi_second_swallowtail}
In the proof of Theorem \ref{THM_Bunge}:
\begin{itemize}
\item $\epsilon$ is colax natural,
\item $\Phi$ is a modification,
\item the second swallowtail identity.
\end{itemize}
\end{lemma}
\begin{proof}
In this proof, we will reference the defining equations for $\gamma', \iota', F\alpha, \epsilon, \Psi$ above the equals sign. In the unlabeled equations we use the middle-four interchange rule combined with the pseudofunctor laws.

\textbf{$\epsilon$ is colax natural}: The \textbf{composition axiom} amounts to proving the equality of the following 2-cells:
\[\begin{tikzcd}
	FUB & FUB & FUB && B & FUB && B \\
	&& FUC && C &&& C \\
	FUD & FUD & FUD && D & FUD && D
	\arrow[""{name=0, anchor=center, inner sep=0}, "{\epsilon_B}", from=1-3, to=1-5]
	\arrow["h", from=1-5, to=2-5]
	\arrow["g", from=2-5, to=3-5]
	\arrow["FUh", from=1-3, to=2-3]
	\arrow["FUg", from=2-3, to=3-3]
	\arrow[""{name=1, anchor=center, inner sep=0}, "{\epsilon_D}"', from=3-3, to=3-5]
	\arrow[""{name=2, anchor=center, inner sep=0}, "{\epsilon_C}"{description}, from=2-3, to=2-5]
	\arrow[""{name=3, anchor=center, inner sep=0}, "{F(Ug \circ Uh)}"{description, pos=0.8}, from=1-2, to=3-2]
	\arrow[""{name=4, anchor=center, inner sep=0}, "{FU(gh)}"', from=1-1, to=3-1]
	\arrow[Rightarrow, no head, from=3-1, to=3-2]
	\arrow[Rightarrow, no head, from=3-2, to=3-3]
	\arrow[Rightarrow, no head, from=1-1, to=1-2]
	\arrow[Rightarrow, no head, from=1-2, to=1-3]
	\arrow["{FU(gh)}"{description}, from=1-6, to=3-6]
	\arrow[""{name=5, anchor=center, inner sep=0}, "{\epsilon_B}", from=1-6, to=1-8]
	\arrow["h", from=1-8, to=2-8]
	\arrow["g", from=2-8, to=3-8]
	\arrow[""{name=6, anchor=center, inner sep=0}, "{\epsilon_D}"', from=3-6, to=3-8]
	\arrow["{\epsilon_h}"', shorten <=4pt, shorten >=4pt, Rightarrow, from=2, to=0]
	\arrow["{\epsilon_g}"', shorten <=4pt, shorten >=4pt, Rightarrow, from=1, to=2]
	\arrow["{\epsilon_{gh}}"', shorten <=9pt, shorten >=9pt, Rightarrow, from=6, to=5]
	\arrow["{\gamma'}", shorten <=5pt, Rightarrow, from=3, to=2-3]
	\arrow["F\gamma", shorten <=7pt, shorten >=7pt, Rightarrow, from=4, to=3]
\end{tikzcd}\]
It is enough to prove these after applying $U(-) \circ \gamma y_{UB} \circ U\epsilon_D \mbbd \circ \Phi_D U(gh)$ on both sides. We then have:
\begin{align*}
& U(LHS) \circ \gamma^{-1} y_{UB} \circ U\epsilon_D \mbbd_{y_{UC} U(gh)} \circ \Phi_D U(gh) =
\phantom{.}
\phantom{.}
\phantom{.}\\
&= U(g\epsilon_h \circ \epsilon_g FUh \circ \epsilon_D \gamma')y_{UB} \circ \gamma^{-1}y_{UB} \circ U\epsilon_D UF\gamma y_{UB} \circ U\epsilon_D \mbbd_{y_{UC} U(gh)} \circ \Phi_D U(gh)
\phantom{.}
\phantom{.}
\phantom{.}\\
&\sstackrel{\eqref{EQP_Falpha}}{=} U(g\epsilon_h \circ \epsilon_g FUh \circ \epsilon_D \gamma')y_{UB} \circ \gamma^{-1} y_{UB} \circ U\epsilon_D \mbbd_{y_{UC} Ug Uh} \circ \Phi_D Ug Uh \circ \gamma
\phantom{.}
\phantom{.}
\phantom{.}\\
&= U(g\epsilon_h \circ \epsilon_g FUh)y_{UB} \circ \gamma^{-1} y_{UB} \circ U\epsilon_D U\gamma' y_{UB} \circ U\epsilon_D \mbbd_{y_{UC} Ug Uh} \circ \Phi_D Ug Uh \circ \gamma
\phantom{.}
\phantom{.}
\phantom{.}\\
&\sstackrel{\eqref{EQP_gamma_prime}}{=} U(g\epsilon_h  \circ \epsilon_g FUh)y_{UB} \circ \gamma^{-1} y_{UB} \circ U\epsilon_D \gamma^{-1} y_{UB} \circ U\epsilon_g UFUg \mbbd_{y_{UC} Uh} \circ\\
&\phantom{=}\circ U\epsilon_D \mbbd_{y_{UD} Ug} Uh \circ \Phi_D Ug Uh \circ \gamma
\phantom{.}
\phantom{.}
\phantom{.}\\
&= U(g\epsilon_h)y_{UB} \circ \gamma^{-1} y_{UB} \circ U\epsilon_g UFUh y_{UB} \circ \gamma^{-1} UFUh y_{UB} \circ U\epsilon_D UFUg \mbbd_{y_{UC} Uh} \circ \\
&\phantom{=}\circ U\epsilon_D \mbbd_{y_{UD} Ug} Uh \circ \Phi_D Ug Uh \circ \gamma 
\phantom{.}
\phantom{.}
\phantom{.}\\
&\sstackrel{\eqref{EQP_Phi_epsilon}}{=} U(g \epsilon_h) y_{UB} \circ \gamma^{-1} y_{UB} \circ \gamma^{-1} UFUh y_{UB} \circ Ug U\epsilon_B \mbbd_{y_{UC} Uh} \circ Ug \Phi_B Uh \circ \gamma
\phantom{.}
\phantom{.}
\phantom{.}\\
&= \gamma^{-1} y_{UB} \circ Ug U\epsilon_h y_{UB} \circ Ug \gamma^{-1} y_{UB} \circ Ug U\epsilon_B \mbbd_{y_{UC} Uh} \circ Ug \Phi_B Uh \circ \gamma
\phantom{.}
\phantom{.}
\phantom{.}\\
&\sstackrel{\eqref{EQP_Phi_epsilon}}{=} \gamma^{-1}y_{UB} \circ Ug \gamma^{-1}y_{UB} \circ Ug Uh \Phi_B \circ \gamma
\phantom{.}
\phantom{.}
\phantom{.}\\
&= \gamma^{-1} y_{UB} \circ \gamma^{-1} U\epsilon_B y_{UB} \circ Ug Uh \Phi_B \circ \gamma
\phantom{.}
\phantom{.}
\phantom{.}\\
&\sstackrel{\eqref{EQP_Phi_epsilon}}{=} U(RHS) \circ \gamma^{-1} y_{UB} \circ U\epsilon_D \mbbd_{y_{UC}U(gh)} \circ \Phi_D U(gh).
\end{align*}

The \textbf{unit axiom} for $\epsilon$ amounts to showing that:
\[\begin{tikzcd}
	FUB && B \\
	&&& {=} & FUB && FUB && B \\
	FUB && B
	\arrow[""{name=0, anchor=center, inner sep=0}, curve={height=-24pt}, Rightarrow, no head, from=2-5, to=2-7]
	\arrow[""{name=1, anchor=center, inner sep=0}, "{FU1_B}"', curve={height=24pt}, from=2-5, to=2-7]
	\arrow[""{name=2, anchor=center, inner sep=0}, "{F1_{UB}}"{description}, from=2-5, to=2-7]
	\arrow["{\epsilon_B}", from=2-7, to=2-9]
	\arrow[""{name=3, anchor=center, inner sep=0}, "{\epsilon_B}", from=1-1, to=1-3]
	\arrow["{FU1_B}"', from=1-1, to=3-1]
	\arrow[Rightarrow, no head, from=1-3, to=3-3]
	\arrow[""{name=4, anchor=center, inner sep=0}, "{\epsilon_B}"', from=3-1, to=3-3]
	\arrow["{\epsilon_{1_B}}"', shorten <=9pt, shorten >=9pt, Rightarrow, from=4, to=3]
	\arrow["F\iota"', shorten <=3pt, shorten >=3pt, Rightarrow, from=1, to=2]
	\arrow["{\iota'}"', shorten <=3pt, shorten >=3pt, Rightarrow, from=2, to=0]
\end{tikzcd}\]
It suffices to prove that these 2-cells are equal after applying the 2-cell\\ $U(-) \circ \gamma^{-1} y_{UB} \circ U\epsilon_D \mbbd_{y_{UB} U1_B} \circ \Phi_D U1_B$ on both sides. This is done as follows:

\[
\adjustbox{scale=0.90,center}{
\begin{tikzcd}
	{U\epsilon_B UFU1_B y_{UB}} & {U(\epsilon_B FU1_B)y_{UB}} && {U(\epsilon_B F1_{UB})y_{UB}} \\
	& {U\epsilon_B UF1_B y_{UB}} & {U\epsilon_B U1_{FUB} y_{UB}} & {U\epsilon_B y_{UB}} \\
	{U\epsilon_B y_{UB}U1_B } &&& {U\epsilon_B y_{UB}} \\
	& {U1_B U\epsilon_B y_{UB}} & {\eqref{EQP_Phi_epsilon}} \\
	{U1_B} & {U\epsilon_B y_{UB} U1_B} & {U\epsilon_B UFU1_B y_{UB}} & {U(\epsilon_B \circ FU1_B)y_{UB}}
	\arrow["{\gamma^{-1}y_{UB}}", from=1-1, to=1-2]
	\arrow["{U\epsilon_B UF\iota y_{UB}}"{description}, from=1-1, to=2-2]
	\arrow["{U(\epsilon_B F\iota)y_{UB}}", from=1-2, to=1-4]
	\arrow["{U(\epsilon_B \iota') y_{UB}}", from=1-4, to=2-4]
	\arrow["{\gamma^{-1}y_{UB}}"{description}, curve={height=-12pt}, from=2-2, to=1-4]
	\arrow["{U\epsilon_B U\iota' y_{UB}}", from=2-2, to=2-3]
	\arrow["{\eqref{EQP_Falpha}}"{description}, draw=none, from=2-2, to=3-1]
	\arrow["{\gamma^{-1}y_{UB}}", from=2-3, to=2-4]
	\arrow["{U\epsilon_B \mbbd_{y_{UB} U1_B}}", from=3-1, to=1-1]
	\arrow["{U\epsilon_B y_{UB} \iota}"{description, pos=0.3}, from=3-1, to=3-4]
	\arrow[""{name=0, anchor=center, inner sep=0}, "{U\epsilon_B \mbbd_{y_{UB}}}"{pos=0.7}, from=3-4, to=2-2]
	\arrow[""{name=1, anchor=center, inner sep=0}, "{U\epsilon_B \iota^{-1}y_{UB}}"', from=3-4, to=2-3]
	\arrow["1"', from=3-4, to=2-4]
	\arrow["{\gamma^{-1} y_{UB}}"{description}, from=4-2, to=3-4]
	\arrow["{\Phi_B U1_B}", from=5-1, to=3-1]
	\arrow["{U1_B \Phi_B}", from=5-1, to=4-2]
	\arrow["{\Phi_B U1_B}"', from=5-1, to=5-2]
	\arrow["{U\epsilon_B \mbbd_{y_{UB} U1_B}}"', shift right=2, from=5-2, to=5-3]
	\arrow["{\gamma^{-1} y_{UB}}"', from=5-3, to=5-4]
	\arrow["{U\epsilon_{1_B}y_{UB}}"', from=5-4, to=3-4]
	\arrow["{\eqref{EQP_iota_prime}}"'{pos=1}, draw=none, from=1, to=0]
\end{tikzcd}
}
\]

The \textbf{local naturality} for $\epsilon$ amounts to showing that the 2-cells below are equal:
\[\begin{tikzcd}
	FUB & FUB && B & FUB && B & B \\
	FUC & FUC && C & FUC && C & C
	\arrow[""{name=0, anchor=center, inner sep=0}, "FUk", from=1-2, to=2-2]
	\arrow[""{name=1, anchor=center, inner sep=0}, "{\epsilon_B}", from=1-2, to=1-4]
	\arrow[""{name=2, anchor=center, inner sep=0}, "{\epsilon_C}"', from=2-2, to=2-4]
	\arrow["k", from=1-4, to=2-4]
	\arrow[Rightarrow, no head, from=1-1, to=1-2]
	\arrow[Rightarrow, no head, from=2-1, to=2-2]
	\arrow[""{name=3, anchor=center, inner sep=0}, "FUh"', from=1-1, to=2-1]
	\arrow[""{name=4, anchor=center, inner sep=0}, "{\epsilon_B}", from=1-5, to=1-7]
	\arrow[""{name=5, anchor=center, inner sep=0}, "h"', from=1-7, to=2-7]
	\arrow["FUh"', from=1-5, to=2-5]
	\arrow[""{name=6, anchor=center, inner sep=0}, "{\epsilon_C}"', from=2-5, to=2-7]
	\arrow[Rightarrow, no head, from=1-7, to=1-8]
	\arrow[Rightarrow, no head, from=2-7, to=2-8]
	\arrow[""{name=7, anchor=center, inner sep=0}, "k", from=1-8, to=2-8]
	\arrow["FU\alpha"', shorten <=7pt, shorten >=7pt, Rightarrow, from=3, to=0]
	\arrow["{\epsilon_h}"', shorten <=4pt, shorten >=4pt, Rightarrow, from=6, to=4]
	\arrow["\alpha", shorten <=6pt, shorten >=6pt, Rightarrow, from=5, to=7]
	\arrow["{\epsilon_k}", shorten <=4pt, shorten >=4pt, Rightarrow, from=2, to=1]
\end{tikzcd}\]

An analogous approach will be done here as well, this time pre-composing with the 2-cell $U(-) \circ \gamma^{-1} y_{UB} \circ U\epsilon_D \mbbd_{y_{UC} Uh} \circ \Phi_D Uh$:
\[
\adjustbox{scale=0.9,center}{
\begin{tikzcd}
	Uh && {U\epsilon_C y_{UC} Uh} & {U\epsilon_C UFUh y_{UB}} & {U(\epsilon_C FUh)y_{UB}} \\
	& Uk & {U\epsilon_C y_{UC} Uk} & {U\epsilon_C UFUky_{UB}} & {U(\epsilon_C FUk)y_{UB}} \\
	{U\epsilon_C y_{UC} Uh} & {Uh U\epsilon_B y_{UB}} && {Uk U\epsilon_B y_{UB}} \\
	{U\epsilon_C UFUh y_{UB}} & {U(\epsilon_C FUh)} & {U(h \epsilon_B) y_{UB}} && {U(k \epsilon_B)y_{UB}}
	\arrow["{\Phi_C Uh}", from=1-1, to=1-3]
	\arrow["U\alpha"{description}, from=1-1, to=2-2]
	\arrow["{\Phi_D Uh}"', from=1-1, to=3-1]
	\arrow["{Uh\Phi_B}"', from=1-1, to=3-2]
	\arrow[""{name=0, anchor=center, inner sep=0}, "{U\epsilon_C \mbbd_{y_{UC} Uh}}", shift left=2, from=1-3, to=1-4]
	\arrow["{U\epsilon_C y_{UC} U\alpha}"', from=1-3, to=2-3]
	\arrow["{\gamma^{-1} y_{UB}}", from=1-4, to=1-5]
	\arrow["{U\epsilon_C UFU\alpha y_{UB}}", from=1-4, to=2-4]
	\arrow["{U(\epsilon_C FU\alpha)y_{UB}}", from=1-5, to=2-5]
	\arrow["{\Phi_C Uk}", from=2-2, to=2-3]
	\arrow["{Uk\Phi_B}"{description}, from=2-2, to=3-4]
	\arrow[""{name=1, anchor=center, inner sep=0}, "{U\epsilon_D \mbbd_{y_{UC} Uk}}"', shift right=2, from=2-3, to=2-4]
	\arrow["{\gamma^{-1} y_{UB}}"', from=2-4, to=2-5]
	\arrow["{\eqref{EQP_Phi_epsilon}}"{description}, draw=none, from=2-4, to=3-4]
	\arrow["{U\epsilon_k y_{UB}}", from=2-5, to=4-5]
	\arrow["{U\epsilon_D \mbbd_{y_{UC} Uh}}"', from=3-1, to=4-1]
	\arrow["{\eqref{EQP_Phi_epsilon}}"{description}, draw=none, from=3-2, to=3-1]
	\arrow["{U\alpha U\epsilon_B y_{UB}}"', from=3-2, to=3-4]
	\arrow["{\gamma^{-1} y_{UB}}"{description, pos=0.4}, from=3-2, to=4-3]
	\arrow["{\gamma^{-1} y_{UB}}"{description}, from=3-4, to=4-5]
	\arrow["{\gamma^{-1}y_{UB}}"', from=4-1, to=4-2]
	\arrow["{U\epsilon_h y_{UB}}"', from=4-2, to=4-3]
	\arrow["{U(\alpha \epsilon_B)y_{UB}}"', from=4-3, to=4-5]
	\arrow["{\eqref{EQP_Falpha}}"{description}, draw=none, from=0, to=1]
\end{tikzcd}
}\]

\textbf{$\Psi$ is a modification}:
This amounts to showing that these 2-cells are equal:
\[\begin{tikzcd}[column sep=scriptsize]
	& {} & {} \\
	FA & FUFA & {} & FA & FA & FB && FA \\
	\\
	& FUFB && FB && FUFB && FB
	\arrow[shift left=3, curve={height=-24pt}, Rightarrow, no head, from=2-1, to=2-4]
	\arrow[curve={height=-24pt}, Rightarrow, no head, from=2-5, to=2-8]
	\arrow["{Fy_A}", from=2-1, to=2-2]
	\arrow["FUFf", from=2-2, to=4-2]
	\arrow[""{name=0, anchor=center, inner sep=0}, "{F(y_Bf)}"', curve={height=30pt}, from=2-1, to=4-2]
	\arrow[""{name=0p, anchor=center, inner sep=0}, phantom, from=2-1, to=4-2, start anchor=center, end anchor=center, curve={height=30pt}]
	\arrow[""{name=1, anchor=center, inner sep=0}, "{\epsilon_{FA}}"{description}, from=2-2, to=2-4]
	\arrow[""{name=2, anchor=center, inner sep=0}, "{\epsilon_{FB}}"', from=4-2, to=4-4]
	\arrow["Ff", from=2-4, to=4-4]
	\arrow[""{name=3, anchor=center, inner sep=0}, "{F(y_B f)}"{description, pos=0.7}, curve={height=30pt}, from=2-5, to=4-6]
	\arrow["Ff"{description}, from=2-5, to=2-6]
	\arrow["{Fy_B}", from=2-6, to=4-6]
	\arrow["{\epsilon_{FB}}"', from=4-6, to=4-8]
	\arrow["Ff", from=2-8, to=4-8]
	\arrow[""{name=4, anchor=center, inner sep=0}, curve={height=-12pt}, Rightarrow, no head, from=2-6, to=4-8]
	\arrow[""{name=5, anchor=center, inner sep=0}, from=2-1, to=4-2]
	\arrow[""{name=5p, anchor=center, inner sep=0}, phantom, from=2-1, to=4-2, start anchor=center, end anchor=center]
	\arrow[""{name=6, anchor=center, inner sep=0}, draw=none, from=1-3, to=1-2]
	\arrow[""{name=6p, anchor=center, inner sep=0}, phantom, from=1-3, to=1-2, start anchor=center, end anchor=center]
	\arrow[""{name=7, anchor=center, inner sep=0}, draw=none, from=2-3, to=2-2]
	\arrow[""{name=7p, anchor=center, inner sep=0}, phantom, from=2-3, to=2-2, start anchor=center, end anchor=center]
	\arrow["{\Psi_B}"', shorten >=7pt, Rightarrow, from=4-6, to=4]
	\arrow["{Fy_f}"', shorten <=5pt, shorten >=5pt, Rightarrow, from=0p, to=5p]
	\arrow["{\gamma'}", shorten <=4pt, Rightarrow, from=5, to=2-2]
	\arrow["{\epsilon_{Ff}}"', shorten <=9pt, shorten >=9pt, Rightarrow, from=2, to=1]
	\arrow["{\gamma'}", shorten <=8pt, Rightarrow, from=3, to=2-6]
	\arrow["{\Psi_A}"'{pos=0.3}, shorten <=4pt, shorten >=9pt, Rightarrow, from=7p, to=6p]
\end{tikzcd}\]

This time we precompose both sides with the 2-cell:
\[
U(-) \circ \gamma^{-1} y_{UB} \circ U\epsilon_{FB} \mbbd_{y_{UFB}y_B f} \circ \Phi_{FB} y_B f.
\]
We obtain:
\begin{align*}
&U(LHS)y_A \circ \gamma^{-1} y_{UB} \circ U\epsilon_{FB} \mbbd_{y_{UFB}y_B f} \circ \Phi_{FB} y_B f =
\phantom{.}
\phantom{.}
\phantom{.}\\
&= U(Ff \Psi_A \circ \epsilon_{Ff} Fy_A)y_A \circ \gamma^{-1}y_A \circ U\epsilon_{FB} U\gamma'y_A \circ U\epsilon_{FB} UF\mbbd_{y_Bf}y_A \circ \\
&\phantom{=} \circ U\epsilon_{FB} \mbbd_{y_{UFB} y_B f} \circ \Phi_{FB} y_B f
\phantom{.}
\phantom{.}
\phantom{.}\\
&\sstackrel{\eqref{EQP_Falpha}}{=} U(Ff \Psi_A \circ \epsilon_{Ff} Fy_A)y_A \circ \gamma^{-1}y_A \circ U\epsilon_{FB} U\gamma'y_A \circ U\epsilon_{FB}\mbbd_{UFf y_A} \circ \\
&\phantom{=} \circ U\epsilon_{FB} y_{UFB} \mbbd_{y_B f} \circ \Phi_{FB} y_B f
\phantom{.}
\phantom{.}
\phantom{.}\\
&\sstackrel{\eqref{EQP_gamma_prime}}{=} U(Ff \Psi_A \circ \epsilon_{Ff} Fy_A)y_A \circ \gamma^{-1}y_A \circ U\epsilon_{FB} \gamma^{-1} y_A \circ U\epsilon_{FB} UFUFf \mbbd_{y_{UFA}y_A} \circ\\
&\phantom{=}\circ U\epsilon_{FB} UFUFf \mbbd_{y_{UFB} UFf} y_A \circ U\epsilon_{FB} y_{UFB} \mbbd_{y_B f} \circ \Phi_{FB} y_B f
\phantom{.}
\phantom{.}
\phantom{.}\\
&= U(Ff \Psi_A)y_A \circ \gamma^{-1}y_A \circ U\epsilon_{Ff} UFy_A y_A \circ U(\epsilon_{FB} FUFf)\mbbd_{y_{UFA}y_A} \circ \\
&\phantom{=}\circ \gamma^{-1} y_{UFA} y_A \circ U\epsilon_{FB} \mbbd_{y_{UFB} UFf} y_A \circ \Phi_{FB} UFf y_A \circ \mbbd_{y_B f}
\phantom{.}
\phantom{.}
\phantom{.}\\
&\sstackrel{\eqref{EQP_Phi_epsilon}}{=} U(Ff \Psi_A)y_A \circ \gamma^{-1}y_A \circ \gamma^{-1} UFy_A y_A \circ UFf U\epsilon_{FA} \mbbd_{y_{UFA} y_A} \circ \\
&\phantom{=} \circ UFf \Phi_A y_A \circ \mbbd_{y_B f}
\phantom{.}
\phantom{.}
\phantom{.}\\
&= \gamma^{-1}y_A \circ UFf U\Psi_A y_A \circ UFf \gamma^{-1} y_A \circ UFf U\epsilon_{FA} \mbbd_{y_{UFA} y_A} \circ \\
&\phantom{=}\circ UFf \Phi_A y_A \circ \mbbd_{y_B f}
\phantom{.}
\phantom{.}
\phantom{.}\\
&\sstackrel{\eqref{EQP_Psi}}{=} \gamma^{-1} y_A \circ UFf \iota^{-1} y_A \circ \mbbd_{y_B f}
\phantom{.}
\phantom{.}
\phantom{.}\\
&\sstackrel{\eqref{EQP_Psi}}{=} \gamma^{-1} y_A \circ U\Psi_B UFf y_A \circ U\epsilon_{FB} \gamma^{-1} y_A \circ U\epsilon_{FB} UFy_B \mbbd_{y_B f}  \circ \\
&\phantom{=} \circ U\epsilon_{FB} \mbbd_{y_{UFB} y_B} f \circ \Phi_{FB} y_B f
\phantom{.}
\phantom{.}
\phantom{.}\\
&= U(\Psi_B Ff) y_A \circ \gamma^{-1}y_A \circ U\epsilon_{FB} \gamma^{-1} y_A \circ U\epsilon_{FB} UFy_B \mbbd_{y_B f} \circ \\
&\phantom{=}\circ U\epsilon_{FB} \mbbd_{y_{UFB} y_B} f \circ \Phi_{FB} y_B f
\phantom{.}
\phantom{.}
\phantom{.}\\
&\sstackrel{\eqref{EQP_gamma_prime}}{=} U(\Psi_B Ff) y_A \circ \gamma^{-1}y_A \circ U\epsilon_{FB} U\gamma' y_A \circ U\epsilon_{FB} \mbbd_{y_{UFB} y_B f} \circ \Phi_{FB} y_B f
\phantom{.}
\phantom{.}
\phantom{.}\\
&= U(RHS)y_A \circ \gamma^{-1} y_{UB} \circ U\epsilon_{FB} \mbbd_{y_{UFB} y_B f} \circ \Phi_{FB} y_B f.
\end{align*}

\textbf{The second swallowtail identity}: This amounts to showing the following equality, which we will again do by an appropriate pre-composition:

\[\begin{tikzcd}
	FUB \\
	&&& {} \\
	{} & {} & FUFU && FUB & {=} & {\epsilon_B \iota} \\
	\\
	&& FUB && B
	\arrow[""{name=0, anchor=center, inner sep=0}, "{Fy_{UB}}"{description}, from=1-1, to=3-3]
	\arrow["{F1_{UB}}"', shift right=5, curve={height=24pt}, from=1-1, to=5-3]
	\arrow["{F(U\epsilon_B y_{UB})}"{description, pos=0.7}, from=1-1, to=5-3]
	\arrow["{(F\Phi)_B}"'{pos=0.8}, shorten <=10pt, Rightarrow, from=3-1, to=3-2]
	\arrow["{\gamma'}", shorten <=4pt, Rightarrow, from=3-2, to=3-3]
	\arrow["{\epsilon_{FUB}}", from=3-3, to=3-5]
	\arrow["{FU \epsilon_B}"{description}, from=3-3, to=5-3]
	\arrow[curve={height=24pt}, Rightarrow, no head, from=3-5, to=1-1]
	\arrow["{\epsilon_B}", from=3-5, to=5-5]
	\arrow["{\epsilon_{\epsilon_B}}"{description}, Rightarrow, from=5-3, to=3-5]
	\arrow["{\epsilon_B}"', from=5-3, to=5-5]
	\arrow["{\Psi_{UB}}"'{pos=0.6}, shift right=3, shorten <=25pt, shorten >=25pt, Rightarrow, from=0, to=2-4]
\end{tikzcd}\]

\begin{align*}
&U(LHS)y_{UB} \circ \gamma^{-1} y_{UB} \circ U\epsilon_B \mbbd_{y_{UB}} \circ \Phi_B = \\
\phantom{.}
\phantom{.}
\phantom{.}
&= U(\epsilon_B \Psi_{UB} \circ \epsilon_{\epsilon_B} Fy_{UB})y_{UB} \circ \gamma^{-1}y_{UB} \circ U\epsilon_B U\gamma' y_{UB} \circ U\epsilon_B UF\Phi_B y_{UB} \circ \\
&\phantom{=}\circ U\epsilon_B \mbbd_{y_{UB}} \circ \Phi_B \\
\phantom{.}
\phantom{.}
\phantom{.}
&\sstackrel{\eqref{EQP_Falpha}}{=} U(\epsilon_B \Psi_{UB} \circ \epsilon_{\epsilon_B}Fy_{UB})y_{UB} \circ \gamma^{-1}y_{UB} \circ U\epsilon_B U\gamma' y_{UB} \circ U\epsilon_B \mbbd_{y_{UB} U\epsilon_B y_{UB}} \circ \\
&\phantom{=}\circ \Phi_B U\epsilon_B y_{UB} \circ \Phi_B\\
\phantom{.}
\phantom{.}
\phantom{.}
&\sstackrel{\eqref{EQP_gamma_prime}}{=} U(\epsilon_B \Psi_{UB} \circ \epsilon_{\epsilon_B}Fy_{UB})y_{UB} \circ \gamma^{-1}y_{UB} \circ U\epsilon_B \gamma^{-1}y_{UB} \circ \\
&\phantom{=}\circ U\epsilon_B UFU\epsilon_B \mbbd_{y_{UFUB}y_{UB}} \circ U\epsilon_B \mbbd_{y_{UB}U\epsilon_B} y_{UB} \circ \Phi_B U\epsilon_B y_{UB} \circ \Phi_B \\
\phantom{.}
\phantom{.}
\phantom{X}
&= U(\epsilon_B \Psi_{UB}) y_{UB} \circ \gamma^{-1} y_{UB} \circ U\epsilon_{\epsilon_B} UFy_{UB} y_{UB} \circ \gamma^{-1} UFy_{UB} y_{UB} \circ \\
&\phantom{=} \circ U\epsilon_B UFU\epsilon_B \mbbd_{y_{UFUB}y_{UB}} \circ U\epsilon_B \mbbd_{y_{UB}\epsilon_B} y_{UB} \circ  U\epsilon_B y_{UB} \Phi_B \circ \Phi_B\\
\phantom{.}
\phantom{.}
\phantom{.}
&\sstackrel{\eqref{EQP_Phi_epsilon}}{=} U(\epsilon_B \Psi_{UB})y_{UB} \circ \gamma^{-1} y_{UB} \circ \gamma^{-1} UFy_{UB} y_{UB} \circ U\epsilon_B U\epsilon_{FUB} \mbbd \circ \\
&\phantom{=}\circ U\epsilon_B \Phi_B y_{UB} \circ \Phi_B      \\
\phantom{.}
\phantom{.}
\phantom{Y}
&= \gamma^{-1}y_{UB} \circ U\epsilon_B U\Psi_B y_{UB} \circ U\epsilon_B U\epsilon_{FUB} \mbbd_{y_{UFUB}y_{UB}} \circ U\epsilon_B \Phi_B y_{UB} \circ \Phi_B      \\
\phantom{.}
\phantom{.}
\phantom{.}
&\sstackrel{\eqref{EQP_Psi}}{=} \gamma^{-1}y_{UB} \circ U\epsilon_B \iota^{-1}y_{UB} \circ \Phi_B \\
\phantom{.}
\phantom{.}
\phantom{.}
&\sstackrel{\eqref{EQP_iota_prime}}{=} \gamma^{-1}y_{UB} \circ U\epsilon_B U\iota' y_{UB} \circ U\epsilon_B \mbbd_{y_{UB}} \circ \Phi_B\\
\phantom{.}
\phantom{.}
\phantom{.}
&= U(RHS) y_{UB} \circ \gamma^{-1}y_{UB} \circ U\epsilon_B \mbbd_{y_{UB}} \circ \Phi_B.
\end{align*}

\end{proof}

\begin{lemma}\label{THML_big_l_adj_thm}
The composite bijection \textbf{(A)+(B)+(C)} in the proof of Theorem \ref{THM_BIG_lax_adj_thm}:
\[ \ck_D(DA,DHL)(f,Dl') \cong \cl(GDA,L)(s_L \circ Gf, s_L \circ GDl'), \]
is given by the assignment:
\[
\alpha \mapsto s_LG\alpha.
\]
\end{lemma}
\begin{proof}

Because of the swallowtail identity for the biadjunction in Proposition \ref{THMP_biadjunkce_ck_ck_D}, it can be seen that the counit of the adjunction \eqref{EQ_hL_yA_yHL_JD} evaluated at $f: A \to HL$ is equal to the following (using notation from Remark \ref{POZN_explicit_form_h_L}):
\[\begin{tikzcd}
	A && DA && DHL && HL \\
	&& HL && {}
	\arrow["{y_A}", from=1-1, to=1-3]
	\arrow["Df", from=1-3, to=1-5]
	\arrow["{h_L}", from=1-5, to=1-7]
	\arrow["f"', from=1-1, to=2-3]
	\arrow["{y_{HL}}"{description}, from=2-3, to=1-5]
	\arrow[curve={height=18pt}, Rightarrow, no head, from=2-3, to=1-7]
	\arrow["{y^{-1}_f}", Rightarrow, from=1-3, to=2-3]
	\arrow["{\epsilon_L}"{pos=0.3}, shorten >=6pt, Rightarrow, from=1-5, to=2-5]
\end{tikzcd}\]
The composite bijection \textbf{(A)+(B)+(C)} is thus the assignment:
\[
\alpha \mapsto s_L \circ GD(\epsilon_L f \circ h_L y^{-1}_{Dl'} \circ h_L \alpha y_A) \circ \beth_f^{-1}.
\]

Unwrapping the definitions of variables $h_L, \epsilon_L, \beth_f^{-1}$, what we need to show that the composite 2-cell below equals $s_L G\alpha$ (note that we use the same convention for the modifications on which a pseudofunctor is applied as in Definition \ref{DEFI_colax_adjunction}):

\[
\adjustbox{scale=0.85,center}{
\begin{tikzcd}[column sep=small]
	&&&& GDA \\
	&&&& {GD^2HL} && GDHL & L \\
	\\
	GDA & {GD^2A} && {GD^2HL} & {GDHGD^2HL} && GDHGDHL & GDHL \\
	\\
	&& GDHL & GDHGDHL
	\arrow[""{name=0, anchor=center, inner sep=0}, Rightarrow, no head, from=4-4, to=2-5]
	\arrow[""{name=0p, anchor=center, inner sep=0}, phantom, from=4-4, to=2-5, start anchor=center, end anchor=center]
	\arrow[""{name=1, anchor=center, inner sep=0}, "{GDc_{DHL}}"', from=4-4, to=4-5]
	\arrow[""{name=1p, anchor=center, inner sep=0}, phantom, from=4-4, to=4-5, start anchor=center, end anchor=center]
	\arrow["{s_{GD^2HL}}"{description}, from=4-5, to=2-5]
	\arrow[""{name=2, anchor=center, inner sep=0}, "{Gp_{DHL}}"{description}, from=2-5, to=2-7]
	\arrow[""{name=2p, anchor=center, inner sep=0}, phantom, from=2-5, to=2-7, start anchor=center, end anchor=center]
	\arrow[""{name=3, anchor=center, inner sep=0}, "{s_L}", from=2-7, to=2-8]
	\arrow[""{name=3p, anchor=center, inner sep=0}, phantom, from=2-7, to=2-8, start anchor=center, end anchor=center]
	\arrow[""{name=4, anchor=center, inner sep=0}, from=4-5, to=4-7]
	\arrow[""{name=4p, anchor=center, inner sep=0}, phantom, from=4-5, to=4-7, start anchor=center, end anchor=center]
	\arrow["{s_{GDHL}}"{description}, from=4-7, to=2-7]
	\arrow[""{name=5, anchor=center, inner sep=0}, "{GDHs_L}"', from=4-7, to=4-8]
	\arrow[""{name=5p, anchor=center, inner sep=0}, phantom, from=4-7, to=4-8, start anchor=center, end anchor=center]
	\arrow["{s_L}"', from=4-8, to=2-8]
	\arrow["Gf", from=1-5, to=2-7]
	\arrow[""{name=6, anchor=center, inner sep=0}, "{Gp_{DA}}"{description, pos=0.6}, from=4-2, to=1-5]
	\arrow[""{name=7, anchor=center, inner sep=0}, "GDf", curve={height=-12pt}, from=4-2, to=4-4]
	\arrow["{GDy_A}"', from=4-1, to=4-2]
	\arrow[""{name=8, anchor=center, inner sep=0}, curve={height=-24pt}, Rightarrow, no head, from=4-1, to=1-5]
	\arrow[""{name=9, anchor=center, inner sep=0}, "{GD^2l'}"', curve={height=12pt}, from=4-2, to=4-4]
	\arrow["{GDl'}"', curve={height=18pt}, from=4-1, to=6-3]
	\arrow["{GDy_{HL}}"{description, pos=0.7}, from=6-3, to=4-4]
	\arrow[from=6-3, to=6-4]
	\arrow["{GDHGDy_{HL}}"{description}, from=6-4, to=4-5]
	\arrow[""{name=10, anchor=center, inner sep=0}, curve={height=18pt}, Rightarrow, no head, from=6-4, to=4-7]
	\arrow[""{name=11, anchor=center, inner sep=0}, shift right=5, curve={height=30pt}, Rightarrow, no head, from=6-3, to=4-8]
	\arrow["{s_{Gp_{DHL}}^{-1}}", shorten <=9pt, shorten >=9pt, Rightarrow, from=2p, to=4p]
	\arrow["{s_{s_L}^{-1}}", shorten <=9pt, shorten >=9pt, Rightarrow, from=3p, to=5p]
	\arrow["GD\alpha", shorten <=3pt, shorten >=3pt, Rightarrow, from=7, to=9]
	\arrow["{(GDy)_{l'}}"'{pos=0.7}, shorten <=12pt, Rightarrow, from=9, to=6-3]
	\arrow["{(GDc)_{y_{HL}}}"'{pos=1}, shift right=3, shorten <=21pt, Rightarrow, from=1, to=6-4]
	\arrow["{(GDHG\Psi)_{HL}}", shorten >=7pt, Rightarrow, from=4-5, to=10]
	\arrow["{(GD\tau)^{-1}_L}", shorten <=3pt, shorten >=3pt, Rightarrow, from=10, to=11]
	\arrow["{(G\Psi)^{-1}_A}"', shorten <=6pt, shorten >=6pt, Rightarrow, from=8, to=6]
	\arrow["{\sigma_{DHL}^{-1}}", shorten <=4pt, shorten >=4pt, Rightarrow, from=0p, to=1p]
	\arrow["{(Gp)_f^{-1}}", shorten <=6pt, shorten >=9pt, Rightarrow, from=6, to=4-4]
\end{tikzcd}
}\]

The diagram below proves this equality:
\[
\adjustbox{scale=0.8,center}{
\begin{tikzcd}[column sep=tiny]
	{s_L Gf} & {s_L Gf Gp_{DA}GDy_A} & {s_L Gp_{DHL} GDf GDy_A} \\
	{s_L GDl'} & {s_L GDl' Gp_{DA}GDy_A} & {s_L Gp_{DHL} GD^2l' GDy_A} \\
	{s_L GDl'} && {s_L Gp_{DHL} GDy_{HL} GDl'} \\
	& {s_L Gp_{DHL} GDy_{HL}s_{GDHL} GDc_{HL} GDl'} & {s_L Gp_{DHL} s_{GD^2HL} GDc_{DHL} GDy_{HL} GDl'} \\
	{s_L GDl'} && {s_L Gp_{DHL} s_{GD^2HL}GDHGDy_{HL} GDc_{HL} GDl'} \\
	{s_L GDHs_L GDc_{HL} GDl'} & {s_L s_{GDHL} GDc_{HL} GDl'} & {s_L s_{GDHL} GDHGp_{DHL} GDHGDy_{HL} GDc_{HL} GDl'}
	\arrow[""{name=0, anchor=center, inner sep=0}, "{s_L Gf (G\Psi)^{-1}_A}", from=1-1, to=1-2]
	\arrow[""{name=0p, anchor=center, inner sep=0}, phantom, from=1-1, to=1-2, start anchor=center, end anchor=center]
	\arrow["{s_L G\alpha}"', from=1-1, to=2-1]
	\arrow["{s_L(Gp)_f^{-1}GDy_A}", from=1-2, to=1-3]
	\arrow["{s_LG\alpha Gp_{DA} GDy_A}", from=1-2, to=2-2]
	\arrow["{s_L Gp_{DHL} GD\alpha GDy_A}", from=1-3, to=2-3]
	\arrow[""{name=1, anchor=center, inner sep=0}, "{s_L GDl'(G\Psi)^{-1}_A}", from=2-1, to=2-2]
	\arrow[""{name=1p, anchor=center, inner sep=0}, phantom, from=2-1, to=2-2, start anchor=center, end anchor=center]
	\arrow["1"', from=2-1, to=3-1]
	\arrow[""{name=2, anchor=center, inner sep=0}, "{s_L (G\Psi)^{-1}_{HL}GDl'}"{description}, from=2-1, to=3-3]
	\arrow[""{name=3, anchor=center, inner sep=0}, "{s_L (Gp)^{-1}_{Dl}GDy_A}"{description}, from=2-2, to=2-3]
	\arrow["{s_L Gp_{DHL} GDy_{l'}^{-1}}", from=2-3, to=3-3]
	\arrow["1"', from=3-1, to=5-1]
	\arrow["{s_L \sigma^{-1}_{HL} GDl'}"{description}, from=3-1, to=6-2]
	\arrow[""{name=4, anchor=center, inner sep=0}, "{s_L (G\Psi)_{HL} GDl'}"{description}, from=3-3, to=3-1]
	\arrow[""{name=5, anchor=center, inner sep=0}, "{s_L Gp_{DHL}GDy_{HL} \sigma^{-1}_{HL} GDl'}"{description}, from=3-3, to=4-2]
	\arrow["{s_L Gp_{DHL} \sigma^{-1}_{DHL} GDy_{HL} GDl'}", from=3-3, to=4-3]
	\arrow["{s_L (G\Psi)_{HL} s_{GDHL} GDc_{HL} GDl'}", from=4-2, to=6-2]
	\arrow["{s_L Gp_{DHL} s_{GD^2HL} GDc_{y_{HL}} GDl'}", from=4-3, to=5-3]
	\arrow[""{name=6, anchor=center, inner sep=0}, "{s_L Gp_{DHL} s_{GDy_{HL}} GDc_{HL} GDl'}"{description}, from=5-3, to=4-2]
	\arrow["{s_L s^{-1}_{Gp_{DHL}} GDHGDy_{HL} GDc_{HL} GDl'}", from=5-3, to=6-3]
	\arrow["{s_L (GD\tau)^{-1}_L GDl'}"{description}, from=6-1, to=5-1]
	\arrow[""{name=7, anchor=center, inner sep=0}, "{s_L \sigma_{HL} GDl'}"{description}, from=6-2, to=5-1]
	\arrow[""{name=8, anchor=center, inner sep=0}, "{s^{-1}_{s_L} GDc_{HL} GDl'}", from=6-2, to=6-1]
	\arrow[""{name=9, anchor=center, inner sep=0}, "{s_L s_{GDHL} (GDHG\Psi)_{HL}GDc_{HL} GDl'}", shift left=3, from=6-3, to=6-2]
	\arrow["{(*)}"{description}, draw=none, from=0p, to=1p]
	\arrow["{(a)}"{description}, draw=none, from=1-3, to=3]
	\arrow["{(b)}"{description}, draw=none, from=3, to=2]
	\arrow["{(*)}"{description}, draw=none, from=4, to=4-2]
	\arrow["{(e)}"{description}, draw=none, from=5, to=6]
	\arrow["{(d)}"{description, pos=0.8}, draw=none, from=6, to=9]
	\arrow["{(c)}"{description}, draw=none, from=7, to=8]
\end{tikzcd}}\]

In this diagram:
\begin{itemize}
\item $(a)$ is the local naturality of $(Gp)^{-1}$,
\item $(b)$ is the modification axiom for $(G\Psi)^{-1}$
\item $(c)$ is the swallowtail identity for $(s,c)$,
\item $(d)$ is the equation derived from the local naturality of $s$,
\item $(e)$ is the modification axiom for $\sigma^{-1}$,
\item $(*)$'s are the middle-four-interchange laws.
\end{itemize}

\end{proof}